\documentclass[sn-mathphys,Numbered]{sn-jnl}

\usepackage{graphicx}
\usepackage{multirow}
\usepackage{amsmath,amssymb,amsfonts}
\usepackage{amsthm}
\usepackage{mathrsfs}
\usepackage[title]{appendix}
\usepackage{xcolor}
\usepackage{textcomp}
\usepackage{manyfoot}
\usepackage{booktabs}
\usepackage{algorithm}
\usepackage{algorithmicx}
\usepackage{algpseudocode}
\usepackage{listings}

\usepackage{subcaption}
\usepackage{url}
\usepackage{adjustbox}
\usepackage{morefloats}

\graphicspath{{figs/}}

\theoremstyle{thmstyleone}

\theoremstyle{thmstyletwo}

\theoremstyle{thmstylethree}

\newcommand{\corrtext}[1]{#1}

\begin{document}

\title[Improving local solution of the ADER-DG method for ODE and DAE systems]{%
	Improving the local solution of the DG predictor of the ADER-DG method for solving
	systems of ordinary differential equations and its applicability to systems of differential-algebraic equations
}

\author*[1]{\fnm{Ivan S.} \sur{Popov}}
\email{diphosgen@mail.ru, popovis@omsu.ru}
\affil*[1]{%
	\orgdiv{Department of Theoretical Physics}, \orgname{Dostoevsky Omsk State University}, 
	\orgaddress{\street{Mira prospekt}, \city{Omsk}, \postcode{644077}, \country{Russia}}
}

\abstract{%
Improved local numerical solution for the ADER-DG numerical method with a local DG predictor for solving the initial value problem for a first-order ODE system is proposed. The improved local numerical solution demonstrates convergence orders of one higher than the convergence order of the local numerical solution of the original ADER-DG numerical method and has the property of continuity at grid nodes. Rigorous proofs of the approximation orders of the local numerical solution and the improved local numerical solution are presented. Obtaining the proposed improved local numerical solution does not require significant changes to the structure of the ADER-DG numerical method. Therefore, all conclusions regarding the convergence orders of the numerical solution at grid nodes, the resulting superconvergence, and the high stability of the ADER-DG numerical method remain unchanged. A wide range of applications of the ADER-DG numerical method is presented for solving specific initial value problems for ODE systems for a wide range of polynomial degrees. The obtained results provide strong confirmation for the developed rigorous theory. The improved local numerical solution is shown to exhibit both higher accuracy and improved smoothness and point-wise comparability. Empirical convergence orders of all individual numerical solutions were calculated for a wide range of error norms, which well agree with the expected convergence orders. The rigorous proof, based on the $\epsilon$-embedding method, of the applicability of the ADER-DG numerical method with a local DG predictor to solving DAE systems is presents.
}



\maketitle

\section*{Introduction}
\addcontentsline{toc}{section}{Introduction}
\label{intro}

In this paper, a study of the arbitrary high order ADER (Arbitrary DERivative) discontinuous Galerkin (DG) method with local DG predictor, which is frequently used for solving problems for partial differential equations, is based on solving the initial value problems for the first-order nonlinear ordinary differential equation (ODE) system chosen in the following form
\begin{equation}\label{eq:ivp_ode_diff_src}
\frac{d\mathbf{u}}{dt} = \mathbf{F}\left(\mathbf{u},\, t\right),\quad  \mathbf{u}(t_{0}) = \mathbf{u}_{0},\quad t\in\Omega = \left\{t\, |\, t \in [t_{0},\ t_{f}]\right\},
\end{equation}
where $\mathbf{u}: \Omega \rightarrow \mathbb{R}^{D}$ is a desired vector-function; the vector-function $\mathbf{F}: \mathbb{R}^{D}\times\Omega \rightarrow \mathbb{R}^{D}$ is a right-side function, which is given; $D$ is the total number of desired functions; $\Omega$ is the domain of definition. The initial condition $\mathbf{u}_{0}$ on the desired vector-function $\mathbf{u}$ was chosen at the point $t_{0}$. The classical theory of ordinary differential equation shows that in case $\mathbf{F} \in C_{1}(\mathbb{R}^{K}\times\Omega)$ the solution of the problem exists and is unique. It is clear that the initial value problem for systems of differential equations containing differential equations of the second and higher orders, uniquely solvable with respect to derivatives of higher orders, can be represented in the chosen form of the first-order nonlinear ordinary differential equations system (\ref{eq:ivp_ode_diff_src}).

Initial value problems for ODE systems and the need to obtain numerical solutions arise in a wide range of scientific and technical problems. Currently, there exist many numerical methods for solving initial value problems for ODE systems~\cite{Butcher_book_2016, Hairer_book_1, Hairer_book_2}. In recent years, significant interest in this area has been associated with the development and study of DG methods for the numerical solution of ODE systems~\cite{Babuska_book_2001, Wahlbin_lectures_1995, dg_ivp_ode_4, dg_ivp_ode_5, dg_ivp_ode_6}. This interest is associated with the possibility of obtaining a highly accurate numerical solution, the occurrence of superconvergence of DG methods when solving ODE systems, and the possibility of obtaining a continuous numerical solution to the problem, not just at a finite set of grid nodes~\cite{dg_ivp_ode_4, dg_ivp_ode_5, dg_ivp_ode_6}. DG methods, in general, occupy a special place in the field of numerical solution of differential equations, due to the high accuracy of the numerical solution obtained using them and their scalability to complex and large problems.

DG methods were initially proposed by Reed and Hill~\cite{lasl_rep_dg_1973} for the numerical simulation of the neutron transport equation. In 1981 Delfour \textit{et al} \cite{Delfour_1981} constructed a Runge-Kutta-like DG method for solving the initial value problem for the ODE system that demonstrates superconvergence up to order $2N+2$ for polynomials with a degree $N$. Later, in 1986 Delfour \textit{et al} \cite{Delfour_1986} constructed a general theory of one-step, hybrid, and multi-step methods using discontinuous polynomial approximations for solving the initial value problem for the ODE system. Cockburn, Shu \textit{et al} in~\cite{Cockburn_base_1, Cockburn_base_2, Cockburn_base_3, Cockburn_base_4, Cockburn_base_5} created an accurate, rigorous and thoroughly developed mathematical basis of DG methods, which stimulated their further development and use for solving a wide class of problems. DG methods are used to solve the initial value problem for ODE systems~\cite{ader_dg_ode_jsc, dg_ivp_ode_5, dg_ivp_ode_6, dg_ivp_ode_1, dg_ivp_ode_2, dg_ivp_ode_3}, the boundary value problem for ODE systems~\cite{dg_ivp_ode_4, dg_bvp_ode_1, dg_bvp_ode_2}, to solve systems of differential-algebraic equations (DAE)~\cite{ader_dg_dae}, stochastic differential equations~\cite{dg_stoch_ode_1, dg_stoch_ode_2}, and especially widely to solve partial differential equations systems~\cite{ader_dg_ideal_flows, ader_dg_diss_flows, ader_dg_ale, ader_dg_grmhd, ader_dg_gr_prd, ader_dg_gr_z4_2024, ader_dg_simple_mod, ader_dg_PNPM, fron_phys, ader_dg_axioms, exahype, PNPM_DG_2009, ader_dg_semiexpl, ader_dg_hyperelastic, ader_dg_seiemic, ader_dg_seiemic_underwater, ader_eno_fv_blood_2022, ader_dg_gr_z4_2024, ader_dg_wb_shwater_2022, ader_eff, ader_eff_blas, dg_entropy, dg_entropy_add, ader_weno_sph, ader_dg_mod_1, ader_dg_mod_2, ader_rev_2024}.

The ADER-DG methods are DG methods based on the ADER paradigm proposed by Titarev and Toro in~\cite{ader_init_1, ader_init_2} in the context of finite volume methods for partial differential equations systems. The ADER paradigm allows the creation of numerical methods of an arbitrarily high order. Among the existing numerical methods for solving ODE systems, it can be considered as somewhat similar to high-order numerical method based on the Taylor expansion for solving ODE systems~\cite{ivp_ode_taylor_series_2017, ivp_ode_taylor_series_soft_2005}, which allow obtaining an arbitrary high order. The modern version of ADER involves the use of a local DG predictor, which was proposed by Dumbser \textit{et al}~\cite{ader_stiff_1, ader_stiff_2}. The use of ADER in DG methods has allowed obtaining unrivaled results in the numerical solution of systems of partial differential equations~\cite{ader_dg_ideal_flows, ader_dg_diss_flows, ader_dg_ale, ader_dg_grmhd, ader_dg_gr_prd, ader_dg_gr_z4_2024, ader_dg_simple_mod, ader_dg_PNPM, fron_phys, ader_dg_axioms, exahype, PNPM_DG_2009, ader_dg_semiexpl, ader_dg_hyperelastic, ader_dg_seiemic, ader_dg_seiemic_underwater, ader_eno_fv_blood_2022, ader_dg_gr_z4_2024, ader_dg_wb_shwater_2022, ader_eff, ader_eff_blas, dg_entropy, dg_entropy_add, ader_weno_sph, ader_dg_mod_1, ader_dg_mod_2, ader_rev_2024, ader_stiff_3, ader_stiff_4}. Han Veiga \textit{et al} in~\cite{dec_vs_ader_2021} and Micalizzi \textit{et al} in~\cite{dec_vs_ader_2023} showed that the methods of the ADER family are significantly interconnected with numerical methods on the deferred correction (DeC) paradigm, which has a long history and its application to initial value problems for ODE systems goes back to Daniel \textit{et al}~\cite{dec_src_1968}, and are effectively used to solve both ODE systems~\cite{dec_dutt_2000, dec_dutt_2000, dec_minion_2003, dec_shu_2008} and partial differential equations systems~\cite{dec_abgrall_2017, dec_abgrall_2019}.

It should be emphasized that an important feature of numerical DG methods, in particular the version of the ADER-DG method with a local DG predictor presented in this paper, is the arbitrarily high convergence order and superconvergence, reaching convergence of order $2N+1$ when using polynomials of degree $N$ in the numerical solution representation. Guaranteed achievement of such high superconvergence is quite difficult and challenging task, especially for stiff problems~\cite{ader_dg_ode_jsc} and significantly nonuniform grids~\cite{siac_rev, siac_ref_1, siac_ref_3}. To solve such problems, post-processing methods for solutions that allow achieving high superconvergence are currently being studied and are being developed. In particular, a relevant approach to increasing the order convergence is the Smoothness-Increasing Accuracy-Conserving (SIAC)~\cite{siac_ref_5} (see details also in~\cite{siac_ref_1, siac_ref_3, siac_ref_5, siac_ref_11, siac_ref_20, siac_ref_12, siac_ref_15, siac_ref_16}).

An important feature of DG methods, especially for solving ODE systems~\cite{ader_dg_ode_jsc, dg_ivp_ode_4, dg_ivp_ode_5, dg_ivp_ode_6, dg_ivp_ode_1, dg_ivp_ode_2, dg_ivp_ode_3}, is the possibility of obtaining a continuous numerical solution to the problem in the domain of definition $\Omega$, and not only at a finite set of grid nodes. In~\cite{ader_dg_ode_jsc} a version of the ADER-DG numerical method with a local DG predictor was presented for solving the initial value problem for a first-order ODE system (\ref{eq:ivp_ode_diff_src}). It was shown that the ADER-DG numerical method with a local DG predictor has a convergence order $2N+1$ for the numerical solution at grid nodes and $N+1$ for the local solution defined over the domain of definition $\Omega$, for polynomial degree $N$. Furthermore, the ADER-DG numerical method with a local DG predictor is a highly stable method, characterized by $A$-stability and $L_{1}$-stability. \corrtext{The theory of ADER-DG methods for solving ODE systems was developed in~\cite{ader_improving_2024, ader_proofs_2025, ader_dg_ode_sinum}.} The work~\cite{ader_dg_dae} presents a version of the ADER-DG numerical method with a local DG predictor for solving a DAE system, based on a modification of the method~\cite{ader_dg_ode_jsc}.

This paper develops a method for improving the local numerical solution of the ADER-DG method with a local DG predictor, defined over the domain of definition $\Omega$, with a higher convergence order $N+2$ and higher accuracy compared to the original local solution of the ADER-DG numerical method with a local DG predictor. The key properties of the original ADER-DG numerical method, such as superconvergence $2N+1$ of the numerical solution at grid nodes and high stability, remain unchanged. A rigorous proof of the approximation orders $N+1$ for the original local solution and $N+2$ for the improved local solution is presented. It is important to note for comparison that the DG methods for solving ODE systems typically have the approximation order and the convergence order $N+1$~\cite{dg_ivp_ode_5} for a continuous solution over the domain of definition $\Omega$, while the improved local numerical solution for the ADER-DG method with a local DG predictor proposed in this paper has the approximation order and the convergence order one unit higher --- $N+2$, which, in particular, was achieved for other DG methods in~\cite{dg_ivp_ode_3}. This paper also presents a rigorous proof, within the framework of the $\epsilon$-embedding method~\cite{Hairer_book_2}, of the applicability of the ADER-DG numerical method with a local DG predictor for solving the initial value problem for a first-order ODE system (\ref{eq:ivp_ode_diff_src}) to solving a DAE system.

Section~\ref{sec:gen_descr} ``General description of the numerical method'' presents a general description of the ADER-DG numerical method with a local DG predictor for solving the initial value problem for a first-order ODE system. Section~\ref{sec:imp_dg} ``Improving the local solution'' presents an improved local solution for the ADER-DG numerical method with a local DG predictor and provides rigorous proofs of the approximation orders $N+1$ for the original local solution (in Subsection~\ref{sec:imp_dg:loc_sol_prop} ``Properties of the local solution'') and $N+2$ for the improved local solution (in Subsection~\ref{sec:imp_dg:imp_loc_sol} ``Improved local solution''). Section~\ref{sec:apps} ``Applications of the numerical method'' presents a wide range of applications of the ADER-DG numerical method with a local DG predictor for solving specific initial value problems for an ODE system for a wide range of polynomial degrees, which covers both the range of sufficiently low polynomial degrees required in applied calculations and the range of very high polynomial degrees that can be used in high-precision calculations. In the final Section~\ref{sec:daes} ``Features of the numerical method for DAE systems'', a rigorous proof is presented, based on the use of the $\epsilon$-embedding method~\cite{Hairer_book_2}, of the applicability of the numerical ADER-DG method with a local DG predictor for solving DAE systems, and it is also revealed that the version of the method proposed in the work can be attributed to the classical version of the ADER-DG method with a local DG predictor for solving the initial value problem for a first-order ODE system.

\section{General description of the numerical method}
\label{sec:gen_descr}

The general structure and implementation features of the ADER-DG numerical method with a local time DG predictor are well described in the works~\cite{ader_dg_ideal_flows, ader_dg_diss_flows, ader_dg_ale, ader_dg_grmhd, ader_dg_gr_prd, ader_dg_gr_z4_2024, ader_dg_simple_mod, ader_dg_PNPM, fron_phys, ader_dg_axioms, exahype, PNPM_DG_2009}. The application of the ADER-DG numerical method to solving the initial value problem for a first-order ODE system is described in detail in~\cite{ader_dg_ode_jsc, ader_dg_dae}, and the features of the internal structure of the method and rigorous proofs are presented in~\cite{ader_dg_ode_sinum}. Therefore, in the following, the presentation of the structure of the ADER-DG numerical method and the local time DG predictor will be limited mainly to the introduction of the concepts necessary for the present work.

The numerical solution of the ODE system (\ref{eq:ivp_ode_diff_src}), obtained by the ADER-DG numerical method, is presented on the discretization of the domain of definition $\Omega$ of the desired functions by a finite number of non-overlapping discretization domains $\Omega_{n} = \{t\, |\, t\in[t_{n},\ t_{n+1}]\}$ that completely cover the domain of definition $\Omega = \cup_{n} \Omega_{n}$, where $t_{n}$ and $t_{n+1} = t_{n} + {\Delta t_{n}}$ are discretization nodes, and $\Delta t_{n}$ is the discretization step, which is not assumed to be constant for different discretization domains $\Omega_{n}$. Therefore, the proposed ADER-DG numerical method is a method that allows implementation for a variable step $\Delta t_{n}$. The set of discretization domains $\Omega_{n}$ covering the domain of definition $\Omega$ of the desired functions is a one-dimensional grid. Further points $t_{n}$ will be denoted by the grid nodes, and the space between nodes $t\in[t_{n},\ t_{n+1}]$ will be denoted by the domain between nodes.

The ADER-DG numerical method for solving the ODE system (\ref{eq:ivp_ode_diff_src}) is based on the use of an integral form on the discretization domains $\Omega_{n}$:
\begin{equation}\label{eq:ivp_ode_int_src}
\mathbf{u}_{n+1} = \mathbf{u}_{n} + \int\limits_{t_{n}}^{t_{n+1}} \mathbf{F}\left(\mathbf{u}(t),\, t\right) dt,
\end{equation}
where $\mathbf{u}_{n}$ and $\mathbf{u}_{n+1}$ are the solution to the problem at the grid nodes $t_{n}$ and $t_{n+1}$, which will hereinafter be called simply the ``solution at the nodes''; the function $\mathbf{u}(t)$ included in the integrand represent the solution to the problem in the space between the nodes $t\in[t_{n},\ t_{n+1}]$. The formula apparatus of the ADER-DG numerical method was formulated in terms of the local variable $\tau$:
\begin{equation}\label{eq:tau_mapping}
t(\tau) = t_{n} + \tau\cdot{\Delta t}_{n}\in\Omega_{n},\quad \tau(t) = \frac{t - t_{n}}{{\Delta t}_{n}}\in\omega,
\end{equation}
using which the discretization domain $\Omega_{n}$ is mapped onto the reference domain $\omega = \{\tau\, |\, \tau\in[0,\ 1]\}$, which is associated with the convenience of determining the basis functions. Using the local variable $\tau$, the integral expression (\ref{eq:ivp_ode_int_src}) is written in the following form:
\begin{equation}\label{eq:ivp_ode_int_src_by_tau}
\mathbf{u}_{n+1} = \mathbf{u}_{n} + {\Delta t}_{n} \int\limits_{0}^{1} \mathbf{F}(\mathbf{u}(t(\tau)),\, t(\tau)) d\tau.
\end{equation}
The solution between the nodes $\mathbf{u}(t(\tau))$ included in the integrand is chosen in the form of a local discrete time solution $\mathbf{q}_{n}: \omega \rightarrow \mathcal{R}^{D}$, which is determined for each individual discretization domain $\Omega_{n}$:
\begin{equation}
\mathbf{u} = \mathbf{u}(t(\tau)),\ t\in\Omega_{n}\ \mapsto\ \mathbf{q}_{n} = \mathbf{q}_{n}(\tau),\ \tau\in\omega.
\end{equation}
The local solution is also defined for the definition domain $\mathbf{u}_{L}: \Omega \rightarrow \mathcal{R}^{D}$ by a piecewise assembly of local solutions $\mathbf{q}_{n}$ for each discretization domain $\Omega_{n}$:
\begin{equation}\label{eq:local_sol_assembly}
\mathbf{u}_{L}(t) = \sum\limits_{n} \chi_{n}(t)\cdot\mathbf{q}_{n}\left(\frac{t - t_{n}}{{\Delta t}_{n}}\right),
\end{equation}
where $\chi_{n}: \Omega\rightarrow\{0,\, 1\}$ is the indicator function of the discretization domain $\Omega_{n}\subseteq\Omega$.

The local discrete time solution $\mathbf{q}_{n}$ in discretization domain $\Omega_{n}$ is represented in the form of an expansion over a set $\{\varphi_{p}\}$ of basis functions $\varphi_{p}: \omega\rightarrow\mathcal{R}$ of the following form:
\begin{equation}\label{eq:local_sol}
\mathbf{q}_{n}(\tau) = \sum\limits_{p} \hat{\mathbf{q}}_{n, p}\varphi_{p}(\tau),
\end{equation}
where $\{\hat{\mathbf{q}}_{n, p}\}$ is the set of expansion coefficients, and is found as a solution to the weak form of the ODE system (\ref{eq:ivp_ode_diff_src}) in $\Omega_{n}$:
\begin{equation}\label{eq:ivp_ode_weak}
\int\limits_{0}^{1}\varphi_{p}(\tau)\left[\frac{d\mathbf{q}_{n}(\tau)}{d\tau} - {\Delta t}_{n}\mathbf{F}(\mathbf{q}_{n}(\tau), t(\tau))\right]d\tau = 0,\quad
\mathbf{q}_{n}(0) = \mathbf{u}_{n},
\end{equation}
where the initial condition is also understood in weak form.

The basis functions are selected in the form of Lagrange interpolation polynomials $\{\varphi_{p}\}_{p = 0}^{N}$ of degree $N$, with nodal points at the roots $\tau_{k}$ of the shifted Legendre polynomials $P_{N+1}: \omega\rightarrow\mathcal{R}$, which are represented in the following form:
\begin{equation}\label{eq:basis_funcs_def}
\varphi_{p}(\tau) = \sum\limits_{k = 0}^{N} \varphi_{p, k}\tau^{k} = \prod\limits_{k \neq p} \frac{\tau - \tau_{k}}{\tau_{p} - \tau_{k}},\quad
\varphi_{p}(\tau_{q}) = \delta_{pq},
\end{equation}
where $\{\tau_{k}\}_{k = 0}^{N}$ is an ascending order set of roots of the shifted Legendre polynomials $P_{N+1}$, $\{\varphi_{p, k}\}$ is a set of coefficients of the basis polynomials $\{\varphi_{p}\}$. The coefficients $\{\varphi_{p, l}\}$ of the Lagrange interpolation polynomials $\varphi_{p}$ are determined from the solution of a system of linear algebraic equations, which is technically implemented in the form of taking the inverse matrix to the Vandermonde matrix $||\tau_{k}^{l}||$: $||\varphi_{p, k}||^{T} = ||\tau_{p}^{k}||^{-1}$. The shifted Legendre polynomials $P_{N+1}$ do not have multiple roots $\tau_{k}$, which proves the existence and uniqueness of the set of basis polynomials $\{\varphi_{p}\}$.

The choice of the set of the basis functions $\{\varphi_{p}\}$ allowed to use the Gauss-Legendre (GL) quadrature formula:
\begin{equation}\label{eq:gl_rule}
\begin{split}
&\int\limits_{0}^{1} f(\tau)d\tau \approx \sum\limits_{p = 0}^{N} w_{p} f(\tau_{p}),\qquad \sum\limits_{p = 0}^{N} w_{p} = 1,\\
&w_{p} = \int\limits_{0}^{1} \varphi_{p}^{2}(\tau) d\tau = \int\limits_{0}^{1} \varphi_{p}(\tau) d\tau > 0,
\end{split}
\end{equation}
where $\{w_{p}\}$ is the weights of GL quadrature formula, to calculate integrals (\ref{eq:ivp_ode_int_src}) and obtain point-wise evaluation
\begin{equation}\label{eq:point_wise_evaluation}
\mathbf{F}(\mathbf{q}_{n}(\tau), t(\tau)) = \sum\limits_{p = 0}^{N} \hat{\mathbf{F}}_{n, p}\varphi_{p}(\tau)\ \mapsto\ 
\sum\limits_{p = 0}^{N} \mathbf{F}(\hat{\mathbf{q}}_{n, p}, t(\tau_{n}))\varphi_{p}(\tau),
\end{equation}
which is based precisely on the use of the GL quadrature formula for calculating the expansion coefficients:
\begin{equation}
\begin{split}
\hat{\mathbf{F}}_{n, p} &=
\left[\int\limits_{0}^{1} \mathbf{F}(\mathbf{q}_{n}(\tau), t(\tau)) \varphi_{p}(\tau) d\tau\right]/\left[\int\limits_{0}^{1} \varphi_{p}^{2}(\tau) d\tau\right]\\ &\approx
\frac{1}{w_{p}} \sum\limits_{q = 0}^{N} w_{q} \mathbf{F}(\mathbf{q}_{n}(\tau_{q}), t(\tau_{q})) \varphi_{p}(\tau_{q}) =
\mathbf{F}(\hat{\mathbf{q}}_{n, p}, t(\tau_{p})).
\end{split}
\end{equation}

Evaluating the weak form (\ref{eq:ivp_ode_weak}) of the ODE system (\ref{eq:ivp_ode_diff_src}), with using (\ref{eq:basis_funcs_def}), (\ref{eq:gl_rule}) and (\ref{eq:point_wise_evaluation}), leads to a system of nonlinear algebraic equations on the expansion coefficients $\{\hat{\mathbf{q}}_{n, p}\}$ of the local discrete time solution $\mathbf{q}_{n}$ (\ref{eq:local_sol}):
\begin{equation}\label{eq:ivp_ode_weak_rewr_dst}
\sum\limits_{q = 0}^{N}\Big[\mathrm{K}_{pq}\hat{\mathbf{q}}_{n, q} - {\Delta t}_{n} \mathrm{M}_{pq}\mathbf{F}(\hat{\mathbf{q}}_{n, q}, t(\tau_{q}))\Big] = 
\varphi_{p}(0)\mathbf{u}_{n}.
\end{equation}
The coefficients can be assembled into matrices $\mathrm{K} = ||\mathrm{K}_{pq}||$, $\mathrm{M} = ||\mathrm{M}_{pq}||$ and vectors $||\varphi_{p}(0)||$, $||\varphi_{p}(1)||$, and can be calculated by the following expressions:
\begin{equation}\label{eq:kappa_mu_defs}
\begin{split}
&\mathrm{K}_{pq} = \varphi_{p}(1)\varphi_{q}(1) - \int\limits_{0}^{1} \frac{d\varphi_{p}(\tau)}{d\tau}\varphi_{q}(\tau)d\tau
			= \varphi_{p}(0)\varphi_{q}(0) + \int\limits_{0}^{1} \varphi_{p}(\tau)\frac{d\varphi_{q}(\tau)}{d\tau}d\tau,\\
&\mathrm{M}_{pq} = \int\limits_{0}^{1} \varphi_{p}(\tau)\varphi_{q}(\tau)d\tau \equiv w_{p} \delta_{pq},\quad
\varphi_{p}(0) = \varphi_{p, 0},\quad \varphi_{p}(1) = \sum\limits_{k = 0}^{N} \varphi_{p, k},
\end{split}
\end{equation}
where the second expression for $\mathrm{K}_{pq}$ is obtained by integrating by parts the integral included in it, and the expression for $\mathrm{M}_{pq}$ explicitly takes into account the orthogonality of the basis $\{\varphi_{p}\}$ and the expression for the weights $w_{p}$ of the GL quadrature formula (\ref{eq:gl_rule}). The coefficients of vectors $||\varphi_{p}(0)||$, $||\varphi_{p}(1)||$ and matrix $\mathrm{K}$ satisfy the following directly verifiable expressions~\cite{ader_dg_ode_sinum}:
\begin{equation}\label{eq:ader_int_prop}
\begin{split}
&\sum\limits_{q = 0}^{N}\mathrm{K}_{pq}  =
\varphi_{p}(0) \sum\limits_{q = 0}^{N}\varphi_{q}(0) + \int\limits_{0}^{1} \varphi_{p}(\tau) \frac{d}{d\tau}\left[\sum\limits_{q = 0}^{N}\varphi_{q}(\tau)\right]d\tau = \varphi_{p}(0),\\
&\sum\limits_{p = 0}^{N}\mathrm{K}_{pq} = \varphi_{q}(1) \sum\limits_{p = 0}^{N}\varphi_{p}(1) - 
\int\limits_{0}^{1} \varphi_{q}(\tau) \frac{d}{d\tau}\left[\sum\limits_{p = 0}^{N}\varphi_{p}(\tau)\right]d\tau = \varphi_{q}(1),\\
&\sum\limits_{q = 0}^{N}\left[\mathrm{K}^{-1}\right]_{pq}\varphi_{q}(0) = 1,\qquad
\sum\limits_{p = 0}^{N}\left[\mathrm{K}^{-1}\right]_{pq}\varphi_{p}(1) = 1,
\end{split}
\end{equation}
where a corollary of the exact point-wise evaluation of the constant-function $f(\tau) \equiv 1$ as expansion over a set of basis functions $\{\varphi_{p}\}$ is taken into account.

Multiplication form (\ref{eq:ivp_ode_weak_rewr_dst}) by $\mathrm{K}^{-1}$ leads to the system of nonlinear algebraic equations of the local DG predictor:
\begin{equation}\label{eq:snae_lstdg}
\hat{\mathbf{q}}_{n, p} - {\Delta t}_{n}\sum\limits_{q = 0}^{N} \mathrm{A}_{pq}\mathbf{F}(\hat{\mathbf{q}}_{n, q}, t(\tau_{q})) = \mathbf{u}_{n},
\end{equation}
where $\mathrm{A}_{pq}$ is elements of the matrix $\mathrm{A} = ||\mathrm{A}_{pq}|| = \mathrm{K}^{-1}\mathrm{M}$, and the property (\ref{eq:ader_int_prop}) is taken into account. This system (\ref{eq:snae_lstdg}) is generally nonlinear and can be solved by Newton's iteration method or Picard's iteration method, as well as other methods for solving systems of nonlinear algebraic equations~\cite{Hairer_book_2, Butcher_book_2016, Ortega_book_1970, Demidovich_1981, Burden_1981, Traub_1982, Evans_1995}.

The expression (\ref{eq:ivp_ode_int_src}) of the solution $\mathbf{u}_{n+1}$ at the grid node $t_{n+1}$ as a result of using the quadrature formula (\ref{eq:gl_rule}) takes the form:
\begin{equation}\label{eq:ader_dg_node_sol}
\mathbf{u}_{n+1} = \mathbf{u}_{n} + {\Delta t}_{n}\sum\limits_{p = 0}^{N} w_{p} \mathbf{F}\left(\hat{\mathbf{q}}_{n, p}, t(\tau_{p})\right).
\end{equation}
The expressions (\ref{eq:snae_lstdg}), (\ref{eq:ader_dg_node_sol}) together constitute the ADER-DG method with the local time DG predictor for solving the initial value problem for ODE system (\ref{eq:ivp_ode_diff_src}).

The work~\cite{ader_dg_ode_jsc} demonstrated that the ADER-DG numerical method for solving an ODE system (\ref{eq:ivp_ode_diff_src}) has a convergence order $p_{\rm G} = 2N+1$ for the solution $\mathbf{u}_{n}$ (\ref{eq:ader_dg_node_sol}) at grid nodes and a convergence order $p_{\rm L} = N+1$ (\ref{eq:local_sol_assembly}) for the local solution $\mathbf{u}_{L}(t)$ in the space between nodes, using polynomials $\{\varphi_{p}\}$ of degree $N$. Furthermore, the ADER-DG numerical method is highly stable and is $A$-stable and $L_{1}$-stable. A rigorous theoretical explanation for these properties of the ADER-DG numerical method is proposed~\cite{ader_dg_ode_sinum}, and several additional properties are identified, such as algebraic stability and $B$-stability.

The software implementation of the ADER-DG numerical method with a local DG predictor is developed using the \texttt{python} programming language. The error of the numerical solution obtained using this numerical method becomes very small even on fairly coarse grids. Therefore, in the software implementation, floating-point numbers of arbitrarily high precision are used within module \texttt{mpmath} (with an additional import of module \texttt{gmpy2}) of the \texttt{python} programming language, the value \texttt{mpmath.mp.dps = 500} is chosen. All calculations presented in this work are performed using this developed software implementation.

\section{Improving the local solution}
\label{sec:imp_dg}

In this Section, an improvement to the local solution is presented and justified, resulting in an improved local solution with higher accuracy and improved smoothness. Before introducing the local solution improvement, the fundamental properties of the original local solution $\mathbf{u}_{L}(t)$, which will be used in particular for the improvement, are presented and rigorously proven.

\subsection{Properties of the local solution}
\label{sec:imp_dg:loc_sol_prop}

The local solution $\mathbf{u}_{L}(t)$, as demonstrated in~\cite{ader_dg_ode_jsc}, possesses a high accuracy, especially for high polynomial degrees $N$. Moreover, at the right-hand point $t_{n+1}$ of each discretization domain $\Omega_{n}$, corresponding to $\tau = 1$, the local solution $\mathbf{q}_{n}(\tau)$ is numerically identical to the solution $\mathbf{u}_{n+1}$ at the grid node $t_{n+1}$:
\begin{equation}\label{eq:local_sol_rp}
\begin{split}
\mathbf{q}_{n}(1)
&= \sum\limits_{p = 0}^{N} \hat{\mathbf{q}}_{n, p} \varphi_{p}(1)
= \mathbf{u}_{n} + {\Delta t}_{n}\sum\limits_{q = 0}^{N}\left[
	\sum\limits_{p = 0}^{N} \mathrm{A}_{pq} \varphi_{p}(1)
\right]\mathbf{F}(\hat{\mathbf{q}}_{n, q}, t(\tau_{q}))\\
&= \mathbf{u}_{n} + {\Delta t}_{n}\sum\limits_{q = 0}^{N}\left[
	\sum\limits_{r = 0}^{N} [\mathrm{K}^{-1}\cdot\mathrm{M}]_{rq} \varphi_{r}(1)
\right]\mathbf{F}(\hat{\mathbf{q}}_{n, q}, t(\tau_{q}))\\
&= \mathbf{u}_{n} + {\Delta t}_{n}\sum\limits_{q = 0}^{N}\left[
	\sum\limits_{r = 0}^{N}\sum\limits_{k = 0}^{N} \left[\mathrm{K}^{-1}\right]_{rk} w_{k} \delta_{kq} \varphi_{r}(1)
\right]\mathbf{F}(\hat{\mathbf{q}}_{n, q}, t(\tau_{q}))\\
&= \mathbf{u}_{n} + {\Delta t}_{n}\sum\limits_{q = 0}^{N}\left[
	\left(\sum\limits_{r = 0}^{N}\left[\mathrm{K}^{-1}\right]_{rq} \varphi_{r}(1)\right) w_{q}
\right]\mathbf{F}(\hat{\mathbf{q}}_{n, q}, t(\tau_{q}))\\
&= \mathbf{u}_{n} + {\Delta t}_{n}\sum\limits_{q = 0}^{N} w_{q}\mathbf{F}(\hat{\mathbf{q}}_{n, q}, t(\tau_{q})) \equiv \mathbf{u}_{n+1},
\end{split}
\end{equation}
where the property (\ref{eq:ader_int_prop}) was used. However, this property does not hold for the local solution $\mathbf{q}_{n}(t)$ at the left-hand point $t_{n}$ of each discretization domain $\Omega_{n}$, corresponding to $\tau = 0$:
\begin{equation}\label{eq:local_sol_lp}
\begin{split}
\mathbf{q}_{n}(0)
&= \sum\limits_{p = 0}^{N} \hat{\mathbf{q}}_{n, p} \varphi_{p}(0)
= \mathbf{u}_{n} + {\Delta t}_{n}\sum\limits_{q = 0}^{N}\left[
	\sum\limits_{p = 0}^{N} \mathrm{A}_{pq} \varphi_{p}(0)
\right]\mathbf{F}(\hat{\mathbf{q}}_{n, q}, t(\tau_{q}))\\
&= \mathbf{u}_{n} + {\Delta t}_{n}\sum\limits_{q = 0}^{N}\left[
	\left(\sum\limits_{r = 0}^{N}\left[\mathrm{K}^{-1}\right]_{rq} \varphi_{r}(0)\right) w_{q}
\right]\mathbf{F}(\hat{\mathbf{q}}_{n, q}, t(\tau_{q}))\\
&= \mathbf{u}_{n} + {\Delta t}_{n}\sum\limits_{q = 0}^{N} \tilde{w}_{q}\mathbf{F}(\hat{\mathbf{q}}_{n, q}, t(\tau_{q})),\quad
\tilde{w}_{q} = w_{q} \sum\limits_{r = 0}^{N}\left[\mathrm{K}^{-1}\right]_{rq} \varphi_{r}(0),
\end{split}
\end{equation}
where the vector $||\tilde{w}_{q}|| \neq 0$ due to $||w_{q}|| \neq 0$, $||\varphi_{p}(0)|| \neq 0$ and non-degeneracy of matrix $\mathrm{K}^{-1}$, therefore $\mathbf{q}_{n}(0) \neq \mathbf{u}_{n}$ in the general case $\mathbf{F} \not\equiv 0$. It turns out that the initial condition in the weak formulation of the problem is not strictly satisfied in the general case $\mathbf{F}$. This property of the local solution leads to discontinuities at the grid nodes $t_{n}$: $\mathbf{u}_{L}(t \rightarrow t_{n}^{+}) \neq \mathbf{u}_{L}(t \rightarrow t_{n}^{-})$, which degrades the local solutions properties and makes it less attractive for use, for example, compared to high-precision reconstructions of the solution. These discontinuities were clearly observed and noted in~\cite{ader_dg_ode_jsc} and will be presented in the Examples in work.

The rigorous proof of the approximation order $p_{\rm L} = N+1$ of the local solution $\mathbf{u}_{L}(t)$ obviously follows from the fact that the weak formulation (\ref{eq:ivp_ode_weak}) with basis functions $\{\varphi_{p}\}$ (\ref{eq:basis_funcs_def}) gives an exact solution $\mathbf{u}$ to the problem (\ref{eq:ivp_ode_diff_src}) for the right-hand side $\mathbf{F} \in \mathcal{P}_{N-1}^{D}(\Omega)$, where $\mathcal{P}_{k}(\Omega)$ is the space of polynomials of degree no greater than $k$ defined on $\Omega$ --- in this case, for a sufficiently smooth functions $\mathbf{u}$ and $\mathbf{F}$, an interpolation representation can be used for solution $\mathbf{u}(t)$ whose error can be estimated as $O({\Delta t}_{n}^{N+1})$~\cite{ader_improving_2024, ader_dg_ode_sinum}. However, let us consider more technical rigorous proof of the approximation order $p_{\rm L} = N+1$ of the local solution $\mathbf{u}_{L}(t)$, since the justification and proof of an improved local solution will be based on this result. The uniform norm will be used to calculate the vector norms below:
\begin{equation}\label{eq:norm_in_rd}
|\mathbf{u}| \equiv ||\mathbf{u}||_{\mathcal{L}_{\infty}} = \max\limits_{1 \leqslant i \leqslant D} |u_{i}|, \quad \mathbf{u} = ||u_{i}|| \in \mathcal{R}^{D},
\end{equation}
which generally does not limit the generality of the results obtained using it due to the equivalence of various possible norms in $\mathcal{R}^{D}$. Rigorous proofs of approximation and the derivation of approximation orders are carried out under the assumption that the right-hand side function $\mathbf{F}(\mathbf{u}, t)$ of the ODE system (\ref{eq:ivp_ode_diff_src}) satisfies the Lipschitz condition for argument $\mathbf{u}$~\cite{dg_ivp_ode_4, dg_ivp_ode_5, dg_ivp_ode_6, dg_ivp_ode_1, dg_ivp_ode_2, dg_ivp_ode_3}, which in this paper is chosen in the following form:
\begin{equation}\label{eq:lip_cond_f}
|\mathbf{F}(\mathbf{u}_{1}, t) - \mathbf{F}(\mathbf{u}_{2}, t)| \leqslant L_{F} |\mathbf{u}_{1} - \mathbf{u}_{2}|,
\end{equation}
for all vectors $\mathbf{u}_{1}, \mathbf{u}_{2} \in \mathcal{R}^{D}$, where $L_{F}$ is the Lipschitz constant for function $\mathbf{F}(\mathbf{u}, t)$ for argument $\mathbf{u}$. A complete proof of the approximation order $p_{\rm L} = N+1$ for a local solution $\mathbf{u}_{L}(t)$ is performed in two stages: first, the approximation order for the local solution at the nodes $\{t_{n} + \tau_{p}{\Delta t}_{n}\}$ of the quadrature formula (\ref{eq:gl_rule}) for all discretization domains $\Omega_{n}$ is proved, and then, based on this, the approximation order for the local solution $\mathbf{u}_{L}(t)$ is proved throughout the domain of definition $\Omega$ of the solution.

The local solution $\mathbf{q}_{n}(\tau)$ over the discretization domain $t\in\Omega_{n}$ of the solution is represented in the form:
\begin{equation}\label{eq:local_sol_on_tau}
\begin{split}
\mathbf{q}_{n}(\tau) = \sum\limits_{p = 0}^{N} \hat{\mathbf{q}}_{n, p}\varphi_{p}(\tau) &=
\sum\limits_{p = 0}^{N} \varphi_{p}(\tau)\left[\mathbf{u}_{n} + {\Delta t}_{n}\sum\limits_{q = 0}^{N} \mathrm{A}_{pq}\mathbf{F}(\hat{\mathbf{q}}_{n, q}, t(\tau_{q}))\right]\\
&= \mathbf{u}_{n} + {\Delta t}_{n} \sum\limits_{p = 0}^{N}\sum\limits_{q = 0}^{N} \varphi_{p}(\tau)\mathrm{A}_{pq}\mathbf{F}(\hat{\mathbf{q}}_{n, q}, t(\tau_{q})),
\end{split}
\end{equation}
where a corollary of the exact point-wise evaluation of the constant-function $f(\tau) \equiv 1$ as expansion over a set of basis functions $\{\varphi_{p}\}$ is taken into account.
The local solution $\mathbf{u}_{L}(t)$ over the domain of definition $t\in\Omega$ of the solution is represented in the form (\ref{eq:local_sol_assembly}):
\begin{equation}\label{eq:local_sol_assembly_by_phi}
\mathbf{u}_{L}(t) = \sum\limits_{n} \chi_{n}(t)\left[
	\mathbf{u}_{n} + {\Delta t}_{n} \sum\limits_{p = 0}^{N}\sum\limits_{q = 0}^{N}
	\varphi_{p}\left(\frac{t - t_{n}}{{\Delta t}_{n}}\right)\mathrm{A}_{pq}\mathbf{F}(\hat{\mathbf{q}}_{n, q}, t(\tau_{q}))
\right],
\end{equation}
therefore, for each discretization domain $\Omega_{n}$, the contribution to the error of the local solution $\mathbf{u}_{L}(t)$ will be determined not only by the error of the local solution $\mathbf{q}_{n}(\tau)$ associated with the limited representation over a finite set of basis functions $\{\varphi_{p}\}$, but also by the error of the solution at the nodes $\mathbf{u}_{n}$. The function $\mathbf{F}(\mathbf{u}, t)$ on the right-hand side of the ODE system (\ref{eq:ivp_ode_diff_src}), with the chosen function $\mathbf{u} = \mathbf{u}(t)$, admits a representation in the form of polynomial interpolation by basis functions $\{\varphi_{p}\}$:
\begin{equation}
\begin{split}
\mathbf{F}(\mathbf{u}(t(\tau)),\, t(\tau)) &=
\mathrm{P}_{\varphi}[\mathbf{F}(\mathbf{u}(t(\tau)),\, t(\tau))] + \mathbf{R}_{F}^{\varphi}(t(\tau))\\ &=
\sum\limits_{p = 0}^{N} \mathbf{F}(\mathbf{u}(t(\tau_{p})),\, t(\tau_{p})) \varphi_{p}(\tau) + \mathbf{R}_{F}^{\varphi}(t(\tau)),\ 
t \in \Omega_{n},
\end{split}
\end{equation}
where $\mathrm{P}_{\varphi}$ is the operator representing the function of argument $\tau\in\omega$ over a set of basis functions $\{\varphi_{p}\}$, and $\mathbf{R}_{F}^{\varphi}(t)$ is the remainder term of the interpolation formula, for which the estimate $\mathbf{R}_{F}^{\varphi} = O\left({\Delta t}_{n}^{N+1}\right)$ is well known from the theory of Lagrangian interpolation polynomials. Using this representation, setting $\mathbf{u}(t)$ as the exact solution of the initial value problem for the ODE system (\ref{eq:ivp_ode_diff_src}), the following representation can be obtained for the exact solution $\mathbf{u}(t)$ over the discretization domain $t\in\Omega_{n}$:
\begin{equation}\label{eq:u_ex_exp_by_int_phi}
\begin{split}
\mathbf{u}(t(\tau)) &= \mathbf{u}(t_{n}) + \int\limits_{t_{n}}^{t(\tau)} \mathbf{F}\left(\mathbf{u}(t),\, t\right) dt\\
&= \mathbf{u}(t_{n}) + \sum\limits_{p = 0}^{N} \mathbf{F}(\mathbf{u}(t(\tau_{p})),\, t(\tau_{p})) \int\limits_{t_{n}}^{t(\tau)} \varphi_{p}(\tau(t)) dt +
\int\limits_{t_{n}}^{t(\tau)} \mathbf{R}_{F}^{\varphi}(t) dt\\
&= \mathbf{u}(t_{n}) + 
{\Delta t}_{n} \sum\limits_{p = 0}^{N} \mathbf{F}(\mathbf{u}(t(\tau_{p})),\, t(\tau_{p})) \int\limits_{0}^{\tau(t)} \varphi_{p}(\xi) d\xi + \mathbf{R}_{u}^{\Phi}(t(\tau)),\ 
t \in \Omega_{n},
\end{split}
\end{equation}
where the left boundary condition of the mapping (\ref{eq:tau_mapping}) is used in the last expression, and $\mathbf{R}_{u}^{\Phi}(t)$ is the remainder term, for which the estimate $\mathbf{R}_{u}^{\Phi} = O\left({\Delta t}_{n}^{N+2}\right)$ is known from the estimate for the remainder term $\mathbf{R}_{F}^{\varphi}(t)$. This representation of the solution $\mathbf{u}$ is an expansion over a set of new linearly independent (but not orthogonal) basis functions:
\begin{equation}\label{eq:u_exp_by_span_phi}
\begin{split}
&\mathbf{u} = \mathrm{P}_{\Phi}[\mathbf{u}] + \mathbf{R}_{u}^{\Phi},\quad
\Phi = \mathrm{span}\left[1,\, \left\{\int_{0}^{\tau}\varphi_{p}(\xi)d\xi\right\}\right],\quad t \in \Omega_{n},\\
&\mathrm{P}_{\Phi}[\mathbf{u}](t(\tau)) = \mathbf{u}(t_{n}) \cdot 1 + 
{\Delta t}_{n} \sum\limits_{p = 0}^{N} \mathbf{F}(\mathbf{u}(t(\tau_{p})),\, t(\tau_{p})) \cdot \int\limits_{0}^{\tau} \varphi_{p}(\xi) d\xi,
\end{split}
\end{equation}
where $\mathrm{P}_{\Phi}$ is the operator representing the function of argument $\tau\in\omega$ on the linear span $\Phi$, which allows an exact representation of the functions $\mathcal{P}_{N+1}(\Omega_{n})$. The difference between the local numerical solution $\mathbf{q}_{n}$ and the exact solution $\mathbf{u}$ in the discretization domain $\Omega_{n}$, expressed as a function of the mapping parameter $\tau$, shows that it is not possible to isolate a common factor from the expression in brackets.

To evaluate the obtained expression (\ref{eq:u_ex_exp_by_int_phi}) and calculate the accuracy order, it is necessary to present both terms in brackets of the obtained expression in the form of an expansion in the same set of basis functions, for which the original set of basis functions $\{\varphi_{p}\}$ of the ADER-DG numerical method is selected, for which the representation of the integration operator is constructed. For this purpose, the following property of integration of a power function $f(\tau) = \tau^{k}$ is defined:
\begin{equation}\label{eq:c_prop_def}
\sum\limits_{q = 0}^{N} \mathrm{A}_{pq}\tau_{q}^{k} = \frac{\tau_{p}^{k+1}}{k+1},\quad 0 \leqslant k \leqslant N-1,
\end{equation}
\corrtext{which corresponds to simplifying condition $C(N)$~\cite{ader_improving_2024, ader_dg_ode_sinum} in the theory of numerical methods for solving ODE systems~\cite{Hairer_book_1, Hairer_book_2, Dekker_Verwer_1984}.} This property can be conveniently rewritten using matrix notation and introducing the vector $\mathbf{t}_{k} = [\tau_{0}^{k} \ldots \tau_{N}^{k}]^{T}$:
\begin{equation}\label{eq:c_prop_matrix}
\mathrm{A}\mathbf{t}_{k} = \frac{\mathbf{t}_{k+1}}{k+1},\ \Rightarrow\ 
\mathrm{M}\mathbf{t}_{k} = \frac{\mathrm{K}\mathbf{t}_{k+1}}{k+1},\quad 0 \leqslant k \leqslant N-1.
\end{equation}
To rigorously prove the property (\ref{eq:c_prop_def}), the right-hand side of the resulting matrix expression (\ref{eq:c_prop_matrix}) is directly calculated, leading to the following chain of equalities:
\begin{equation}\label{eq:c_prop_proof}
\begin{split}
\sum\limits_{q = 0}^{N} &\frac{\mathrm{K}_{pq} \tau_{q}^{k+1}}{k+1} =
\sum\limits_{q = 0}^{N} \left[\varphi_{p}(0)\varphi_{q}(0) + \int\limits_{0}^{1} \varphi_{p}(\tau)\frac{d\varphi_{q}(\tau)}{d\tau}d\tau\right]\frac{\tau_{q}^{k+1}}{k+1}\\
&= \frac{\varphi_{p}(0)}{k+1} \sum\limits_{q = 0}^{N} \tau_{q}^{k+1}\varphi_{q}(0) +
\frac{1}{k+1} \int\limits_{0}^{1} \varphi_{p}(\tau) \frac{d}{d\tau}\left[\sum\limits_{q = 0}^{N}\tau_{q}^{k+1}\varphi_{q}(\tau)\right]d\tau\\
&= \frac{\varphi_{p}(0)}{k+1} (0)^{k+1} + \frac{1}{k+1} \int\limits_{0}^{1} \varphi_{p}(\tau) \frac{d\tau^{k+1}}{d\tau} d\tau
= \frac{1}{k+1} \int\limits_{0}^{1} \frac{d\tau^{k+1}}{d\tau}\varphi_{p}(\tau)d\tau\\
&= \int\limits_{0}^{1} \tau^{k}\varphi_{p}(\tau)d\tau = \sum\limits_{q = 0}^{N} w_{q} \tau_{q}^{k} \varphi_{q}(\tau_{p}) =
\sum\limits_{q = 0}^{N} w_{q} \tau_{q}^{k} \delta_{qp} = w_{p}\tau_{p}^{k},
\end{split}
\end{equation}
which leads to the left-hand side of the resulting matrix expression (\ref{eq:c_prop_matrix}), thus proving the property (\ref{eq:c_prop_def}). It is important to note that this set of transformations (\ref{eq:c_prop_proof}) only holds for the case $k \leqslant N-1$ --- in the case $k \geqslant N$, the function $\tau^{k+1}$ cannot be exactly represented as an expansion in basis polynomials $\{\varphi_{p}(\tau)\}$ of degree $N$. Multiplying the resulting expression (\ref{eq:c_prop_def}) by the basis function $\varphi_{p}(\tau)$ and summing over index $0 \leqslant p \leqslant N$ yields the following expression for the integral of the power function $f(\tau) = \tau^{k}$:
\begin{equation}
\int\limits_{0}^{\tau} \tau^{k} d\tau =
\frac{\tau^{k+1}}{k+1} = \sum\limits_{p = 0}^{N} \frac{\tau_{p}^{k+1}\varphi_{p}(\tau)}{k+1} =
\sum\limits_{p = 0}^{N}\sum\limits_{q = 0}^{N} \varphi_{p}(\tau)\mathrm{A}_{pq}\tau_{q}^{k},\ 0 \leqslant k \leqslant N-1,
\end{equation}
where, again, the degree $k \leqslant N-1$. Based on this representation for the integral of a power function $\tau^{k}$, the representation is constructed for the integral of basis polynomials $\{\varphi_{p}(\tau)\}$, in which all degrees except the highest $N$ are accurately taken into account:
\begin{equation}
\begin{split}
&\int\limits_{0}^{\tau} \varphi_{q}(\xi) d\xi = \sum\limits_{k = 0}^{N} \varphi_{q, k} \frac{\tau^{k+1}}{k+1} \mapsto
\sum\limits_{k = 0}^{N-1} \varphi_{q, k} \frac{\tau^{k+1}}{k+1} = 
\sum\limits_{k = 0}^{N-1} \varphi_{q, k} \sum\limits_{p = 0}^{N}\sum\limits_{r = 0}^{N} \varphi_{p}(\tau)\mathrm{A}_{pr}\tau_{r}^{k}\\
&= \sum\limits_{p = 0}^{N}\sum\limits_{r = 0}^{N} \varphi_{p}(\tau)\mathrm{A}_{pr} \left[\sum\limits_{k = 0}^{N} \varphi_{q, k}\tau_{r}^{k} - \varphi_{q, N}\tau_{r}^{N}\right] =
\sum\limits_{p = 0}^{N}\sum\limits_{r = 0}^{N} \varphi_{p}(\tau) \mathrm{A}_{pr} \left[\varphi_{q}(\tau_{r}) - \varphi_{q, N}\tau_{r}^{N}\right]\\
&= \sum\limits_{p = 0}^{N}\sum\limits_{r = 0}^{N} \varphi_{p}(\tau) \mathrm{A}_{pr} \left[\delta_{q, r} - \varphi_{q, N}\tau_{r}^{N}\right] =
\sum\limits_{p = 0}^{N} \varphi_{p}(\tau)\mathrm{A}_{pq} - \sum\limits_{p = 0}^{N}\sum\limits_{r = 0}^{N} \varphi_{p}(\tau) \mathrm{A}_{pr} \varphi_{q, N}\tau_{r}^{N}\\
&\phantom{= \sum\limits_{p = 0}^{N}\sum\limits_{r = 0}^{N} \varphi_{p}(\tau) \mathrm{A}_{pr} \left[\delta_{q, r} - \varphi_{q, N}\tau_{r}^{N}\right]\,} =
\sum\limits_{p = 0}^{N} \varphi_{p}(\tau)\mathrm{A}_{pq} - \varphi_{q, N}\sum\limits_{p = 0}^{N} \sum\limits_{r = 0}^{N} \varphi_{p}(\tau) \mathrm{A}_{pr} \tau_{r}^{N},
\end{split}
\end{equation}
where the appearance of the last term is due to the inaccuracy of the representation of a polynomial of degree $N+1$ using polynomials $\{\varphi_{p}(\tau)\}$ of degree $N$. Then the operator for representing integrals of basis polynomials $\{\varphi_{p}(\tau)\}$ in the form of an expansion in basis polynomials $\{\varphi_{p}(\tau)\}$ is chosen in the following form:
\begin{equation}
\mathrm{P}_{\varphi}\left[\int\limits_{0}^{\tau} \varphi_{q}(\xi) d\xi\right] =
\sum\limits_{p = 0}^{N} \varphi_{p}(\tau)\mathrm{A}_{pq},
\end{equation}
which defines a multiplicative action, in the sense of matrix multiplication, on the vector of coefficients of the function representation chosen in the form of expansion in basis polynomials $\{\varphi_{p}(\tau)\}$, and has an error of $O({\Delta t}_{n}^{N+1})$ due to the use of polynomials $\{\varphi_{p}(\tau)\}$ of degree $N$ for the representation. In this case, the integral of function $f = f(\tau)$, defined by the interpolation representation over a set of basis polynomials $\{\varphi_{p}(\tau)\}$ with a remainder term $R_{f}^{\varphi}$, takes the following form:
\begin{equation}
\begin{split}
&f(\tau) = \sum\limits_{p = 0}^{N} f_{p} \varphi_{p}(\tau) + R_{f}^{\varphi} = \sum\limits_{p = 0}^{N} f(\tau_{p}) \varphi_{p}(\tau) + R_{f}^{\varphi},\\
&\int\limits_{0}^{\tau} f(\xi) d\xi = \sum\limits_{q = 0}^{N} f(\tau_{q}) \mathrm{P}_{\varphi}\left[\int\limits_{0}^{\tau} \varphi_{q}(\xi) d\xi\right] =
\sum\limits_{p = 0}^{N} \sum\limits_{q = 0}^{N} \varphi_{p}(\tau)\mathrm{A}_{pq} f(\tau_{q}) + R_{if}^{\varphi},
\end{split}
\end{equation}
where the remainder term $R_{if}^{\varphi}$ (if the integral of the function $f(\tau)$ is calculated on a grid with partitioning $t_{n}$, and the function $f(\tau)$ is a slice of the function $f(t)$ on the discretization domain $\Omega_{n}$) has the same order as the remainder term $R_{f}^{\varphi}$ of the original expansion and not one unit higher, as in case (\ref{eq:u_ex_exp_by_int_phi}): $R_{f}^{\varphi} = O({\Delta t}_{n}^{N+1})$ and $R_{if}^{\varphi} = O({\Delta t}_{n}^{N+1})$, since interpolation is also carried out using polynomials $\{\varphi_{p}(\tau)\}$ of degree $N$. Using the resulting representation, the expression (\ref{eq:u_ex_exp_by_int_phi}) for the exact solution is the following:
\begin{equation}
\begin{split}
\mathbf{u}(t(\tau)) &= \mathrm{P}_{\varphi}\left[\mathbf{u}(t(\tau)) + \mathbf{R}_{u}^{\Phi}(t(\tau))\right] =
\mathrm{P}_{\varphi}\left[\mathbf{u}(t(\tau))\right] + \mathbf{R}_{u}^{\varphi}(t(\tau))\\
&= \mathbf{u}(t_{n}) + 
{\Delta t}_{n} \sum\limits_{p = 0}^{N} \mathbf{F}(\mathbf{u}(t(\tau_{p})),\, t(\tau_{p}))\,
\mathrm{P}_{\varphi}\left[\int\limits_{0}^{\tau(t)} \varphi_{p}(\xi) d\xi\right] + \mathbf{R}_{u}^{\varphi}(t(\tau))\\ &= \mathbf{u}(t_{n}) + 
{\Delta t}_{n} \sum\limits_{p = 0}^{N}\sum\limits_{q = 0}^{N}
\varphi_{p}(\tau)\mathrm{A}_{pq} \mathbf{F}(\mathbf{u}(t(\tau_{q})),\, t(\tau_{q})) + \mathbf{R}_{u}^{\varphi}(t(\tau)),
\end{split}
\end{equation}
where the remainder term has the estimate $\mathbf{R}_{u}^{\varphi}(t(\tau)) = O({\Delta t}_{n}^{N+1})$ instead of $\mathbf{R}_{u}^{\Phi}(t(\tau)) = O({\Delta t}_{n}^{N+2})$. 

Subtracting the resulting expression for representing the exact solution from the local numerical solution $\mathbf{q}_{n}(\tau)$ (\ref{eq:local_sol_on_tau}) yields the following expression:
\begin{equation}
\begin{split}
&\mathbf{q}_{n}(\tau) - \mathbf{u}(t(\tau)) = \mathbf{u}_{n} - \mathbf{u}(t_{n})\\
&+ {\Delta t}_{n} \sum\limits_{p = 0}^{N} \sum\limits_{q = 0}^{N}
\varphi_{p}(\tau)\mathrm{A}_{pq}\left[\mathbf{F}(\hat{\mathbf{q}}_{n, q}, t(\tau_{q})) - \mathbf{F}(\mathbf{u}(t(\tau_{q})),\, t(\tau_{q}))\right]
- \mathbf{R}_{u}^{\varphi}(t(\tau)),
\end{split}
\end{equation}
where (\ref{eq:lip_cond_f}) could conveniently be used. However, the presence of double summation on the right-hand side does not allow the desired estimate to be expressed directly, so the estimate for the local numerical solution $\mathbf{q}_{n}(\tau)$ at the quadrature points $\{\tau_{p}\}$ is first calculated. Direct calculation yields result:
\begin{equation}
\begin{split}
&\mathbf{q}_{n}(\tau_{p}) - \mathbf{u}(t(\tau_{p})) = \mathbf{q}_{n, p} - \mathbf{u}(t(\tau_{p})) = \mathbf{u}_{n} - \mathbf{u}(t_{n})\\
&+ {\Delta t}_{n} \sum\limits_{q = 0}^{N}
\mathrm{A}_{pq}\left[\mathbf{F}(\hat{\mathbf{q}}_{n, q}, t(\tau_{q})) - \mathbf{F}(\mathbf{u}(t(\tau_{q})),\, t(\tau_{q}))\right]
- \mathbf{R}_{u}^{\varphi}(t(\tau_{p})),
\end{split}
\end{equation}
resulting in the following error estimates:
\begin{equation}
\begin{split}
|\mathbf{q}_{n, p} &- \mathbf{u}(t(\tau_{p}))| \leqslant |\mathbf{u}_{n} - \mathbf{u}(t_{n})| + |\mathbf{R}_{u}^{\varphi}(t(\tau_{p}))|\\
&+ {\Delta t}_{n} \left|\sum\limits_{q = 0}^{N}
\mathrm{A}_{pq}\left[\mathbf{F}(\hat{\mathbf{q}}_{n, q}, t(\tau_{q})) - \mathbf{F}(\mathbf{u}(t(\tau_{q})),\, t(\tau_{q}))\right]\right|\\
&\leqslant |\mathbf{u}_{n} - \mathbf{u}(t_{n})| + |\mathbf{R}_{u}^{\varphi}(t(\tau_{p}))|
+ {\Delta t}_{n} \sum\limits_{q = 0}^{N} |\mathrm{A}_{pq}| |\mathbf{F}(\hat{\mathbf{q}}_{n, q}, t(\tau_{q})) - \mathbf{F}(\mathbf{u}(t(\tau_{q})),\, t(\tau_{q}))|\\
&\leqslant |\mathbf{u}_{n} - \mathbf{u}(t_{n})| + |\mathbf{R}_{u}^{\varphi}(t(\tau_{p}))|
+ {\Delta t}_{n} L_{F} \sum\limits_{q = 0}^{N} |\mathrm{A}_{pq}| |\hat{\mathbf{q}}_{n, q} -\mathbf{u}(t(\tau_{q}))|,
\end{split}
\end{equation}
where (\ref{eq:lip_cond_f}) was used, and the following estimates can be used:
\begin{equation}
|\mathbf{u}_{n} - \mathbf{u}(t_{n})| \leqslant C_{u} {\Delta t}_{n}^{2N+1},\qquad
|\mathbf{R}_{u}^{\varphi}(t(\tau_{p}))| \leqslant C_{R}^{\varphi} {\Delta t}_{n}^{N+1},
\end{equation}
where the first estimate is a consequence of the accuracy order $p_{\rm G} = 2N+1$~\cite{ader_dg_ode_jsc} for the numerical solution at nodes $\mathbf{u}_{n}$, presented for the numerical calculation at several steps simultaneously (in the case of a single step, estimate $\propto {\Delta t}_{n}^{2N+2}$ could be used). As a result of using the obtained expressions and estimates, a system of linear algebraic equations of the following form can be obtained for the error of the local numerical solution $\mathbf{q}_{n, p}$ at quadrature points $\{\tau_{p}\}$:
\begin{equation}
\begin{split}
\sum\limits_{q = 0}^{N} [ \delta_{pq} - {\Delta t}_{n} L_{F} |\mathrm{A}_{pq}| ] |\mathbf{q}_{n, q} &- \mathbf{u}(t(\tau_{q}))|
\leqslant |\mathbf{u}_{n} - \mathbf{u}(t_{n})| + |\mathbf{R}_{u}^{\varphi}(t(\tau_{p}))|\\
&\leqslant C_{u} {\Delta t}_{n}^{2N+1} + C_{R}^{\varphi} {\Delta t}_{n}^{N+1} \simeq C_{L, p}^{q} {\Delta t}_{n}^{N+1},
\end{split}
\end{equation}
which, in general, is uniquely solvable, and the constant $C_{L, p}^{q} > 0$ can be chosen to depend on the index $p$ of the quadrature point. It is important to note three details of the structure of this system of linear equations: the system matrix $[\mathrm{I} - {\Delta t}_{n} L_{F} \cdot \mathrm{A}_{\rm abs}]$, where $\mathrm{A}_{\rm abs} = ||\,|A_{pq}|\,||$, is nonsingular (except in the unique cases in which ${\Delta t}_{n} L_{F}$ is the characteristic number of the matrix $\mathrm{A}_{\rm abs}$, in which case a different Lipschitz constant $L_{F}$ can always be chosen for the estimate); the vector on the right-hand side is proportional to ${\Delta t}_{n}^{N+1}$ --- all components of the vector have a common factor ${\Delta t}_{n}^{N+1}$; and the asymptotics of the solution ${\Delta t}_{n} \rightarrow 0^{+}$ imply that the inverse matrix has no rows with a common factor ${\Delta t}_{n}$ (otherwise, the matrix $[\mathrm{I} - {\Delta t}_{n} L_{F} \cdot \mathrm{A}_{\rm abs}]$ would be singular in the limit ${\Delta t}_{n} \rightarrow 0^{+}$, but it is the identity matrix $[\mathrm{I} - {\Delta t}_{n} L_{F} \cdot \mathrm{A}_{\rm abs}] \rightarrow \mathrm{I}$, and therefore invertible). Therefore, the resolution of this system of linear equations allows us to obtain an estimate of the error of the local numerical solution $\mathbf{q}_{n, p}$ at quadrature points $\{\tau_{p}\}$:
\begin{equation}\label{eq:local_sol_err_qpnts}
\begin{split}
|\mathbf{q}_{n, p} - \mathbf{u}(t(\tau_{p}))| &\leqslant
\max\limits_{p} \left|\sum\limits_{q = 0}^{N}\left[\left(\mathrm{I} - {\Delta t}_{n} L_{F} \cdot \mathrm{A}_{\rm abs}\right)^{-1}\right]_{pq} C_{L, q}^{q} {\Delta t}_{n}^{N+1}\right|\\
&= \max\limits_{p} \left|\sum\limits_{q = 0}^{N}\left[\left(\mathrm{I} - {\Delta t}_{n} L_{F} \cdot \mathrm{A}_{\rm abs}\right)^{-1}\right]_{pq} C_{L, q}^{q}\right| {\Delta t}_{n}^{N+1}
\leqslant C_{L}^{q} {\Delta t}_{n}^{N+1},
\end{split}
\end{equation}
where the existence of the constant $C_{L}^{q} > 0$ follows from the asymptotics of the system matrix in the limit ${\Delta t}_{n} \rightarrow 0^{+}$.

Using the obtained error estimate of the local numerical solution $\mathbf{q}_{n, q}$ at quadrature points $\{\tau_{p}\}$ (\ref{eq:local_sol_err_qpnts}), the error estimate for the local solution $\mathbf{q}_{n}(\tau)$ in the norm $\mathcal{L}_{\infty}$ is obtained in the following form:
\begin{equation}\label{eq:local_sol_errs_L_inf_est}
\begin{split}
||\mathbf{q}_{n}(\tau) &- \mathbf{u}(t(\tau))||_{\mathcal{L}_{\infty}} =
\operatorname{ess}\sup\limits_{\hspace{-5mm}\tau\in[0,1]}|\mathbf{q}_{n}(\tau) - \mathbf{u}(t(\tau))| =
\max\limits_{\tau\in[0,1]}|\mathbf{q}_{n}(\tau) - \mathbf{u}(t(\tau))|\\
&\leqslant |\mathbf{R}_{u}^{\varphi}(t(\tau_{p}))|
+ \max\limits_{\tau\in[0,1]}\left|\sum\limits_{p = 0}^{N}\varphi_{p}(\tau)[\mathbf{q}_{n, p} - \mathbf{u}(t(\tau_{p}))]\right|\\
&\leqslant C_{R}^{\varphi} {\Delta t}^{N+1} + \Lambda_{\varphi} C_{L}^{q} {\Delta t}^{N+1} =
(C_{R}^{\varphi} + \Lambda_{\varphi} C_{L}^{q}) {\Delta t}^{N+1} =
C_{L}^{\mathcal{L}_{\infty}} {\Delta t}^{N+1},
\end{split}
\end{equation}
where $\Lambda_{\varphi}$ is the Lebesgue constant of the set of basis functions $\{\varphi_{p}\}$, and the characteristic grid discretization step ${\Delta t}$ is chosen, which can be ${\Delta t} = \max_{n} {\Delta t}_{n}$. From this error estimate, a similar error estimate for the local solution $\mathbf{u}_{L}(t)$ (\ref{eq:local_sol_assembly_by_phi}) follows immediately:
\begin{equation}
\begin{split}
||\mathbf{u}_{L}(t) &- \mathbf{u}(t)||_{\mathcal{L}_{\infty}} =
\operatorname{ess}\sup\limits_{\hspace{-5mm}t\in\Omega}\left|\mathbf{u}_{L}(t) - \mathbf{u}(t)\right|\\
&=\max\limits_{n}\sup\limits_{t\in\Omega_{n}}\left|\mathbf{u}_{L}(t) - \mathbf{u}(t)\right| =
\max\limits_{n}\sup\limits_{t\in\Omega_{n}}\left|\mathbf{q}_{n}(t) - \mathbf{u}(t)\right|\\
&= \max\limits_{n} ||\mathbf{q}_{n}(\tau) - \mathbf{u}(t(\tau))||_{\mathcal{L}_{\infty}} \leqslant C_{L}^{\mathcal{L}_{\infty}} {\Delta t}^{N+1}.
\end{split}
\end{equation}
Similar error estimates for the local numerical solution $\mathbf{u}_{L}(t)$ (\ref{eq:local_sol_assembly_by_phi}) are also obtained for other classical error norms in $\mathcal{L}_{p}$ ($p \geqslant 1$):
\begin{equation}
\begin{split}
||\mathbf{u}_{L}(t) - \mathbf{u}(t)||_{\mathcal{L}_{p}} &=
\left[\int\limits_{\Omega} |\mathbf{u}_{L}(t) - \mathbf{u}(t)|^{p} dt \right]^{1/p} \leqslant
\left[|\Omega| \cdot \left(\operatorname{ess}\sup\limits_{\hspace{-5mm}t\in\Omega}\left|\mathbf{u}_{L}(t) - \mathbf{u}(t)\right|\right)^{p}\,\right]^{1/p}\\
&= |\Omega|^{1/p} ||\mathbf{u}_{L}(t) - \mathbf{u}(t)||_{\mathcal{L}_{\infty}} \leqslant
|\Omega|^{1/p} C_{L}^{\mathcal{L}_{\infty}} {\Delta t}^{N+1} \equiv C_{L}^{\mathcal{L}_{p}} {\Delta t}^{N+1}.
\end{split}
\end{equation}
The obtained rigorous proofs substantiates the empirical result~\cite{ader_dg_ode_jsc}, that the local solution $\mathbf{u}_{L}(t)$ (\ref{eq:local_sol_assembly}) obtained using the local DG predictor of the ADER-DG numerical method has an approximation order $p_{\rm L} = N+1$, in norms $\mathcal{L}_{\infty}$ and $\mathcal{L}_{p}$ ($p \geqslant 1$). However, a rigorous proof relies heavily on the Lipschitz condition (\ref{eq:lip_cond_f}) for the function $\mathbf{F}(\mathbf{u},\, t)$ of the right-hand side of the original system of equations (\ref{eq:ivp_ode_diff_src}), which is not uncommon in the case of strictly studying the approximation and convergence of DG methods.

\subsection{Improved local solution}
\label{sec:imp_dg:imp_loc_sol}

A significant problem of the local solution is the approximation order $p_{\rm L} = N+1$ and, especially, the discontinuity of the assembled local solution $\mathbf{u}_{L}(t)$ (\ref{eq:local_sol_assembly}) at the left point $\tau = 0$ of each discretization domain $\Omega_{n}$: $\mathbf{q}_{n}(0) \neq \mathbf{u}_{n}$ (\ref{eq:local_sol_lp}), $\mathbf{q}_{n}(1) = \mathbf{u}_{n+1}$ (\ref{eq:local_sol_rp}), including at the point $t_{0}$ where the initial condition is specified: $\mathbf{q}_{0}(0) \neq \mathbf{u}_{0}$ in (\ref{eq:ivp_ode_diff_src}).

In this paper, the following form of an improved local solution obtained by a local DG predictor, based on the following relation:
\begin{equation}\label{eq:imp_local_sol}
\mathbf{q}_{n}^{\rm IL}(\tau) = \mathbf{u}_{n} + {\Delta t}_{n}\int\limits_{0}^{\tau} \mathbf{F}\left(\mathbf{q}_{n}(\xi), t(\xi)\right) d\xi,\ t\in\Omega_{n},
\end{equation}
is proposed. The improved local solution $\mathbf{q}_{n}^{\rm IL}(\tau)$ is also defined for the definition domain $\mathbf{u}_{L}: \Omega \rightarrow \mathcal{R}^{D}$ by a piecewise assembly of local solutions $\mathbf{q}_{n}$ for each discretization domain $\Omega_{n}$:
\begin{equation}\label{eq:imp_local_sol_assembly}
\mathbf{u}_{\rm IL}(t) = \sum\limits_{n} \chi_{n}(t)\cdot\mathbf{q}_{n}^{\rm IL}\left(\frac{t - t_{n}}{{\Delta t}_{n}}\right),
\end{equation}
where $\chi_{n}: \Omega\rightarrow\{0,\, 1\}$ is the indicator function of the discretization domain $\Omega_{n}\subseteq\Omega$ as before (\ref{eq:local_sol_assembly}). Using a point-wise evaluation (\ref{eq:point_wise_evaluation}) of function $\mathbf{F}$ when calculating the integral leads to the following expression for evaluation the improved local solution (\ref{eq:imp_local_sol}):
\begin{equation}\label{eq:imp_local_sol_by_phi}
\mathbf{q}_{n}^{\rm IL}(\tau) = \mathbf{u}_{n} + 
{\Delta t}_{n} \sum\limits_{p = 0}^{N} \mathbf{F}\left(\mathbf{q}_{n, p}, t(\tau_{p})\right) \int\limits_{0}^{\tau} \varphi_{p}(\xi) d\xi,\ 
t\in\Omega_{n},
\end{equation}
which is a representation of the solution in the form of an expansion in basis functions of the linear span $\Phi$ (\ref{eq:u_exp_by_span_phi}). At the boundary points $\tau = 0$ and $\tau = 1$ of the discretization domain $\Omega_{n}$, the improved local solution $\mathbf{q}_{n}^{\rm IL}(\tau)$ coincides with the numerical solution $\mathbf{u}_{n}$ and $\mathbf{u}_{n+1}$ at the grid nodes $t_{n}$ and $t_{n+1}$:
\begin{equation}\label{eq:imp_local_bound_pnts}
\mathbf{q}_{n}^{\rm IL}(0) = \mathbf{u}_{n},\quad
\mathbf{q}_{n}^{\rm IL}(\tau) = \mathbf{u}_{n} + {\Delta t}_{n} \sum\limits_{p = 0}^{N} w_{p}\mathbf{F}\left(\mathbf{q}_{n, p}, t(\tau_{p})\right) \equiv \mathbf{u}_{n+1},
\end{equation}
where the expression for the weights $w_{p}$ of quadrature formula (\ref{eq:gl_rule}) was used, resulting in the continuity of the improved local solution $\mathbf{u}_{\rm IL}(t)$ defined over the definition domain $\mathbf{u}_{L}$: $\mathbf{u}_{\rm IL}(t \rightarrow t_{n}^{+}) = \mathbf{u}_{\rm IL}(t \rightarrow t_{n}^{-}) = \mathbf{u}_{\rm IL}(t_{n}) = \mathbf{u}_{n}$. This approach to constructing a local solution is similar to obtaining a continuous solution for collocation implicit Runge-Kutta methods, which is well known in the literature~\cite{Hairer_book_1}. However, in this case, the main difference is that the ADER-DG numerical method with a local DG predictor is not initially collocation-based; otherwise, relation $\int_{0}^{\tau_{p}} \varphi_{q}(\tau) d\tau = \mathrm{A}_{pq}$ would hold, and the matrix $\mathrm{A}$ would be identical to the coefficient matrix of the classical Gauss-Legendre method.

Determining the approximation order $p_{\rm IL}$ of the improved local solution $\mathbf{u}_{\rm IL}(t)$ (\ref{eq:imp_local_sol_assembly}), in contrast to the original local solution $\mathbf{u}_{L}(t)$ obtained by the local DG predictor, is based on direct use of the representation of the exact solution $\mathbf{u}(t)$ (\ref{eq:u_ex_exp_by_int_phi}). The error estimate of the improved local solution $\mathbf{q}_{n}^{\rm IL}(\tau)$ in the discretization domain $t\in\Omega_{n}$ is obtained by subtracting expression (\ref{eq:u_ex_exp_by_int_phi}) from expression (\ref{eq:imp_local_sol_by_phi}), which leads to the following expression:
\begin{equation}
\begin{split}
\mathbf{q}_{n}^{\rm IL}(\tau) &- \mathbf{u}(t(\tau)) = \mathbf{u}_{n} - \mathbf{u}(t_{n}) - \mathbf{R}_{u}^{\Phi}(t(\tau))\\
&+ {\Delta t}_{n} \sum\limits_{p = 0}^{N} \left[
	\mathbf{F}\left(\mathbf{q}_{n, p}, t(\tau_{p})\right) -
	\mathbf{F}\left(\mathbf{u}(t(\tau_{p})), t(\tau_{p})\right)
\right]\int\limits_{0}^{\tau} \varphi_{p}(\tau) d\tau,
\end{split}
\end{equation}
using the absolute value of which leads to the following chain of inequalities:
\begin{equation}
\begin{split}
&|\mathbf{q}_{n}^{\rm IL}(\tau) - \mathbf{u}(t(\tau))| \leqslant
\left|\mathbf{u}_{n} - \mathbf{u}(t_{n})\right| + |\mathbf{R}_{u}^{\Phi}(t(\tau))|\\
&+ \left|{\Delta t}_{n} \sum\limits_{p = 0}^{N}
	\mathbf{F}\left(\mathbf{q}_{n, p}, t(\tau_{p})\right) -
	\mathbf{F}\left(\mathbf{u}(t(\tau_{p})), t(\tau_{p})\right)
\int\limits_{0}^{\tau} \varphi_{p}(\tau) d\tau\right|\\
&\leqslant \left|\mathbf{u}_{n} - \mathbf{u}(t_{n})\right| + |\mathbf{R}_{u}^{\Phi}(t(\tau))| +
{\Delta t}_{n} \sum\limits_{p = 0}^{N} \left|
	\mathbf{F}\left(\mathbf{q}_{n, p}, t(\tau_{p})\right) -
	\mathbf{F}\left(\mathbf{u}(t(\tau_{p})), t(\tau_{p})\right)
\right|\left|\int\limits_{0}^{\tau} \varphi_{p}(\tau) d\tau\right|\\
&\leqslant \left|\mathbf{u}_{n} - \mathbf{u}(t_{n})\right| + |\mathbf{R}_{u}^{\Phi}(t(\tau))| +
{\Delta t}_{n} \sum\limits_{p = 0}^{N} \left|\mathbf{q}_{n, p} - \mathbf{u}(t(\tau_{p}), t(\tau_{p}))\right|
\left|\int\limits_{0}^{\tau} \varphi_{p}(\tau) d\tau\right|\\
&\leqslant C_{u} {\Delta t}_{n}^{2N+1} + C_{R}^{\Phi} {\Delta t}_{n}^{N+2} + {\Delta t}_{n} \cdot \left[\sum\limits_{p = 0}^{N} B_{p} C_{L, p}^{q} {\Delta t}_{n}^{N+1}\right] \simeq
C_{\rm IL}^{q} {\Delta t}^{N+2},
\end{split}
\end{equation}
where estimate (\ref{eq:local_sol_errs_L_inf_est}) was used for the local solution $\mathbf{q}_{n}(\tau)$, and where the following notations and known estimates were used:
\begin{equation}
B_{p} = \max\limits_{\tau\in[0,1]} \left|\int\limits_{0}^{\tau} \varphi_{p}(\tau) d\tau\right|,\quad
|\mathbf{R}_{u}^{\Phi}(t(\tau))| \leqslant C_{R}^{\Phi} {\Delta t}_{n}^{N+2}.
\end{equation}
This allows us to conclude that the approximation order $p_{\rm IL}$ for the improved local solution in the discretization domain is $N+2$. From this error estimate, a similar error estimate for the local solution $\mathbf{u}_{\rm IL}(t)$ (\ref{eq:imp_local_sol_assembly}) follows immediately:
\begin{equation}
\begin{split}
||\mathbf{u}_{\rm IL}(t) &- \mathbf{u}(t)||_{\mathcal{L}_{\infty}} =
\operatorname{ess}\sup\limits_{\hspace{-5mm}t\in\Omega}\left|\mathbf{u}_{\rm IL}(t) - \mathbf{u}(t)\right|\\
&= \max\limits_{n}\sup\limits_{t\in\Omega_{n}}\left|\mathbf{u}_{\rm IL}(t) - \mathbf{u}(t)\right| =
\max\limits_{n}\sup\limits_{t\in\Omega_{n}}\left|\mathbf{q}_{n}^{\rm IL}(t) - \mathbf{u}(t)\right|\\
&= \max\limits_{n}\max\limits_{t\in\Omega_{n}}|\mathbf{q}_{n}^{\rm IL}(\tau) - \mathbf{u}(t(\tau))| \leqslant C_{\rm IL}^{\mathcal{L}_{\infty}} {\Delta t}^{N+2}.
\end{split}
\end{equation}
Similar error estimates for the improved local numerical solution $\mathbf{u}_{\rm IL}(t)$ (\ref{eq:imp_local_sol_assembly}) are also obtained for other classical error norms in $\mathcal{L}_{p}$ ($p \geqslant 1$):
\begin{equation}
\begin{split}
||\mathbf{u}_{\rm IL}(t) - \mathbf{u}(t)||_{\mathcal{L}_{p}} &=
\left[\int\limits_{\Omega} |\mathbf{u}_{\rm IL}(t) - \mathbf{u}(t)|^{p} dt \right]^{1/p} \leqslant
\left[|\Omega| \cdot \left(\operatorname{ess}\sup\limits_{\hspace{-5mm}t\in\Omega}\left|\mathbf{u}_{\rm IL}(t) - \mathbf{u}(t)\right|\right)^{p}\,\right]^{1/p}\\
&= |\Omega|^{1/p} ||\mathbf{u}_{\rm IL}(t) - \mathbf{u}(t)||_{\mathcal{L}_{\infty}} \leqslant
|\Omega|^{1/p} C_{\rm IL}^{\mathcal{L}_{\infty}} {\Delta t}^{N+2} \equiv C_{\rm IL}^{\mathcal{L}_{p}} {\Delta t}^{N+2}.
\end{split}
\end{equation}
The obtained rigorous proofs show that the improved local solution $\mathbf{u}_{\rm IL}(t)$ (\ref{eq:local_sol_assembly}), obtained based on the use of the local solution $\mathbf{q}_{n}(\tau)$ obtained by the local DG predictor of the ADER-DG numerical method, has the approximation order $p_{\rm IL} = N+2$ in the norms of $\mathcal{L}_{\infty}$ and $\mathcal{L}_{p}$ ($p \geqslant 1$).

It is important to note that the resulting improved local solution $\mathbf{u}_{\rm IL}(t)$ (\ref{eq:local_sol_assembly}) can be used to obtain a numerical solution in the space between grid nodes $t\in\Omega_{n}$, and its use does not affect the process of obtaining a numerical solution $\mathbf{u}_{n}$ at the grid nodes $t_{n}$. Therefore, the fundamental properties of the ADER-DG numerical method with a local DG predictor, identified in the previous studies~\cite{ader_dg_ode_jsc, ader_dg_ode_sinum}, remain unchanged; in particular, the method retains its high stability.

It should also be noted that calculating an improved numerical solution $\mathbf{u}_{\rm IL}(t)$ (\ref{eq:local_sol_assembly}) requires no additional computational costs other than the computation of the improved numerical solution values $\mathbf{u}_{\rm IL}(t_{i})$ at the required points $\{t_{n, i} = t_{n} + \tau_{i}{\Delta t}_{n}\}$. Technically, this can be accomplished using classical matrix-matrix operations, which can be used to calculate an improved numerical solution $\mathbf{u}_{\rm IL}(t)$ at a predetermined set of points $\{\tau_{i}\}$ in the discretization domain $\Omega_{n}$. This requires pre-computing the basis function $\int_{0}^{\tau} \varphi_{p}(\xi) d\xi$ values at the given set of points $\{\tau_{i}\}$, then multiplying the resulting matrix by the coefficient matrix $\hat{\mathbf{q}}_{p}$ and adding $\mathbf{u}_{n}$.

A demonstration of the qualitative properties of the improved local solution is performed based on the numerical solution of a simple initial value problem for a second-order ODE describing a one-dimensional linear harmonic oscillator:
\begin{equation}\label{eq:demo_ode}
\begin{split}
&\ddot{x} + x = 0,\quad
x(0) = 1,\quad \dot{x}(0) = 0,\quad
t \in [0,\, 8\pi],\\
&x^{\rm ex}(t) = \cos(t),\quad
\dot{x}^{\rm ex}(t) = -\sin(t),
\end{split}
\end{equation}
for which the solution vectors in the original notation of the ODE system (\ref{eq:ivp_ode_diff_src}), with $D = 2$, are chosen in the form $\mathbf{u}(t) = [x(t)\ \dot{x}(t)]^{T}$, $\mathbf{u}^{\rm ex}(t) = [x^{\rm ex}(t)\ \dot{x}^{\rm ex}(t)]^{T} = [\cos(t)\ -\sin(t)]^{T}$.

\begin{figure}[h!]
\captionsetup[subfigure]{%
	position=bottom,
	font+=smaller,
	textfont=normalfont,
	singlelinecheck=off,
	justification=raggedright
}
\centering
\begin{subfigure}{0.24\textwidth}
\includegraphics[width=\textwidth]{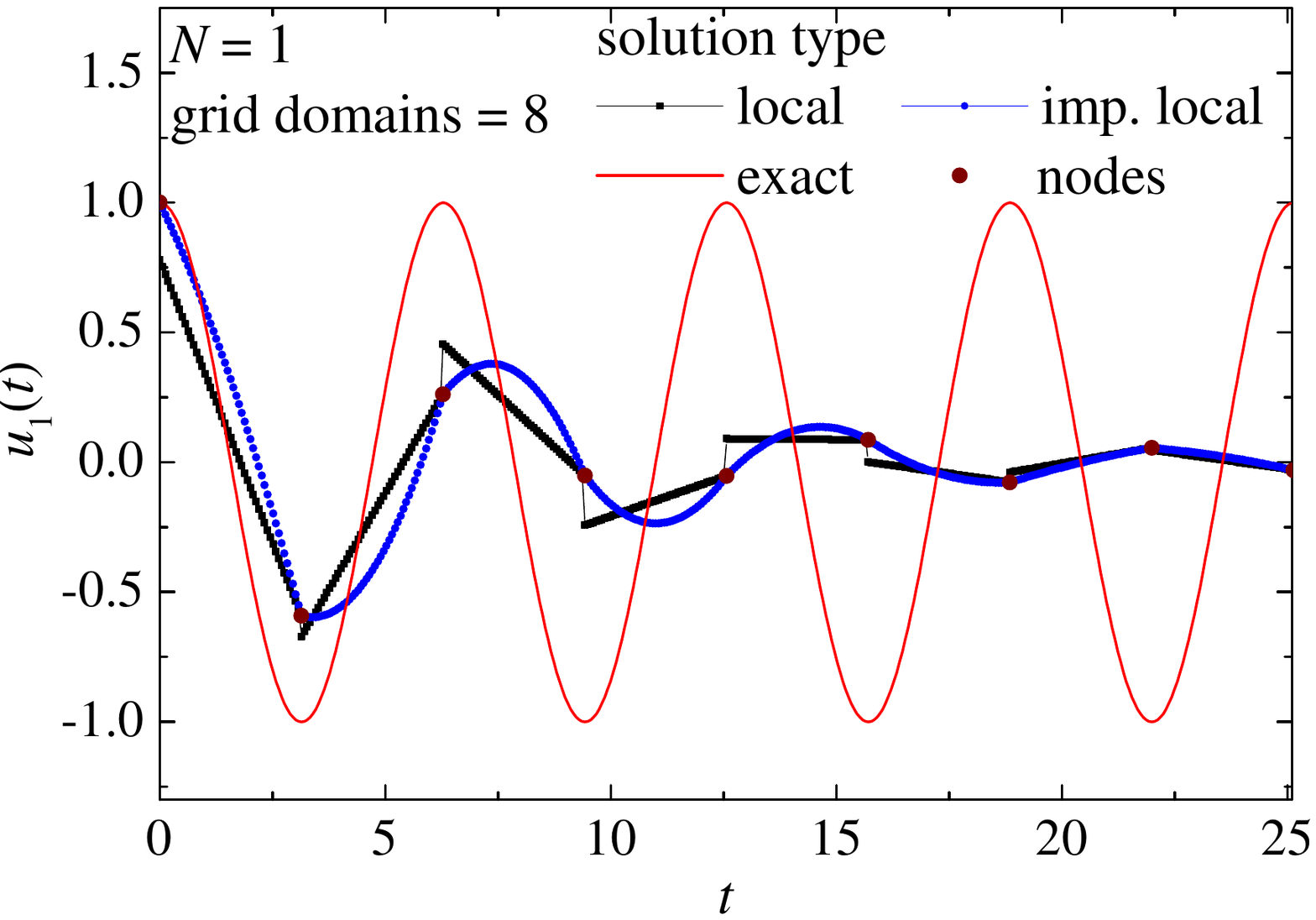}
\vspace{-8mm}\caption{\label{fig:demo_nodes_8:a1}}
\end{subfigure}
\begin{subfigure}{0.24\textwidth}
\includegraphics[width=\textwidth]{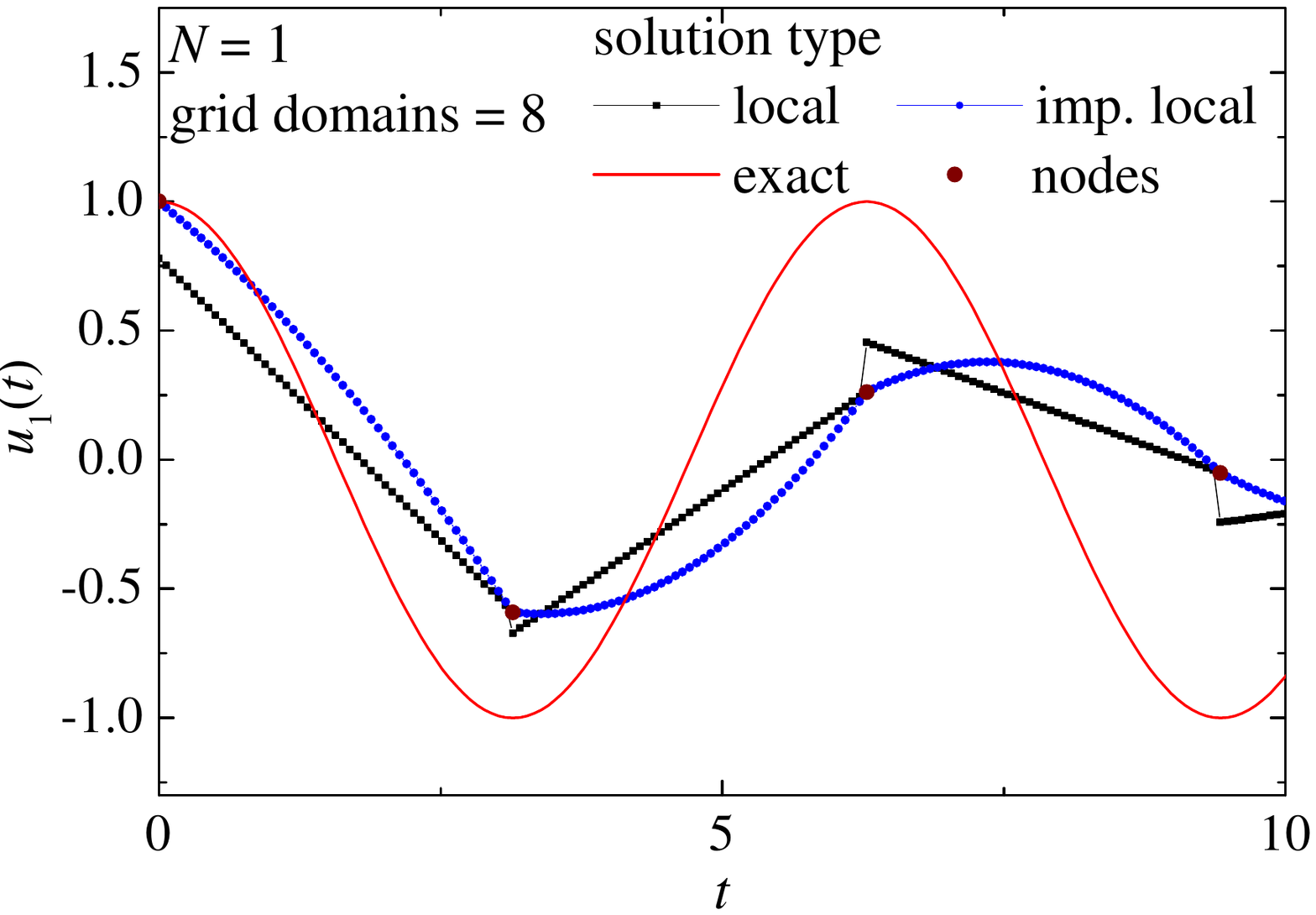}
\vspace{-8mm}\caption{\label{fig:demo_nodes_8:a2}}
\end{subfigure}
\begin{subfigure}{0.24\textwidth}
\includegraphics[width=\textwidth]{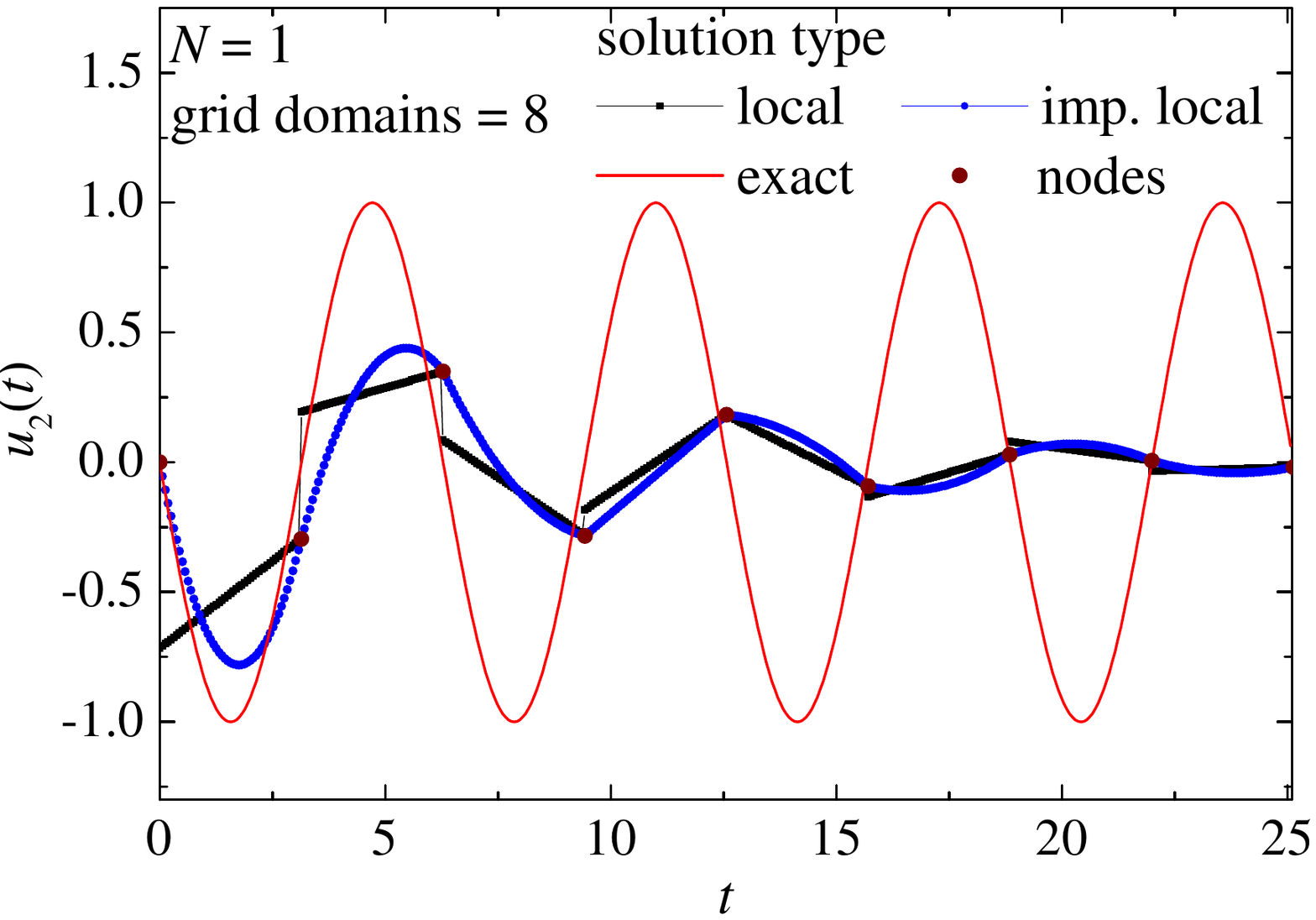}
\vspace{-8mm}\caption{\label{fig:demo_nodes_8:a3}}
\end{subfigure}
\begin{subfigure}{0.24\textwidth}
\includegraphics[width=\textwidth]{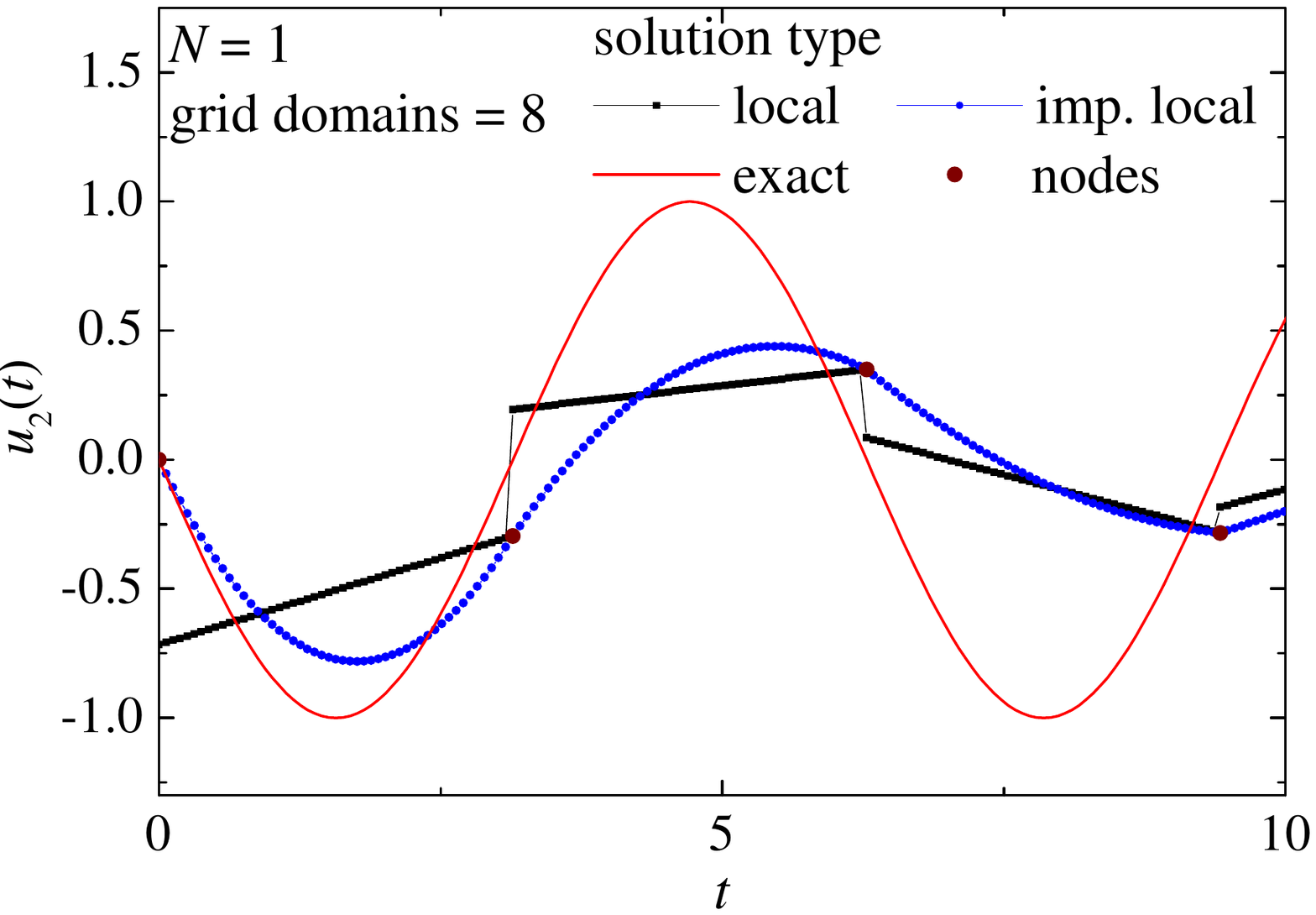}
\vspace{-8mm}\caption{\label{fig:demo_nodes_8:a4}}
\end{subfigure}\\[-2mm]
\begin{subfigure}{0.24\textwidth}
\includegraphics[width=\textwidth]{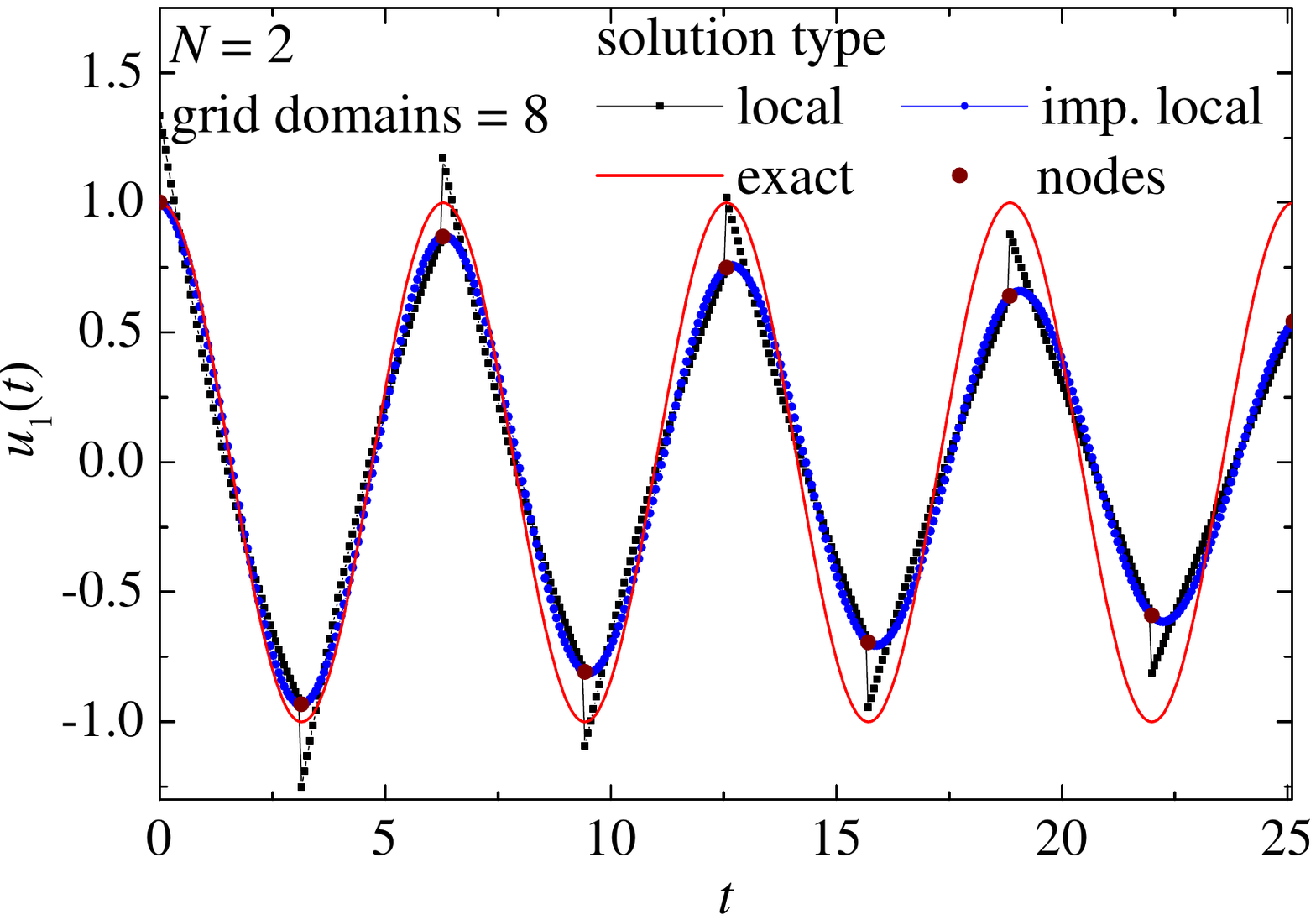}
\vspace{-8mm}\caption{\label{fig:demo_nodes_8:b1}}
\end{subfigure}
\begin{subfigure}{0.24\textwidth}
\includegraphics[width=\textwidth]{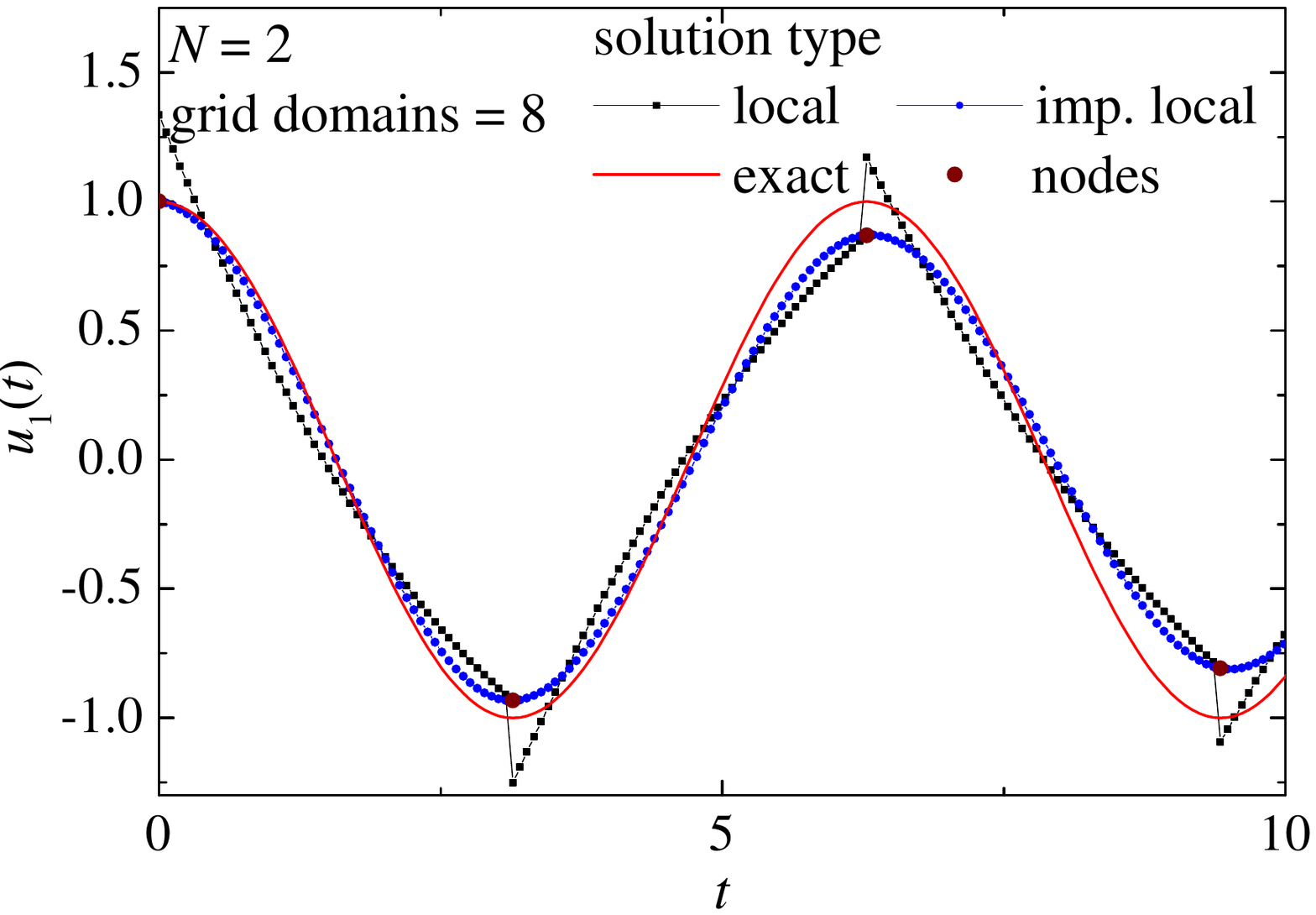}
\vspace{-8mm}\caption{\label{fig:demo_nodes_8:b2}}
\end{subfigure}
\begin{subfigure}{0.24\textwidth}
\includegraphics[width=\textwidth]{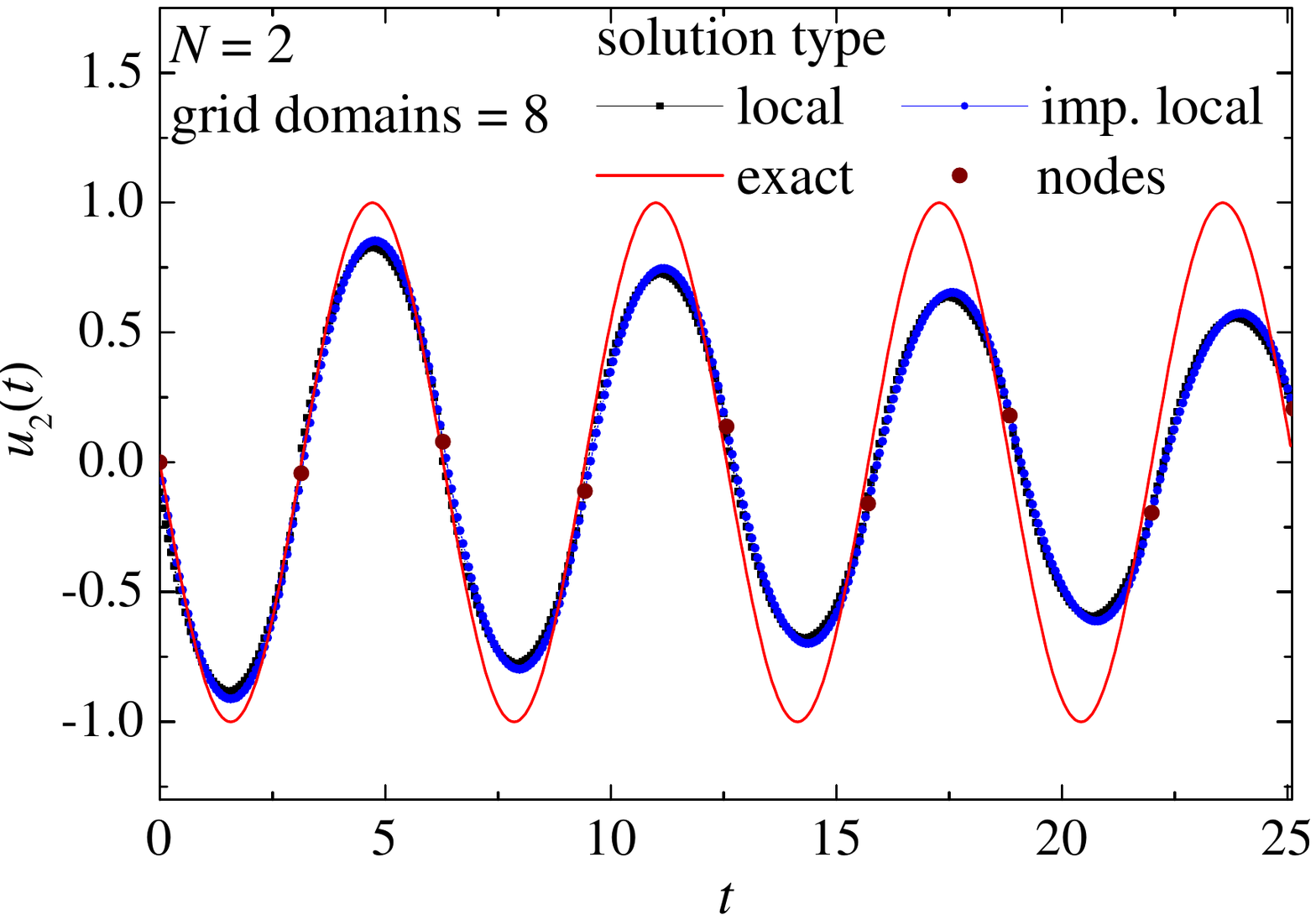}
\vspace{-8mm}\caption{\label{fig:demo_nodes_8:b3}}
\end{subfigure}
\begin{subfigure}{0.24\textwidth}
\includegraphics[width=\textwidth]{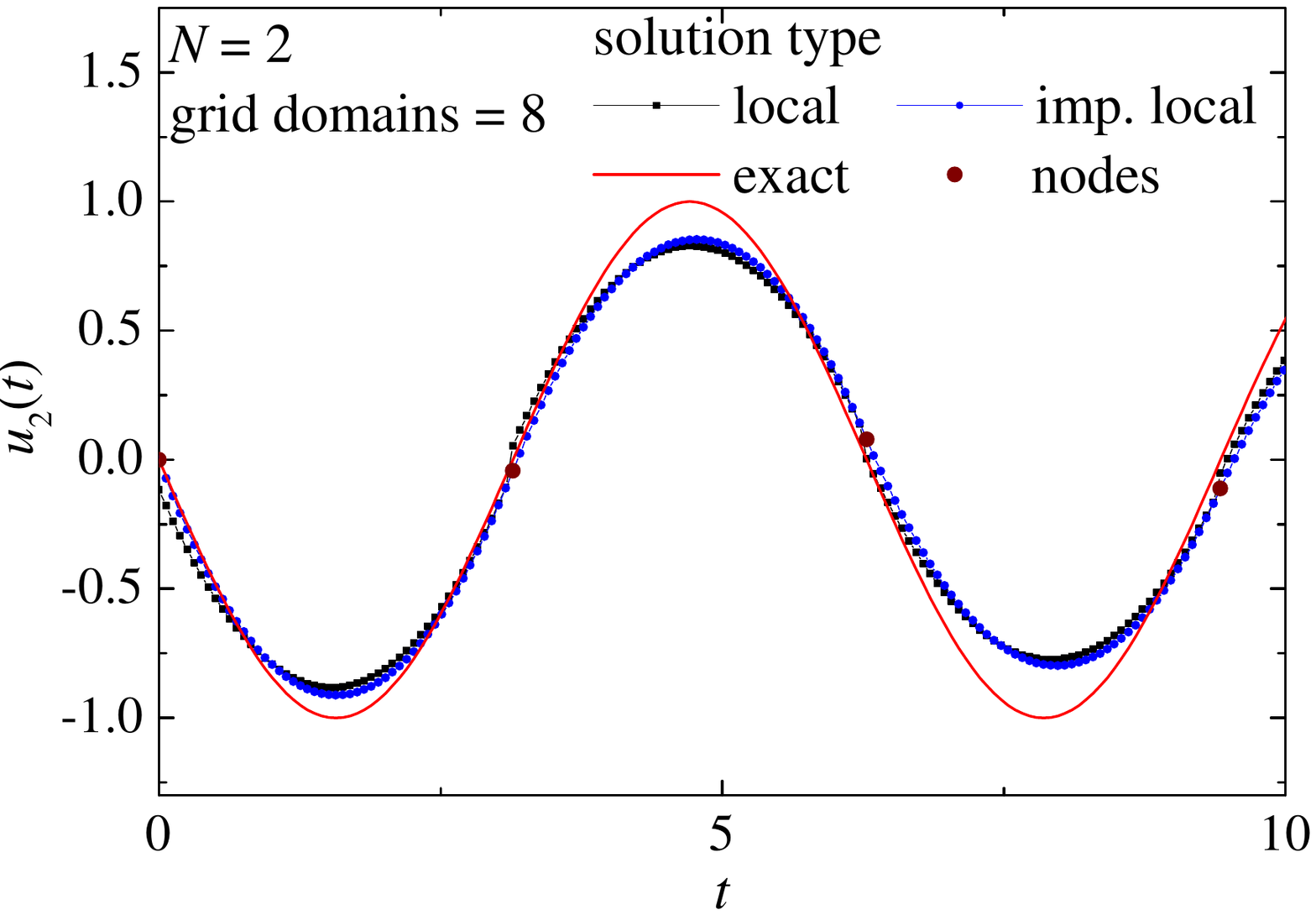}
\vspace{-8mm}\caption{\label{fig:demo_nodes_8:b4}}
\end{subfigure}\\[-2mm]
\begin{subfigure}{0.24\textwidth}
\includegraphics[width=\textwidth]{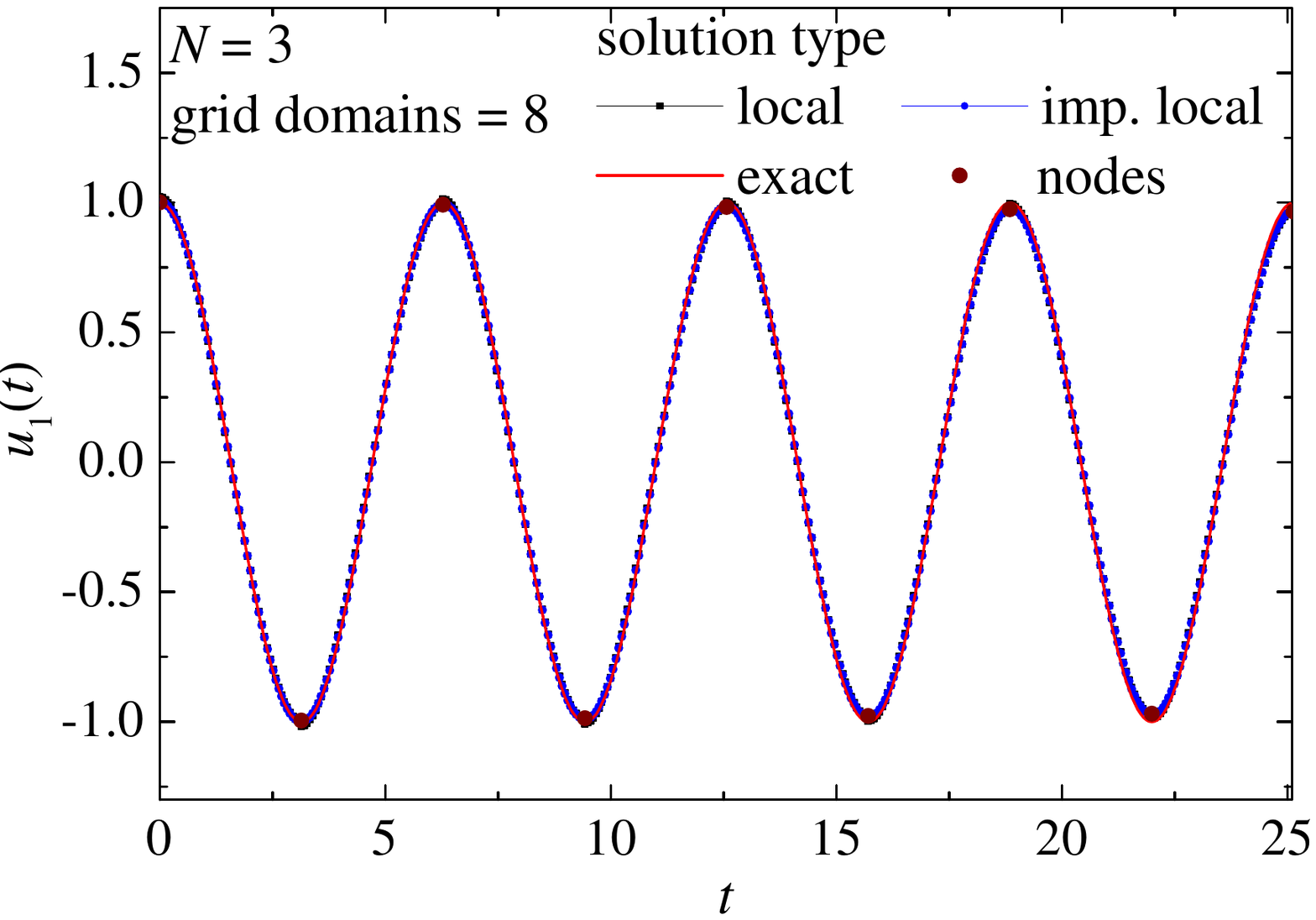}
\vspace{-8mm}\caption{\label{fig:demo_nodes_8:c1}}
\end{subfigure}
\begin{subfigure}{0.24\textwidth}
\includegraphics[width=\textwidth]{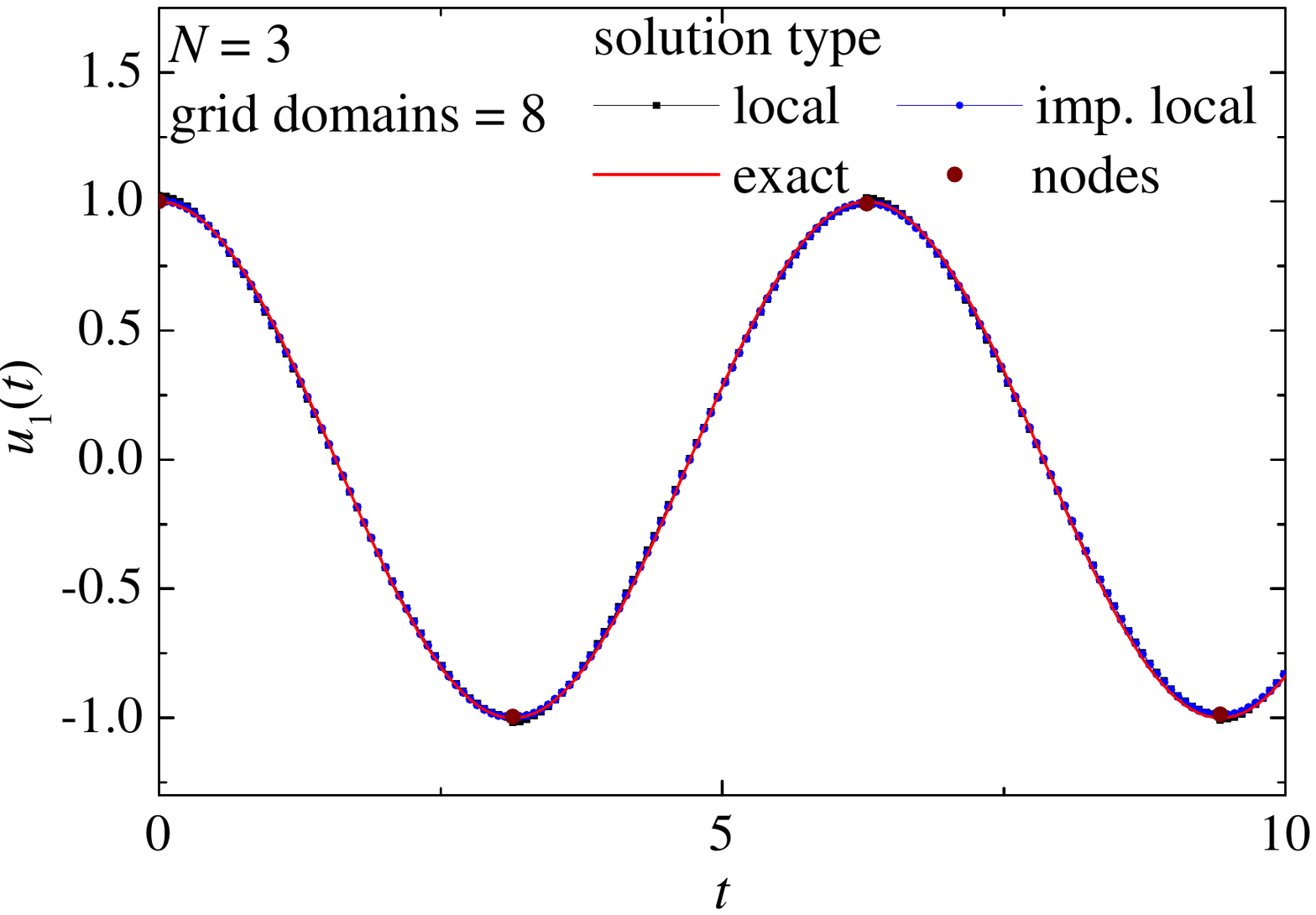}
\vspace{-8mm}\caption{\label{fig:demo_nodes_8:c2}}
\end{subfigure}
\begin{subfigure}{0.24\textwidth}
\includegraphics[width=\textwidth]{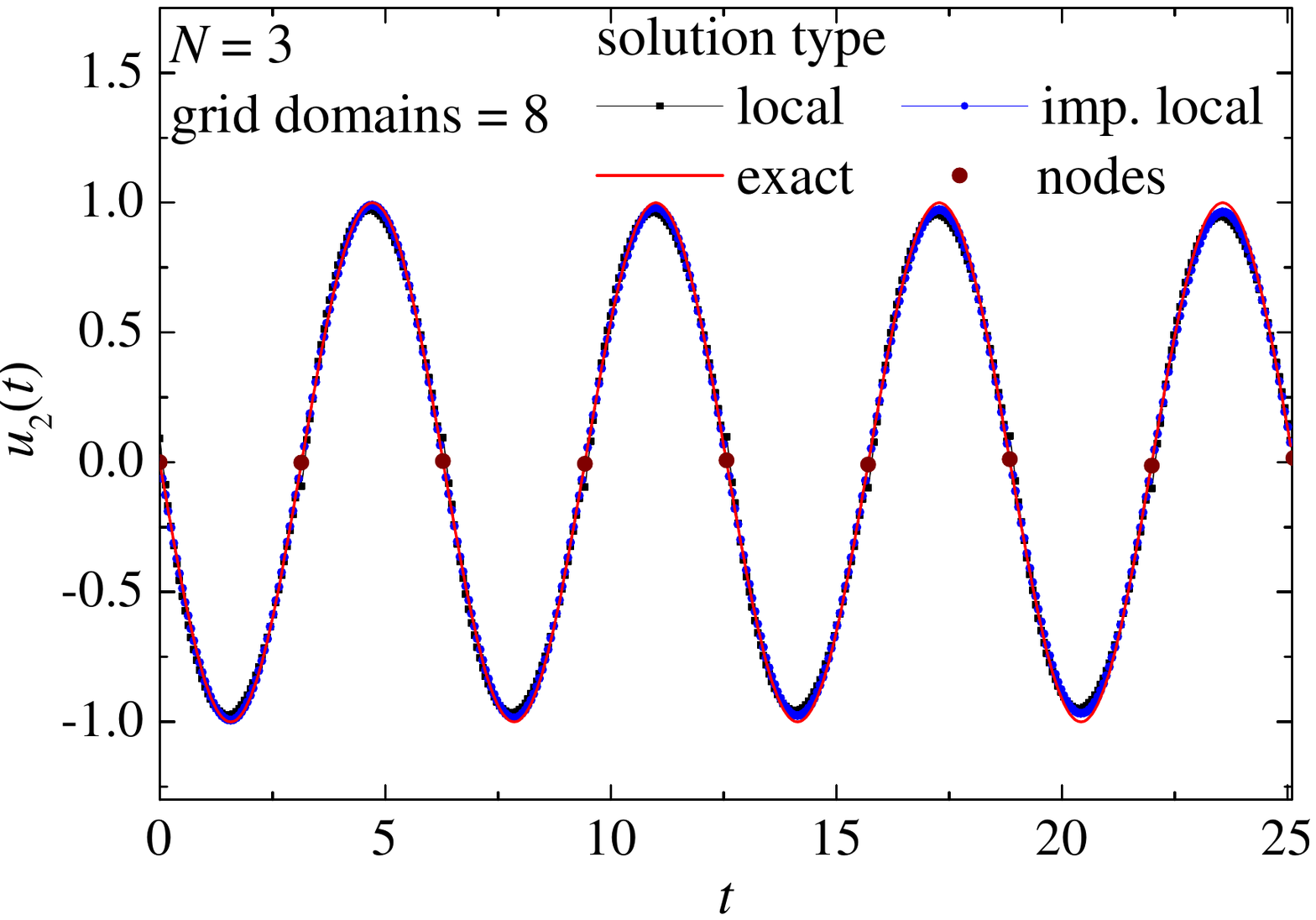}
\vspace{-8mm}\caption{\label{fig:demo_nodes_8:c3}}
\end{subfigure}
\begin{subfigure}{0.24\textwidth}
\includegraphics[width=\textwidth]{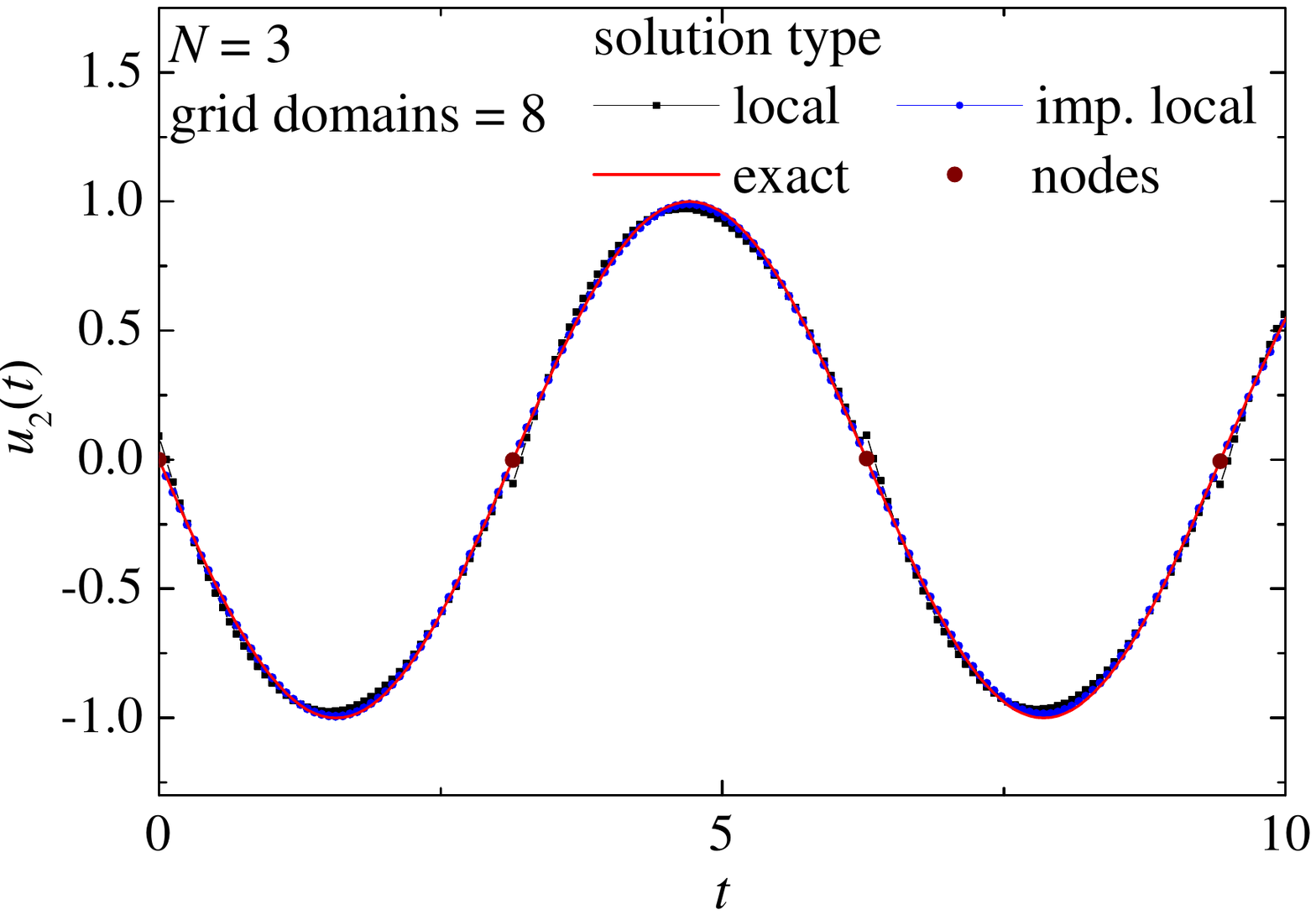}
\vspace{-8mm}\caption{\label{fig:demo_nodes_8:c4}}
\end{subfigure}\\[-2mm]
\caption{%
Numerical solution of the problem (\ref{eq:demo_ode}), in the domain $0 \leqslant t \leqslant 8\pi$ with step $\mathrm{\Delta}t = \pi$. Comparison of the solution at nodes $\mathbf{u}_{n}$, the local solution $\mathbf{u}_{L}(t)$, the improved local solution $\mathbf{u}_{IL}(t)$ and the exact solution $\mathbf{u}^{\rm ex}(t)$ for components $u_{1} \equiv x$ (\subref{fig:demo_nodes_8:a1}, \subref{fig:demo_nodes_8:a2}, \subref{fig:demo_nodes_8:b1}, \subref{fig:demo_nodes_8:b2}, \subref{fig:demo_nodes_8:c1}, \subref{fig:demo_nodes_8:c2}) and $u_{2} \equiv \dot{x}$ (\subref{fig:demo_nodes_8:a3}, \subref{fig:demo_nodes_8:a4}, \subref{fig:demo_nodes_8:b3}, \subref{fig:demo_nodes_8:b4}, \subref{fig:demo_nodes_8:c3}, \subref{fig:demo_nodes_8:c4}), obtained using polynomials with degrees $N = 1$ (\subref{fig:demo_nodes_8:a1}, \subref{fig:demo_nodes_8:a2}, \subref{fig:demo_nodes_8:a3}, \subref{fig:demo_nodes_8:a4}), $N = 2$ (\subref{fig:demo_nodes_8:b1}, \subref{fig:demo_nodes_8:b2}, \subref{fig:demo_nodes_8:b3}, \subref{fig:demo_nodes_8:b4}) and $N = 3$ (\subref{fig:demo_nodes_8:c1}, \subref{fig:demo_nodes_8:c2}, \subref{fig:demo_nodes_8:c3}, \subref{fig:demo_nodes_8:c4}). Zoomed domain $0 \leqslant t \leqslant 10$ is presented on (\subref{fig:demo_nodes_8:a2}, \subref{fig:demo_nodes_8:a4}, \subref{fig:demo_nodes_8:b2}, \subref{fig:demo_nodes_8:b4}, \subref{fig:demo_nodes_8:c2}, \subref{fig:demo_nodes_8:c4}).
}
\label{fig:demo_nodes_8}
\end{figure}

\begin{figure}[h!]
\captionsetup[subfigure]{%
	position=bottom,
	font+=smaller,
	textfont=normalfont,
	singlelinecheck=off,
	justification=raggedright
}
\centering
\begin{subfigure}{0.24\textwidth}
\includegraphics[width=\textwidth]{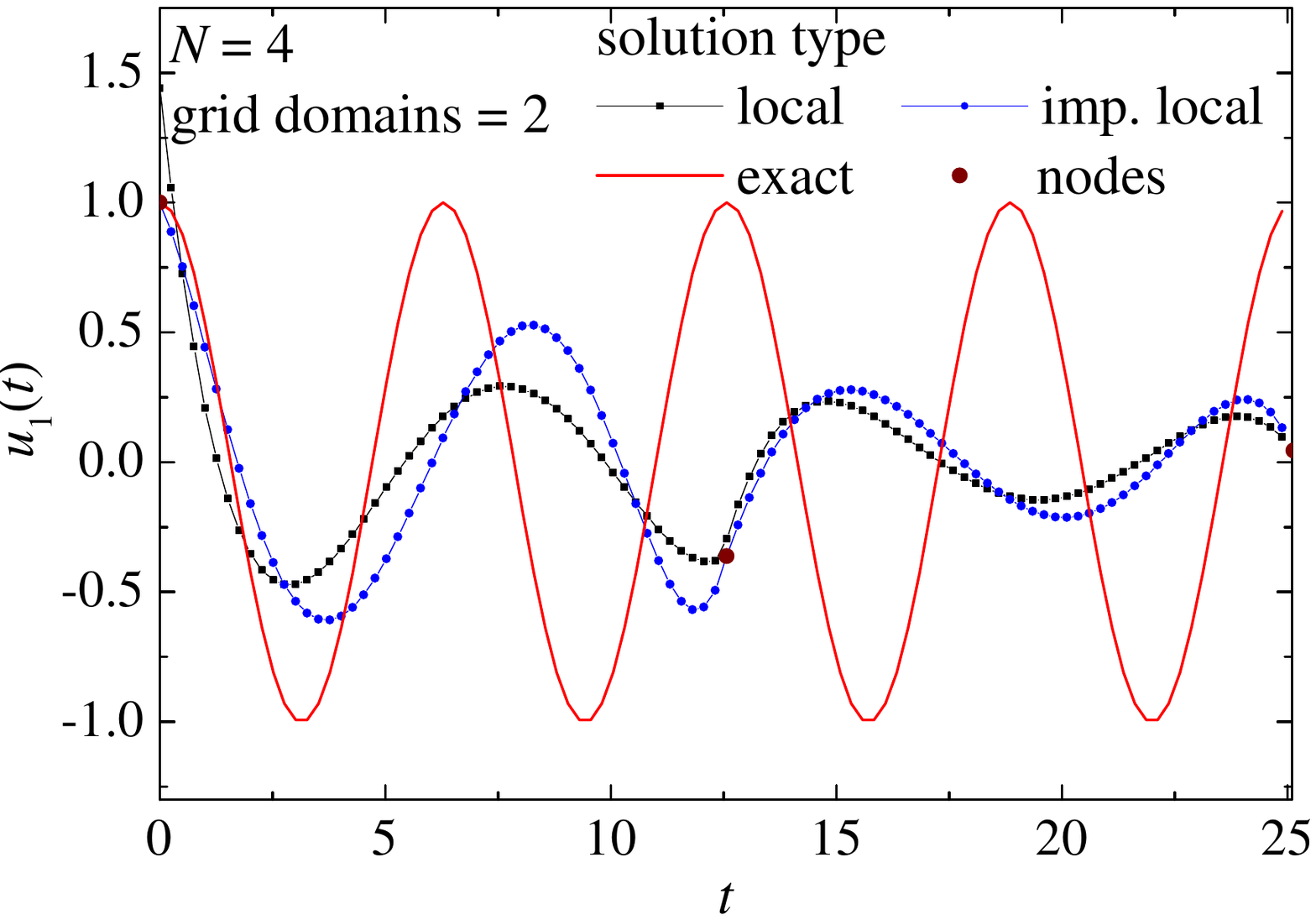}
\vspace{-8mm}\caption{\label{fig:demo_nodes_2:a1}}
\end{subfigure}
\begin{subfigure}{0.24\textwidth}
\includegraphics[width=\textwidth]{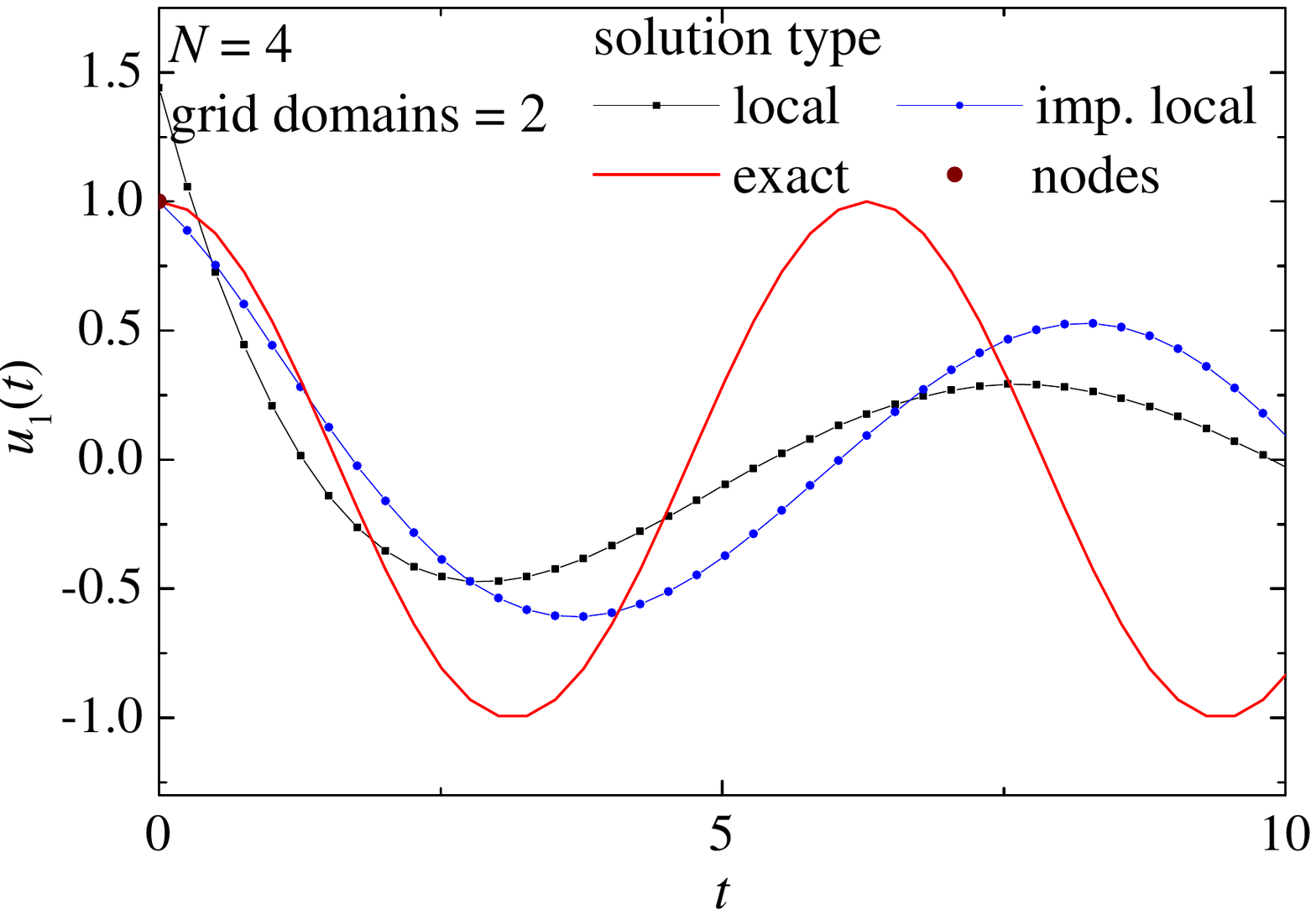}
\vspace{-8mm}\caption{\label{fig:demo_nodes_2:a2}}
\end{subfigure}
\begin{subfigure}{0.24\textwidth}
\includegraphics[width=\textwidth]{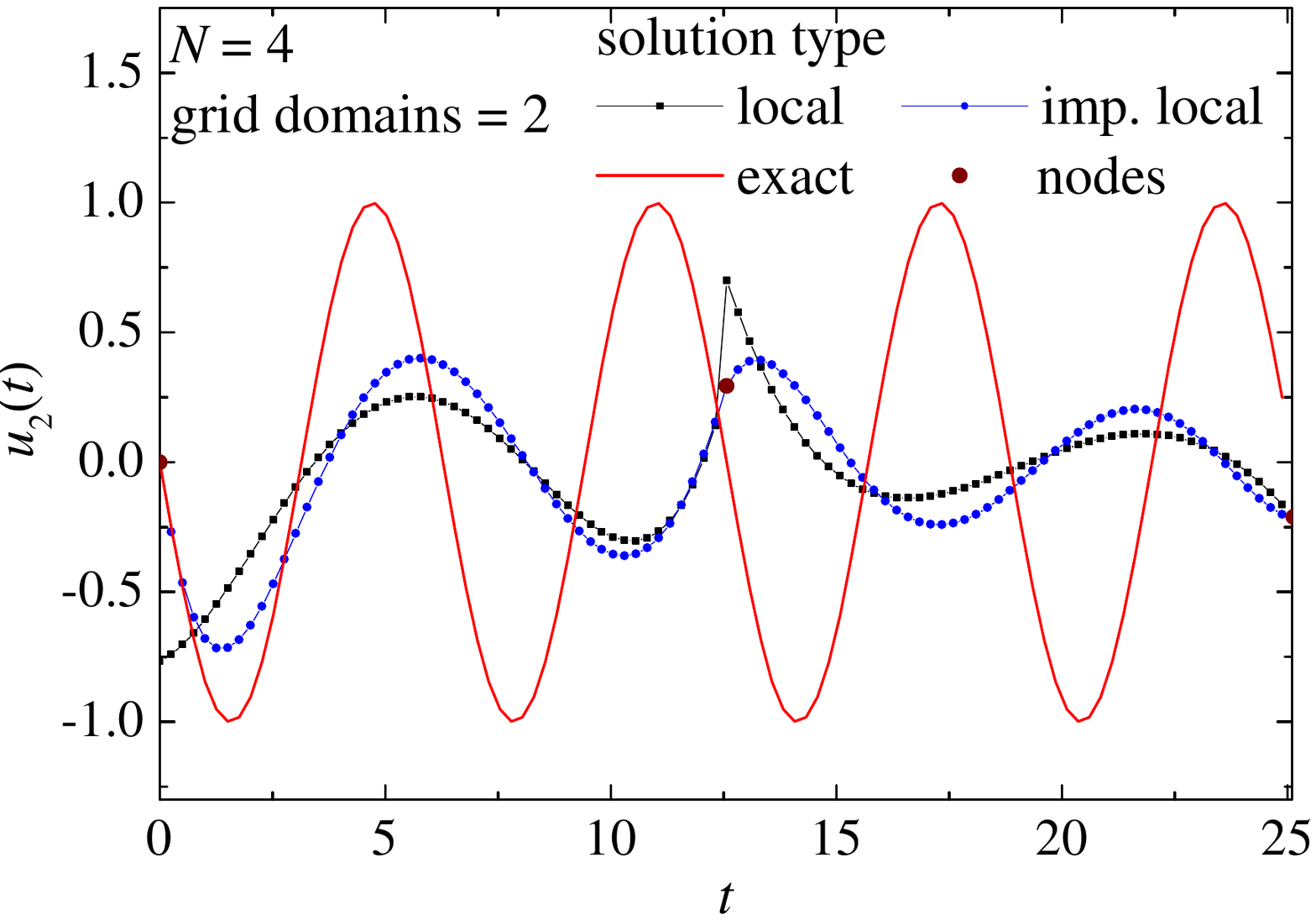}
\vspace{-8mm}\caption{\label{fig:demo_nodes_2:a3}}
\end{subfigure}
\begin{subfigure}{0.24\textwidth}
\includegraphics[width=\textwidth]{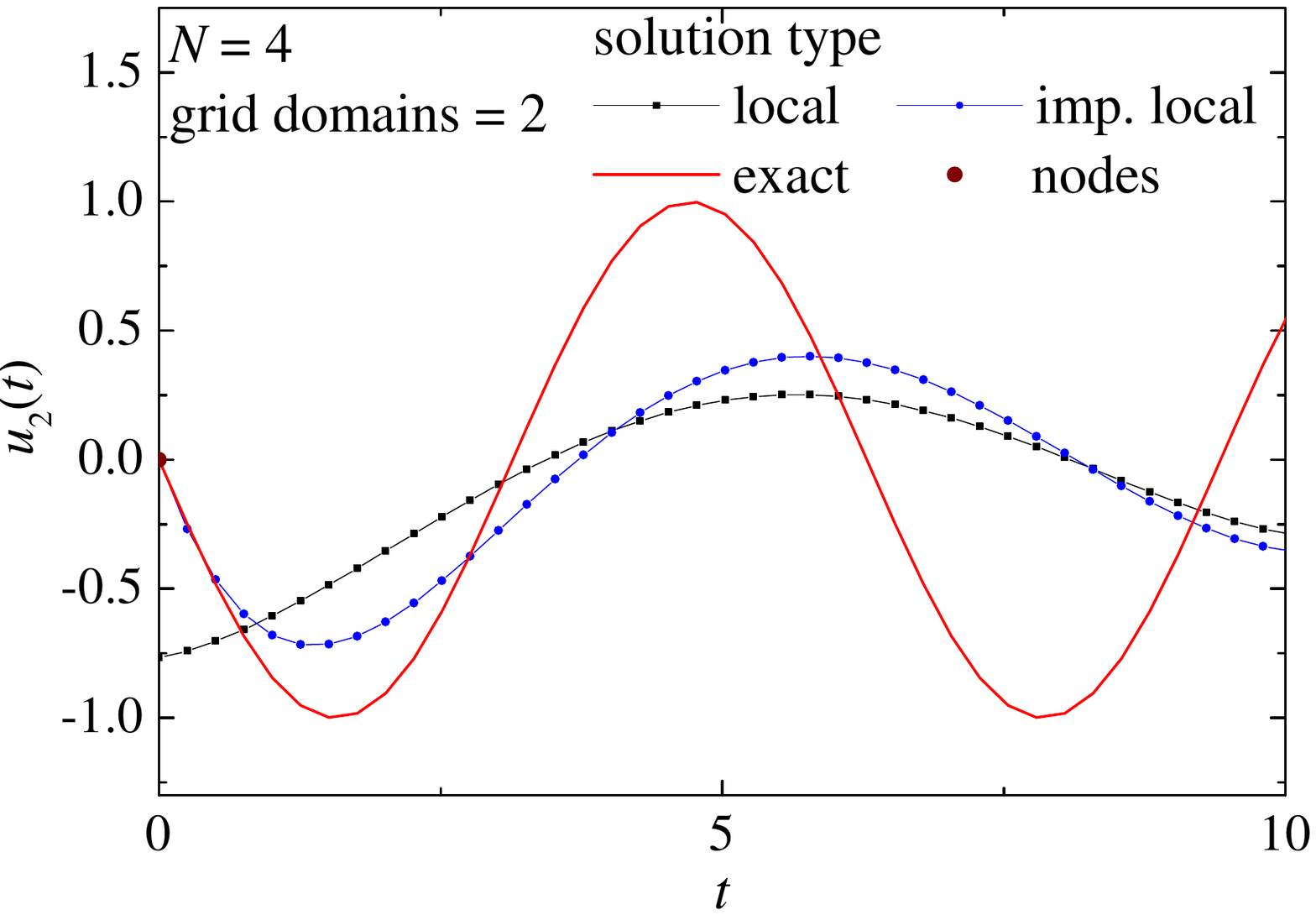}
\vspace{-8mm}\caption{\label{fig:demo_nodes_2:a4}}
\end{subfigure}\\[-2mm]
\begin{subfigure}{0.24\textwidth}
\includegraphics[width=\textwidth]{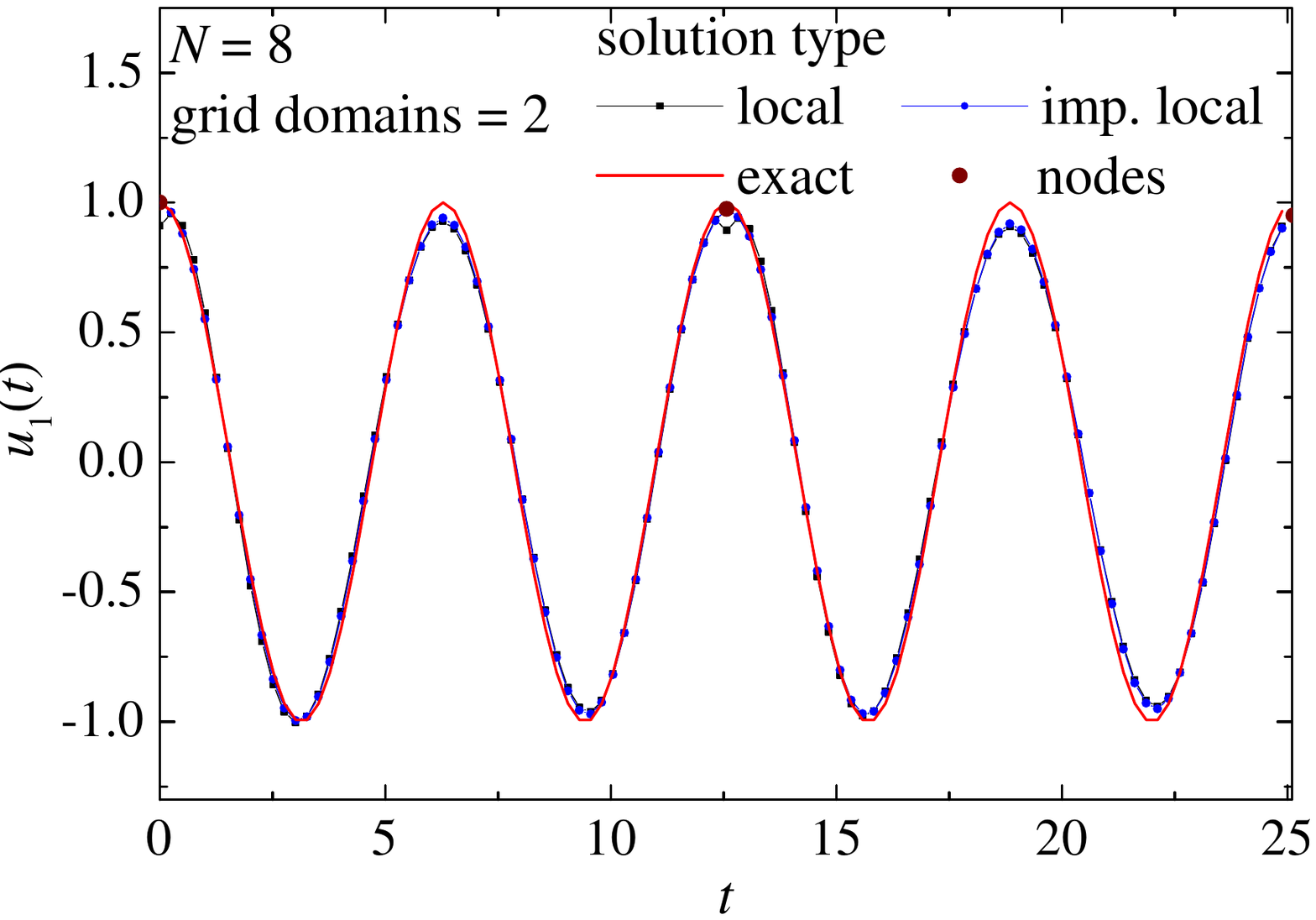}
\vspace{-8mm}\caption{\label{fig:demo_nodes_2:b1}}
\end{subfigure}
\begin{subfigure}{0.24\textwidth}
\includegraphics[width=\textwidth]{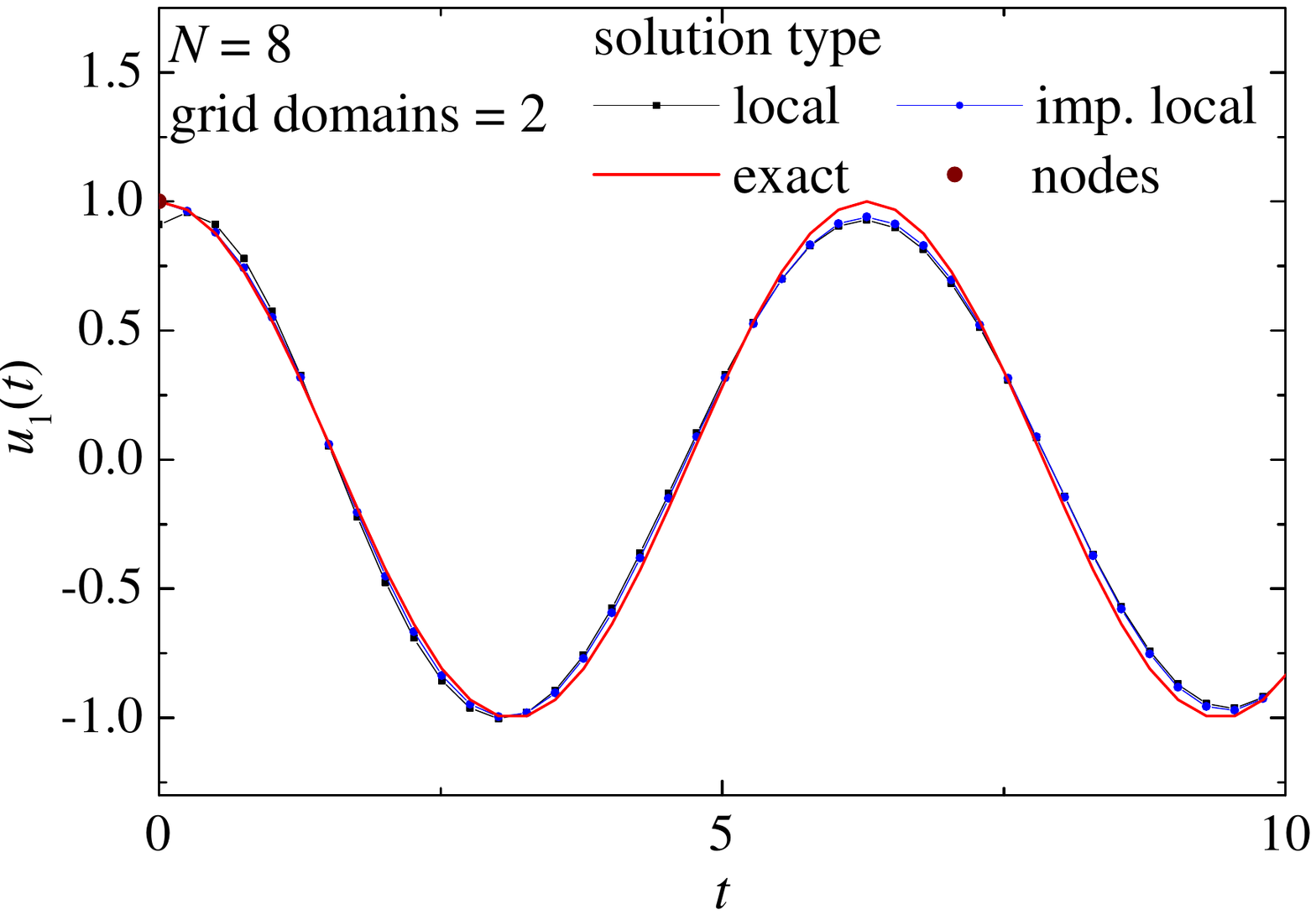}
\vspace{-8mm}\caption{\label{fig:demo_nodes_2:b2}}
\end{subfigure}
\begin{subfigure}{0.24\textwidth}
\includegraphics[width=\textwidth]{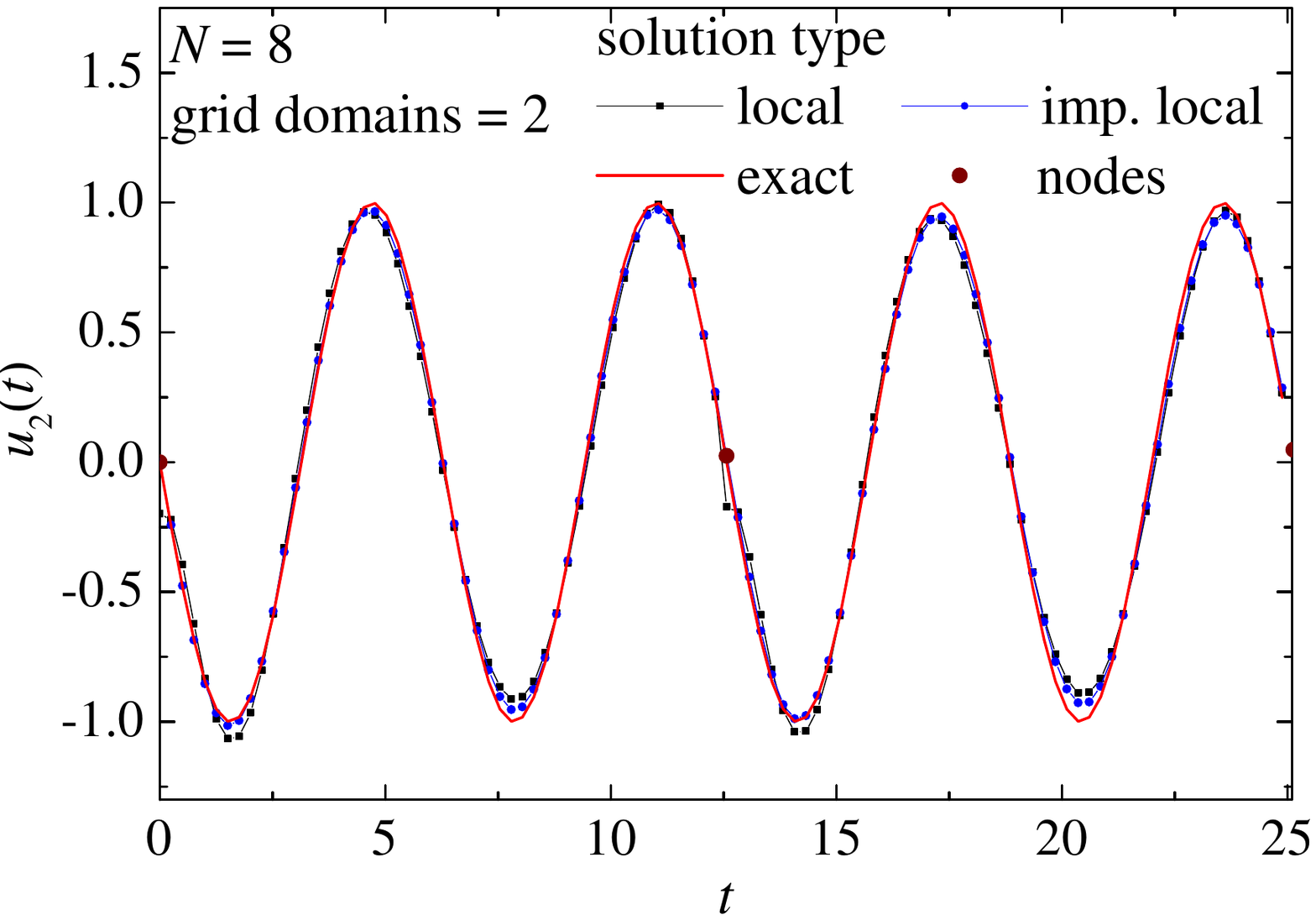}
\vspace{-8mm}\caption{\label{fig:demo_nodes_2:b3}}
\end{subfigure}
\begin{subfigure}{0.24\textwidth}
\includegraphics[width=\textwidth]{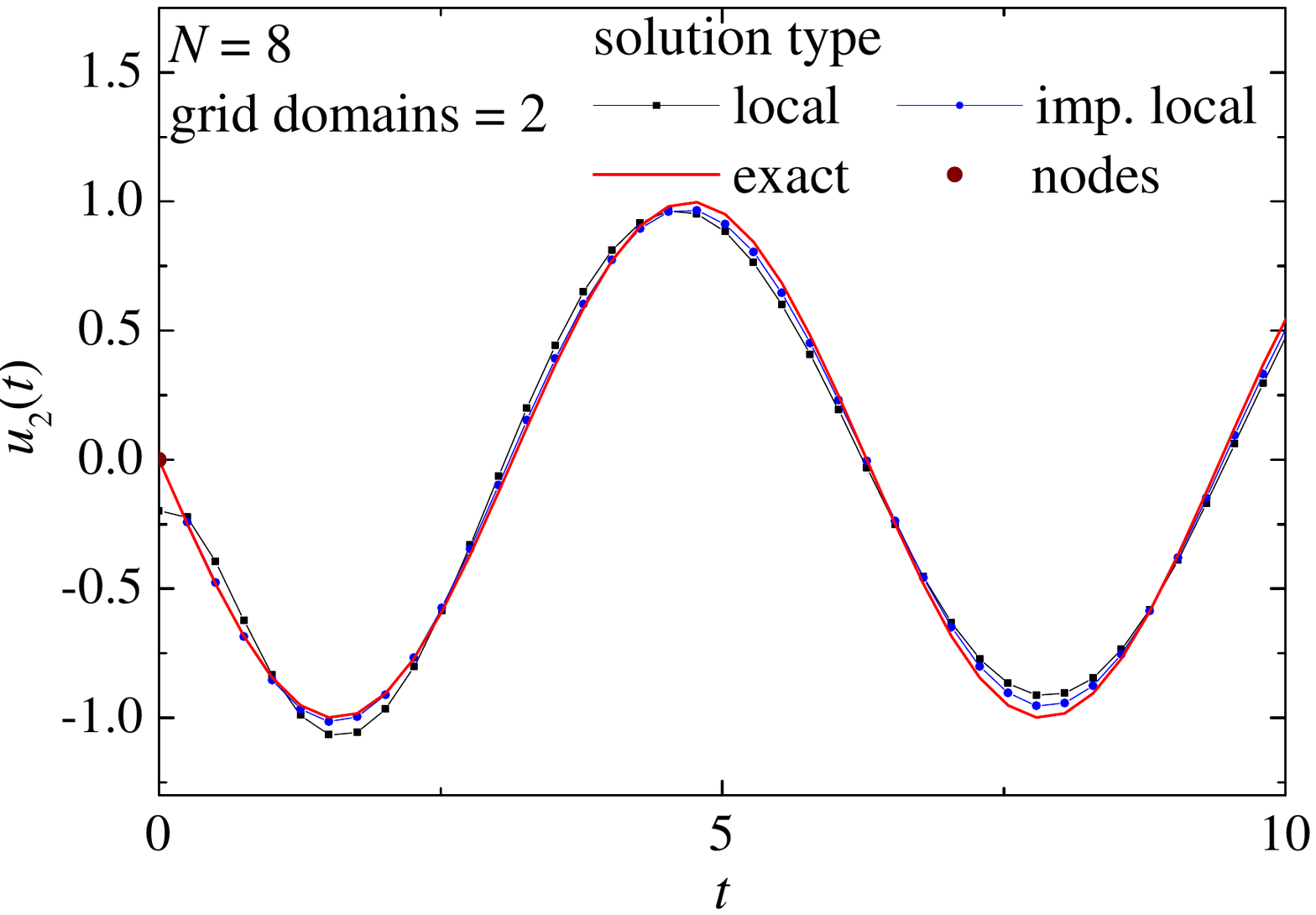}
\vspace{-8mm}\caption{\label{fig:demo_nodes_2:b4}}
\end{subfigure}\\[-2mm]
\begin{subfigure}{0.24\textwidth}
\includegraphics[width=\textwidth]{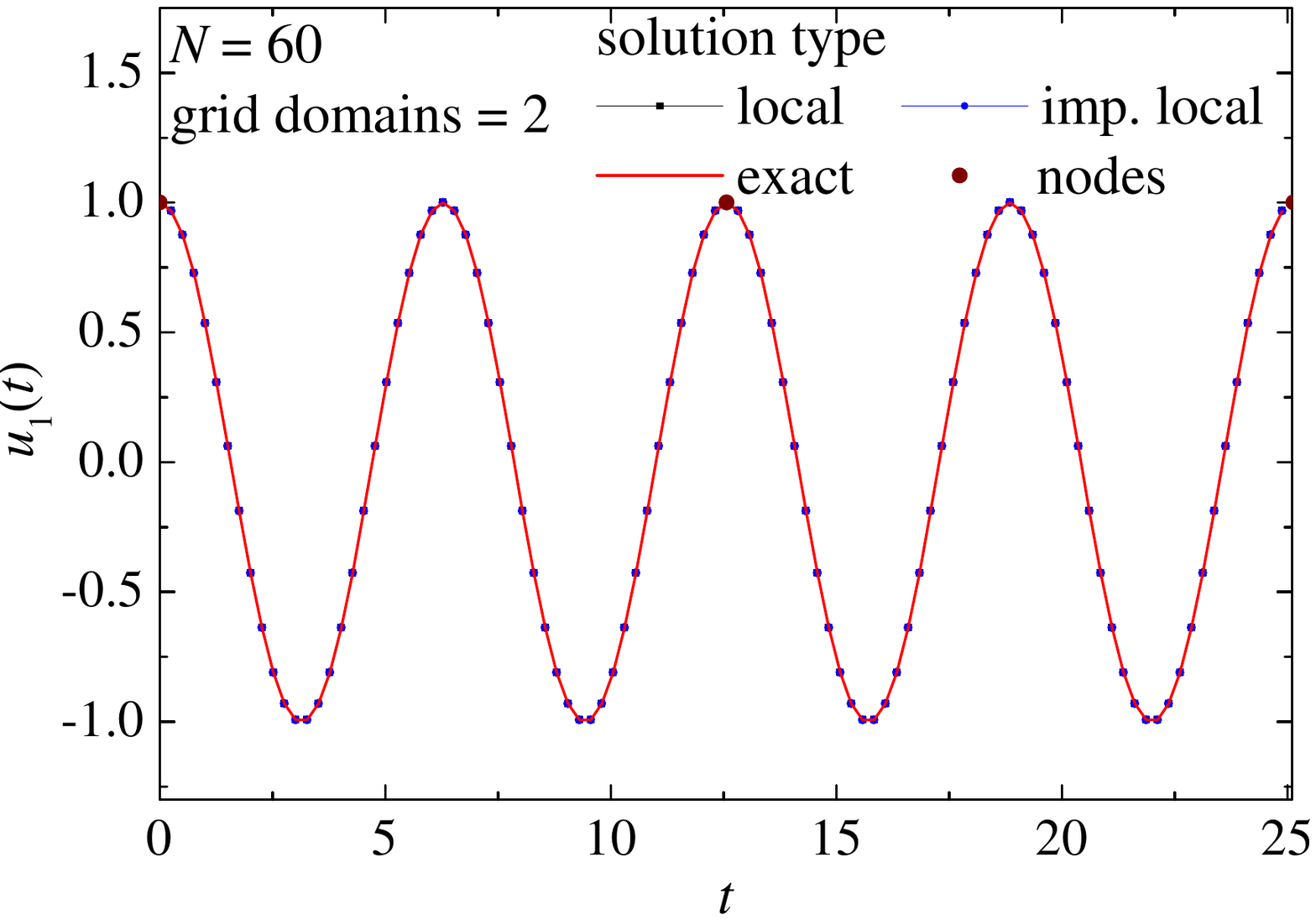}
\vspace{-8mm}\caption{\label{fig:demo_nodes_2:c1}}
\end{subfigure}
\begin{subfigure}{0.24\textwidth}
\includegraphics[width=\textwidth]{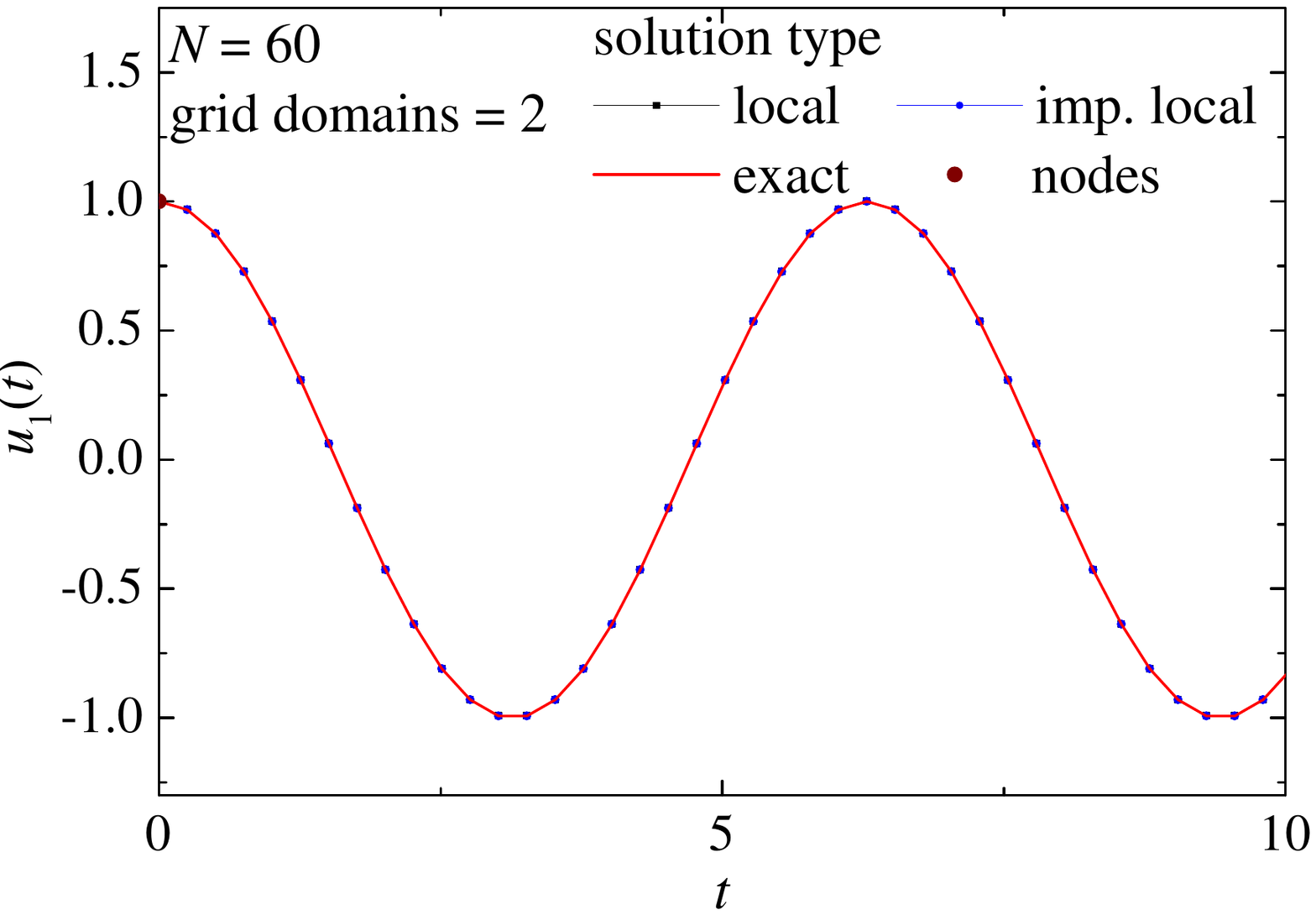}
\vspace{-8mm}\caption{\label{fig:demo_nodes_2:c2}}
\end{subfigure}
\begin{subfigure}{0.24\textwidth}
\includegraphics[width=\textwidth]{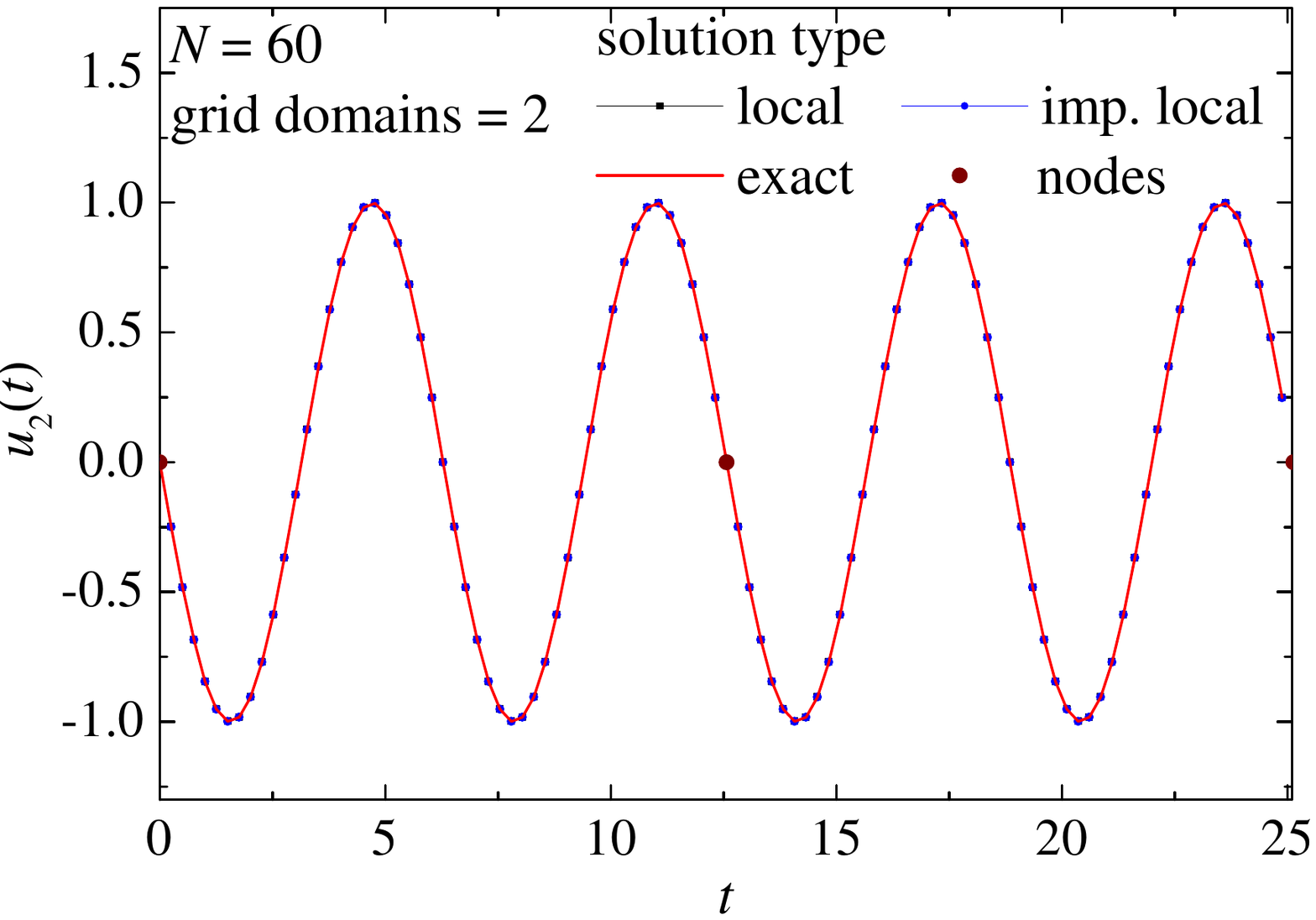}
\vspace{-8mm}\caption{\label{fig:demo_nodes_2:c3}}
\end{subfigure}
\begin{subfigure}{0.24\textwidth}
\includegraphics[width=\textwidth]{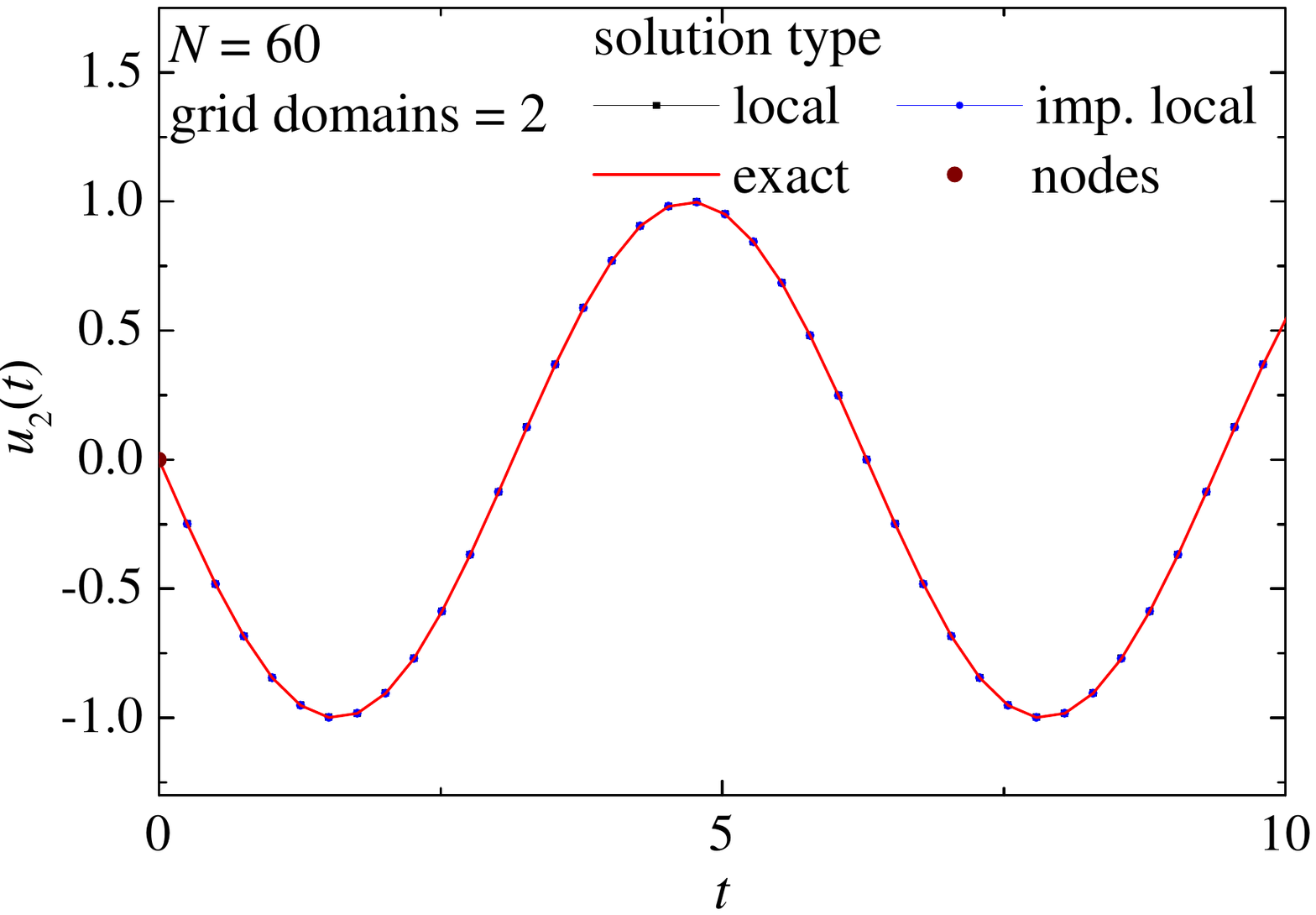}
\vspace{-8mm}\caption{\label{fig:demo_nodes_2:c4}}
\end{subfigure}\\[-2mm]
\caption{%
Numerical solution of the problem (\ref{eq:demo_ode}), in the domain $0 \leqslant t \leqslant 8\pi$ with step $\mathrm{\Delta}t = 4\pi$. Comparison of the solution at nodes $\mathbf{u}_{n}$, the local solution $\mathbf{u}_{L}(t)$, the improved local solution $\mathbf{u}_{\rm IL}(t)$ and the exact solution $\mathbf{u}^{\rm ex}(t)$ for components $u_{1} \equiv x$ (\subref{fig:demo_nodes_2:a1}, \subref{fig:demo_nodes_2:a2}, \subref{fig:demo_nodes_2:b1}, \subref{fig:demo_nodes_2:b2}, \subref{fig:demo_nodes_2:c1}, \subref{fig:demo_nodes_2:c2}) and $u_{2} \equiv \dot{x}$ (\subref{fig:demo_nodes_2:a3}, \subref{fig:demo_nodes_2:a4}, \subref{fig:demo_nodes_2:b3}, \subref{fig:demo_nodes_2:b4}, \subref{fig:demo_nodes_2:c3}, \subref{fig:demo_nodes_2:c4}), obtained using polynomials with degrees $N = 4$ (\subref{fig:demo_nodes_2:a1}, \subref{fig:demo_nodes_2:a2}, \subref{fig:demo_nodes_2:a3}, \subref{fig:demo_nodes_2:a4}), $N = 8$ (\subref{fig:demo_nodes_2:b1}, \subref{fig:demo_nodes_2:b2}, \subref{fig:demo_nodes_2:b3}, \subref{fig:demo_nodes_2:b4}) and $N = 60$ (\subref{fig:demo_nodes_2:c1}, \subref{fig:demo_nodes_2:c2}, \subref{fig:demo_nodes_2:c3}, \subref{fig:demo_nodes_2:c4}). Zoomed domain $0 \leqslant t \leqslant 10$ is presented on (\subref{fig:demo_nodes_2:a2}, \subref{fig:demo_nodes_2:a4}, \subref{fig:demo_nodes_2:b2}, \subref{fig:demo_nodes_2:b4}, \subref{fig:demo_nodes_2:c2}, \subref{fig:demo_nodes_2:c4}).
}
\label{fig:demo_nodes_2}
\end{figure}

The obtained results are presented in Fig.~\ref{fig:demo_nodes_8} and~\ref{fig:demo_nodes_2}. Comparisons are presented between the exact analytical solution $\mathbf{u}^{\rm ex}(t)$, the numerical solution at grid nodes $\mathbf{u}_{n}$ (\ref{eq:ader_dg_node_sol}), the local numerical solution $\mathbf{u}_{L}$ (\ref{eq:local_sol_assembly}), and the improved local numerical solution $\mathbf{u}_{\rm IL}$ (\ref{eq:imp_local_sol_assembly}). Fig.~\ref{fig:demo_nodes_8} presents the results for basis polynomials of degrees $N = 1$, $2$ and $3$, obtained with a discretization step ${\Delta t} = \pi$ ($8$ discretization domains on the definition domain $[0,\, 8\pi]$). Fig.~\ref{fig:demo_nodes_2} presents the results for higher degrees of basis polynomials $N = 4$, $8$ and $60$, obtained with a larger discretization step ${\Delta t} = 4\pi$ (only $2$ discretization domains on the definition domain $[0,\, 8\pi]$). For a more visual and detailed comparison of the exact analytical solution and the different numerical solutions, Fig.~\ref{fig:demo_nodes_8} and~\ref{fig:demo_nodes_2} also separately show zoomed domain $0 \leqslant \tau \leqslant 10$. The exact analytical solution $\mathbf{u}^{\rm ex}(t)$, local numerical solution $\mathbf{u}_{L}(t)$, and improved local numerical solution $\mathbf{u}_{\rm IL}(t)$ for plotting in the Figures were tabulated with $50$ equally spaced points per discretization domain $\Omega_{n}$.

The dependencies of the numerical solutions $\mathbf{u}_{L}$, $\mathbf{u}_{\rm IL}$, $\mathbf{u}_{n}$ for polynomial degree $N = 1$ shown in Fig.~\ref{fig:demo_nodes_8} (\subref{fig:demo_nodes_8:a1}, \subref{fig:demo_nodes_8:a2}, \subref{fig:demo_nodes_8:a3}, \subref{fig:demo_nodes_8:a4}) demonstrate significant dissipation, which is not characteristic of the exact analytical solution of the original conservative problem. This is a numerical artifact and is associated with the large discretization step ${\Delta t} = \pi$ and the very high stability~\cite{ader_dg_ode_jsc} of the ADER-DG numerical method with a local DG predictor. In this case, the discontinuity of the local numerical solution $\mathbf{u}_{L}$ at the left ends $t_{n}$ of the discretization domains $\Omega_{n}$ is clearly evident, where the values of the numerical solution $\mathbf{u}_{L}(t_{n}^{+})$ do not coincide with the values of the numerical solution $\mathbf{u}_{n}$ at the grid nodes $t_{n}$. These features are demonstrated even more clearly in Fig.~\ref{fig:demo_nodes_8} (\subref{fig:demo_nodes_8:a2}, \subref{fig:demo_nodes_8:a4}) with zoomed domain $0 \leqslant \tau \leqslant 10$. These effects occur for both solution component $u_{1}$ and component $u_{2}$, with no significant differences in the behavior of the individual components. The presented dependencies of the improved local numerical solution $\mathbf{u}_{\rm IL}$, shown in the same Fig.~\ref{fig:demo_nodes_8} for the case of polynomials of degree $N = 1$, demonstrate a significantly higher quality of the numerical solution. In particular, the observed discontinuities in the local solution $\mathbf{u}_{L}(t)$ at the grid nodes $t_{n}$ are eliminated, while at the left $t_{n}$ and right $t_{n+1}$ ends of the discretization domain $\Omega_{n}$, the improved numerical solution $\mathbf{u}_{\rm IL}(t_{n})$ coincides with the numerical solution $\mathbf{u}_{n}$ at the grid nodes $t_{n}$. It is also clear that in the case of the improved local numerical solution $\mathbf{u}_{\rm IL}$, the functional dependence is represented by quadratic functions of the mapping coordinate $\tau$, while for the local numerical solution $\mathbf{u}_{L}(t_{n})$, it is represented only by linear functions of the mapping coordinate $\tau$, which significantly affects the qualitative properties of the numerical solution.


The dependencies of the numerical solutions $\mathbf{u}_{L}$, $\mathbf{u}_{\rm IL}$, $\mathbf{u}_{n}$ for polynomial degree $N = 2$ presented in Fig.~\ref{fig:demo_nodes_8} (\subref{fig:demo_nodes_8:b1}, \subref{fig:demo_nodes_8:b2}, \subref{fig:demo_nodes_8:b3}, \subref{fig:demo_nodes_8:b4}) no longer demonstrate the strong numerical dissipation characteristic of a large discretization step ${\Delta t} = \pi$, taking into account the high stability of the ADER-DG numerical method with a local DG predictor. The selected discretization step ${\Delta t} = \pi$ also allows one to clearly observe the discontinuity of the local numerical solution $\mathbf{u}_{L}$ at the left ends $t_{n}$ of the discretization domains $\Omega_{n}$, when the numerical solution $\mathbf{u}_{L}(t_{n}^{+})$ do not coincide with the values of the numerical solution $\mathbf{u}_{n}$ at the grid nodes $t_{n}$. The dependencies of the improved local numerical solution $\mathbf{u}_{\rm IL}$, presented in the same Fig.~\ref{fig:demo_nodes_8} (\subref{fig:demo_nodes_8:b1}, \subref{fig:demo_nodes_8:b2}, \subref{fig:demo_nodes_8:b3}, \subref{fig:demo_nodes_8:b4}) for the case of polynomial degree $N = 2$, also demonstrate a significantly higher quality of the numerical solution, where the observed discontinuities in the local solution $\mathbf{u}_{L}$ at the grid nodes $t_{n}$ are eliminated. At the left $t_{n}$ and right $t_{n+1}$ ends of the discretization domain $\Omega_{n}$, the improved numerical solution $\mathbf{u}_{\rm IL}$ coincides with the $\mathbf{u}_{n}$ at the grid nodes $t_{n}$. These differences are especially clearly presented in Fig.~\ref{fig:demo_nodes_8} (\subref{fig:demo_nodes_8:b2}, \subref{fig:demo_nodes_8:b4}) with zoomed domain $0 \leqslant \tau \leqslant 10$. The dependencies of the numerical solutions $\mathbf{u}_{L}$, $\mathbf{u}_{\rm IL}$, $\mathbf{u}_{n}$ for the degree of polynomials of $N = 3$ for the selected value of the discretization step ${\Delta t} = \pi$ presented in Fig.~\ref{fig:demo_nodes_8} (\subref{fig:demo_nodes_8:c1}, \subref{fig:demo_nodes_8:c2}, \subref{fig:demo_nodes_8:c3}, \subref{fig:demo_nodes_8:c4}) practically do not demonstrate numerical dissipation, as well as discontinuities of the local numerical solution at the left ends of the discretization domains --- partial discontinuities of the local numerical solution can be observed in Fig.~\ref{fig:demo_nodes_8} (\subref{fig:demo_nodes_8:c4}) with zoomed domain $0 \leqslant \tau \leqslant 10$ for component $u_{2}$ of the solution.

Fig.~\ref{fig:demo_nodes_2} presents numerical solutions $\mathbf{u}_{L}$, $\mathbf{u}_{\rm IL}$, $\mathbf{u}_{n}$ of the problem (\ref{eq:demo_ode}) for significantly higher polynomial degrees $N = 4$, $8$, $60$ and a significantly larger discretization step ${\Delta t} = 4\pi$, compared to $N = 1$, $2$, $3$ and ${\Delta t} = \pi$ for the results presented in Fig.~\ref{fig:demo_nodes_8}. This is done to demonstrate the qualitative properties of the local $\mathbf{u}_{L}$ and improved local $\mathbf{u}_{\rm IL}$ numerical solutions and the possibilities of their application. The dependencies of the local numerical solution $\mathbf{u}_{L}(t)$ for polynomial degree $N = 4$ presented in Fig.~\ref{fig:demo_nodes_2} (\subref{fig:demo_nodes_2:a1}, \subref{fig:demo_nodes_2:a2}, \subref{fig:demo_nodes_2:a3}, \subref{fig:demo_nodes_2:a4}) again demonstrate significant dissipation, which is associated with a significant fourfold increase in the discretization step ${\Delta t}$. Since the degree $N = 4$ is already quite high, discontinuities in the dependence of the local numerical solution $\mathbf{u}_{L}(t)$, as in Fig.~\ref{fig:demo_nodes_8} (\subref{fig:demo_nodes_8:a1}, \subref{fig:demo_nodes_8:a2}, \subref{fig:demo_nodes_8:a3}, \subref{fig:demo_nodes_8:a4}), are practically not observed, with the exception of the discontinuities of the local solution in Fig.~\ref{fig:demo_nodes_2} (\subref{fig:demo_nodes_2:a3}). However, in this case, the dependencies of the improved local numerical solution $\mathbf{u}_{\rm IL}(t)$ also demonstrate a significantly higher quality of the numerical solution, especially in terms of the agreement between the improved local solution and the exact analytical solution at the point $t_{0}$. The results presented in Figs.~\ref{fig:demo_nodes_2} (\subref{fig:demo_nodes_2:b1}, \subref{fig:demo_nodes_2:b2}, \subref{fig:demo_nodes_2:b3}, \subref{fig:demo_nodes_2:b4}) and~\ref{fig:demo_nodes_2} (\subref{fig:demo_nodes_2:c1}, \subref{fig:demo_nodes_2:c2}, \subref{fig:demo_nodes_2:c3}, \subref{fig:demo_nodes_2:c4}) for polynomial degrees $N = 8$ and $60$ again demonstrate excellent dependencies of the numerical solutions $\mathbf{u}_{L}$, $\mathbf{u}_{\rm IL}$, $\mathbf{u}_{n}$, especially considering such a large discretization step ${\Delta t} = 4\pi$, which accommodates two full periods of oscillations of the harmonic oscillator (\ref{eq:demo_ode}).

\section{Applications of the numerical method}
\label{sec:apps}

This Section presents four examples of applying the ADER-DG numerical method with a local DG predictor to solving an initial value problem for the ODE system (\ref{eq:ivp_ode_diff_src}). The examples are presented separately in four Subsections: Subsection~\ref{sec:apps:exp_diss} ``Dahlquist's test equation'', Subsection~\ref{sec:apps:lin_diss} ``Linear exp--test'', Subsection~\ref{sec:apps:harm_osc} ``Harmonic oscillator'', and Subsection~\ref{sec:apps:pend} ``Nonlinear mathematical pendulum''.

The qualitative and quantitative properties of the numerical solutions are studied, and the dependencies of the local error $\varepsilon(t)$ of the numerical solution on the argument $t\in\Omega$ are calculated. The local error $\varepsilon(t)$ of the numerical solutions was selected based on a point estimate and is calculated based on the following relationships:
\begin{equation}\label{eq:eps_local_def}
\begin{split}
\varepsilon_{L}(t) = \left|\mathbf{u}_{L}(t) - \mathbf{u}^{\rm ex}(t)\right|,\ 
\varepsilon_{\rm IL}(t) = \left|\mathbf{u}_{\rm IL}(t) - \mathbf{u}^{\rm ex}(t)\right|,\ 
\varepsilon_{n}(t_{n}) = \left|\mathbf{u}_{n} - \mathbf{u}^{\rm ex}(t_{n})\right|,
\end{split}
\end{equation}
where $|\ldots| \equiv ||\ldots||_{\mathcal{L}_{\infty}}$ is defined by (\ref{eq:norm_in_rd}), $\mathbf{u}^{\rm ex}(t)$ is the exact analytical solution of the problem (or the reference solution, but for all presented examples, the exact analytical solution is known, which is used as the reference), used for the local numerical solution $\mathbf{u}_{L}(t)$, the improved local numerical solution $\mathbf{u}_{\rm IL}(t)$, and the numerical solution $\mathbf{u}_{n}$ at grid nodes, respectively. The errors $\varepsilon(t_{n, s})$ of the local solution $\mathbf{u}_{L}$ and improved local solution $\mathbf{u}_{\rm IL}$ will be represented in sub-nodes $t_{n, s}$ located between the grid nodes $t_{n}$. In this work, $S = 50$ points are uniformly added in the discretization domains $\Omega_{n}$, at which estimates of the local numerical solutions are made.

The dependencies of the global error $e({\Delta t})$ on the discretization step ${\Delta t}$ are calculated, based on which empirical convergence orders $p$ for various types of numerical solutions are determined. The domain of definition $\Omega$ is discretized into $M$ discretization domains $\{\Omega_{n}\}$, with a discretization step ${\Delta t} = (t_{f} - t_{0})/M$. In this paper, $10$ types of global errors $e$ of numerical solutions are used: $4$ types of global errors for the numerical solution $\mathbf{u}_{n}$ at grid nodes and $3$ types of global errors each for the local numerical solution $\mathbf{u}_{L}$ and the improved local numerical solution $\mathbf{u}_{\rm IL}$. The global errors $e^{n}$ for the numerical solution $\mathbf{u}_{n}$ at nodes are selected as follows:
\begin{equation}\label{eq:eps_un_global_def}
\begin{array}{ll}
e^{n}_{f} = \left|\mathbf{u}_{M} - \mathbf{u}^{\rm ex}(t_{f})\right|,&\quad
e^{n}_{L_{\infty}} = \max\limits_{0 \leqslant n \leqslant M}\left|\mathbf{u}_{n} - \mathbf{u}^{\rm ex}(t_{n})\right|,\\
e^{n}_{L_{1}} = {\displaystyle\sum\limits_{n = 0}^{M}} {\Delta t}_{n}\left|\mathbf{u}_{n} - \mathbf{u}^{\rm ex}(t_{n})\right|,&\quad
e^{n}_{L_{2}} = \left[{\displaystyle\sum\limits_{n = 0}^{M}} {\Delta t}_{n}\left|\mathbf{u}_{n} - \mathbf{u}^{\rm ex}(t_{n})\right|^{2}\right]^{1/2}.
\end{array}
\end{equation}
The global errors $e^{l}$ and $e^{\rm imp}$ for the local numerical solution $\mathbf{u}_{L}$ and the improved local numerical solution $\mathbf{u}_{\rm IL}$, respectively, are selected as follows:
\begin{equation}\label{eq:eps_ul_global_def}
\begin{array}{ll}
e^{l}_{L_{\infty}} = \max\limits_{0 \leqslant n \leqslant M}\sup\limits_{t\in\Omega_{n}}\left|\mathbf{u}_{L}(t) - \mathbf{u}^{\rm ex}(t)\right|,&\quad
e^{\rm imp}_{L_{\infty}} = \max\limits_{0 \leqslant n \leqslant M}\sup\limits_{t\in\Omega_{n}}\left|\mathbf{u}_{\rm IL}(t) - \mathbf{u}^{\rm ex}(t)\right|,\\
e^{l}_{L_{1}} = {\displaystyle\sum\limits_{n = 0}^{M-1}}\ {\displaystyle\int\limits_{\Omega_{n}}}\left|\mathbf{u}_{L}(t) - \mathbf{u}^{\rm ex}(t)\right| dt,&\quad
e^{\rm imp}_{L_{1}} = {\displaystyle\sum\limits_{n = 0}^{M-1}}\ {\displaystyle\int\limits_{\Omega_{n}}}\left|\mathbf{u}_{\rm IL}(t) - \mathbf{u}^{\rm ex}(t)\right| dt,\\
e^{l}_{L_{2}} = \left[{\displaystyle\sum\limits_{n = 0}^{M-1}}\ {\displaystyle\int\limits_{\Omega_{n}}}\left|\mathbf{u}_{L}(t) - \mathbf{u}^{\rm ex}(t)\right|^{2} dt\right]^{1/2},&\quad
e^{\rm imp}_{L_{2}} = \left[{\displaystyle\sum\limits_{n = 0}^{M-1}}\ {\displaystyle\int\limits_{\Omega_{n}}}\left|\mathbf{u}_{\rm IL}(t) - \mathbf{u}^{\rm ex}(t)\right|^{2} dt\right]^{1/2}.
\end{array}
\end{equation}
The calculation of global errors (\ref{eq:eps_ul_global_def}) for the local solution $\mathbf{u}_{L}(t)$ and improved local solution $\mathbf{u}_{\rm IL}(t)$ in the discretization domain $\Omega_{n}$ is carried out on the basis of the replacement of integrals by finite sums over sub-nodes $t_{n, s}$, multiplied by the discretization steps between sub-nodes $\Delta t_{n}/S$. The operation $\sup$ by $t\in\Omega_{n}$ in (\ref{eq:eps_ul_global_def}) has been replaced by the operation $\max$ by sub-nodes $t_{n, s}$.

According to these $10$ types of global error $e$ (\ref{eq:eps_un_global_def}), (\ref{eq:eps_ul_global_def}), $10$ empirical convergence orders $p$ for the numerical solutions $\mathbf{u}_{L}$, $\mathbf{u}_{\rm IL}$, $\mathbf{u}_{n}$ are selected. The calculation of empirical convergence orders $p$ is carried out separately for each presented type of numerical solution, and is based on a power approximation of the dependence of global errors $e({\Delta t}) \propto {\Delta t}^{p}$ on the discretization step ${\Delta t}$. In all examples, the number of discretization domains $M = 10$, $12$, $14$, $16$, $18$, $20$, $22$, $24$ is chosen. Based on which set of the numbers of discretization domains $M$ the empirical convergence orders $p$ are calculated by least squares approximation in log-log scale: $\lg{e({\Delta t})} \propto p\cdot\lg{{\Delta t}}$, at $8$ data points. The empirical convergence orders $p$ are compared with the expected theoretical values $p_{\rm th.}$:
\begin{equation}\label{eq:conv_ords_exp}
p_{\rm th.}^{n} = 2N+1,\quad
p_{\rm th.}^{l} = N+1,\quad
p_{\rm th.}^{\rm imp} = N+2,
\end{equation}
which coincide with the approximation orders defined above of numerical solutions $\mathbf{u}_{L}$, $\mathbf{u}_{\rm IL}$, $\mathbf{u}_{n}$, respectively, for the ADER-DG numerical method with a local DG predictor.

\subsection{Example 1: Dahlquist's test equation}
\label{sec:apps:exp_diss}

\begin{figure}[h!]
\captionsetup[subfigure]{%
	position=bottom,
	font+=smaller,
	textfont=normalfont,
	singlelinecheck=off,
	justification=raggedright
}
\centering
\begin{subfigure}{0.24\textwidth}
\includegraphics[width=\textwidth]{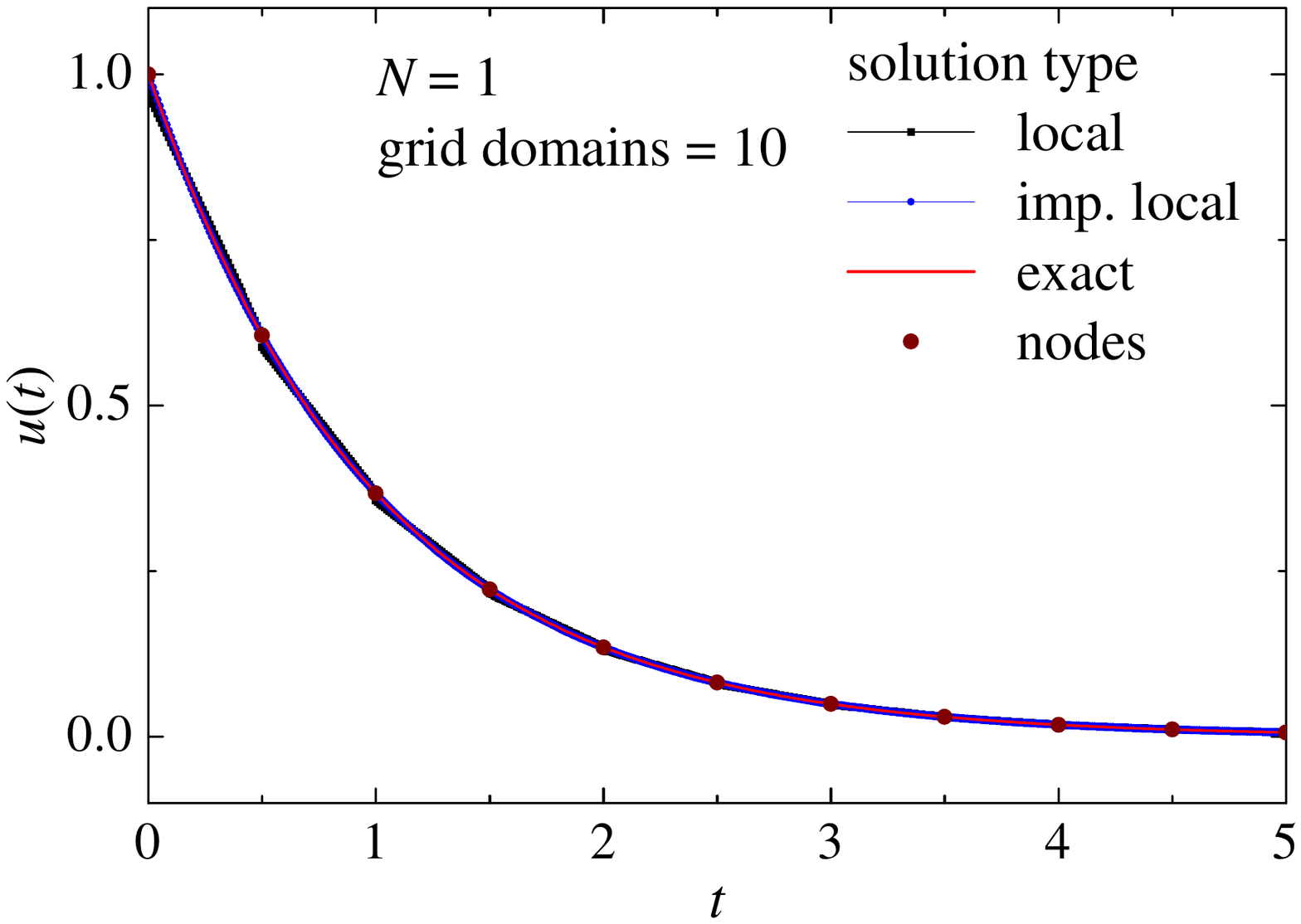}
\vspace{-8mm}\caption{\label{fig:exp_diss:a1}}
\end{subfigure}
\begin{subfigure}{0.24\textwidth}
\includegraphics[width=\textwidth]{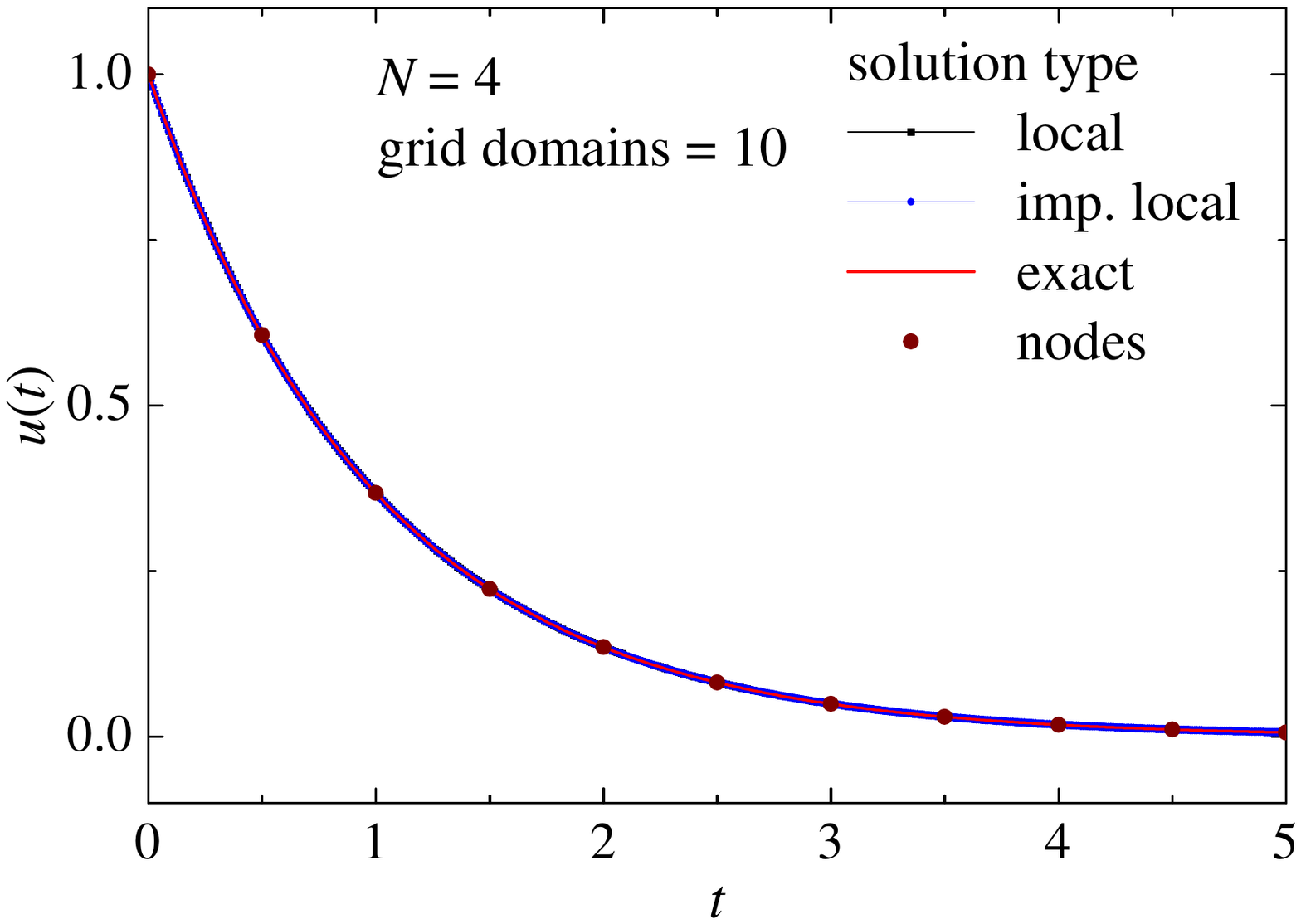}
\vspace{-8mm}\caption{\label{fig:exp_diss:a2}}
\end{subfigure}
\begin{subfigure}{0.24\textwidth}
\includegraphics[width=\textwidth]{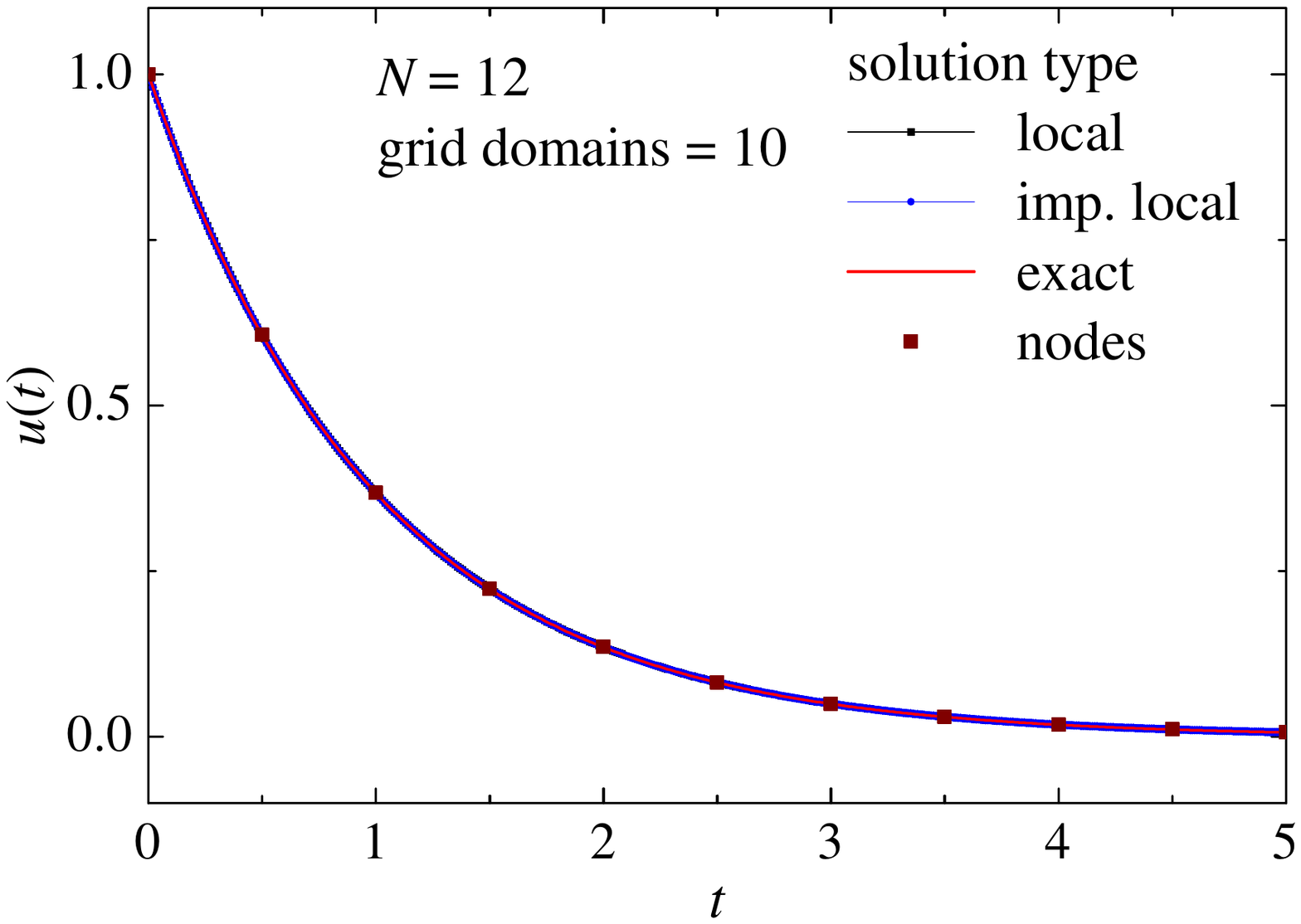}
\vspace{-8mm}\caption{\label{fig:exp_diss:a3}}
\end{subfigure}
\begin{subfigure}{0.24\textwidth}
\includegraphics[width=\textwidth]{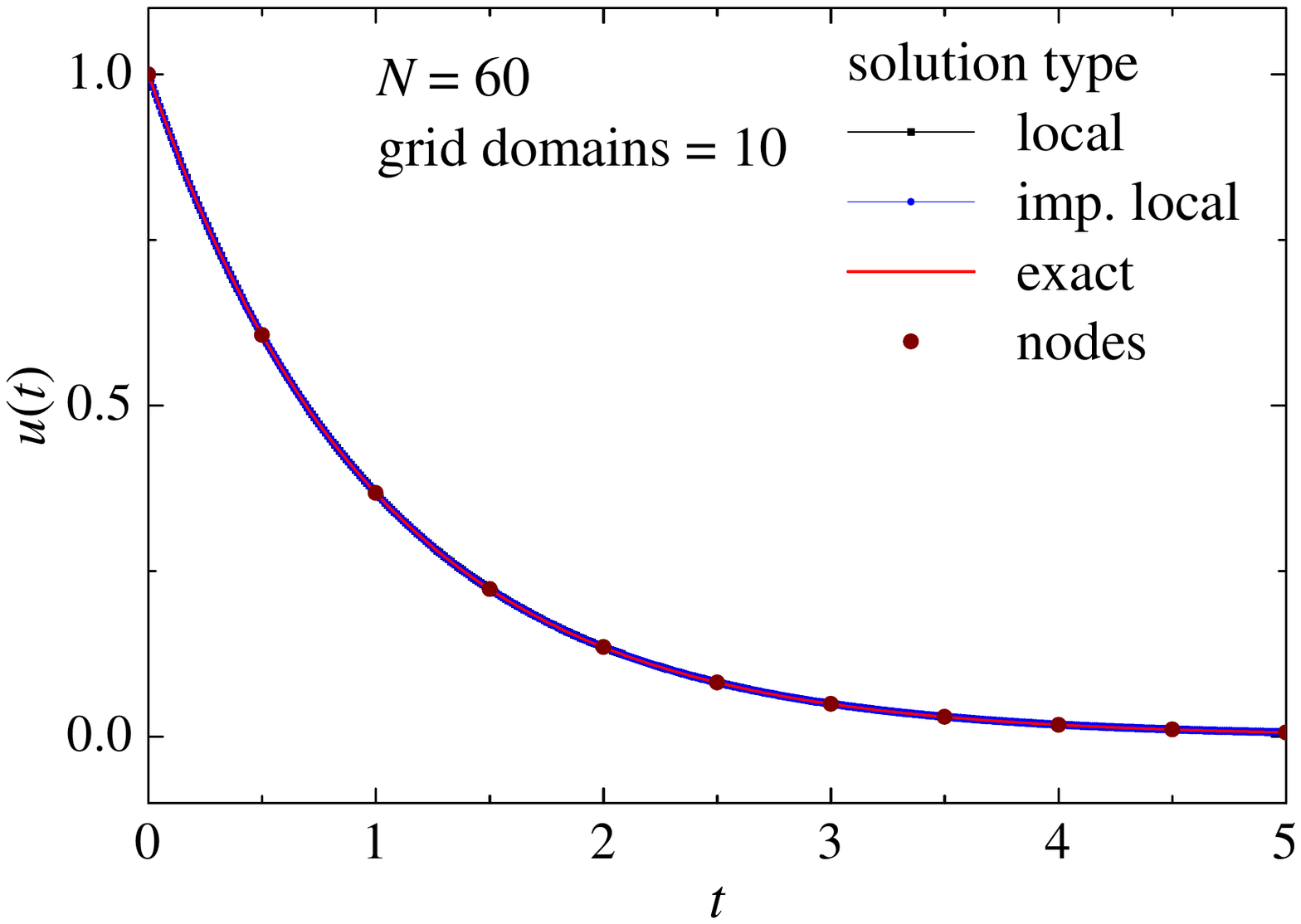}
\vspace{-8mm}\caption{\label{fig:exp_diss:a4}}
\end{subfigure}\\
\begin{subfigure}{0.24\textwidth}
\includegraphics[width=\textwidth]{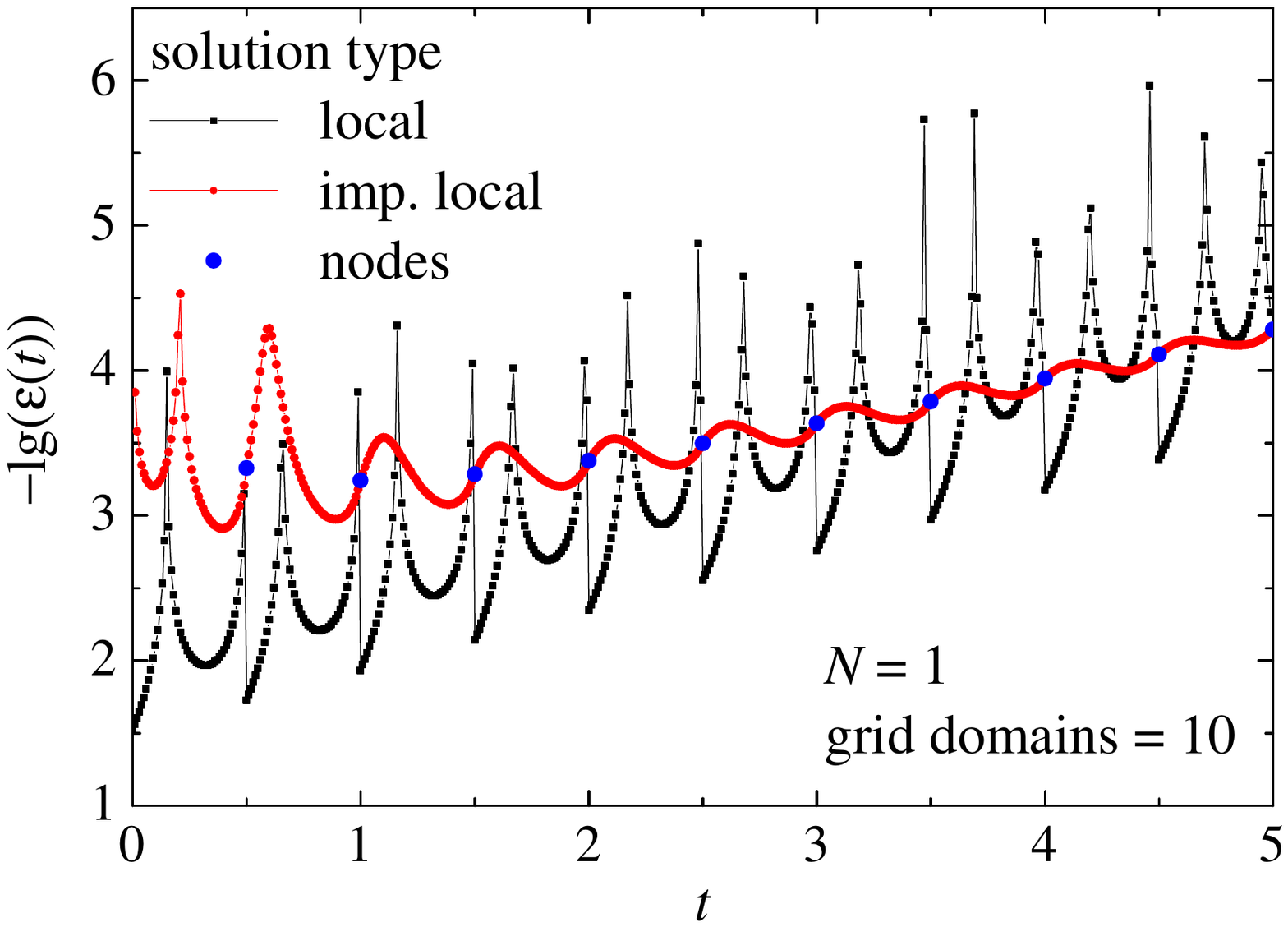}
\vspace{-8mm}\caption{\label{fig:exp_diss:b1}}
\end{subfigure}
\begin{subfigure}{0.24\textwidth}
\includegraphics[width=\textwidth]{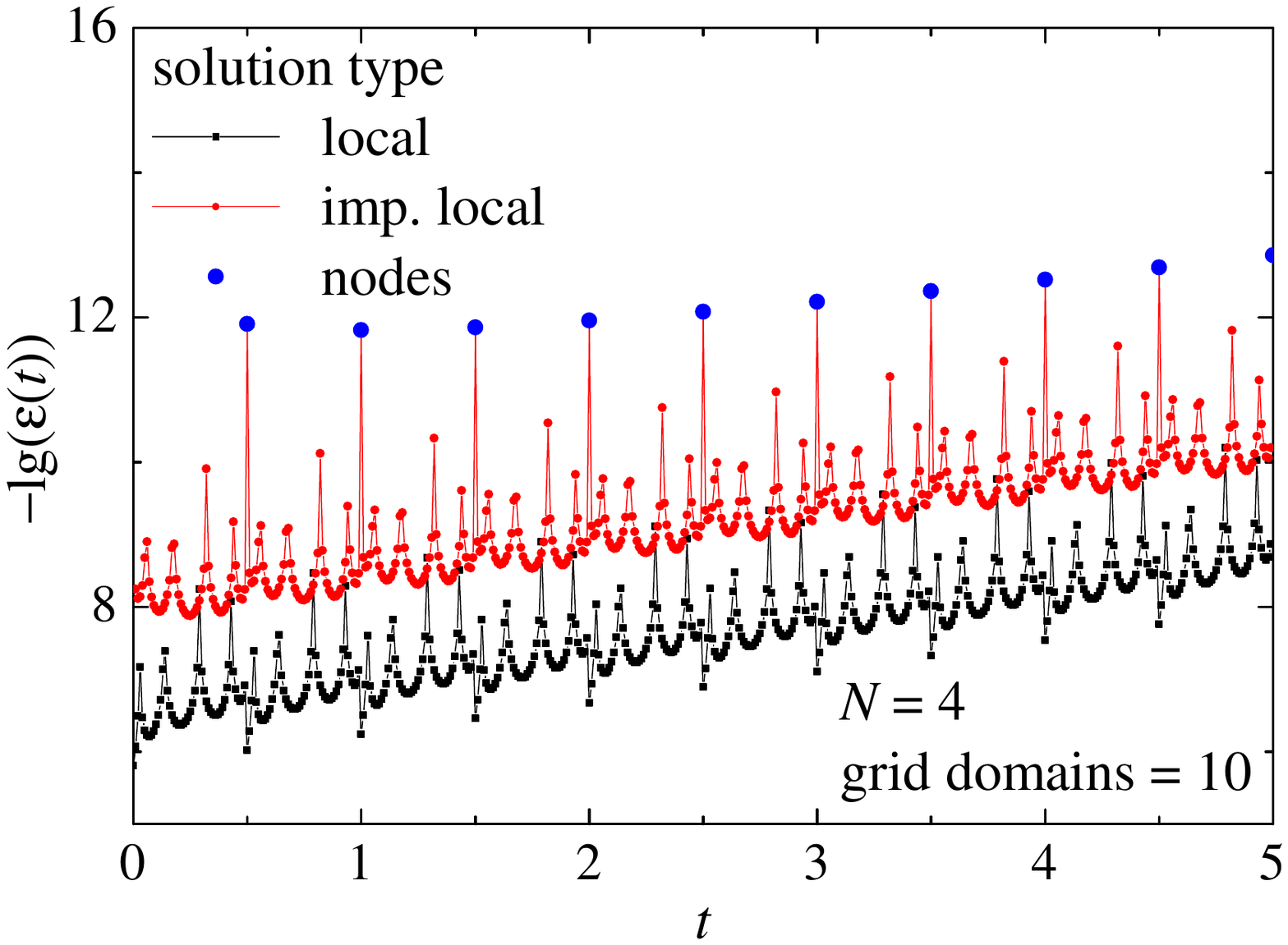}
\vspace{-8mm}\caption{\label{fig:exp_diss:b2}}
\end{subfigure}
\begin{subfigure}{0.24\textwidth}
\includegraphics[width=\textwidth]{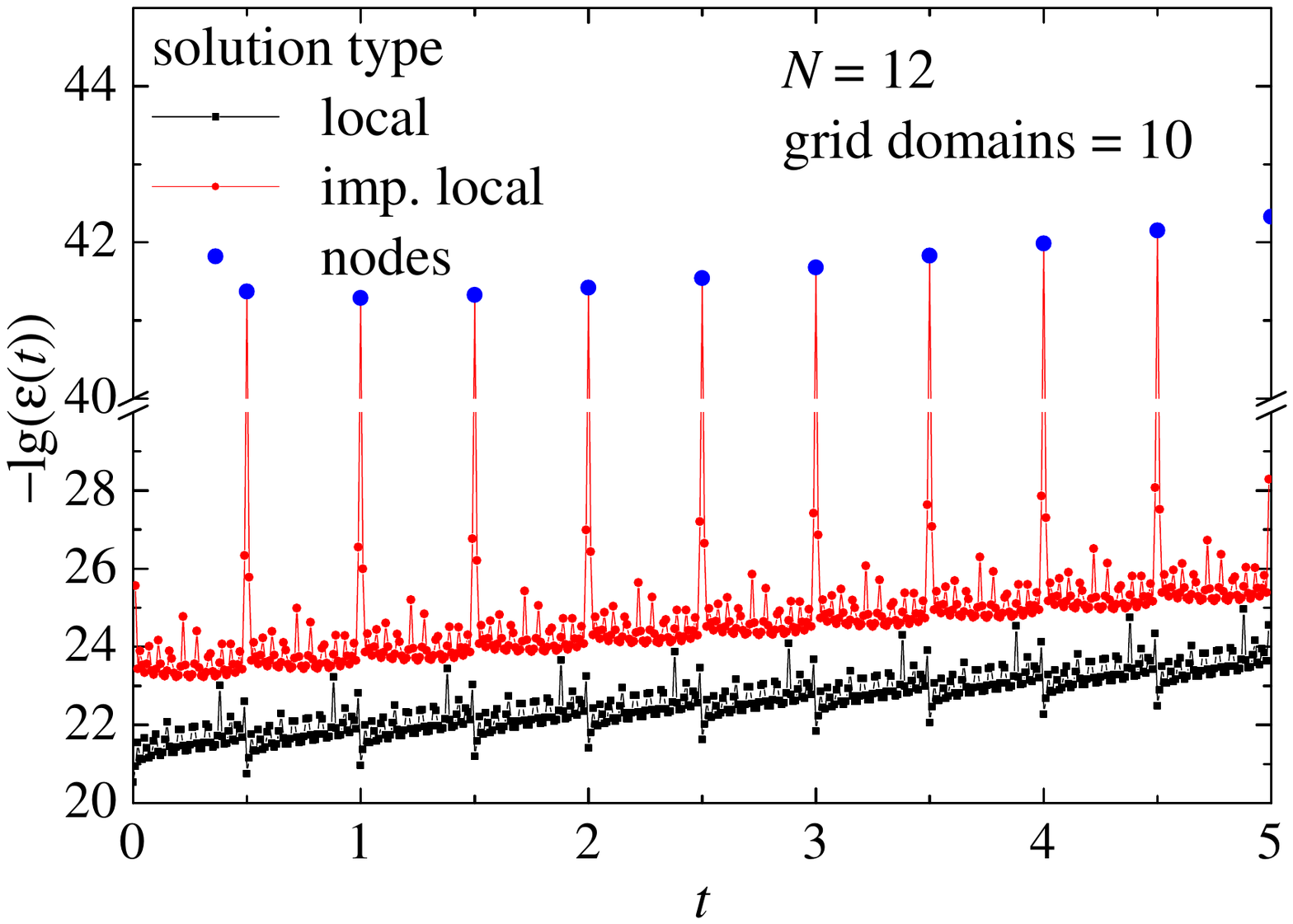}
\vspace{-8mm}\caption{\label{fig:exp_diss:b3}}
\end{subfigure}
\begin{subfigure}{0.24\textwidth}
\includegraphics[width=\textwidth]{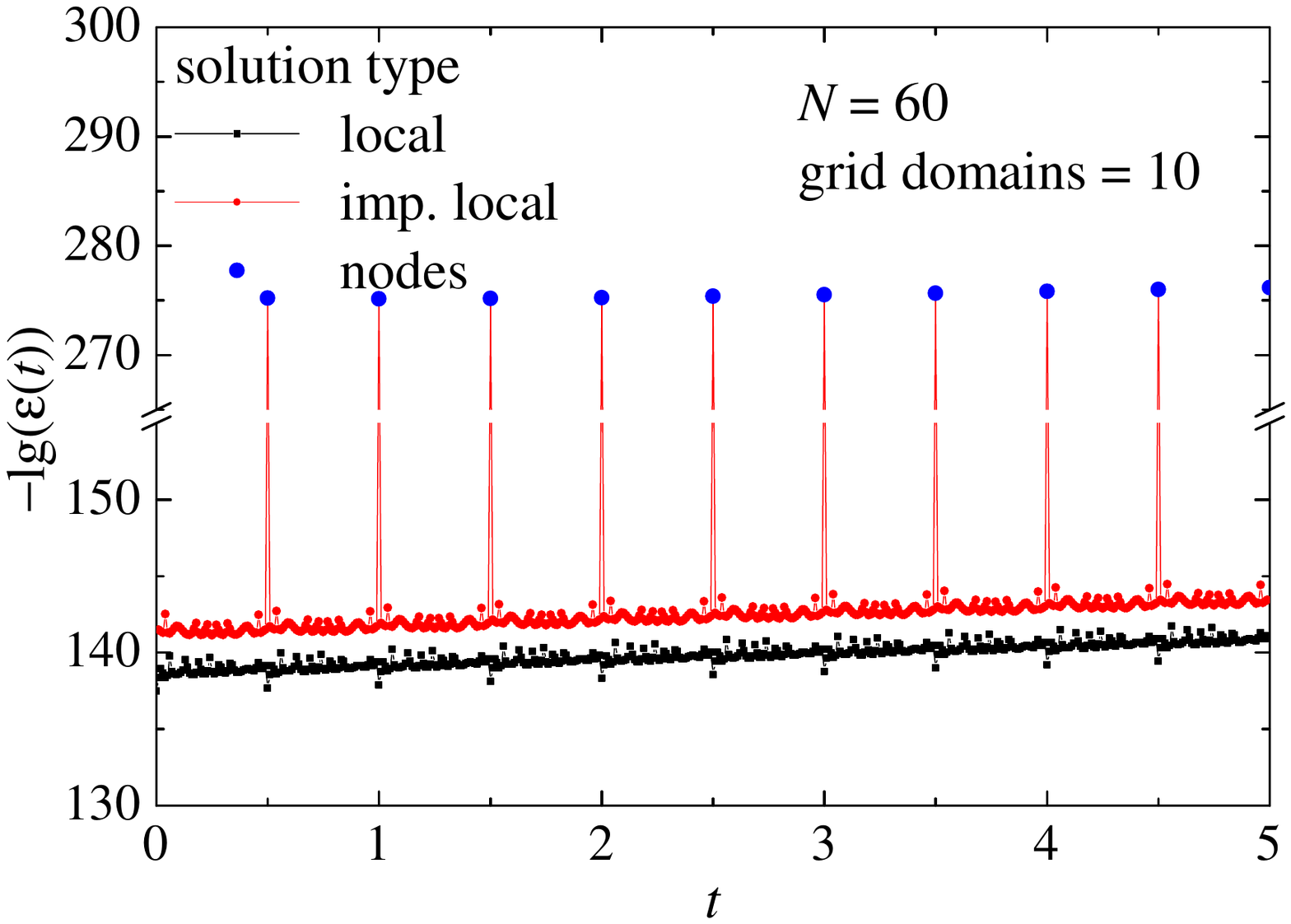}
\vspace{-8mm}\caption{\label{fig:exp_diss:b4}}
\end{subfigure}\\
\begin{subfigure}{0.24\textwidth}
\includegraphics[width=\textwidth]{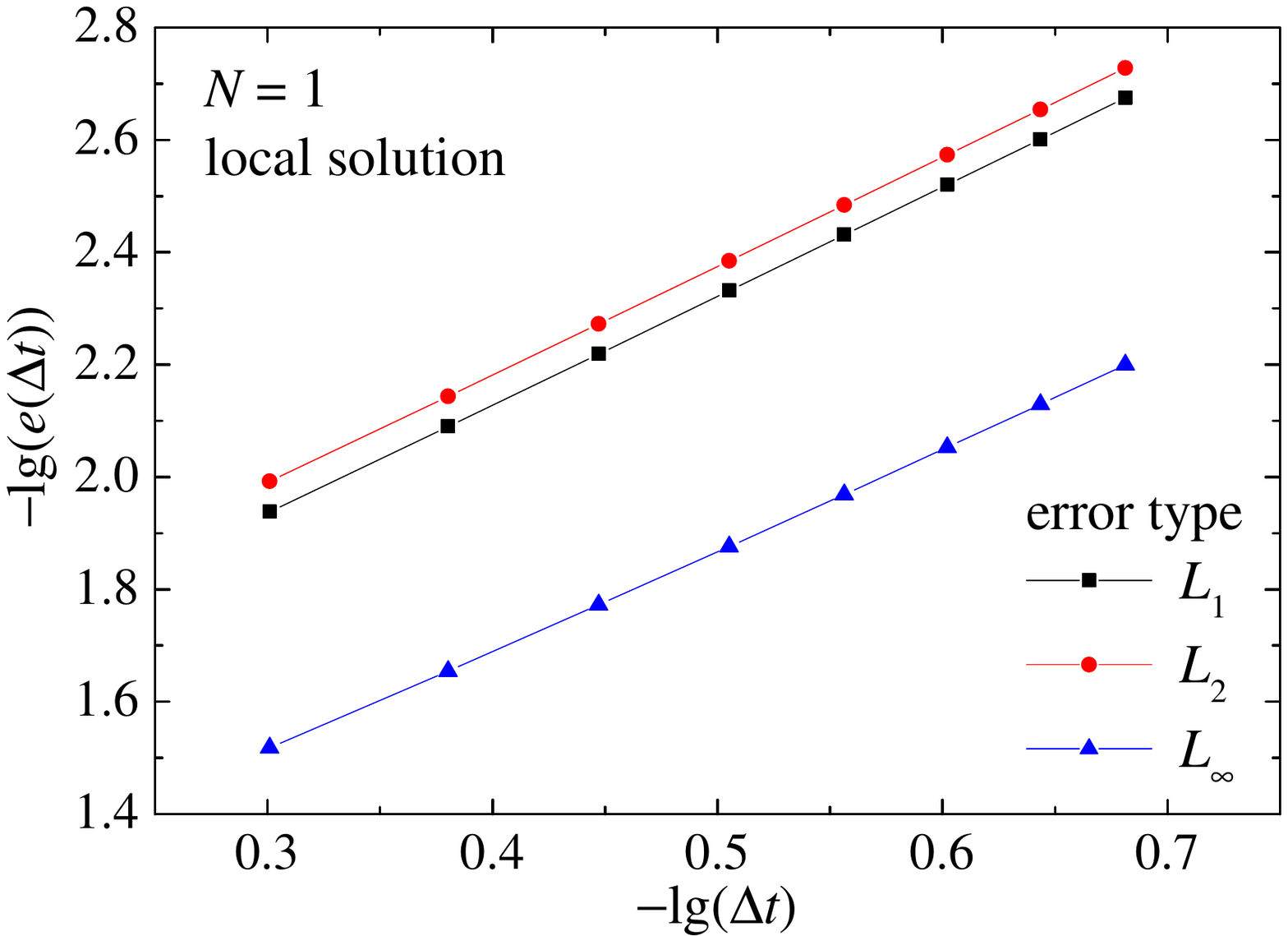}
\vspace{-8mm}\caption{\label{fig:exp_diss:c1}}
\end{subfigure}
\begin{subfigure}{0.24\textwidth}
\includegraphics[width=\textwidth]{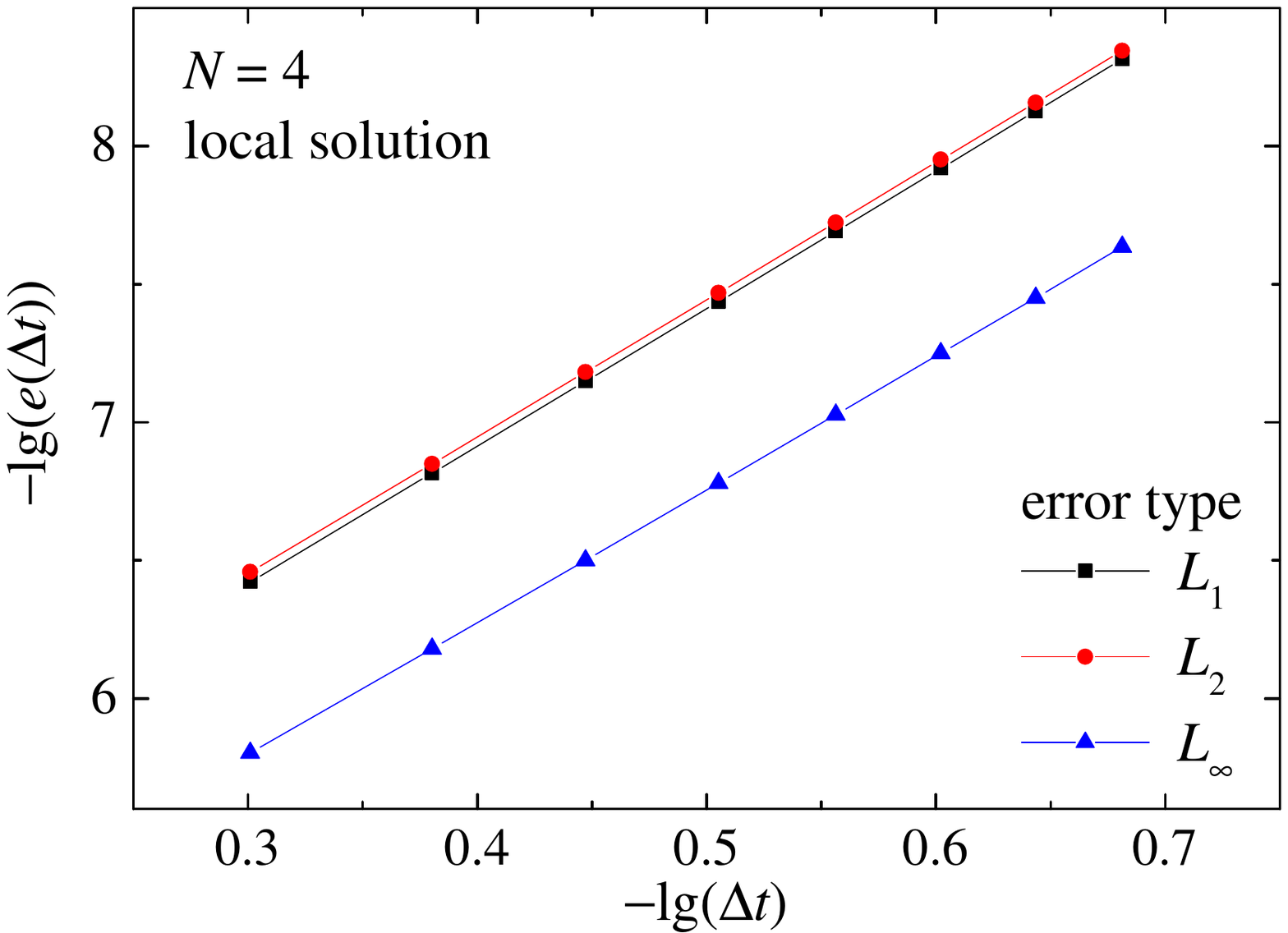}
\vspace{-8mm}\caption{\label{fig:exp_diss:c2}}
\end{subfigure}
\begin{subfigure}{0.24\textwidth}
\includegraphics[width=\textwidth]{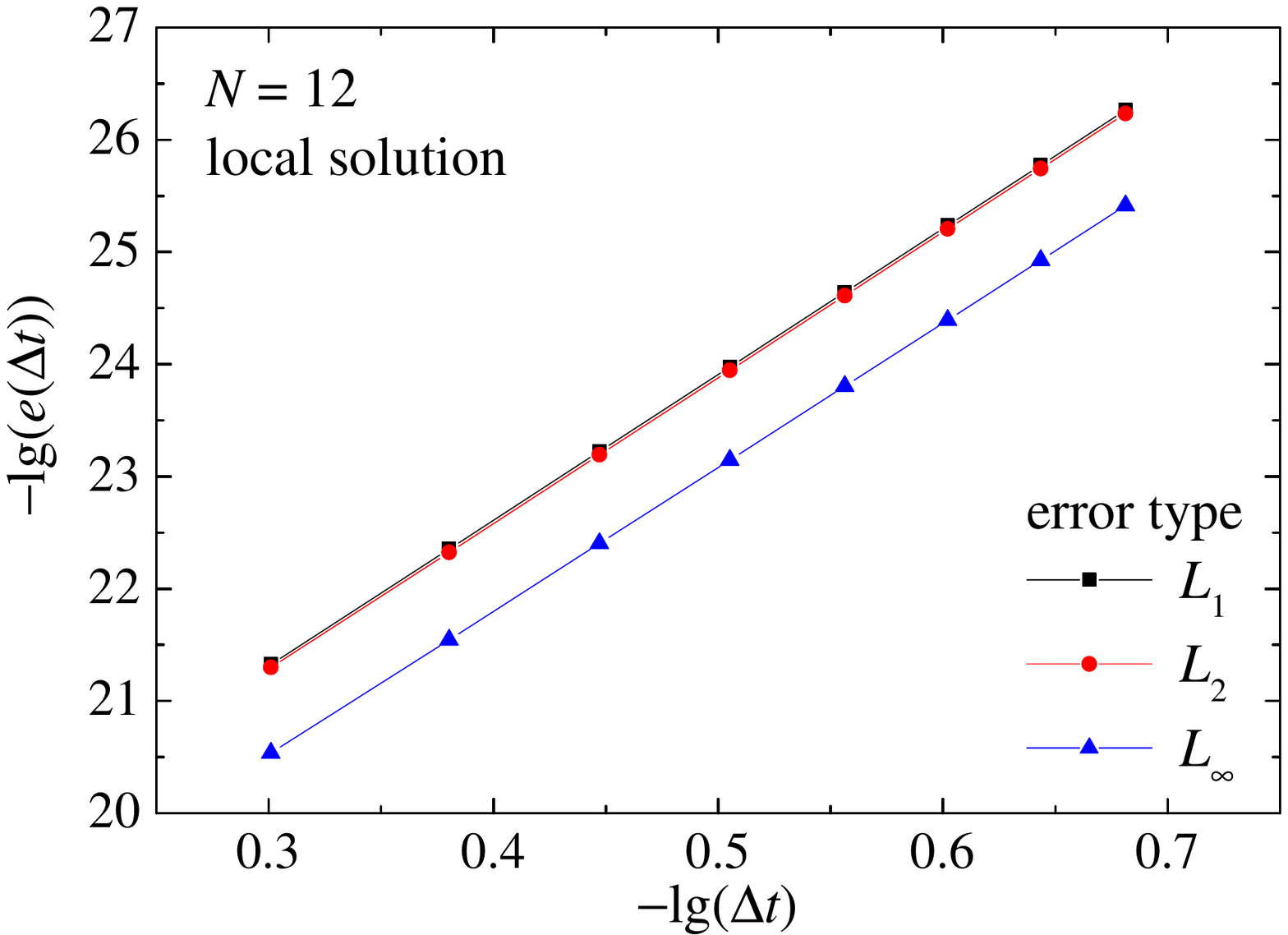}
\vspace{-8mm}\caption{\label{fig:exp_diss:c3}}
\end{subfigure}
\begin{subfigure}{0.24\textwidth}
\includegraphics[width=\textwidth]{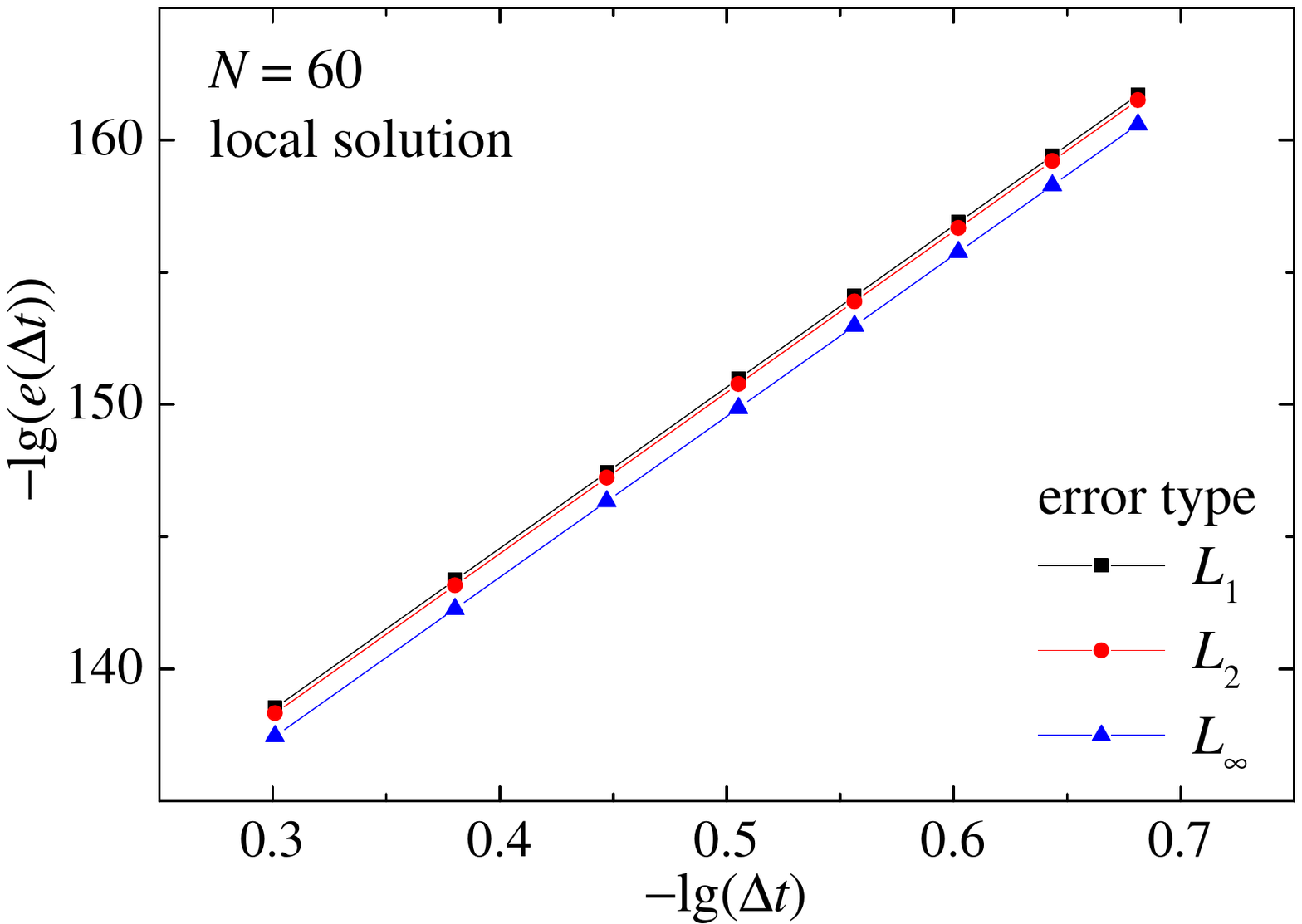}
\vspace{-8mm}\caption{\label{fig:exp_diss:c4}}
\end{subfigure}\\
\begin{subfigure}{0.24\textwidth}
\includegraphics[width=\textwidth]{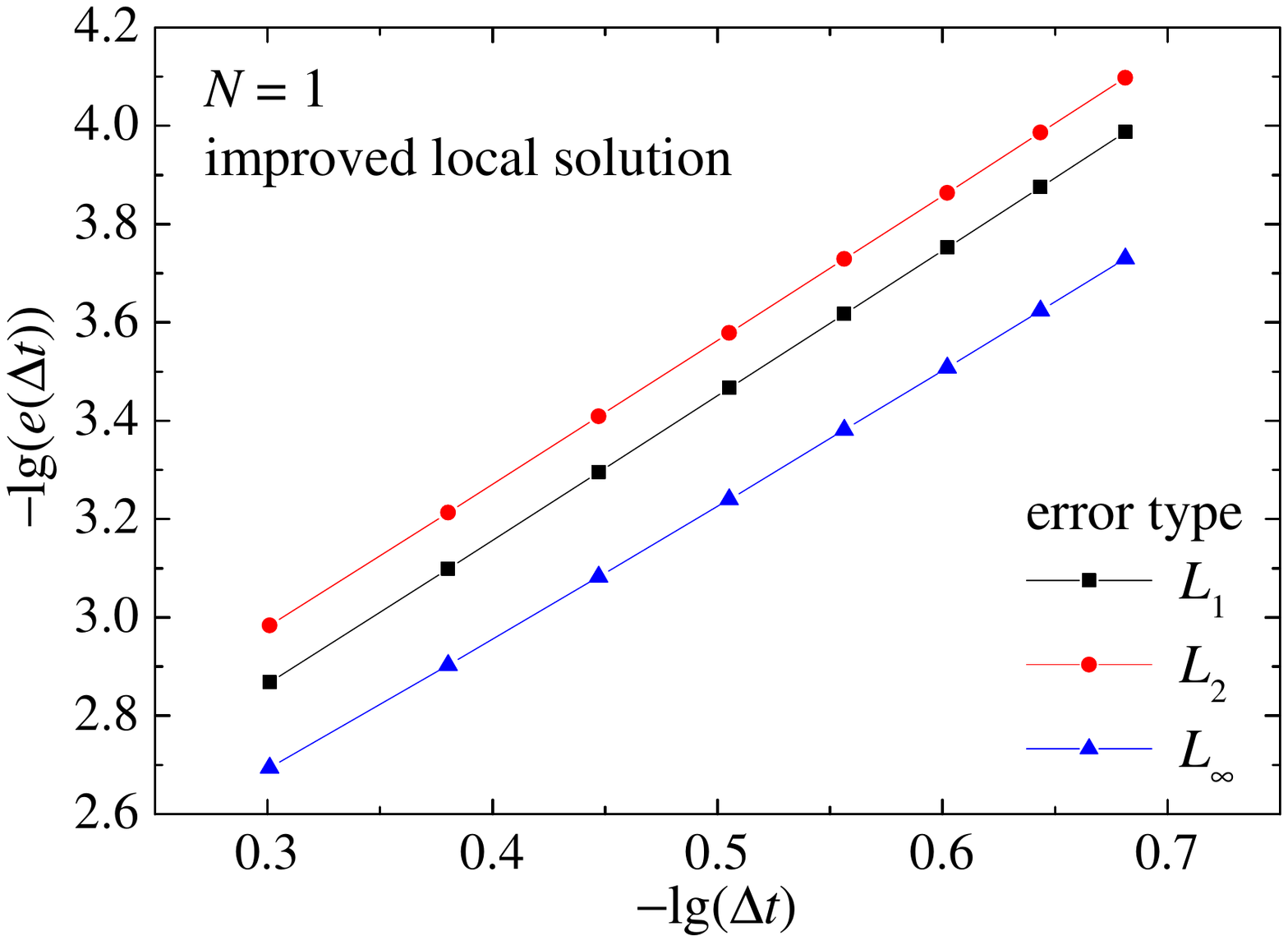}
\vspace{-8mm}\caption{\label{fig:exp_diss:d1}}
\end{subfigure}
\begin{subfigure}{0.24\textwidth}
\includegraphics[width=\textwidth]{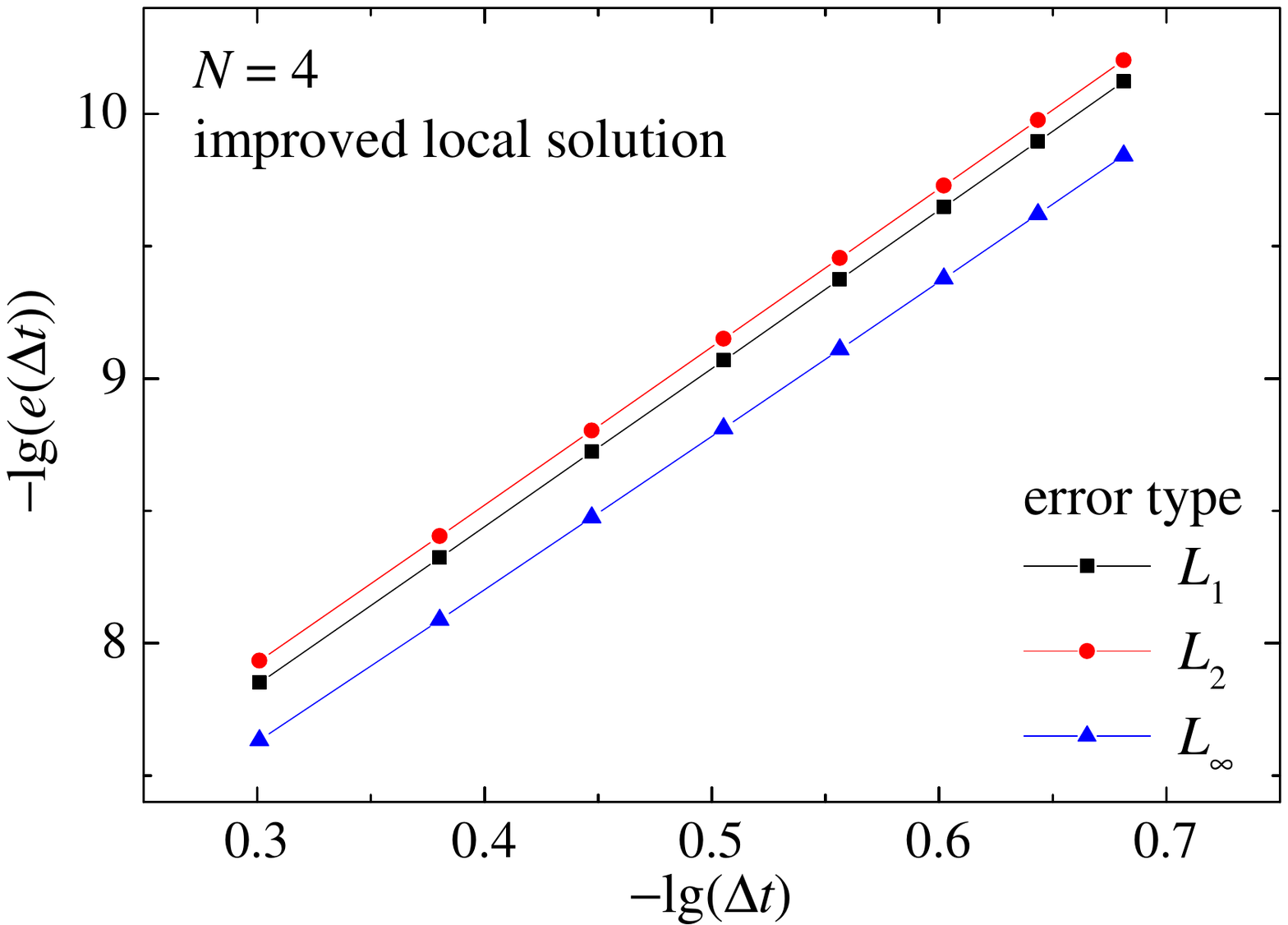}
\vspace{-8mm}\caption{\label{fig:exp_diss:d2}}
\end{subfigure}
\begin{subfigure}{0.24\textwidth}
\includegraphics[width=\textwidth]{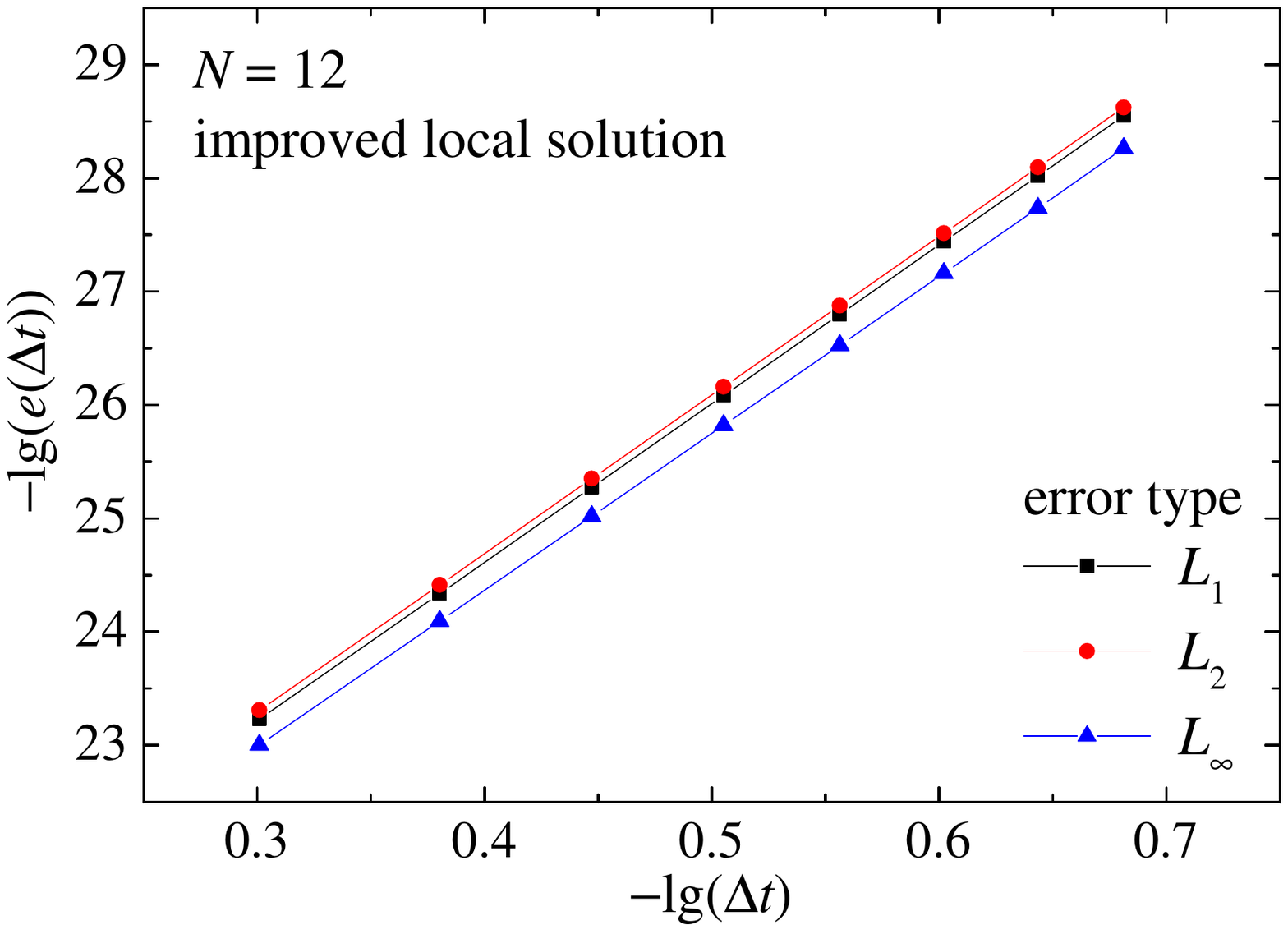}
\vspace{-8mm}\caption{\label{fig:exp_diss:d3}}
\end{subfigure}
\begin{subfigure}{0.24\textwidth}
\includegraphics[width=\textwidth]{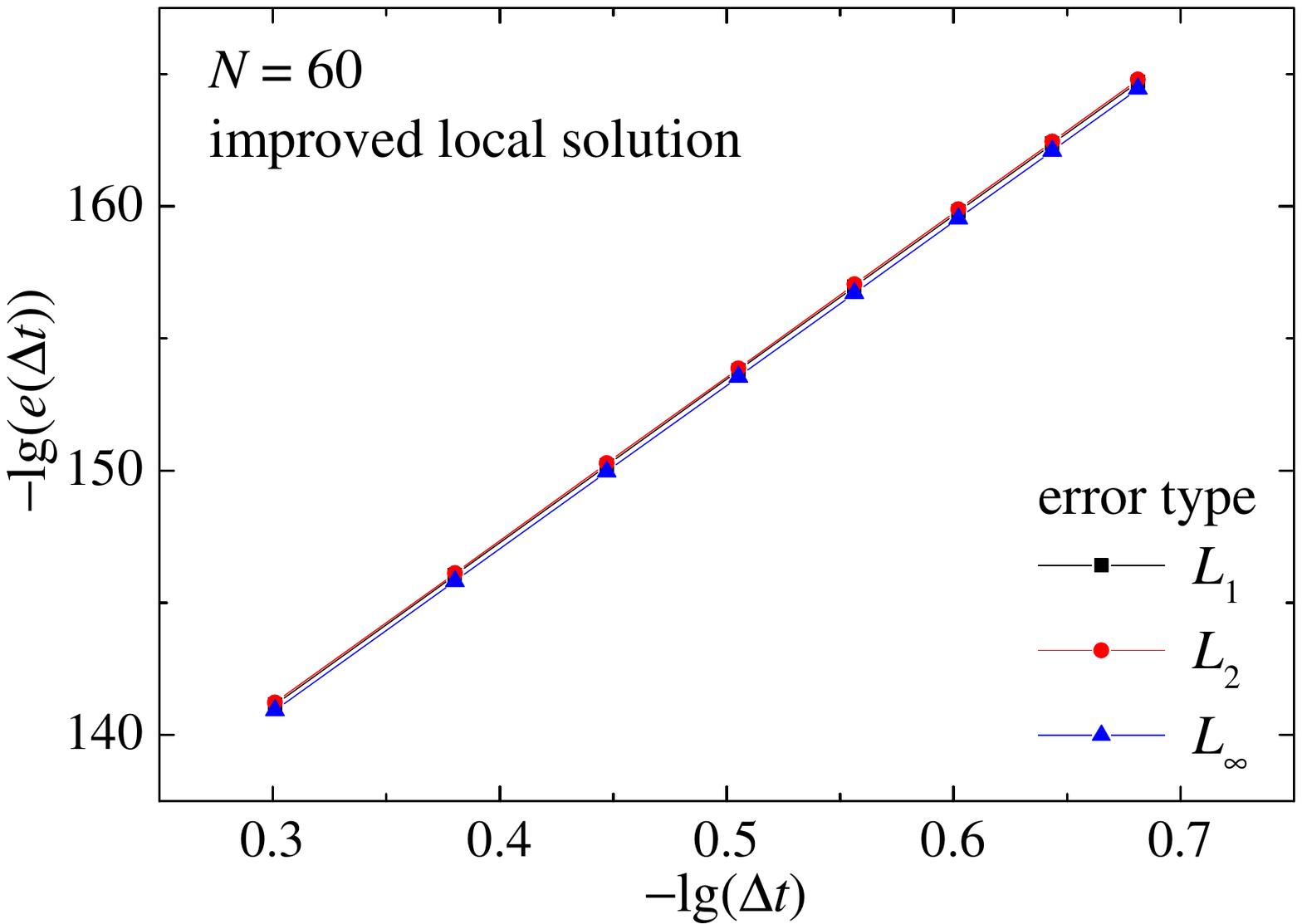}
\vspace{-8mm}\caption{\label{fig:exp_diss:d4}}
\end{subfigure}\\
\begin{subfigure}{0.24\textwidth}
\includegraphics[width=\textwidth]{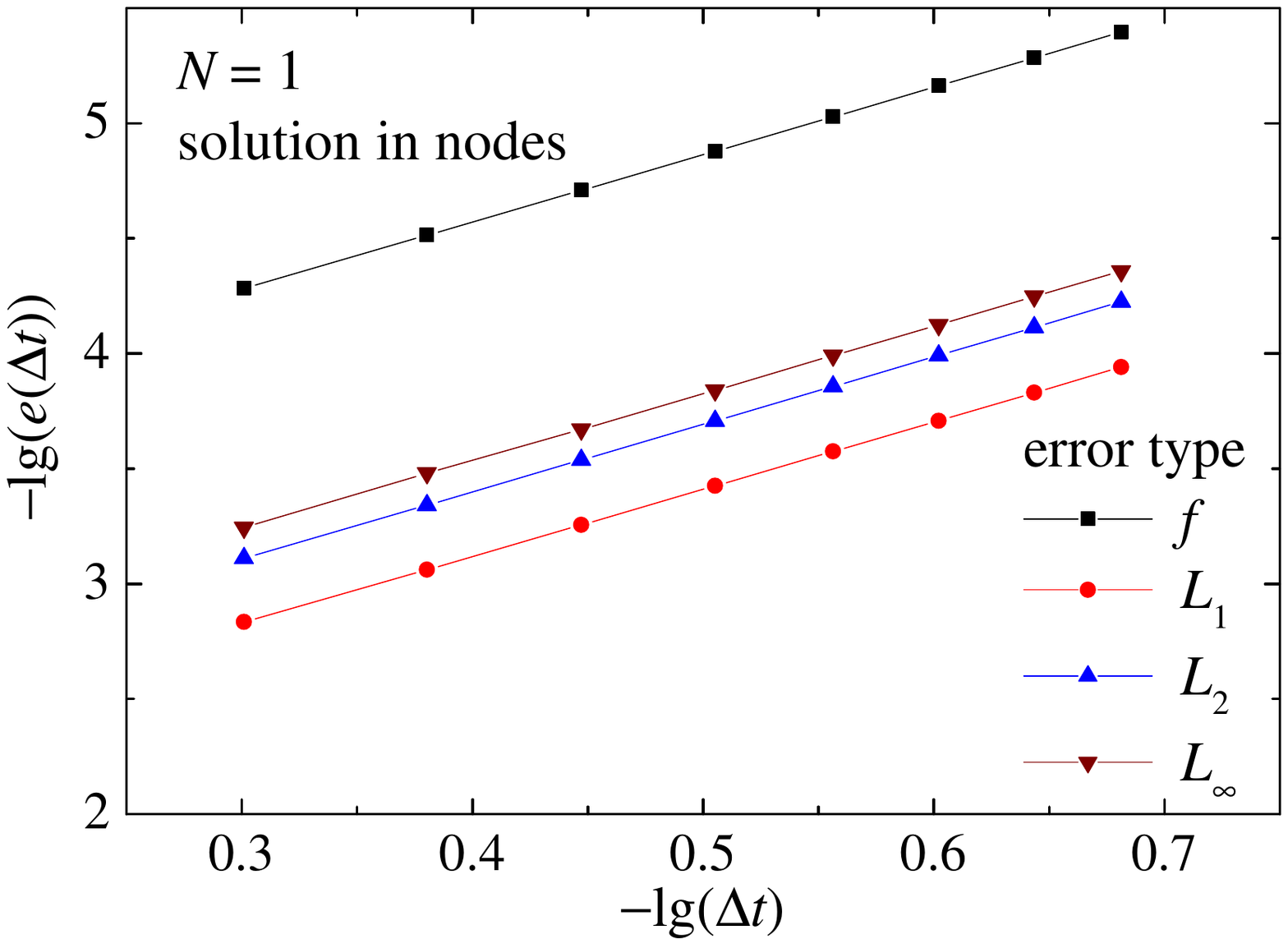}
\vspace{-8mm}\caption{\label{fig:exp_diss:e1}}
\end{subfigure}
\begin{subfigure}{0.24\textwidth}
\includegraphics[width=\textwidth]{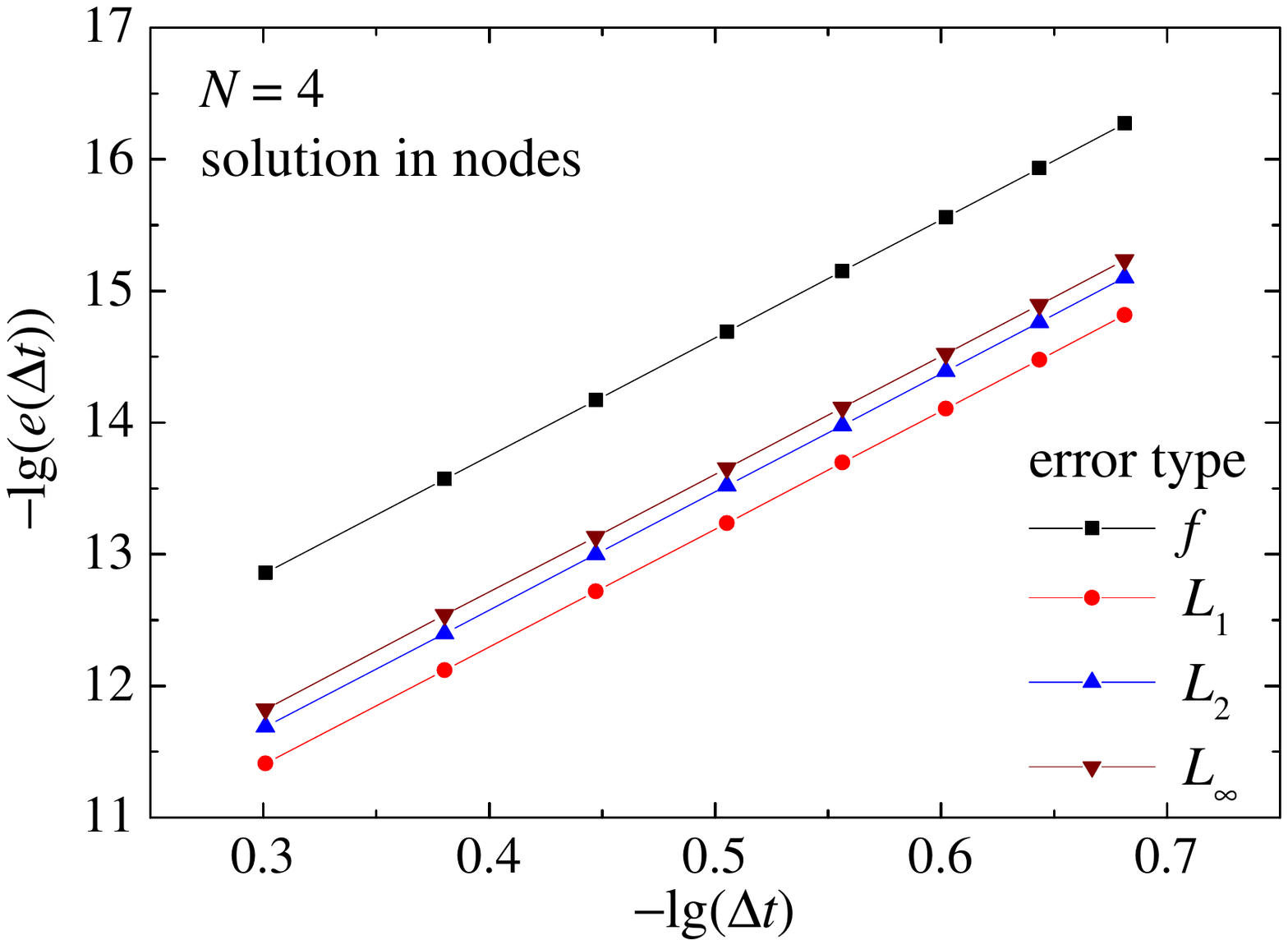}
\vspace{-8mm}\caption{\label{fig:exp_diss:e2}}
\end{subfigure}
\begin{subfigure}{0.24\textwidth}
\includegraphics[width=\textwidth]{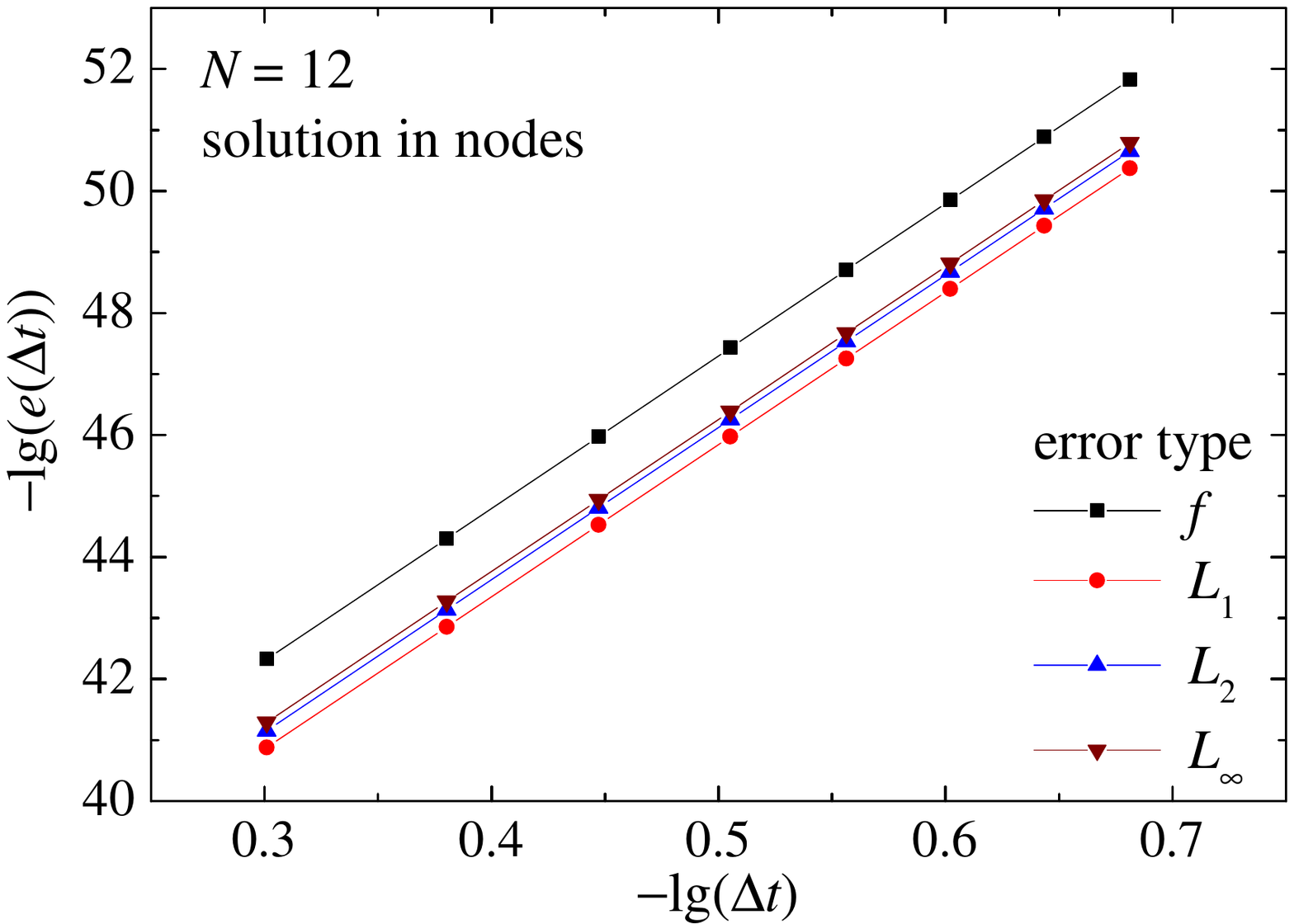}
\vspace{-8mm}\caption{\label{fig:exp_diss:e3}}
\end{subfigure}
\begin{subfigure}{0.24\textwidth}
\includegraphics[width=\textwidth]{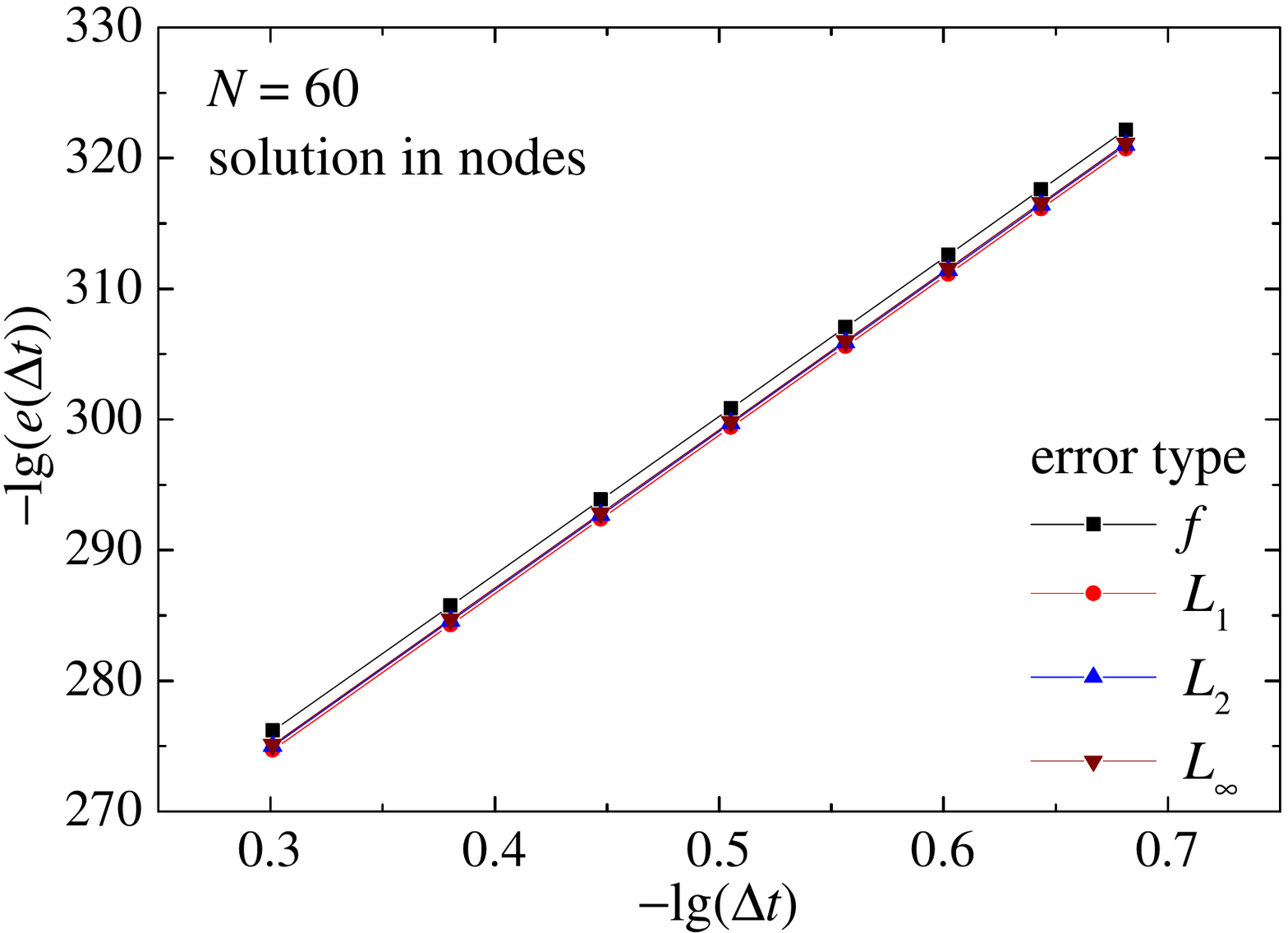}
\vspace{-8mm}\caption{\label{fig:exp_diss:e4}}
\end{subfigure}\\
\caption{%
Numerical solution of the system (\ref{eq:exp_diss_ode}). Comparison of the solution at nodes $u_{n}$, the local solution $u_{L}(t)$, the improved local solution $u_{\rm IL}(t)$ and the exact solution $u^{\rm ex}(t)$ (\subref{fig:exp_diss:a1}, \subref{fig:exp_diss:a2}, \subref{fig:exp_diss:a3}, \subref{fig:exp_diss:a4}) and the errors $\varepsilon(t)$ (\subref{fig:harm_osc:b1}, \subref{fig:harm_osc:b2}, \subref{fig:harm_osc:b3}, \subref{fig:harm_osc:b4}), obtained using polynomials with degrees $N = 1$ (\subref{fig:exp_diss:a1}, \subref{fig:exp_diss:b1}), $N = 4$ (\subref{fig:exp_diss:a2}, \subref{fig:exp_diss:b2}), $N = 12$ (\subref{fig:exp_diss:a3}, \subref{fig:exp_diss:b3}) and $N = 60$ (\subref{fig:exp_diss:a4}, \subref{fig:exp_diss:b4}). Log-log plot of the dependence of the global error for the local solution $e^{l}$ (\subref{fig:exp_diss:c1}, \subref{fig:exp_diss:c2}, \subref{fig:exp_diss:c3}, \subref{fig:exp_diss:c4}), the improved local solution $e^{\rm imp}$ (\subref{fig:exp_diss:d1}, \subref{fig:exp_diss:d2}, \subref{fig:exp_diss:d3}, \subref{fig:exp_diss:d4}) and the solution at nodes $e^{n}$ (\subref{fig:exp_diss:e1}, \subref{fig:exp_diss:e2}, \subref{fig:exp_diss:e3}, \subref{fig:exp_diss:e4}) on the discretization step $\mathrm{\Delta}t$, obtained in the $f$-norm and norms $L_{1}$, $L_{2}$ and $L_{\infty}$, obtained using polynomials with degrees $N = 1$ (\subref{fig:exp_diss:c1}, \subref{fig:exp_diss:d1}, \subref{fig:exp_diss:e1}), $N = 4$ (\subref{fig:exp_diss:c2}, \subref{fig:exp_diss:d2}, \subref{fig:exp_diss:e2}), $N = 12$ (\subref{fig:exp_diss:c3}, \subref{fig:exp_diss:d3}, \subref{fig:exp_diss:e3}) and $N = 60$ (\subref{fig:exp_diss:c4}, \subref{fig:exp_diss:d4}, \subref{fig:exp_diss:e4}).
}
\label{fig:exp_diss}
\end{figure}

The first example of applying the ADER-DG numerical method with a local DG predictor to solving the initial value problem for the ODE system (\ref{eq:ivp_ode_diff_src}) is a trivial case corresponding to the Dahlquist's test equation:
\begin{equation}\label{eq:exp_diss_ode}
\dot{u} + u = 0,\quad
u(0) = 1,\quad t \in [0,\, 5],\quad
u^{\rm ex}(t) = \exp(-t).
\end{equation}
The obtained results are presented in Fig.~\ref{fig:exp_diss} and Table~\ref{tab:conv_orders_exp_diss}.

\begin{table*}[h!]
\centering
\normalsize
\caption{%
Convergence orders $p_{f}$, $p_{L_{1}}$, $p_{L_{2}}$, $p_{L_{\infty}}$, calculated in $f$-norm and norms $L_{1}$, $L_{2}$, $L_{\infty}$, of the numerical solution of the ADER-DG method for the problem (\ref{eq:exp_diss_ode}); $N$ is the degree of the basis polynomials $\varphi_{p}$. Orders $p^{n}$ are calculated for \textit{the numerical solution at the nodes} $\mathbf{u}_{n}$; orders $p^{\rm imp}$ --- for \textit{the improved local solution} $\mathbf{u}_{\rm IL}$; orders $p^{l}$ --- for \textit{the local solution} $\mathbf{u}_{L}$. The theoretical values of convergence order $p_{\rm th.}^{n} = 2N+1$, $p_{\rm th.}^{l} = N+1$ and $p^{\rm imp}_{\rm th.} = N+2$ are presented for comparison.
}
\label{tab:conv_orders_exp_diss}
\setlength{\tabcolsep}{3.5pt}
\begin{tabular}{@{}|l|llll|c|lll|c|lll|c|@{}}
\toprule
$N$ & $p^{n}_{f}$ &
$p^{n}_{L_{1}}$ & $p^{n}_{L_{2}}$ & $p^{n}_{L_{\infty}}$ & $p^{n}_{\rm th.}$ &
$p^{l}_{L_{1}}$ & $p^{l}_{L_{2}}$ & $p^{l}_{L_{\infty}}$ & $p^{l}_{\rm th.}$ &
$p^{\rm imp}_{L_{1}}$ & $p^{\rm imp}_{L_{2}}$ & $p^{\rm imp}_{L_{\infty}}$ & $p^{\rm imp}_{\rm th.}$\\
\midrule
$1$ & $2.93$ & $2.92$ & $2.93$ & $2.93$ & $3$ & $1.94$ & $1.94$ & $1.79$ & $2$ & $2.94$ & $2.93$ & $2.73$ & $3$\\
$2$ & $4.95$ & $4.94$ & $4.95$ & $4.95$ & $5$ & $2.97$ & $2.95$ & $2.81$ & $3$ & $3.96$ & $3.95$ & $3.78$ & $4$\\
$3$ & $6.97$ & $6.95$ & $6.96$ & $6.96$ & $7$ & $3.97$ & $3.96$ & $3.82$ & $4$ & $4.97$ & $4.96$ & $4.82$ & $5$\\
$4$ & $8.97$ & $8.96$ & $8.97$ & $8.97$ & $9$ & $4.98$ & $4.97$ & $4.82$ & $5$ & $5.98$ & $5.97$ & $5.81$ & $6$\\
$5$ & $11.0$ & $11.0$ & $11.0$ & $11.0$ & $11$ & $5.98$ & $5.97$ & $5.82$ & $6$ & $6.98$ & $6.97$ & $6.82$ & $7$\\
$6$ & $13.0$ & $13.0$ & $13.0$ & $13.0$ & $13$ & $6.98$ & $6.97$ & $6.83$ & $7$ & $7.98$ & $7.97$ & $7.82$ & $8$\\
$7$ & $15.0$ & $15.0$ & $15.0$ & $15.0$ & $15$ & $7.98$ & $7.97$ & $7.83$ & $8$ & $8.98$ & $8.97$ & $8.83$ & $9$\\
$8$ & $17.0$ & $17.0$ & $17.0$ & $17.0$ & $17$ & $8.98$ & $8.97$ & $8.83$ & $9$ & $9.98$ & $9.97$ & $9.82$ & $10$\\
$9$ & $19.0$ & $19.0$ & $19.0$ & $19.0$ & $19$ & $9.98$ & $9.97$ & $9.83$ & $10$ & $11.0$ & $11.0$ & $10.8$ & $11$\\
$10$ & $21.0$ & $21.0$ & $21.0$ & $21.0$ & $21$ & $11.0$ & $11.0$ & $10.8$ & $11$ & $12.0$ & $12.0$ & $11.8$ & $12$\\
\midrule
$11$ & $23.0$ & $23.0$ & $23.0$ & $23.0$ & $23$ & $12.0$ & $12.0$ & $11.8$ & $12$ & $13.0$ & $13.0$ & $12.8$ & $13$\\
$12$ & $25.0$ & $25.0$ & $25.0$ & $25.0$ & $25$ & $13.0$ & $13.0$ & $12.8$ & $13$ & $14.0$ & $14.0$ & $13.8$ & $14$\\
$13$ & $27.0$ & $27.0$ & $27.0$ & $27.0$ & $27$ & $14.0$ & $14.0$ & $13.8$ & $14$ & $15.0$ & $15.0$ & $14.8$ & $15$\\
$14$ & $29.0$ & $29.0$ & $29.0$ & $29.0$ & $29$ & $15.0$ & $15.0$ & $14.8$ & $15$ & $16.0$ & $16.0$ & $15.8$ & $16$\\
$15$ & $31.0$ & $31.0$ & $31.0$ & $31.0$ & $31$ & $16.0$ & $16.0$ & $15.8$ & $16$ & $17.0$ & $17.0$ & $16.8$ & $17$\\
$16$ & $33.0$ & $33.0$ & $33.0$ & $33.0$ & $33$ & $17.0$ & $17.0$ & $16.8$ & $17$ & $18.0$ & $18.0$ & $17.8$ & $18$\\
$17$ & $35.0$ & $35.0$ & $35.0$ & $35.0$ & $35$ & $18.0$ & $18.0$ & $17.8$ & $18$ & $19.0$ & $19.0$ & $18.8$ & $19$\\
$18$ & $37.0$ & $37.0$ & $37.0$ & $37.0$ & $37$ & $19.0$ & $19.0$ & $18.8$ & $19$ & $20.0$ & $20.0$ & $19.8$ & $20$\\
$19$ & $39.0$ & $39.0$ & $39.0$ & $39.0$ & $39$ & $20.0$ & $20.0$ & $19.8$ & $20$ & $21.0$ & $21.0$ & $20.8$ & $21$\\
$20$ & $41.0$ & $41.0$ & $41.0$ & $41.0$ & $41$ & $21.0$ & $21.0$ & $20.8$ & $21$ & $22.0$ & $22.0$ & $21.8$ & $22$\\
\midrule
$21$ & $43.0$ & $43.0$ & $43.0$ & $43.0$ & $43$ & $22.0$ & $22.0$ & $21.8$ & $22$ & $23.0$ & $23.0$ & $22.8$ & $23$\\
$22$ & $45.0$ & $45.0$ & $45.0$ & $45.0$ & $45$ & $23.0$ & $23.0$ & $22.8$ & $23$ & $24.0$ & $24.0$ & $23.8$ & $24$\\
$23$ & $47.0$ & $47.0$ & $47.0$ & $47.0$ & $47$ & $24.0$ & $24.0$ & $23.8$ & $24$ & $25.0$ & $25.0$ & $24.8$ & $25$\\
$24$ & $49.0$ & $49.0$ & $49.0$ & $49.0$ & $49$ & $25.0$ & $25.0$ & $24.8$ & $25$ & $26.0$ & $26.0$ & $25.8$ & $26$\\
$25$ & $51.0$ & $51.0$ & $51.0$ & $51.0$ & $51$ & $26.0$ & $26.0$ & $25.8$ & $26$ & $27.0$ & $27.0$ & $26.8$ & $27$\\
$26$ & $53.0$ & $53.0$ & $53.0$ & $53.0$ & $53$ & $27.0$ & $27.0$ & $26.8$ & $27$ & $28.0$ & $28.0$ & $27.8$ & $28$\\
$27$ & $55.0$ & $55.0$ & $55.0$ & $55.0$ & $55$ & $28.0$ & $28.0$ & $27.8$ & $28$ & $29.0$ & $29.0$ & $28.8$ & $29$\\
$28$ & $57.0$ & $57.0$ & $57.0$ & $57.0$ & $57$ & $29.0$ & $29.0$ & $28.8$ & $29$ & $30.0$ & $30.0$ & $29.8$ & $30$\\
$29$ & $59.0$ & $59.0$ & $59.0$ & $59.0$ & $59$ & $30.0$ & $30.0$ & $29.8$ & $30$ & $31.0$ & $31.0$ & $30.8$ & $31$\\
$30$ & $61.0$ & $61.0$ & $61.0$ & $61.0$ & $61$ & $31.0$ & $31.0$ & $30.8$ & $31$ & $32.0$ & $32.0$ & $31.8$ & $32$\\
\midrule
$35$ & $71.0$ & $71.0$ & $71.0$ & $71.0$ & $71$ & $36.0$ & $36.0$ & $35.8$ & $36$ & $37.0$ & $37.0$ & $36.8$ & $37$\\
$40$ & $81.0$ & $81.0$ & $81.0$ & $81.0$ & $81$ & $41.0$ & $41.0$ & $40.8$ & $41$ & $42.0$ & $42.0$ & $41.8$ & $42$\\
$45$ & $91.0$ & $91.0$ & $91.0$ & $91.0$ & $91$ & $46.0$ & $46.0$ & $45.8$ & $46$ & $47.0$ & $47.0$ & $46.8$ & $47$\\
$50$ & $101.0$ & $101.0$ & $101.0$ & $101.0$ & $101$ & $51.0$ & $51.0$ & $50.8$ & $51$ & $52.0$ & $52.0$ & $51.8$ & $52$\\
$55$ & $111.0$ & $111.0$ & $111.0$ & $111.0$ & $111$ & $56.0$ & $56.0$ & $55.8$ & $56$ & $57.0$ & $57.0$ & $56.8$ & $57$\\
$60$ & $121.0$ & $121.0$ & $121.0$ & $121.0$ & $121$ & $61.0$ & $61.0$ & $60.8$ & $61$ & $62.0$ & $62.0$ & $61.8$ & $62$\\
\bottomrule
\end{tabular}
\end{table*}

Fig.~\ref{fig:exp_diss} shows the dependencies of the numerical solutions $u_{L}$, $u_{\rm IL}$, $u_{n}$ and the exact analytical solution $u^{\rm ex}$, the dependencies of the local error $\varepsilon$ (\ref{eq:eps_local_def}) of the numerical solutions, and the dependence of the global error $e$ (\ref{eq:eps_un_global_def}), (\ref{eq:eps_ul_global_def}) of the numerical solutions on the discretization step ${\Delta t}$, for polynomial degrees $N = 1$, $4$, $12$ and $60$. A comparison of the obtained dependencies of the numerical solutions $u_{L}(t)$, $u_{\rm IL}(t)$, $u_{n}$ with the exact analytical solution $u^{\rm ex}(t)$, presented in Fig.~\ref{fig:exp_diss} (\subref{fig:exp_diss:a1}, \subref{fig:exp_diss:a2}, \subref{fig:exp_diss:a3}, \subref{fig:exp_diss:a4}), demonstrates a high-quality agreement. The local error $\varepsilon$ (\ref{eq:eps_local_def}) of the numerical solutions shows that for polynomial degree $N = 1$, the errors of the local solution $\mathbf{u}_{L}$ and the improved local solution $\mathbf{u}_{\rm IL}$ are approximately comparable throughout the entire domain, but in the initial region, the improved local solution exhibits higher accuracy. It is evident that the local error $\varepsilon(t)$ decreases with increasing $t$, which is related to the exponential decrease in the solution $u(t)$ itself. With increasing polynomial degree $N$, in particular, in the cases of polynomial degrees $N = 4$, $12$, $60$, the local error $\varepsilon$ of the improved local solution $u_{\rm IL}$ is significantly lower than the local error $\varepsilon$ of the local solution $u_{L}$, with the differences reaching $2$--$4$ orders of magnitude. Moreover, the local error $\varepsilon$ of the numerical solution $u_{n}$ at the grid nodes $t_{n}$ is significantly smaller than the local errors $\varepsilon$ of the local solution $u_{L}$ and the improved local solution $u_{\rm IL}$, amounting to approximately $1$--$3$ orders of magnitude in the case of polynomial degree $N = 1$, approximately $4$--$6$ orders of magnitude in the case of polynomial degree $N = 4$, approximately $20$--$22$ orders of magnitude in the case of polynomial degree $N = 12$ and approximately $130$--$138$ orders of magnitude in the case of polynomial degree $N = 60$ (for this purpose, breaks in the graph along the vertical axis are inserted in Fig.~\ref{fig:exp_diss} (\subref{fig:exp_diss:b3}, \subref{fig:exp_diss:b4})).

The resulting dependencies of the global error $e$ (\ref{eq:eps_un_global_def}), (\ref{eq:eps_ul_global_def}) on the discretization step ${\Delta t}$, presented in Fig.~\ref{fig:exp_diss} (\subref{fig:exp_diss:c1}, \subref{fig:exp_diss:c2}, \subref{fig:exp_diss:c3}, \subref{fig:exp_diss:c4}, \subref{fig:exp_diss:d1}, \subref{fig:exp_diss:d2}, \subref{fig:exp_diss:d3}, \subref{fig:exp_diss:d4}, \subref{fig:exp_diss:e1}, \subref{fig:exp_diss:e2}, \subref{fig:exp_diss:e3}, \subref{fig:exp_diss:e4}), demonstrate a very high quality of linear approximation for all studied polynomial degrees $N$ (Fig.~\ref{fig:exp_diss} only shows results for polynomial degrees $N = 1$, $4$, $12$ and $60$). Based on the approximation of the obtained dependencies $e({\Delta t})$ in log-log scale by a linear function $\lg{e({\Delta t})} \propto p\cdot\lg{{\Delta t}}$, empirical convergence orders $p$ are calculated and presented in Table~\ref{tab:conv_orders_exp_diss} for all polynomial degrees $N = 1, \ldots, 30$ and polynomial degrees $N$ up to $60$ with a step of $5$. A comparison of the obtained empirical convergence orders $p$ with the expected theoretical values $p_{\rm th.}$ (\ref{eq:conv_ords_exp}) shows excellent agreement. It is particularly noticeable that the empirical convergence orders $p^{\rm imp}$ of the improved local numerical solution $u_{\rm IL}$ are one unit higher than the empirical convergence orders $p^{l}$ of the local numerical solution $u_{L}$.

The obtained results allowed to conclude that the ADER-DG numerical method with a local DG predictor provides a highly accurate numerical solution to the initial value problem for the ODE system (\ref{eq:exp_diss_ode}) presented in this example. The obtained results are in good agreement with the theory developed above. The improved local numerical solution $u_{\rm IL}$ indeed demonstrates higher accuracy and a higher convergence order compared to the local numerical solution $u_{L}$.

\subsection{Example 2: Linear exp--test}
\label{sec:apps:lin_diss}

\begin{figure}[h!]
\captionsetup[subfigure]{%
	position=bottom,
	font+=smaller,
	textfont=normalfont,
	singlelinecheck=off,
	justification=raggedright
}
\centering
\begin{subfigure}{0.24\textwidth}
\includegraphics[width=\textwidth]{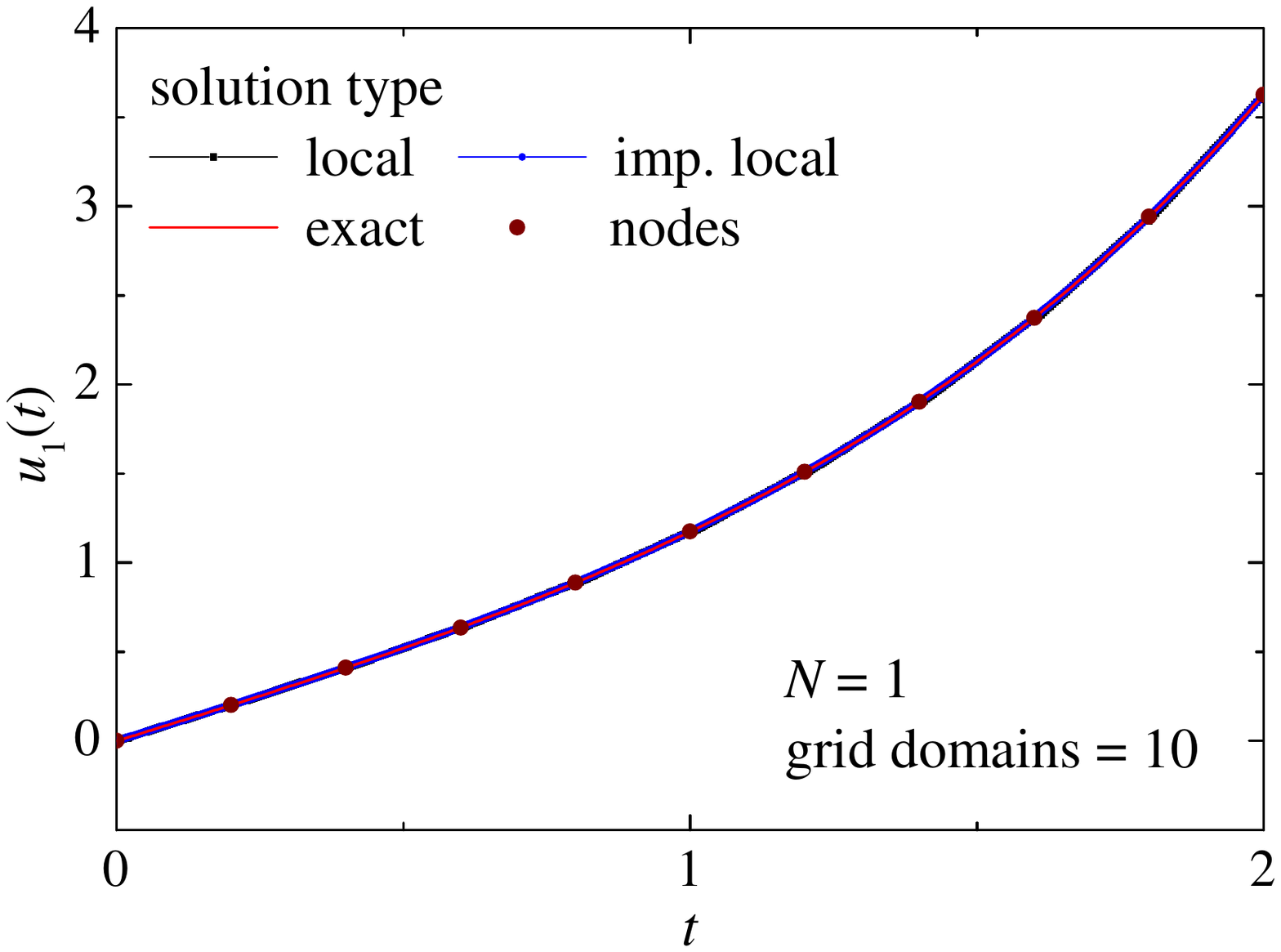}
\vspace{-8mm}\caption{\label{fig:lin_diss:a1}}
\end{subfigure}
\begin{subfigure}{0.24\textwidth}
\includegraphics[width=\textwidth]{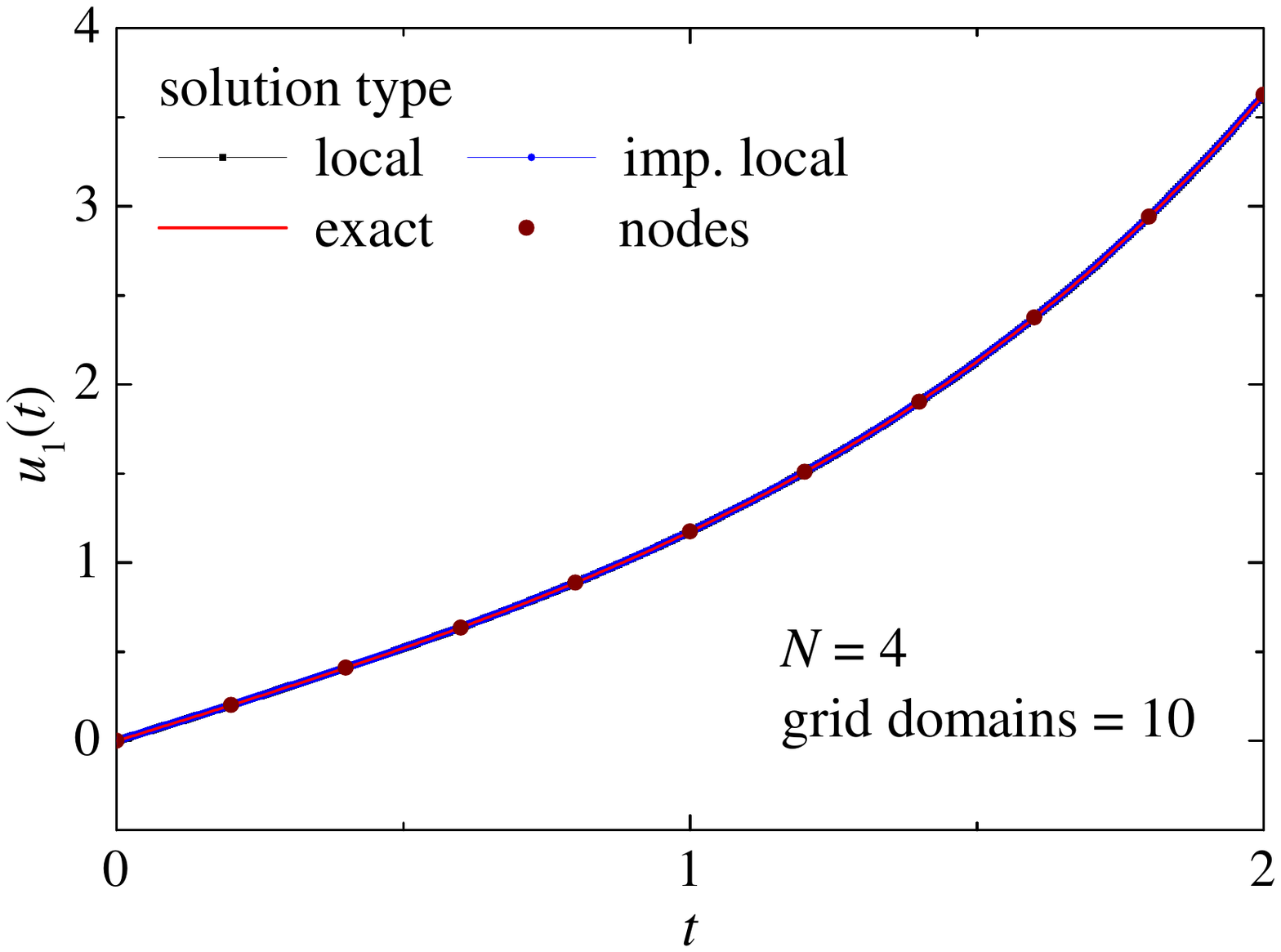}
\vspace{-8mm}\caption{\label{fig:lin_diss:a2}}
\end{subfigure}
\begin{subfigure}{0.24\textwidth}
\includegraphics[width=\textwidth]{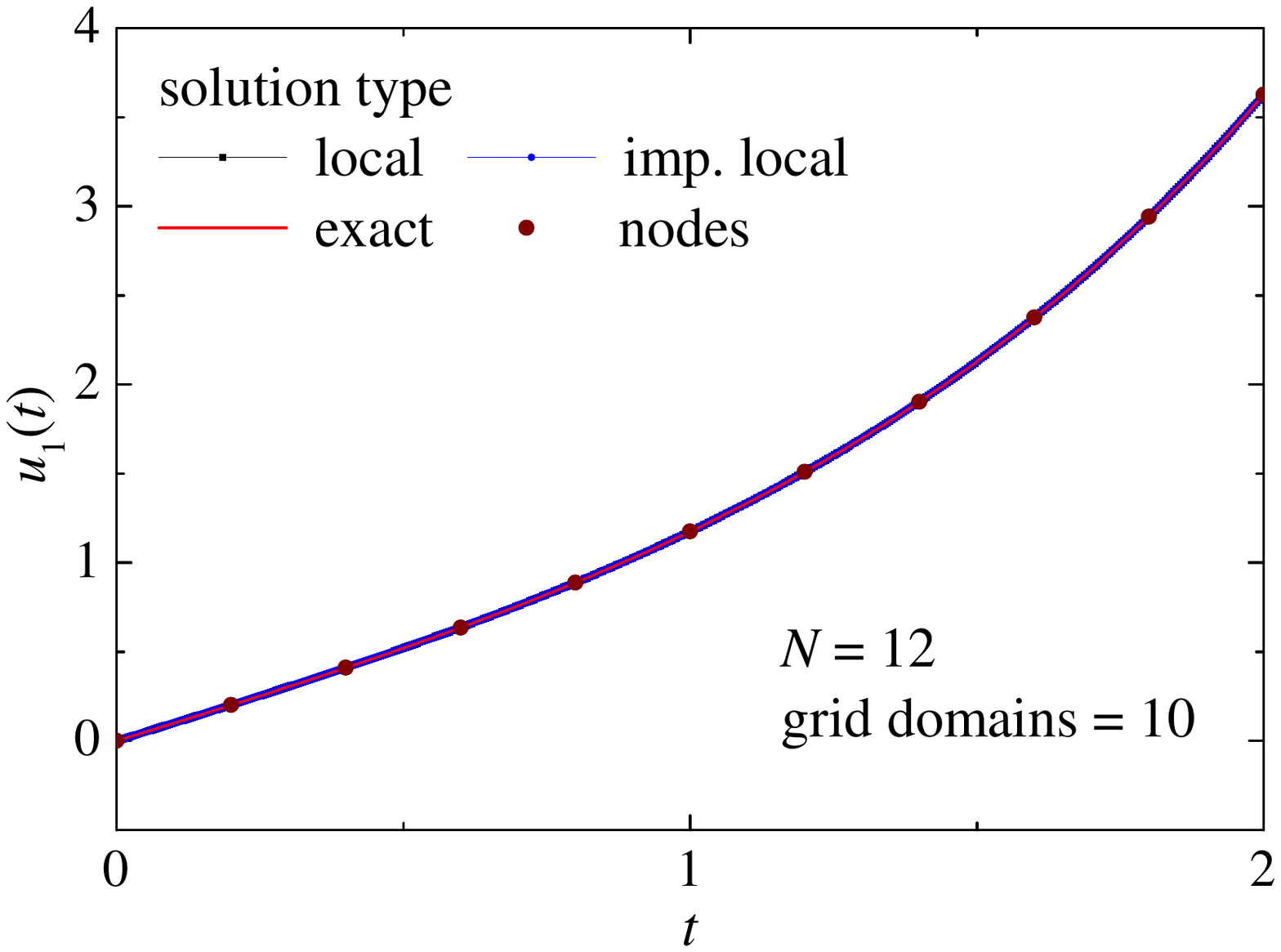}
\vspace{-8mm}\caption{\label{fig:lin_diss:a3}}
\end{subfigure}
\begin{subfigure}{0.24\textwidth}
\includegraphics[width=\textwidth]{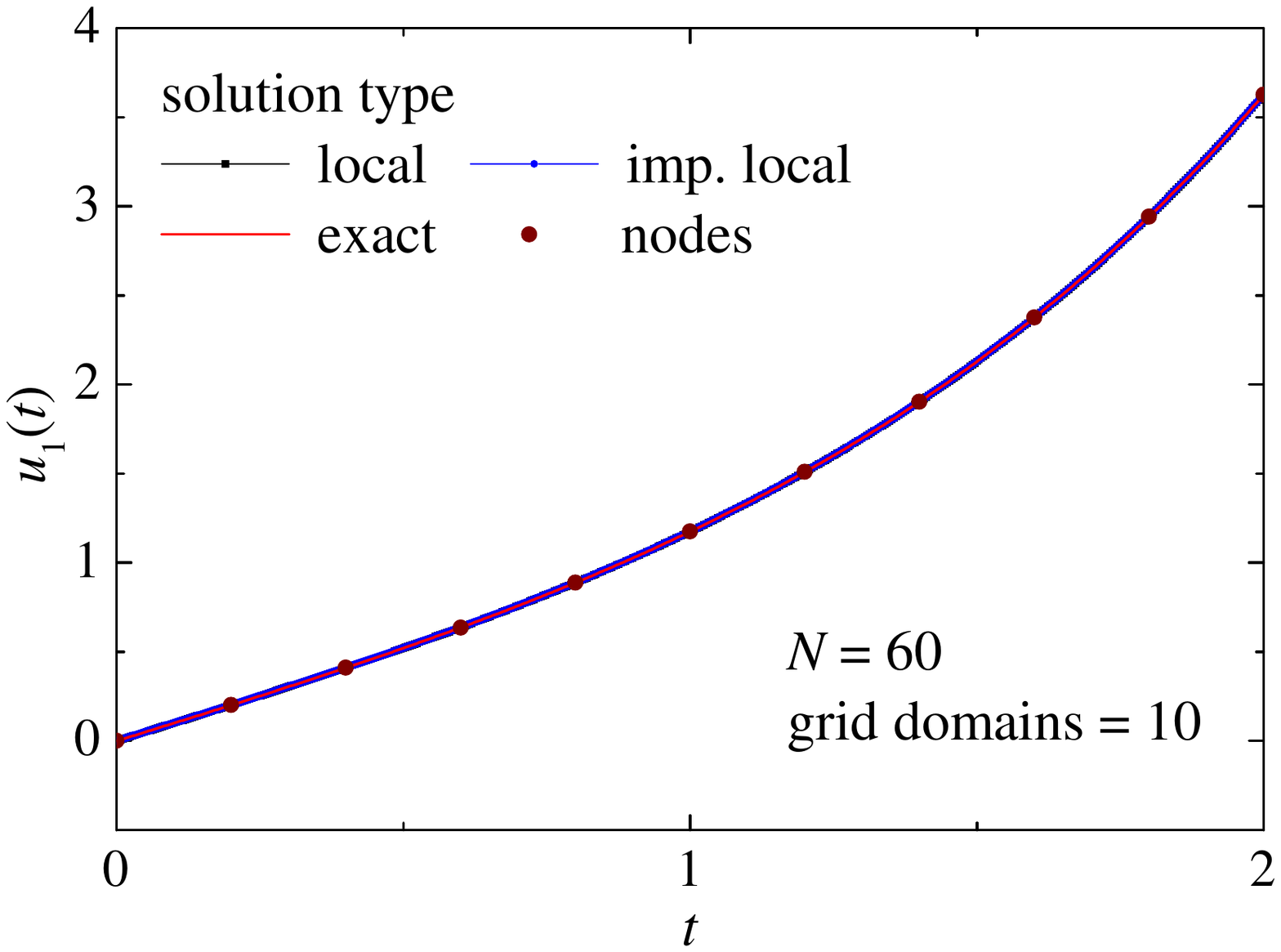}
\vspace{-8mm}\caption{\label{fig:lin_diss:a4}}
\end{subfigure}\\
\begin{subfigure}{0.24\textwidth}
\includegraphics[width=\textwidth]{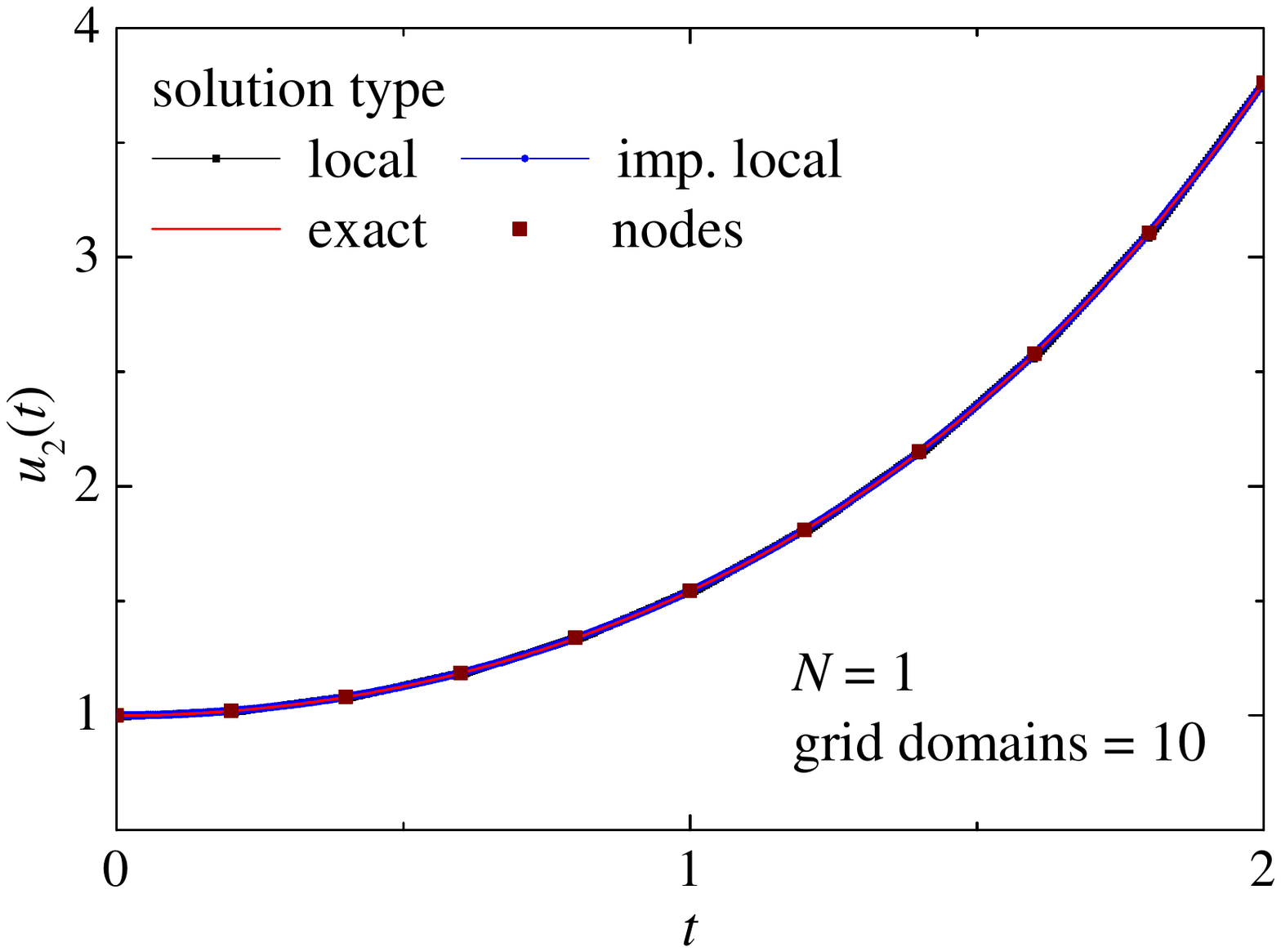}
\vspace{-8mm}\caption{\label{fig:lin_diss:b1}}
\end{subfigure}
\begin{subfigure}{0.24\textwidth}
\includegraphics[width=\textwidth]{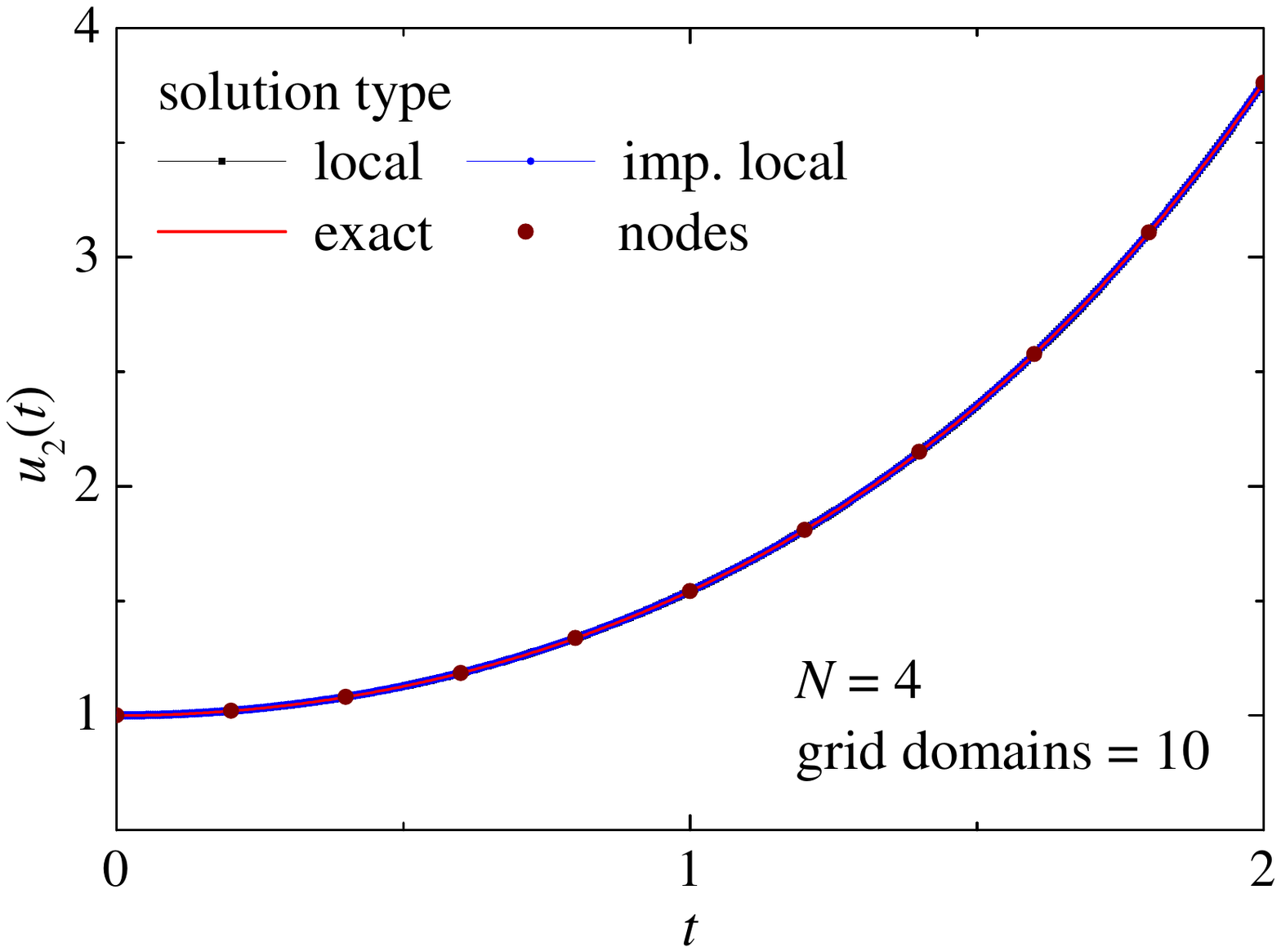}
\vspace{-8mm}\caption{\label{fig:lin_diss:b2}}
\end{subfigure}
\begin{subfigure}{0.24\textwidth}
\includegraphics[width=\textwidth]{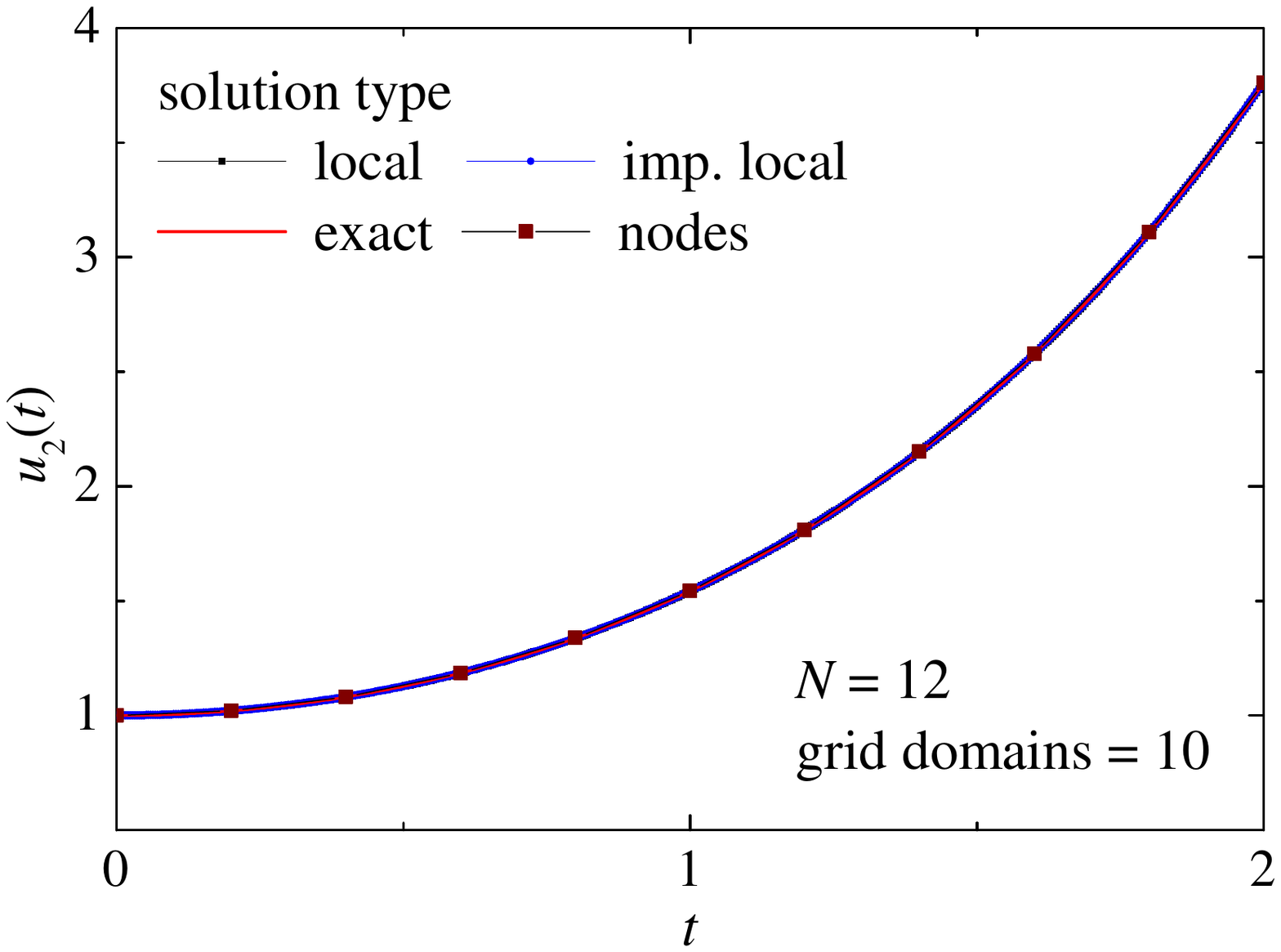}
\vspace{-8mm}\caption{\label{fig:lin_diss:b3}}
\end{subfigure}
\begin{subfigure}{0.24\textwidth}
\includegraphics[width=\textwidth]{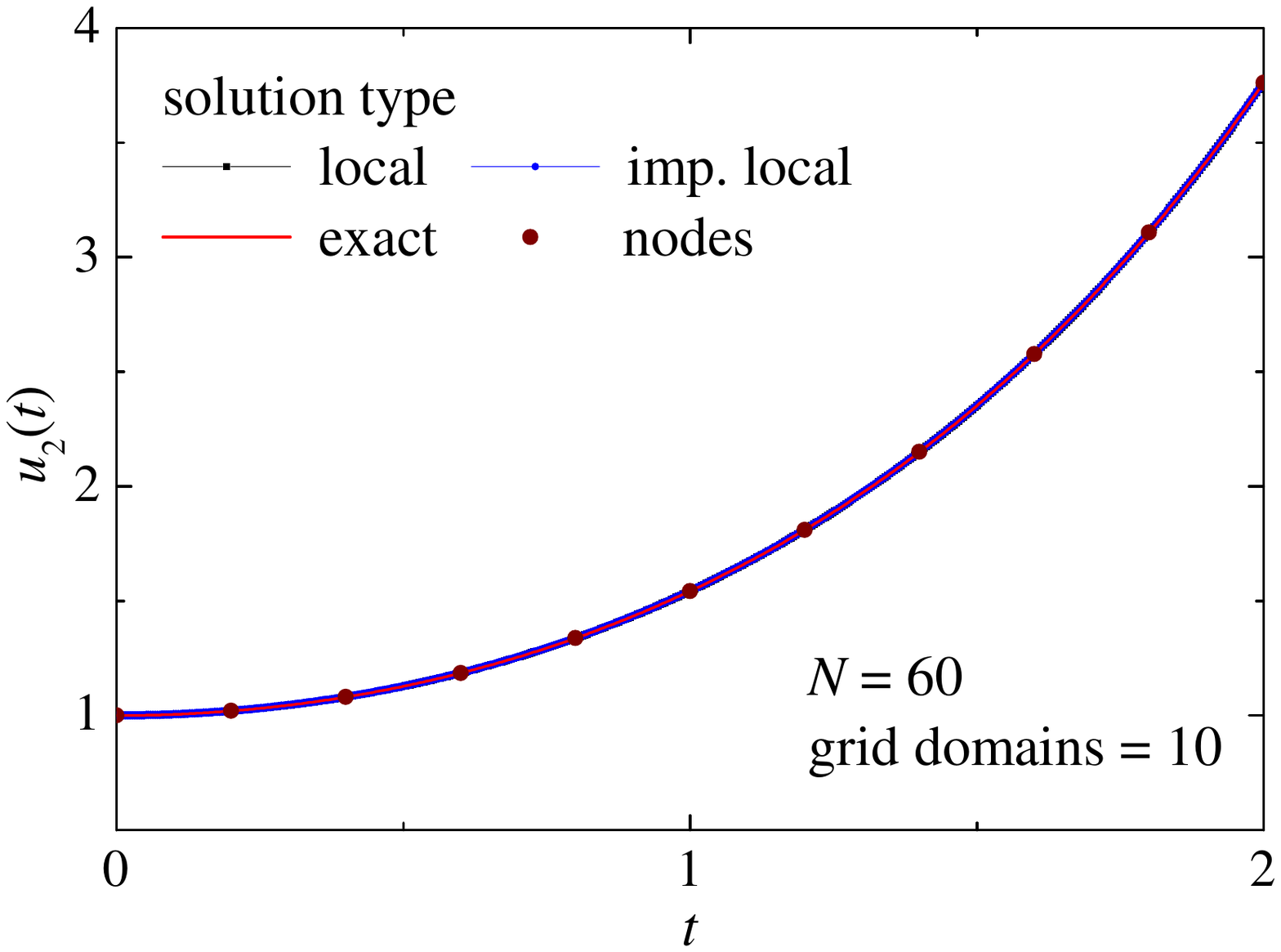}
\vspace{-8mm}\caption{\label{fig:lin_diss:b4}}
\end{subfigure}\\
\begin{subfigure}{0.24\textwidth}
\includegraphics[width=\textwidth]{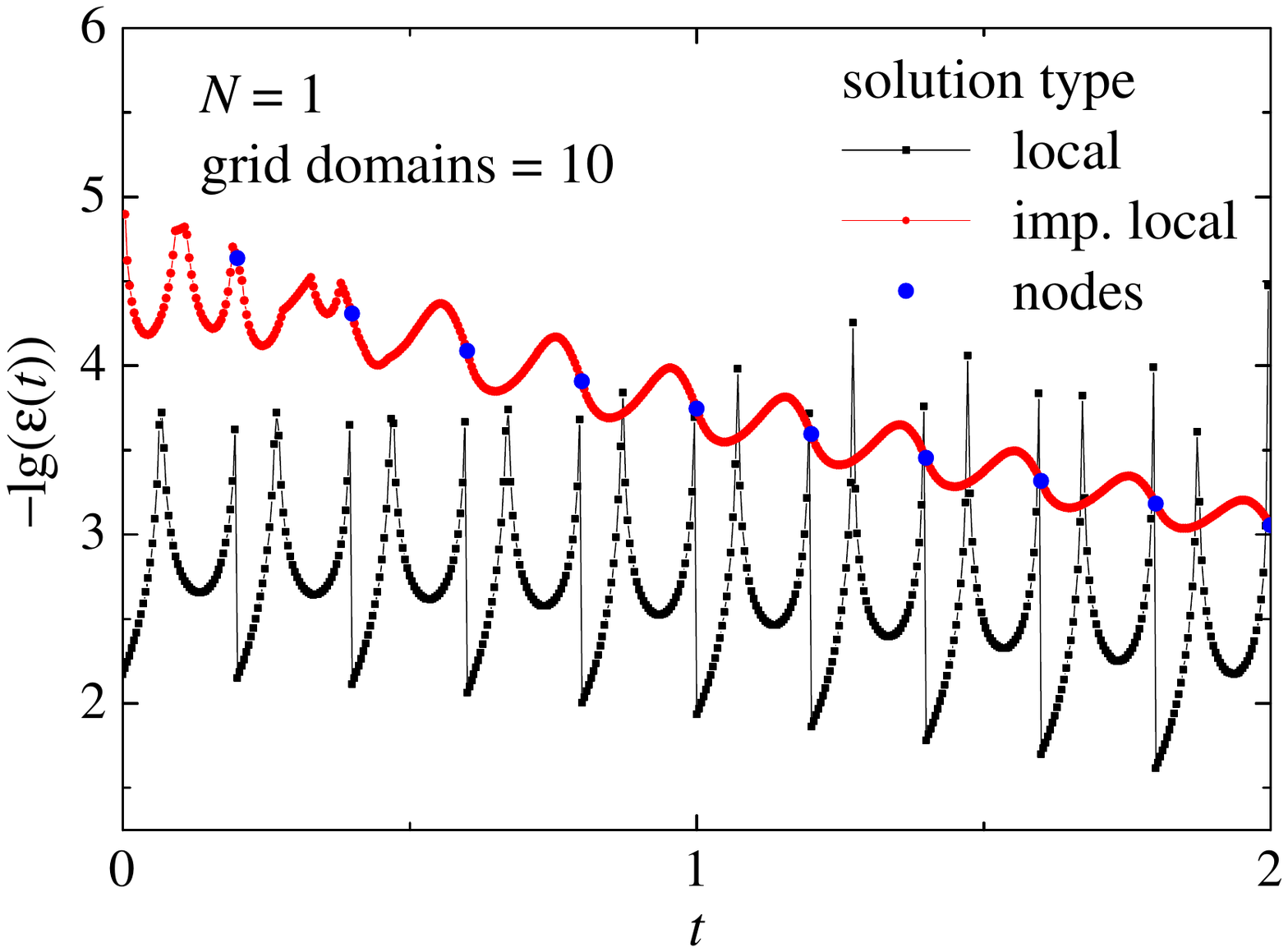}
\vspace{-8mm}\caption{\label{fig:lin_diss:c1}}
\end{subfigure}
\begin{subfigure}{0.24\textwidth}
\includegraphics[width=\textwidth]{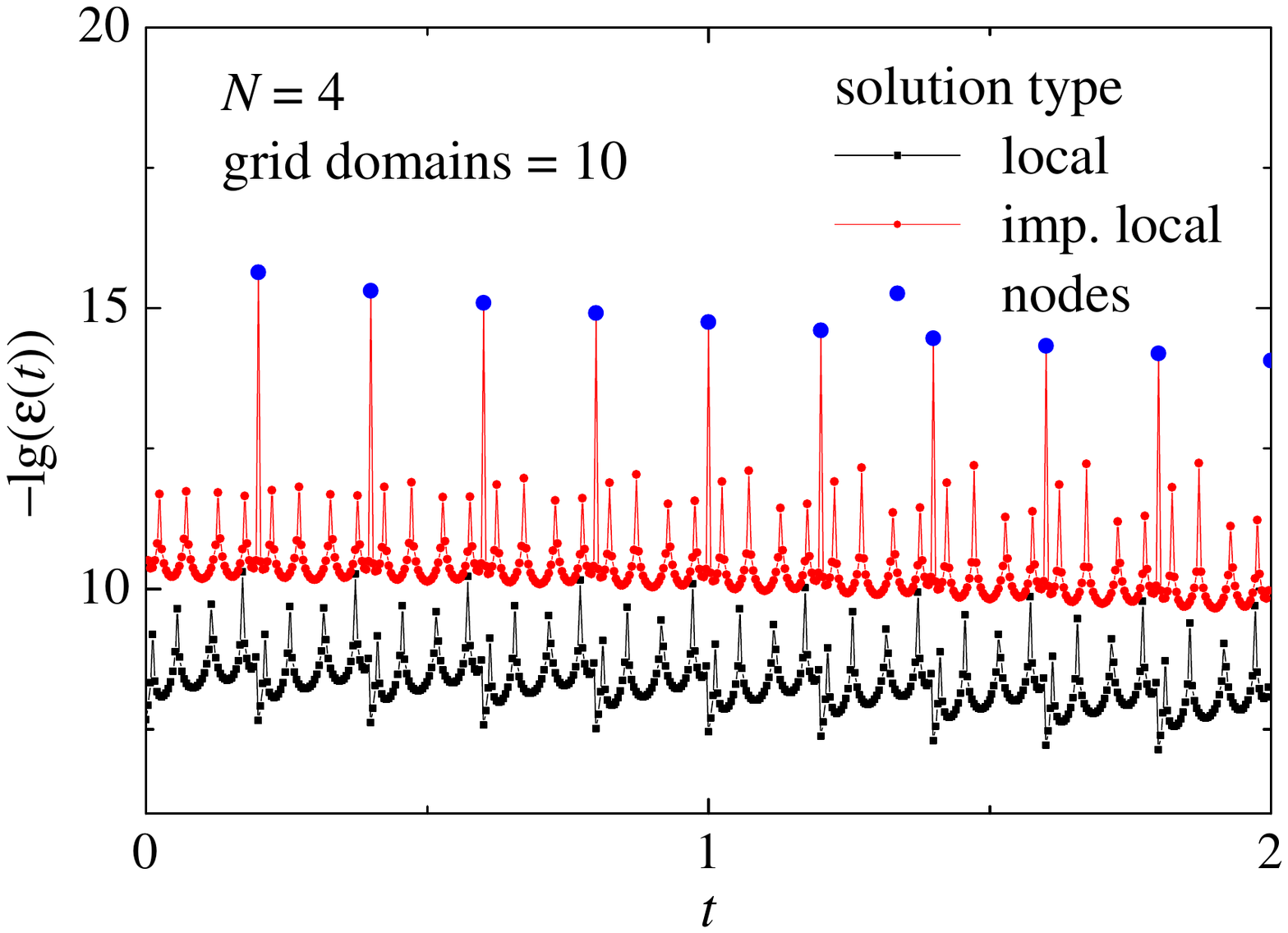}
\vspace{-8mm}\caption{\label{fig:lin_diss:c2}}
\end{subfigure}
\begin{subfigure}{0.24\textwidth}
\includegraphics[width=\textwidth]{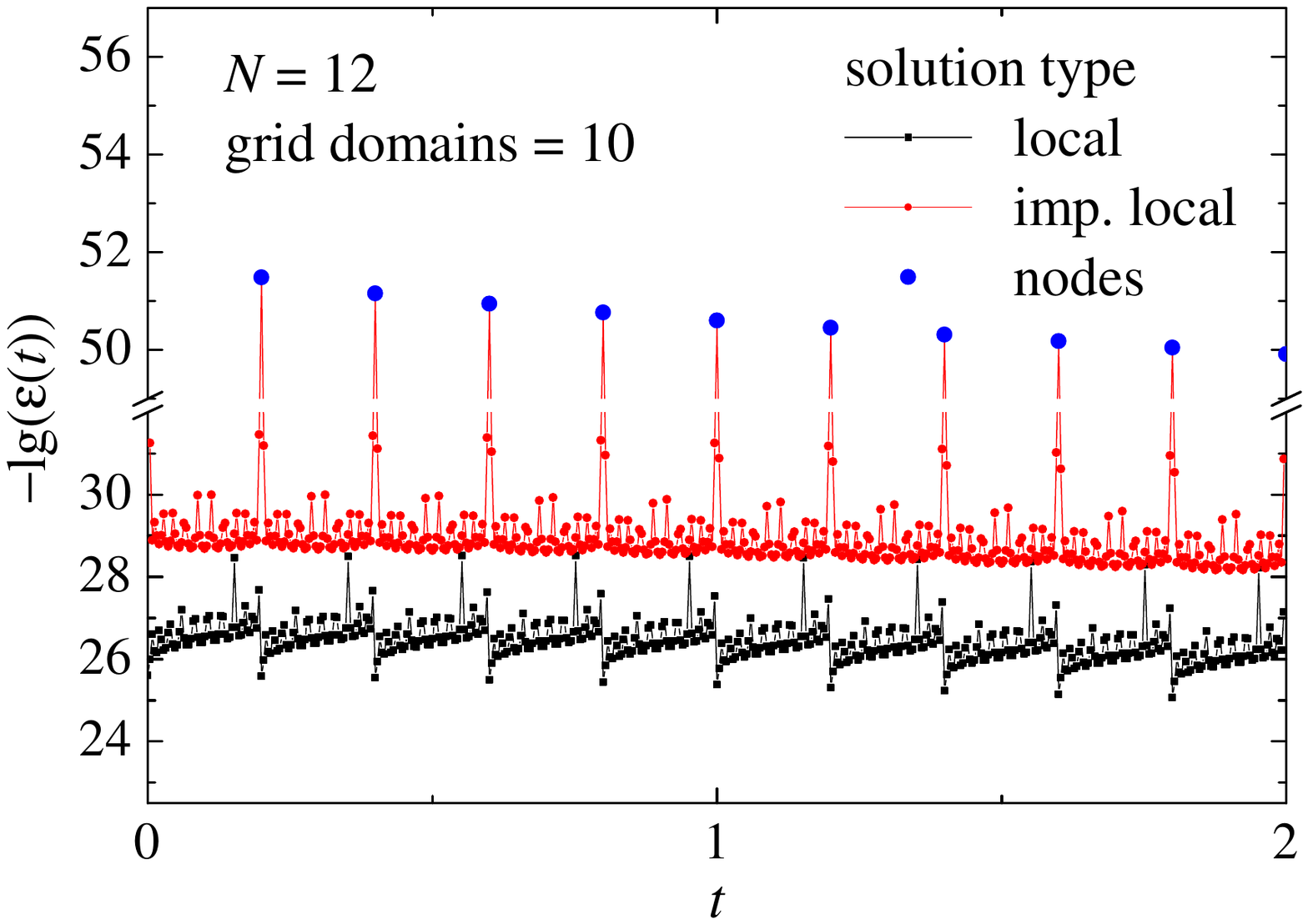}
\vspace{-8mm}\caption{\label{fig:lin_diss:c3}}
\end{subfigure}
\begin{subfigure}{0.24\textwidth}
\includegraphics[width=\textwidth]{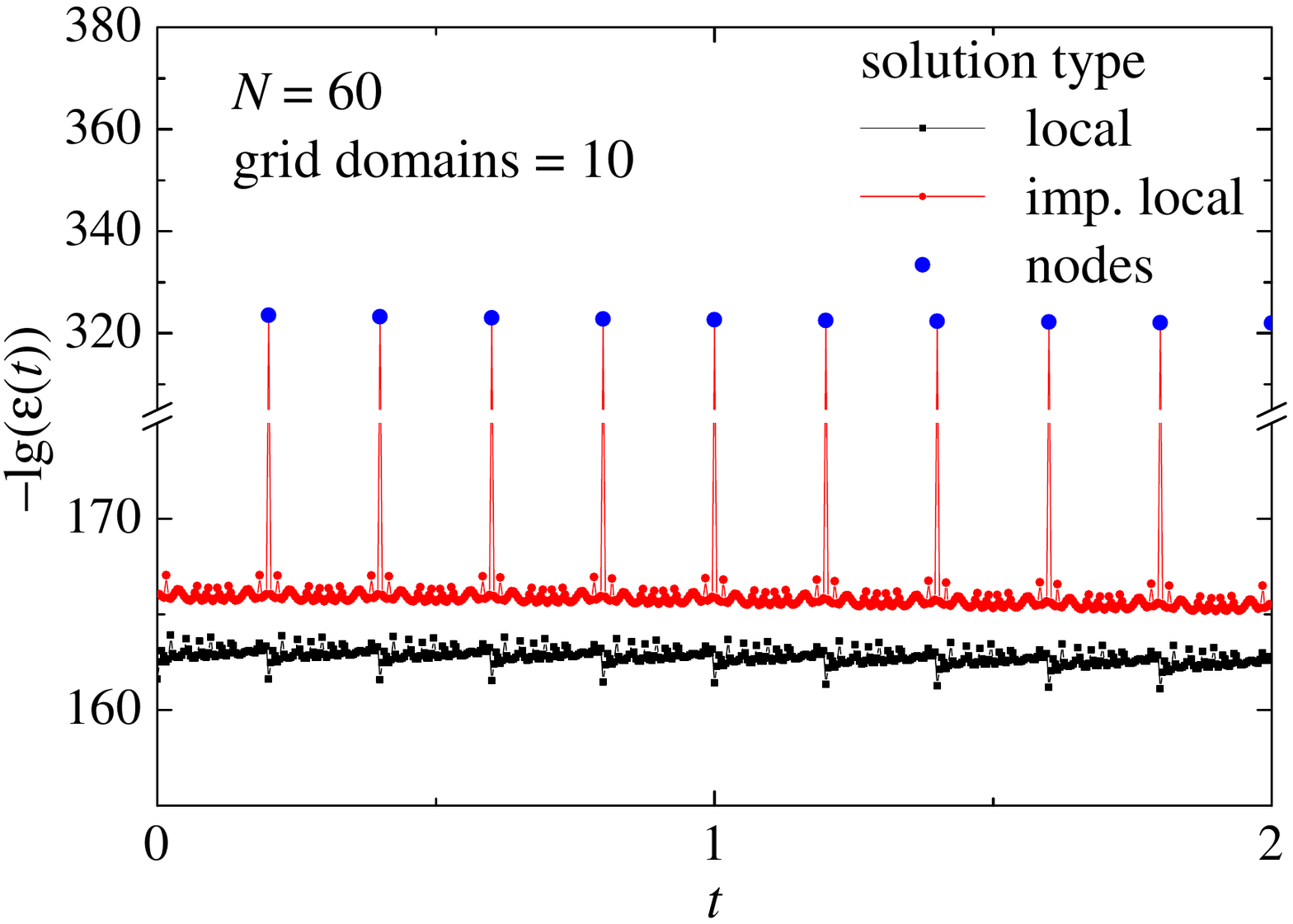}
\vspace{-8mm}\caption{\label{fig:lin_diss:c4}}
\end{subfigure}\\
\begin{subfigure}{0.24\textwidth}
\includegraphics[width=\textwidth]{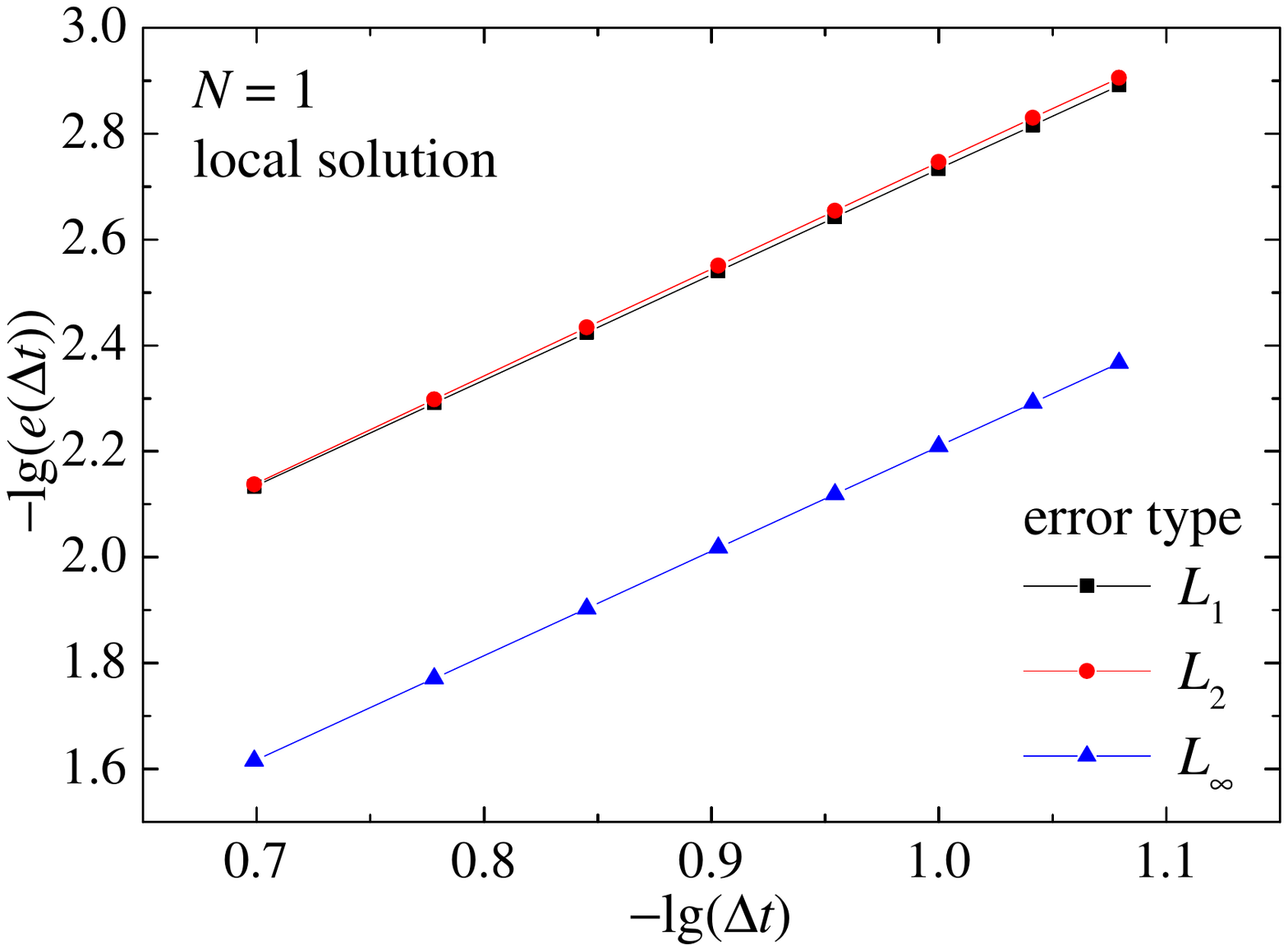}
\vspace{-8mm}\caption{\label{fig:lin_diss:d1}}
\end{subfigure}
\begin{subfigure}{0.24\textwidth}
\includegraphics[width=\textwidth]{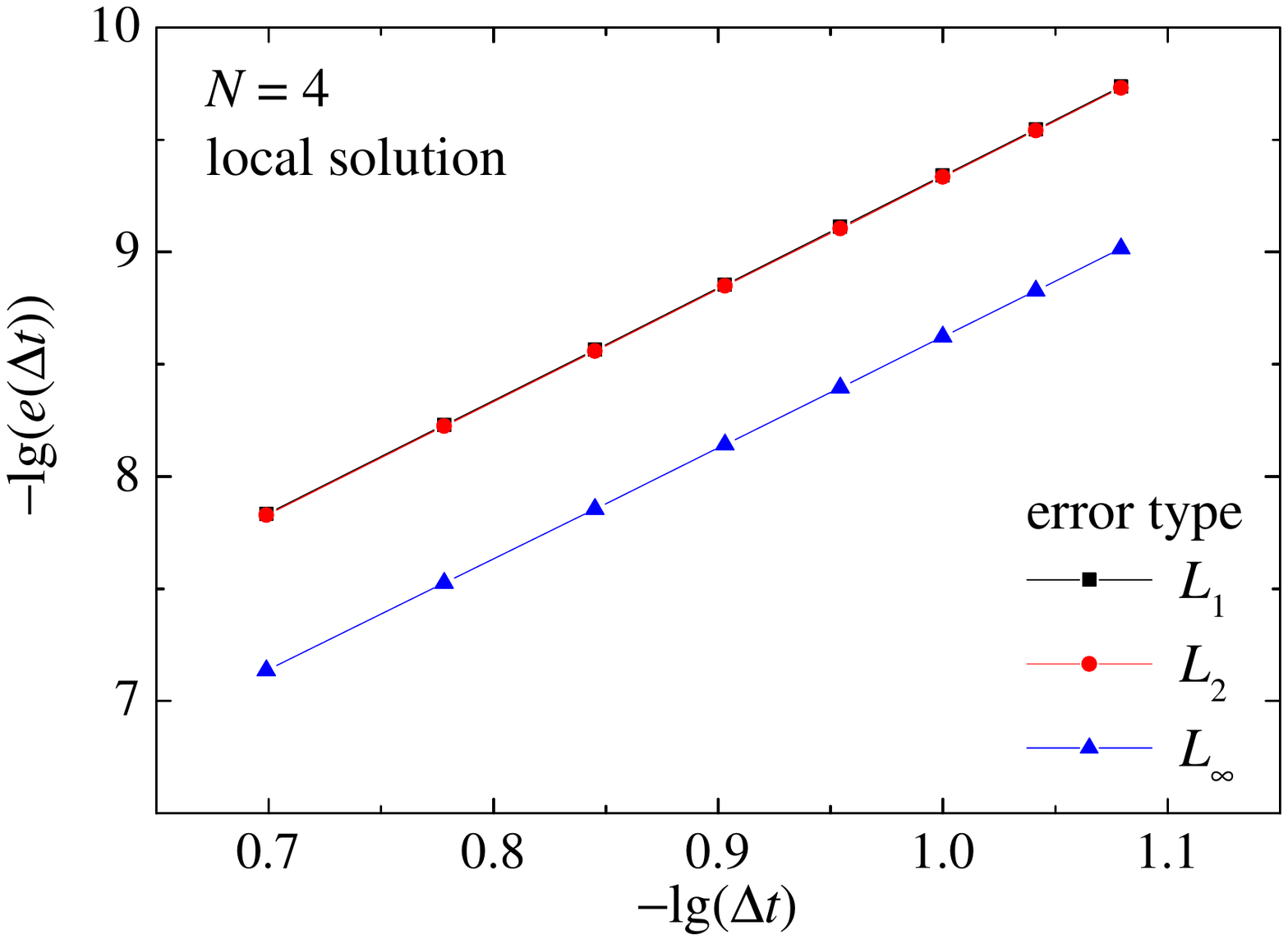}
\vspace{-8mm}\caption{\label{fig:lin_diss:d2}}
\end{subfigure}
\begin{subfigure}{0.24\textwidth}
\includegraphics[width=\textwidth]{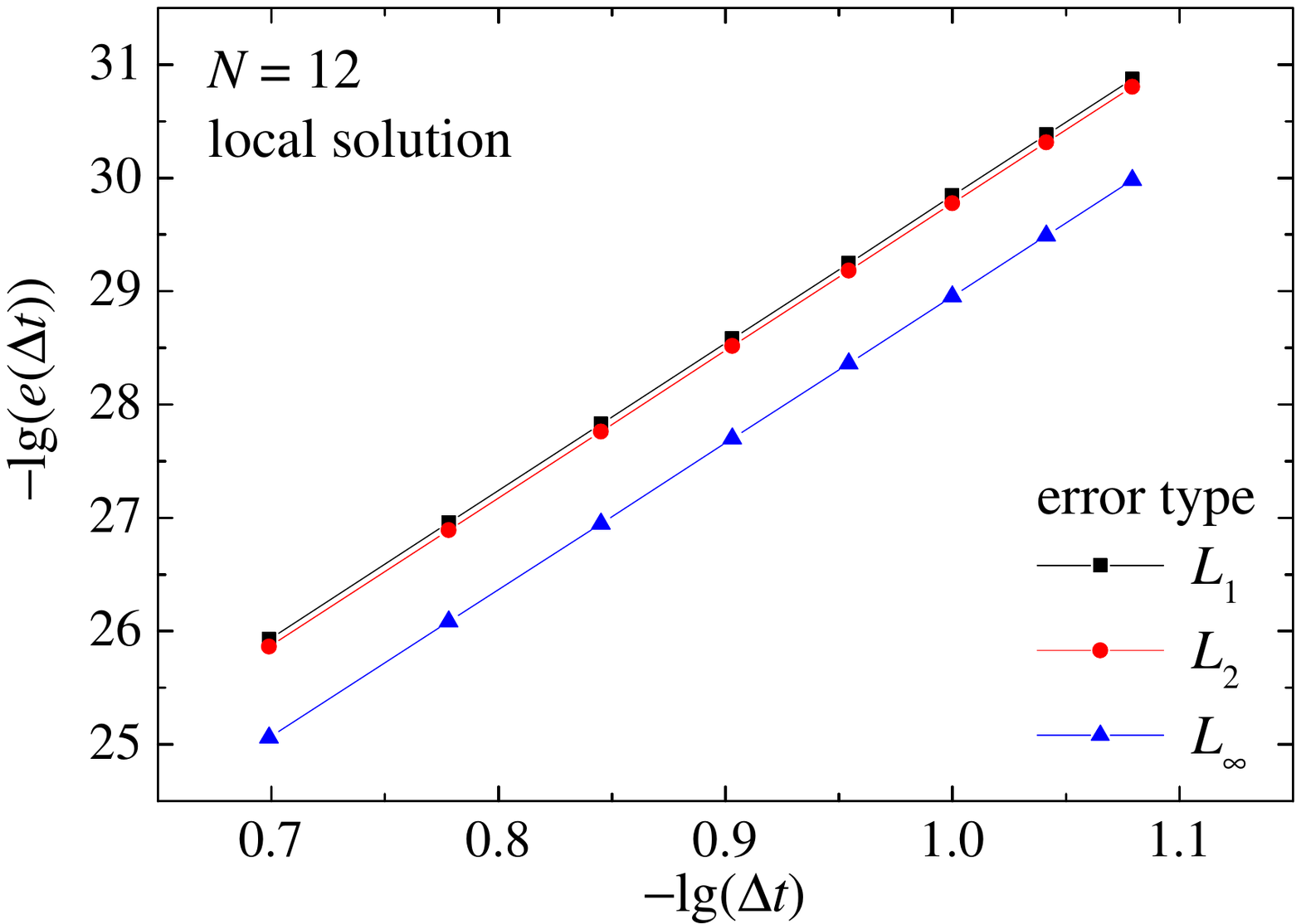}
\vspace{-8mm}\caption{\label{fig:lin_diss:d3}}
\end{subfigure}
\begin{subfigure}{0.24\textwidth}
\includegraphics[width=\textwidth]{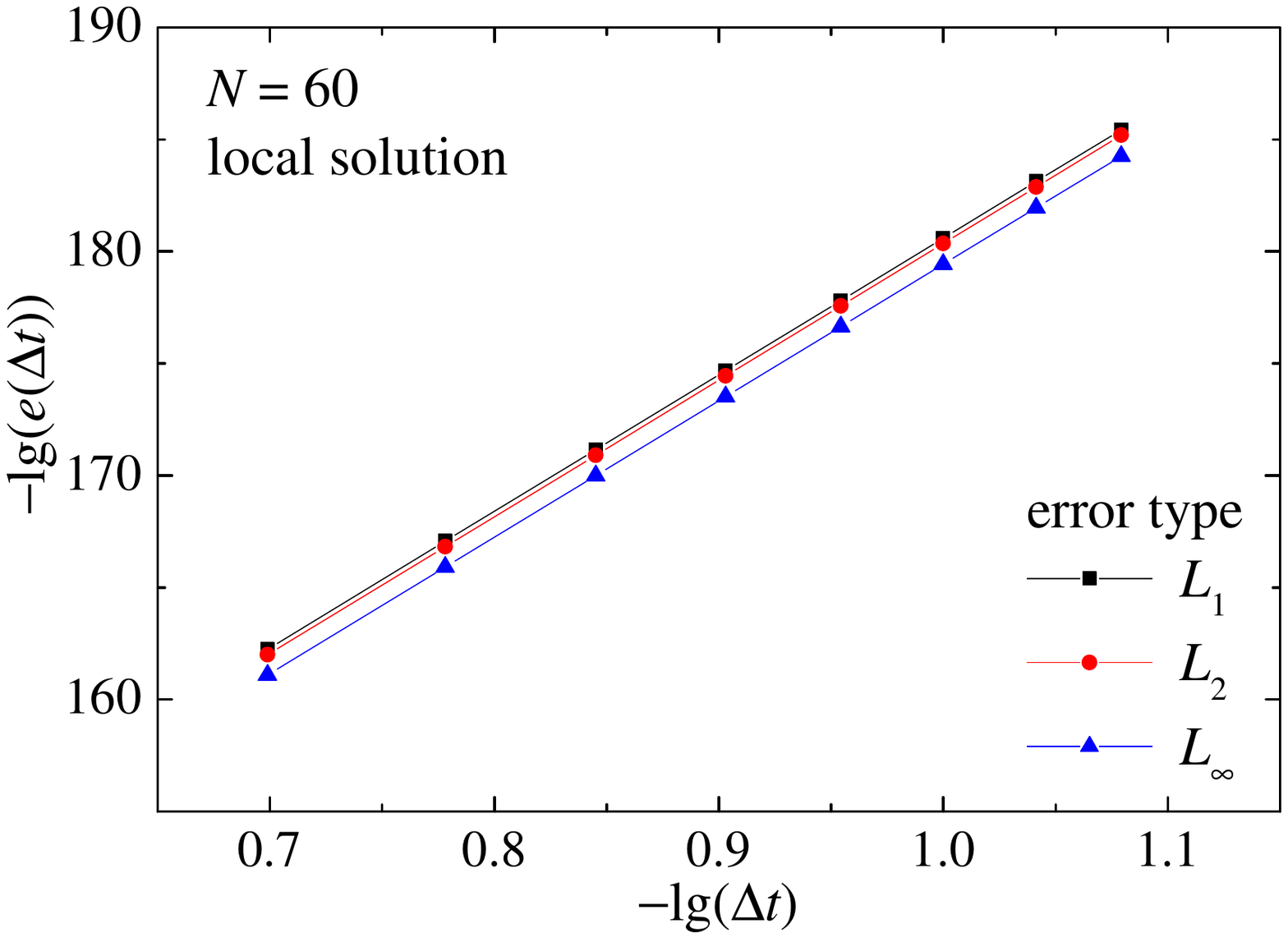}
\vspace{-8mm}\caption{\label{fig:lin_diss:d4}}
\end{subfigure}\\
\begin{subfigure}{0.24\textwidth}
\includegraphics[width=\textwidth]{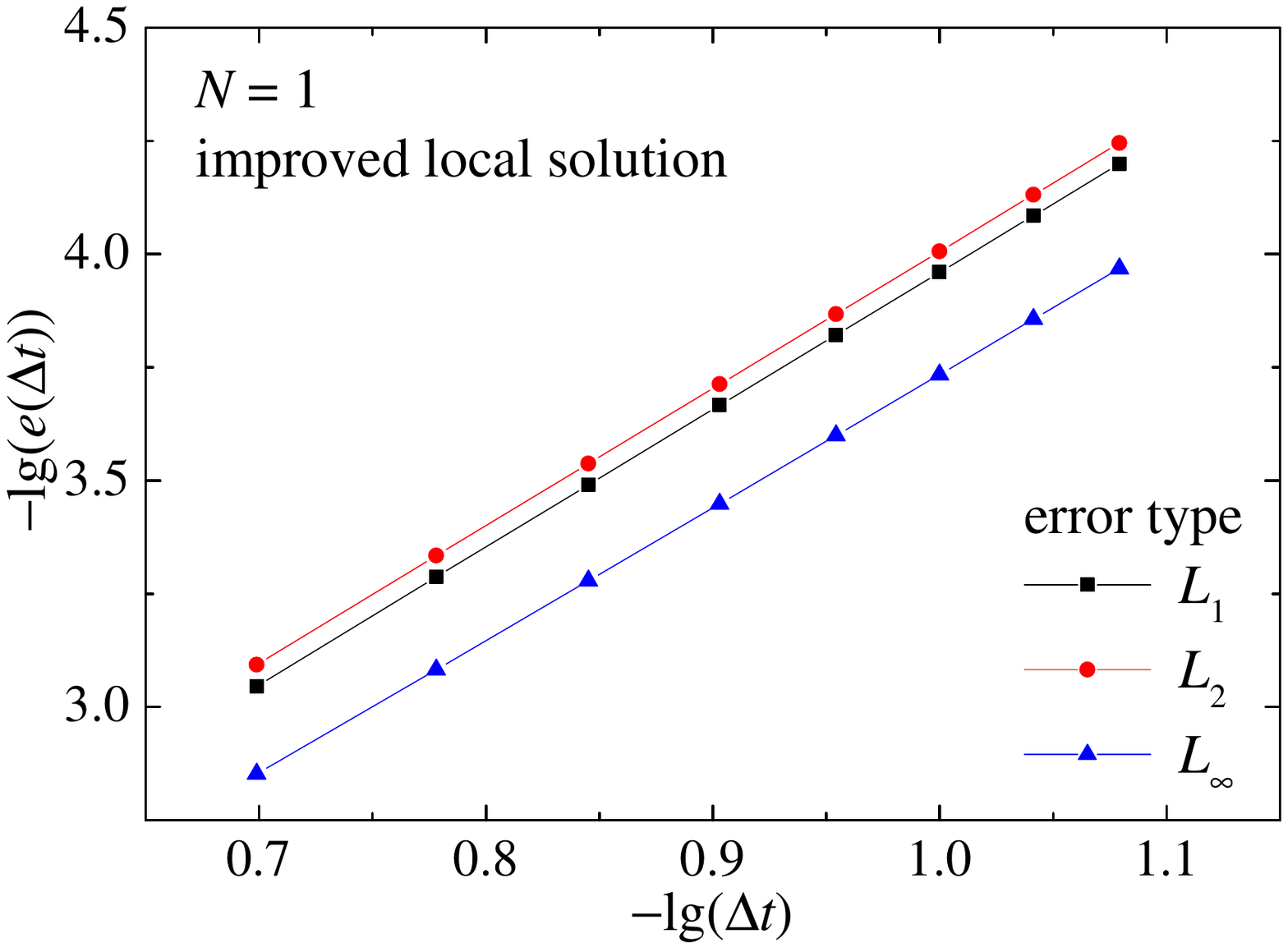}
\vspace{-8mm}\caption{\label{fig:lin_diss:e1}}
\end{subfigure}
\begin{subfigure}{0.24\textwidth}
\includegraphics[width=\textwidth]{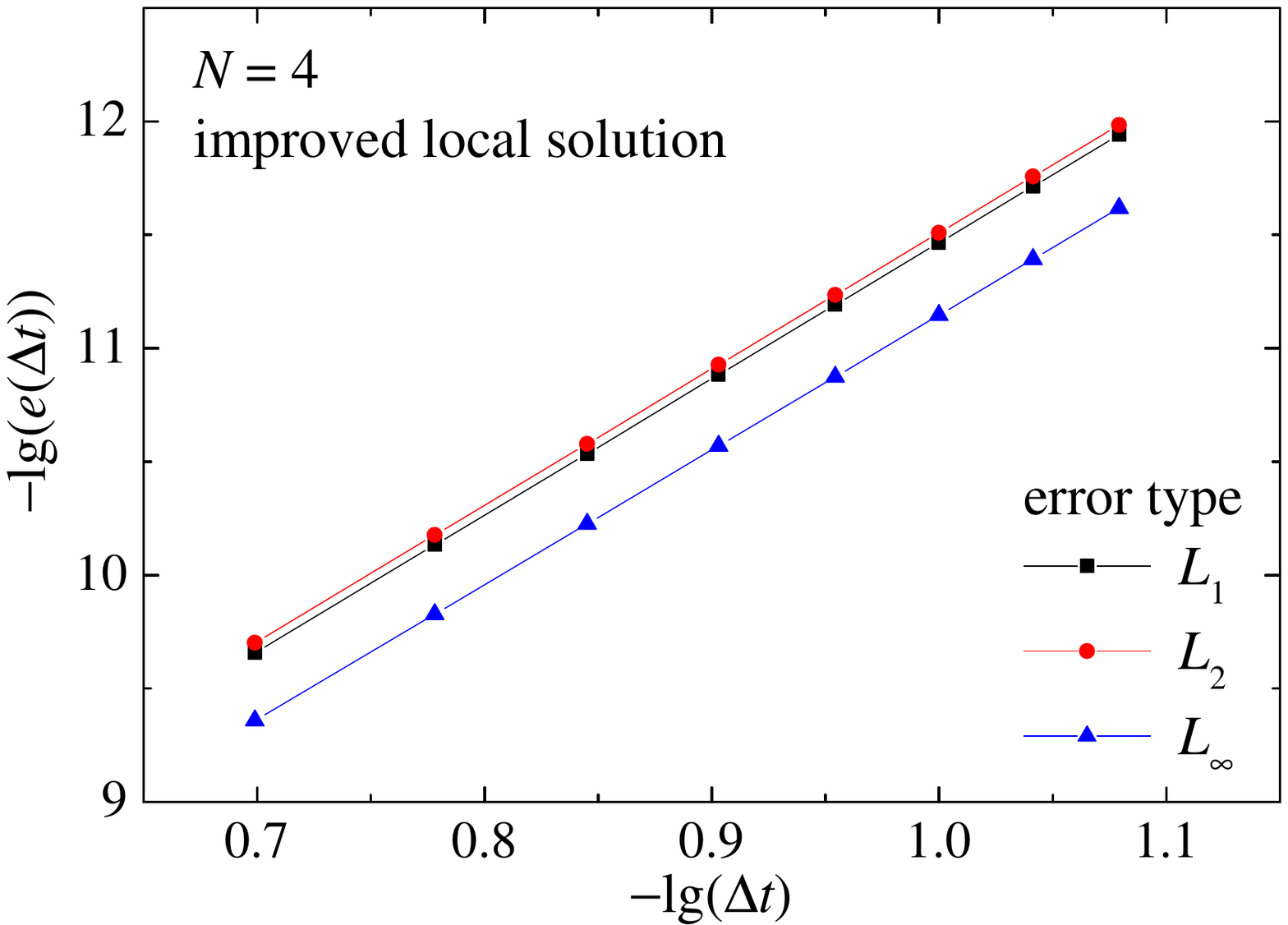}
\vspace{-8mm}\caption{\label{fig:lin_diss:e2}}
\end{subfigure}
\begin{subfigure}{0.24\textwidth}
\includegraphics[width=\textwidth]{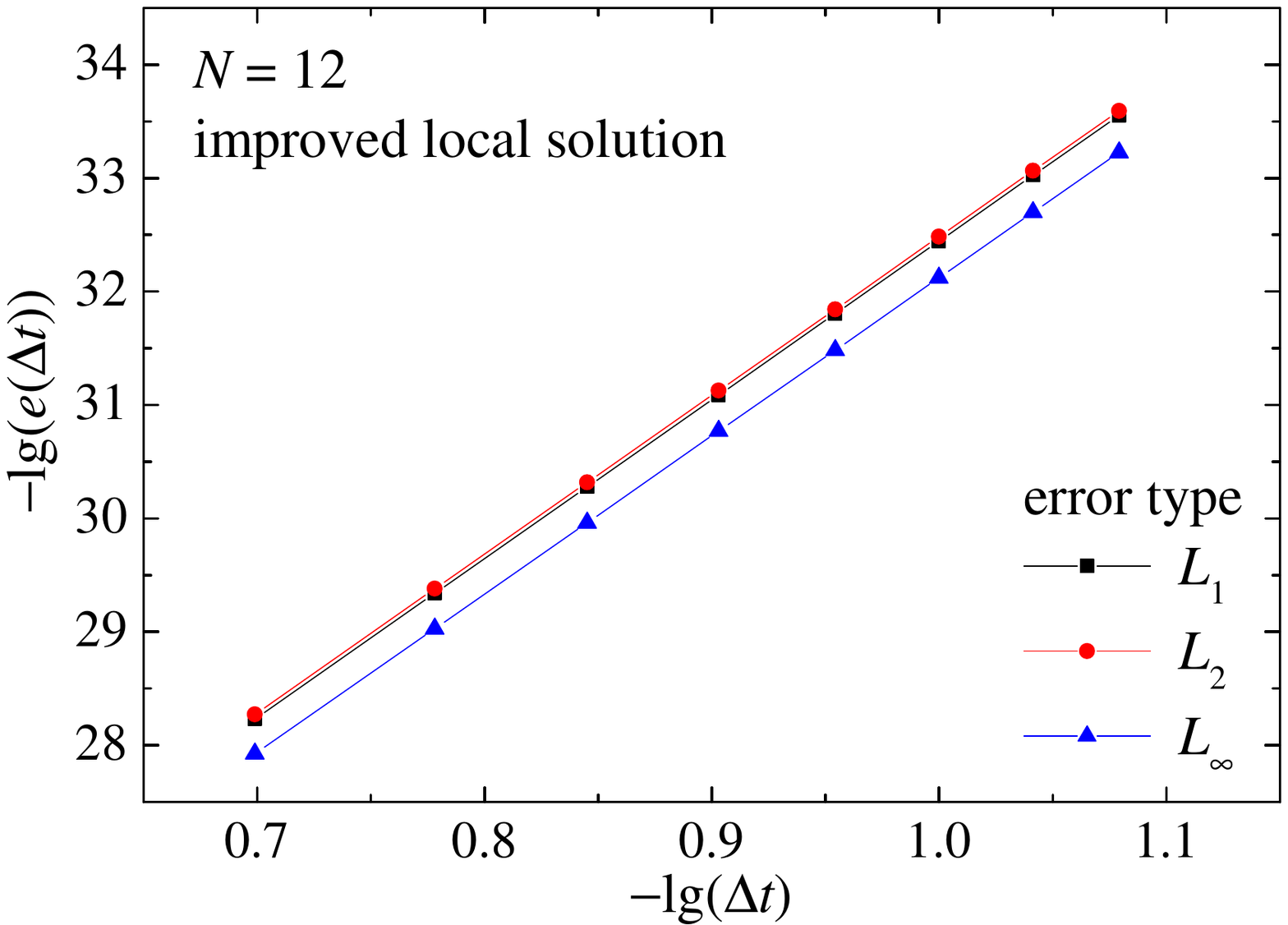}
\vspace{-8mm}\caption{\label{fig:lin_diss:e3}}
\end{subfigure}
\begin{subfigure}{0.24\textwidth}
\includegraphics[width=\textwidth]{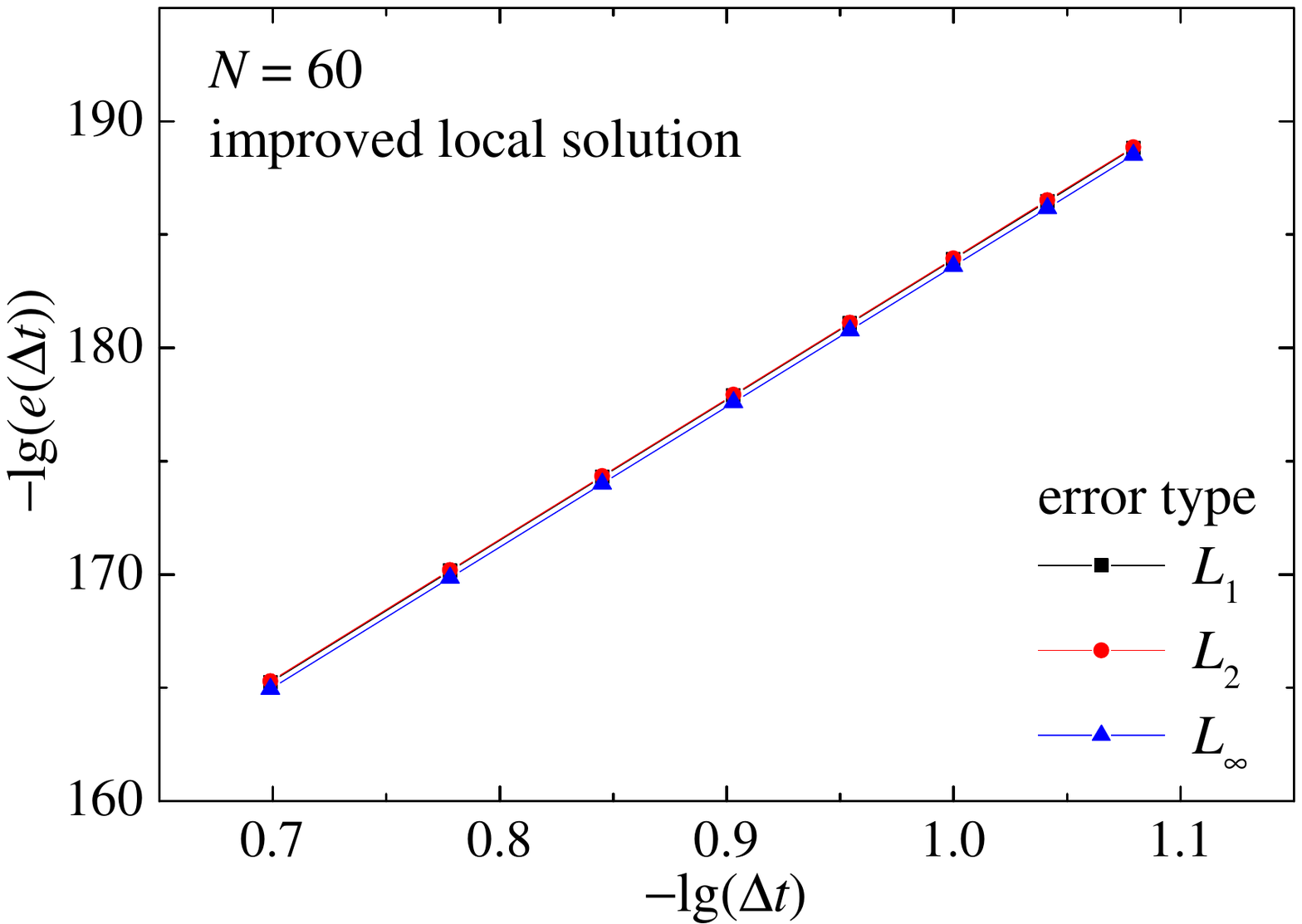}
\vspace{-8mm}\caption{\label{fig:lin_diss:e4}}
\end{subfigure}\\
\begin{subfigure}{0.24\textwidth}
\includegraphics[width=\textwidth]{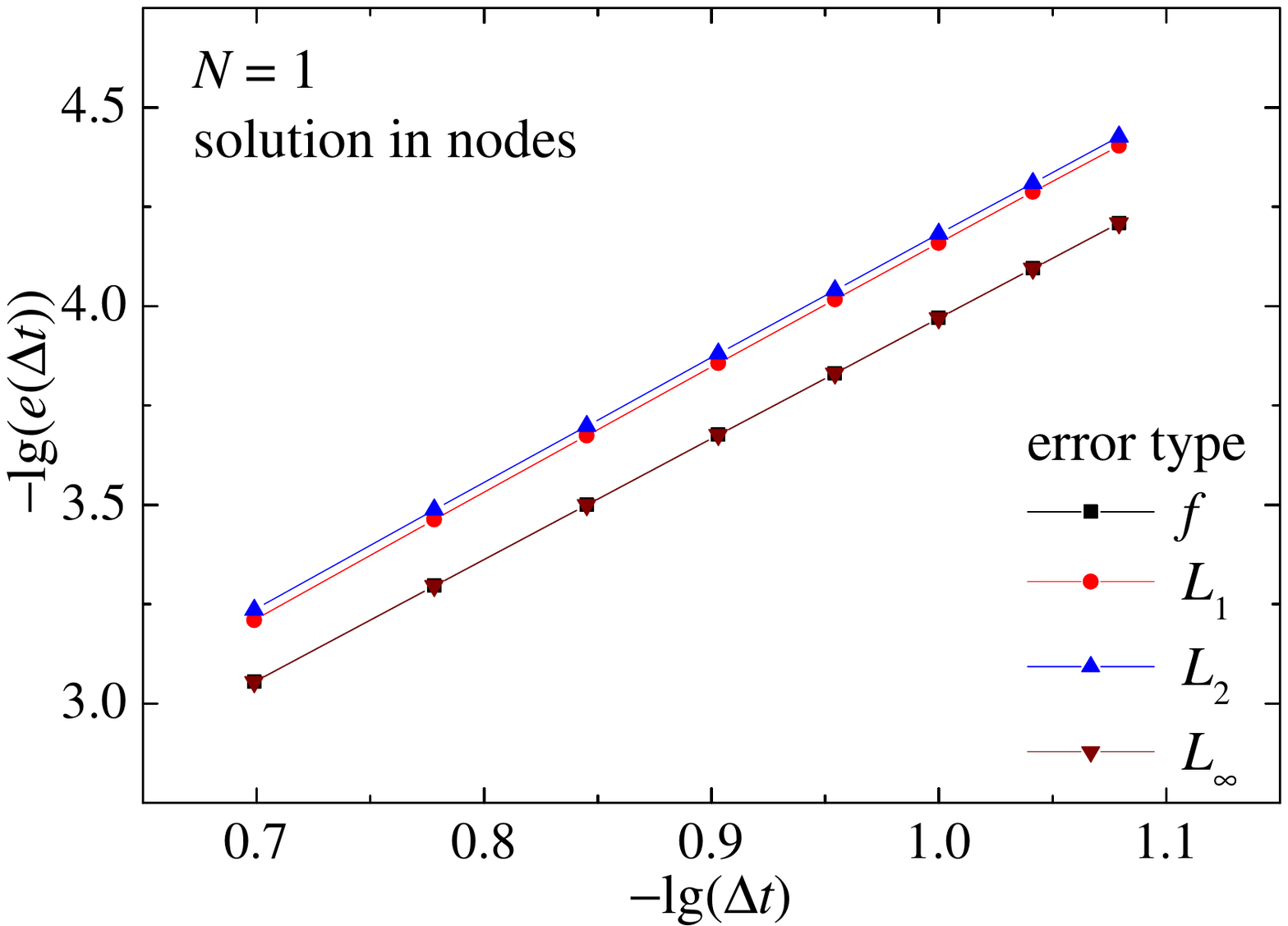}
\vspace{-8mm}\caption{\label{fig:lin_diss:f1}}
\end{subfigure}
\begin{subfigure}{0.24\textwidth}
\includegraphics[width=\textwidth]{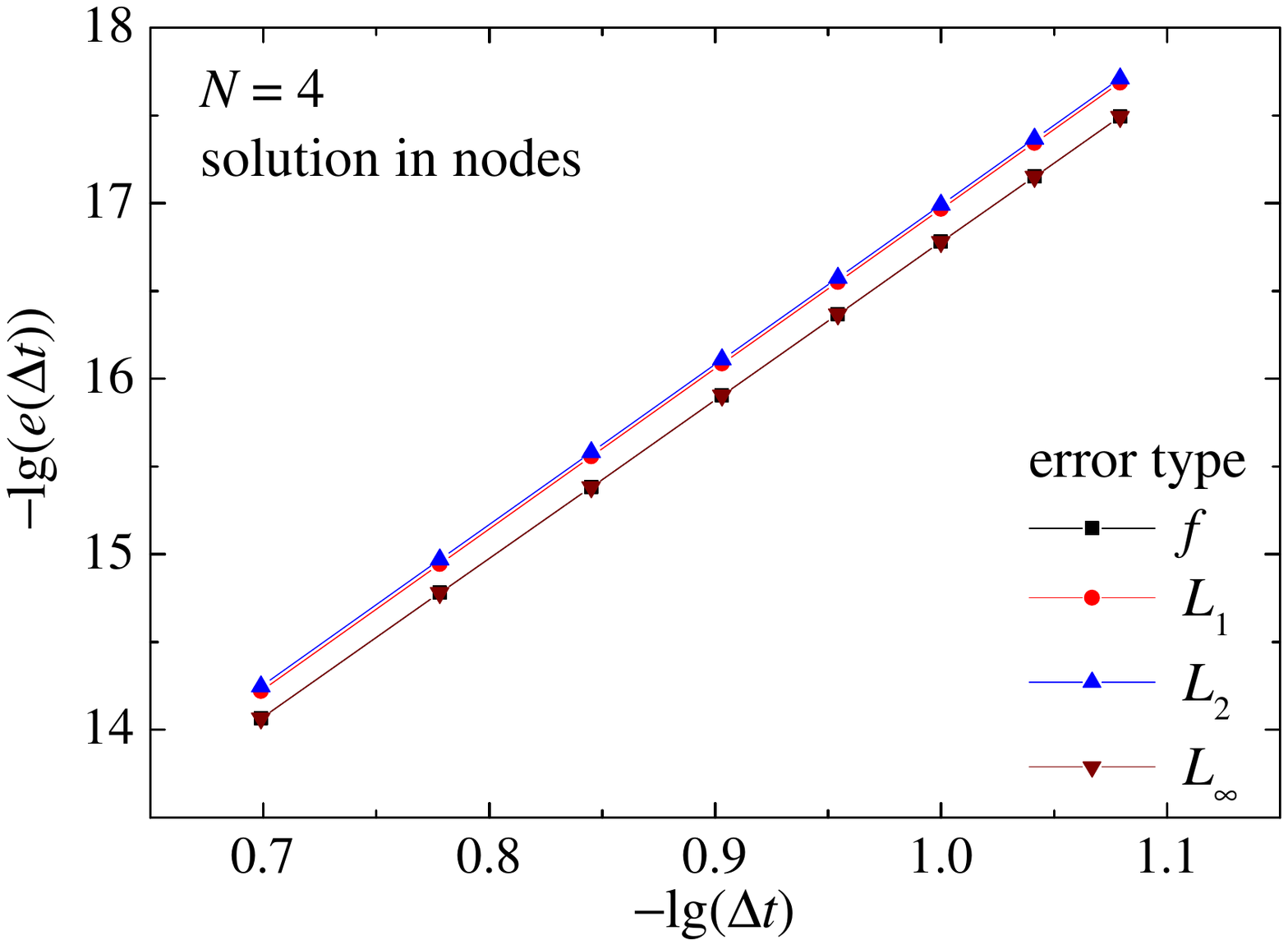}
\vspace{-8mm}\caption{\label{fig:lin_diss:f2}}
\end{subfigure}
\begin{subfigure}{0.24\textwidth}
\includegraphics[width=\textwidth]{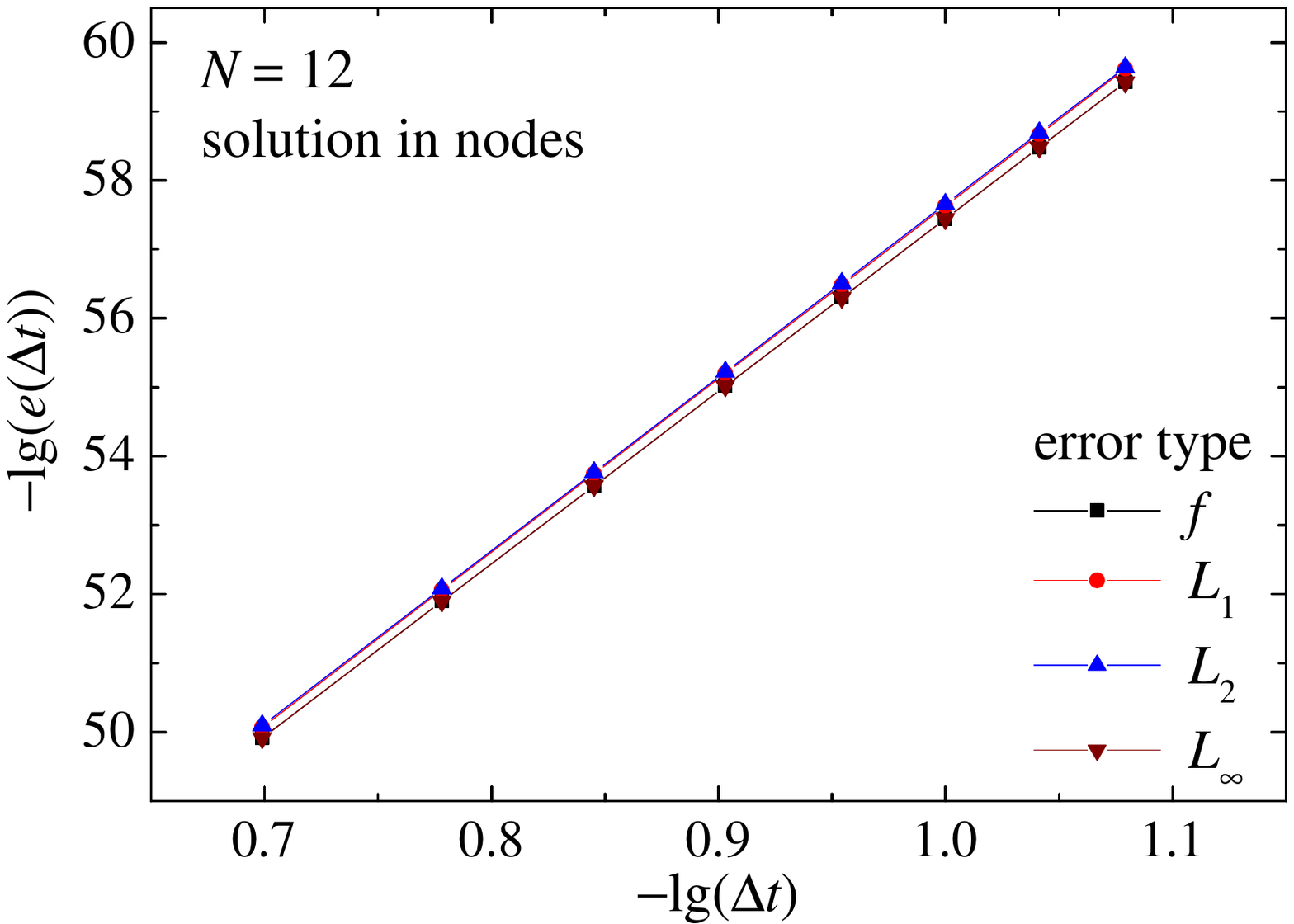}
\vspace{-8mm}\caption{\label{fig:lin_diss:f3}}
\end{subfigure}
\begin{subfigure}{0.24\textwidth}
\includegraphics[width=\textwidth]{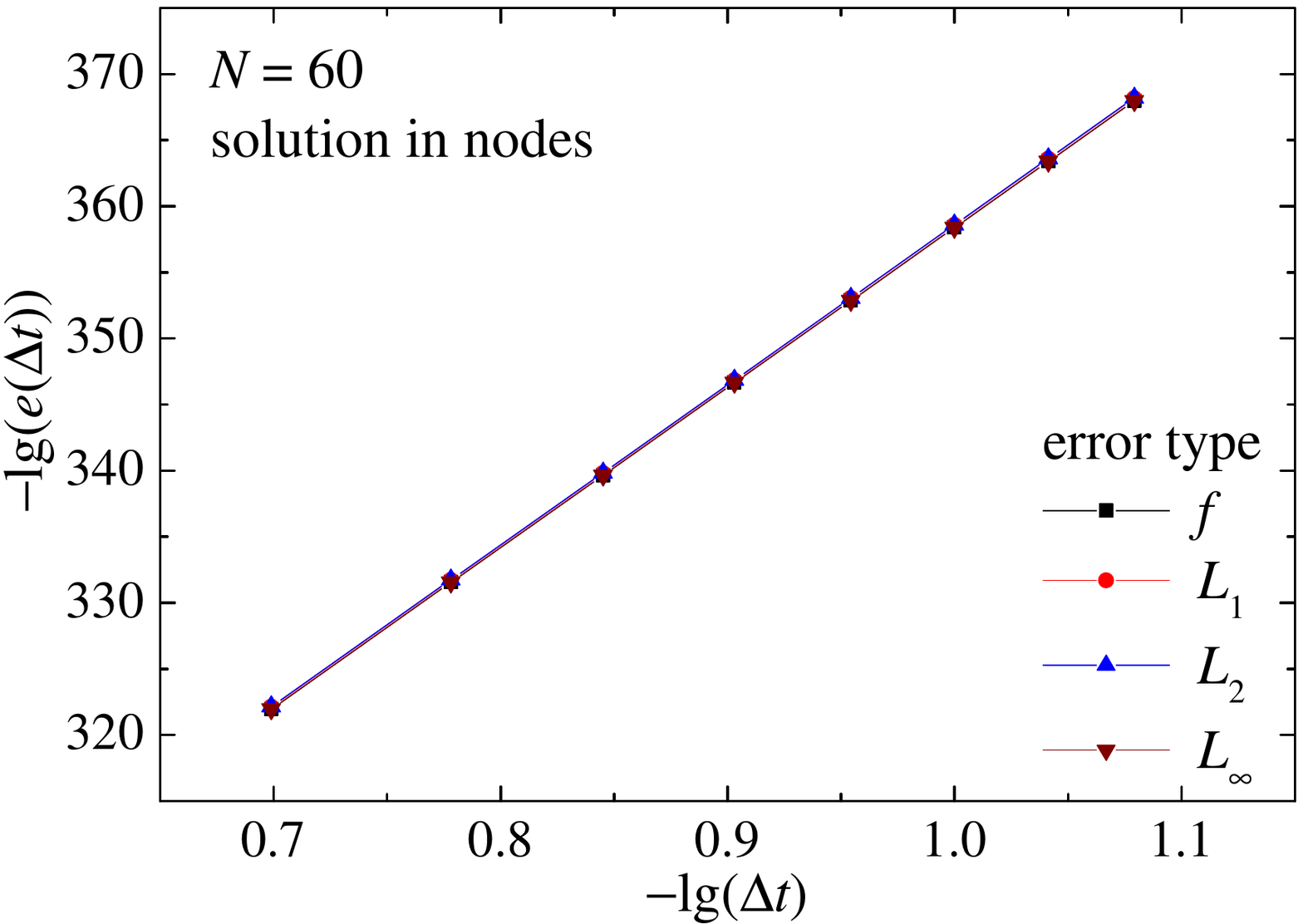}
\vspace{-8mm}\caption{\label{fig:lin_diss:f4}}
\end{subfigure}\\
\caption{%
Numerical solution of the system (\ref{eq:lin_diss_ode}). Comparison of the solution at nodes $\mathbf{u}_{n}$, the local solution $\mathbf{u}_{L}(t)$, the improved local solution $\mathbf{u}_{\rm IL}(t)$ and the exact solution $\mathbf{u}^{\rm ex}(t)$ (\ref{eq:lin_diss_sol_ex}) for components $u_{1} \equiv x$ (\subref{fig:lin_diss:a1}, \subref{fig:lin_diss:a2}, \subref{fig:lin_diss:a3}, \subref{fig:lin_diss:a4}) and $u_{2} \equiv \dot{x}$ (\subref{fig:lin_diss:b1}, \subref{fig:lin_diss:b2}, \subref{fig:lin_diss:b3}, \subref{fig:lin_diss:b4}), the errors $\varepsilon(t)$ (\subref{fig:lin_diss:c1}, \subref{fig:lin_diss:c2}, \subref{fig:lin_diss:c3}, \subref{fig:lin_diss:c4}), obtained using polynomials with degrees $N = 1$ (\subref{fig:lin_diss:a1}, \subref{fig:lin_diss:b1}, \subref{fig:lin_diss:c1}), $N = 4$ (\subref{fig:lin_diss:a2}, \subref{fig:lin_diss:b2}, \subref{fig:lin_diss:c2}), $N = 12$ (\subref{fig:lin_diss:a3}, \subref{fig:lin_diss:b3}, \subref{fig:lin_diss:c3}) and $N = 60$ (\subref{fig:lin_diss:a4}, \subref{fig:lin_diss:b4}, \subref{fig:lin_diss:c4}). Log-log plot of the dependence of the global error for the local solution $e^{l}$ (\subref{fig:lin_diss:d1}, \subref{fig:lin_diss:d2}, \subref{fig:lin_diss:d3}, \subref{fig:lin_diss:d4}), the improved local solution $e^{\rm imp}$ (\subref{fig:lin_diss:e1}, \subref{fig:lin_diss:e2}, \subref{fig:lin_diss:e3}, \subref{fig:lin_diss:e4}) and the solution at nodes $e^{n}$ (\subref{fig:lin_diss:f1}, \subref{fig:lin_diss:f2}, \subref{fig:lin_diss:f3}, \subref{fig:lin_diss:f4}) on the discretization step $\mathrm{\Delta}t$, obtained in the $f$-norm and norms $L_{1}$, $L_{2}$ and $L_{\infty}$, obtained using polynomials with degrees $N = 1$ (\subref{fig:lin_diss:d1}, \subref{fig:lin_diss:e1}, \subref{fig:lin_diss:f1}), $N = 4$ (\subref{fig:lin_diss:d2}, \subref{fig:lin_diss:e2}, \subref{fig:lin_diss:f2}), $N = 12$ (\subref{fig:lin_diss:d3}, \subref{fig:lin_diss:e3}, \subref{fig:lin_diss:f3}) and $N = 60$ (\subref{fig:lin_diss:d4}, \subref{fig:lin_diss:e4}, \subref{fig:lin_diss:f4}).
}
\label{fig:lin_diss}
\end{figure}

The second example of application of the ADER-DG numerical method with a local DG predictor to solving the initial value problem for the ODE system (\ref{eq:ivp_ode_diff_src}) is the following problem:
\begin{equation}\label{eq:lin_diss_ode}
\ddot{x} - x = 0,\quad
x(0) = 0,\quad \dot{x}(0) = 1,\quad
t \in [0,\, 2],
\end{equation}
for which the solution vectors in the original notation of the ODE system (\ref{eq:ivp_ode_diff_src}), with $D = 2$, are chosen in the form $\mathbf{u}(t) = [x(t)\ \dot{x}(t)]^{T}$. The exact analytical solution is as follows:
\begin{equation}\label{eq:lin_diss_sol_ex}
x^{\rm ex}(t) = \sinh(t),\quad
\dot{x}^{\rm ex}(t) = \cosh(t).
\end{equation}
The obtained results are presented in Fig.~\ref{fig:lin_diss} and Table~\ref{tab:conv_orders_lin_diss}.

\begin{table*}[h!]
\centering
\normalsize
\caption{%
Convergence orders $p_{f}$, $p_{L_{1}}$, $p_{L_{2}}$, $p_{L_{\infty}}$, calculated in $f$-norm and norms $L_{1}$, $L_{2}$, $L_{\infty}$, of the numerical solution of the ADER-DG method for the problem (\ref{eq:lin_diss_ode}); $N$ is the degree of the basis polynomials $\varphi_{p}$. Orders $p^{n}$ are calculated for \textit{the numerical solution at the nodes} $\mathbf{u}_{n}$; orders $p^{\rm imp}$ --- for \textit{the improved local solution} $\mathbf{u}_{\rm IL}$; orders $p^{l}$ --- for \textit{the local solution} $\mathbf{u}_{L}$. The theoretical values of convergence order $p_{\rm th.}^{n} = 2N+1$, $p_{\rm th.}^{l} = N+1$ and $p^{\rm imp}_{\rm th.} = N+2$ are presented for comparison.
}
\label{tab:conv_orders_lin_diss}
\setlength{\tabcolsep}{3.5pt}
\begin{tabular}{@{}|l|llll|c|lll|c|lll|c|@{}}
\toprule
$N$ & $p^{n}_{f}$ &
$p^{n}_{L_{1}}$ & $p^{n}_{L_{2}}$ & $p^{n}_{L_{\infty}}$ & $p^{n}_{\rm th.}$ &
$p^{l}_{L_{1}}$ & $p^{l}_{L_{2}}$ & $p^{l}_{L_{\infty}}$ & $p^{l}_{\rm th.}$ &
$p^{\rm imp}_{L_{1}}$ & $p^{\rm imp}_{L_{2}}$ & $p^{\rm imp}_{L_{\infty}}$ & $p^{\rm imp}_{\rm th.}$\\
\midrule
$1$ & $3.04$ & $3.14$ & $3.13$ & $3.04$ & $3$ & $2.00$ & $2.02$ & $1.98$ & $2$ & $3.03$ & $3.03$ & $2.94$ & $3$\\
$2$ & $5.02$ & $5.13$ & $5.12$ & $5.02$ & $5$ & $3.01$ & $3.01$ & $2.95$ & $3$ & $4.00$ & $4.01$ & $3.99$ & $4$\\
$3$ & $7.02$ & $7.12$ & $7.11$ & $7.02$ & $7$ & $4.01$ & $4.01$ & $3.95$ & $4$ & $5.01$ & $5.01$ & $4.94$ & $5$\\
$4$ & $9.01$ & $9.12$ & $9.11$ & $9.01$ & $9$ & $5.00$ & $5.00$ & $4.94$ & $5$ & $6.00$ & $6.00$ & $5.95$ & $6$\\
$5$ & $11.0$ & $11.1$ & $11.1$ & $11.0$ & $11$ & $6.00$ & $6.00$ & $5.94$ & $6$ & $7.00$ & $7.00$ & $6.94$ & $7$\\
$6$ & $13.0$ & $13.1$ & $13.1$ & $13.0$ & $13$ & $7.00$ & $7.00$ & $6.94$ & $7$ & $8.00$ & $8.00$ & $7.94$ & $8$\\
$7$ & $15.0$ & $15.1$ & $15.1$ & $15.0$ & $15$ & $8.00$ & $8.00$ & $7.94$ & $8$ & $9.00$ & $9.00$ & $8.94$ & $9$\\
$8$ & $17.0$ & $17.1$ & $17.1$ & $17.0$ & $17$ & $9.00$ & $9.00$ & $8.94$ & $9$ & $10.0$ & $10.0$ & $9.94$ & $10$\\
$9$ & $19.0$ & $19.1$ & $19.1$ & $19.0$ & $19$ & $10.0$ & $10.0$ & $9.94$ & $10$ & $11.0$ & $11.0$ & $10.9$ & $11$\\
$10$ & $21.0$ & $21.1$ & $21.1$ & $21.0$ & $21$ & $11.0$ & $11.0$ & $10.9$ & $11$ & $12.0$ & $12.0$ & $11.9$ & $12$\\
\midrule
$11$ & $23.0$ & $23.1$ & $23.1$ & $23.0$ & $23$ & $12.0$ & $12.0$ & $11.9$ & $12$ & $13.0$ & $13.0$ & $12.9$ & $13$\\
$12$ & $25.0$ & $25.1$ & $25.1$ & $25.0$ & $25$ & $13.0$ & $13.0$ & $12.9$ & $13$ & $14.0$ & $14.0$ & $13.9$ & $14$\\
$13$ & $27.0$ & $27.1$ & $27.1$ & $27.0$ & $27$ & $14.0$ & $14.0$ & $13.9$ & $14$ & $15.0$ & $15.0$ & $14.9$ & $15$\\
$14$ & $29.0$ & $29.1$ & $29.1$ & $29.0$ & $29$ & $15.0$ & $15.0$ & $14.9$ & $15$ & $16.0$ & $16.0$ & $15.9$ & $16$\\
$15$ & $31.0$ & $31.1$ & $31.1$ & $31.0$ & $31$ & $16.0$ & $16.0$ & $15.9$ & $16$ & $17.0$ & $17.0$ & $16.9$ & $17$\\
$16$ & $33.0$ & $33.1$ & $33.1$ & $33.0$ & $33$ & $17.0$ & $17.0$ & $16.9$ & $17$ & $18.0$ & $18.0$ & $17.9$ & $18$\\
$17$ & $35.0$ & $35.1$ & $35.1$ & $35.0$ & $35$ & $18.0$ & $18.0$ & $17.9$ & $18$ & $19.0$ & $19.0$ & $18.9$ & $19$\\
$18$ & $37.0$ & $37.1$ & $37.1$ & $37.0$ & $37$ & $19.0$ & $19.0$ & $18.9$ & $19$ & $20.0$ & $20.0$ & $19.9$ & $20$\\
$19$ & $39.0$ & $39.1$ & $39.1$ & $39.0$ & $39$ & $20.0$ & $20.0$ & $19.9$ & $20$ & $21.0$ & $21.0$ & $20.9$ & $21$\\
$20$ & $41.0$ & $41.1$ & $41.1$ & $41.0$ & $41$ & $21.0$ & $21.0$ & $20.9$ & $21$ & $22.0$ & $22.0$ & $21.9$ & $22$\\
\midrule
$21$ & $43.0$ & $43.1$ & $43.1$ & $43.0$ & $43$ & $22.0$ & $22.0$ & $21.9$ & $22$ & $23.0$ & $23.0$ & $22.9$ & $23$\\
$22$ & $45.0$ & $45.1$ & $45.1$ & $45.0$ & $45$ & $23.0$ & $23.0$ & $22.9$ & $23$ & $24.0$ & $24.0$ & $23.9$ & $24$\\
$23$ & $47.0$ & $47.1$ & $47.1$ & $47.0$ & $47$ & $24.0$ & $24.0$ & $23.9$ & $24$ & $25.0$ & $25.0$ & $24.9$ & $25$\\
$24$ & $49.0$ & $49.1$ & $49.1$ & $49.0$ & $49$ & $25.0$ & $25.0$ & $24.9$ & $25$ & $26.0$ & $26.0$ & $25.9$ & $26$\\
$25$ & $51.0$ & $51.1$ & $51.1$ & $51.0$ & $51$ & $26.0$ & $26.0$ & $25.9$ & $26$ & $27.0$ & $27.0$ & $26.9$ & $27$\\
$26$ & $53.0$ & $53.1$ & $53.1$ & $53.0$ & $53$ & $27.0$ & $27.0$ & $26.9$ & $27$ & $28.0$ & $28.0$ & $27.9$ & $28$\\
$27$ & $55.0$ & $55.1$ & $55.1$ & $55.0$ & $55$ & $28.0$ & $28.0$ & $27.9$ & $28$ & $29.0$ & $29.0$ & $28.9$ & $29$\\
$28$ & $57.0$ & $57.1$ & $57.1$ & $57.0$ & $57$ & $29.0$ & $29.0$ & $28.9$ & $29$ & $30.0$ & $30.0$ & $29.9$ & $30$\\
$29$ & $59.0$ & $59.1$ & $59.1$ & $59.0$ & $59$ & $30.0$ & $30.0$ & $29.9$ & $30$ & $31.0$ & $31.0$ & $30.9$ & $31$\\
$30$ & $61.0$ & $61.1$ & $61.1$ & $61.0$ & $61$ & $31.0$ & $31.0$ & $30.9$ & $31$ & $32.0$ & $32.0$ & $31.9$ & $32$\\
\midrule
$35$ & $71.0$ & $71.1$ & $71.1$ & $71.0$ & $71$ & $36.0$ & $36.0$ & $35.9$ & $36$ & $37.0$ & $37.0$ & $36.9$ & $37$\\
$40$ & $81.0$ & $81.1$ & $81.1$ & $81.0$ & $81$ & $41.0$ & $41.0$ & $40.9$ & $41$ & $42.0$ & $42.0$ & $41.9$ & $42$\\
$45$ & $91.0$ & $91.1$ & $91.1$ & $91.0$ & $91$ & $46.0$ & $46.0$ & $45.9$ & $46$ & $47.0$ & $47.0$ & $46.9$ & $47$\\
$50$ & $101.0$ & $101.1$ & $101.1$ & $101.0$ & $101$ & $51.0$ & $51.0$ & $50.9$ & $51$ & $52.0$ & $52.0$ & $51.9$ & $52$\\
$55$ & $111.0$ & $111.1$ & $111.1$ & $111.0$ & $111$ & $56.0$ & $56.0$ & $55.9$ & $56$ & $57.0$ & $57.0$ & $56.9$ & $57$\\
$60$ & $121.0$ & $121.1$ & $121.1$ & $121.0$ & $121$ & $61.0$ & $61.0$ & $60.9$ & $61$ & $62.0$ & $62.0$ & $61.9$ & $62$\\
\bottomrule
\end{tabular}
\end{table*}

Fig.~\ref{fig:lin_diss} shows the dependencies of the numerical solutions $\mathbf{u}_{L}$, $\mathbf{u}_{\rm IL}$, $\mathbf{u}_{n}$ and the exact analytical solution $\mathbf{u}^{\rm ex}$, the dependencies of the local error $\varepsilon$ (\ref{eq:eps_local_def}) of the numerical solutions, and the dependence of the global error $e$ (\ref{eq:eps_un_global_def}), (\ref{eq:eps_ul_global_def}) of the numerical solutions on the discretization step ${\Delta t}$, for polynomial degrees $N = 1$, $4$, $12$ and $60$. A comparison of the obtained dependencies of the numerical solutions $\mathbf{u}_{L}(t)$, $\mathbf{u}_{\rm IL}(t)$, $\mathbf{u}_{n}$ with the exact analytical solution $\mathbf{u}^{\rm ex}(t)$, presented in Fig.~\ref{fig:lin_diss} (\subref{fig:lin_diss:a1}, \subref{fig:lin_diss:a2}, \subref{fig:lin_diss:a3}, \subref{fig:lin_diss:a4}) for the component $u_{1}$ and in Fig.~\ref{fig:lin_diss} (\subref{fig:lin_diss:b1}, \subref{fig:lin_diss:b2}, \subref{fig:lin_diss:b3}, \subref{fig:lin_diss:b4}) for the component $u_{2}$, demonstrates a high-quality agreement.

\begin{figure}[h!]
\captionsetup[subfigure]{%
	position=bottom,
	font+=smaller,
	textfont=normalfont,
	singlelinecheck=off,
	justification=raggedright
}
\centering
\begin{subfigure}{0.24\textwidth}
\includegraphics[width=\textwidth]{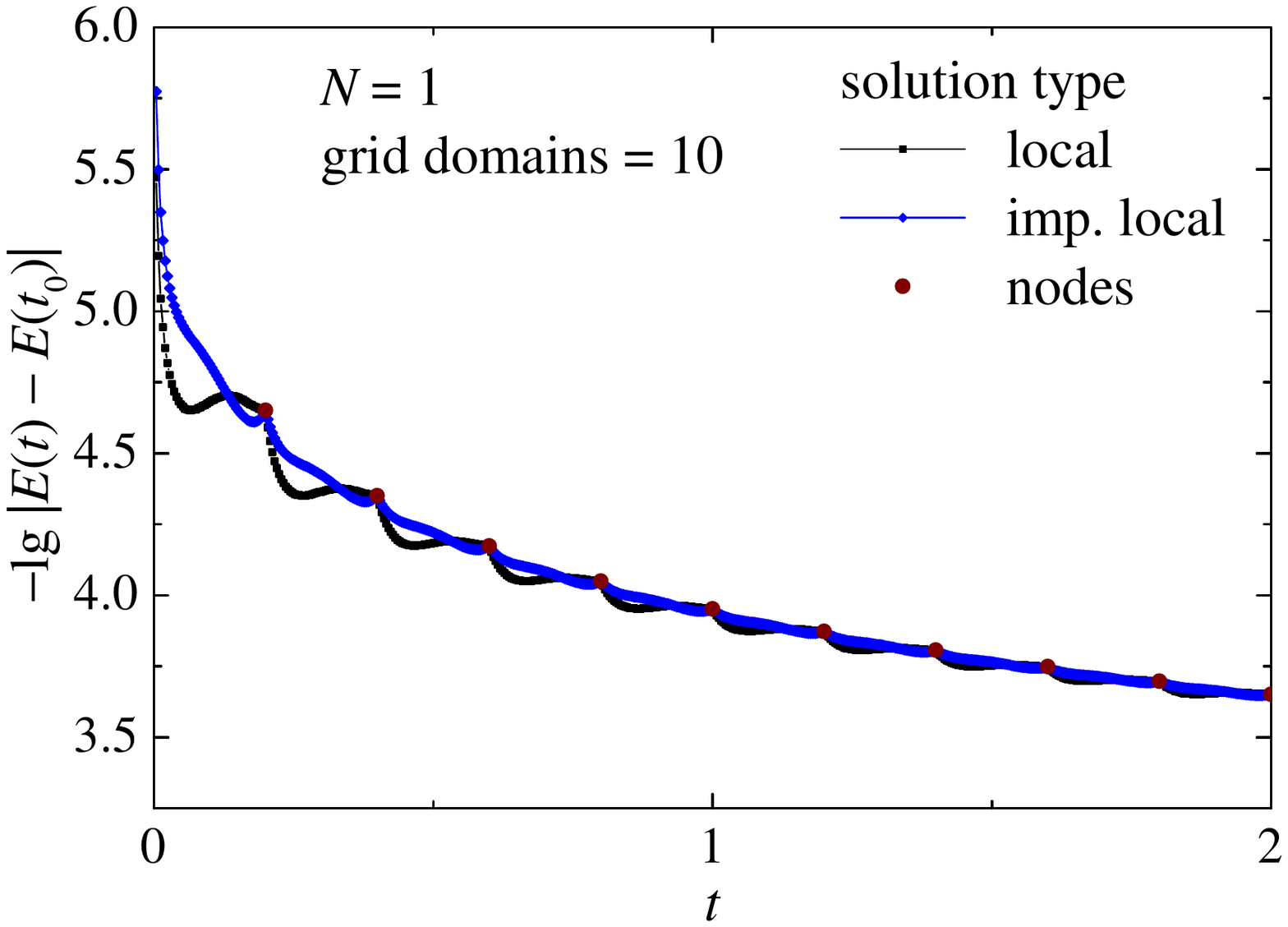}
\vspace{-8mm}\caption{\label{fig:econs_lin_diss:a1}}
\end{subfigure}
\begin{subfigure}{0.24\textwidth}
\includegraphics[width=\textwidth]{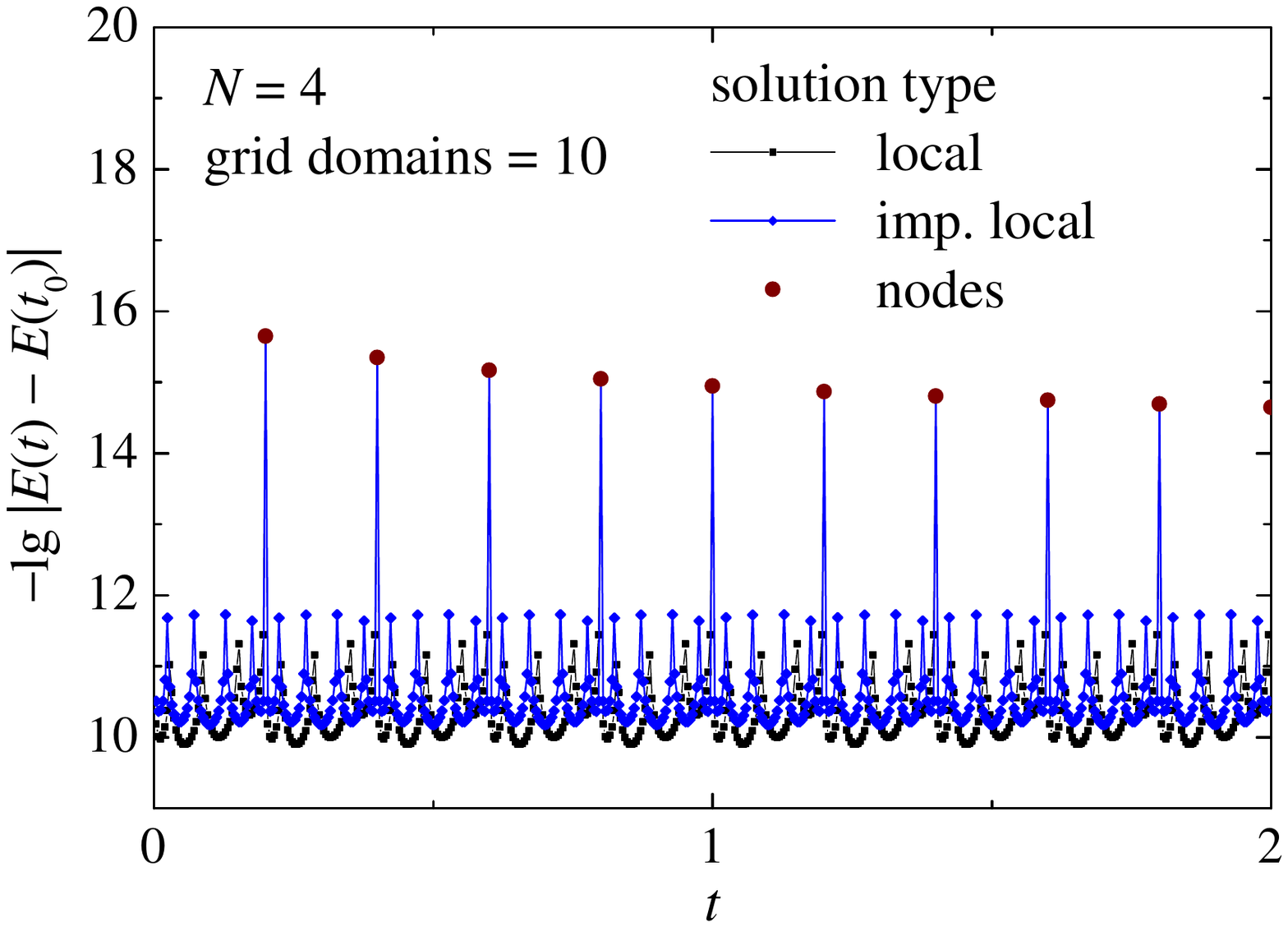}
\vspace{-8mm}\caption{\label{fig:econs_lin_diss:a2}}
\end{subfigure}
\begin{subfigure}{0.24\textwidth}
\includegraphics[width=\textwidth]{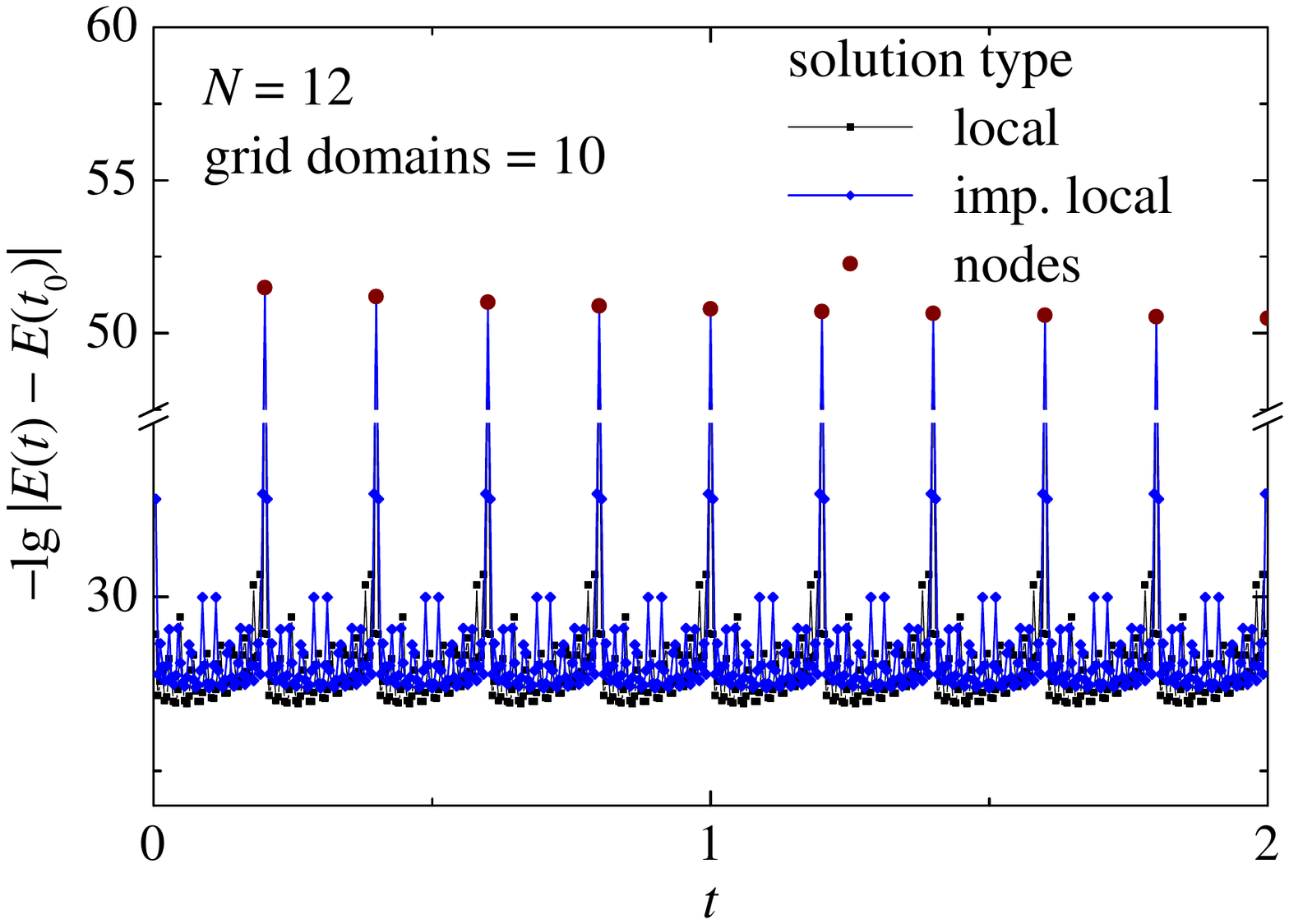}
\vspace{-8mm}\caption{\label{fig:econs_lin_diss:a3}}
\end{subfigure}
\begin{subfigure}{0.24\textwidth}
\includegraphics[width=\textwidth]{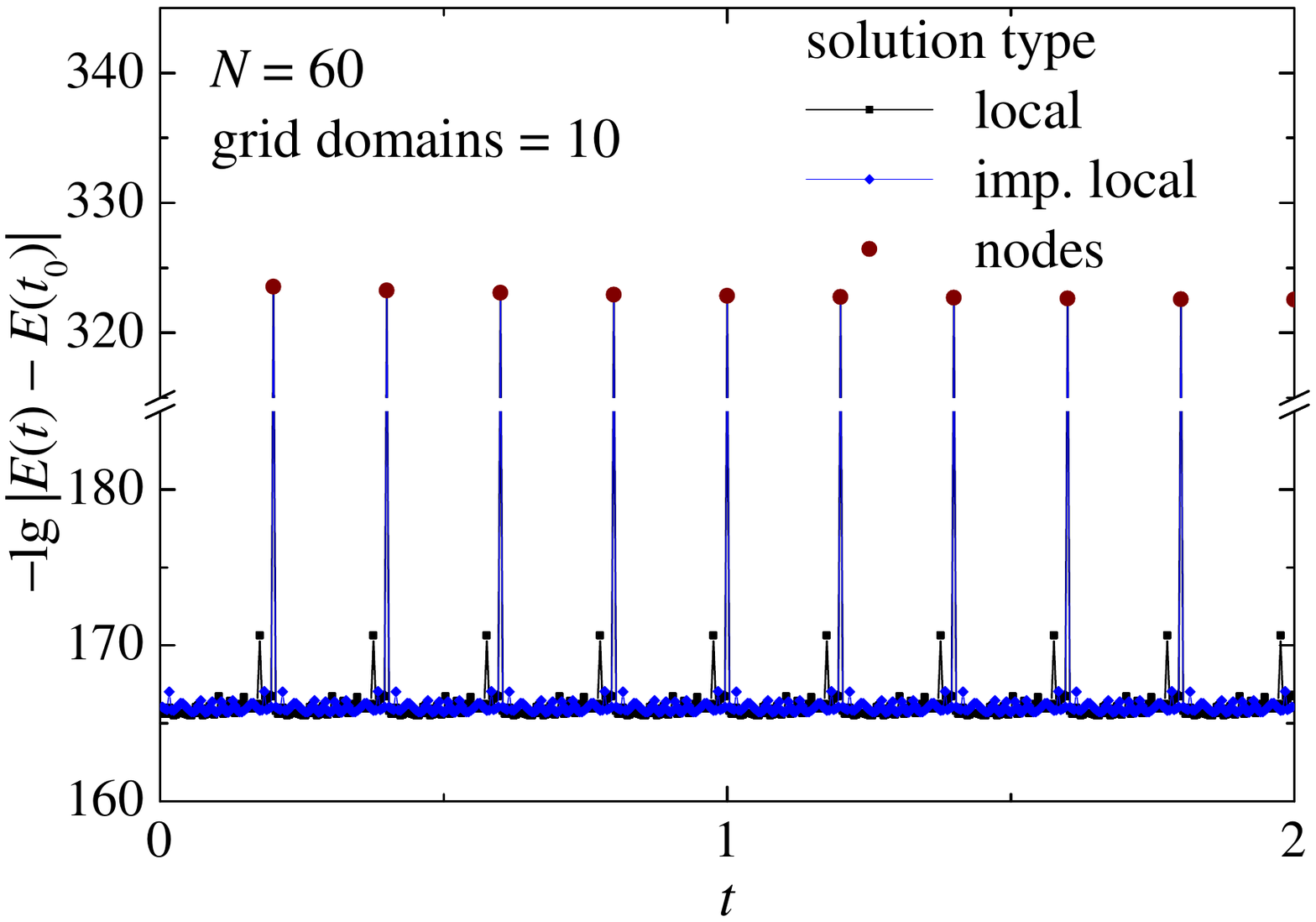}
\vspace{-8mm}\caption{\label{fig:econs_lin_diss:a4}}
\end{subfigure}
\caption{%
Dependence of the negative logarithm of the error $-\lg|E(t)-E(t_{0})|$ of the energy conservation law (\ref{eq:econs_lin_diss}) of the numerical solution of the system (\ref{eq:lin_diss_ode}) on the argument $t$ for solution at nodes $\mathbf{u}_{n}$, the local solution $\mathbf{u}_{L}(t)$ and the improved local solution $\mathbf{u}_{\rm IL}(t)$, obtained using polynomials with degrees $N = 1$~(\subref{fig:econs_lin_diss:a1}), $N = 4$~(\subref{fig:econs_lin_diss:a4}), $N = 12$~(\subref{fig:econs_lin_diss:a3}), $N = 60$~(\subref{fig:econs_lin_diss:a4}).
}
\label{fig:econs_lin_diss}
\end{figure}

The local error $\varepsilon$ (\ref{eq:eps_local_def}) of the numerical solutions, witch is presented in Fig.~\ref{fig:lin_diss} (\subref{fig:lin_diss:c1}, \subref{fig:lin_diss:c2}, \subref{fig:lin_diss:c3}, \subref{fig:lin_diss:c4}), shows that for polynomial degree $N = 1$, the errors of the improved local solution $\mathbf{u}_{\rm IL}$ are lower than those of the local solution $\mathbf{u}_{L}$ across the domain of definition $\Omega$, but this is particularly evident in the initial domain, where the improved local solution $\mathbf{u}_{\rm IL}$ demonstrates significantly higher accuracy. With increasing polynomial degree $N$, in particular, in the cases of polynomial degrees $N = 4$, $12$, $60$, the local error $\varepsilon$ of the improved local solution $\mathbf{u}_{\rm IL}$ is significantly lower than the local error $\varepsilon$ of the local solution $\mathbf{u}_{L}$, with the differences reaching $2$--$4$ orders of magnitude. Moreover, the local error $\varepsilon$ of the numerical solution $\mathbf{u}_{n}$ at the grid nodes $t_{n}$ is significantly smaller than the local errors $\varepsilon$ of the local solution $\mathbf{u}_{L}$ and the improved local solution $\mathbf{u}_{\rm IL}$, amounting to approximately $1$--$2$ orders of magnitude in the case of polynomial degree $N = 1$, approximately $5$--$8$ orders of magnitude in the case of polynomial degree $N = 4$, approximately $22$--$24$ orders of magnitude in the case of polynomial degree $N = 12$ and approximately $146$--$148$ orders of magnitude in the case of polynomial degree $N = 60$ (for this purpose, breaks in the graph along the vertical axis are inserted in Fig.~\ref{fig:lin_diss} (\subref{fig:lin_diss:c3}, \subref{fig:lin_diss:c4})).

The resulting dependencies of the global error $e$ (\ref{eq:eps_un_global_def}), (\ref{eq:eps_ul_global_def}) on the discretization step ${\Delta t}$, presented in Fig.~\ref{fig:lin_diss} (\subref{fig:lin_diss:d1}, \subref{fig:lin_diss:d2}, \subref{fig:lin_diss:d3}, \subref{fig:lin_diss:d4}, \subref{fig:lin_diss:e1}, \subref{fig:lin_diss:e2}, \subref{fig:lin_diss:e3}, \subref{fig:lin_diss:e4}, \subref{fig:lin_diss:f1}, \subref{fig:lin_diss:f2}, \subref{fig:lin_diss:f3}, \subref{fig:lin_diss:f4}), demonstrate a very high quality of linear approximation for all studied polynomial degrees $N$ (Fig.~\ref{fig:lin_diss} only shows results for polynomial degrees $N = 1$, $4$, $12$ and $60$). In the case of a numerical solution $\mathbf{u}_{n}$ in grid nodes for polynomials of degree $N = 60$, the global error $e$ for which is shown in Fig.~\ref{fig:lin_diss} (\subref{fig:lin_diss:f4}), a decrease in the discretization step ${\Delta t}$ by a factor of $2$ leads to a decrease in global errors $e$ by $40$--$50$ orders of magnitude.

Based on the approximation of the obtained dependencies $e({\Delta t})$ in log-log scale by a linear function $\lg{e({\Delta t})} \propto p\cdot\lg{{\Delta t}}$, empirical convergence orders $p$ are calculated and presented in Table~\ref{tab:conv_orders_lin_diss} for all polynomial degrees $N = 1, \ldots, 30$ and polynomial degrees $N$ up to $60$ with a step of $5$. A comparison of the obtained empirical convergence orders $p$ with the expected theoretical values $p_{\rm th.}$ (\ref{eq:conv_ords_exp}) shows excellent agreement. It is particularly noticeable that the empirical convergence orders $p^{\rm imp}$ of the improved local numerical solution $\mathbf{u}_{\rm IL}$ are one unit higher than the empirical convergence orders $p^{l}$ of the local numerical solution $\mathbf{u}_{L}$.

\corrtext{The system of equations (\ref{eq:lin_diss_ode}) considered in this Example can be presented in a conservative form, which corresponds to the mechanical problem of obtaining the law of motion of a particle of unit mass in an external field with potential energy $U(x)$~\cite{LandauMechanics, GoldsteinMechanics}:
\begin{equation}
\ddot{x} = -\frac{dU(x)}{dx},\qquad U(x) = -\frac{x^{2}}{2},
\end{equation}
which, together with the initial conditions (\ref{eq:lin_diss_ode}), shows that for the exact solution the energy conservation law $E(t) = \mathrm{const}$ of the following form holds
\begin{equation}\label{eq:econs_lin_diss}
E(t) = \frac{\dot{x}^{2}(t)}{2} + U(x) = \frac12\left[\dot{x}^{2}(t) - x^{2}(t)\right] = \mathrm{const} \equiv 
E(t_{0}) = \frac12\left[\dot{x}^{2}(0) - x^{2}(0)\right] = \frac12,
\end{equation}
which represents the integral of motion. It is known~\cite{Butcher_book_2016, Hairer_book_1, Hairer_book_2} that satisfying the energy conservation law, as well as other symmetries, in a numerical solution in some form of difference analog is an important property of the numerical method for solving the initial value problem for an ODE system (\ref{eq:ivp_ode_diff_src}).}

\corrtext{Fig.~\ref{fig:econs_lin_diss} shows the dependence of the negative logarithm of the error $-\lg|E(t)-E(t_{0})|$ of the energy conservation law (\ref{eq:econs_lin_diss}) of the numerical solution of the system (\ref{eq:lin_diss_ode}) on the argument $t$ for solution at nodes $\mathbf{u}_{n}$, the local solution $\mathbf{u}_{L}(t)$ and the improved local solution $\mathbf{u}_{\rm IL}(t)$, obtained using polynomials with degrees $N = 1$~(\subref{fig:econs_lin_diss:a1}), $N = 4$~(\subref{fig:econs_lin_diss:a4}), $N = 12$~(\subref{fig:econs_lin_diss:a3}), $N = 60$~(\subref{fig:econs_lin_diss:a4}). The obtained results clearly demonstrate that the energy conservation law (\ref{eq:lin_diss_ode}) is not strictly satisfied in the numerical solution, which is due to the dissipative nature of the ADER-DG numerical method with local DG predictor~\cite{ader_dg_ode_jsc, ader_dg_ode_sinum, ader_improving_2024, ader_proofs_2025}. However, due to the possibility of achieving a high order $p$, even in the case of degree $N = 4$, the error $-\lg|E(t)-E(t_{0})|$ in satisfying the energy conservation law (\ref{eq:econs_lin_diss}) over the entire solution domain becomes less than the characteristic rounding error of single-precision floating-point numbers $\sim 10^{-7}$-$10^{-9}$. For degrees $N = 12$ and $60$, the error in the energy conservation law's fulfillment over the entire solution domain becomes smaller than the characteristic rounding error of double-precision floating-point numbers $\sim 10^{-15}$-$10^{-17}$. The presented results also clearly demonstrate that the accuracy of the fulfillment of the energy conservation law (\ref{eq:econs_lin_diss}) for the improved local solution $\mathbf{u}_{\rm IL}(t)$ is significantly higher than for the local solution $\mathbf{u}_{L}(t)$. Therefore, it can be concluded that, despite the dissipative nature of the ADER-DG numerical method with local DG predictor, a sufficiently high degree $N$ can be chosen such that the accuracy of the fulfillment of the energy conservation law will be at or below the characteristic error of representing real numbers as floating-point numbers.}

The obtained results allowed to conclude that the ADER-DG numerical method with a local DG predictor provides a highly accurate numerical solution to the initial value problem for the ODE system (\ref{eq:lin_diss_ode}) presented in this example. The obtained results are in good agreement with the theory developed above. The improved local numerical solution $\mathbf{u}_{\rm IL}$ indeed demonstrates higher accuracy and a higher convergence order compared to the local numerical solution $\mathbf{u}_{L}$.

\subsection{Example 3: Harmonic oscillator}
\label{sec:apps:harm_osc}

\begin{figure}[h!]
\captionsetup[subfigure]{%
	position=bottom,
	font+=smaller,
	textfont=normalfont,
	singlelinecheck=off,
	justification=raggedright
}
\centering
\begin{subfigure}{0.24\textwidth}
\includegraphics[width=\textwidth]{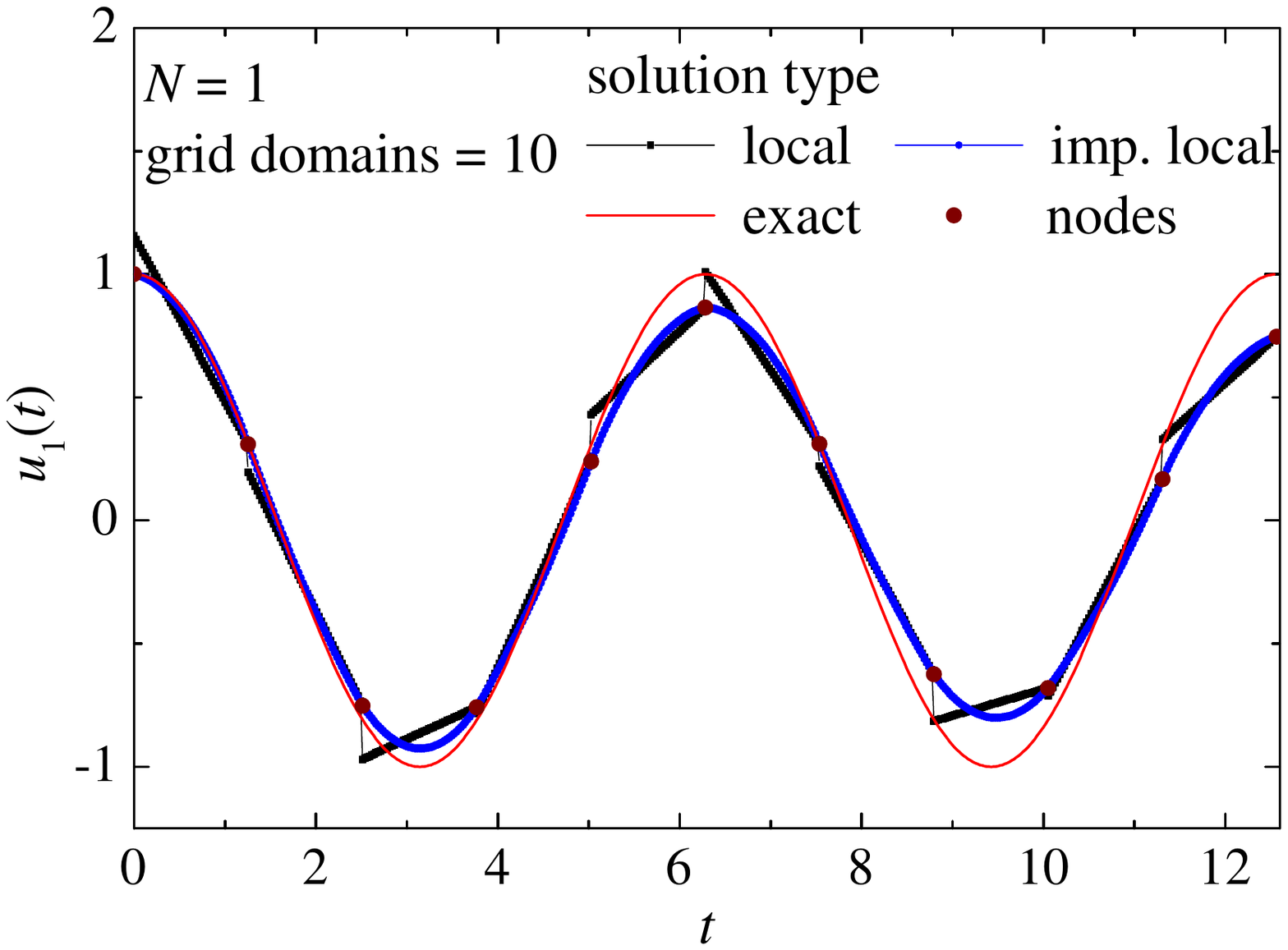}
\vspace{-8mm}\caption{\label{fig:harm_osc:a1}}
\end{subfigure}
\begin{subfigure}{0.24\textwidth}
\includegraphics[width=\textwidth]{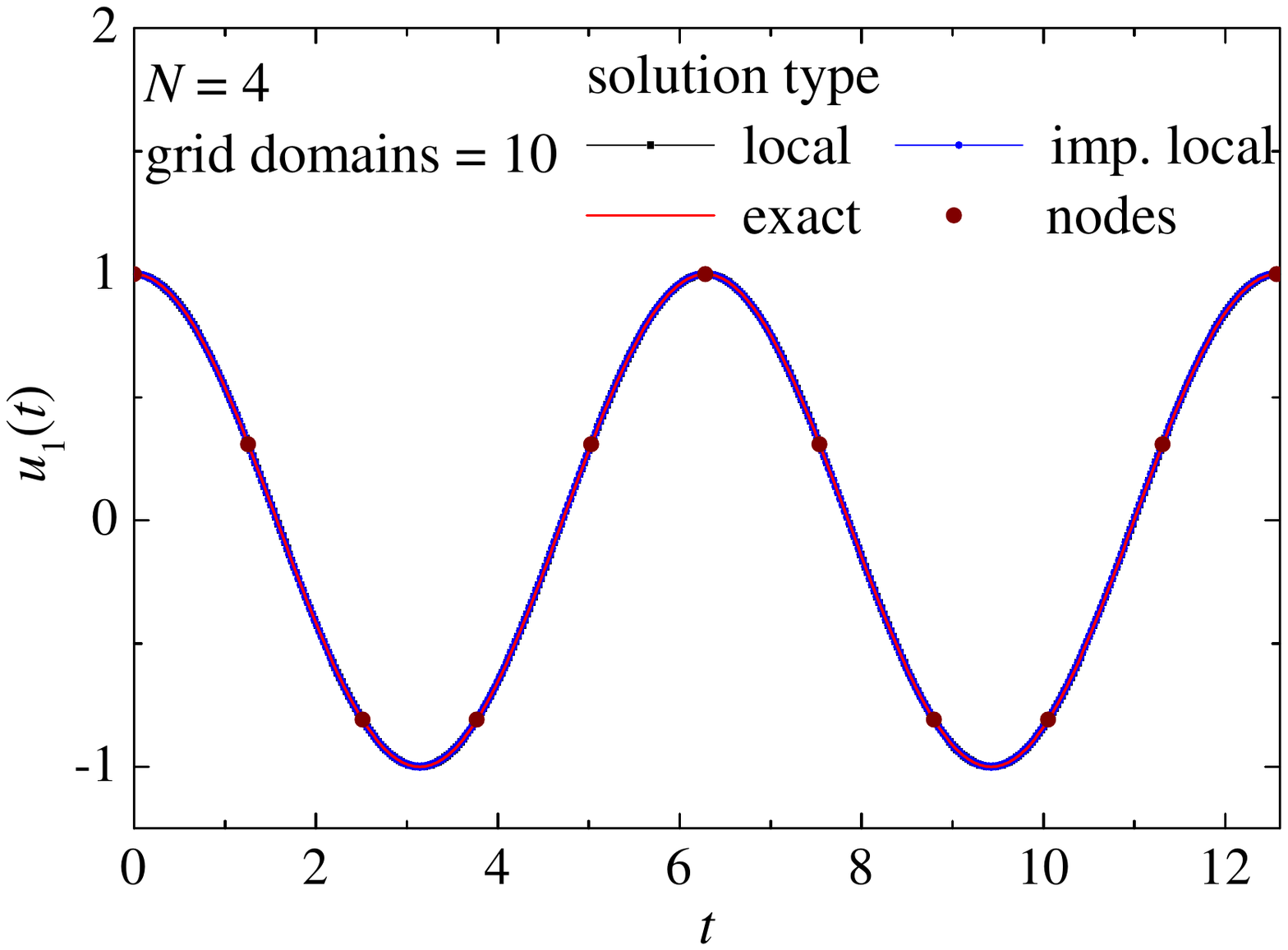}
\vspace{-8mm}\caption{\label{fig:harm_osc:a2}}
\end{subfigure}
\begin{subfigure}{0.24\textwidth}
\includegraphics[width=\textwidth]{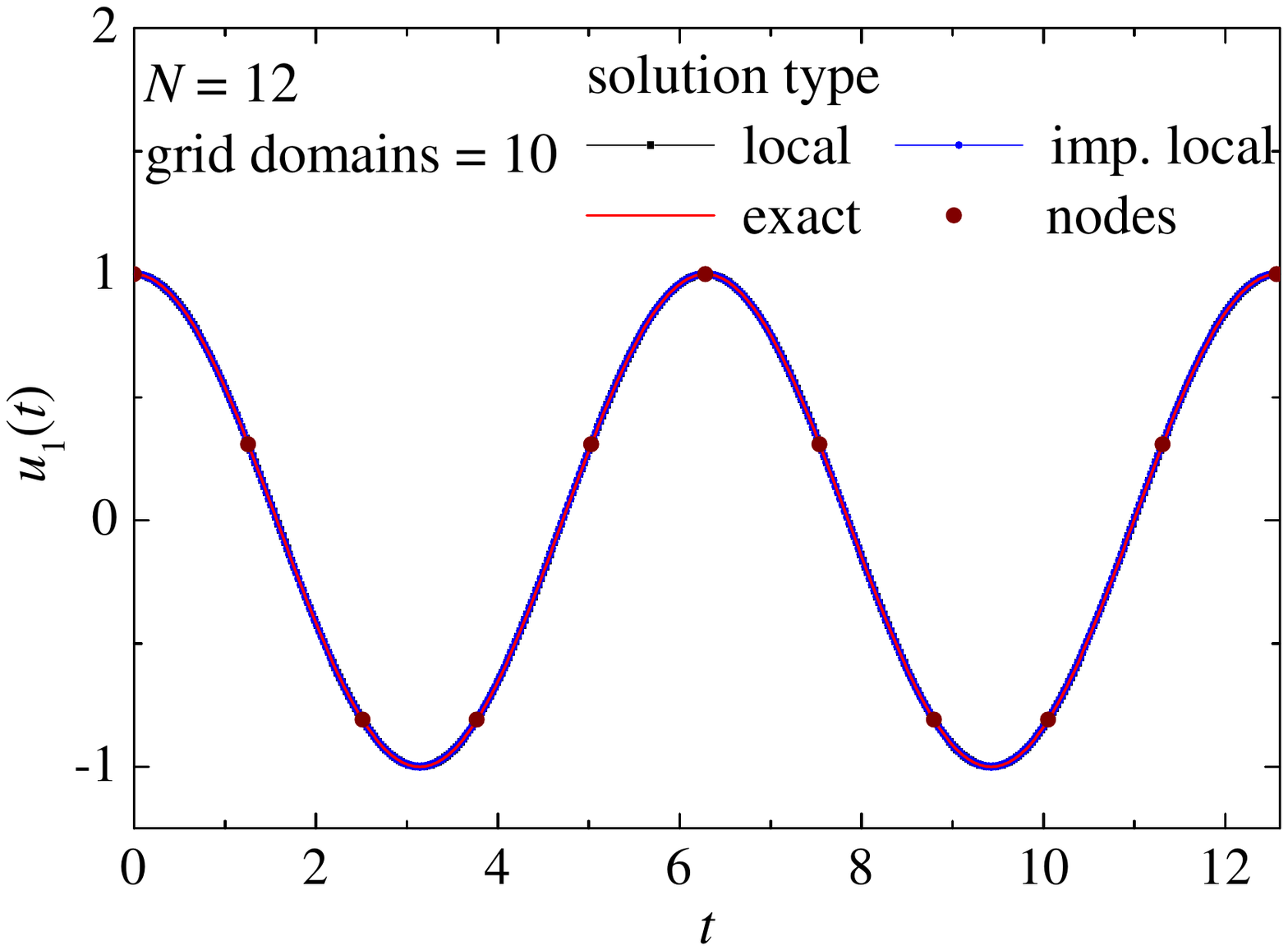}
\vspace{-8mm}\caption{\label{fig:harm_osc:a3}}
\end{subfigure}
\begin{subfigure}{0.24\textwidth}
\includegraphics[width=\textwidth]{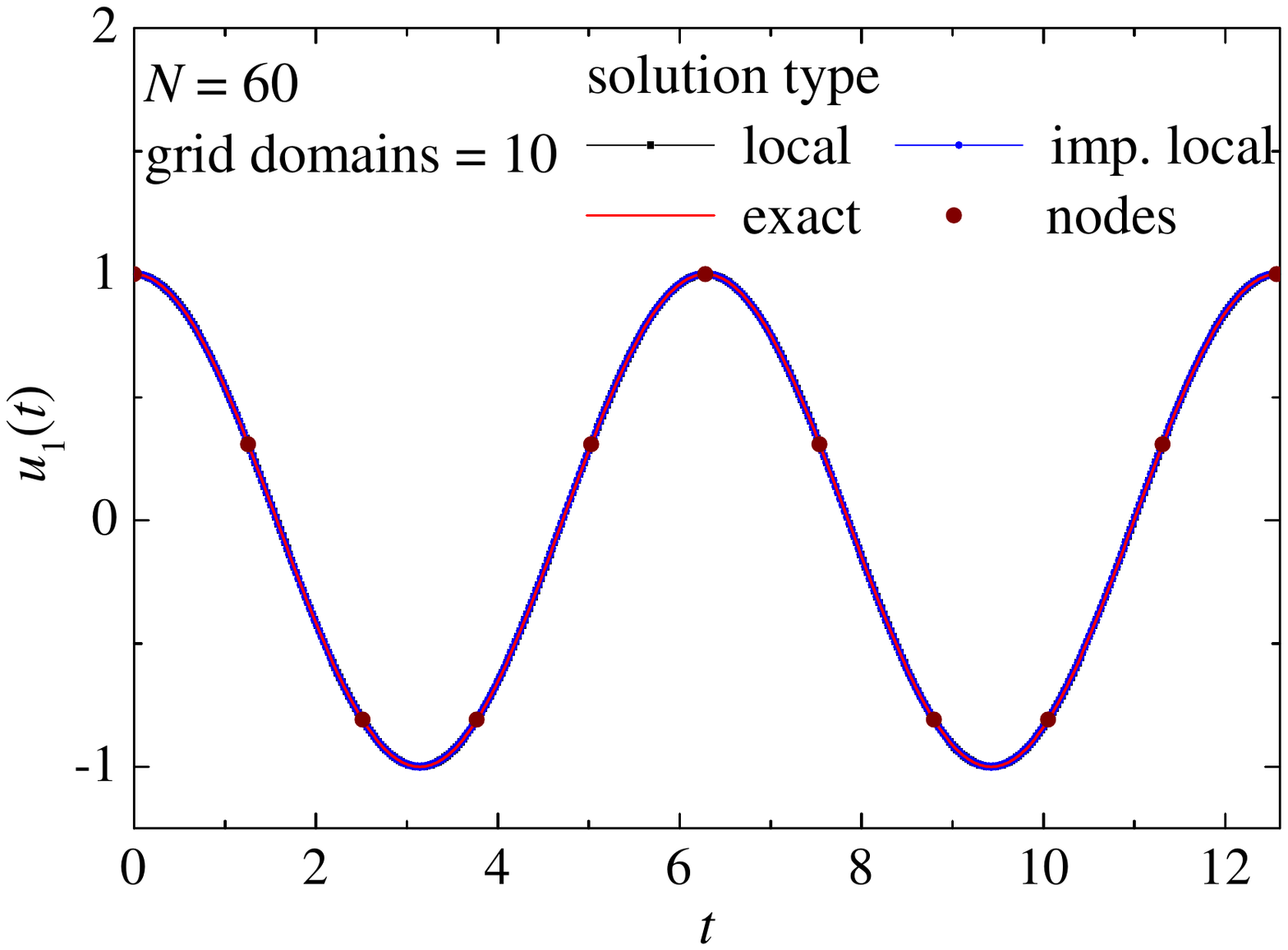}
\vspace{-8mm}\caption{\label{fig:harm_osc:a4}}
\end{subfigure}\\
\begin{subfigure}{0.24\textwidth}
\includegraphics[width=\textwidth]{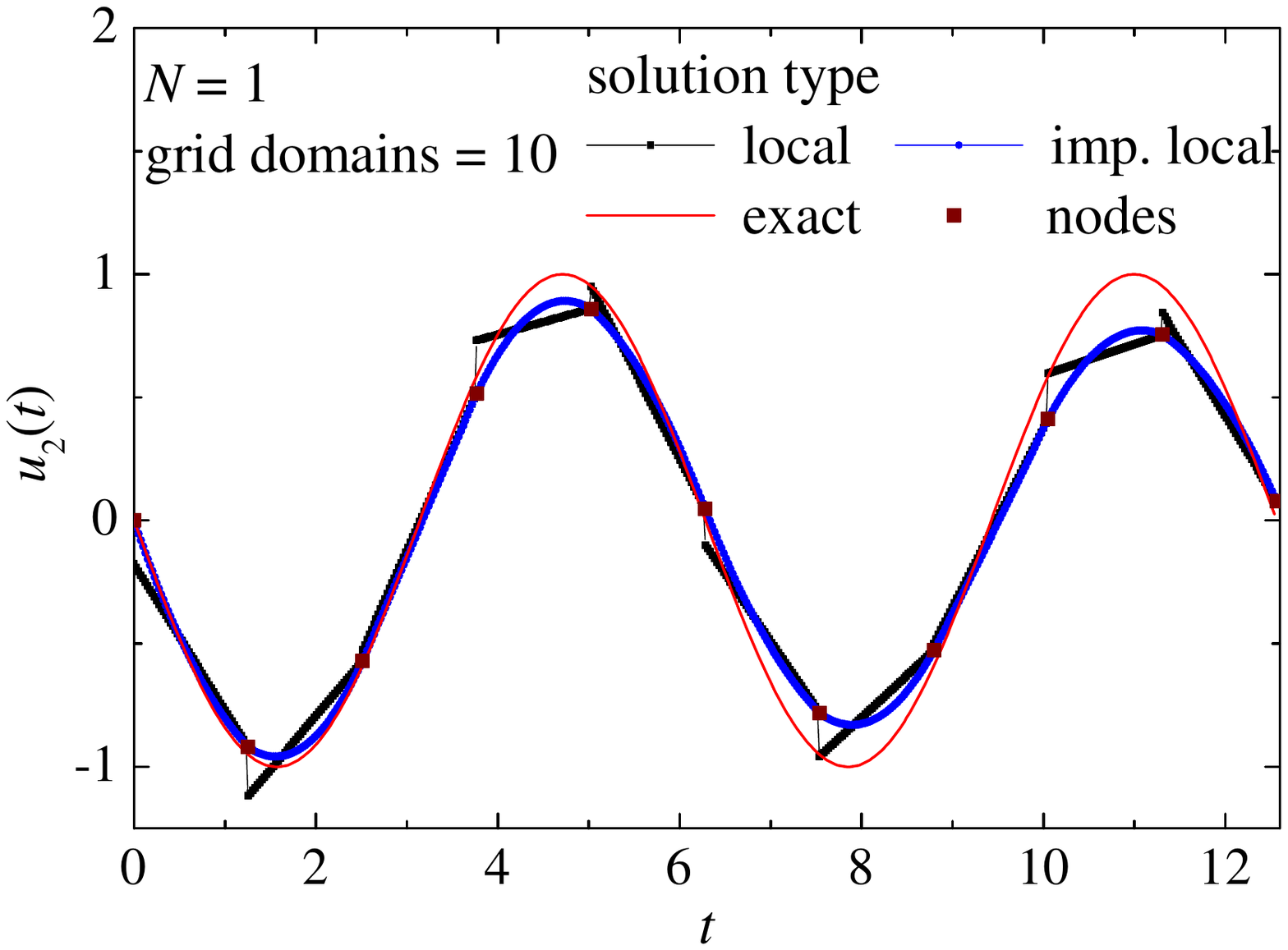}
\vspace{-8mm}\caption{\label{fig:harm_osc:b1}}
\end{subfigure}
\begin{subfigure}{0.24\textwidth}
\includegraphics[width=\textwidth]{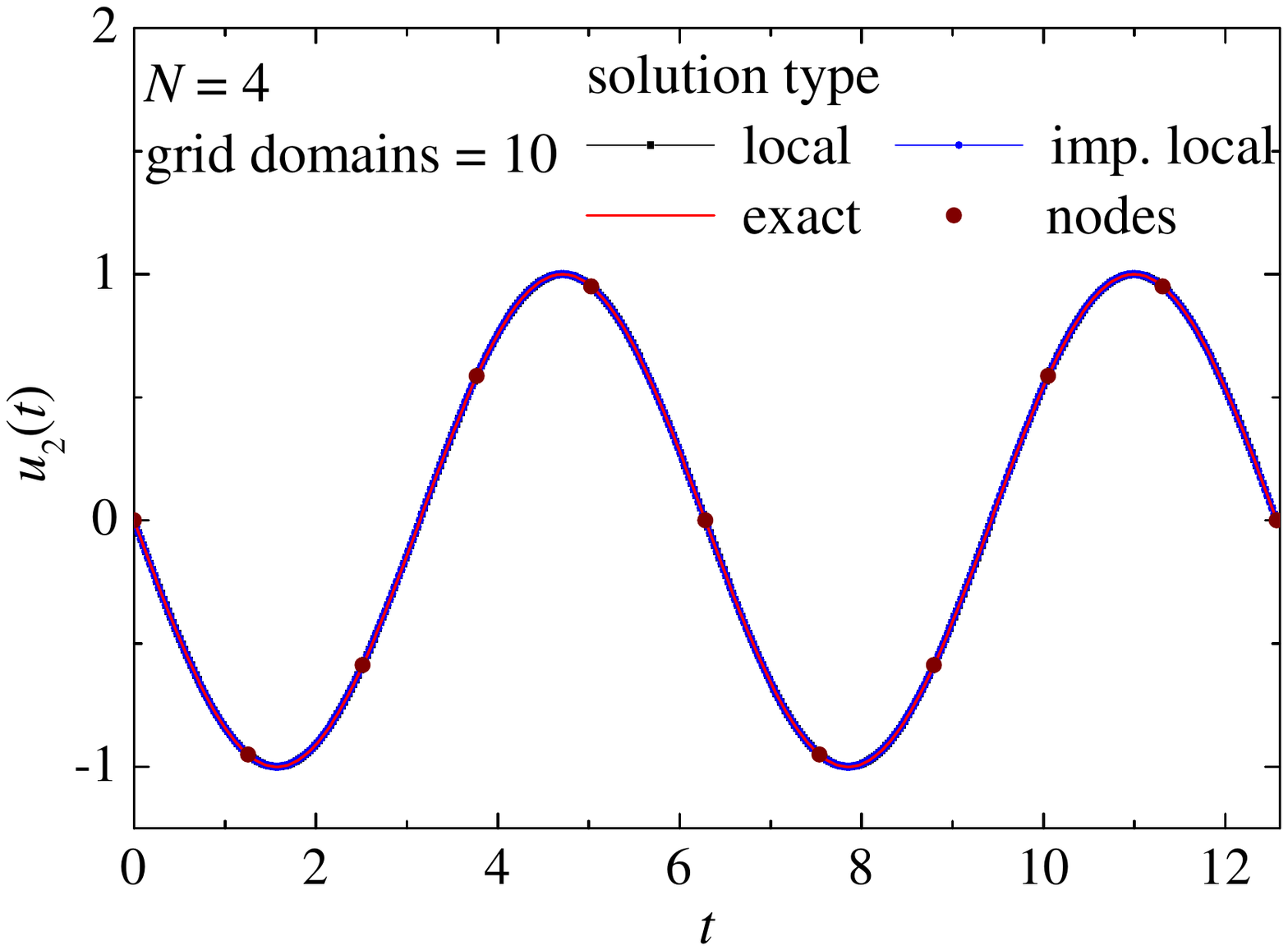}
\vspace{-8mm}\caption{\label{fig:harm_osc:b2}}
\end{subfigure}
\begin{subfigure}{0.24\textwidth}
\includegraphics[width=\textwidth]{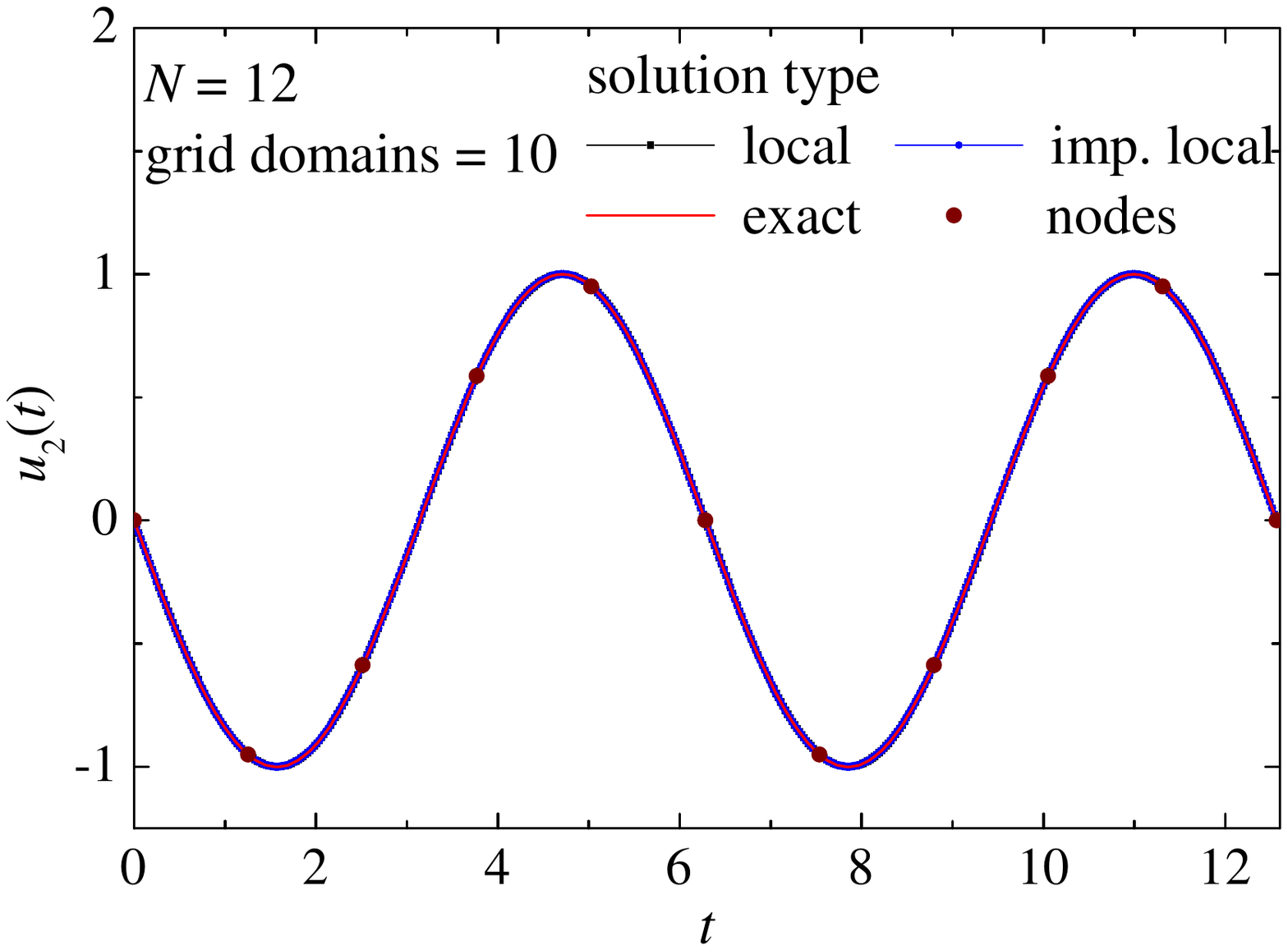}
\vspace{-8mm}\caption{\label{fig:harm_osc:b3}}
\end{subfigure}
\begin{subfigure}{0.24\textwidth}
\includegraphics[width=\textwidth]{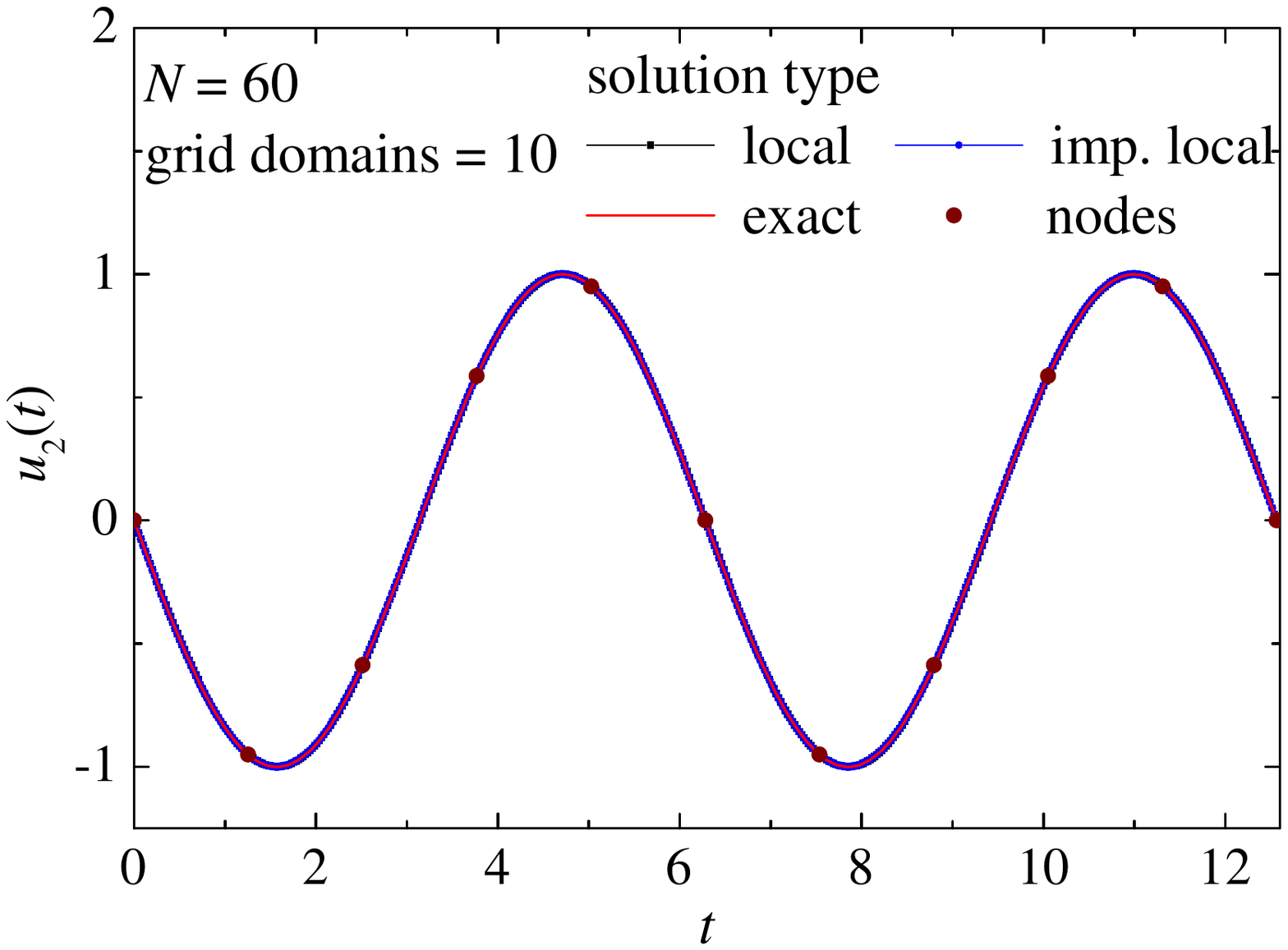}
\vspace{-8mm}\caption{\label{fig:harm_osc:b4}}
\end{subfigure}\\
\begin{subfigure}{0.24\textwidth}
\includegraphics[width=\textwidth]{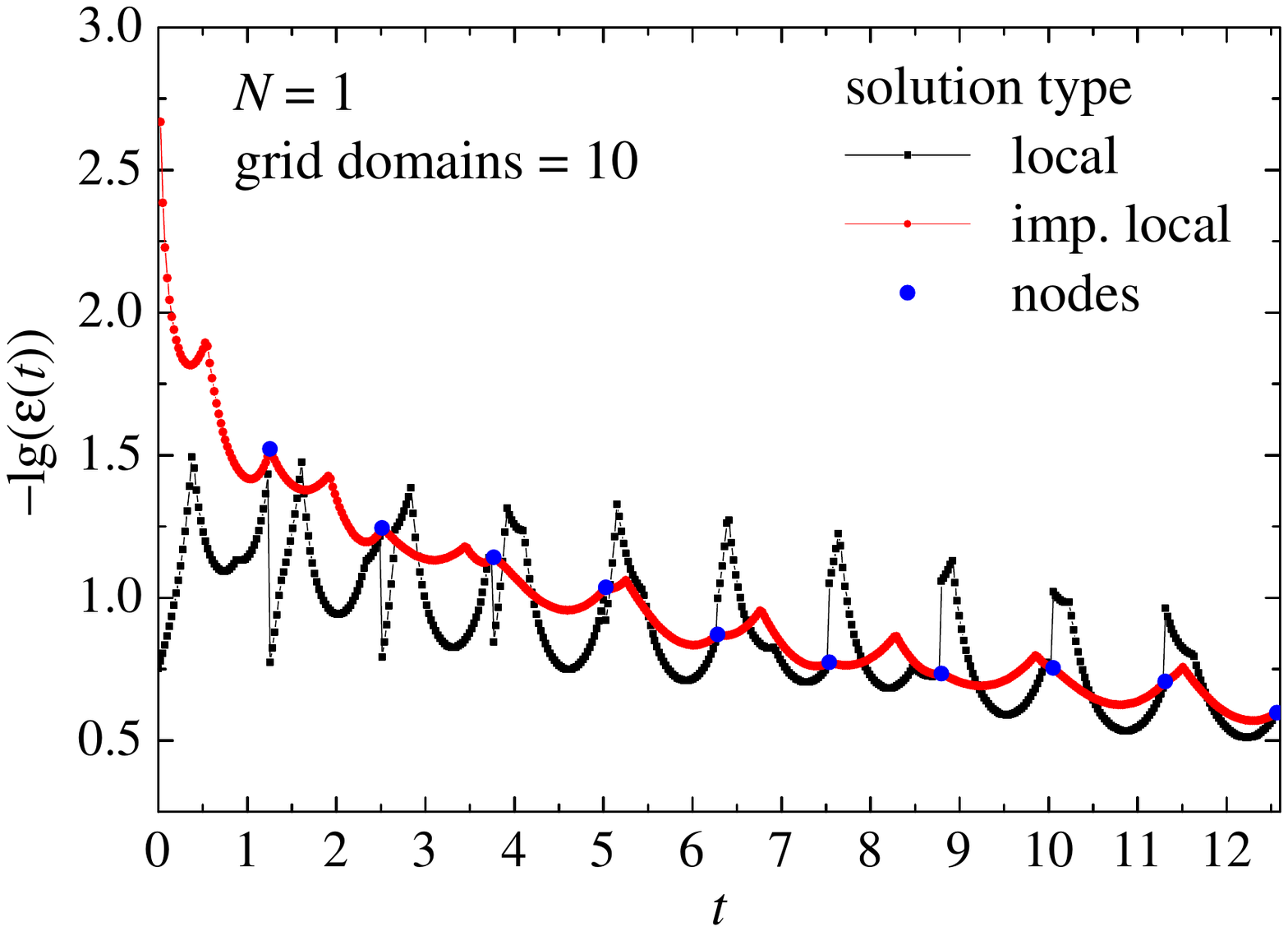}
\vspace{-8mm}\caption{\label{fig:harm_osc:c1}}
\end{subfigure}
\begin{subfigure}{0.24\textwidth}
\includegraphics[width=\textwidth]{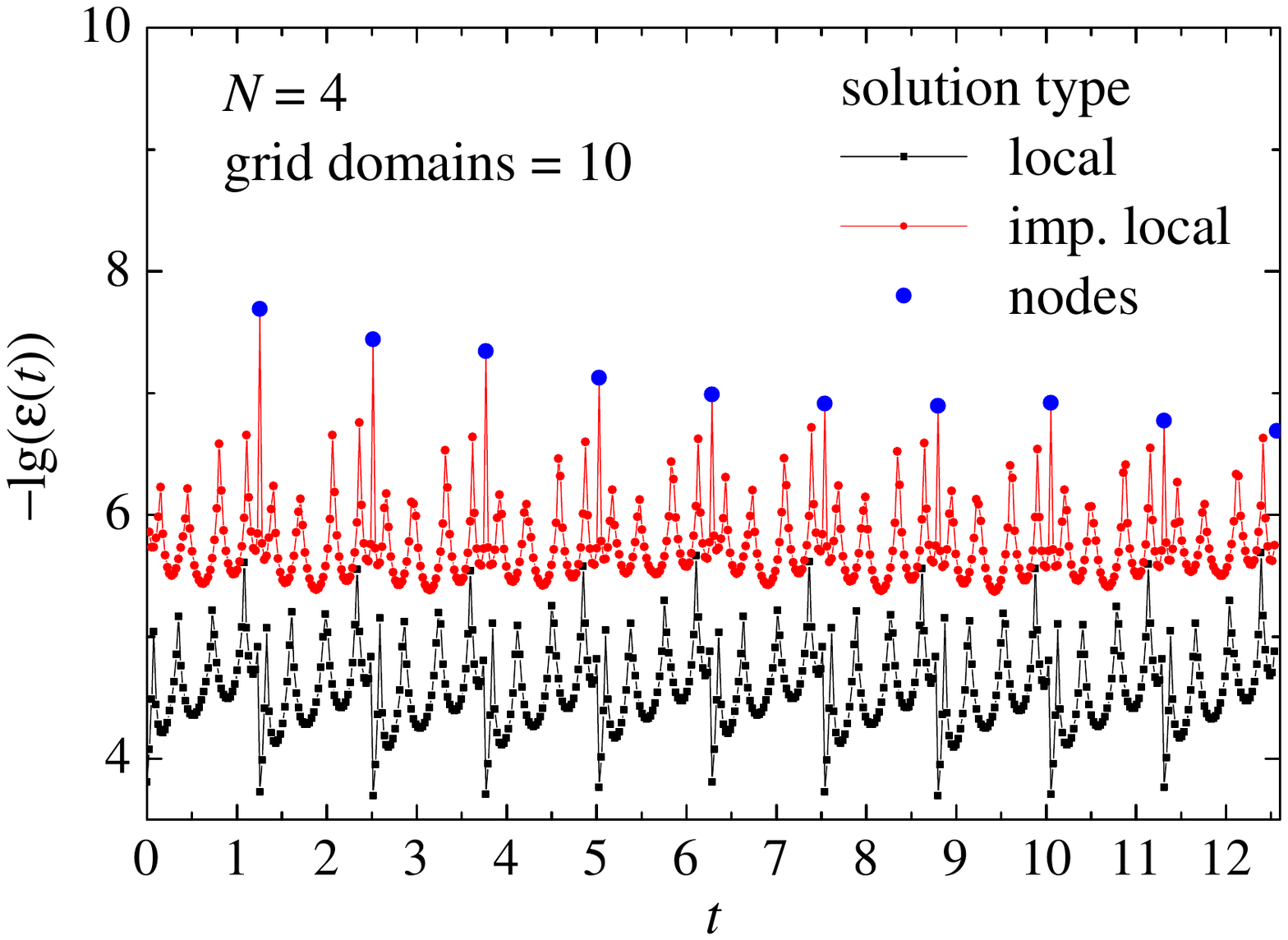}
\vspace{-8mm}\caption{\label{fig:harm_osc:c2}}
\end{subfigure}
\begin{subfigure}{0.24\textwidth}
\includegraphics[width=\textwidth]{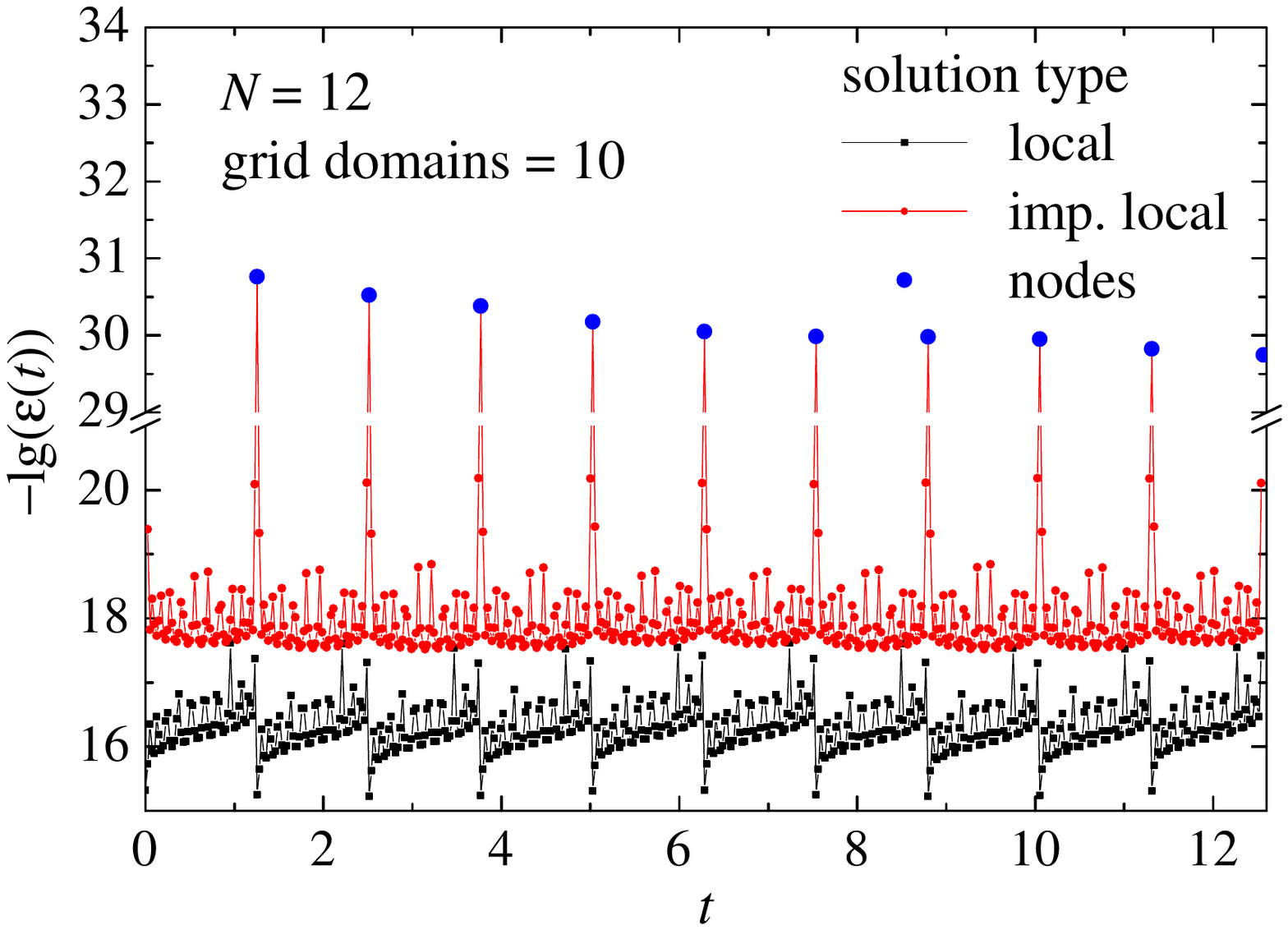}
\vspace{-8mm}\caption{\label{fig:harm_osc:c3}}
\end{subfigure}
\begin{subfigure}{0.24\textwidth}
\includegraphics[width=\textwidth]{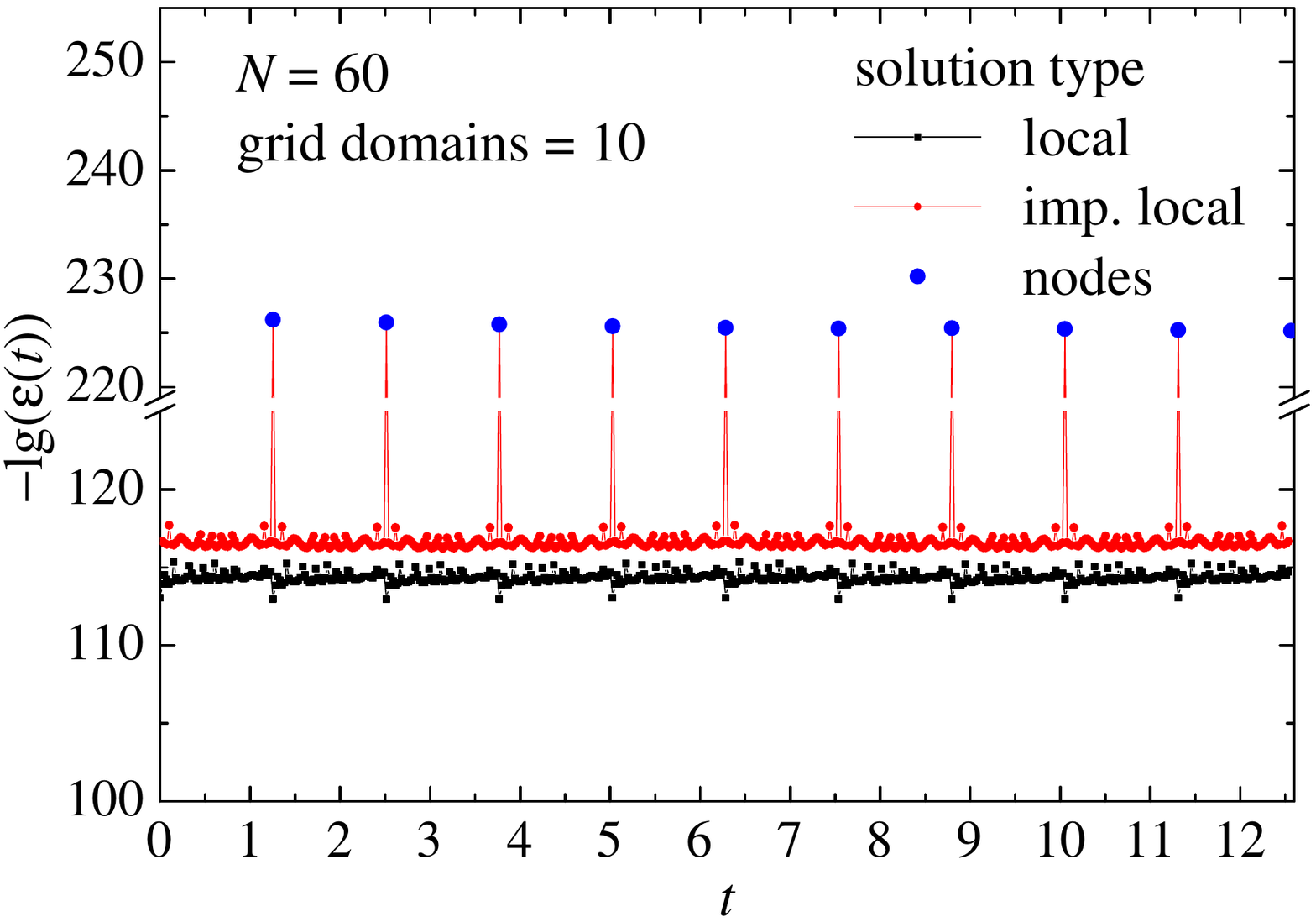}
\vspace{-8mm}\caption{\label{fig:harm_osc:c4}}
\end{subfigure}\\
\begin{subfigure}{0.24\textwidth}
\includegraphics[width=\textwidth]{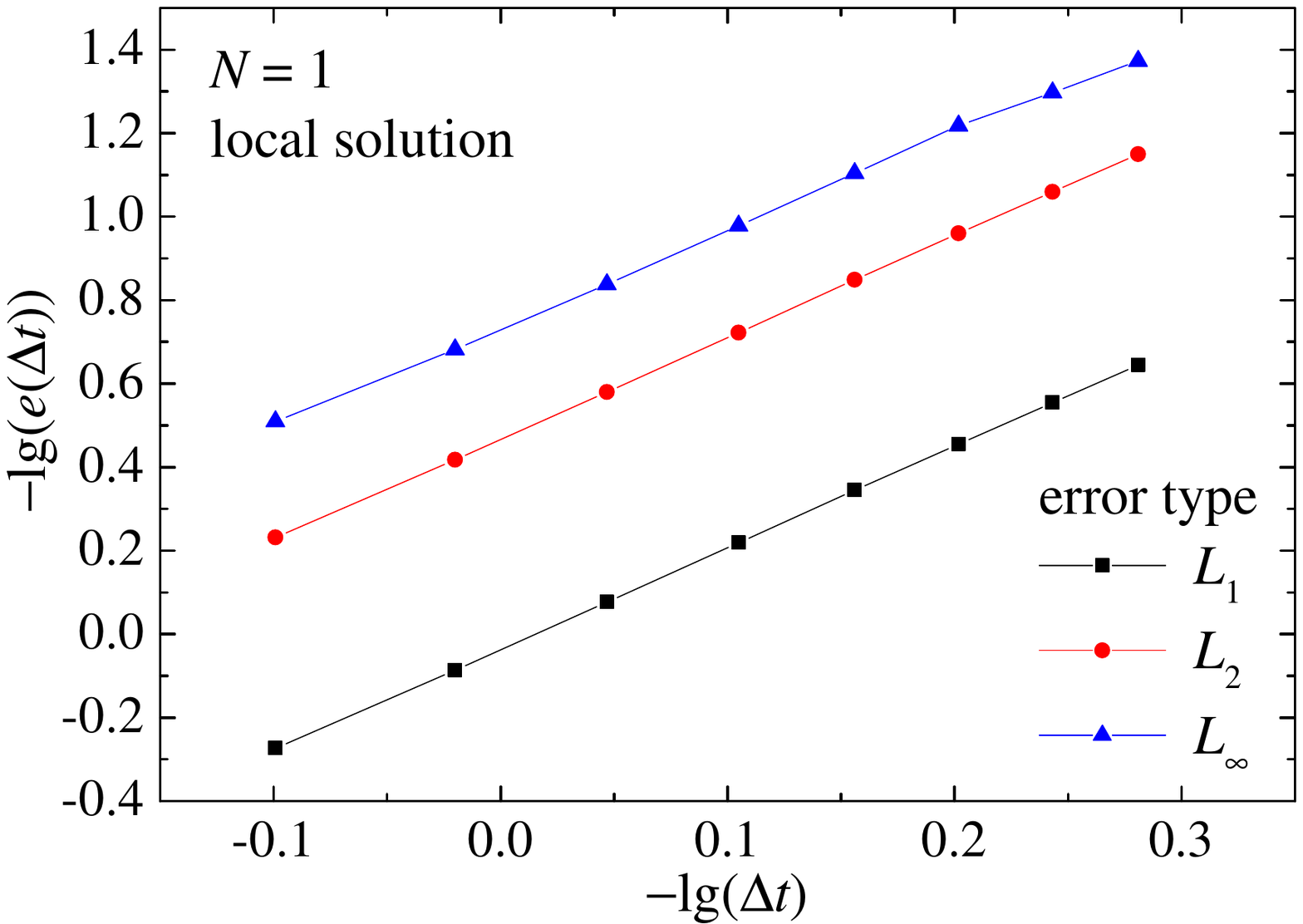}
\vspace{-8mm}\caption{\label{fig:harm_osc:d1}}
\end{subfigure}
\begin{subfigure}{0.24\textwidth}
\includegraphics[width=\textwidth]{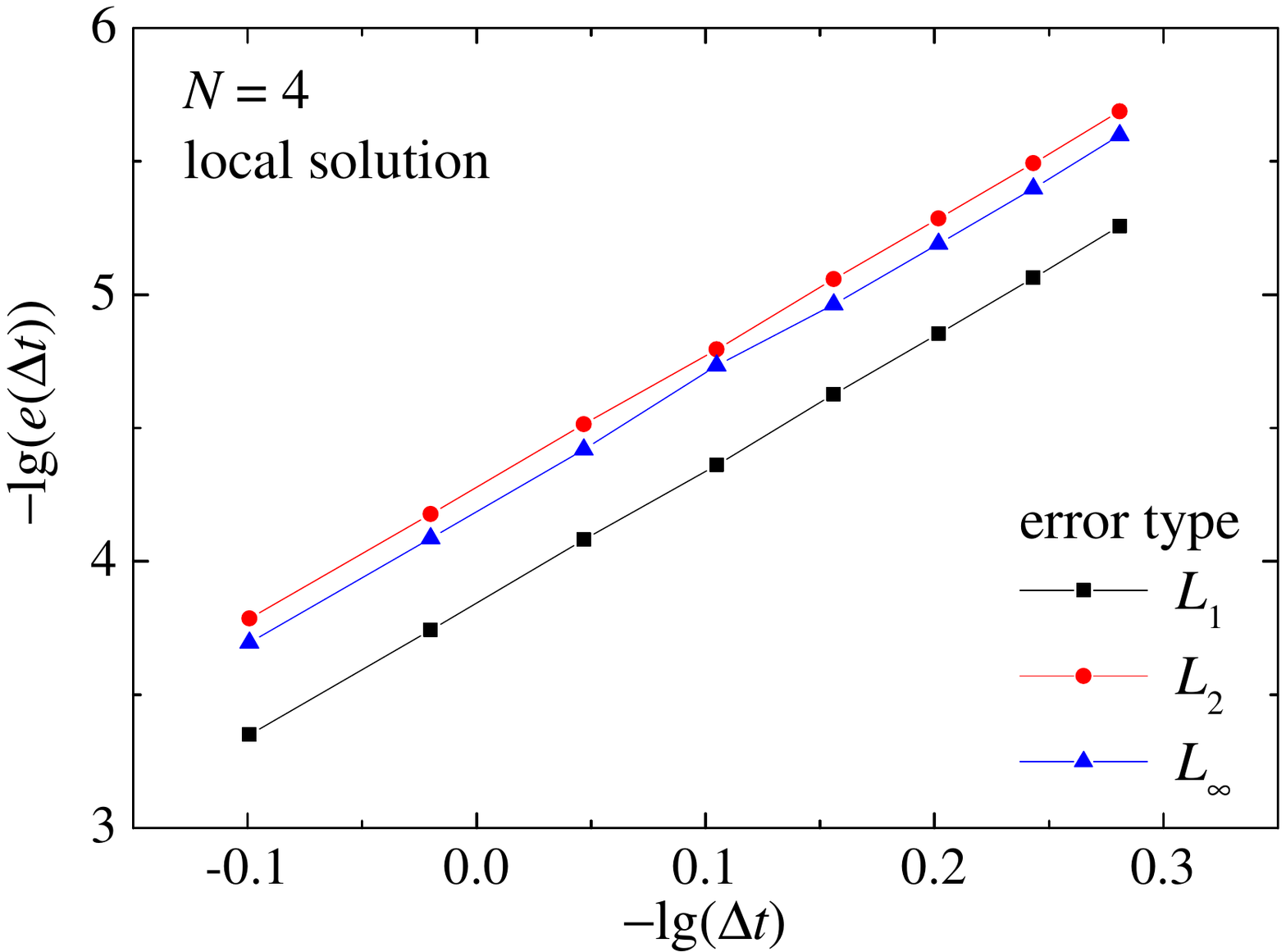}
\vspace{-8mm}\caption{\label{fig:harm_osc:d2}}
\end{subfigure}
\begin{subfigure}{0.24\textwidth}
\includegraphics[width=\textwidth]{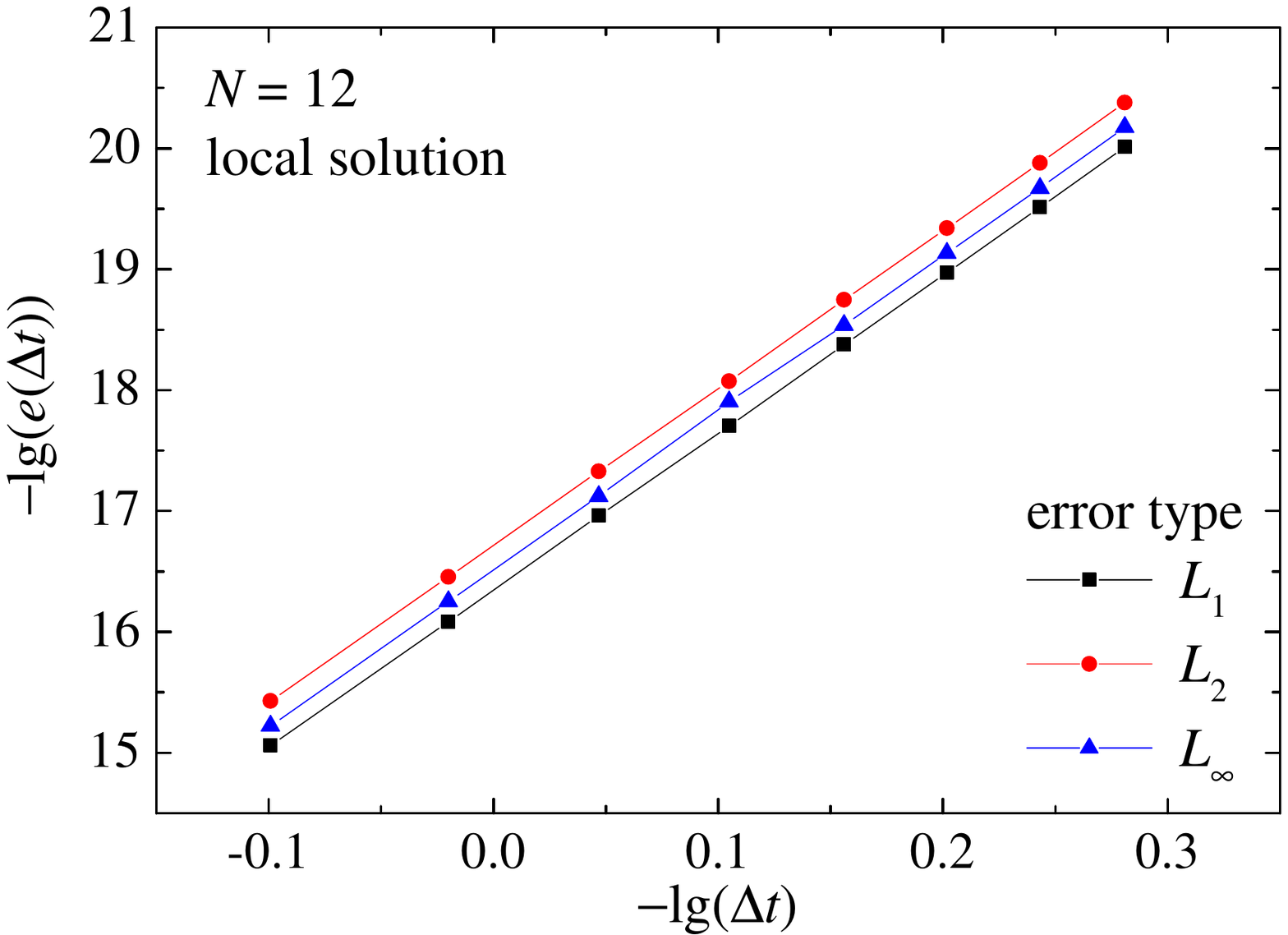}
\vspace{-8mm}\caption{\label{fig:harm_osc:d3}}
\end{subfigure}
\begin{subfigure}{0.24\textwidth}
\includegraphics[width=\textwidth]{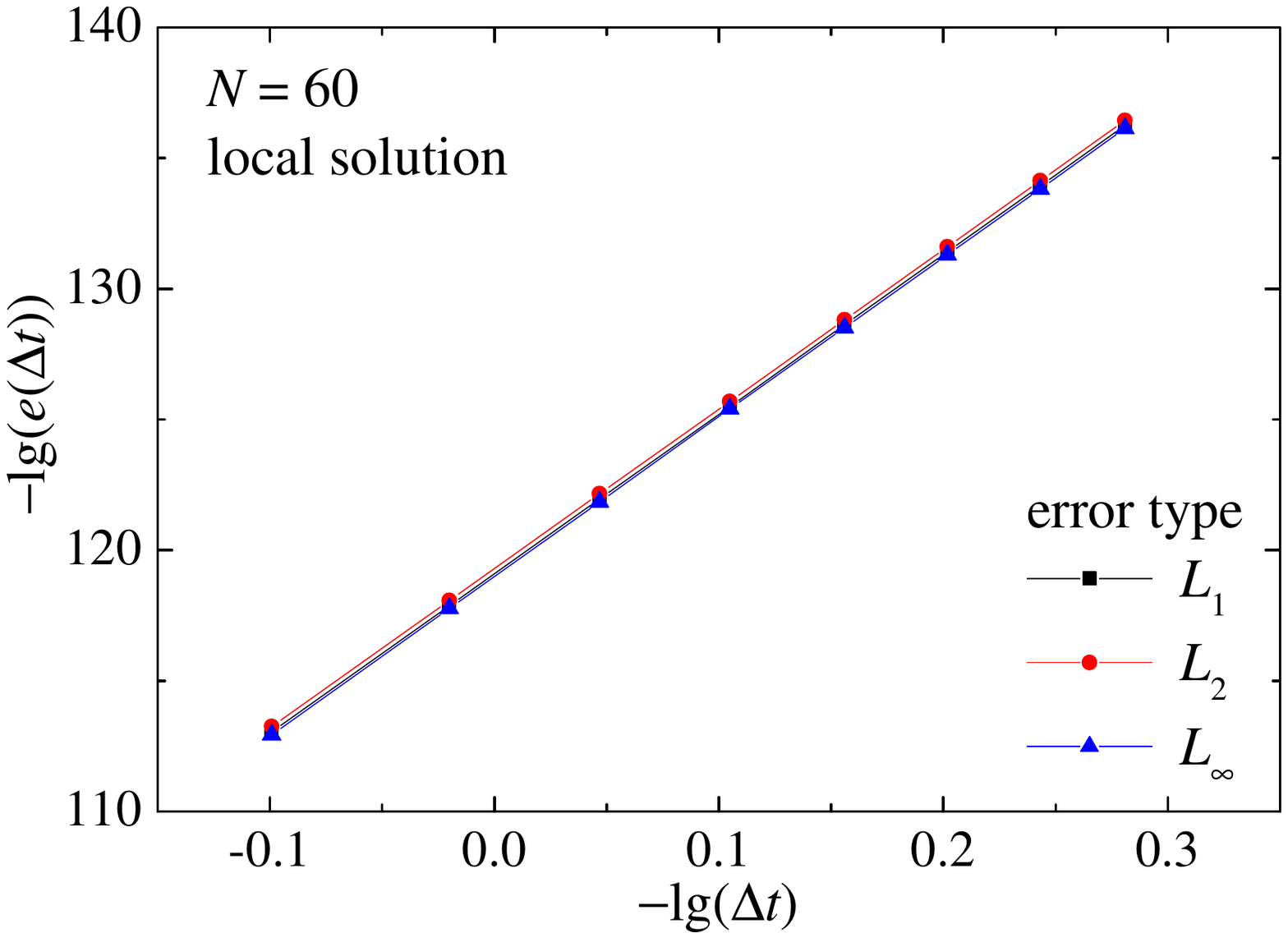}
\vspace{-8mm}\caption{\label{fig:harm_osc:d4}}
\end{subfigure}\\
\begin{subfigure}{0.24\textwidth}
\includegraphics[width=\textwidth]{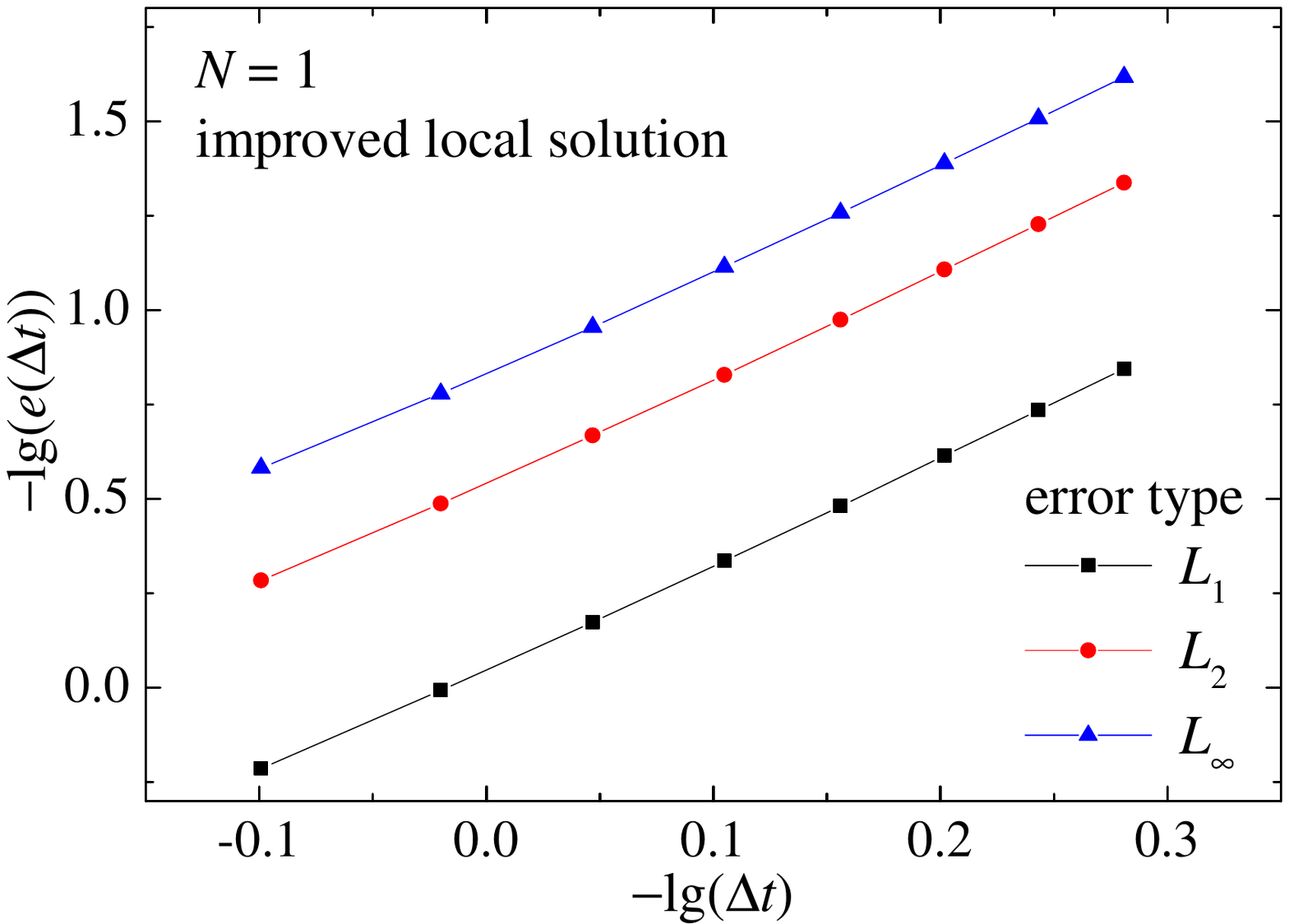}
\vspace{-8mm}\caption{\label{fig:harm_osc:e1}}
\end{subfigure}
\begin{subfigure}{0.24\textwidth}
\includegraphics[width=\textwidth]{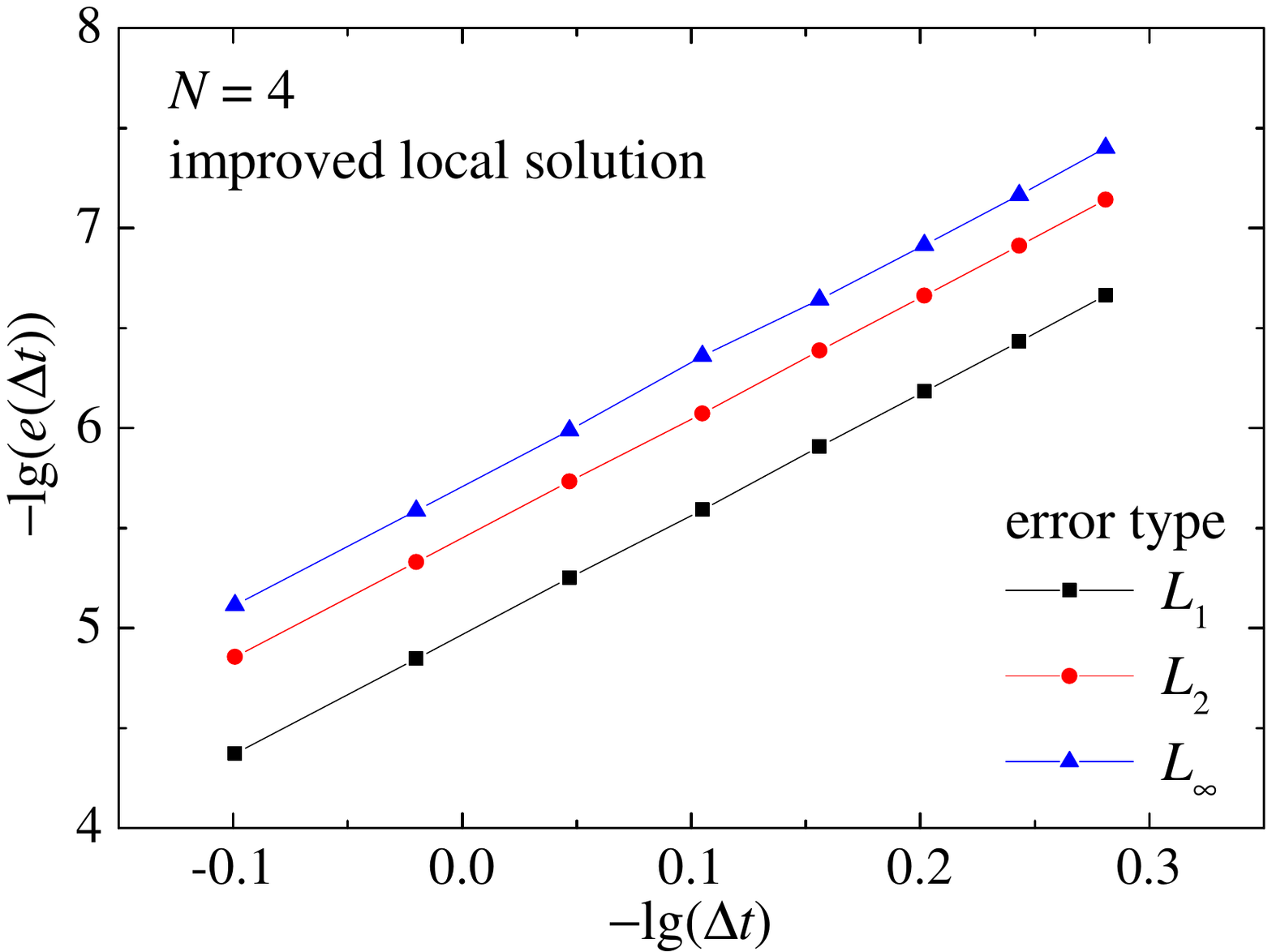}
\vspace{-8mm}\caption{\label{fig:harm_osc:e2}}
\end{subfigure}
\begin{subfigure}{0.24\textwidth}
\includegraphics[width=\textwidth]{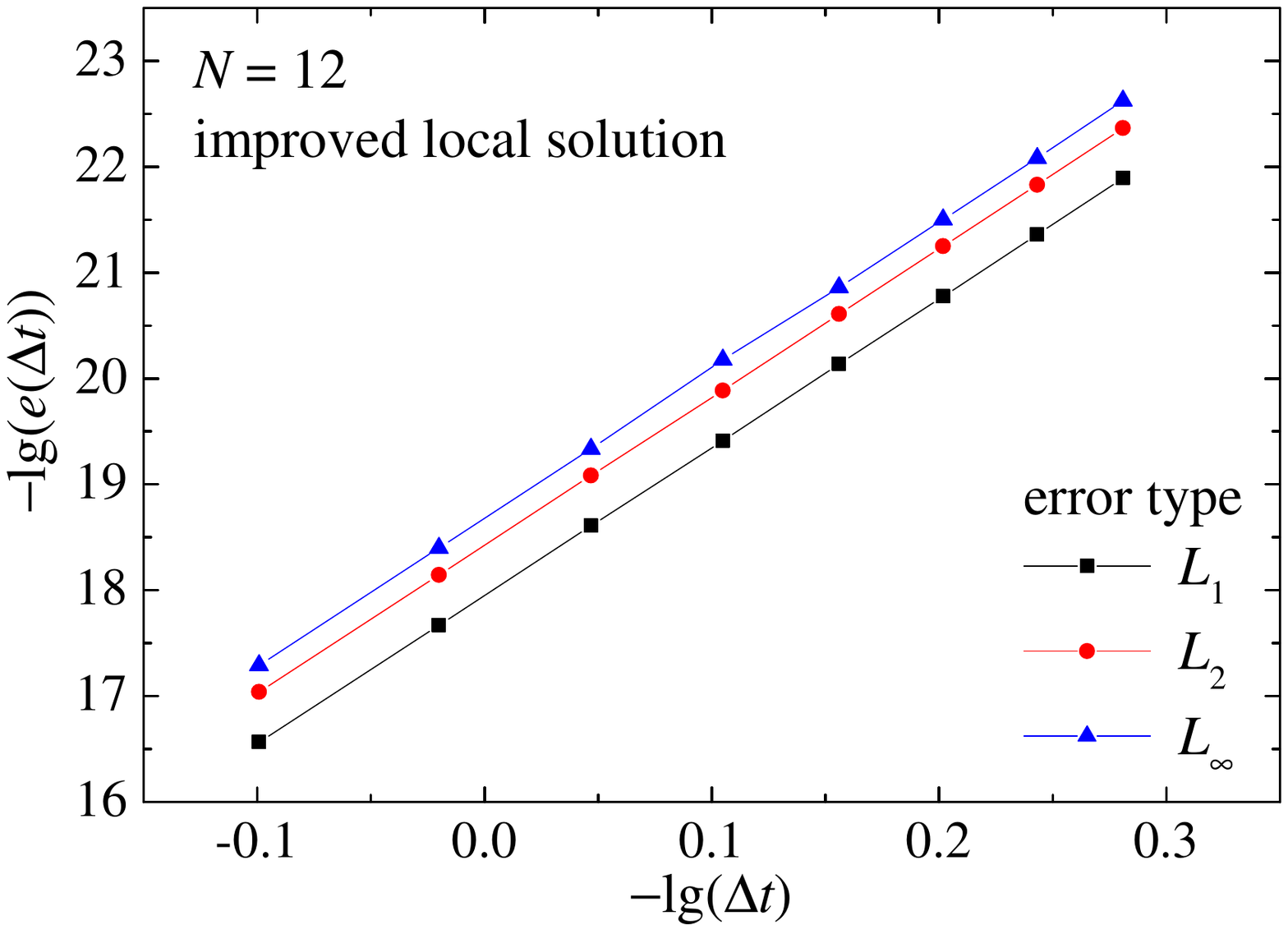}
\vspace{-8mm}\caption{\label{fig:harm_osc:e3}}
\end{subfigure}
\begin{subfigure}{0.24\textwidth}
\includegraphics[width=\textwidth]{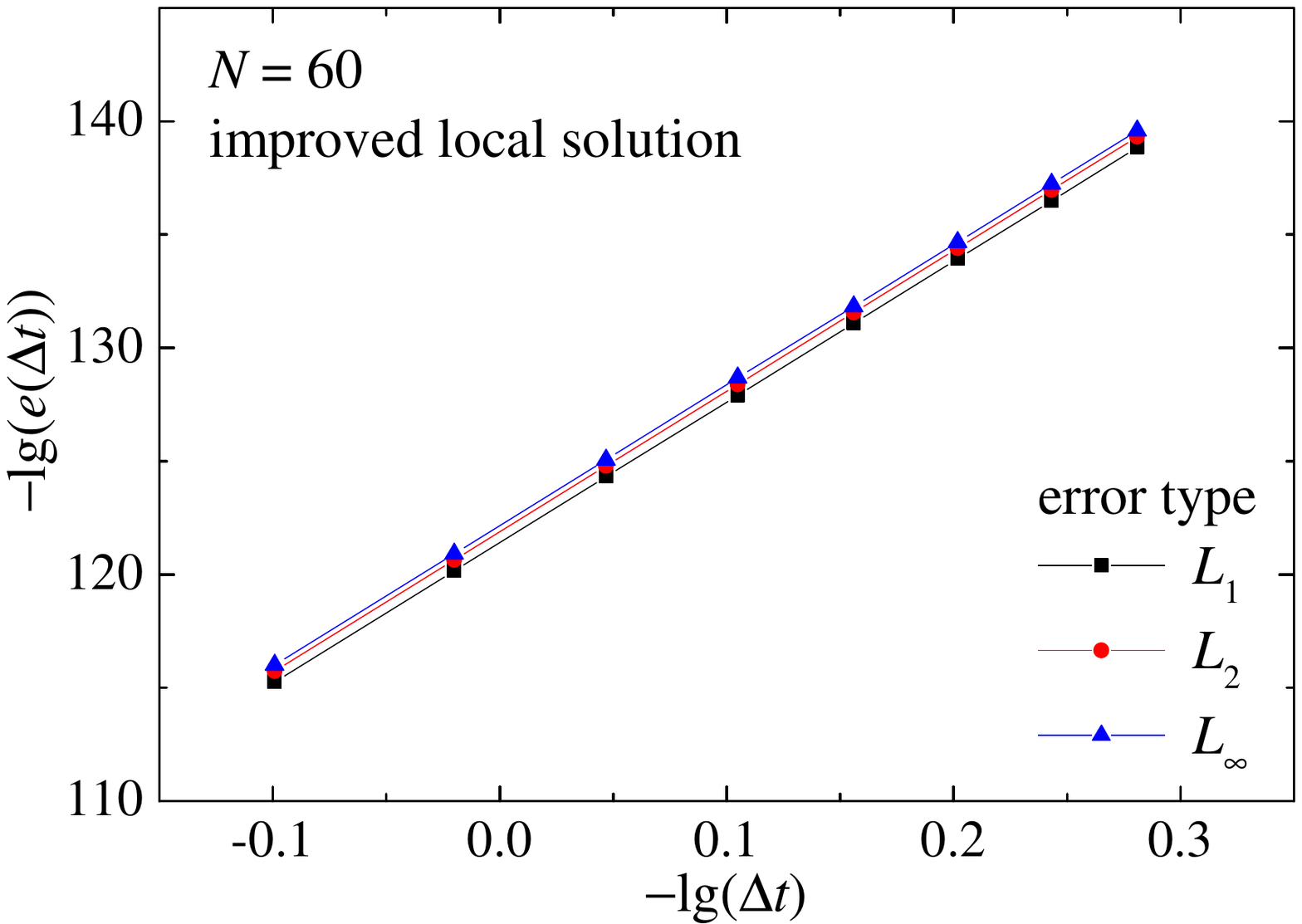}
\vspace{-8mm}\caption{\label{fig:harm_osc:e4}}
\end{subfigure}\\
\begin{subfigure}{0.24\textwidth}
\includegraphics[width=\textwidth]{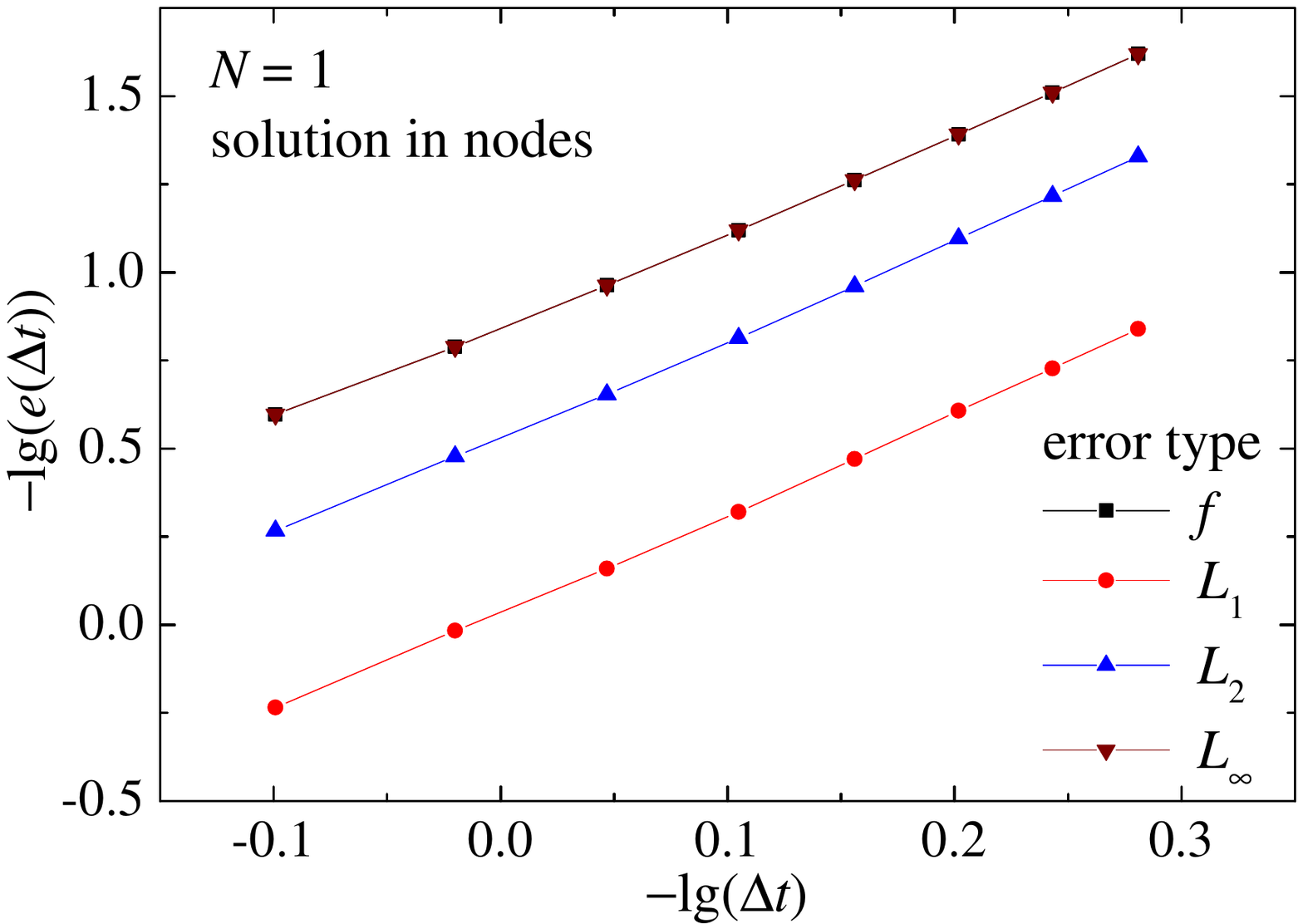}
\vspace{-8mm}\caption{\label{fig:harm_osc:f1}}
\end{subfigure}
\begin{subfigure}{0.24\textwidth}
\includegraphics[width=\textwidth]{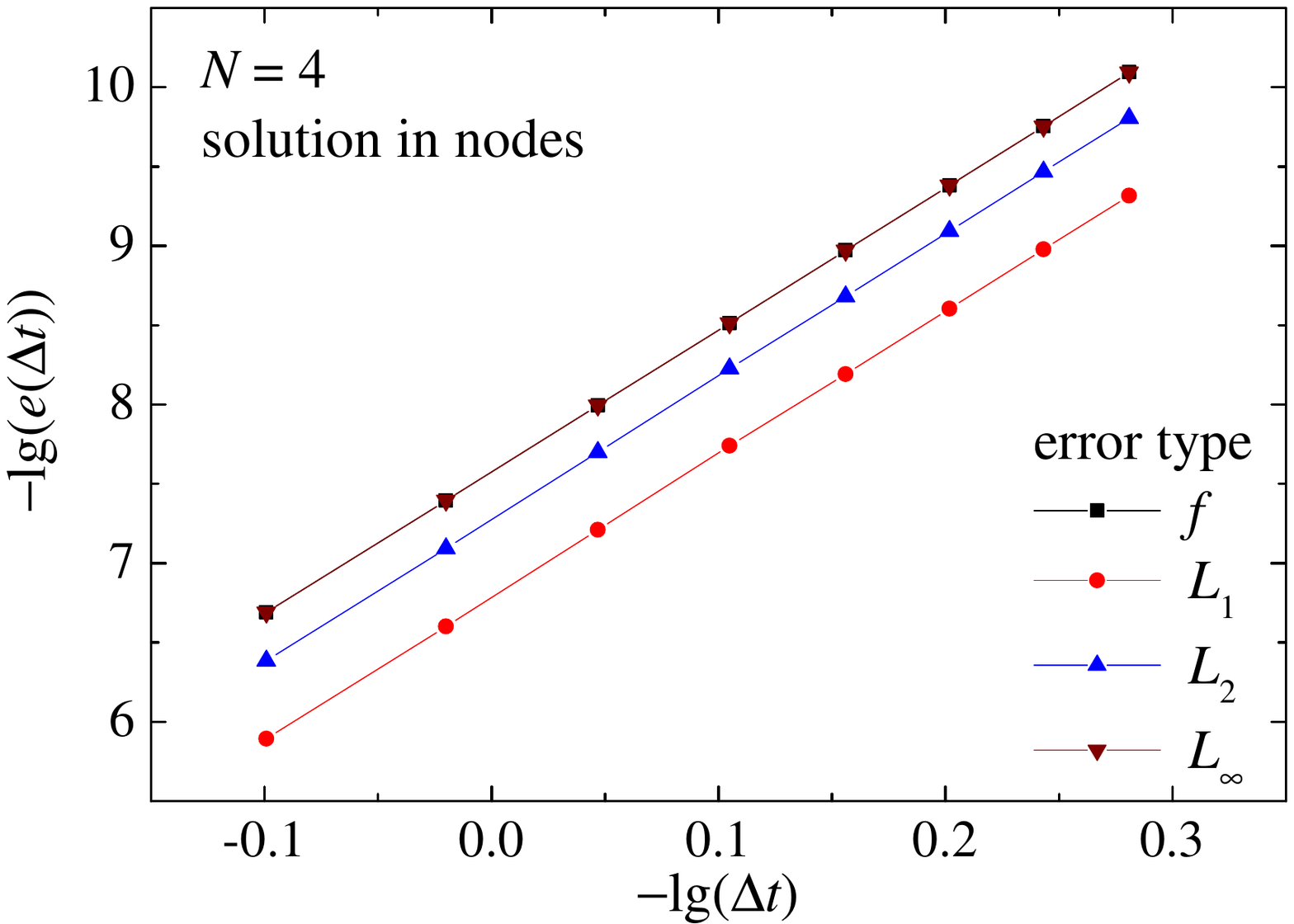}
\vspace{-8mm}\caption{\label{fig:harm_osc:f2}}
\end{subfigure}
\begin{subfigure}{0.24\textwidth}
\includegraphics[width=\textwidth]{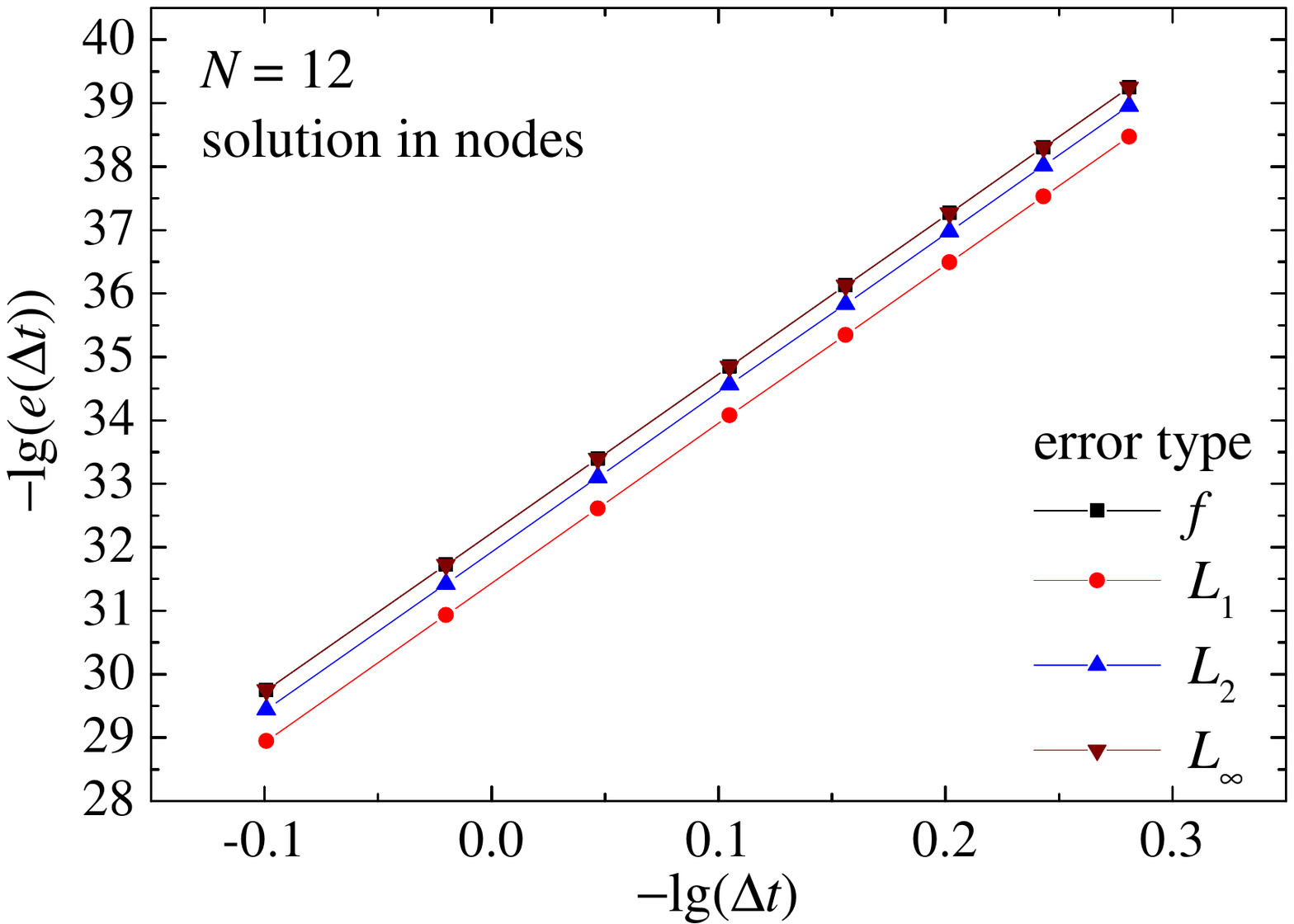}
\vspace{-8mm}\caption{\label{fig:harm_osc:f3}}
\end{subfigure}
\begin{subfigure}{0.24\textwidth}
\includegraphics[width=\textwidth]{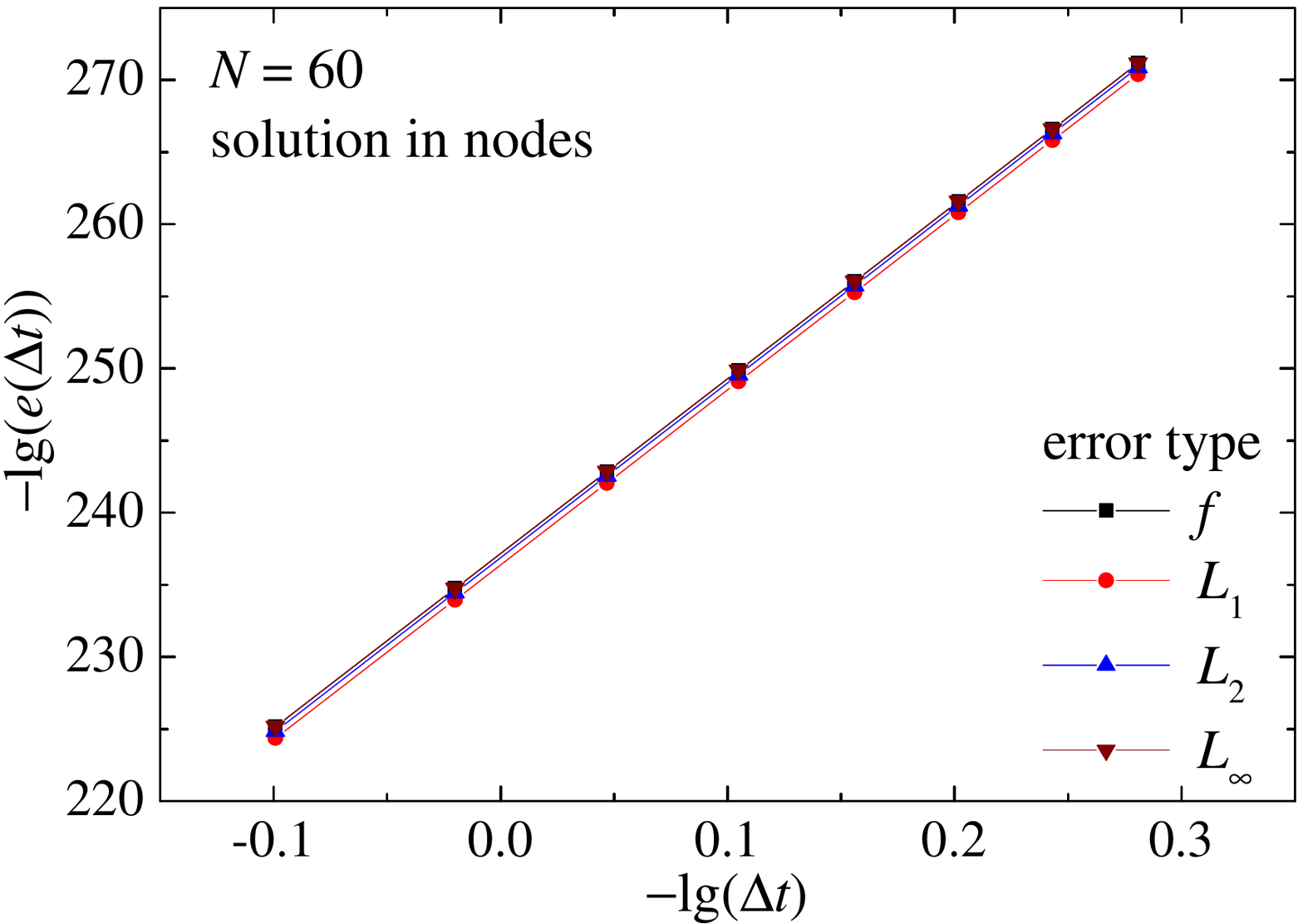}
\vspace{-8mm}\caption{\label{fig:harm_osc:f4}}
\end{subfigure}\\
\caption{%
Numerical solution of the system (\ref{eq:harm_osc_ode}). Comparison of the solution at nodes $\mathbf{u}_{n}$, the local solution $\mathbf{u}_{L}(t)$, the improved local solution $\mathbf{u}_{\rm IL}(t)$ and the exact solution $\mathbf{u}^{\rm ex}(t)$ (\ref{eq:harm_osc_sol_ex}) for components $u_{1} \equiv x$ (\subref{fig:harm_osc:a1}, \subref{fig:harm_osc:a2}, \subref{fig:harm_osc:a3}, \subref{fig:harm_osc:a4}) and $u_{2} \equiv \dot{x}$ (\subref{fig:harm_osc:b1}, \subref{fig:harm_osc:b2}, \subref{fig:harm_osc:b3}, \subref{fig:harm_osc:b4}), the errors $\varepsilon(t)$ (\subref{fig:harm_osc:c1}, \subref{fig:harm_osc:c2}, \subref{fig:harm_osc:c3}, \subref{fig:harm_osc:c4}), obtained using polynomials with degrees $N = 1$ (\subref{fig:harm_osc:a1}, \subref{fig:harm_osc:b1}, \subref{fig:harm_osc:c1}), $N = 4$ (\subref{fig:harm_osc:a2}, \subref{fig:harm_osc:b2}, \subref{fig:harm_osc:c2}), $N = 12$ (\subref{fig:harm_osc:a3}, \subref{fig:harm_osc:b3}, \subref{fig:harm_osc:c3}) and $N = 60$ (\subref{fig:harm_osc:a4}, \subref{fig:harm_osc:b4}, \subref{fig:harm_osc:c4}). Log-log plot of the dependence of the global error for the local solution $e^{l}$ (\subref{fig:harm_osc:d1}, \subref{fig:harm_osc:d2}, \subref{fig:harm_osc:d3}, \subref{fig:harm_osc:d4}), the improved local solution $e^{\rm imp}$ (\subref{fig:harm_osc:e1}, \subref{fig:harm_osc:e2}, \subref{fig:harm_osc:e3}, \subref{fig:harm_osc:e4}) and the solution at nodes $e^{n}$ (\subref{fig:harm_osc:f1}, \subref{fig:harm_osc:f2}, \subref{fig:harm_osc:f3}, \subref{fig:harm_osc:f4}) on the discretization step $\mathrm{\Delta}t$, obtained in the $f$-norm and norms $L_{1}$, $L_{2}$ and $L_{\infty}$, obtained using polynomials with degrees $N = 1$ (\subref{fig:harm_osc:d1}, \subref{fig:harm_osc:e1}, \subref{fig:harm_osc:f1}), $N = 4$ (\subref{fig:harm_osc:d2}, \subref{fig:harm_osc:e2}, \subref{fig:harm_osc:f2}), $N = 12$ (\subref{fig:harm_osc:d3}, \subref{fig:harm_osc:e3}, \subref{fig:harm_osc:f3}) and $N = 60$ (\subref{fig:harm_osc:d4}, \subref{fig:harm_osc:e4}, \subref{fig:harm_osc:f4}).
}
\label{fig:harm_osc}
\end{figure}

The third example of applying the ADER-DG numerical method with a local DG predictor to solving the initial value problem for the ODE system (\ref{eq:ivp_ode_diff_src}) is the problem of describing the motion of a linear one-dimensional harmonic oscillator:
\begin{equation}\label{eq:harm_osc_ode}
\ddot{x} + x = 0,\quad
x(0) = 1,\quad \dot{x}(0) = 0,\quad
t \in [0,\, 4\pi],
\end{equation}
for which the solution vectors in the original notation of the ODE system (\ref{eq:ivp_ode_diff_src}), with $D = 2$, are chosen in the form $\mathbf{u}(t) = [x(t)\ \dot{x}(t)]^{T}$. The exact analytical solution is as follows:
\begin{equation}\label{eq:harm_osc_sol_ex}
x^{\rm ex}(t) = \cos(t),\quad
\dot{x}^{\rm ex}(t) = -\sin(t),
\end{equation}
This example corresponds to the demonstration example (\ref{eq:demo_ode}) presented above for representing the local solution $\mathbf{u}_{L}$ and the improved solution $\mathbf{u}_{\rm IL}$; however, in this Subsection, this example is studied from a quantitative perspective, in terms of the accuracy and convergence of the numerical solutions $\mathbf{u}_{L}$, $\mathbf{u}_{\rm IL}$, $\mathbf{u}_{n}$. The obtained results are presented in Fig.~\ref{fig:harm_osc} and Table~\ref{tab:conv_orders_harm_osc}.

\begin{table*}[h!]
\centering
\normalsize
\caption{%
Convergence orders $p_{f}$, $p_{L_{1}}$, $p_{L_{2}}$, $p_{L_{\infty}}$, calculated in $f$-norm and norms $L_{1}$, $L_{2}$, $L_{\infty}$, of the numerical solution of the ADER-DG method for the problem (\ref{eq:harm_osc_ode}); $N$ is the degree of the basis polynomials $\varphi_{p}$. Orders $p^{n}$ are calculated for \textit{the numerical solution at the nodes} $\mathbf{u}_{n}$; orders $p^{\rm imp}$ --- for \textit{the improved local solution} $\mathbf{u}_{\rm IL}$; orders $p^{l}$ --- for \textit{the local solution} $\mathbf{u}_{L}$. The theoretical values of convergence order $p_{\rm th.}^{n} = 2N+1$, $p_{\rm th.}^{l} = N+1$ and $p^{\rm imp}_{\rm th.} = N+2$ are presented for comparison.
}
\label{tab:conv_orders_harm_osc}
\setlength{\tabcolsep}{3.5pt}
\begin{tabular}{@{}|l|llll|c|lll|c|lll|c|@{}}
\toprule
$N$ & $p^{n}_{f}$ &
$p^{n}_{L_{1}}$ & $p^{n}_{L_{2}}$ & $p^{n}_{L_{\infty}}$ & $p^{n}_{\rm th.}$ &
$p^{l}_{L_{1}}$ & $p^{l}_{L_{2}}$ & $p^{l}_{L_{\infty}}$ & $p^{l}_{\rm th.}$ &
$p^{\rm imp}_{L_{1}}$ & $p^{\rm imp}_{L_{2}}$ & $p^{\rm imp}_{L_{\infty}}$ & $p^{\rm imp}_{\rm th.}$\\
\midrule
$1$ & $2.70$ & $2.83$ & $2.79$ & $2.70$ & $3$ & $2.42$ & $2.42$ & $2.32$ & $2$ & $2.79$ & $2.78$ & $2.73$ & $3$\\
$2$ & $4.91$ & $4.97$ & $4.96$ & $4.91$ & $5$ & $3.12$ & $3.03$ & $2.95$ & $3$ & $4.56$ & $4.56$ & $4.50$ & $4$\\
$3$ & $6.94$ & $7.00$ & $6.98$ & $6.94$ & $7$ & $4.02$ & $4.00$ & $3.99$ & $4$ & $5.10$ & $5.05$ & $5.05$ & $5$\\
$4$ & $8.96$ & $9.01$ & $9.00$ & $8.96$ & $9$ & $5.01$ & $5.00$ & $4.99$ & $5$ & $6.02$ & $6.01$ & $6.00$ & $6$\\
$5$ & $11.0$ & $11.0$ & $11.0$ & $11.0$ & $11$ & $6.01$ & $6.00$ & $6.00$ & $6$ & $7.02$ & $7.01$ & $7.00$ & $7$\\
$6$ & $13.0$ & $13.0$ & $13.0$ & $13.0$ & $13$ & $7.01$ & $7.00$ & $7.00$ & $7$ & $8.01$ & $8.01$ & $8.01$ & $8$\\
$7$ & $15.0$ & $15.0$ & $15.0$ & $15.0$ & $15$ & $8.01$ & $8.00$ & $8.00$ & $8$ & $9.01$ & $9.01$ & $9.01$ & $9$\\
$8$ & $17.0$ & $17.0$ & $17.0$ & $17.0$ & $17$ & $9.01$ & $9.01$ & $9.00$ & $9$ & $10.0$ & $10.0$ & $10.0$ & $10$\\
$9$ & $19.0$ & $19.0$ & $19.0$ & $19.0$ & $19$ & $10.0$ & $10.0$ & $10.0$ & $10$ & $11.0$ & $11.0$ & $11.0$ & $11$\\
$10$ & $21.0$ & $21.0$ & $21.0$ & $21.0$ & $21$ & $11.0$ & $11.0$ & $11.0$ & $11$ & $12.0$ & $12.0$ & $12.0$ & $12$\\
\midrule
$11$ & $23.0$ & $23.0$ & $23.0$ & $23.0$ & $23$ & $12.0$ & $12.0$ & $12.0$ & $12$ & $13.0$ & $13.0$ & $13.0$ & $13$\\
$12$ & $25.0$ & $25.0$ & $25.0$ & $25.0$ & $25$ & $13.0$ & $13.0$ & $13.0$ & $13$ & $14.0$ & $14.0$ & $14.0$ & $14$\\
$13$ & $27.0$ & $27.0$ & $27.0$ & $27.0$ & $27$ & $14.0$ & $14.0$ & $14.0$ & $14$ & $15.0$ & $15.0$ & $15.0$ & $15$\\
$14$ & $29.0$ & $29.1$ & $29.0$ & $29.0$ & $29$ & $15.0$ & $15.0$ & $15.0$ & $15$ & $16.0$ & $16.0$ & $16.0$ & $16$\\
$15$ & $31.0$ & $31.1$ & $31.0$ & $31.0$ & $31$ & $16.0$ & $16.0$ & $16.0$ & $16$ & $17.0$ & $17.0$ & $17.0$ & $17$\\
$16$ & $33.0$ & $33.1$ & $33.0$ & $33.0$ & $33$ & $17.0$ & $17.0$ & $17.0$ & $17$ & $18.0$ & $18.0$ & $18.0$ & $18$\\
$17$ & $35.0$ & $35.1$ & $35.0$ & $35.0$ & $35$ & $18.0$ & $18.0$ & $18.0$ & $18$ & $19.0$ & $19.0$ & $19.0$ & $19$\\
$18$ & $37.0$ & $37.1$ & $37.0$ & $37.0$ & $37$ & $19.0$ & $19.0$ & $19.0$ & $19$ & $20.0$ & $20.0$ & $20.0$ & $20$\\
$19$ & $39.0$ & $39.1$ & $39.0$ & $39.0$ & $39$ & $20.0$ & $20.0$ & $20.0$ & $20$ & $21.0$ & $21.0$ & $21.0$ & $21$\\
$20$ & $41.0$ & $41.1$ & $41.0$ & $41.0$ & $41$ & $21.0$ & $21.0$ & $21.0$ & $21$ & $22.0$ & $22.0$ & $22.0$ & $22$\\
\midrule
$21$ & $43.0$ & $43.1$ & $43.0$ & $43.0$ & $43$ & $22.0$ & $22.0$ & $22.0$ & $22$ & $23.0$ & $23.0$ & $23.0$ & $23$\\
$22$ & $45.0$ & $45.1$ & $45.0$ & $45.0$ & $45$ & $23.0$ & $23.0$ & $23.0$ & $23$ & $24.0$ & $24.0$ & $24.0$ & $24$\\
$23$ & $47.0$ & $47.1$ & $47.0$ & $47.0$ & $47$ & $24.0$ & $24.0$ & $24.0$ & $24$ & $25.0$ & $25.0$ & $25.0$ & $25$\\
$24$ & $49.0$ & $49.1$ & $49.0$ & $49.0$ & $49$ & $25.0$ & $25.0$ & $25.0$ & $25$ & $26.0$ & $26.0$ & $26.0$ & $26$\\
$25$ & $51.0$ & $51.1$ & $51.0$ & $51.0$ & $51$ & $26.0$ & $26.0$ & $26.0$ & $26$ & $27.0$ & $27.0$ & $27.0$ & $27$\\
$26$ & $53.0$ & $53.1$ & $53.0$ & $53.0$ & $53$ & $27.0$ & $27.0$ & $27.0$ & $27$ & $28.0$ & $28.0$ & $28.0$ & $28$\\
$27$ & $55.0$ & $55.1$ & $55.0$ & $55.0$ & $55$ & $28.0$ & $28.0$ & $28.0$ & $28$ & $29.0$ & $29.0$ & $29.0$ & $29$\\
$28$ & $57.0$ & $57.1$ & $57.0$ & $57.0$ & $57$ & $29.0$ & $29.0$ & $29.0$ & $29$ & $30.0$ & $30.0$ & $30.0$ & $30$\\
$29$ & $59.0$ & $59.1$ & $59.0$ & $59.0$ & $59$ & $30.0$ & $30.0$ & $30.0$ & $30$ & $31.0$ & $31.0$ & $31.0$ & $31$\\
$30$ & $61.0$ & $61.1$ & $61.0$ & $61.0$ & $61$ & $31.0$ & $31.0$ & $31.0$ & $31$ & $32.0$ & $32.0$ & $32.0$ & $32$\\
\midrule
$35$ & $71.0$ & $71.1$ & $71.1$ & $71.0$ & $71$ & $36.0$ & $36.0$ & $36.0$ & $36$ & $37.0$ & $37.0$ & $37.0$ & $37$\\
$40$ & $81.0$ & $81.1$ & $81.1$ & $81.0$ & $81$ & $41.0$ & $41.0$ & $41.0$ & $41$ & $42.0$ & $42.0$ & $42.0$ & $42$\\
$45$ & $91.0$ & $91.1$ & $91.1$ & $91.0$ & $91$ & $46.0$ & $46.0$ & $46.0$ & $46$ & $47.0$ & $47.0$ & $47.0$ & $47$\\
$50$ & $101.0$ & $101.1$ & $101.1$ & $101.0$ & $101$ & $51.0$ & $51.0$ & $51.0$ & $51$ & $52.0$ & $52.0$ & $52.0$ & $52$\\
$55$ & $111.0$ & $111.1$ & $111.1$ & $111.0$ & $111$ & $56.0$ & $56.0$ & $56.0$ & $56$ & $57.0$ & $57.0$ & $57.0$ & $57$\\
$60$ & $121.0$ & $121.1$ & $121.1$ & $121.0$ & $121$ & $61.0$ & $61.0$ & $61.0$ & $61$ & $62.0$ & $62.0$ & $62.0$ & $62$\\
\bottomrule
\end{tabular}
\end{table*}

The dependencies of the numerical solutions $\mathbf{u}_{L}$, $\mathbf{u}_{\rm IL}$, $\mathbf{u}_{n}$ and the exact analytical solution $\mathbf{u}^{\rm ex}$, the dependencies of the local error $\varepsilon$ (\ref{eq:eps_local_def}) of the numerical solutions, and the dependence of the global error $e$ (\ref{eq:eps_un_global_def}), (\ref{eq:eps_ul_global_def}) of the numerical solutions on the discretization step ${\Delta t}$, for polynomial degrees $N = 1$, $4$, $12$ and $60$, are shown in Fig.~\ref{fig:harm_osc}. The obtained dependencies of the local numerical solution $\mathbf{u}_{L}$ and the improved local numerical solution $\mathbf{u}_{\rm IL}$ for polynomials degree $N = 1$ exhibit features similar to those presented in the demonstration example (\ref{eq:demo_ode}), as can be seen from a comparison of Fig.~\ref{fig:harm_osc} (\subref{fig:harm_osc:a1}, \subref{fig:harm_osc:b1}) with Fig.~\ref{fig:demo_nodes_8} (\subref{fig:demo_nodes_8:b1}, \subref{fig:demo_nodes_8:b2}, \subref{fig:demo_nodes_8:b3}, \subref{fig:demo_nodes_8:b4}) (however, due to the smaller discretization step ${\Delta t}$, the numerical dissipation effect is significantly less than in Fig.~\ref{fig:demo_nodes_8} (\subref{fig:demo_nodes_8:a1}, \subref{fig:demo_nodes_8:a2}, \subref{fig:demo_nodes_8:a3}, \subref{fig:demo_nodes_8:a4})). A comparison of the obtained dependencies of the numerical solutions $\mathbf{u}_{L}(t)$, $\mathbf{u}_{\rm IL}(t)$, $\mathbf{u}_{n}$ for polynomial degrees $N = 4$, $12$ and $60$ with the exact analytical solution $\mathbf{u}^{\rm ex}(t)$, presented in Fig.~\ref{fig:harm_osc} (\subref{fig:harm_osc:a2}, \subref{fig:harm_osc:a3}, \subref{fig:harm_osc:a4}) for the component $u_{1}$ and in Fig.~\ref{fig:harm_osc} (\subref{fig:harm_osc:b2}, \subref{fig:harm_osc:b3}, \subref{fig:harm_osc:b4}) for the component $u_{2}$, demonstrates a high-quality agreement.

\begin{figure}[h!]
\captionsetup[subfigure]{%
	position=bottom,
	font+=smaller,
	textfont=normalfont,
	singlelinecheck=off,
	justification=raggedright
}
\centering
\begin{subfigure}{0.24\textwidth}
\includegraphics[width=\textwidth]{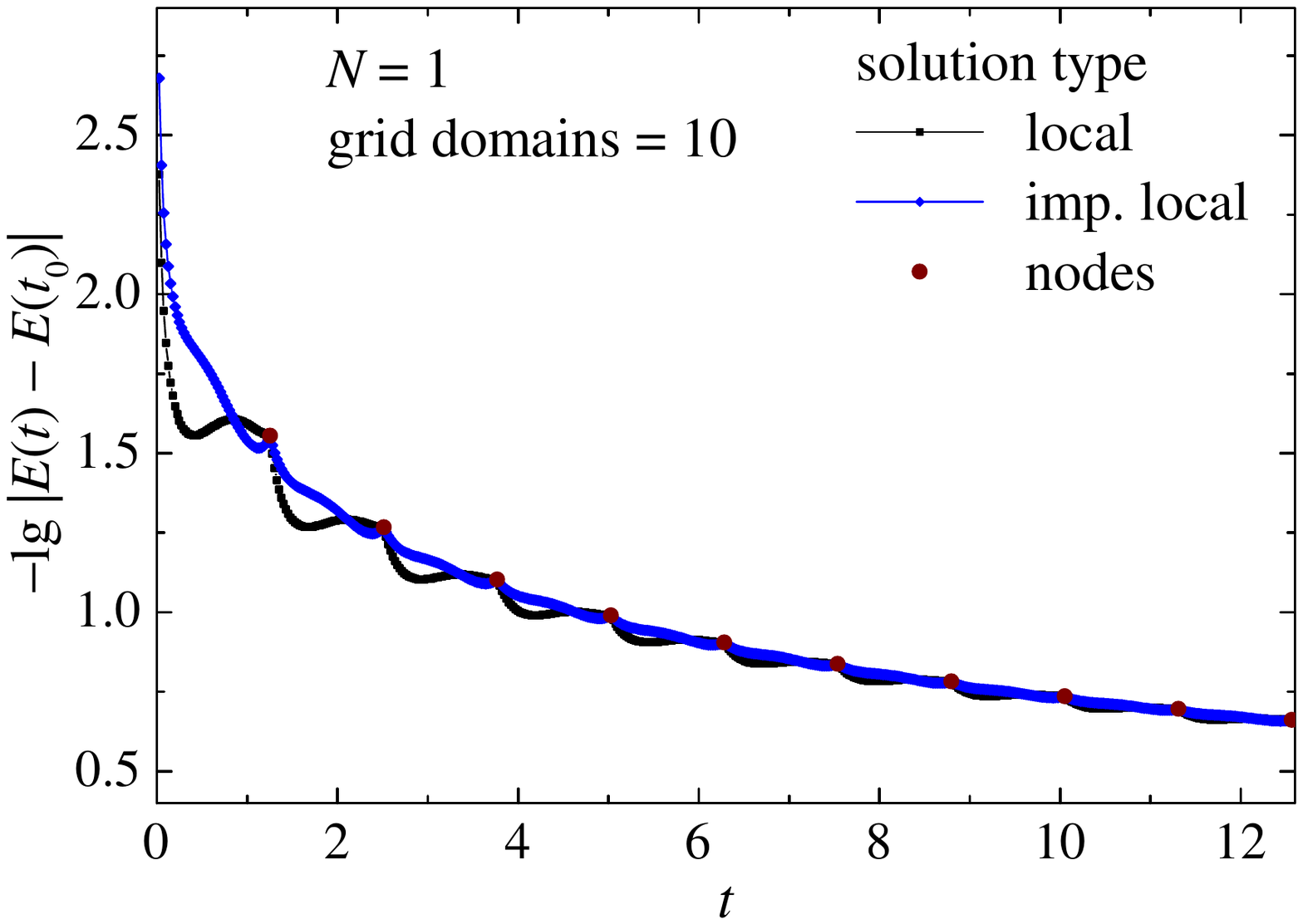}
\vspace{-8mm}\caption{\label{fig:econs_harm_osc:a1}}
\end{subfigure}
\begin{subfigure}{0.24\textwidth}
\includegraphics[width=\textwidth]{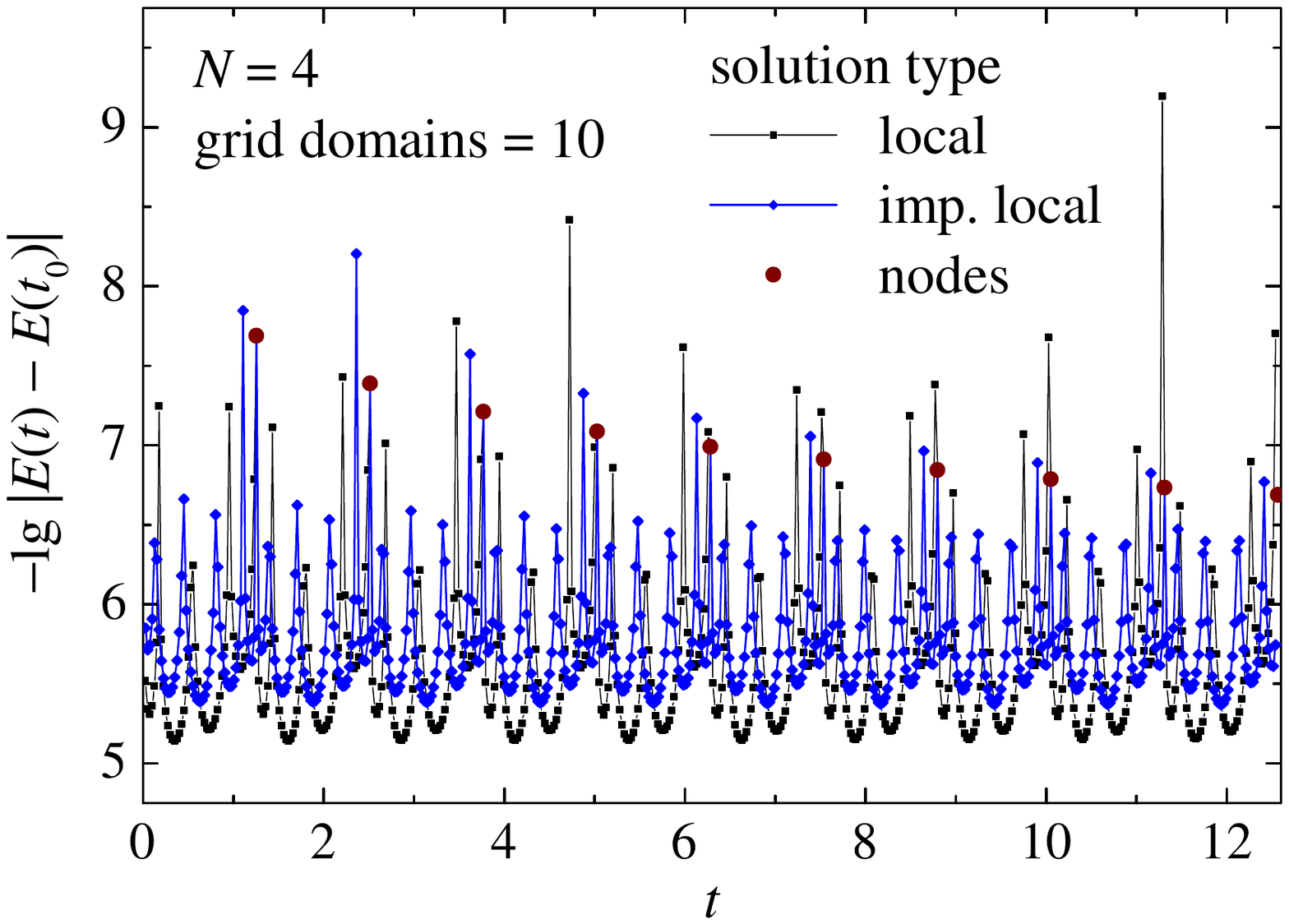}
\vspace{-8mm}\caption{\label{fig:econs_harm_osc:a2}}
\end{subfigure}
\begin{subfigure}{0.24\textwidth}
\includegraphics[width=\textwidth]{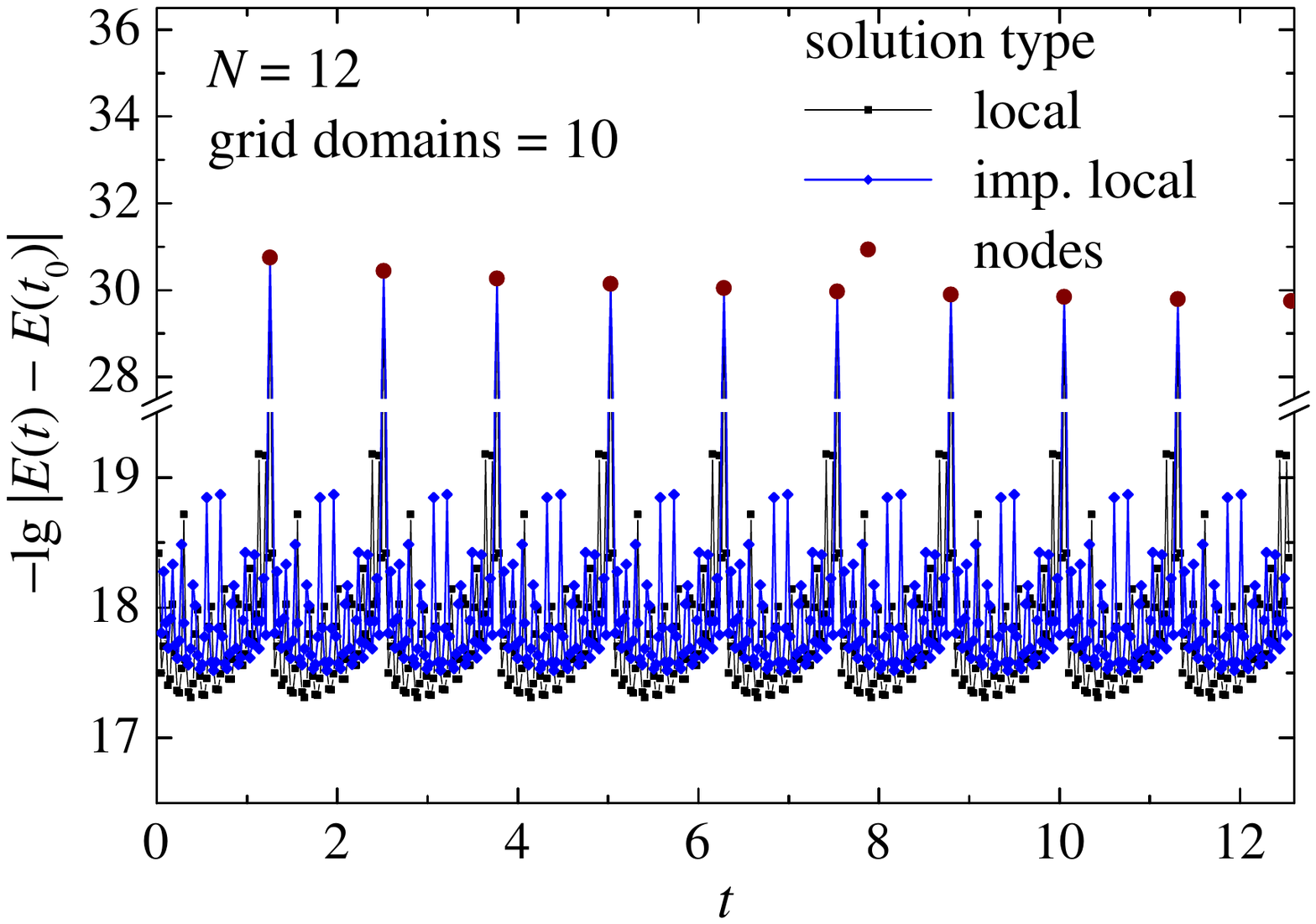}
\vspace{-8mm}\caption{\label{fig:econs_harm_osc:a3}}
\end{subfigure}
\begin{subfigure}{0.24\textwidth}
\includegraphics[width=\textwidth]{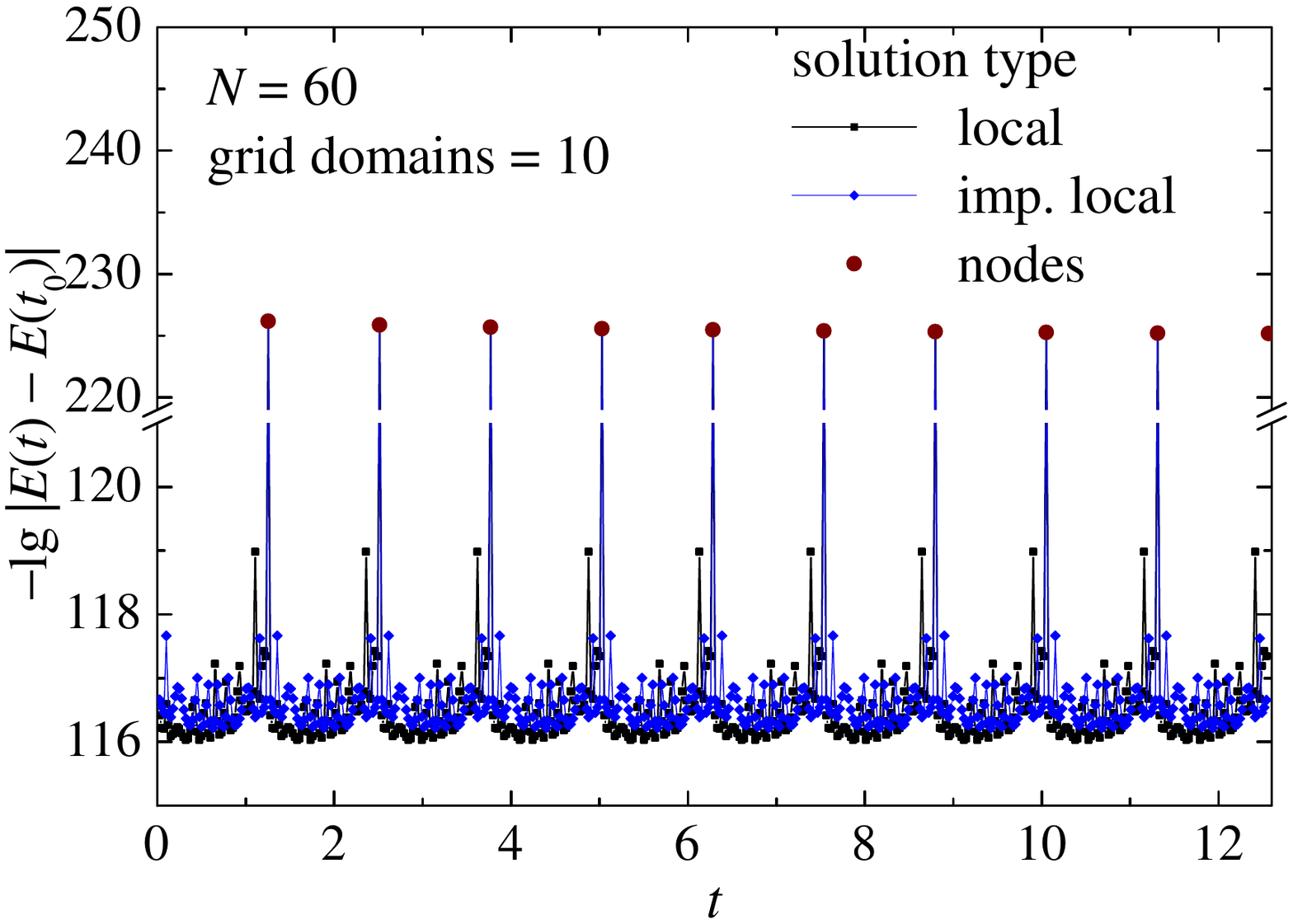}
\vspace{-8mm}\caption{\label{fig:econs_harm_osc:a4}}
\end{subfigure}
\caption{%
Dependence of the negative logarithm of the error $-\lg|E(t)-E(t_{0})|$ of the energy conservation law of the numerical solution of the system (\ref{eq:harm_osc_ode}) on the argument $t$ for solution at nodes $\mathbf{u}_{n}$, the local solution $\mathbf{u}_{L}(t)$ and the improved local solution $\mathbf{u}_{\rm IL}(t)$, obtained using polynomials with degrees $N = 1$~(\subref{fig:econs_harm_osc:a1}), $N = 4$~(\subref{fig:econs_harm_osc:a4}), $N = 12$~(\subref{fig:econs_harm_osc:a3}), $N = 60$~(\subref{fig:econs_harm_osc:a4}).
}
\label{fig:econs_harm_osc}
\end{figure}

The local error $\varepsilon$ (\ref{eq:eps_local_def}) of the numerical solutions shows that for polynomial degree $N = 1$, the errors of the local solution $\mathbf{u}_{L}$ and the improved local solution $\mathbf{u}_{\rm IL}$ are approximately comparable throughout the entire domain, but in the initial region, the improved local solution exhibits higher accuracy. With increasing polynomial degree $N$, in particular, in the cases of polynomial degrees $N = 4$, $12$, $60$, the local error $\varepsilon$ of the improved local solution $\mathbf{u}_{\rm IL}$ is significantly lower than the local error $\varepsilon$ of the local solution $\mathbf{u}_{L}$, with the differences reaching $2$--$3$ orders of magnitude. The local error $\varepsilon$ of the numerical solution $\mathbf{u}_{n}$ at the grid nodes $t_{n}$ is significantly smaller than the local errors $\varepsilon$ of the local solution $\mathbf{u}_{L}$ and the improved local solution $\mathbf{u}_{\rm IL}$, amounting to approximately $1$ orders of magnitude in the case of polynomial degree $N = 1$, approximately $2$--$4$ orders of magnitude in the case of polynomial degree $N = 4$, approximately $12$--$15$ orders of magnitude in the case of polynomial degree $N = 12$ and approximately $110$--$115$ orders of magnitude in the case of polynomial degree $N = 60$ (for this purpose, breaks in the graph along the vertical axis are inserted in Fig.~\ref{fig:harm_osc} (\subref{fig:harm_osc:c3}, \subref{fig:harm_osc:c4})).

The resulting dependencies of the global error $e$ (\ref{eq:eps_un_global_def}), (\ref{eq:eps_ul_global_def}) on the discretization step ${\Delta t}$, presented in Fig.~\ref{fig:harm_osc} (\subref{fig:harm_osc:d1}, \subref{fig:harm_osc:d2}, \subref{fig:harm_osc:d3}, \subref{fig:harm_osc:d4}, \subref{fig:harm_osc:e1}, \subref{fig:harm_osc:e2}, \subref{fig:harm_osc:e3}, \subref{fig:harm_osc:e4}, \subref{fig:harm_osc:f1}, \subref{fig:harm_osc:f2}, \subref{fig:harm_osc:f3}, \subref{fig:harm_osc:f4}), demonstrate a very high quality of linear approximation for all studied polynomial degrees $N$ (Fig.~\ref{fig:harm_osc} only shows results for polynomial degrees $N = 1$, $4$, $12$ and $60$). Based on the approximation of the obtained dependencies $e({\Delta t})$ in log-log scale by a linear function $\lg{e({\Delta t})} \propto p\cdot\lg{{\Delta t}}$, empirical convergence orders $p$ are calculated and presented in Table~\ref{tab:conv_orders_harm_osc} for all polynomial degrees $N = 1, \ldots, 30$ and polynomial degrees $N$ up to $60$ with a step of $5$. A comparison of the obtained empirical convergence orders $p$ with the expected theoretical values $p_{\rm th.}$ (\ref{eq:conv_ords_exp}) shows excellent agreement. It is particularly noticeable that the empirical convergence orders $p^{\rm imp}$ of the improved local numerical solution $\mathbf{u}_{\rm IL}$ are one unit higher than the empirical convergence orders $p^{l}$ of the local numerical solution $\mathbf{u}_{L}$.

\corrtext{The system of equations (\ref{eq:harm_osc_ode}) considered in this Example, similar to the previous Example in Section~\ref{sec:apps:lin_diss}, can be presented in a conservative form, which corresponds to the mechanical problem of obtaining the law of motion of a particle of unit mass in an external field with potential energy $U(x) = x^{2}/2$~\cite{LandauMechanics, GoldsteinMechanics}, which, together with the initial conditions (\ref{eq:harm_osc_ode}), shows that for the exact solution the energy conservation law $E(t) = \mathrm{const}$ of the following form holds
\begin{equation}\label{eq:econs_harm_osc}
E(t) = \frac{\dot{x}^{2}(t)}{2} + U(x) = \frac12\left[\dot{x}^{2}(t) + x^{2}(t)\right] = \mathrm{const} \equiv 
E(t_{0}) = \frac12\left[\dot{x}^{2}(0) + x^{2}(0)\right] = \frac12,
\end{equation}
which represents the integral of motion.}

\corrtext{Fig.~\ref{fig:econs_harm_osc} shows the dependence of the negative logarithm of the error $-\lg|E(t)-E(t_{0})|$ of the energy conservation law (\ref{eq:econs_harm_osc}) of the numerical solution of the system (\ref{eq:harm_osc_ode}) on the argument $t$ for solution at nodes $\mathbf{u}_{n}$, the local solution $\mathbf{u}_{L}(t)$ and the improved local solution $\mathbf{u}_{\rm IL}(t)$, obtained using polynomials with degrees $N = 1$~(\subref{fig:econs_harm_osc:a1}), $N = 4$~(\subref{fig:econs_harm_osc:a4}), $N = 12$~(\subref{fig:econs_harm_osc:a3}), $N = 60$~(\subref{fig:econs_harm_osc:a4}). The obtained results clearly demonstrate that the energy conservation law (\ref{eq:harm_osc_ode}) is not strictly satisfied in the numerical solution, which is due to the dissipative nature of the ADER-DG numerical method with local DG predictor~\cite{ader_dg_ode_jsc, ader_dg_ode_sinum, ader_improving_2024, ader_proofs_2025}. However, due to the possibility of achieving a high order $p$, even in the cases of degrees $N = 12$ and $60$, the error in the energy conservation law's fulfillment over the entire solution domain becomes smaller than the characteristic rounding error of double-precision floating-point numbers $\sim 10^{-15}$-$10^{-17}$. The presented results also clearly demonstrate that the accuracy of the fulfillment of the energy conservation law (\ref{eq:econs_harm_osc}) for the improved local solution $\mathbf{u}_{\rm IL}(t)$ is significantly higher than for the local solution $\mathbf{u}_{L}(t)$. Therefore, it can be concluded that, despite the dissipative nature of the ADER-DG numerical method with local DG predictor, a sufficiently high degree $N$ can be chosen such that the accuracy of the fulfillment of the energy conservation law will be at or below the characteristic error of representing real numbers as floating-point numbers. The results obtained are qualitatively consistent with the results presented previously in Example in Section~\ref{sec:apps:lin_diss}.}

The obtained results allowed to conclude that the ADER-DG numerical method with a local DG predictor provides a highly accurate numerical solution to the initial value problem for the ODE system (\ref{eq:harm_osc_ode}) presented in this example. The obtained results are in good agreement with the theory developed above. The improved local numerical solution $\mathbf{u}_{\rm IL}$ indeed demonstrates higher accuracy and a higher convergence order compared to the local numerical solution $\mathbf{u}_{L}$.

\subsection{Example 4: Nonlinear mathematical pendulum}
\label{sec:apps:pend}

\begin{figure}[h!]
\captionsetup[subfigure]{%
	position=bottom,
	font+=smaller,
	textfont=normalfont,
	singlelinecheck=off,
	justification=raggedright
}
\centering
\begin{subfigure}{0.24\textwidth}
\includegraphics[width=\textwidth]{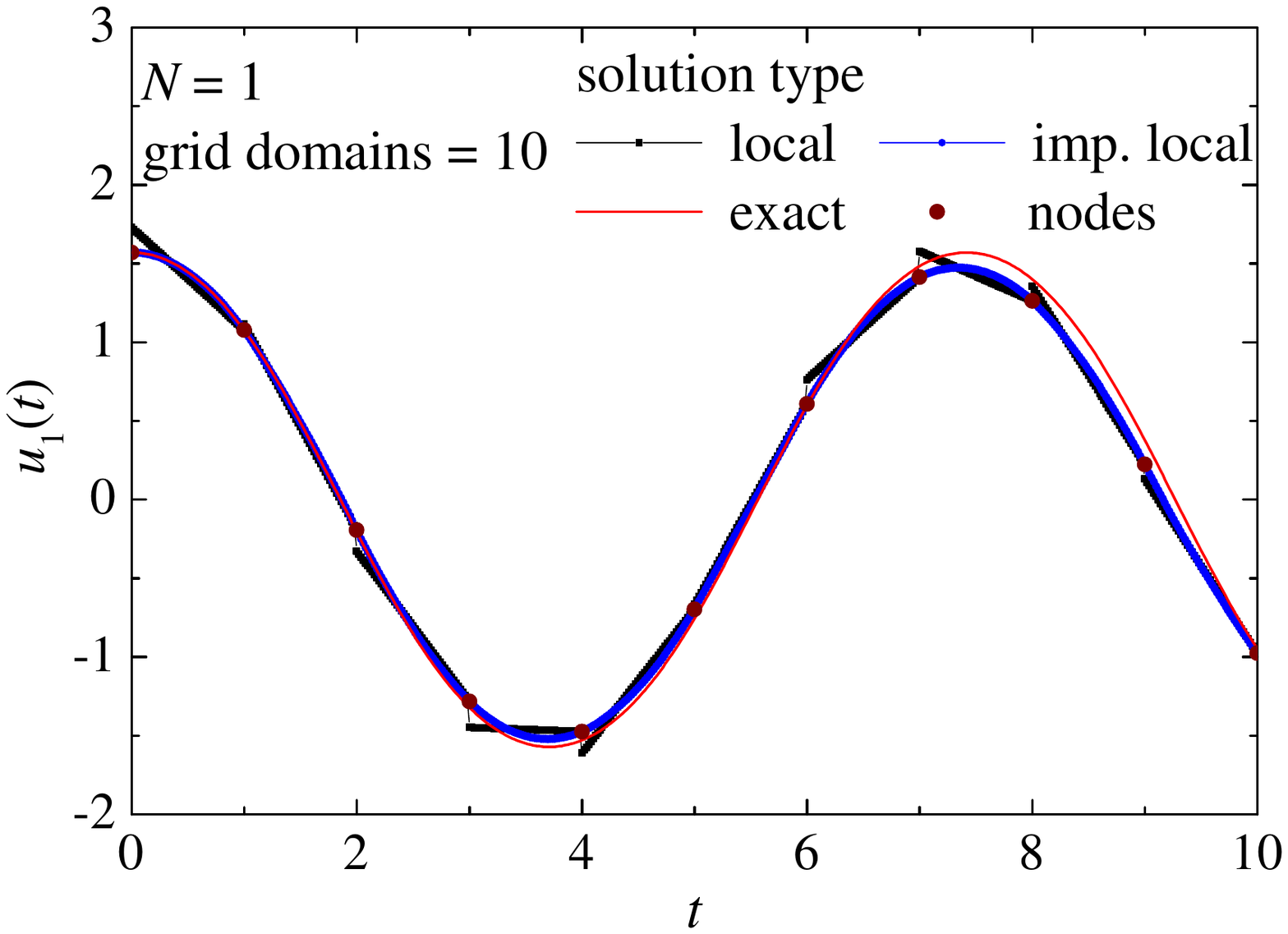}
\vspace{-8mm}\caption{\label{fig:pend:a1}}
\end{subfigure}
\begin{subfigure}{0.24\textwidth}
\includegraphics[width=\textwidth]{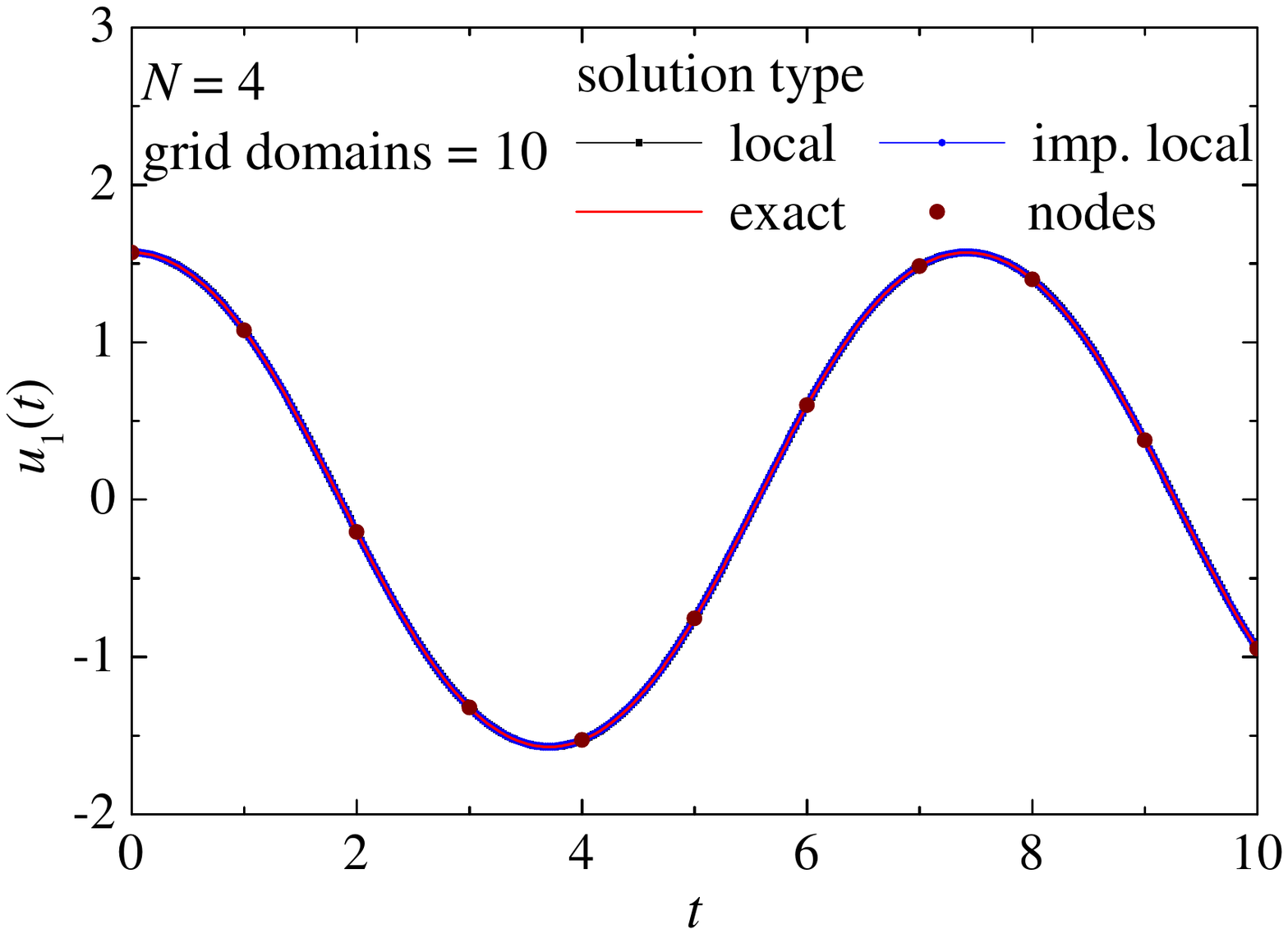}
\vspace{-8mm}\caption{\label{fig:pend:a2}}
\end{subfigure}
\begin{subfigure}{0.24\textwidth}
\includegraphics[width=\textwidth]{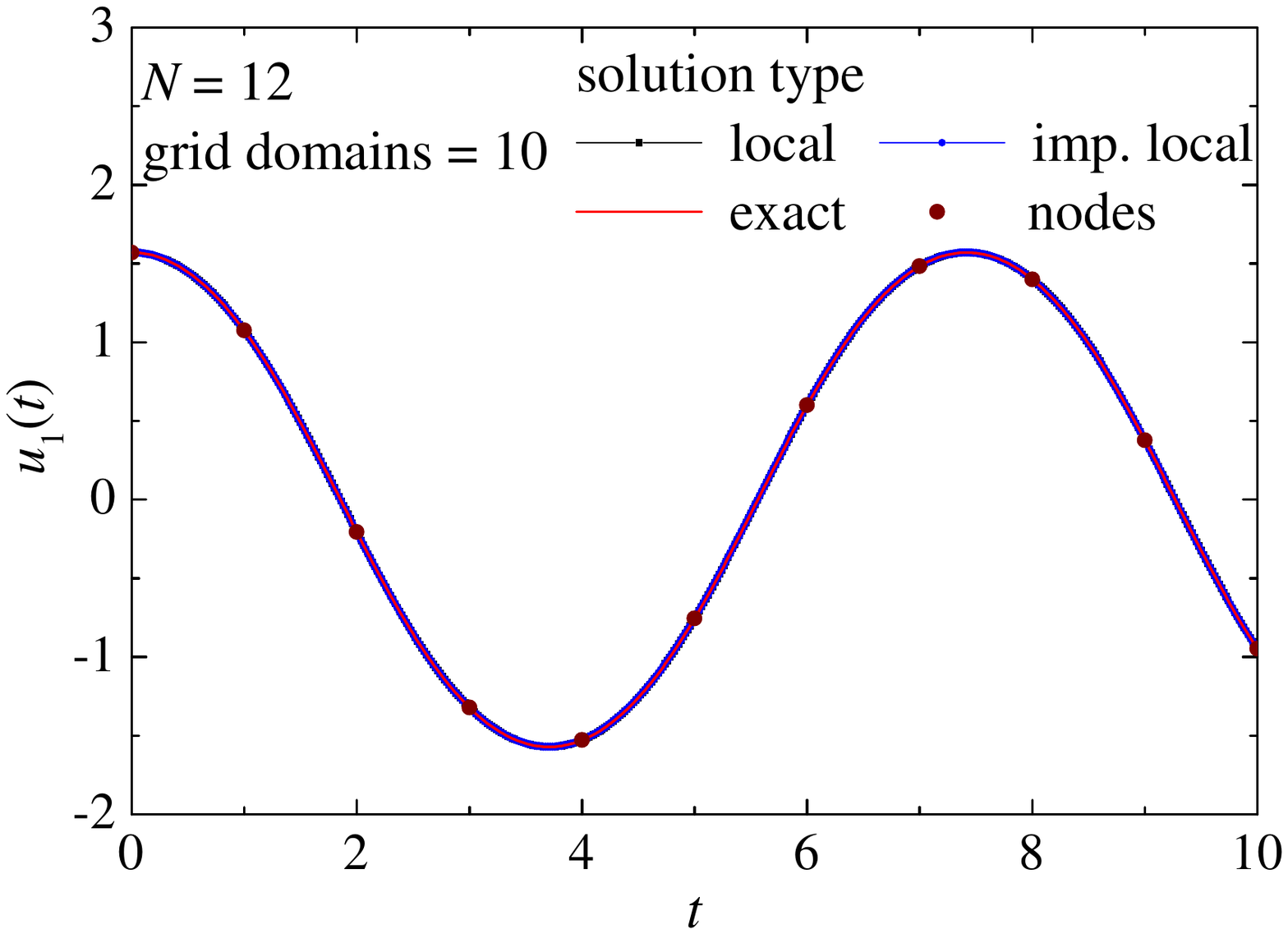}
\vspace{-8mm}\caption{\label{fig:pend:a3}}
\end{subfigure}
\begin{subfigure}{0.24\textwidth}
\includegraphics[width=\textwidth]{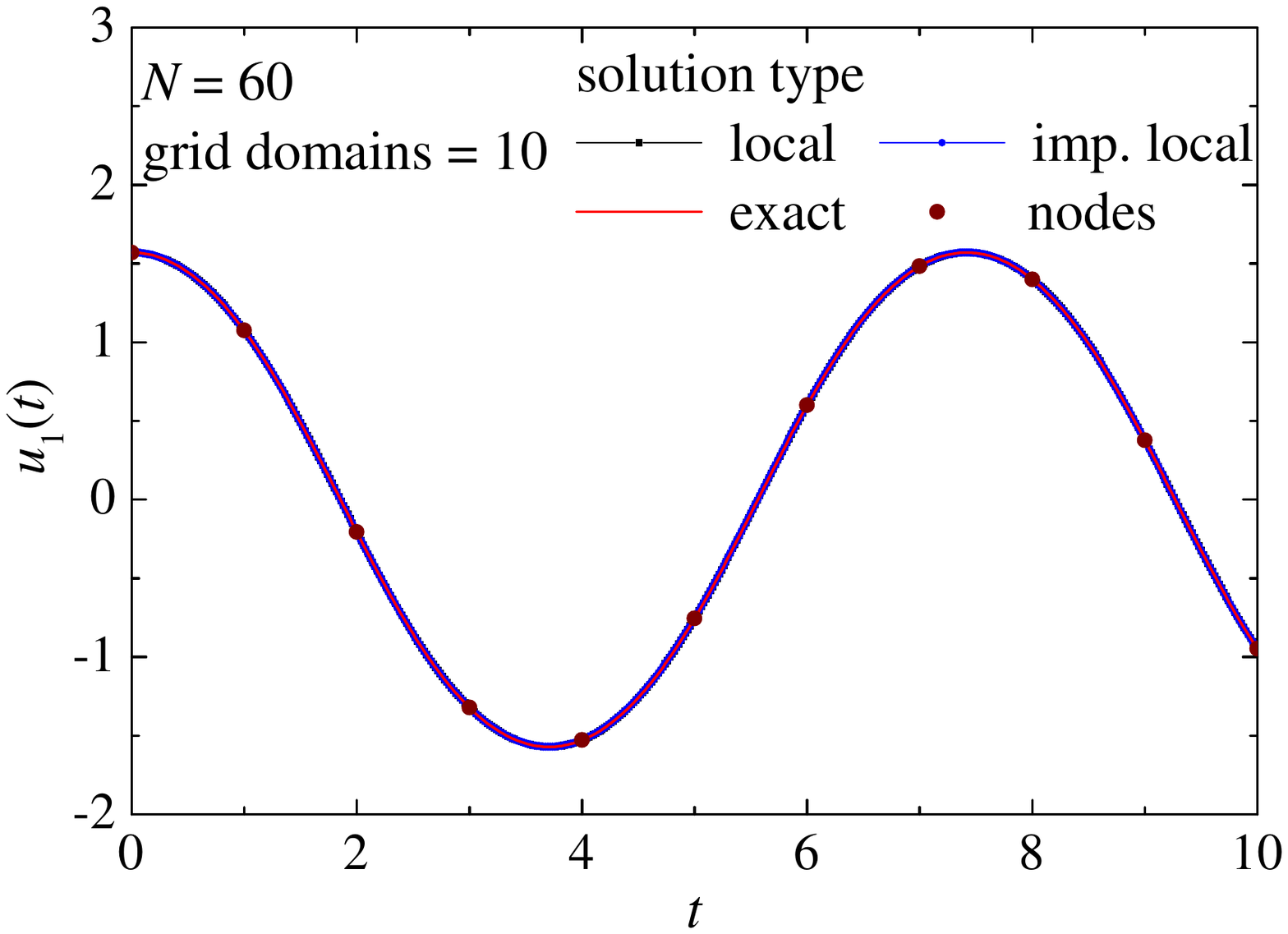}
\vspace{-8mm}\caption{\label{fig:pend:a4}}
\end{subfigure}\\
\begin{subfigure}{0.24\textwidth}
\includegraphics[width=\textwidth]{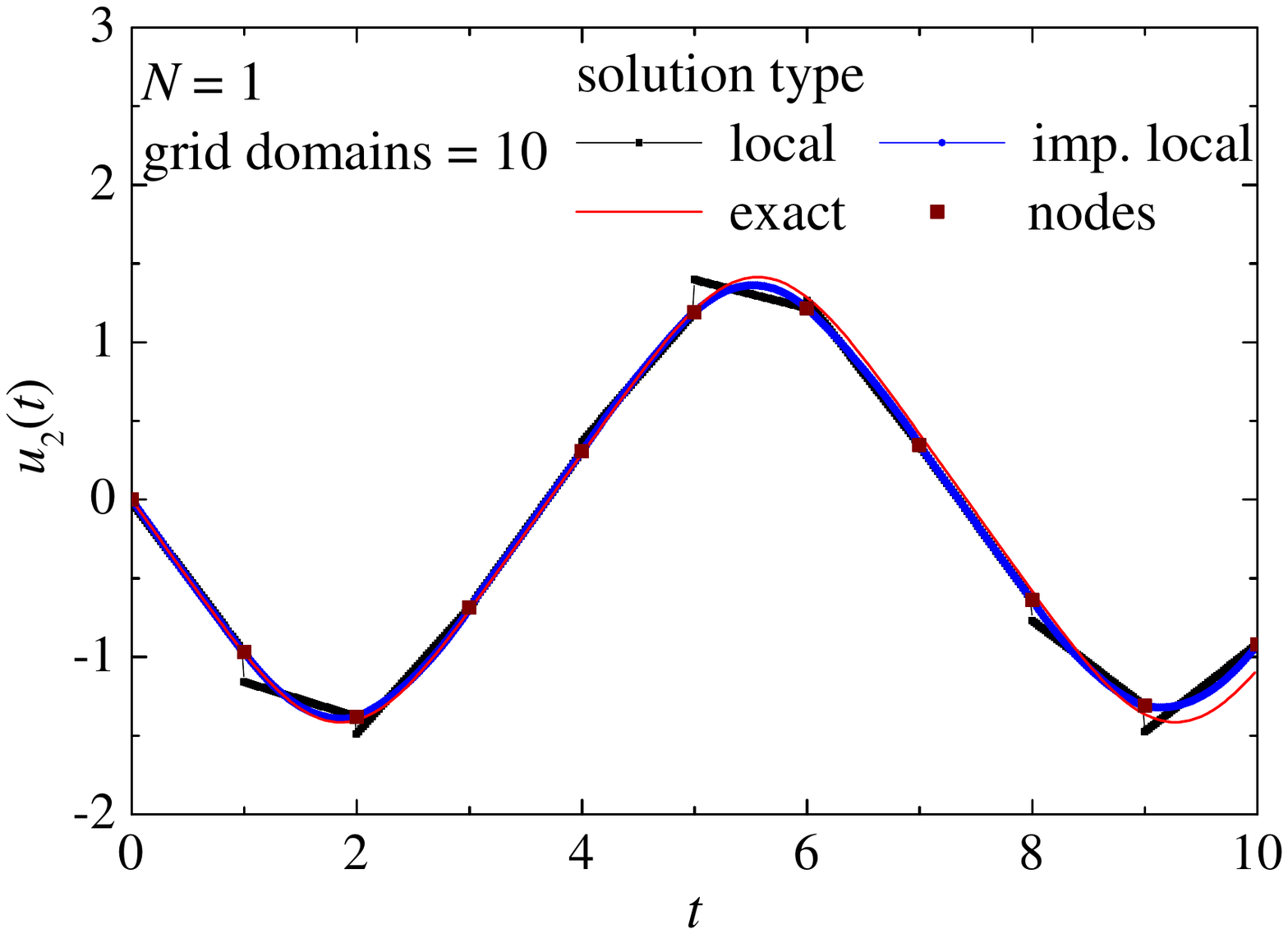}
\vspace{-8mm}\caption{\label{fig:pend:b1}}
\end{subfigure}
\begin{subfigure}{0.24\textwidth}
\includegraphics[width=\textwidth]{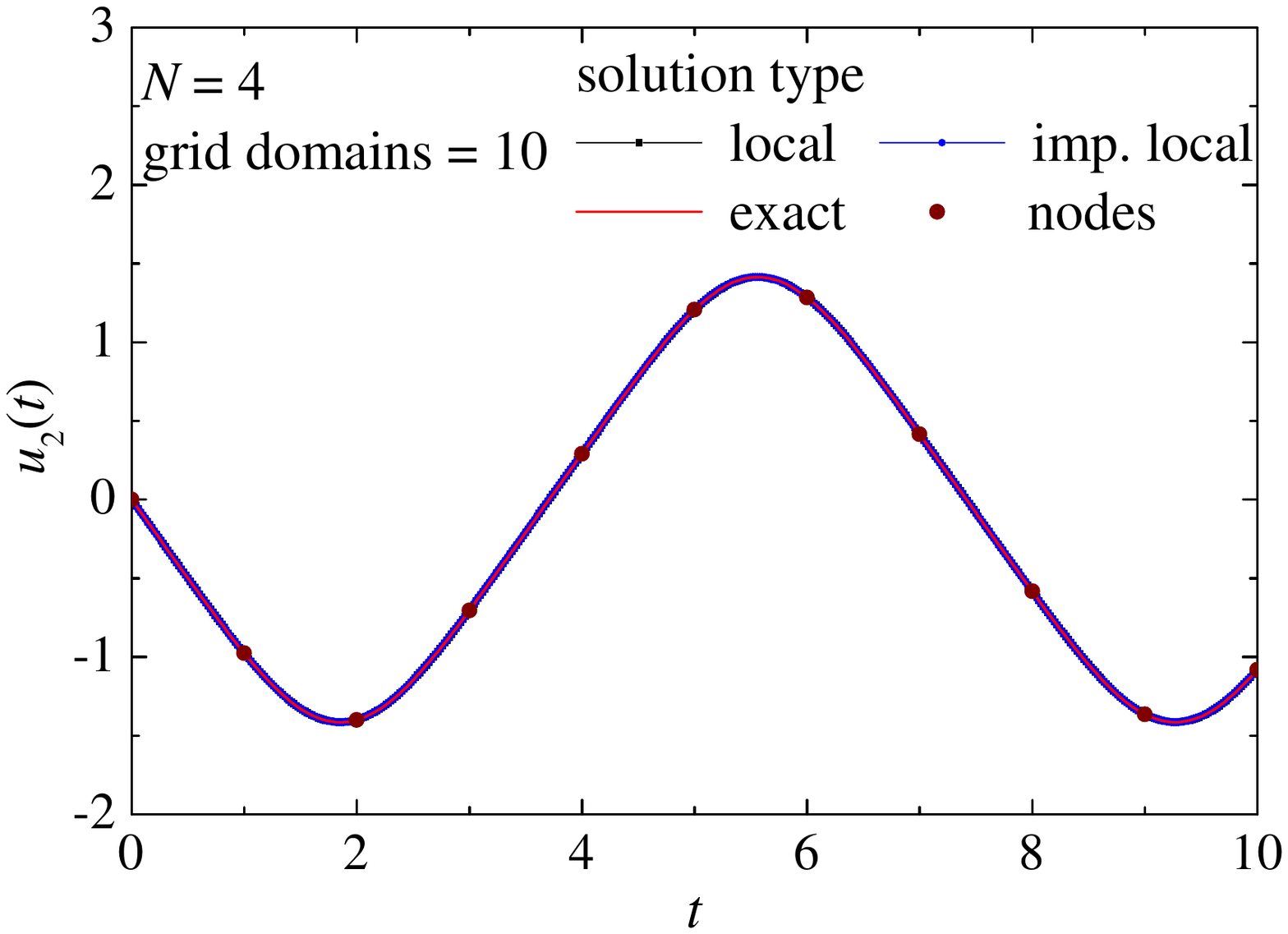}
\vspace{-8mm}\caption{\label{fig:pend:b2}}
\end{subfigure}
\begin{subfigure}{0.24\textwidth}
\includegraphics[width=\textwidth]{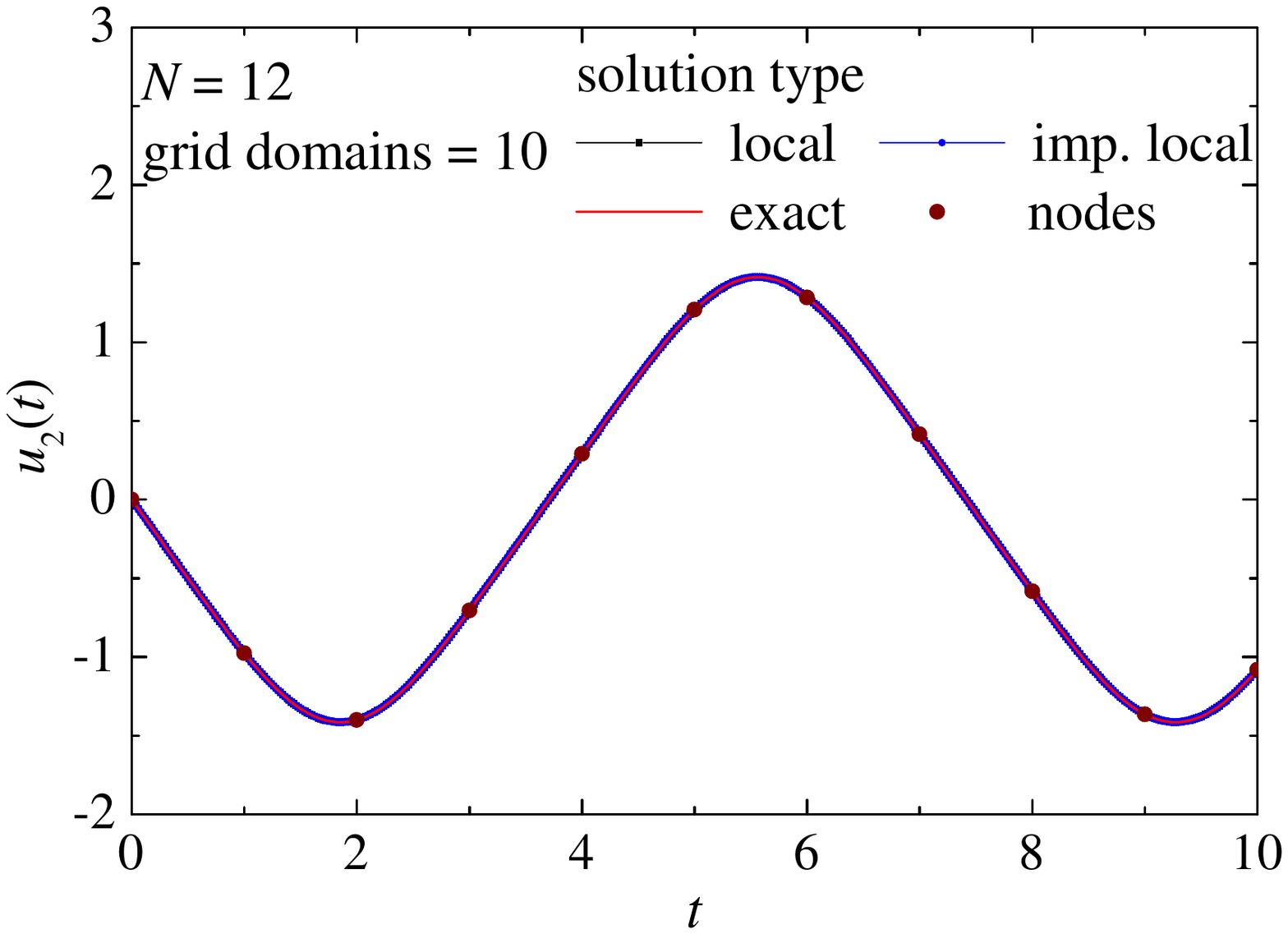}
\vspace{-8mm}\caption{\label{fig:pend:b3}}
\end{subfigure}
\begin{subfigure}{0.24\textwidth}
\includegraphics[width=\textwidth]{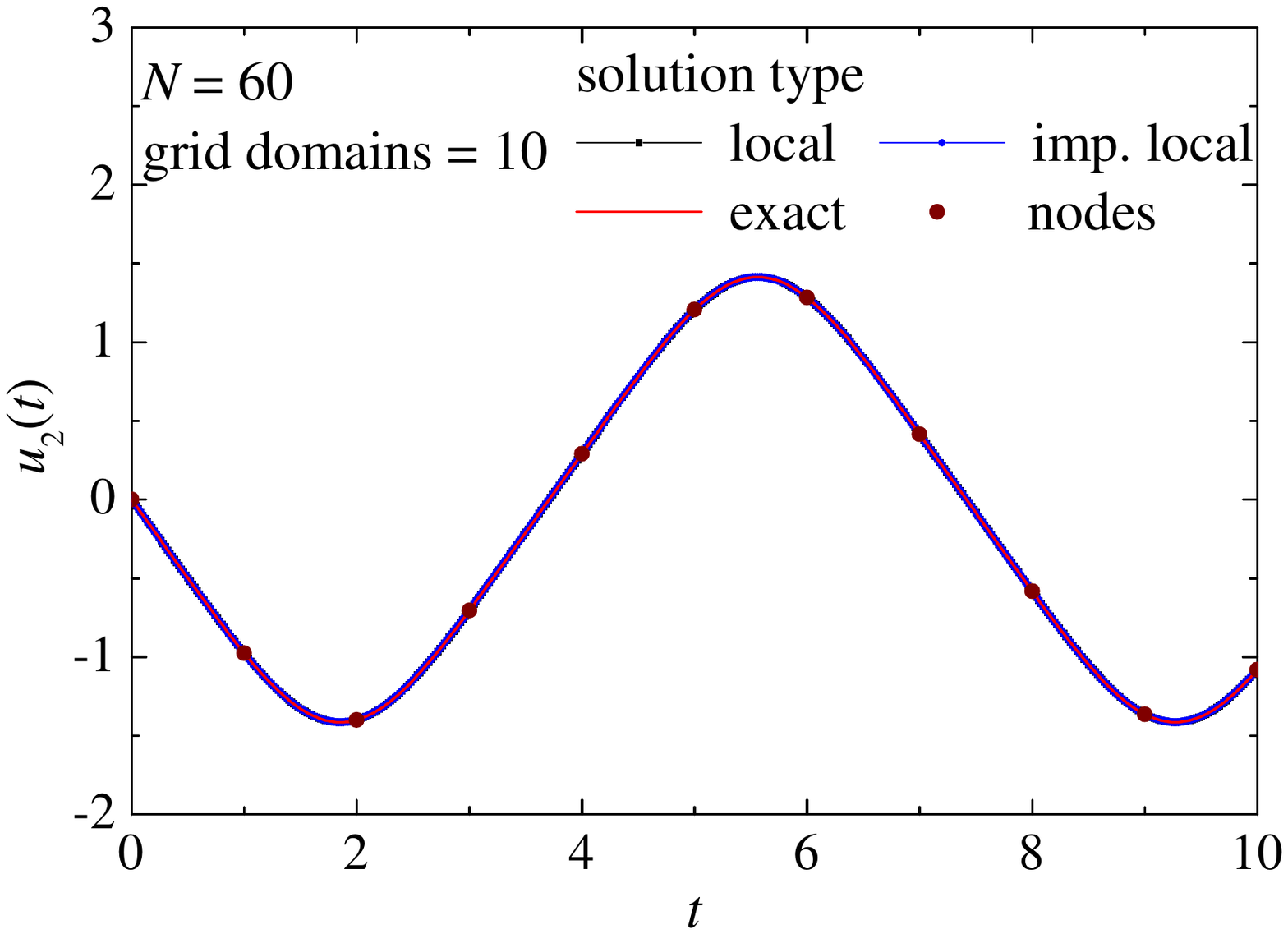}
\vspace{-8mm}\caption{\label{fig:pend:b4}}
\end{subfigure}\\
\begin{subfigure}{0.24\textwidth}
\includegraphics[width=\textwidth]{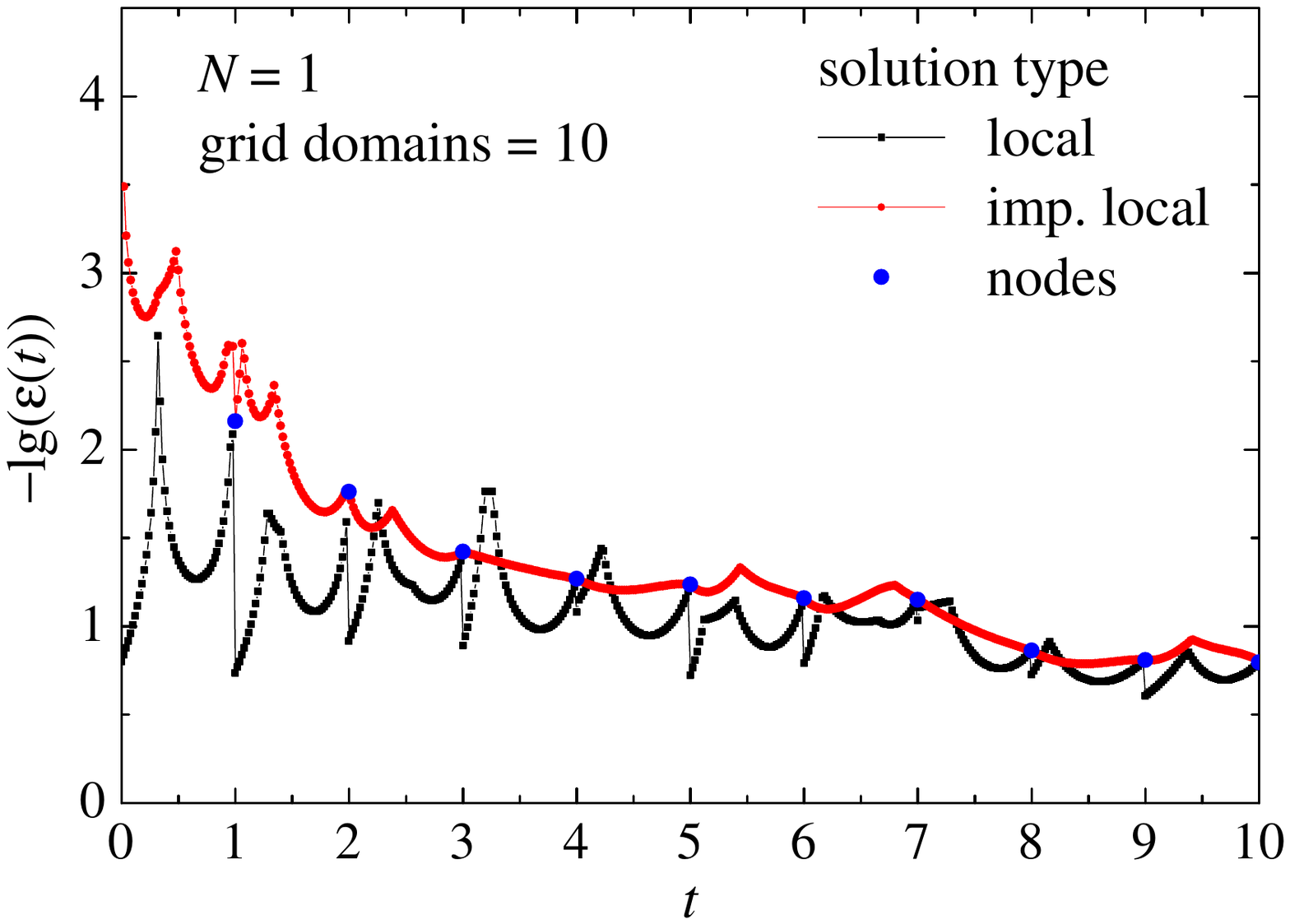}
\vspace{-8mm}\caption{\label{fig:pend:c1}}
\end{subfigure}
\begin{subfigure}{0.24\textwidth}
\includegraphics[width=\textwidth]{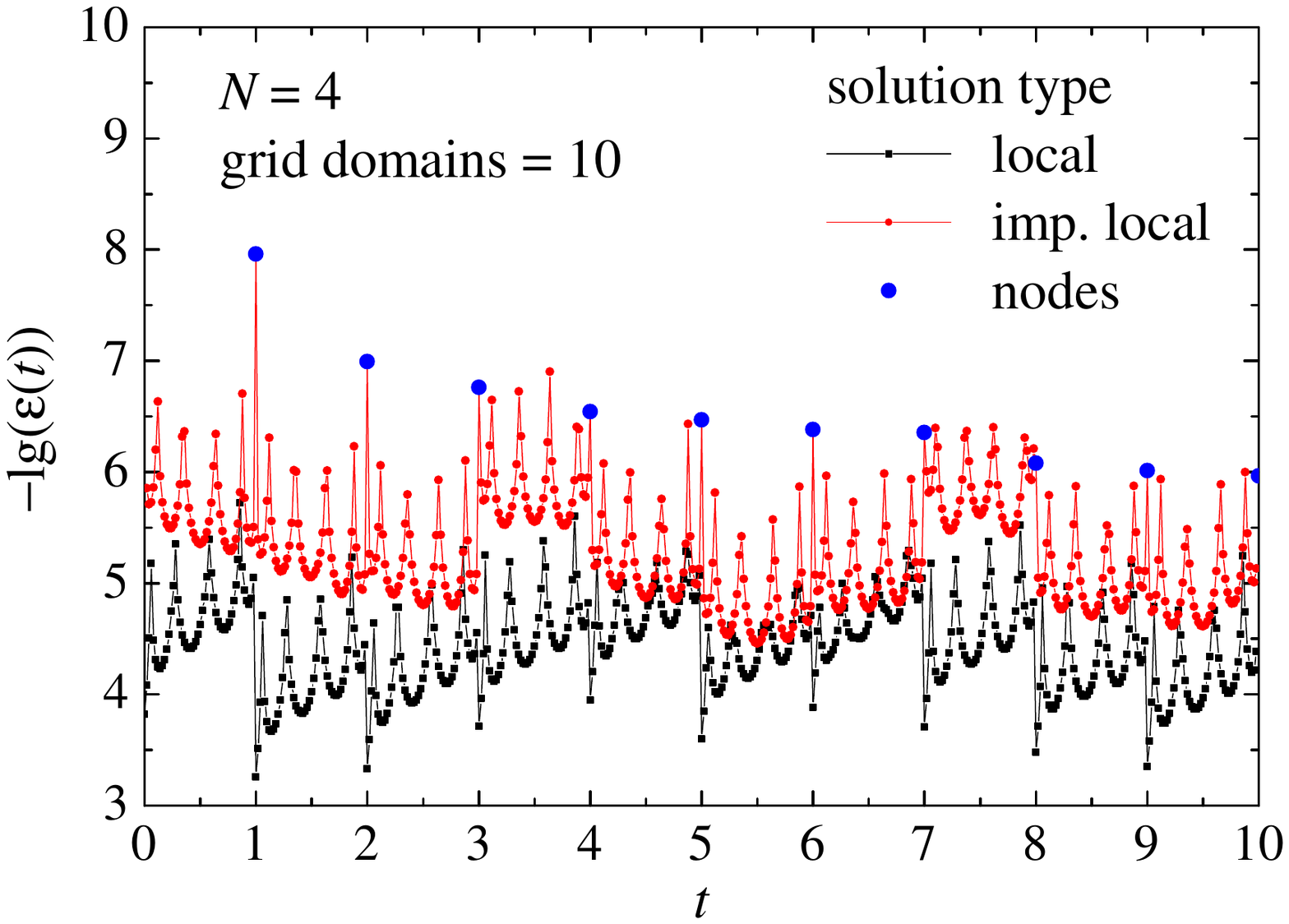}
\vspace{-8mm}\caption{\label{fig:pend:c2}}
\end{subfigure}
\begin{subfigure}{0.24\textwidth}
\includegraphics[width=\textwidth]{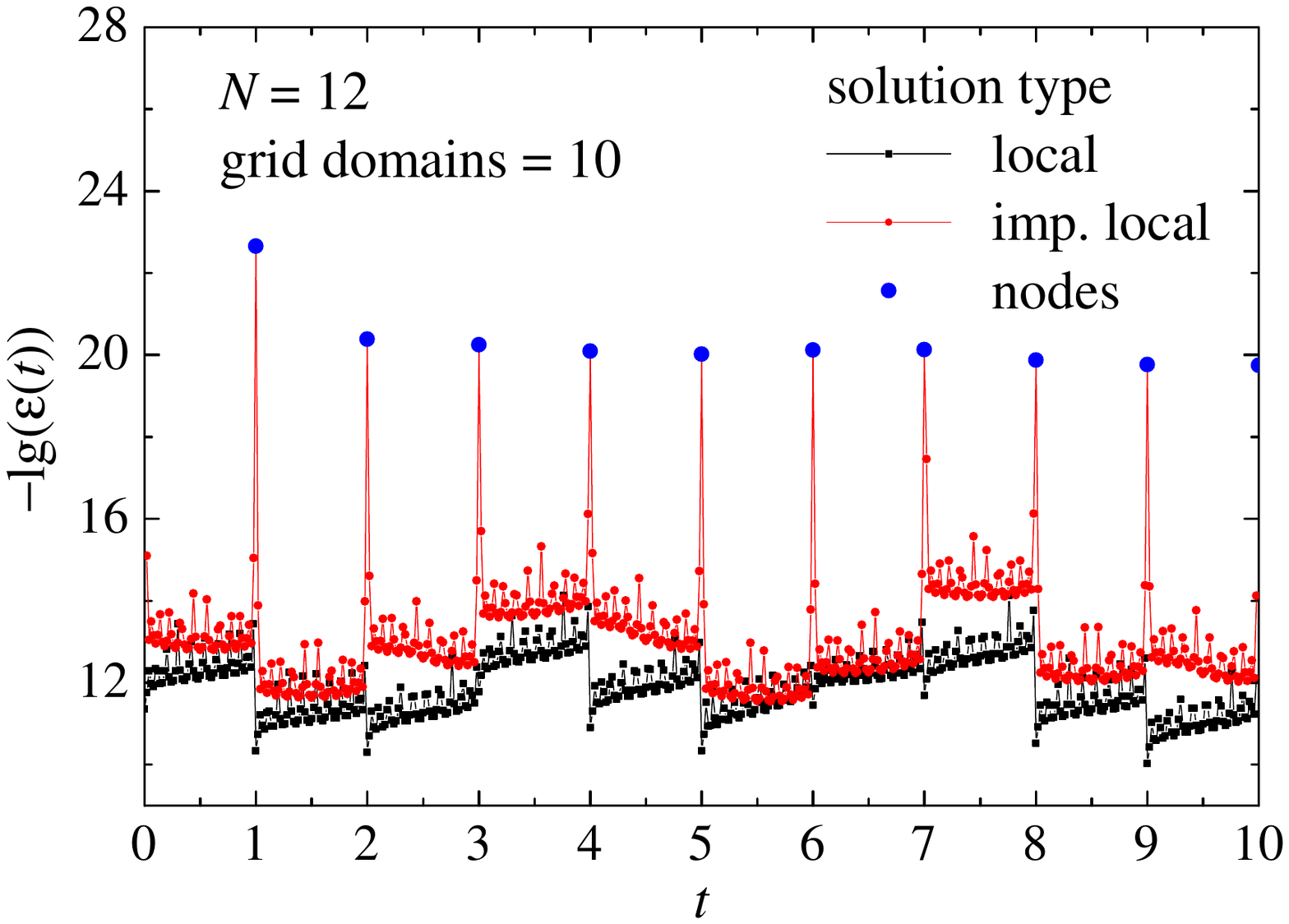}
\vspace{-8mm}\caption{\label{fig:pend:c3}}
\end{subfigure}
\begin{subfigure}{0.24\textwidth}
\includegraphics[width=\textwidth]{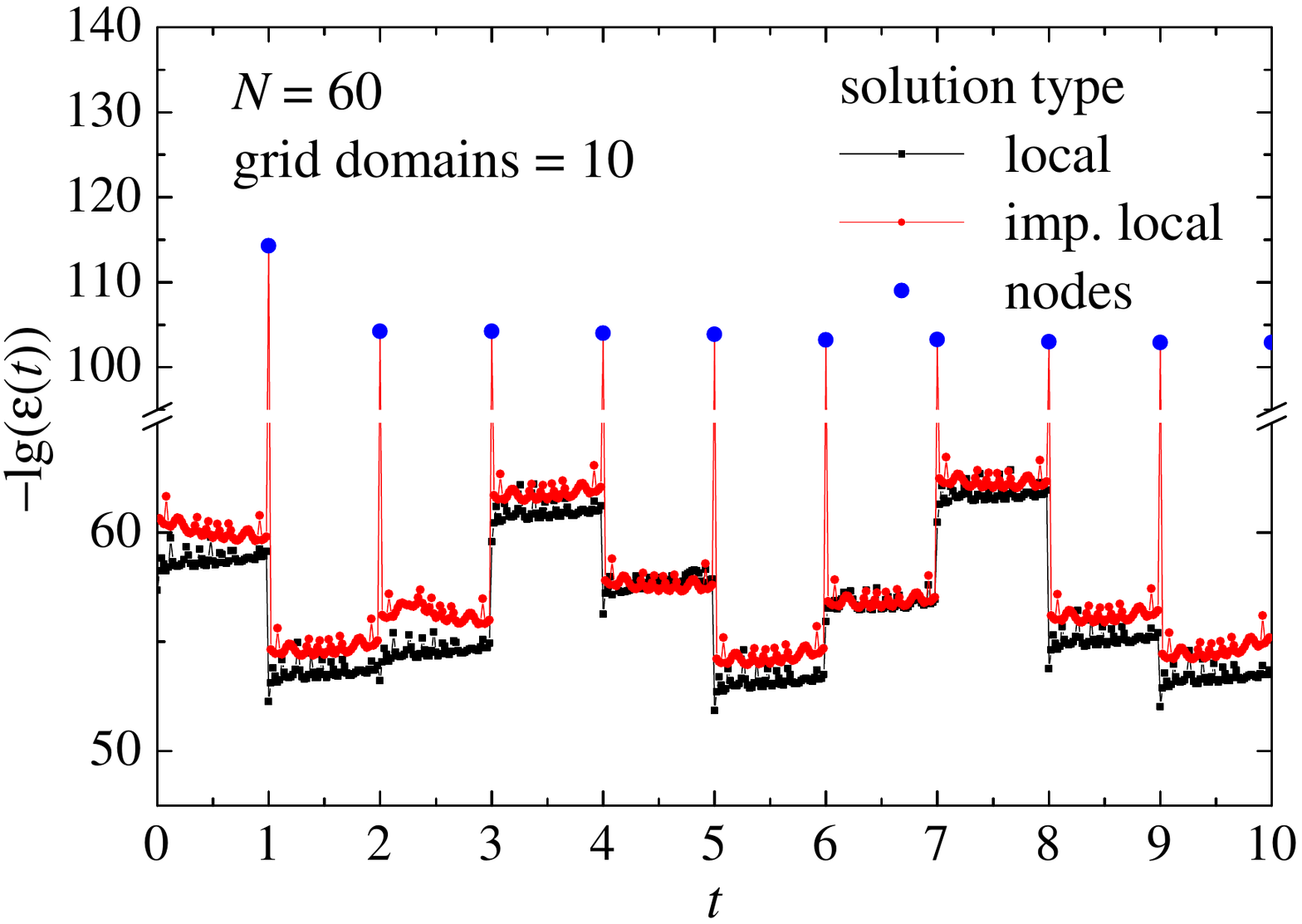}
\vspace{-8mm}\caption{\label{fig:pend:c4}}
\end{subfigure}\\
\begin{subfigure}{0.24\textwidth}
\includegraphics[width=\textwidth]{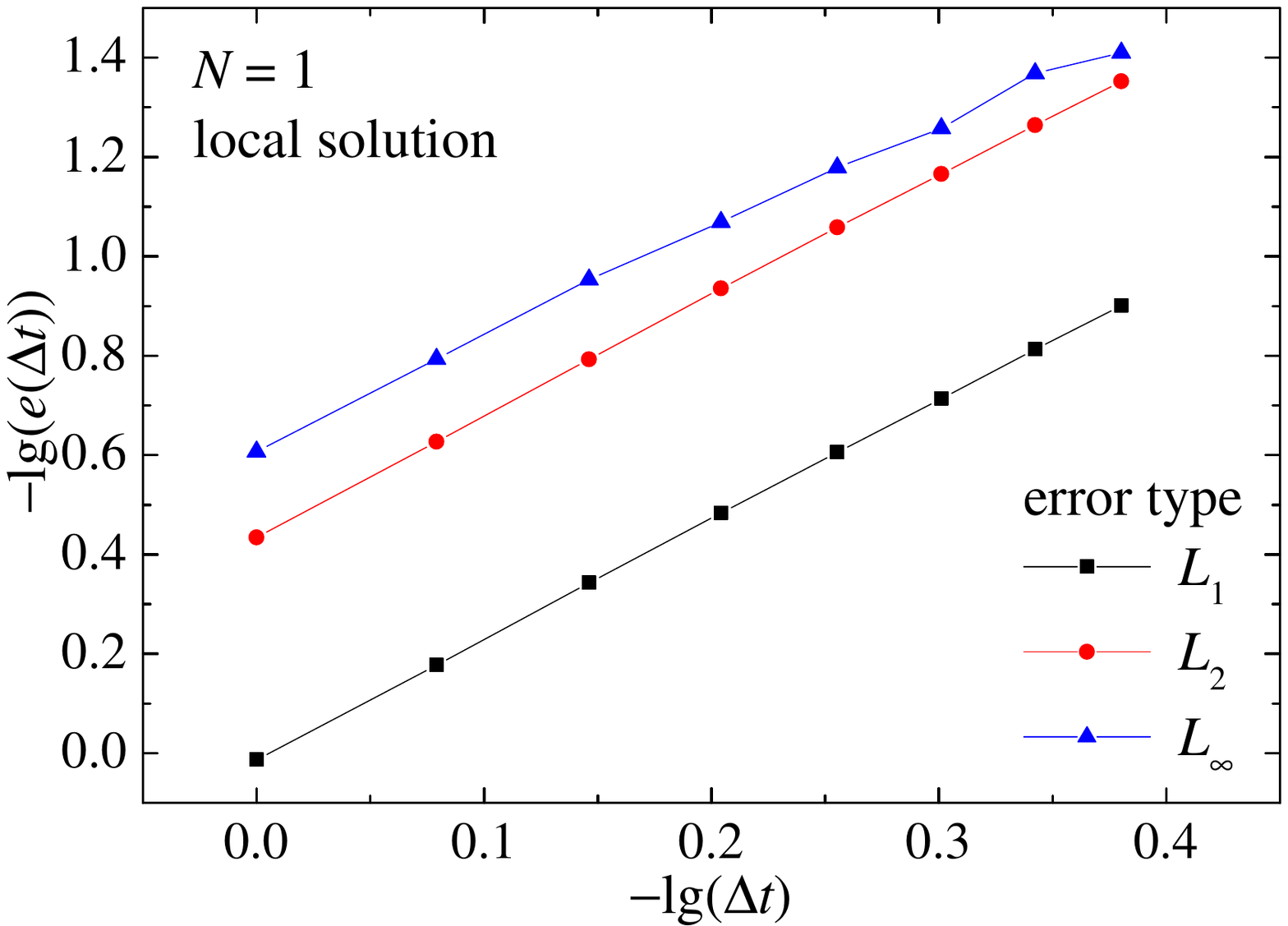}
\vspace{-8mm}\caption{\label{fig:pend:d1}}
\end{subfigure}
\begin{subfigure}{0.24\textwidth}
\includegraphics[width=\textwidth]{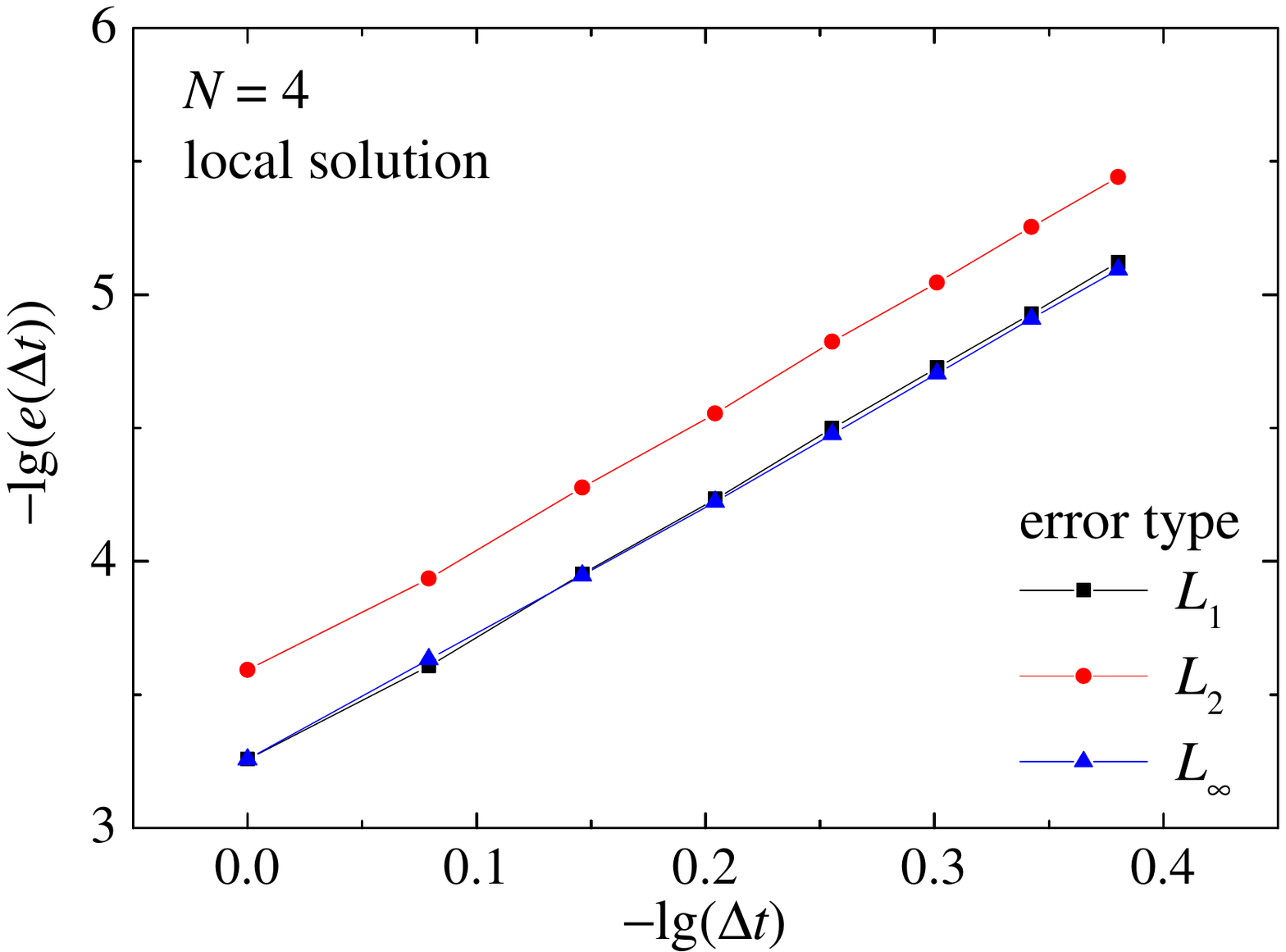}
\vspace{-8mm}\caption{\label{fig:pend:d2}}
\end{subfigure}
\begin{subfigure}{0.24\textwidth}
\includegraphics[width=\textwidth]{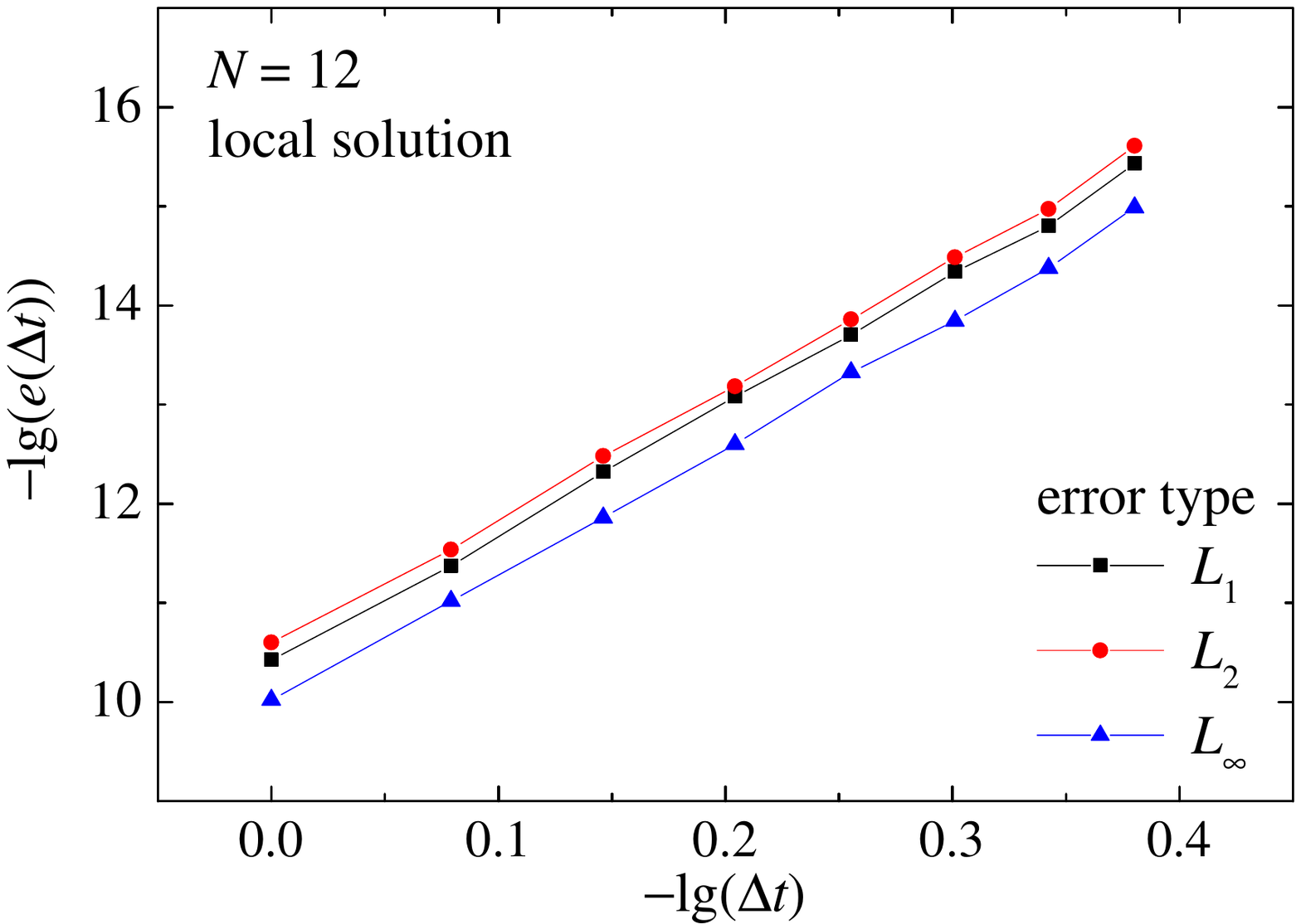}
\vspace{-8mm}\caption{\label{fig:pend:d3}}
\end{subfigure}
\begin{subfigure}{0.24\textwidth}
\includegraphics[width=\textwidth]{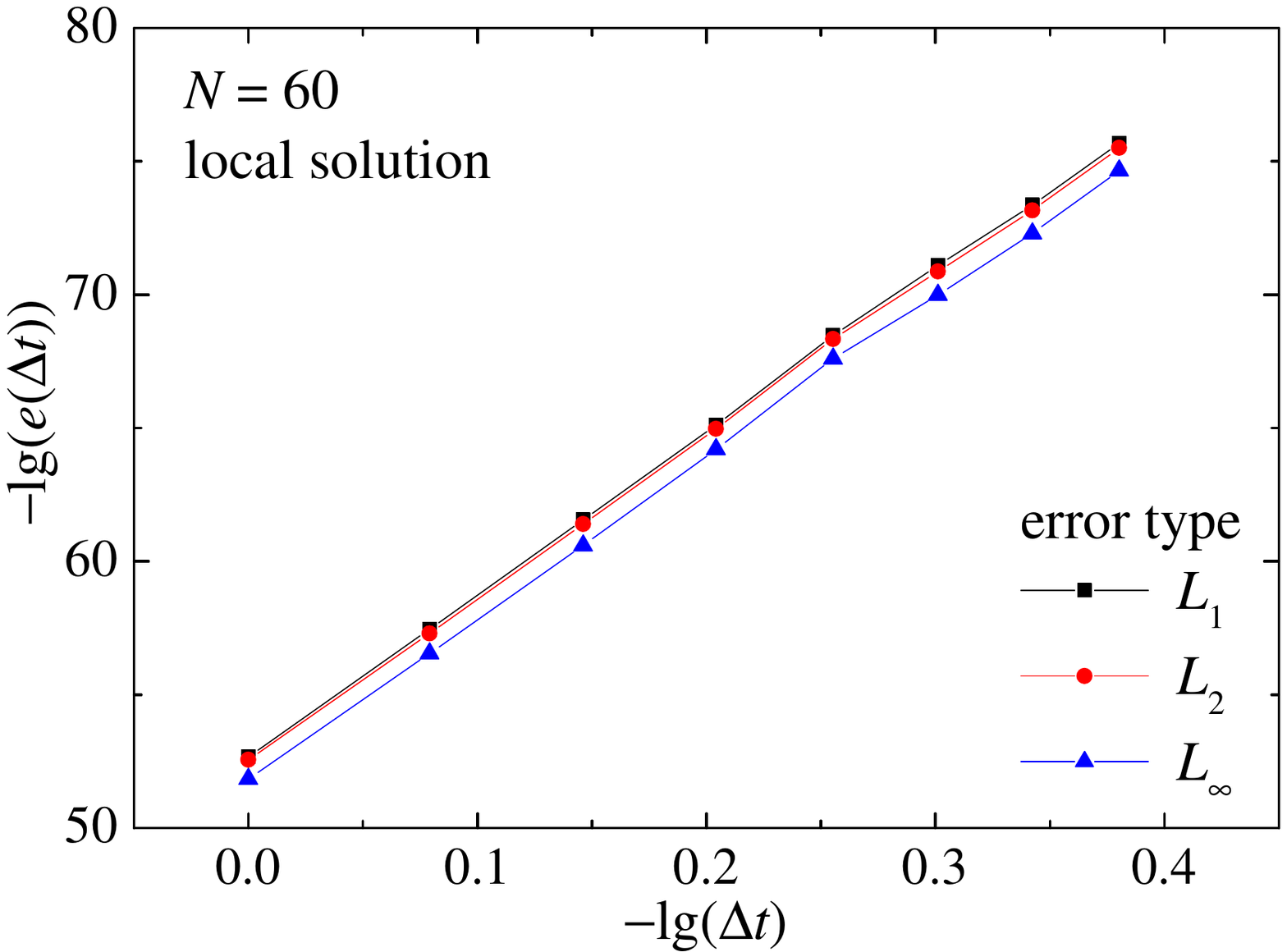}
\vspace{-8mm}\caption{\label{fig:pend:d4}}
\end{subfigure}\\
\begin{subfigure}{0.24\textwidth}
\includegraphics[width=\textwidth]{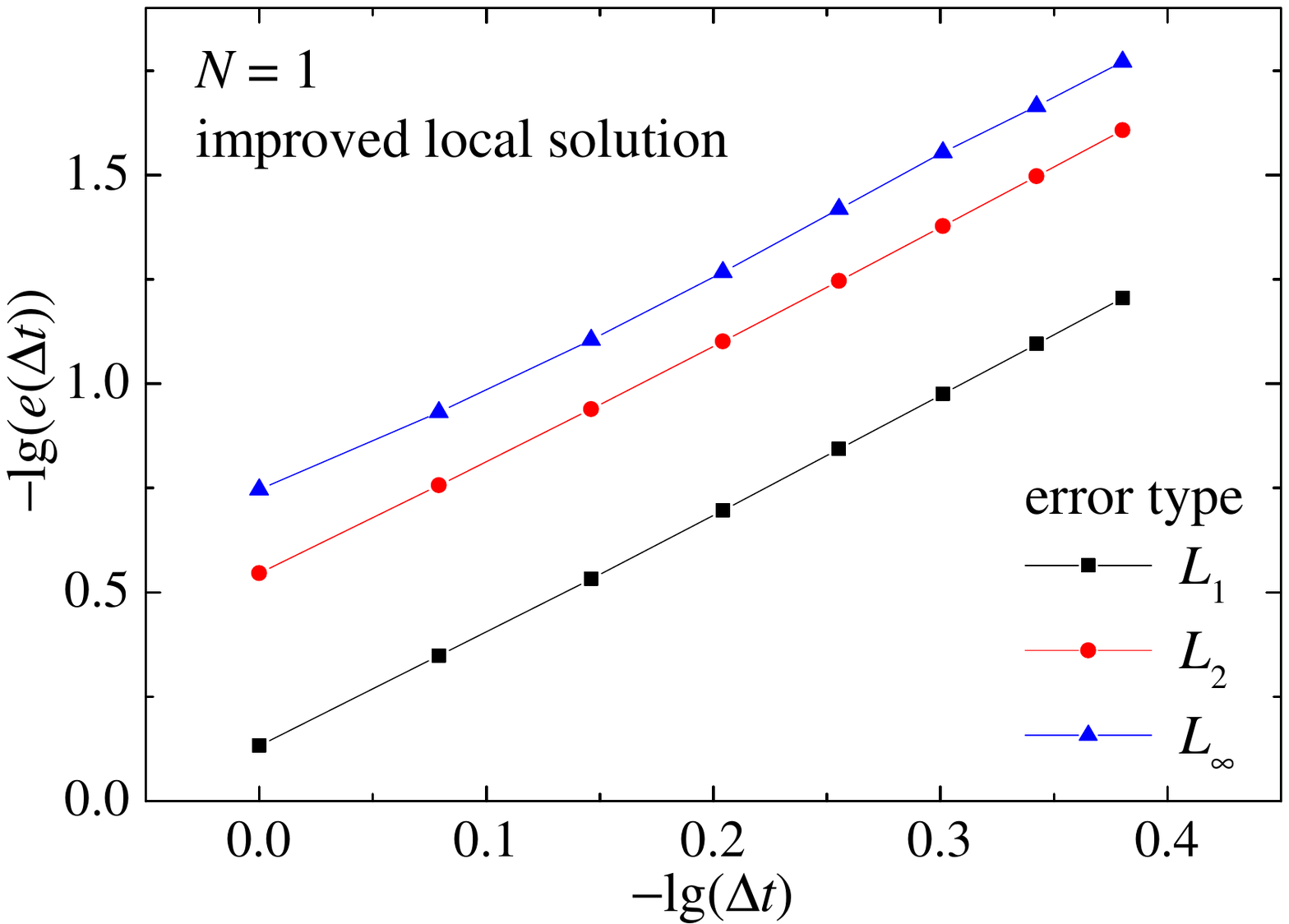}
\vspace{-8mm}\caption{\label{fig:pend:e1}}
\end{subfigure}
\begin{subfigure}{0.24\textwidth}
\includegraphics[width=\textwidth]{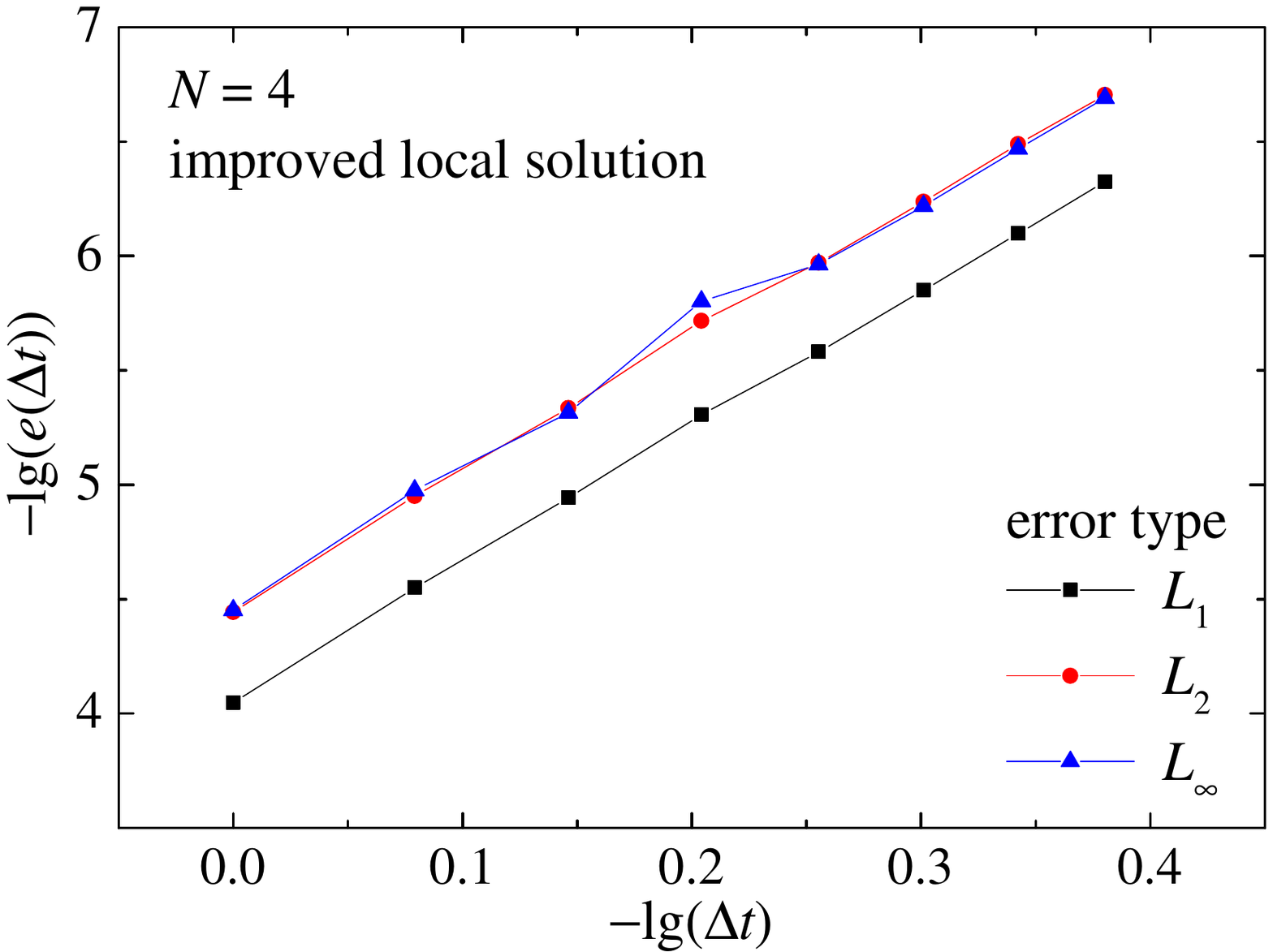}
\vspace{-8mm}\caption{\label{fig:pend:e2}}
\end{subfigure}
\begin{subfigure}{0.24\textwidth}
\includegraphics[width=\textwidth]{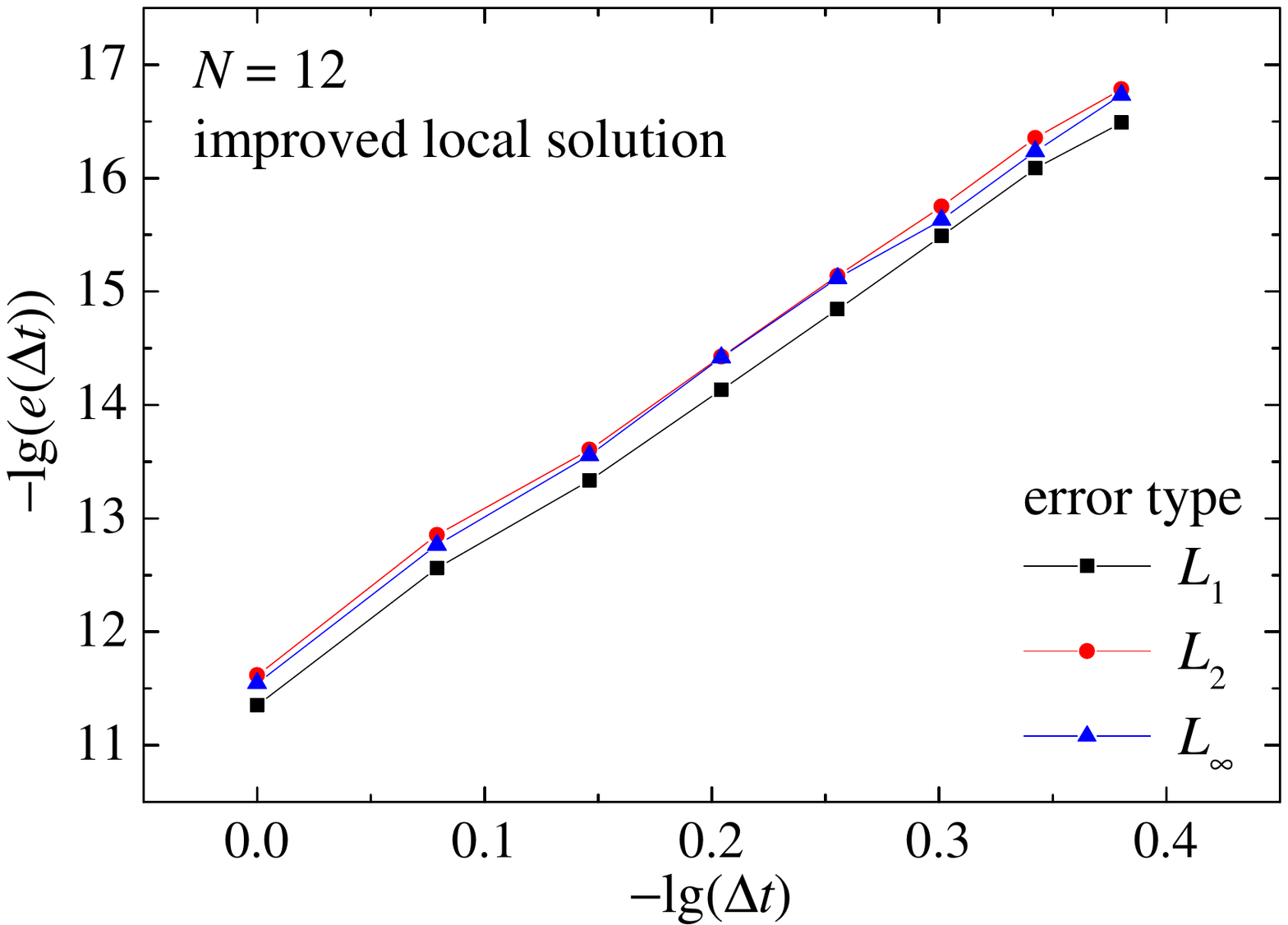}
\vspace{-8mm}\caption{\label{fig:pend:e3}}
\end{subfigure}
\begin{subfigure}{0.24\textwidth}
\includegraphics[width=\textwidth]{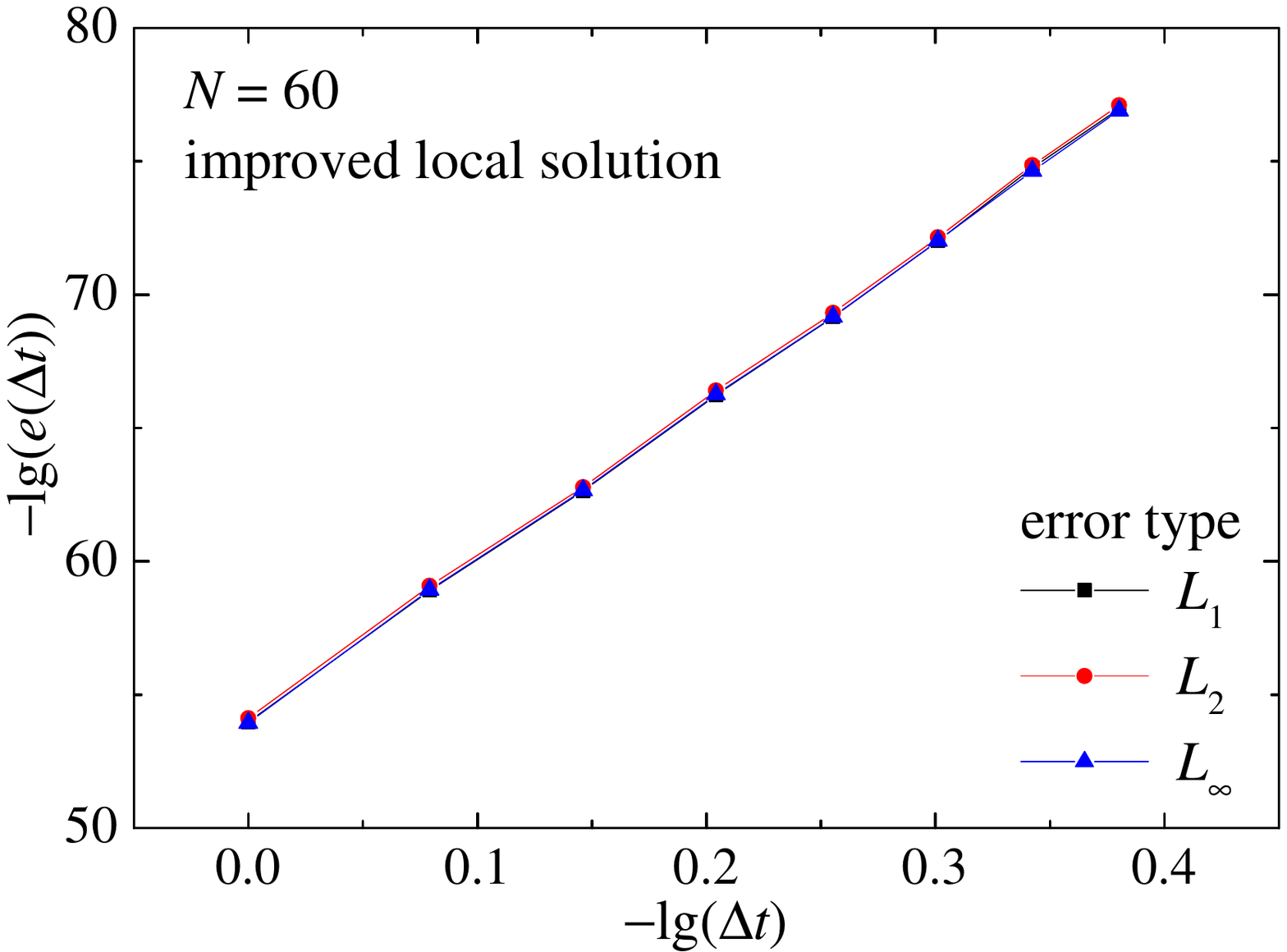}
\vspace{-8mm}\caption{\label{fig:pend:e4}}
\end{subfigure}\\
\begin{subfigure}{0.24\textwidth}
\includegraphics[width=\textwidth]{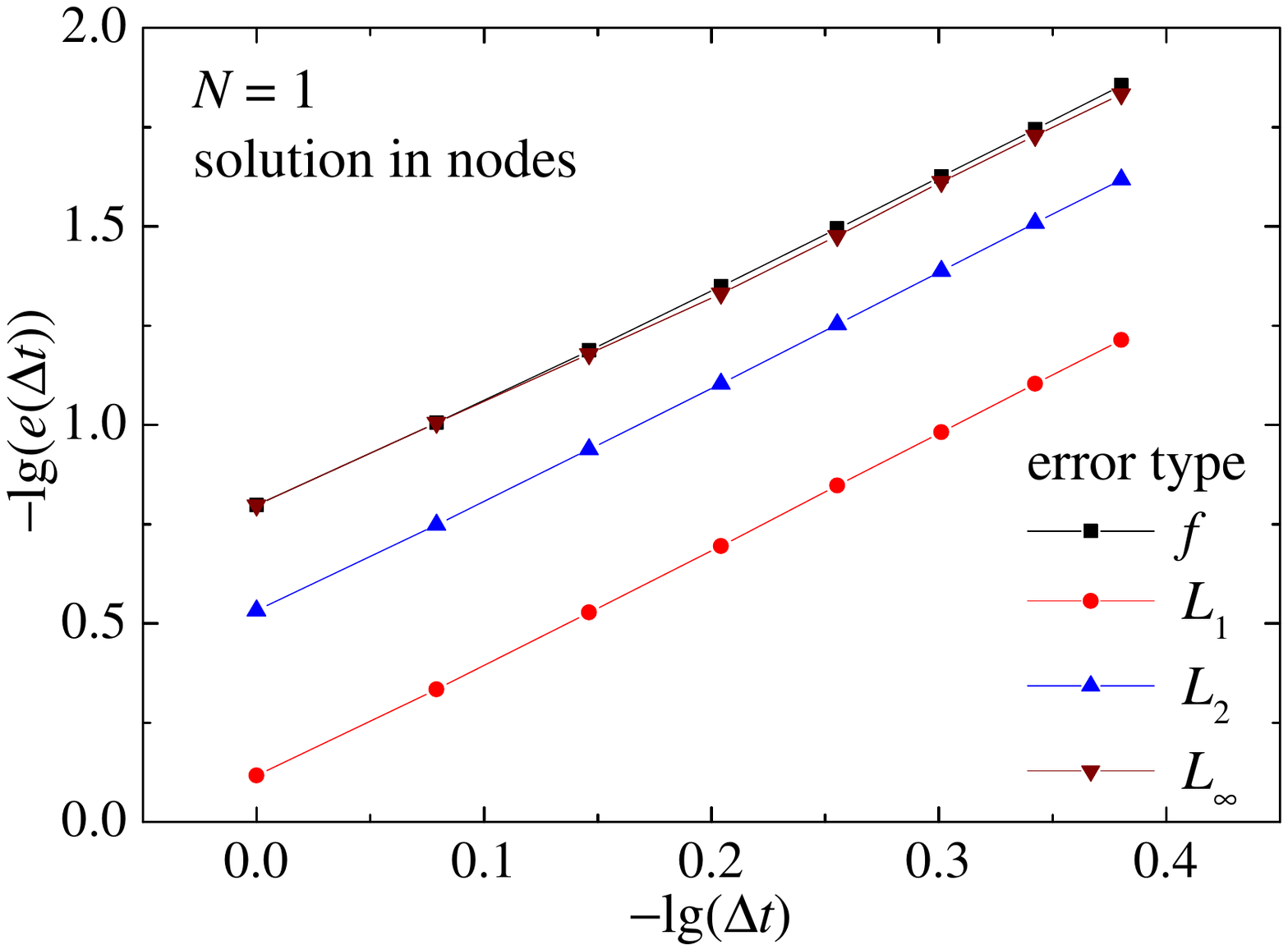}
\vspace{-8mm}\caption{\label{fig:pend:f1}}
\end{subfigure}
\begin{subfigure}{0.24\textwidth}
\includegraphics[width=\textwidth]{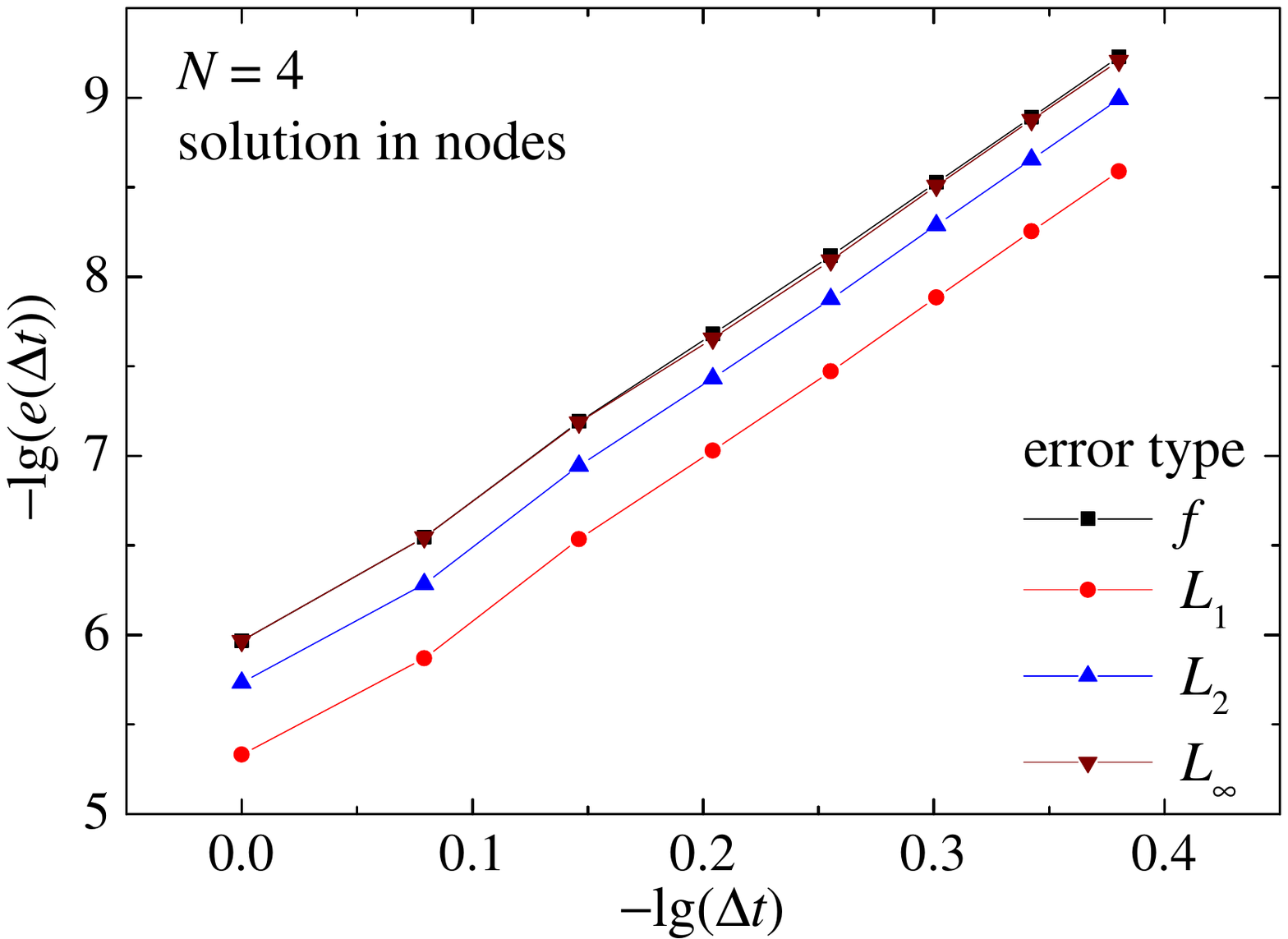}
\vspace{-8mm}\caption{\label{fig:pend:f2}}
\end{subfigure}
\begin{subfigure}{0.24\textwidth}
\includegraphics[width=\textwidth]{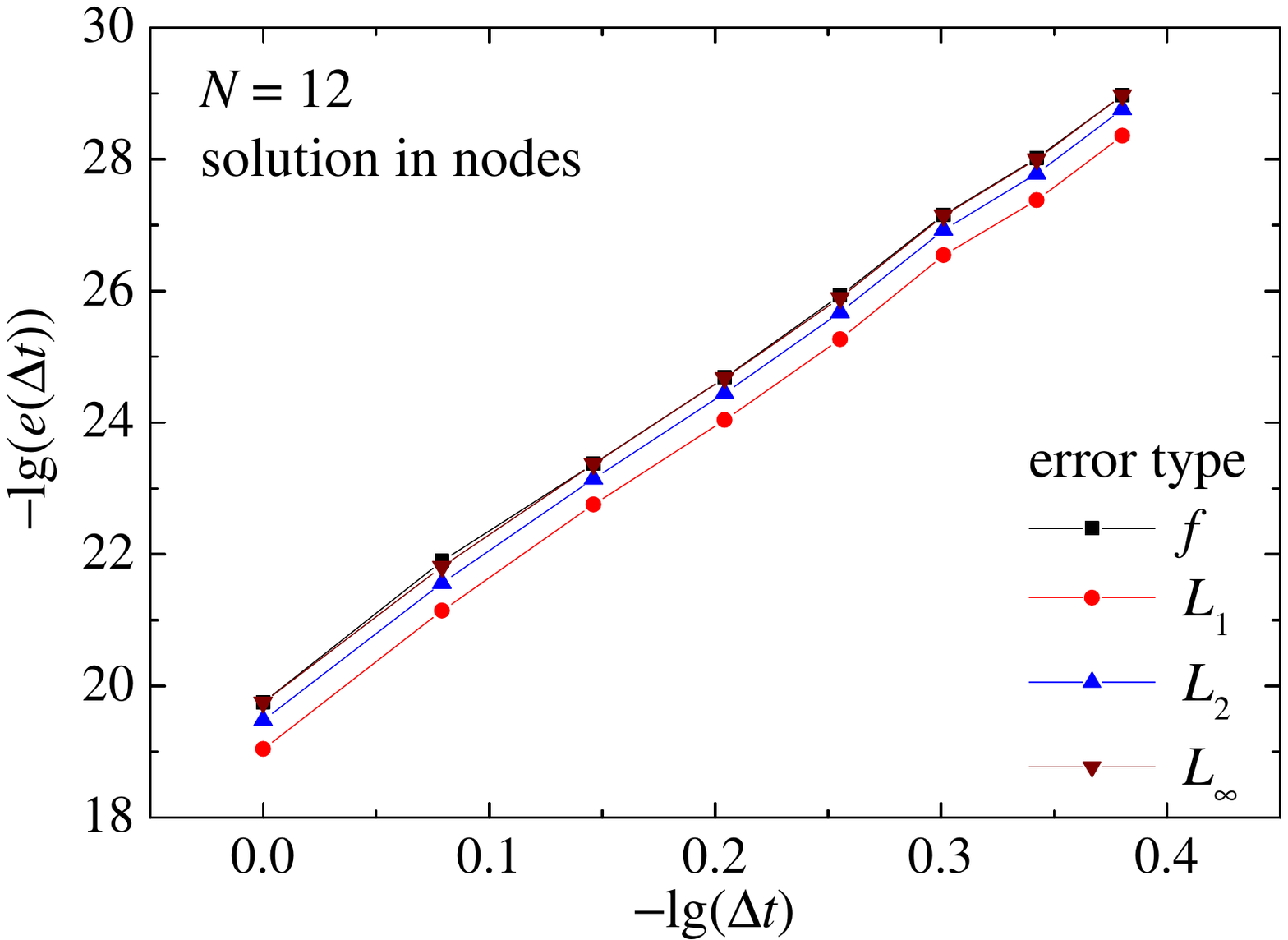}
\vspace{-8mm}\caption{\label{fig:pend:f3}}
\end{subfigure}
\begin{subfigure}{0.24\textwidth}
\includegraphics[width=\textwidth]{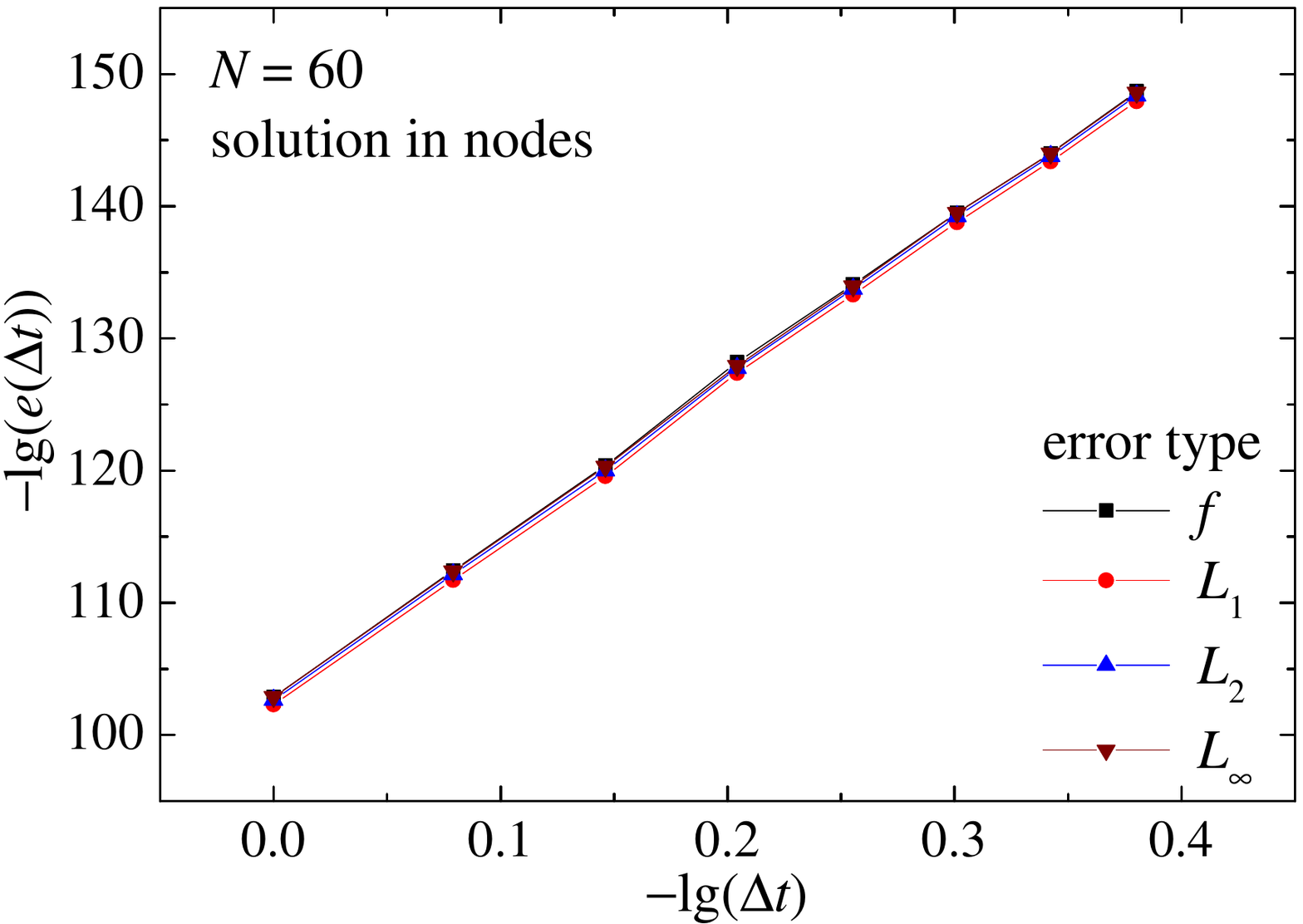}
\vspace{-8mm}\caption{\label{fig:pend:f4}}
\end{subfigure}\\
\caption{%
Numerical solution of the system (\ref{eq:pend_ode}). Comparison of the solution at nodes $\mathbf{u}_{n}$, the local solution $\mathbf{u}_{L}(t)$, the improved local solution $\mathbf{u}_{\rm IL}(t)$ and the exact solution $\mathbf{u}^{\rm ex}(t)$ (\ref{eq:pend_sol_ex}) for components $u_{1} \equiv \phi$ (\subref{fig:pend:a1}, \subref{fig:pend:a2}, \subref{fig:pend:a3}, \subref{fig:pend:a4}) and $u_{2} \equiv \dot{\phi}$ (\subref{fig:pend:b1}, \subref{fig:pend:b2}, \subref{fig:pend:b3}, \subref{fig:pend:b4}), the errors $\varepsilon(t)$ (\subref{fig:pend:c1}, \subref{fig:pend:c2}, \subref{fig:pend:c3}, \subref{fig:pend:c4}), obtained using polynomials with degrees $N = 1$ (\subref{fig:pend:a1}, \subref{fig:pend:b1}, \subref{fig:pend:c1}), $N = 4$ (\subref{fig:pend:a2}, \subref{fig:pend:b2}, \subref{fig:pend:c2}), $N = 12$ (\subref{fig:pend:a3}, \subref{fig:pend:b3}, \subref{fig:pend:c3}) and $N = 60$ (\subref{fig:pend:a4}, \subref{fig:pend:b4}, \subref{fig:pend:c4}). Log-log plot of the dependence of the global error for the local solution $e^{l}$ (\subref{fig:pend:d1}, \subref{fig:pend:d2}, \subref{fig:pend:d3}, \subref{fig:pend:d4}), the improved local solution $e^{\rm imp}$ (\subref{fig:pend:e1}, \subref{fig:pend:e2}, \subref{fig:pend:e3}, \subref{fig:pend:e4}) and the solution at nodes $e^{n}$ (\subref{fig:pend:f1}, \subref{fig:pend:f2}, \subref{fig:pend:f3}, \subref{fig:pend:f4}) on the discretization step $\mathrm{\Delta}t$, obtained in the $f$-norm and norms $L_{1}$, $L_{2}$ and $L_{\infty}$, obtained using polynomials with degrees $N = 1$ (\subref{fig:pend:d1}, \subref{fig:pend:e1}, \subref{fig:pend:f1}), $N = 4$ (\subref{fig:pend:d2}, \subref{fig:pend:e2}, \subref{fig:pend:f2}), $N = 12$ (\subref{fig:pend:d3}, \subref{fig:pend:e3}, \subref{fig:pend:f3}) and $N = 60$ (\subref{fig:pend:d4}, \subref{fig:pend:e4}, \subref{fig:pend:f4}).
}
\label{fig:pend}
\end{figure}

The fourth example of the application of the ADER-DG numerical method with a local DG predictor to solving the problem of the initial value for the ODE system (\ref{eq:ivp_ode_diff_src}) is the problem of describing the motion of a mathematical pendulum in a complete nonlinear formulation:
\begin{equation}\label{eq:pend_ode}
\ddot{\phi} + \omega_{0}^{2}\sin(\phi) = 0,\quad
\phi(0) = \phi_{0},\quad \dot{\phi}(0) = 0,\quad
t \in [0,\, t_{f}].
\end{equation}
The main difference from the previous examples presented in Subsections~\ref{sec:apps:exp_diss},~\ref{sec:apps:lin_diss} and~\ref{sec:apps:harm_osc} is the nonlinearity of the problem formulation. Exact analytical solution for the function $\phi(t)$ was obtained by trivial integration (see, for example, interesting derivation and detailed analysis of the analytical solution in the work~\cite{math_pend_exact_sol_ref}):
\begin{equation}\label{eq:pend_sol_ex}
\begin{split}
&\phi^{\rm ex}(t) = 2\,\arcsin\left\{\Upsilon_{0}\cdot
	\mathrm{sn}\left[K\left(\Upsilon_{0}\right) - \omega_{0}t,\, \Upsilon_{0}\right]\right\},\quad
\Upsilon_{0} = \sin\left(\frac{\phi_{0}}{2}\right),\\
&\dot{\phi}^{\rm ex}(t) = -\frac{\displaystyle
	2\omega_{0}\Upsilon_{0}\cdot
	\mathrm{cn}\left[K\left(\Upsilon_{0}\right) - \omega_{0}t,\, \Upsilon_{0}\right]\cdot
	\mathrm{dn}\left[K\left(\Upsilon_{0}\right) - \omega_{0}t,\, \Upsilon_{0}\right]
}{\displaystyle
	\sqrt{1 - \Upsilon_{0}^{2}\cdot\mathrm{sn}^{2}\left[K\left(\Upsilon_{0}\right) - \omega_{0}t,\, \Upsilon_{0}\right]}
},
\end{split}
\end{equation}
where $K$ is the complete elliptical integral of the first kind, $(\mathrm{sn},\, \mathrm{cn},\, \mathrm{dn})$ are the set of Jacobi elliptic functions: sine, cosine, delta amplitude, respectively. It should be noted that the problem solved in this Subsection is nonlinear --- the presented solution defines the law of motion of a mathematical pendulum for an arbitrary amplitude of oscillations (limited only by the zero value of the initial velocity $\dot{\phi}(0) = 0$, therefore, it is impossible to obtain a full rotation in this formulation of the problem), and not only for small values of amplitude, when it would be possible to limit oneself only to the harmonic approximation of small oscillations, as in the problem (\ref{eq:harm_osc_ode}). The values $\omega_{0} = 1$, $\phi_{0} = \pi/2$ and $t_{f} = 10$ are chosen in the calculations. The solution vectors in the original notation of the ODE system (\ref{eq:ivp_ode_diff_src}), with $D = 2$, are chosen in the form $\mathbf{u}(t) = [\phi(t)\ \dot{\phi}(t)]^{T}$. The obtained results are presented in Fig.~\ref{fig:pend} and Table~\ref{tab:conv_orders_pend}.

\begin{table*}[h!]
\centering
\normalsize
\caption{%
Convergence orders $p_{f}$, $p_{L_{1}}$, $p_{L_{2}}$, $p_{L_{\infty}}$, calculated in $f$-norm and norms $L_{1}$, $L_{2}$, $L_{\infty}$, of the numerical solution of the ADER-DG method for the problem (\ref{eq:pend_ode}); $N$ is the degree of the basis polynomials $\varphi_{p}$. Orders $p^{n}$ are calculated for \textit{the numerical solution at the nodes} $\mathbf{u}_{n}$; orders $p^{\rm imp}$ --- for \textit{the improved local solution} $\mathbf{u}_{\rm IL}$; orders $p^{l}$ --- for \textit{the local solution} $\mathbf{u}_{L}$. The theoretical values of convergence order $p_{\rm th.}^{n} = 2N+1$, $p_{\rm th.}^{l} = N+1$ and $p^{\rm imp}_{\rm th.} = N+2$ are presented for comparison.
}
\label{tab:conv_orders_pend}
\setlength{\tabcolsep}{3.5pt}
\begin{tabular}{@{}|l|llll|c|lll|c|lll|c|@{}}
\toprule
$N$ & $p^{n}_{f}$ &
$p^{n}_{L_{1}}$ & $p^{n}_{L_{2}}$ & $p^{n}_{L_{\infty}}$ & $p^{n}_{\rm th.}$ &
$p^{l}_{L_{1}}$ & $p^{l}_{L_{2}}$ & $p^{l}_{L_{\infty}}$ & $p^{l}_{\rm th.}$ &
$p^{\rm imp}_{L_{1}}$ & $p^{\rm imp}_{L_{2}}$ & $p^{\rm imp}_{L_{\infty}}$ & $p^{\rm imp}_{\rm th.}$\\
\midrule
$1$ & $2.79$ & $2.90$ & $2.87$ & $2.73$ & $3$ & $2.41$ & $2.42$ & $2.13$ & $2$ & $2.82$ & $2.80$ & $2.74$ & $3$\\
$2$ & $4.78$ & $4.86$ & $4.84$ & $4.76$ & $5$ & $3.09$ & $2.99$ & $2.93$ & $3$ & $4.41$ & $4.43$ & $4.21$ & $4$\\
$3$ & $6.81$ & $6.96$ & $6.93$ & $6.82$ & $7$ & $3.98$ & $3.96$ & $3.86$ & $4$ & $4.99$ & $4.96$ & $4.90$ & $5$\\
$4$ & $8.66$ & $8.70$ & $8.69$ & $8.60$ & $9$ & $4.94$ & $4.91$ & $4.84$ & $5$ & $5.95$ & $5.90$ & $5.81$ & $6$\\
$5$ & $10.8$ & $11.0$ & $11.0$ & $10.9$ & $11$ & $5.97$ & $5.91$ & $5.78$ & $6$ & $6.91$ & $6.90$ & $6.88$ & $7$\\
$6$ & $12.6$ & $12.7$ & $12.6$ & $12.6$ & $13$ & $6.92$ & $6.91$ & $6.87$ & $7$ & $7.83$ & $7.82$ & $7.70$ & $8$\\
$7$ & $14.7$ & $14.8$ & $14.8$ & $14.7$ & $15$ & $7.83$ & $7.81$ & $7.68$ & $8$ & $8.97$ & $8.99$ & $8.83$ & $9$\\
$8$ & $16.6$ & $16.7$ & $16.7$ & $16.6$ & $17$ & $8.97$ & $8.99$ & $8.81$ & $9$ & $9.71$ & $9.67$ & $9.60$ & $10$\\
$9$ & $18.6$ & $18.7$ & $18.7$ & $18.6$ & $19$ & $9.71$ & $9.66$ & $9.58$ & $10$ & $11.1$ & $11.1$ & $10.9$ & $11$\\
$10$ & $20.4$ & $20.6$ & $20.6$ & $20.4$ & $21$ & $11.1$ & $11.1$ & $10.9$ & $11$ & $11.6$ & $11.6$ & $11.5$ & $12$\\
\midrule
$11$ & $22.8$ & $22.8$ & $22.9$ & $22.8$ & $23$ & $11.6$ & $11.6$ & $11.5$ & $12$ & $13.1$ & $13.1$ & $13.0$ & $13$\\
$12$ & $24.0$ & $24.3$ & $24.2$ & $24.1$ & $25$ & $13.1$ & $13.1$ & $12.9$ & $13$ & $13.5$ & $13.5$ & $13.5$ & $14$\\
$13$ & $27.1$ & $27.1$ & $27.2$ & $27.1$ & $27$ & $13.6$ & $13.6$ & $13.5$ & $14$ & $15.0$ & $15.0$ & $14.8$ & $15$\\
$14$ & $27.3$ & $27.6$ & $27.5$ & $27.3$ & $29$ & $15.0$ & $14.9$ & $14.8$ & $15$ & $15.7$ & $15.7$ & $15.5$ & $16$\\
$15$ & $31.3$ & $31.3$ & $31.3$ & $31.2$ & $31$ & $15.7$ & $15.7$ & $15.5$ & $16$ & $16.7$ & $16.7$ & $16.7$ & $17$\\
$16$ & $30.5$ & $31.2$ & $31.1$ & $30.9$ & $33$ & $16.6$ & $16.7$ & $16.6$ & $17$ & $17.9$ & $17.8$ & $17.5$ & $18$\\
$17$ & $35.2$ & $35.4$ & $35.3$ & $35.2$ & $35$ & $18.0$ & $17.8$ & $17.5$ & $18$ & $18.6$ & $18.6$ & $18.4$ & $19$\\
$18$ & $35.4$ & $35.3$ & $35.4$ & $35.4$ & $37$ & $18.5$ & $18.6$ & $18.4$ & $19$ & $19.9$ & $19.8$ & $19.4$ & $20$\\
$19$ & $39.1$ & $39.3$ & $39.3$ & $39.2$ & $39$ & $19.9$ & $19.8$ & $19.4$ & $20$ & $20.7$ & $20.6$ & $20.3$ & $21$\\
$20$ & $39.7$ & $39.4$ & $39.5$ & $39.6$ & $41$ & $20.7$ & $20.6$ & $20.3$ & $21$ & $21.7$ & $21.7$ & $21.4$ & $22$\\
\midrule
$21$ & $43.0$ & $43.2$ & $43.1$ & $43.0$ & $43$ & $21.8$ & $21.7$ & $21.4$ & $22$ & $22.7$ & $22.6$ & $22.4$ & $23$\\
$22$ & $43.7$ & $43.7$ & $43.6$ & $43.6$ & $45$ & $22.6$ & $22.6$ & $22.3$ & $23$ & $23.7$ & $23.7$ & $23.5$ & $24$\\
$23$ & $46.7$ & $46.9$ & $46.9$ & $46.8$ & $47$ & $23.7$ & $23.7$ & $23.5$ & $24$ & $24.6$ & $24.6$ & $24.4$ & $25$\\
$24$ & $46.9$ & $47.9$ & $47.8$ & $47.6$ & $49$ & $24.5$ & $24.5$ & $24.3$ & $25$ & $25.7$ & $25.7$ & $25.5$ & $26$\\
$25$ & $50.1$ & $50.2$ & $50.3$ & $50.3$ & $51$ & $25.7$ & $25.7$ & $25.5$ & $26$ & $26.5$ & $26.5$ & $26.4$ & $27$\\
$26$ & $51.6$ & $52.3$ & $52.1$ & $51.7$ & $53$ & $26.5$ & $26.5$ & $26.3$ & $27$ & $27.8$ & $27.7$ & $27.4$ & $28$\\
$27$ & $54.5$ & $54.6$ & $54.7$ & $54.5$ & $55$ & $27.8$ & $27.7$ & $27.4$ & $28$ & $28.5$ & $28.5$ & $28.3$ & $29$\\
$28$ & $55.8$ & $56.1$ & $56.0$ & $55.7$ & $57$ & $28.6$ & $28.5$ & $28.2$ & $29$ & $29.8$ & $29.7$ & $29.3$ & $30$\\
$29$ & $58.6$ & $58.9$ & $58.9$ & $58.8$ & $59$ & $29.9$ & $29.8$ & $29.3$ & $30$ & $30.5$ & $30.5$ & $30.3$ & $31$\\
$30$ & $59.7$ & $59.8$ & $59.7$ & $59.6$ & $61$ & $30.5$ & $30.5$ & $30.2$ & $31$ & $31.8$ & $31.7$ & $31.2$ & $32$\\
\midrule
$35$ & $70.5$ & $71.2$ & $71.1$ & $70.7$ & $71$ & $35.3$ & $35.2$ & $35.0$ & $36$ & $36.6$ & $36.6$ & $36.4$ & $37$\\
$40$ & $79.3$ & $79.3$ & $79.3$ & $79.2$ & $81$ & $40.6$ & $40.6$ & $40.3$ & $41$ & $40.5$ & $40.4$ & $40.3$ & $42$\\
$45$ & $90.7$ & $91.2$ & $91.1$ & $90.7$ & $91$ & $45.2$ & $44.9$ & $44.7$ & $46$ & $46.4$ & $46.4$ & $46.1$ & $47$\\
$50$ & $100.1$ & $100.1$ & $100.2$ & $100.1$ & $101$ & $50.1$ & $50.3$ & $50.1$ & $51$ & $51.2$ & $51.1$ & $50.7$ & $52$\\
$55$ & $109.2$ & $109.6$ & $109.6$ & $109.6$ & $111$ & $54.9$ & $54.8$ & $54.5$ & $56$ & $56.5$ & $56.5$ & $56.2$ & $57$\\
$60$ & $120.7$ & $120.6$ & $120.7$ & $120.7$ & $121$ & $60.7$ & $60.6$ & $60.2$ & $61$ & $60.3$ & $60.2$ & $60.1$ & $62$\\
\bottomrule
\end{tabular}
\end{table*}

The dependencies of the numerical solutions $\mathbf{u}_{L}$, $\mathbf{u}_{\rm IL}$, $\mathbf{u}_{n}$ and the exact analytical solution $\mathbf{u}^{\rm ex}$, the dependencies of the local error $\varepsilon$ (\ref{eq:eps_local_def}) of the numerical solutions, and the dependence of the global error $e$ (\ref{eq:eps_un_global_def}), (\ref{eq:eps_ul_global_def}) of the numerical solutions on the discretization step ${\Delta t}$, for polynomial degrees $N = 1$, $4$, $12$ and $60$, are shown in Fig.~\ref{fig:pend}. The obtained dependencies of the local numerical solution $\mathbf{u}_{L}$ and the improved local numerical solution $\mathbf{u}_{\rm IL}$ for polynomials degree $N = 1$ in Fig.~\ref{fig:pend} (\subref{fig:pend:a1}, \subref{fig:pend:b1}) demonstrate effects similar to those in the case of a harmonic oscillator in Fig.~\ref{fig:harm_osc} (\subref{fig:harm_osc:a1}, \subref{fig:harm_osc:b1}). A comparison of the obtained dependencies of the numerical solutions $\mathbf{u}_{L}(t)$, $\mathbf{u}_{\rm IL}(t)$, $\mathbf{u}_{n}$ for polynomial degrees $N = 4$, $12$ and $60$ with the exact analytical solution $\mathbf{u}^{\rm ex}(t)$, presented in Fig.~\ref{fig:pend} (\subref{fig:pend:a2}, \subref{fig:pend:a3}, \subref{fig:pend:a4}) for the component $u_{1}$ and in Fig.~\ref{fig:pend} (\subref{fig:pend:b2}, \subref{fig:pend:b3}, \subref{fig:pend:b4}) for the component $u_{2}$, demonstrates a high-quality comparison.

\begin{figure}[h!]
\captionsetup[subfigure]{%
	position=bottom,
	font+=smaller,
	textfont=normalfont,
	singlelinecheck=off,
	justification=raggedright
}
\centering
\begin{subfigure}{0.24\textwidth}
\includegraphics[width=\textwidth]{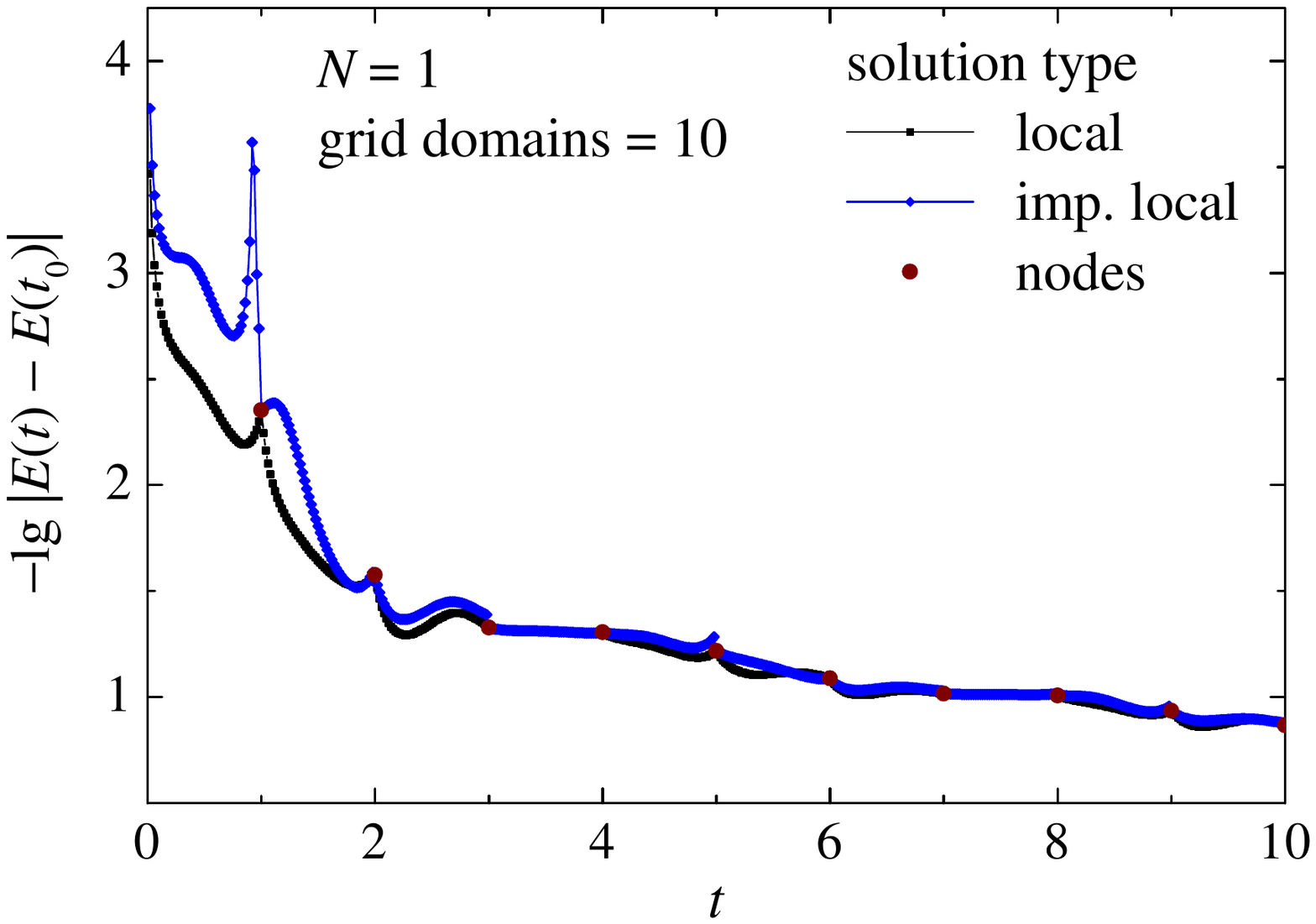}
\vspace{-8mm}\caption{\label{fig:econs_pend:a1}}
\end{subfigure}
\begin{subfigure}{0.24\textwidth}
\includegraphics[width=\textwidth]{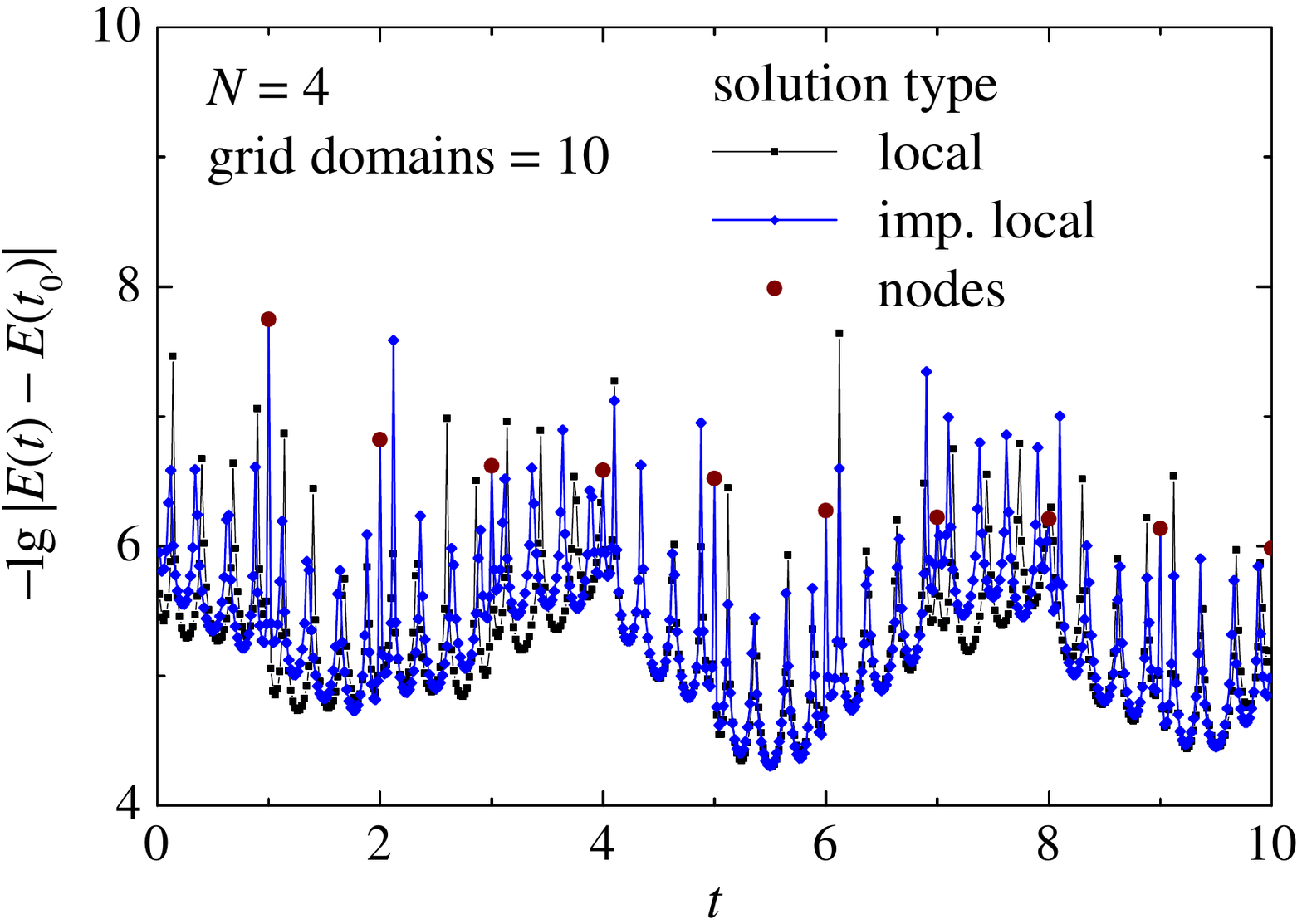}
\vspace{-8mm}\caption{\label{fig:econs_pend:a2}}
\end{subfigure}
\begin{subfigure}{0.24\textwidth}
\includegraphics[width=\textwidth]{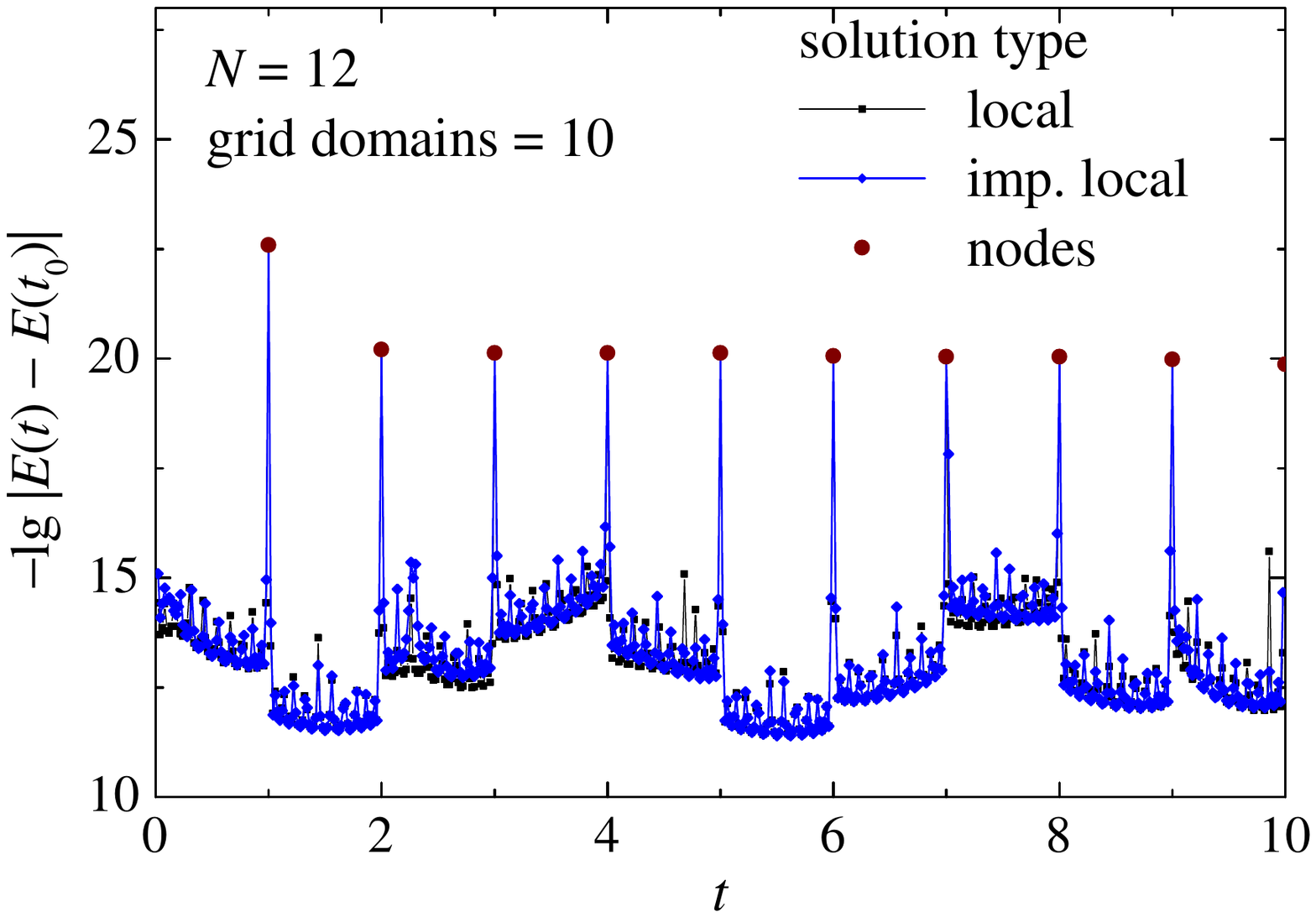}
\vspace{-8mm}\caption{\label{fig:econs_pend:a3}}
\end{subfigure}
\begin{subfigure}{0.24\textwidth}
\includegraphics[width=\textwidth]{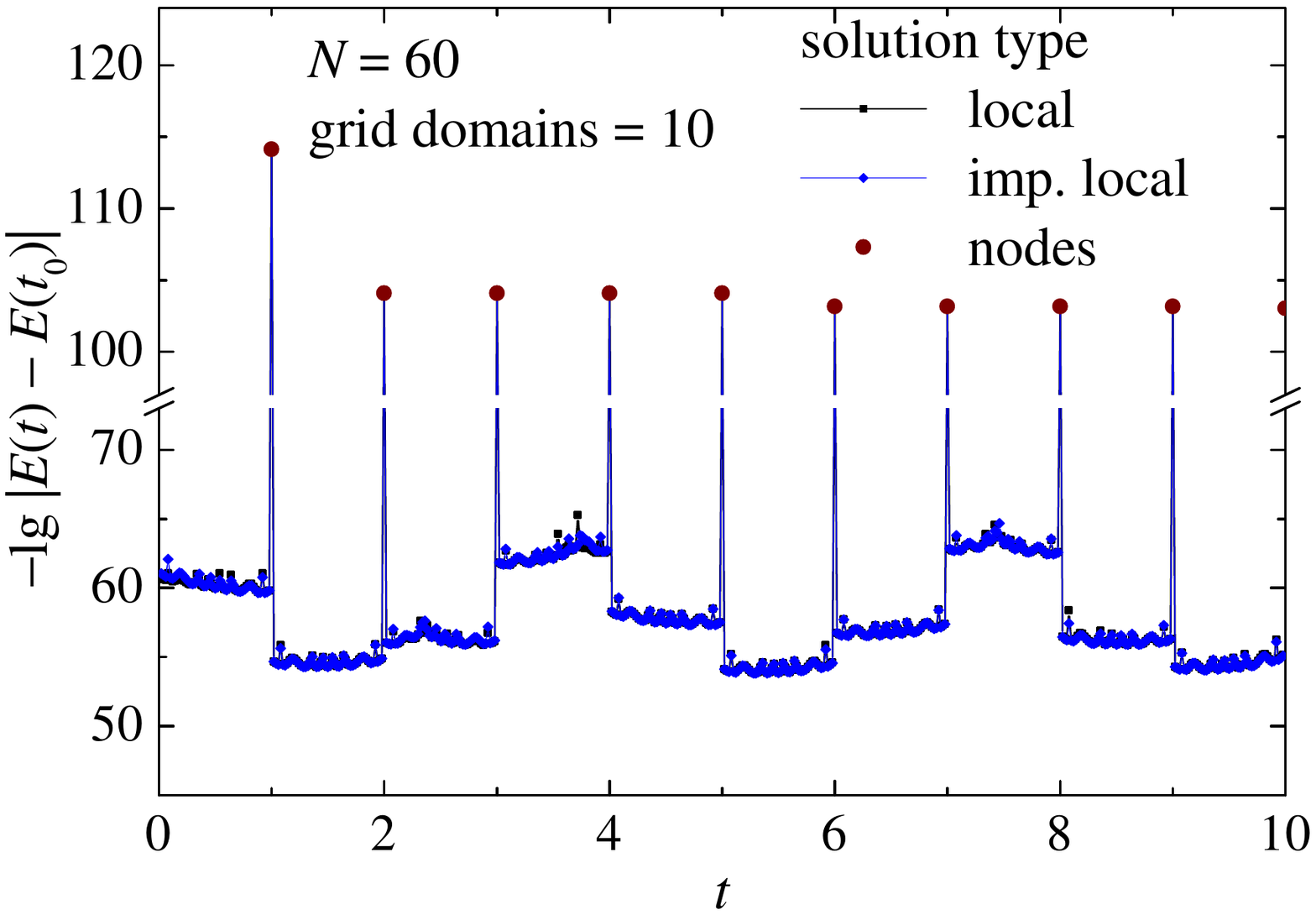}
\vspace{-8mm}\caption{\label{fig:econs_pend:a4}}
\end{subfigure}
\caption{%
Dependence of the negative logarithm of the error $-\lg|E(t)-E(t_{0})|$ of the energy conservation law of the numerical solution of the system (\ref{eq:pend_ode}) on the argument $t$ for solution at nodes $\mathbf{u}_{n}$, the local solution $\mathbf{u}_{L}(t)$ and the improved local solution $\mathbf{u}_{\rm IL}(t)$, obtained using polynomials with degrees $N = 1$~(\subref{fig:econs_pend:a1}), $N = 4$~(\subref{fig:econs_pend:a4}), $N = 12$~(\subref{fig:econs_pend:a3}), $N = 60$~(\subref{fig:econs_pend:a4}).
}
\label{fig:econs_pend}
\end{figure}

Analysis of the local error $\varepsilon$ (\ref{eq:eps_local_def}) presented in Fig.~\ref{fig:pend} (\subref{fig:pend:c1}, \subref{fig:pend:c2}, \subref{fig:pend:c3}, \subref{fig:pend:c4}) shows that it also exhibits behavior similar to the case of a harmonic oscillator in Fig.~\ref{fig:harm_osc} (\subref{fig:harm_osc:c1}, \subref{fig:harm_osc:c2}, \subref{fig:harm_osc:c3}, \subref{fig:harm_osc:c4}), with the difference in the local error between the local numerical solution $\mathbf{u}_{L}$ and the improved local numerical solution $\mathbf{u}_{\rm IL}$ being $1$--$3$ orders of magnitude. The local error $\varepsilon$ of the numerical solution $\mathbf{u}_{n}$ at the grid nodes $t_{n}$ is significantly smaller than the local errors $\varepsilon$ of the local solution $\mathbf{u}_{L}$ and the improved local solution $\mathbf{u}_{\rm IL}$, but not as significant as in the previous examples with linear problems --- amounting to approximately $1$ orders of magnitude in the case of polynomial degree $N = 1$, approximately $1$--$4$ orders of magnitude in the case of polynomial degree $N = 4$, approximately $11$--$13$ orders of magnitude in the case of polynomial degree $N = 12$ and approximately $50$--$60$ orders of magnitude in the case of polynomial degree $N = 60$ (for this purpose, breaks in the graph along the vertical axis are inserted in Fig.~\ref{fig:pend} (\subref{fig:pend:c4})).

The resulting dependencies of the global error $e$ (\ref{eq:eps_un_global_def}), (\ref{eq:eps_ul_global_def}) on the discretization step ${\Delta t}$, presented in Fig.~\ref{fig:pend} (\subref{fig:pend:d1}, \subref{fig:pend:d2}, \subref{fig:pend:d3}, \subref{fig:pend:d4}, \subref{fig:pend:e1}, \subref{fig:pend:e2}, \subref{fig:pend:e3}, \subref{fig:pend:e4}, \subref{fig:pend:f1}, \subref{fig:pend:f2}, \subref{fig:pend:f3}, \subref{fig:pend:f4}), demonstrate a quite acceptable quality of linear approximation for all studied polynomial degrees $N$ (Fig.~\ref{fig:pend} only shows results for polynomial degrees $N = 1$, $4$, $12$ and $60$). Some deviations and anomalies are observed for error $e^{l}_{L_{\infty}}$ in Fig.~\ref{fig:pend} (\subref{fig:pend:d1}) for polynomial degree $N = 1$ and for error $e^{\rm imp}_{L_{\infty}}$ in Fig.~\ref{fig:pend} (\subref{fig:pend:d2}) for polynomial degree $N = 4$. These are related to the approximate method for calculating the maximum local error $\varepsilon$ in the discretization ranges $\Omega_{n}$. In the remaining cases of the global errors $e$, small deviations from a highly accurate linear dependence are observed.

Based on the approximation of the obtained dependencies $e({\Delta t})$ in log-log scale by a linear function $\lg{e({\Delta t})} \propto p\cdot\lg{{\Delta t}}$, empirical convergence orders $p$ are calculated and presented in Table~\ref{tab:conv_orders_pend} for all polynomial degrees $N = 1, \ldots, 30$ and polynomial degrees $N$ up to $60$ with a step of $5$. Unlike the previous examples involving solving linear problems, in this example, the comparison of empirical convergence orders $p$ with expected values $p_{\rm th.}$ (\ref{eq:conv_ords_exp}) is not as good. For small polynomial degrees $N \leqslant 10$, the empirical convergence orders $p$ are somewhat lower than expected values $p_{\rm th.}$, but no more than $0.5$--$1.0$ (however, for polynomial degrees $N = 2$, the empirical convergence orders $p$ of the improved local solution $\mathbf{u}_{\rm IL}$ are even slightly higher than the expected value of $p_{\rm th.}^{\rm imp}(2) = 5$). As the polynomial degree $N$ increases, the difference often becomes more significant, sometimes reaching values greater than $1.0$. It should be noted that the empirical convergence orders $p$ for the improved local numerical solution $\mathbf{u}_{\rm IL}$ are usually higher than for the local numerical solution $\mathbf{u}_{L}$, with the difference sometimes exceeding $1.0$ (except in some cases of very high polynomial degrees $N$). These behavioral characteristics of the empirical convergence orders $p$ can be considered expected, and they were encountered and studied in~\cite{ader_dg_ode_jsc}.

\corrtext{Similar to the two previous Examples in Sections~\ref{sec:apps:lin_diss} and~\ref{sec:apps:harm_osc}, the system of equations (\ref{eq:pend_ode}) considered in this Example can be presented in a conservative form, which is characterized by the energy conservation law $E(t) = \mathrm{const}$~\cite{LandauMechanics, GoldsteinMechanics} of the following form
\begin{equation}\label{eq:econs_pend}
E(t) = \frac{\dot{\phi}^{2}(t)}{2} - \omega_{0}^{2}\cos\phi(t)= \mathrm{const} \equiv 
E(t_{0}) = \frac{\dot{\phi}^{2}(0)}{2} - \omega_{0}^{2}\cos\phi(0) = -\omega_{0}^{2}\cos\phi_{0},
\end{equation}
which represents the integral of motion.}

\corrtext{Fig.~\ref{fig:econs_pend} shows the dependence of the negative logarithm of the error $-\lg|E(t)-E(t_{0})|$ of the energy conservation law (\ref{eq:econs_pend}) of the numerical solution of the system (\ref{eq:pend_ode}) on the argument $t$ for solution at nodes $\mathbf{u}_{n}$, the local solution $\mathbf{u}_{L}(t)$ and the improved local solution $\mathbf{u}_{\rm IL}(t)$, obtained using polynomials with degrees $N = 1$~(\subref{fig:econs_pend:a1}), $N = 4$~(\subref{fig:econs_pend:a4}), $N = 12$~(\subref{fig:econs_pend:a3}), $N = 60$~(\subref{fig:econs_pend:a4}). The obtained results clearly demonstrate that the energy conservation law (\ref{eq:pend_ode}) is not strictly satisfied in the numerical solution. However, due to the possibility of achieving a high order $p$, even in the cases of degrees $N = 12$ and $60$, the error in the energy conservation law's fulfillment over the entire solution domain becomes smaller than the characteristic rounding error of double-precision floating-point numbers $\sim 10^{-15}$-$10^{-17}$. The presented results also clearly demonstrate that the accuracy of the fulfillment of the energy conservation law (\ref{eq:econs_pend}) for the improved local solution $\mathbf{u}_{\rm IL}(t)$ is significantly higher than for the local solution $\mathbf{u}_{L}(t)$. Therefore, it can be concluded that, similar to the two previous Examples in Sections~\ref{sec:apps:lin_diss} and~\ref{sec:apps:harm_osc}, despite the dissipative nature of the ADER-DG numerical method with local DG predictor, a sufficiently high degree $N$ can be chosen such that the accuracy of the fulfillment of the energy conservation law will be at or below the characteristic error of representing real numbers as floating-point numbers. The results obtained are qualitatively consistent with the results presented previously in Examples in Sections~\ref{sec:apps:lin_diss} and~\ref{sec:apps:harm_osc}. However, the quantitatively large error in the satisfiability of the energy conservation law in this example is due to the nonlinearity of the ODE system under study.}

The obtained results allowed to conclude that the ADER-DG numerical method with a local DG predictor provides a highly accurate numerical solution to the initial value problem for the nonlinear ODE system (\ref{eq:pend_ode}) presented in this example. The obtained results are in quite acceptable agreement with the theory developed above. The improved local numerical solution $\mathbf{u}_{\rm IL}$ demonstrates higher accuracy and a higher convergence order compared to the local numerical solution $\mathbf{u}_{L}$.

\subsection{Example 5: Bratu problem}
\label{sec:apps:bratu}

\begin{figure}[h!]
\captionsetup[subfigure]{%
	position=bottom,
	font+=smaller,
	textfont=normalfont,
	singlelinecheck=off,
	justification=raggedright
}
\centering
\begin{subfigure}{0.24\textwidth}
\includegraphics[width=\textwidth]{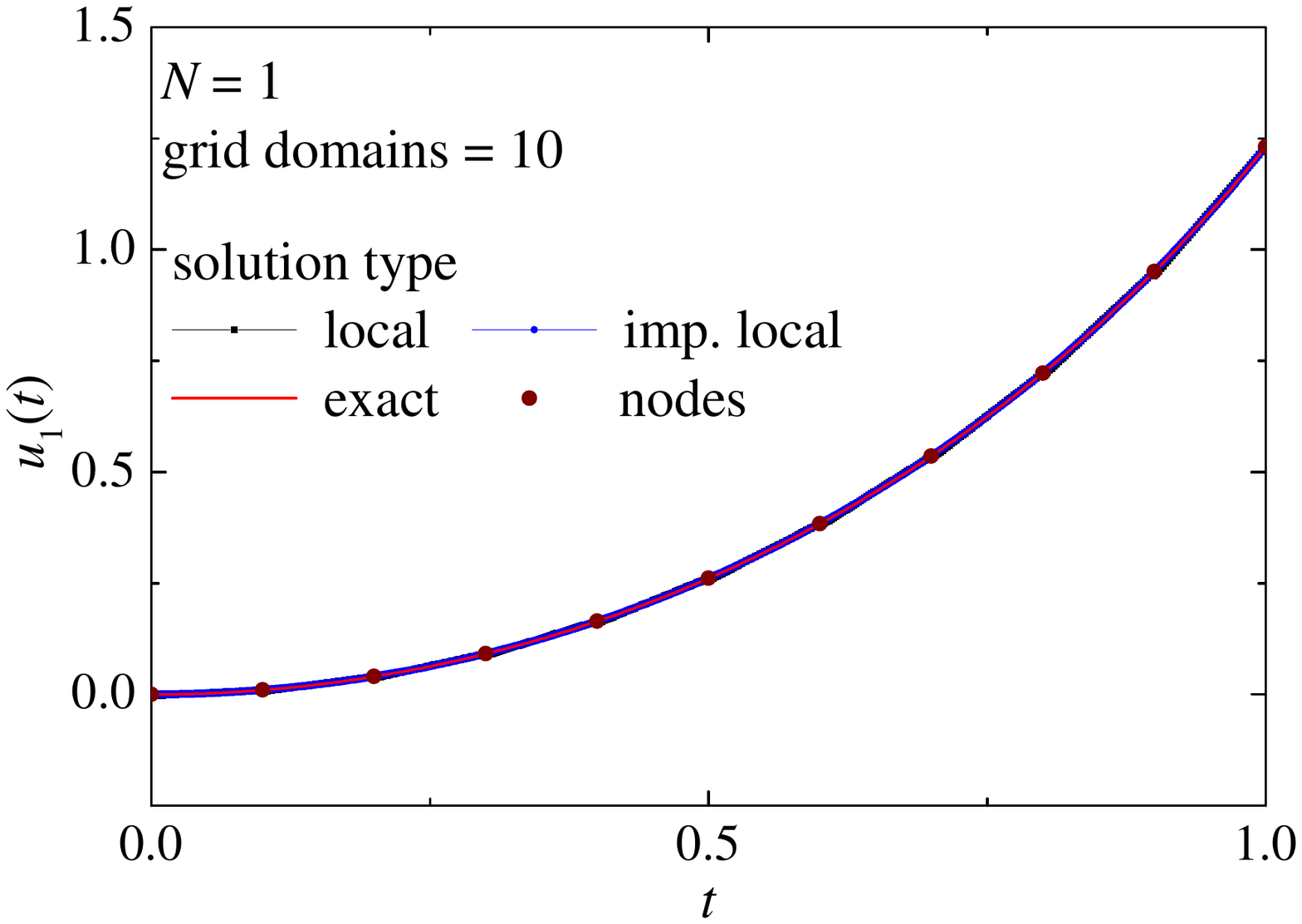}
\vspace{-8mm}\caption{\label{fig:bratu:a1}}
\end{subfigure}
\begin{subfigure}{0.24\textwidth}
\includegraphics[width=\textwidth]{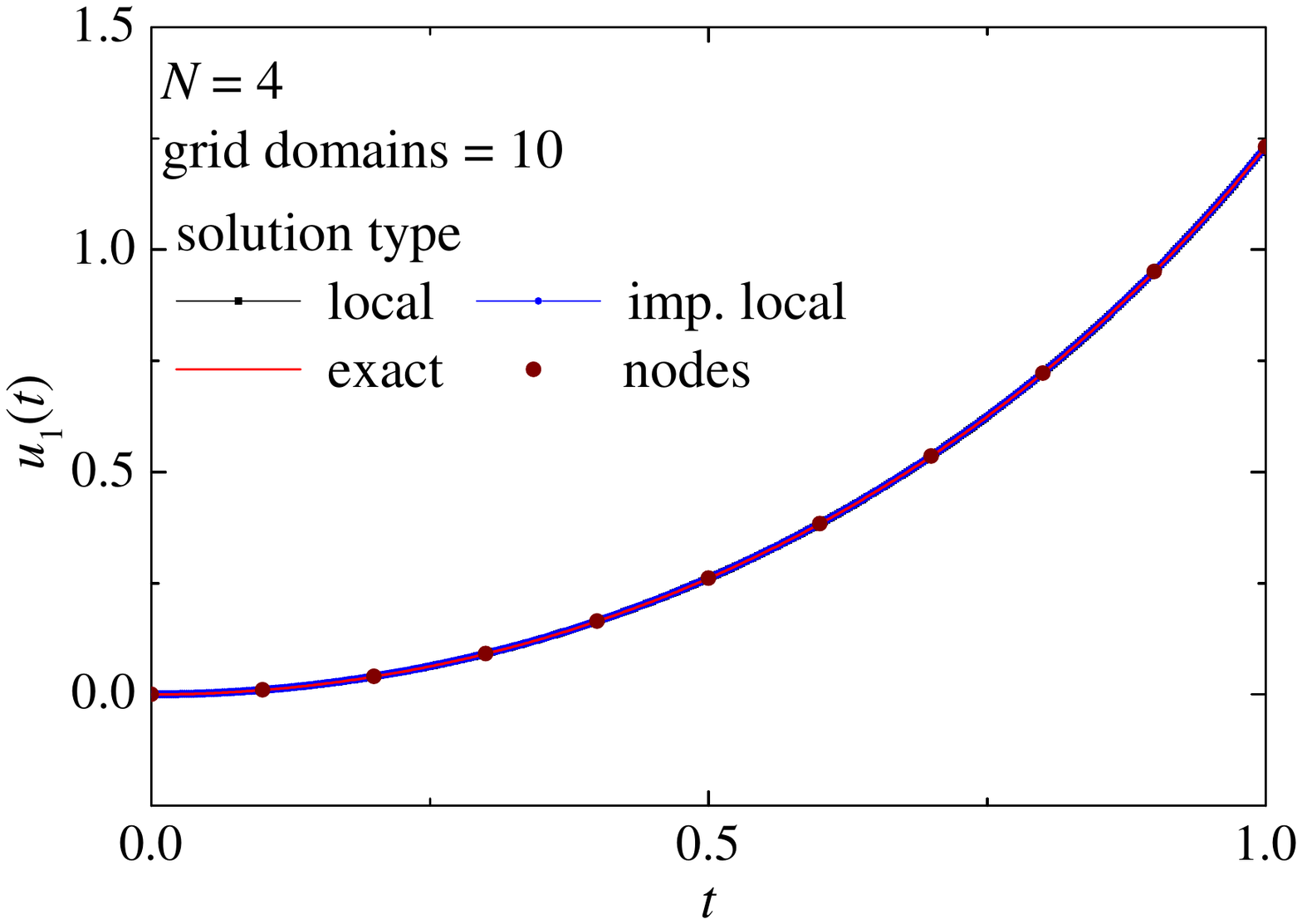}
\vspace{-8mm}\caption{\label{fig:bratu:a2}}
\end{subfigure}
\begin{subfigure}{0.24\textwidth}
\includegraphics[width=\textwidth]{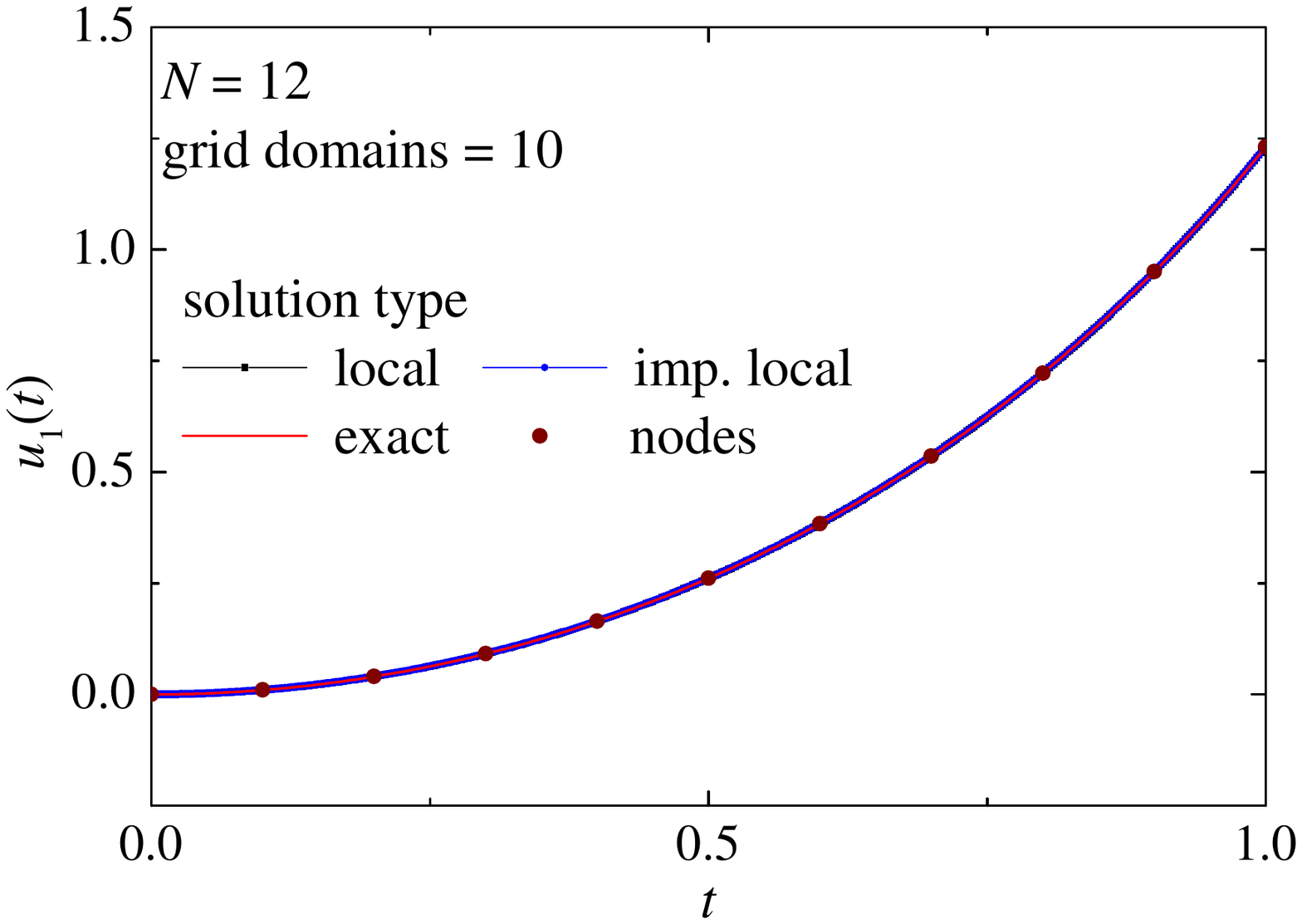}
\vspace{-8mm}\caption{\label{fig:bratu:a3}}
\end{subfigure}
\begin{subfigure}{0.24\textwidth}
\includegraphics[width=\textwidth]{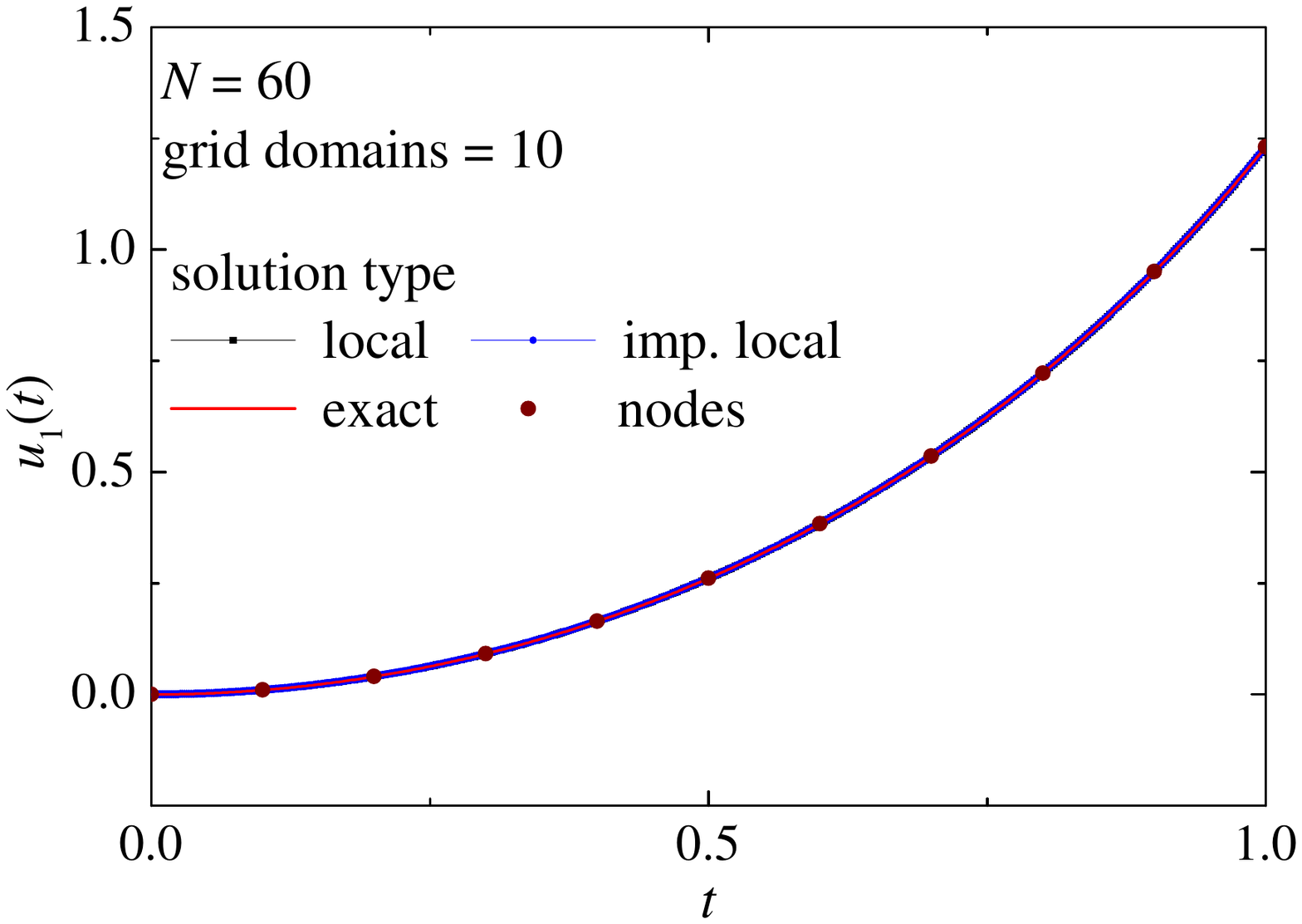}
\vspace{-8mm}\caption{\label{fig:bratu:a4}}
\end{subfigure}\\
\begin{subfigure}{0.24\textwidth}
\includegraphics[width=\textwidth]{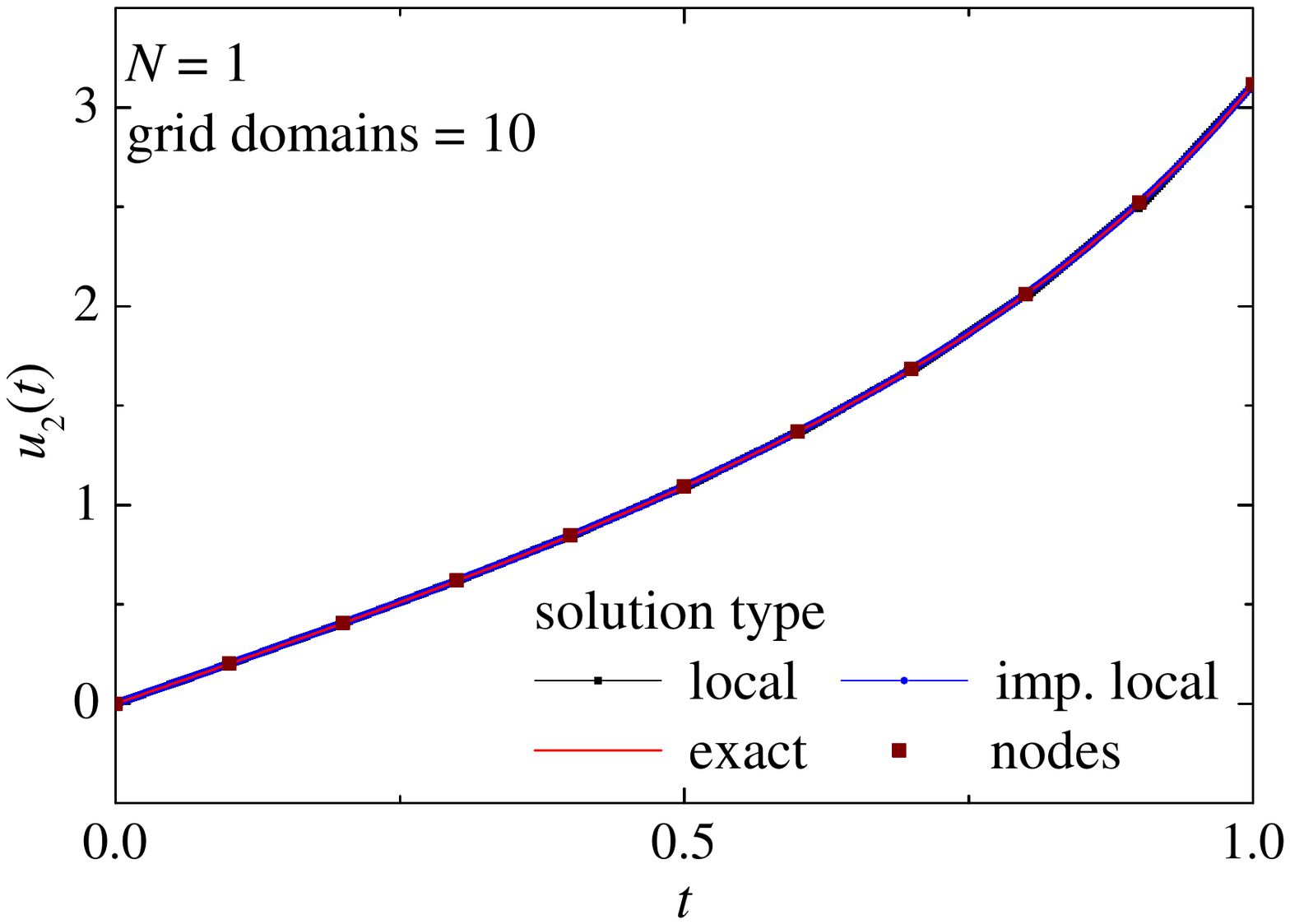}
\vspace{-8mm}\caption{\label{fig:bratu:b1}}
\end{subfigure}
\begin{subfigure}{0.24\textwidth}
\includegraphics[width=\textwidth]{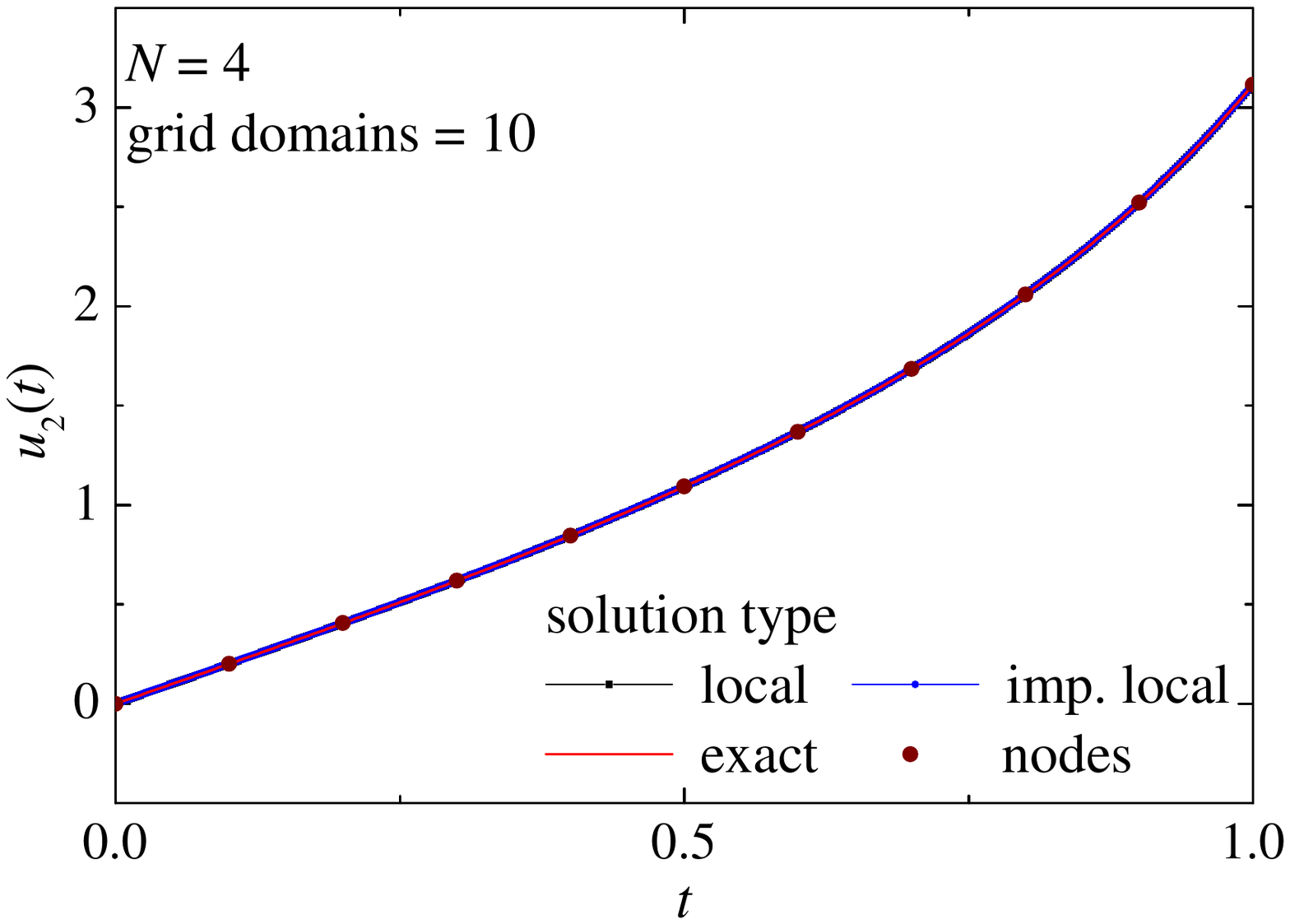}
\vspace{-8mm}\caption{\label{fig:bratu:b2}}
\end{subfigure}
\begin{subfigure}{0.24\textwidth}
\includegraphics[width=\textwidth]{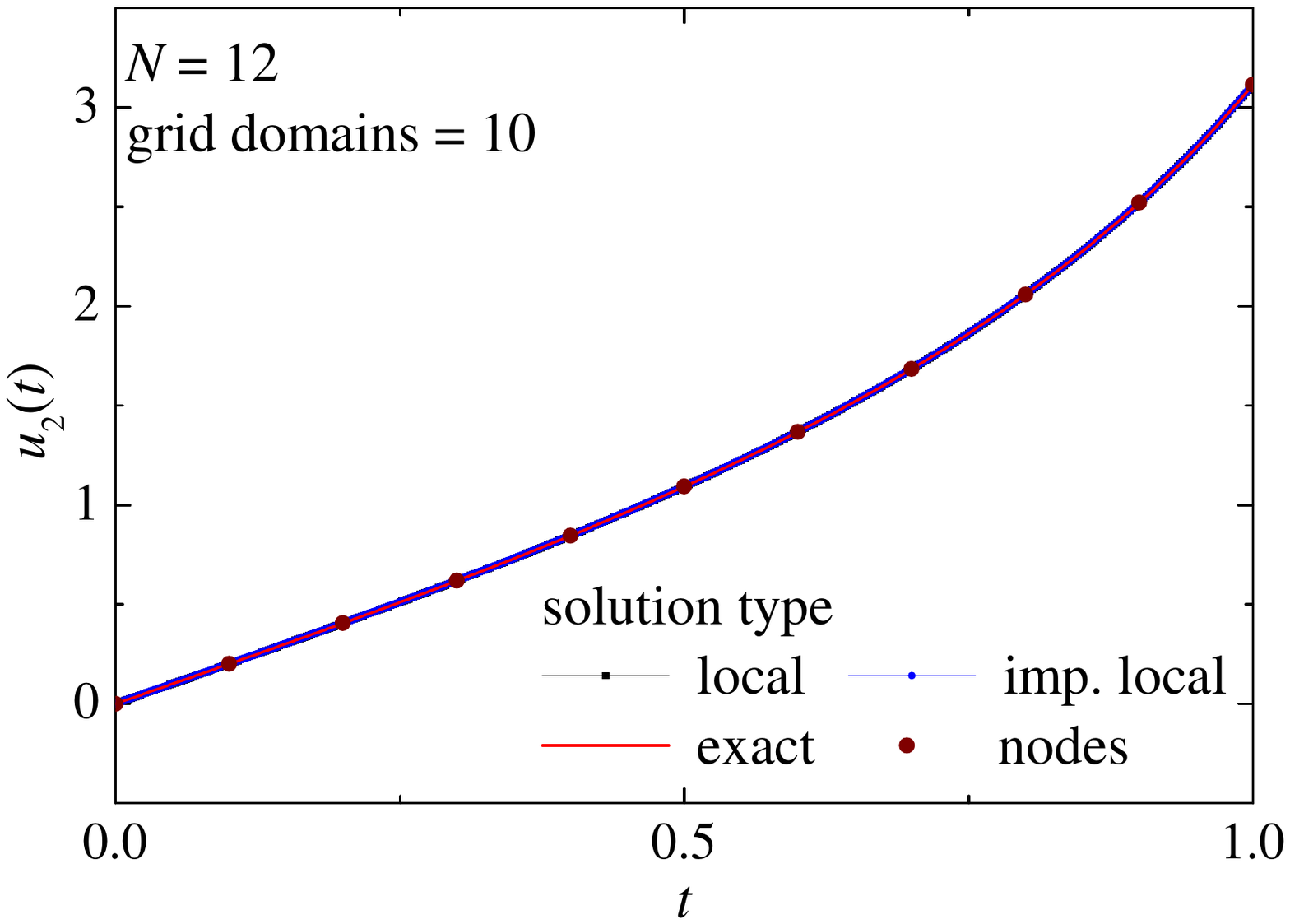}
\vspace{-8mm}\caption{\label{fig:bratu:b3}}
\end{subfigure}
\begin{subfigure}{0.24\textwidth}
\includegraphics[width=\textwidth]{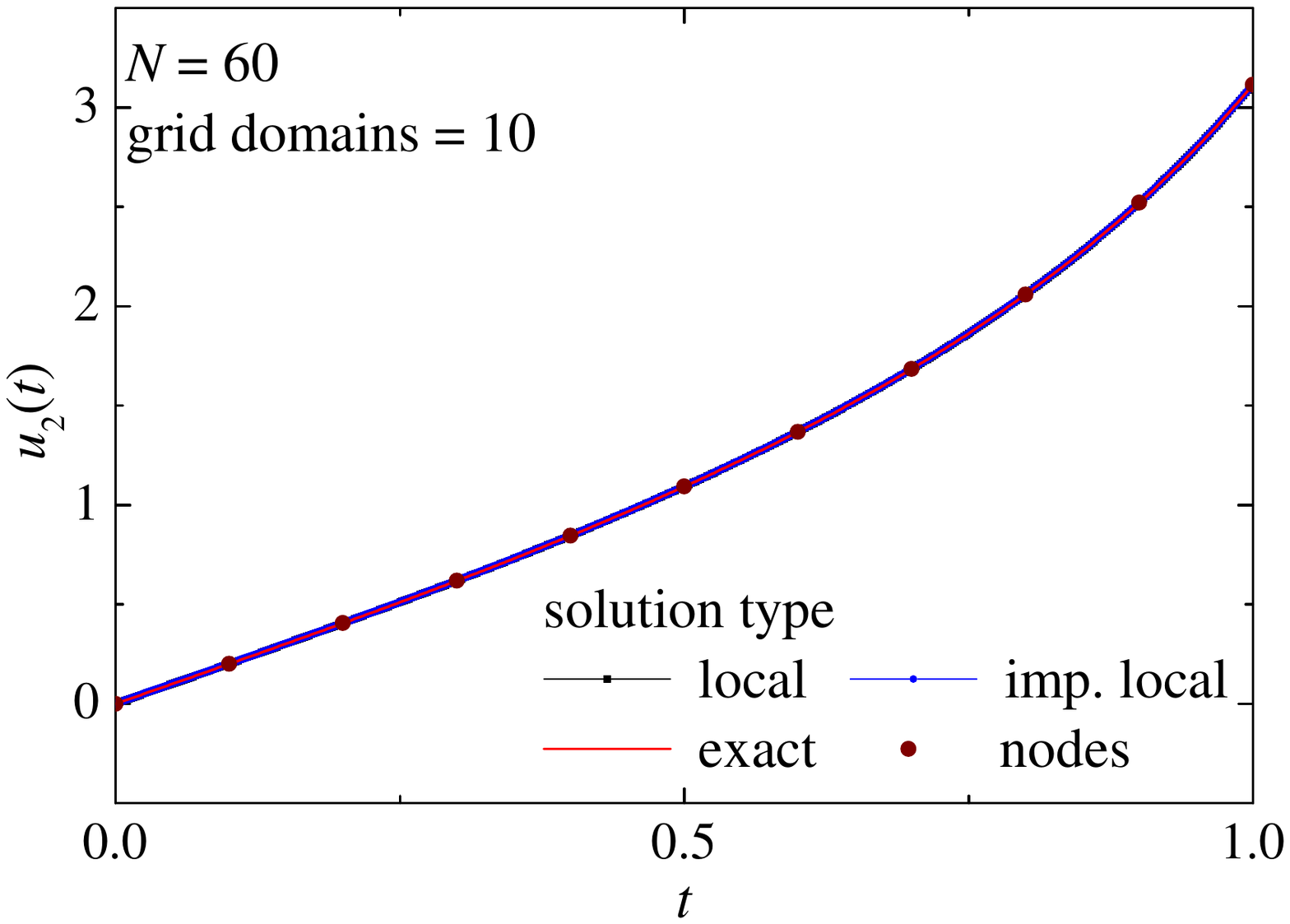}
\vspace{-8mm}\caption{\label{fig:bratu:b4}}
\end{subfigure}\\
\begin{subfigure}{0.24\textwidth}
\includegraphics[width=\textwidth]{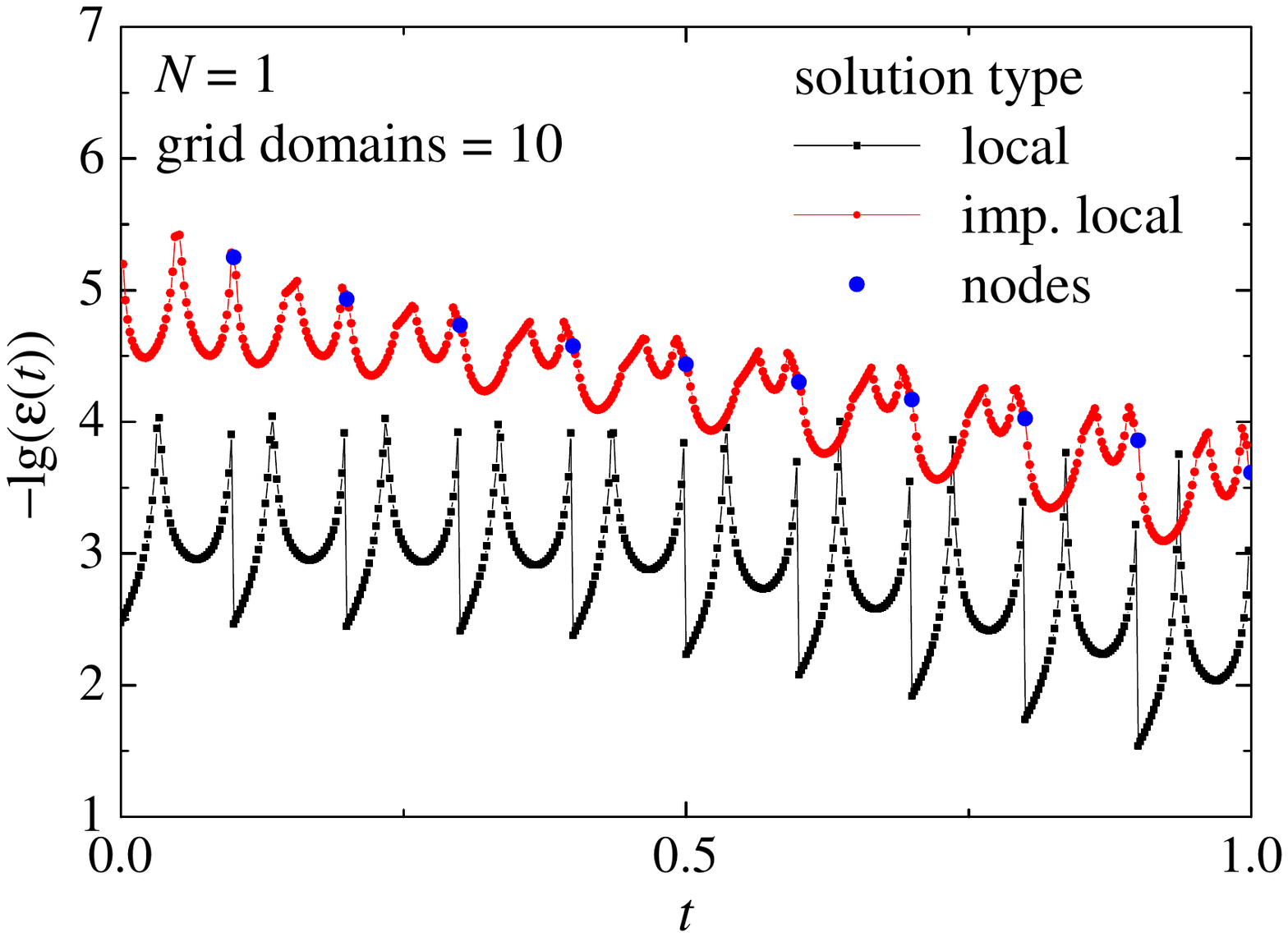}
\vspace{-8mm}\caption{\label{fig:bratu:c1}}
\end{subfigure}
\begin{subfigure}{0.24\textwidth}
\includegraphics[width=\textwidth]{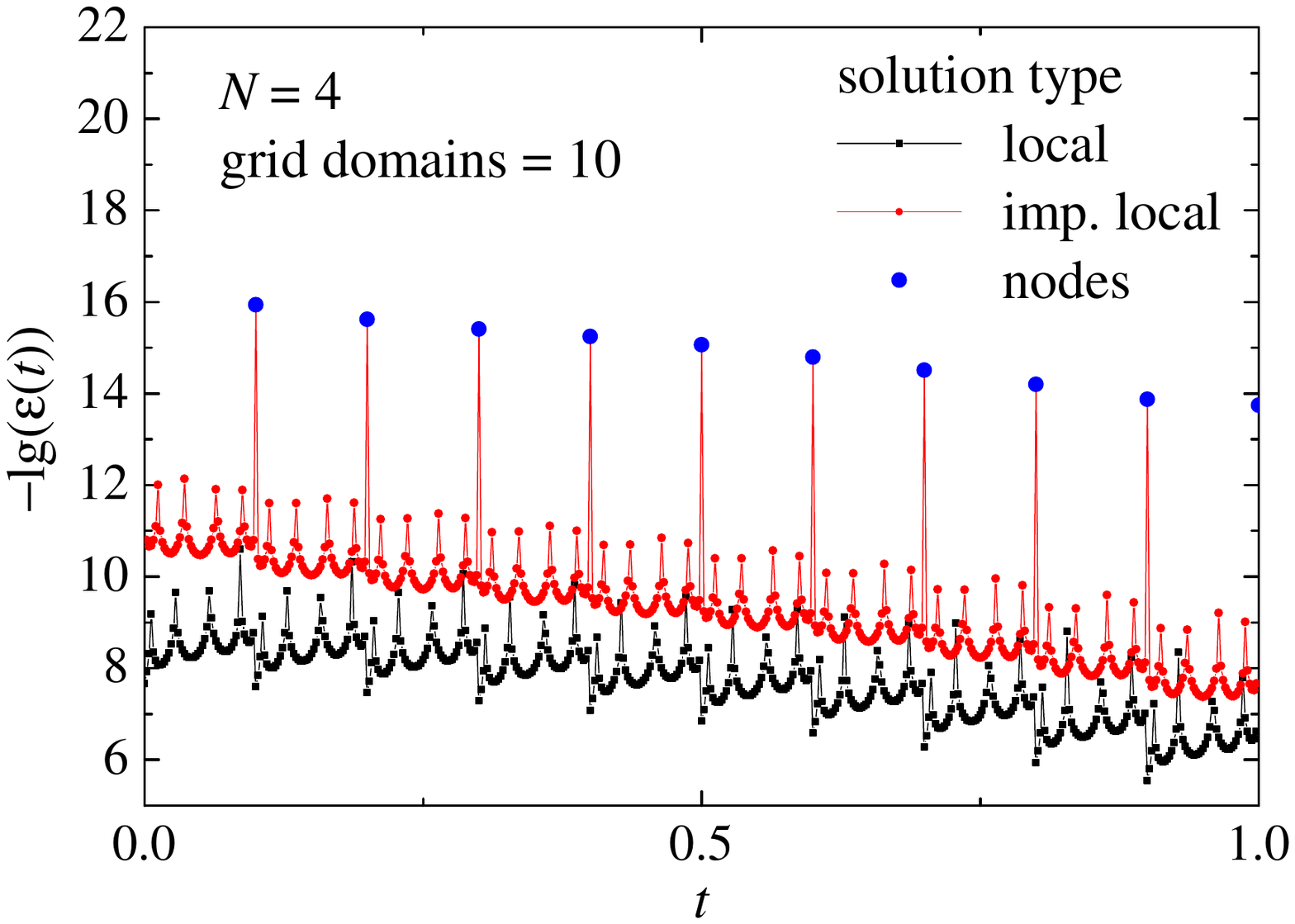}
\vspace{-8mm}\caption{\label{fig:bratu:c2}}
\end{subfigure}
\begin{subfigure}{0.24\textwidth}
\includegraphics[width=\textwidth]{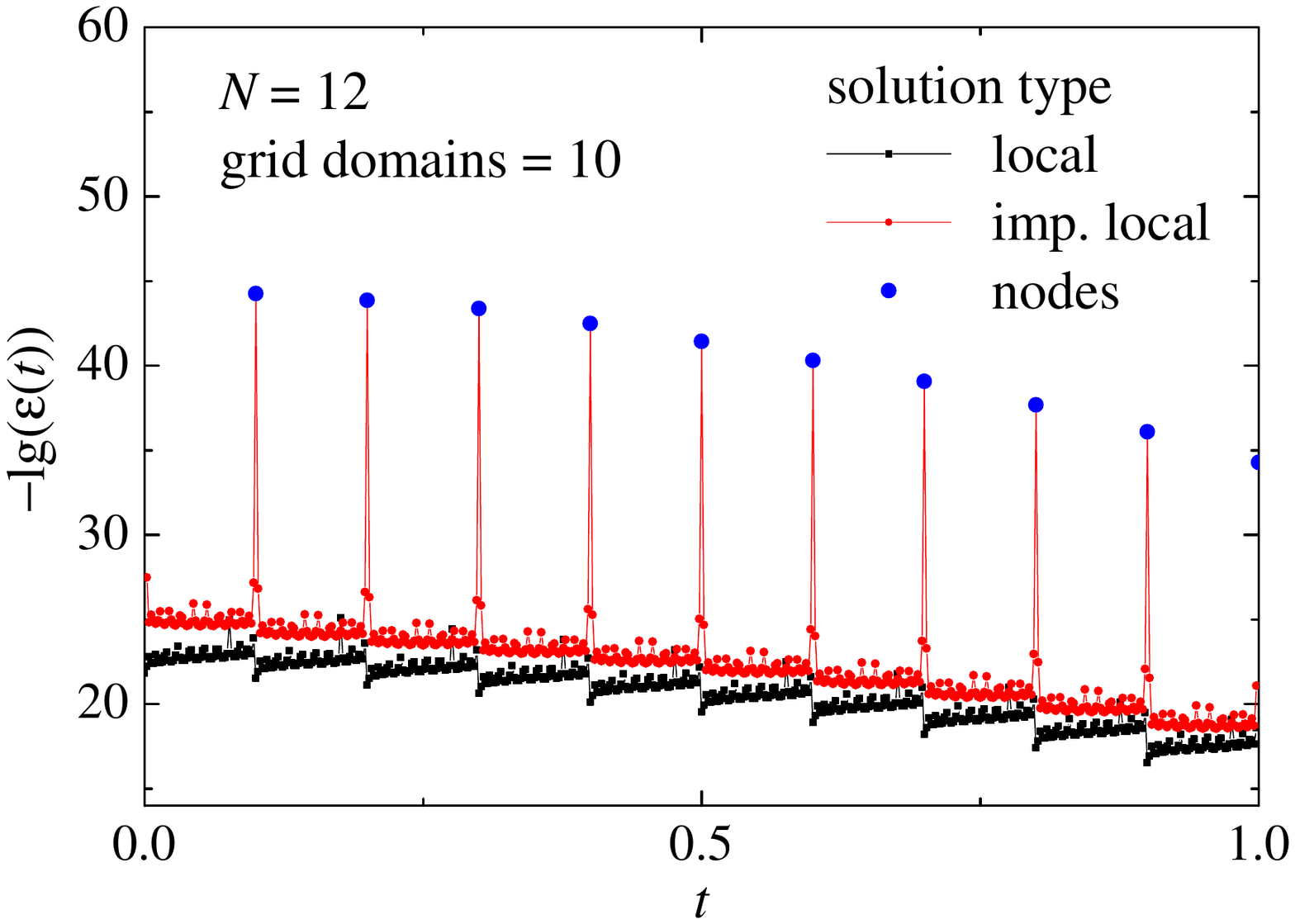}
\vspace{-8mm}\caption{\label{fig:bratu:c3}}
\end{subfigure}
\begin{subfigure}{0.24\textwidth}
\includegraphics[width=\textwidth]{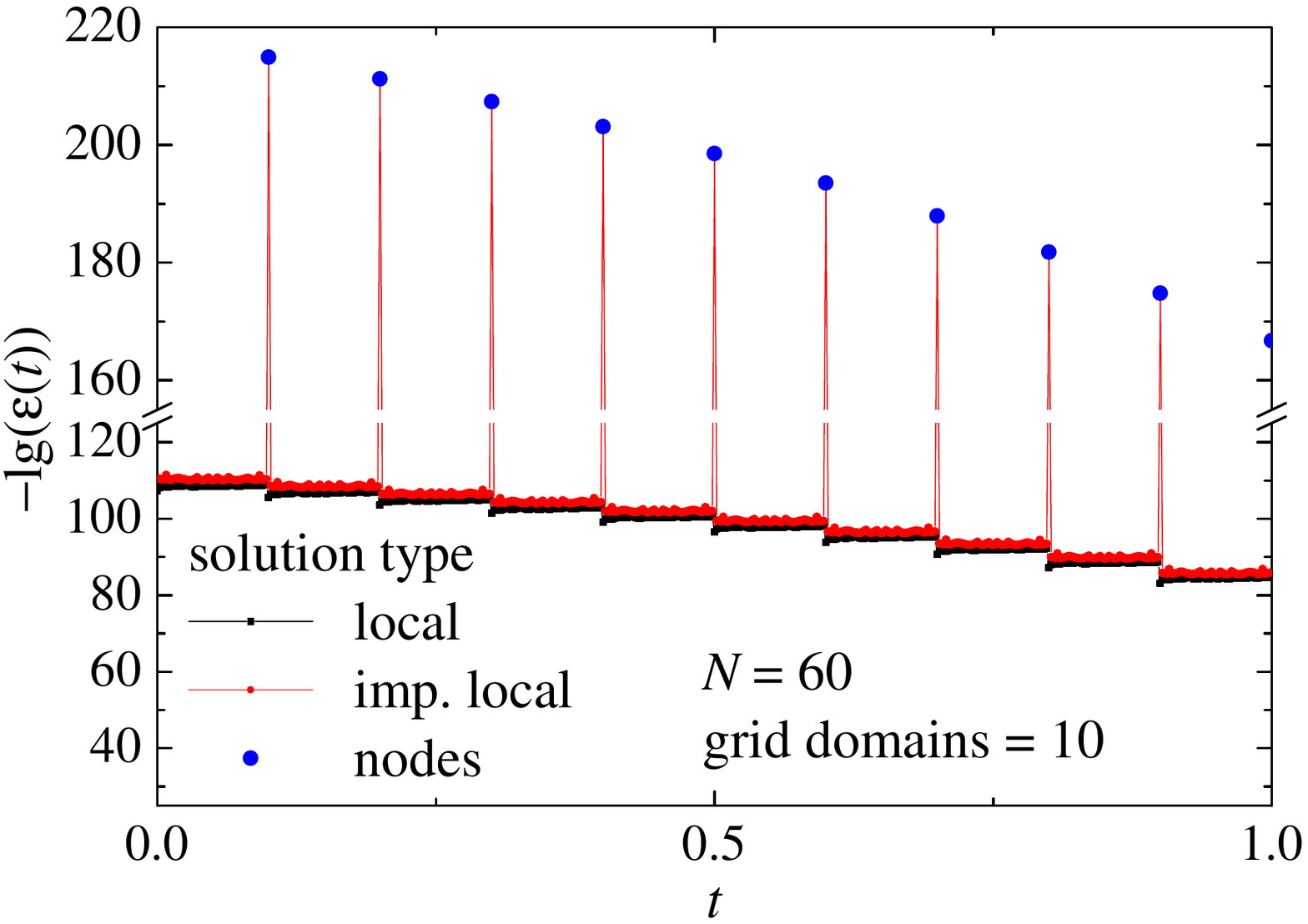}
\vspace{-8mm}\caption{\label{fig:bratu:c4}}
\end{subfigure}\\
\begin{subfigure}{0.24\textwidth}
\includegraphics[width=\textwidth]{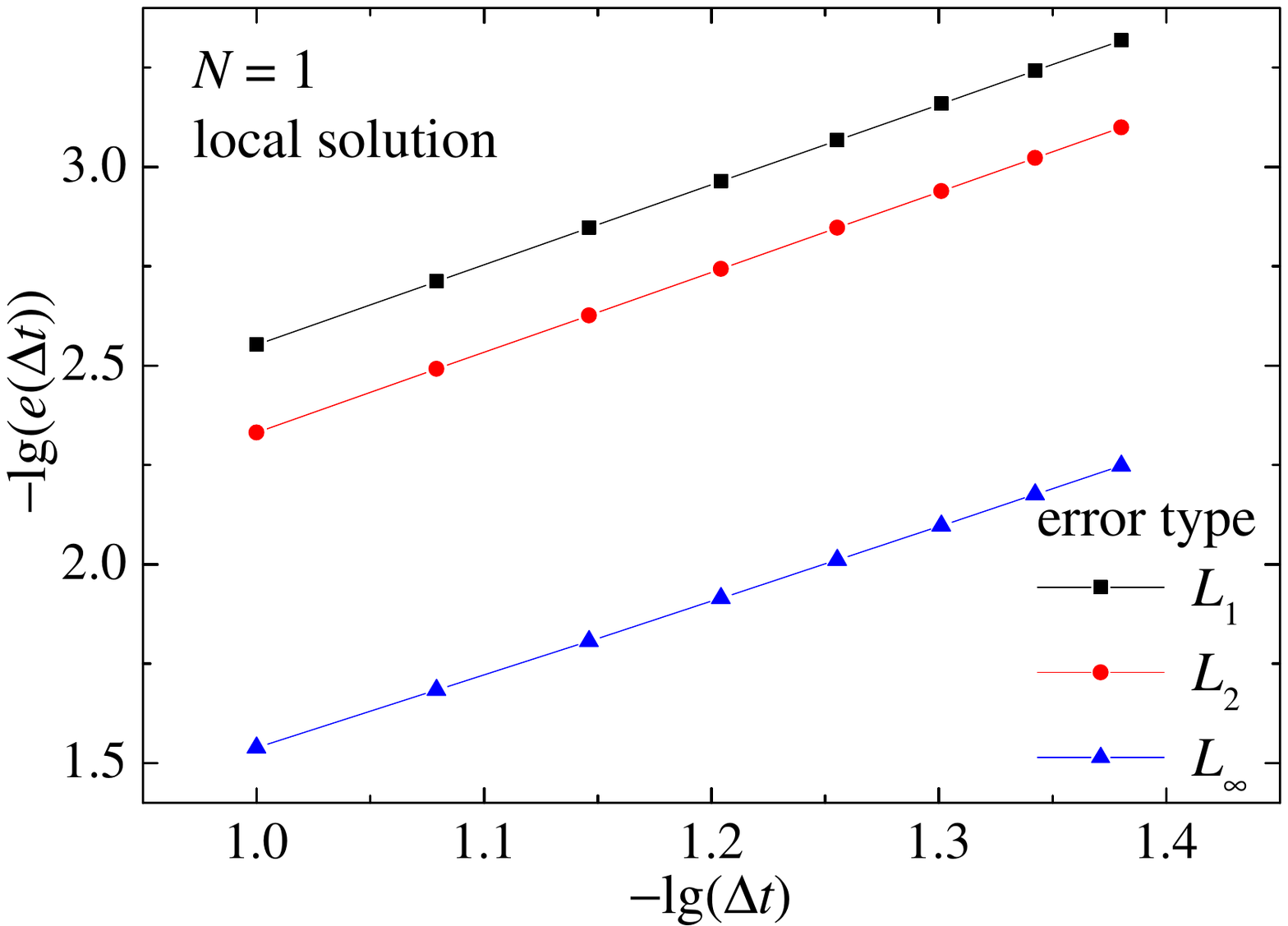}
\vspace{-8mm}\caption{\label{fig:bratu:d1}}
\end{subfigure}
\begin{subfigure}{0.24\textwidth}
\includegraphics[width=\textwidth]{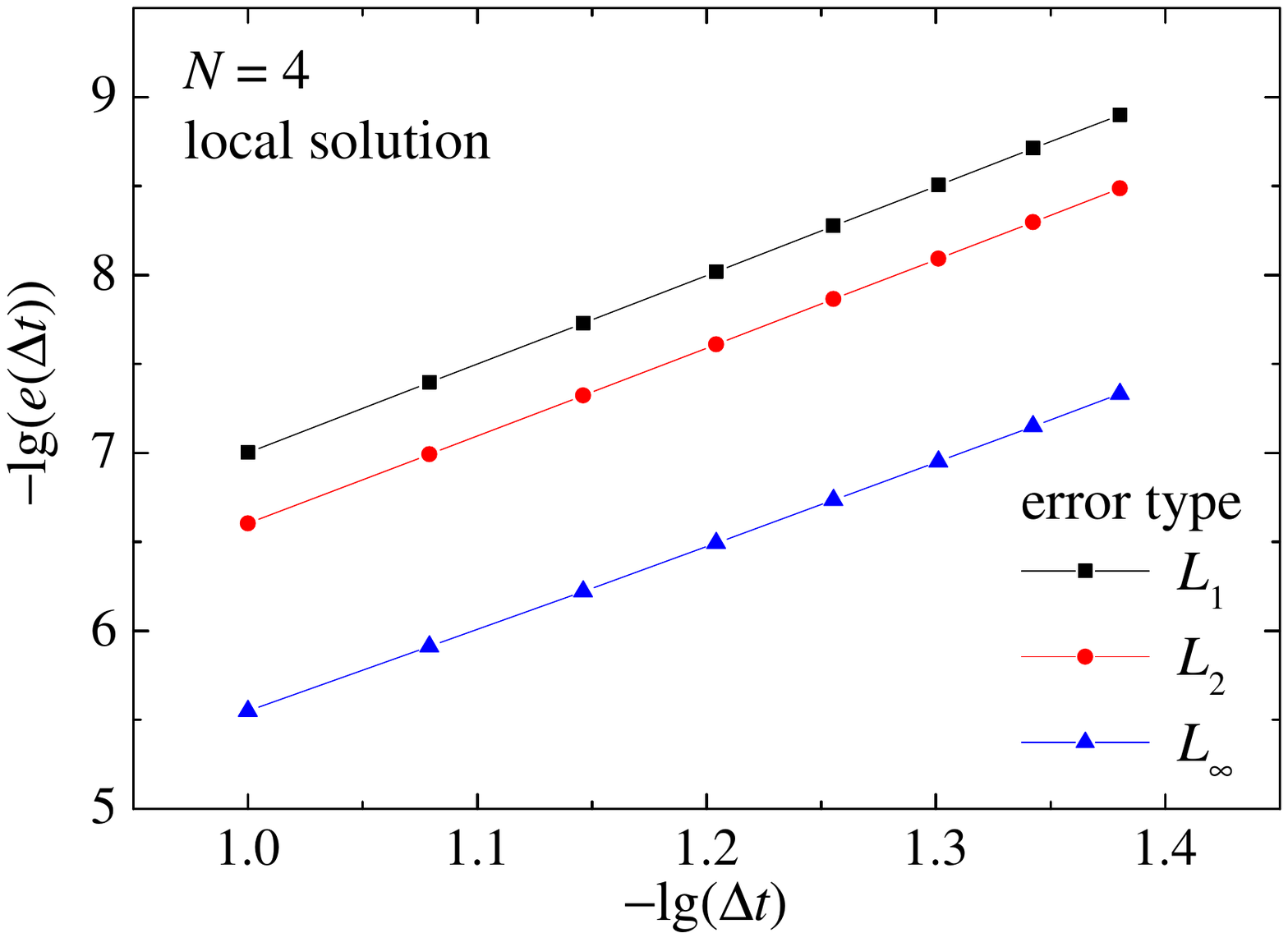}
\vspace{-8mm}\caption{\label{fig:bratu:d2}}
\end{subfigure}
\begin{subfigure}{0.24\textwidth}
\includegraphics[width=\textwidth]{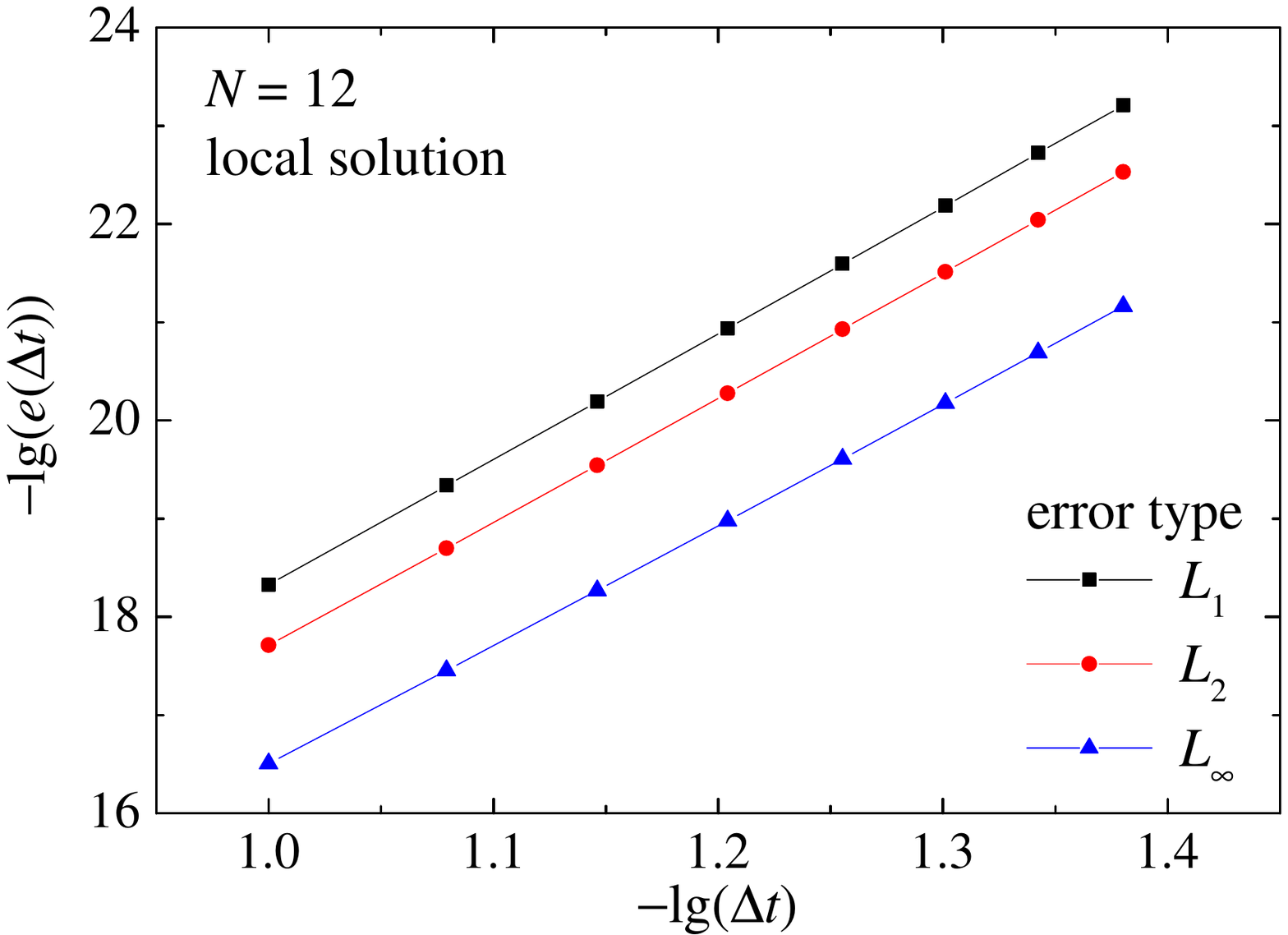}
\vspace{-8mm}\caption{\label{fig:bratu:d3}}
\end{subfigure}
\begin{subfigure}{0.24\textwidth}
\includegraphics[width=\textwidth]{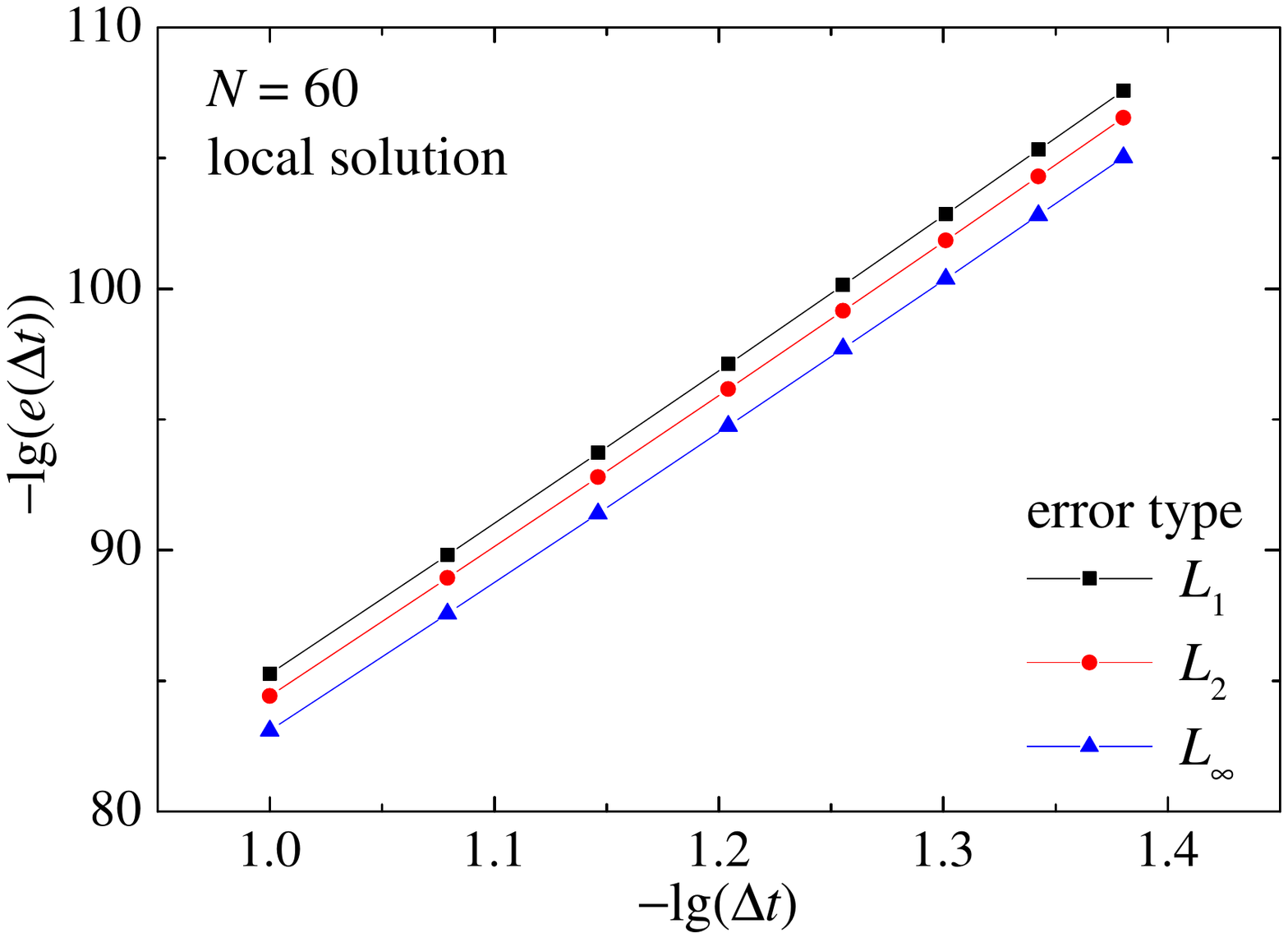}
\vspace{-8mm}\caption{\label{fig:bratu:d4}}
\end{subfigure}\\
\begin{subfigure}{0.24\textwidth}
\includegraphics[width=\textwidth]{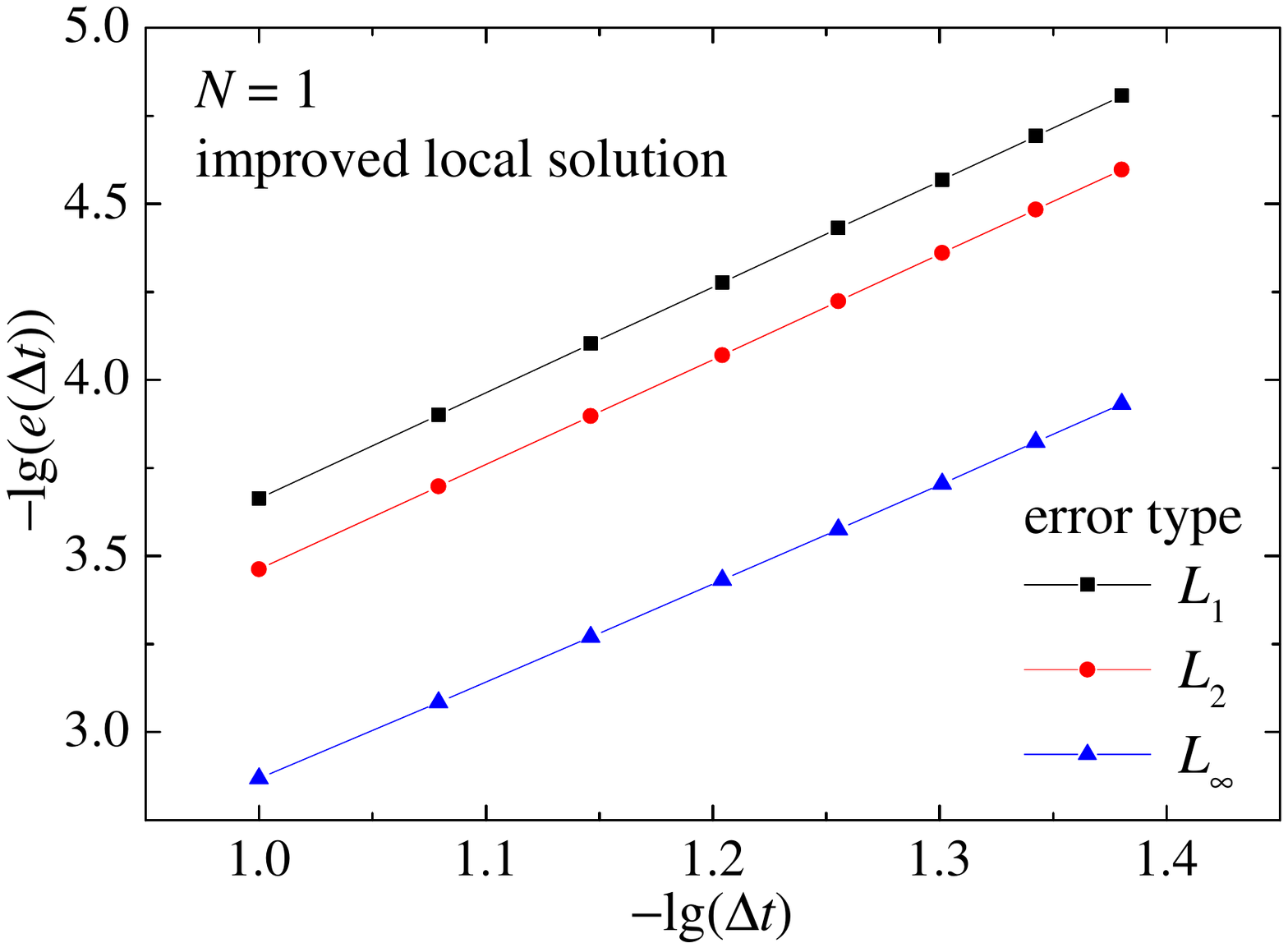}
\vspace{-8mm}\caption{\label{fig:bratu:e1}}
\end{subfigure}
\begin{subfigure}{0.24\textwidth}
\includegraphics[width=\textwidth]{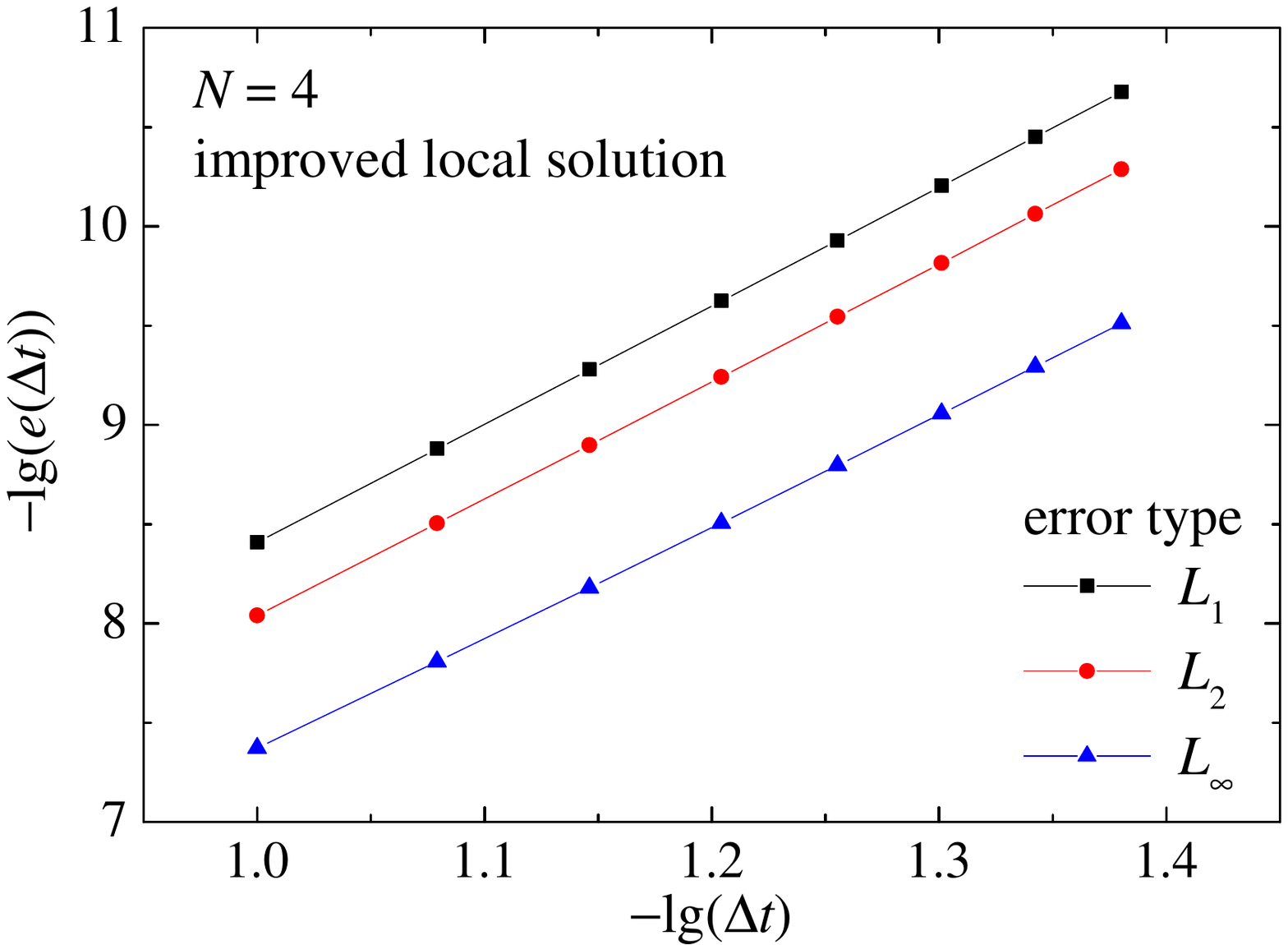}
\vspace{-8mm}\caption{\label{fig:bratu:e2}}
\end{subfigure}
\begin{subfigure}{0.24\textwidth}
\includegraphics[width=\textwidth]{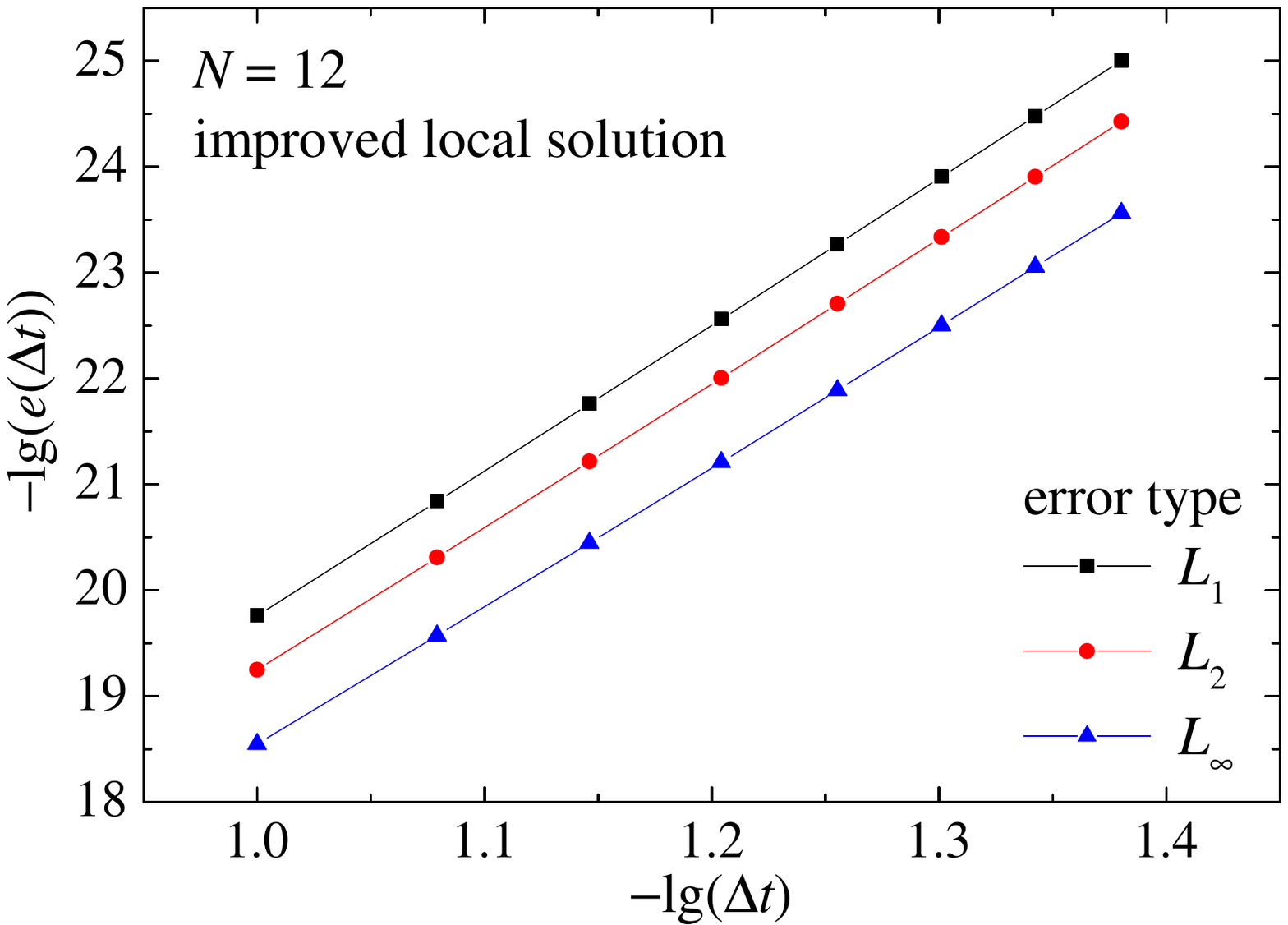}
\vspace{-8mm}\caption{\label{fig:bratu:e3}}
\end{subfigure}
\begin{subfigure}{0.24\textwidth}
\includegraphics[width=\textwidth]{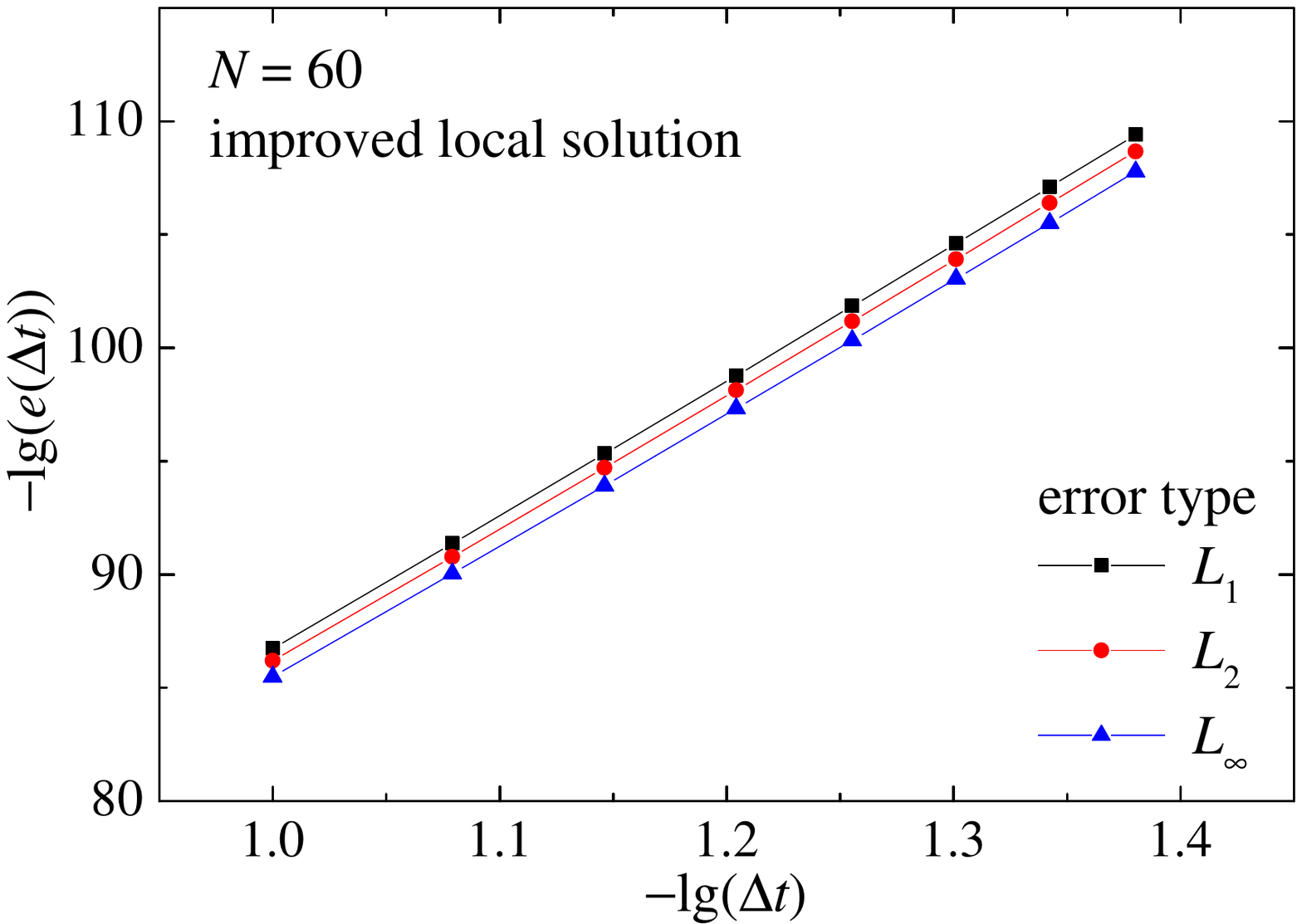}
\vspace{-8mm}\caption{\label{fig:bratu:e4}}
\end{subfigure}\\
\begin{subfigure}{0.24\textwidth}
\includegraphics[width=\textwidth]{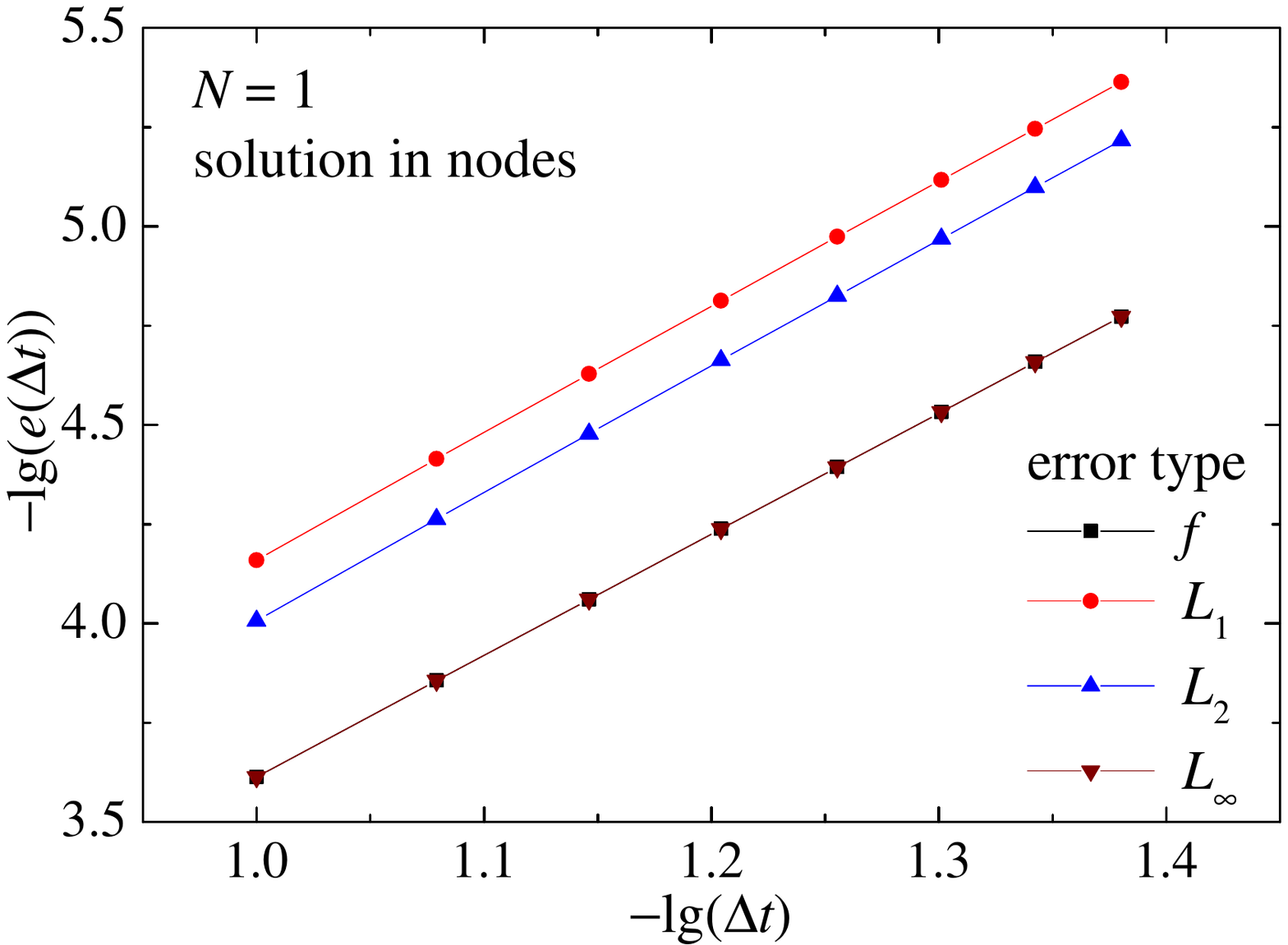}
\vspace{-8mm}\caption{\label{fig:bratu:f1}}
\end{subfigure}
\begin{subfigure}{0.24\textwidth}
\includegraphics[width=\textwidth]{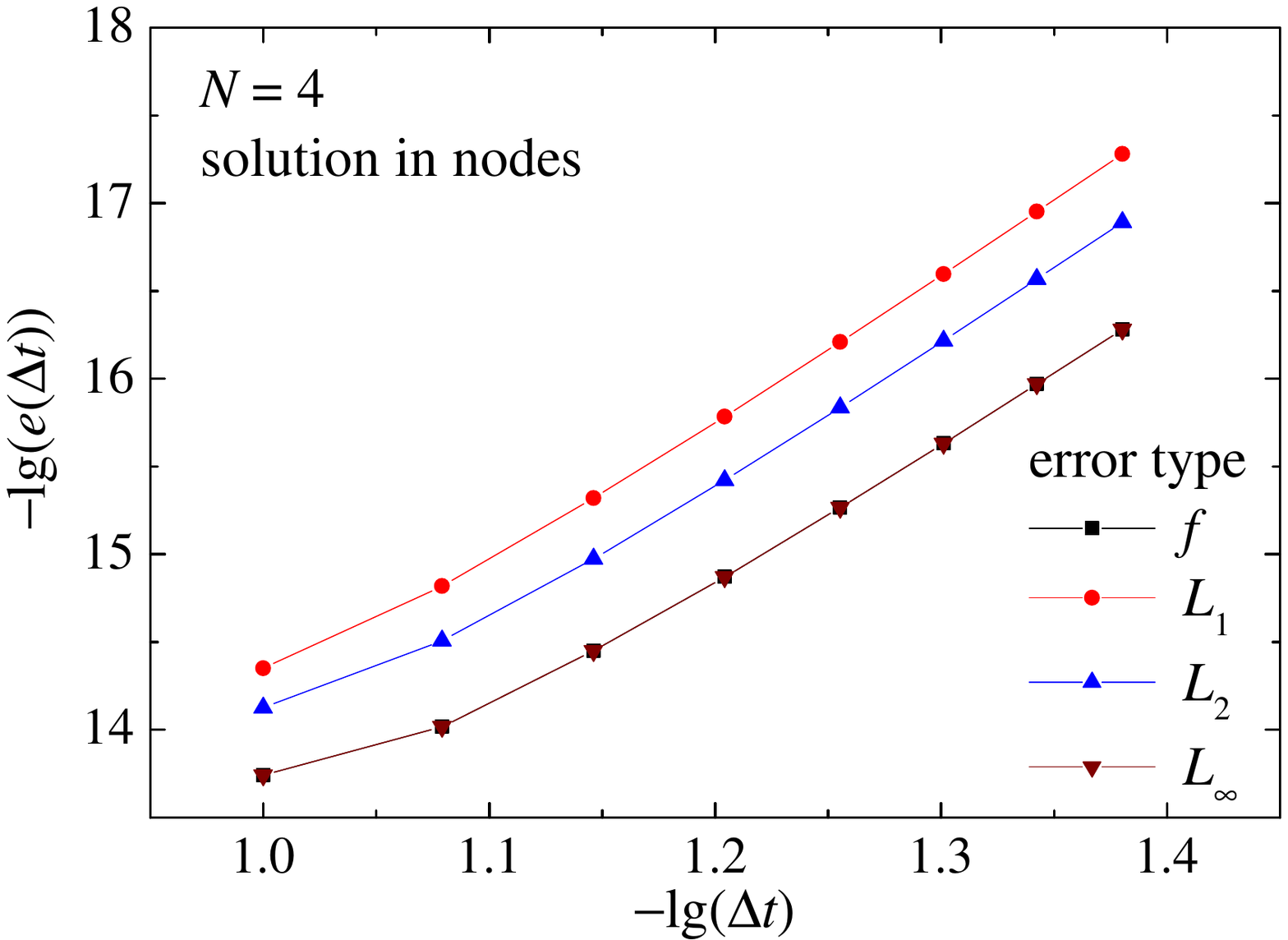}
\vspace{-8mm}\caption{\label{fig:bratu:f2}}
\end{subfigure}
\begin{subfigure}{0.24\textwidth}
\includegraphics[width=\textwidth]{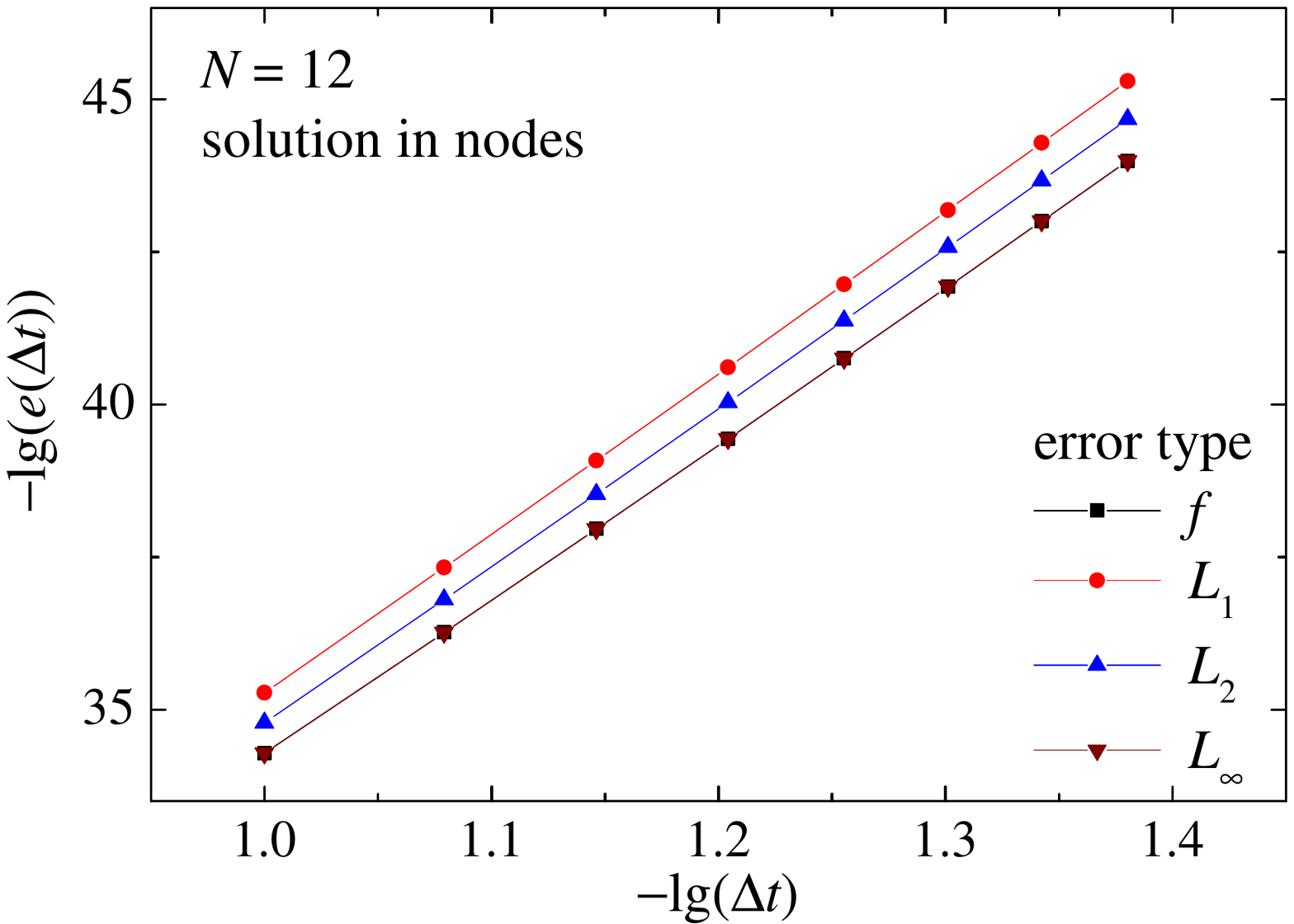}
\vspace{-8mm}\caption{\label{fig:bratu:f3}}
\end{subfigure}
\begin{subfigure}{0.24\textwidth}
\includegraphics[width=\textwidth]{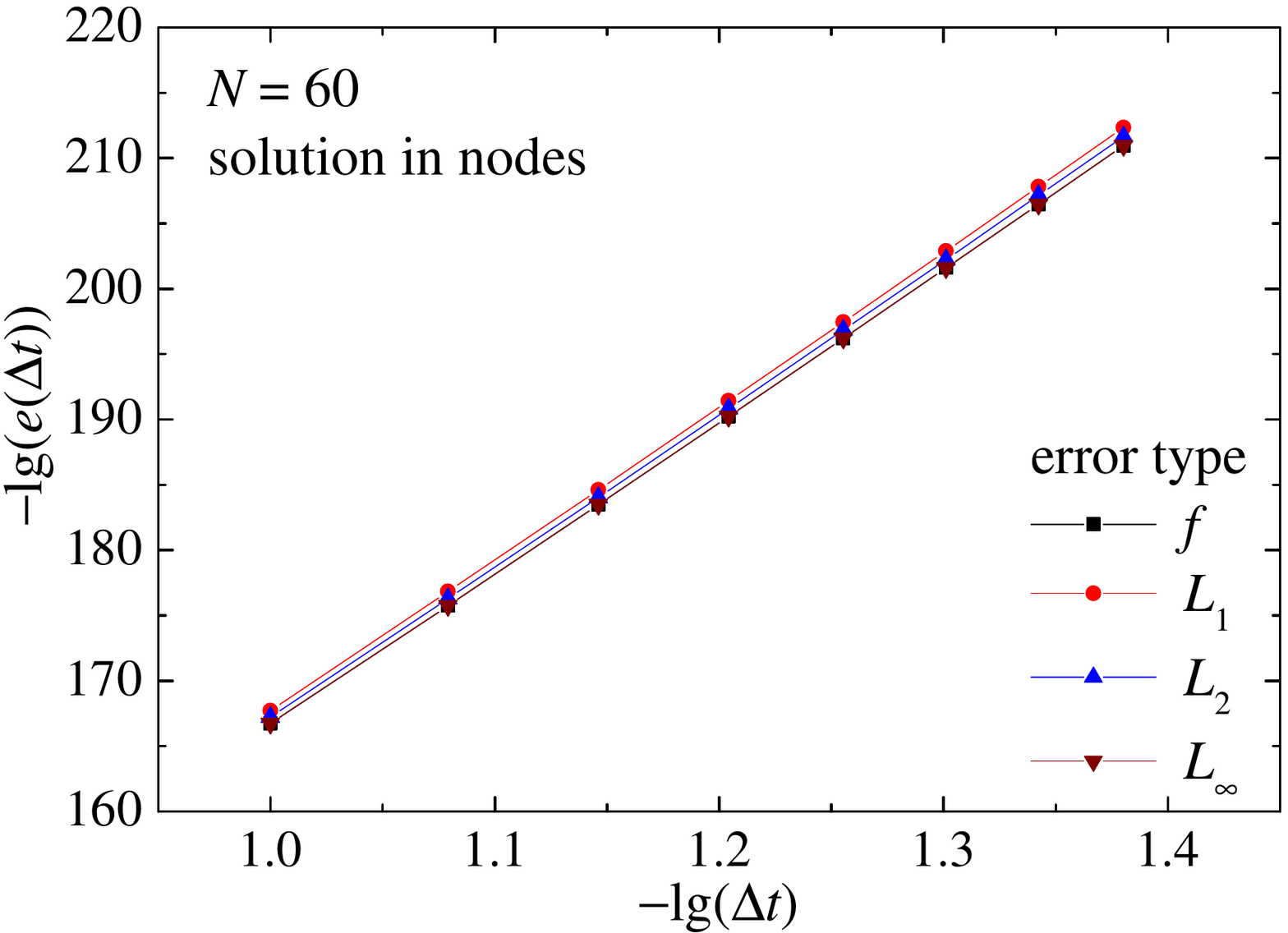}
\vspace{-8mm}\caption{\label{fig:bratu:f4}}
\end{subfigure}\\
\caption{%
Numerical solution of the system (\ref{eq:bratu_ode}). Comparison of the solution at nodes $\mathbf{u}_{n}$, the local solution $\mathbf{u}_{L}(t)$, the improved local solution $\mathbf{u}_{\rm IL}(t)$ and the exact solution $\mathbf{u}^{\rm ex}(t)$ (\ref{eq:bratu_sol_ex}) for components $u_{1} \equiv x$ (\subref{fig:bratu:a1}, \subref{fig:bratu:a2}, \subref{fig:bratu:a3}, \subref{fig:bratu:a4}) and $u_{2} \equiv \dot{x}$ (\subref{fig:bratu:b1}, \subref{fig:bratu:b2}, \subref{fig:bratu:b3}, \subref{fig:bratu:b4}), the errors $\varepsilon(t)$ (\subref{fig:bratu:c1}, \subref{fig:bratu:c2}, \subref{fig:bratu:c3}, \subref{fig:bratu:c4}), obtained using polynomials with degrees $N = 1$ (\subref{fig:bratu:a1}, \subref{fig:bratu:b1}, \subref{fig:bratu:c1}), $N = 4$ (\subref{fig:bratu:a2}, \subref{fig:bratu:b2}, \subref{fig:bratu:c2}), $N = 12$ (\subref{fig:bratu:a3}, \subref{fig:bratu:b3}, \subref{fig:bratu:c3}) and $N = 60$ (\subref{fig:bratu:a4}, \subref{fig:bratu:b4}, \subref{fig:bratu:c4}). Log-log plot of the debratuence of the global error for the local solution $e^{l}$ (\subref{fig:bratu:d1}, \subref{fig:bratu:d2}, \subref{fig:bratu:d3}, \subref{fig:bratu:d4}), the improved local solution $e^{\rm imp}$ (\subref{fig:bratu:e1}, \subref{fig:bratu:e2}, \subref{fig:bratu:e3}, \subref{fig:bratu:e4}) and the solution at nodes $e^{n}$ (\subref{fig:bratu:f1}, \subref{fig:bratu:f2}, \subref{fig:bratu:f3}, \subref{fig:bratu:f4}) on the discretization step $\mathrm{\Delta}t$, obtained in the $f$-norm and norms $L_{1}$, $L_{2}$ and $L_{\infty}$, obtained using polynomials with degrees $N = 1$ (\subref{fig:bratu:d1}, \subref{fig:bratu:e1}, \subref{fig:bratu:f1}), $N = 4$ (\subref{fig:bratu:d2}, \subref{fig:bratu:e2}, \subref{fig:bratu:f2}), $N = 12$ (\subref{fig:bratu:d3}, \subref{fig:bratu:e3}, \subref{fig:bratu:f3}) and $N = 60$ (\subref{fig:bratu:d4}, \subref{fig:bratu:e4}, \subref{fig:bratu:f4}).
}
\label{fig:bratu}
\end{figure}

The fifth and final example of applying the ADER-DG numerical method with a local DG predictor to solving the initial value problem for the ODE system (\ref{eq:ivp_ode_diff_src}) is the Bratu problem:
\begin{equation}\label{eq:bratu_ode}
\ddot{x} - 2\exp(x) = 0,\quad
x(0) = 0,\quad \dot{x}(0) = 0,\quad
t \in [0,\, 1],
\end{equation}
for which the solution vectors in the original notation of the ODE system (\ref{eq:ivp_ode_diff_src}), with $D = 2$, are chosen in the form $\mathbf{u}(t) = [x(t)\ \dot{x}(t)]^{T}$. The exact analytical solution is as follows:
\begin{equation}\label{eq:bratu_sol_ex}
x^{\rm ex}(t) = -2\ln(\cos(t)),\quad
\dot{x}^{\rm ex}(t) = 2\tan(t).
\end{equation}
The main difference between this example and the first three examples in the Subsections~\ref{sec:apps:exp_diss},~\ref{sec:apps:lin_diss} and~\ref{sec:apps:harm_osc}, which studied linear ODE systems, and the fourth example in the Subsection~\ref{sec:apps:pend}, which studied a nonlinear ODE system (\ref{eq:pend_ode}), is the violation of the Lipschitz condition (\ref{eq:lip_cond_f}) in the vector argument $\mathbf{u}$ for the function $\mathbf{F}(\mathbf{u},\, t)$ on the right-hand side of the ODE system (\ref{eq:ivp_ode_diff_src}). Therefore, the results and proofs presented above for approximation orders should generally not hold. A rigorous derivation of the approximation order $N+1$ of the local solution $\mathbf{u}_{L}$ obtained by the local DG predictor essentially relies on the assumption that the function $\mathbf{F}(\mathbf{u},\, t)$ satisfies the Lipschitz condition (\ref{eq:lip_cond_f}) in the vector argument $\mathbf{u}$, which casts doubt on the expected convergence order $p_{\rm th.}^{l} = N+1$ for the local solution $\mathbf{u}_{L}$. However, a rigorous derivation of the approximation order $N+2$ of the improved local solution $\mathbf{u}_{\rm IL}$ is not related to this assumption regarding function $\mathbf{F}(\mathbf{u},\, t)$, but is related to the approximation order $N+1$ of the local solution (and, of course, the approximation order $2N+1$ of the numerical solution at the nodes $\mathbf{u}_{n}$). It is expected that, regardless of a possible approximation order other than $N+1$ the local solution $\mathbf{u}_{L}$, the approximation order of the improved local solution $\mathbf{u}_{\rm IL}$ will still be one higher than the approximation order of the local solution $\mathbf{u}_{L}$. Similar conclusions are expected for the expected convergence orders $p_{\rm th.}^{\rm imp} = p_{\rm th.}^{l} + 1$. The obtained results are presented in Fig.~\ref{fig:bratu} and Table~\ref{tab:conv_orders_bratu}.

\begin{table*}[h!]
\centering
\normalsize
\caption{%
Convergence orders $p_{f}$, $p_{L_{1}}$, $p_{L_{2}}$, $p_{L_{\infty}}$, calculated in $f$-norm and norms $L_{1}$, $L_{2}$, $L_{\infty}$, of the numerical solution of the ADER-DG method for the problem (\ref{eq:bratu_ode}); $N$ is the degree of the basis polynomials $\varphi_{p}$. Orders $p^{n}$ are calculated for \textit{the numerical solution at the nodes} $\mathbf{u}_{n}$; orders $p^{\rm imp}$ --- for \textit{the improved local solution} $\mathbf{u}_{\rm IL}$; orders $p^{l}$ --- for \textit{the local solution} $\mathbf{u}_{L}$. The theoretical values of convergence order $p_{\rm th.}^{n} = 2N+1$, $p_{\rm th.}^{l} = N+1$ and $p^{\rm imp}_{\rm th.} = N+2$ are presented for comparison.
}
\label{tab:conv_orders_bratu}
\setlength{\tabcolsep}{3.5pt}
\begin{tabular}{@{}|l|llll|c|lll|c|lll|c|@{}}
\toprule
$N$ & $p^{n}_{f}$ &
$p^{n}_{L_{1}}$ & $p^{n}_{L_{2}}$ & $p^{n}_{L_{\infty}}$ & $p^{n}_{\rm th.}$ &
$p^{l}_{L_{1}}$ & $p^{l}_{L_{2}}$ & $p^{l}_{L_{\infty}}$ & $p^{l}_{\rm th.}$ &
$p^{\rm imp}_{L_{1}}$ & $p^{\rm imp}_{L_{2}}$ & $p^{\rm imp}_{L_{\infty}}$ & $p^{\rm imp}_{\rm th.}$\\
\midrule
$1$ & $3.05$ & $3.16$ & $3.18$ & $3.05$ & $3$ & $2.02$ & $2.02$ & $1.87$ & $2$ & $3.01$ & $2.98$ & $2.80$ & $3$\\
$2$ & $4.90$ & $5.14$ & $5.12$ & $4.90$ & $5$ & $3.01$ & $2.99$ & $2.80$ & $3$ & $3.99$ & $3.96$ & $3.74$ & $4$\\
$3$ & $6.37$ & $6.91$ & $6.76$ & $6.37$ & $7$ & $4.00$ & $3.97$ & $3.74$ & $4$ & $4.98$ & $4.94$ & $4.67$ & $5$\\
$4$ & $6.91$ & $7.83$ & $7.45$ & $6.91$ & $9$ & $4.99$ & $4.95$ & $4.69$ & $5$ & $5.96$ & $5.91$ & $5.63$ & $6$\\
$5$ & $11.0$ & $11.3$ & $11.4$ & $11.0$ & $11$ & $5.98$ & $5.93$ & $5.63$ & $6$ & $6.95$ & $6.89$ & $6.57$ & $7$\\
$6$ & $16.7$ & $16.5$ & $16.7$ & $16.5$ & $13$ & $6.97$ & $6.90$ & $6.58$ & $7$ & $7.93$ & $7.86$ & $7.52$ & $8$\\
$7$ & $17.1$ & $17.9$ & $17.6$ & $17.1$ & $15$ & $7.95$ & $7.87$ & $7.53$ & $8$ & $8.92$ & $8.83$ & $8.48$ & $9$\\
$8$ & $18.4$ & $19.2$ & $18.9$ & $18.4$ & $17$ & $8.93$ & $8.84$ & $8.47$ & $9$ & $9.90$ & $9.79$ & $9.42$ & $10$\\
$9$ & $20.1$ & $20.8$ & $20.5$ & $20.1$ & $19$ & $9.91$ & $9.80$ & $9.42$ & $10$ & $10.9$ & $10.8$ & $10.4$ & $11$\\
$10$ & $21.9$ & $22.6$ & $22.3$ & $21.9$ & $21$ & $10.9$ & $10.8$ & $10.4$ & $11$ & $11.9$ & $11.7$ & $11.3$ & $12$\\
\midrule
$11$ & $23.7$ & $24.5$ & $24.2$ & $23.7$ & $23$ & $11.9$ & $11.7$ & $11.3$ & $12$ & $12.8$ & $12.7$ & $12.3$ & $13$\\
$12$ & $25.6$ & $26.4$ & $26.0$ & $25.6$ & $25$ & $12.8$ & $12.7$ & $12.3$ & $13$ & $13.8$ & $13.6$ & $13.2$ & $14$\\
$13$ & $27.4$ & $28.3$ & $27.9$ & $27.4$ & $27$ & $13.8$ & $13.6$ & $13.2$ & $14$ & $14.8$ & $14.6$ & $14.2$ & $15$\\
$14$ & $29.3$ & $30.2$ & $29.8$ & $29.3$ & $29$ & $14.8$ & $14.6$ & $14.2$ & $15$ & $15.7$ & $15.5$ & $15.1$ & $16$\\
$15$ & $31.2$ & $32.1$ & $31.7$ & $31.2$ & $31$ & $15.8$ & $15.6$ & $15.1$ & $16$ & $16.7$ & $16.5$ & $16.1$ & $17$\\
$16$ & $33.1$ & $34.0$ & $33.6$ & $33.1$ & $33$ & $16.7$ & $16.5$ & $16.1$ & $17$ & $17.7$ & $17.5$ & $17.0$ & $18$\\
$17$ & $35.0$ & $35.9$ & $35.5$ & $35.0$ & $35$ & $17.7$ & $17.5$ & $17.0$ & $18$ & $18.7$ & $18.4$ & $17.9$ & $19$\\
$18$ & $36.9$ & $37.8$ & $37.4$ & $36.9$ & $37$ & $18.7$ & $18.4$ & $17.9$ & $19$ & $19.6$ & $19.4$ & $18.9$ & $20$\\
$19$ & $38.8$ & $39.7$ & $39.3$ & $38.8$ & $39$ & $19.6$ & $19.4$ & $18.9$ & $20$ & $20.6$ & $20.3$ & $19.8$ & $21$\\
$20$ & $40.7$ & $41.6$ & $41.1$ & $40.7$ & $41$ & $20.6$ & $20.3$ & $19.8$ & $21$ & $21.5$ & $21.3$ & $20.8$ & $22$\\
\midrule
$21$ & $42.5$ & $43.5$ & $43.0$ & $42.5$ & $43$ & $21.6$ & $21.3$ & $20.8$ & $22$ & $22.5$ & $22.2$ & $21.7$ & $23$\\
$22$ & $44.4$ & $45.4$ & $44.9$ & $44.4$ & $45$ & $22.5$ & $22.2$ & $21.7$ & $23$ & $23.5$ & $23.2$ & $22.7$ & $24$\\
$23$ & $46.3$ & $47.3$ & $46.8$ & $46.3$ & $47$ & $23.5$ & $23.2$ & $22.7$ & $24$ & $24.4$ & $24.1$ & $23.6$ & $25$\\
$24$ & $48.2$ & $49.2$ & $48.7$ & $48.2$ & $49$ & $24.4$ & $24.1$ & $23.6$ & $25$ & $25.4$ & $25.1$ & $24.6$ & $26$\\
$25$ & $50.1$ & $51.1$ & $50.6$ & $50.1$ & $51$ & $25.4$ & $25.1$ & $24.6$ & $26$ & $26.4$ & $26.0$ & $25.5$ & $27$\\
$26$ & $52.0$ & $53.0$ & $52.5$ & $52.0$ & $53$ & $26.4$ & $26.0$ & $25.5$ & $27$ & $27.3$ & $27.0$ & $26.5$ & $28$\\
$27$ & $53.9$ & $54.9$ & $54.4$ & $53.9$ & $55$ & $27.3$ & $27.0$ & $26.5$ & $28$ & $28.3$ & $27.9$ & $27.4$ & $29$\\
$28$ & $55.8$ & $56.8$ & $56.3$ & $55.8$ & $57$ & $28.3$ & $27.9$ & $27.4$ & $29$ & $29.2$ & $28.9$ & $28.4$ & $30$\\
$29$ & $57.7$ & $58.7$ & $58.2$ & $57.7$ & $59$ & $29.2$ & $28.9$ & $28.4$ & $30$ & $30.2$ & $29.8$ & $29.3$ & $31$\\
$30$ & $59.6$ & $60.6$ & $60.1$ & $59.6$ & $61$ & $30.2$ & $29.8$ & $29.3$ & $31$ & $31.1$ & $30.8$ & $30.3$ & $32$\\
\midrule
$35$ & $69.1$ & $70.1$ & $69.6$ & $69.1$ & $71$ & $35.0$ & $34.5$ & $34.0$ & $36$ & $35.9$ & $35.5$ & $35.0$ & $37$\\
$40$ & $78.5$ & $79.5$ & $79.0$ & $78.5$ & $81$ & $39.7$ & $39.3$ & $38.8$ & $41$ & $40.7$ & $40.2$ & $39.7$ & $42$\\
$45$ & $88.0$ & $89.0$ & $88.5$ & $88.0$ & $91$ & $44.5$ & $44.0$ & $43.5$ & $46$ & $45.4$ & $45.0$ & $44.5$ & $47$\\
$50$ & $97.5$ & $98.5$ & $98.0$ & $97.5$ & $101$ & $49.2$ & $48.8$ & $48.3$ & $51$ & $50.2$ & $49.7$ & $49.2$ & $52$\\
$55$ & $106.9$ & $107.9$ & $107.4$ & $106.9$ & $111$ & $54.0$ & $53.5$ & $53.0$ & $56$ & $54.9$ & $54.4$ & $53.9$ & $57$\\
$60$ & $116.4$ & $117.4$ & $116.9$ & $116.4$ & $121$ & $58.7$ & $58.2$ & $57.7$ & $61$ & $59.7$ & $59.2$ & $58.7$ & $62$\\
\bottomrule
\end{tabular}
\end{table*}

The dependencies of the numerical solutions $\mathbf{u}_{L}$, $\mathbf{u}_{\rm IL}$, $\mathbf{u}_{n}$ and the exact analytical solution $\mathbf{u}^{\rm ex}$, the dependencies of the local error $\varepsilon$ (\ref{eq:eps_local_def}) of the numerical solutions, and the dependence of the global error $e$ (\ref{eq:eps_un_global_def}), (\ref{eq:eps_ul_global_def}) of the numerical solutions on the discretization step ${\Delta t}$, for polynomial degrees $N = 1$, $4$, $12$ and $60$, are shown in Fig.~\ref{fig:bratu}. A comparison of the obtained dependencies of the numerical solutions $\mathbf{u}_{L}(t)$, $\mathbf{u}_{\rm IL}(t)$, $\mathbf{u}_{n}$ with the exact analytical solution $\mathbf{u}^{\rm ex}(t)$, presented in Fig.~\ref{fig:bratu} (\subref{fig:bratu:a1}, \subref{fig:bratu:a2}, \subref{fig:bratu:a3}, \subref{fig:bratu:a4}) for the component $u_{1}$ and in Fig.~\ref{fig:bratu} (\subref{fig:bratu:b1}, \subref{fig:bratu:b2}, \subref{fig:bratu:b3}, \subref{fig:bratu:b4}) for the component $u_{2}$, demonstrates a high-quality agreement. The local error $\varepsilon$ (\ref{eq:eps_local_def}) of the numerical solutions, witch is presented in Fig.~\ref{fig:bratu} (\subref{fig:bratu:c1}, \subref{fig:bratu:c2}, \subref{fig:bratu:c3}, \subref{fig:bratu:c4}), shows that the local error $\varepsilon$ of the improved local solution $\mathbf{u}_{\rm IL}$ is significantly lower than the local error $\varepsilon$ of the local solution $\mathbf{u}_{L}$, with the differences reaching $2$--$4$ orders of magnitude. Moreover, the local error $\varepsilon$ of the numerical solution $\mathbf{u}_{n}$ at the grid nodes $t_{n}$ is significantly smaller than the local errors $\varepsilon$ of the local solution $\mathbf{u}_{L}$ and the improved local solution $\mathbf{u}_{\rm IL}$, amounting to approximately $1$--$2$ orders of magnitude in the case of polynomial degree $N = 1$, approximately $6$--$8$ orders of magnitude in the case of polynomial degree $N = 4$, approximately $20$--$24$ orders of magnitude in the case of polynomial degree $N = 12$ and approximately $102$--$114$ orders of magnitude in the case of polynomial degree $N = 60$ (for this purpose, breaks in the graph along the vertical axis are inserted in Fig.~\ref{fig:bratu} (\subref{fig:bratu:c4})).

\begin{figure}[h!]
\captionsetup[subfigure]{%
	position=bottom,
	font+=smaller,
	textfont=normalfont,
	singlelinecheck=off,
	justification=raggedright
}
\centering
\begin{subfigure}{0.24\textwidth}
\includegraphics[width=\textwidth]{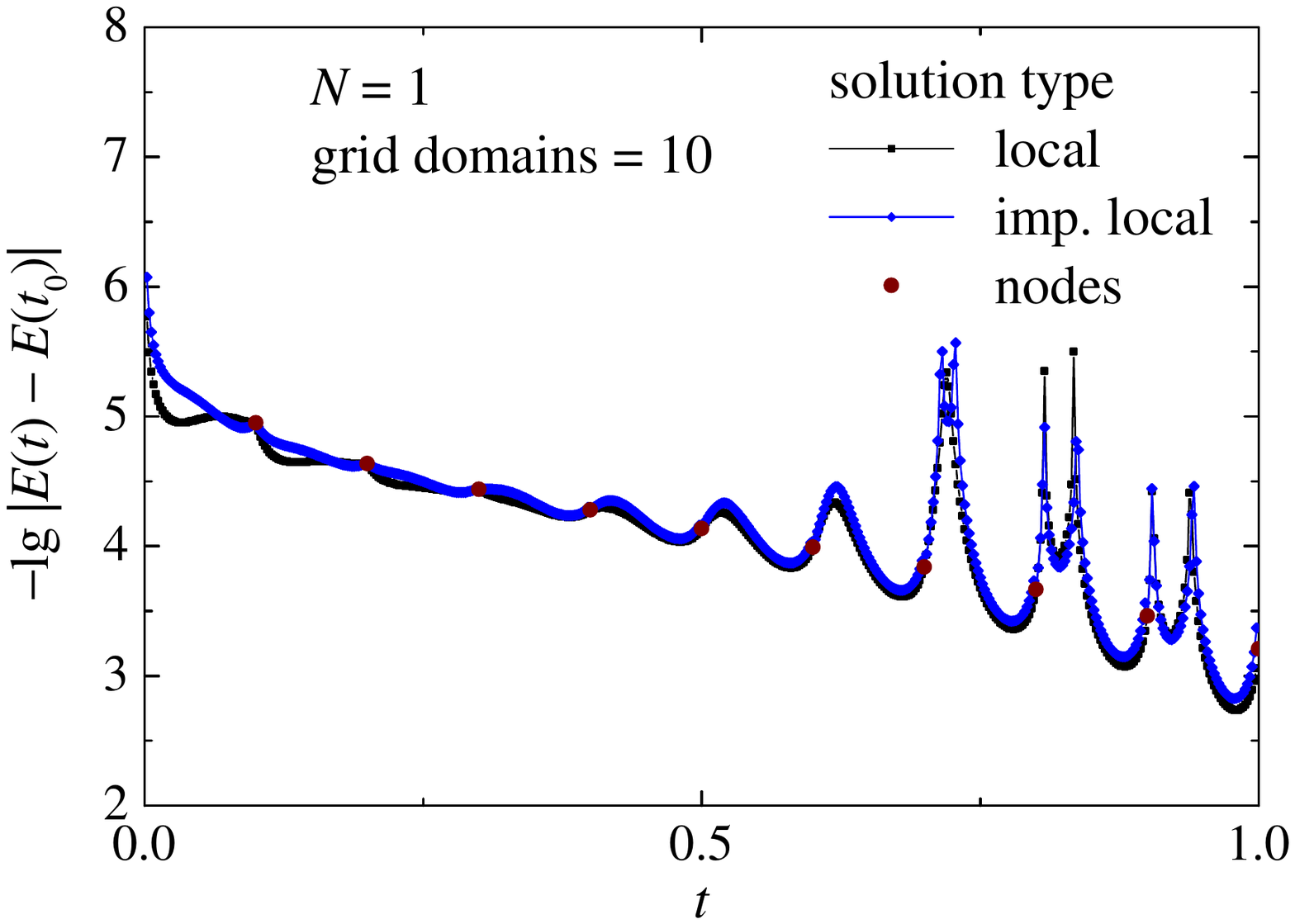}
\vspace{-8mm}\caption{\label{fig:econs_bratu:a1}}
\end{subfigure}
\begin{subfigure}{0.24\textwidth}
\includegraphics[width=\textwidth]{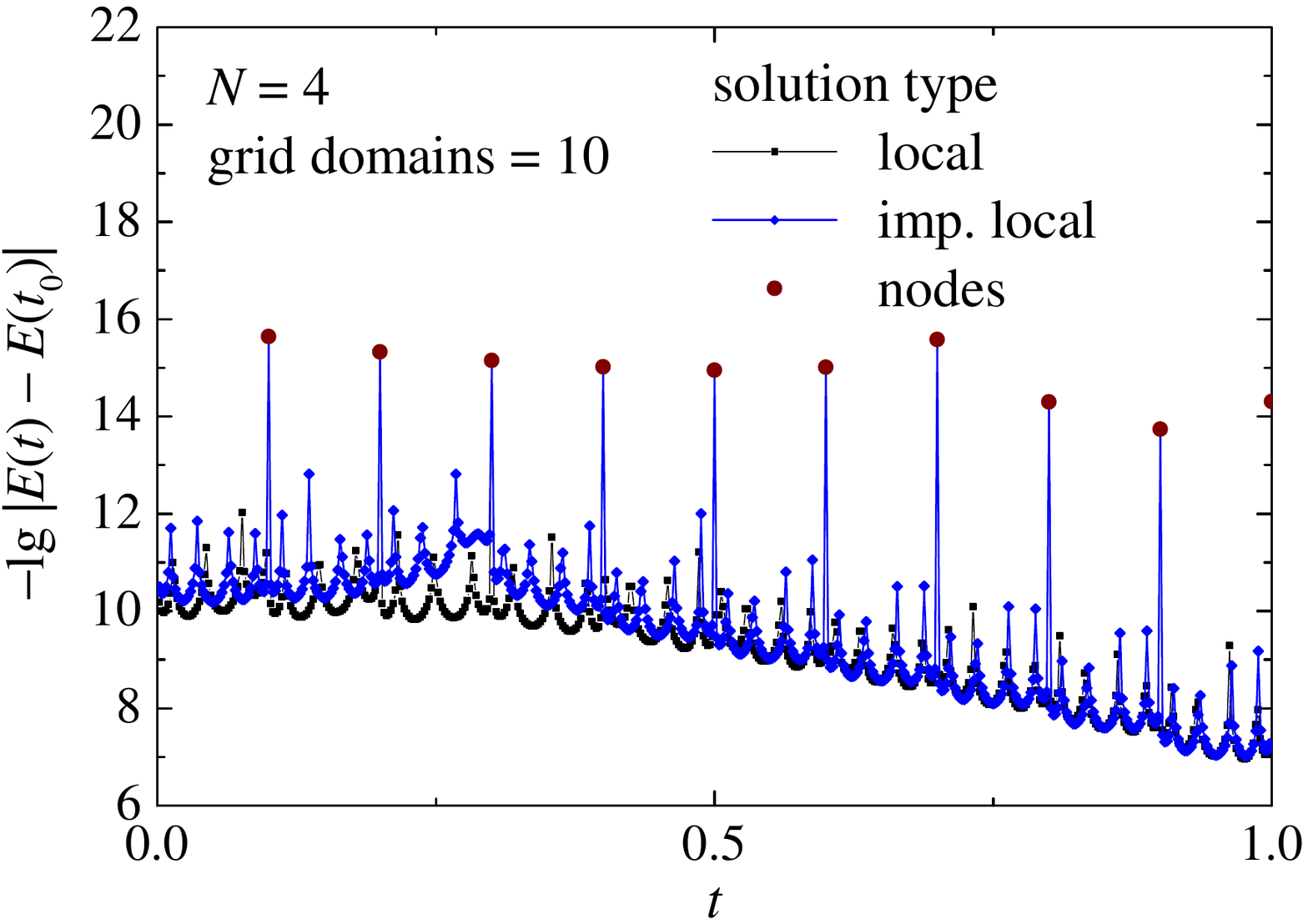}
\vspace{-8mm}\caption{\label{fig:econs_bratu:a2}}
\end{subfigure}
\begin{subfigure}{0.24\textwidth}
\includegraphics[width=\textwidth]{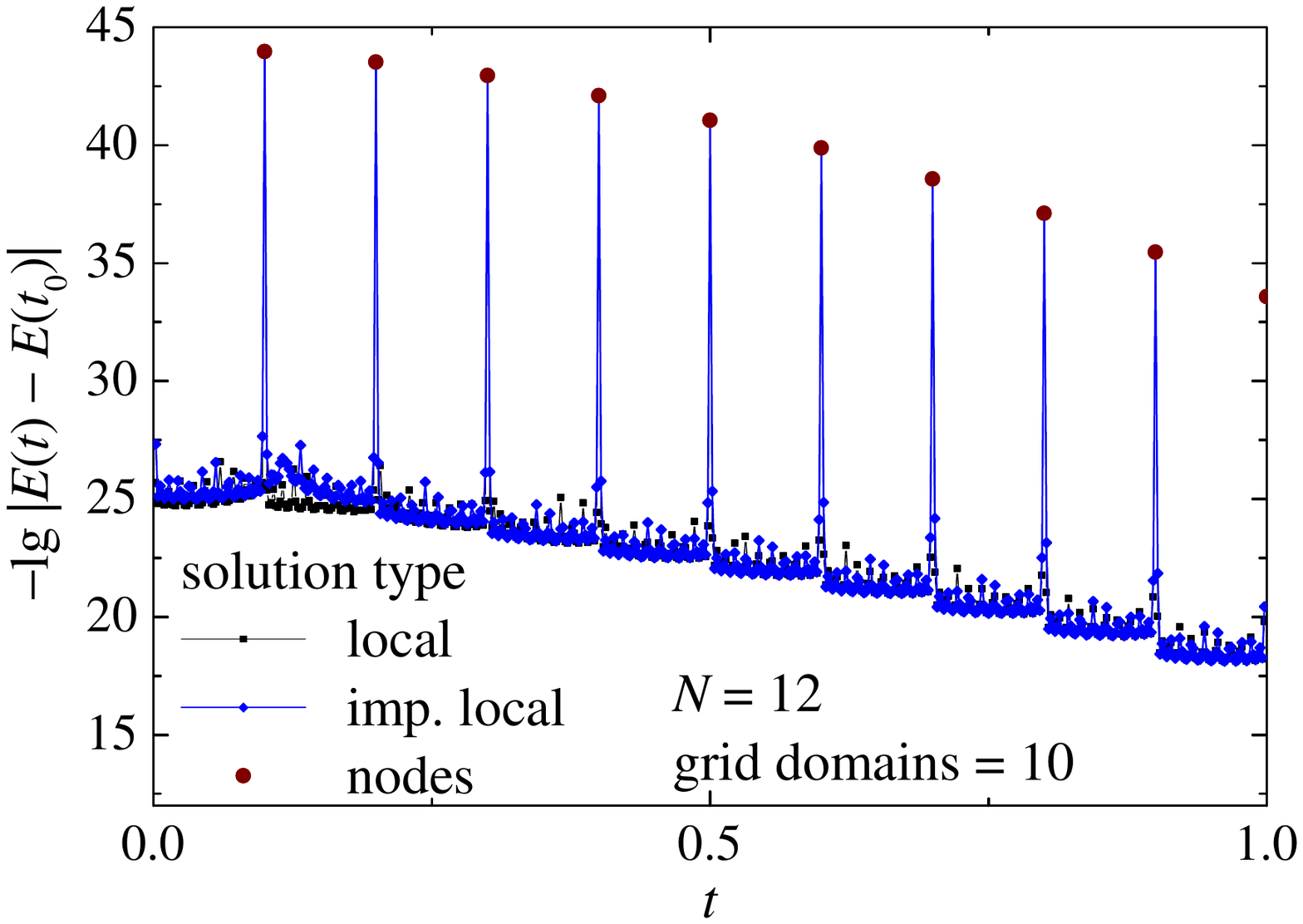}
\vspace{-8mm}\caption{\label{fig:econs_bratu:a3}}
\end{subfigure}
\begin{subfigure}{0.24\textwidth}
\includegraphics[width=\textwidth]{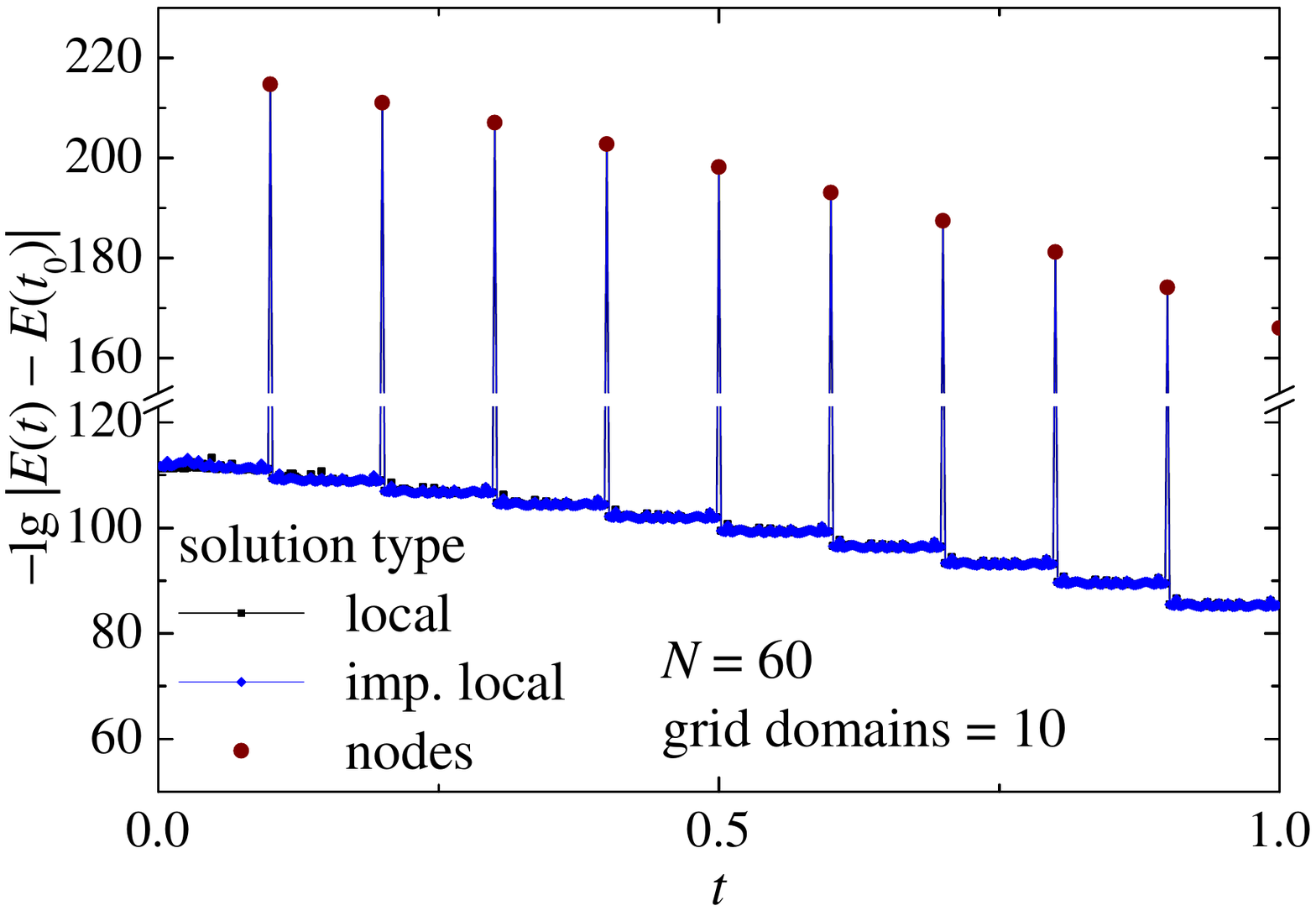}
\vspace{-8mm}\caption{\label{fig:econs_bratu:a4}}
\end{subfigure}
\caption{%
Dependence of the negative logarithm of the error $-\lg|E(t)-E(t_{0})|$ of the energy conservation law of the numerical solution of the system (\ref{eq:bratu_ode}) on the argument $t$ for solution at nodes $\mathbf{u}_{n}$, the local solution $\mathbf{u}_{L}(t)$ and the improved local solution $\mathbf{u}_{\rm IL}(t)$, obtained using polynomials with degrees $N = 1$~(\subref{fig:econs_bratu:a1}), $N = 4$~(\subref{fig:econs_bratu:a4}), $N = 12$~(\subref{fig:econs_bratu:a3}), $N = 60$~(\subref{fig:econs_bratu:a4}).
}
\label{fig:econs_bratu}
\end{figure}

The resulting dependencies of the global error $e$ (\ref{eq:eps_un_global_def}), (\ref{eq:eps_ul_global_def}) on the discretization step ${\Delta t}$, presented in Fig.~\ref{fig:bratu} (\subref{fig:bratu:d1}, \subref{fig:bratu:d2}, \subref{fig:bratu:d3}, \subref{fig:bratu:d4}, \subref{fig:bratu:e1}, \subref{fig:bratu:e2}, \subref{fig:bratu:e3}, \subref{fig:bratu:e4}, \subref{fig:bratu:f1}, \subref{fig:bratu:f2}, \subref{fig:bratu:f3}, \subref{fig:bratu:f4}), demonstrate a high quality of linear approximation for all studied polynomial degrees $N$ (Fig.~\ref{fig:bratu} only shows results for polynomial degrees $N = 1$, $4$, $12$ and $60$). In the case of a numerical solution $\mathbf{u}_{n}$ in grid nodes for polynomials of degree $N = 60$, the global error $e$ for which is shown in Fig.~\ref{fig:bratu} (\subref{fig:bratu:f4}), a decrease in the discretization step ${\Delta t}$ by a factor of $2$ leads to a decrease in global errors $e$ by $40$--$42$ orders of magnitude. Some deviations and anomalies are observed for errors $e^{n}$ in Fig.~\ref{fig:bratu} (\subref{fig:bratu:f2}) for polynomial degree $N = 4$.

Based on the approximation of the obtained dependencies $e({\Delta t})$ in log-log scale by a linear function $\lg{e({\Delta t})} \propto p\cdot\lg{{\Delta t}}$, empirical convergence orders $p$ are calculated and presented in Table~\ref{tab:conv_orders_bratu} for all polynomial degrees $N = 1, \ldots, 30$ and polynomial degrees $N$ up to $60$ with a step of $5$. The obtained results show that the convergence orders $p$ for the numerical solution $\mathbf{u}_{n}$ at grid nodes in many cases of polynomial degrees $N$ correspond well to the expected values $p_{\rm th.}^{n} = 2N+1$ (\ref{eq:conv_ords_exp}), or even exceed them, sometimes even by $1.0$--$3.0$ (especially in the case of polynomial degrees $N = 6$--$8$), and in some cases $N$ are lower than the expected values $p_{\rm th.}^{n}$, especially in cases of large polynomial degrees $N \geqslant 30$--$35$. The convergence orders $p$ for the local numerical solution $\mathbf{u}_{L}$ correspond well to the expected values $p_{\rm th.}^{l} = N+1$ (\ref{eq:conv_ords_exp}) for polynomial degrees $N \leqslant 25$--$30$, and with increasing of polynomial degree $N$ the convergence orders $p$ become lower than the expected values $p_{\rm th.}^{l}$ by $1.0$--$3.0$. These behavioral characteristics of the empirical convergence orders $p$ can be considered expected, and they were encountered and studied in~\cite{ader_dg_ode_jsc}. It is important to note that the empirical convergence orders $p^{\rm imp}$ of the improved local numerical solution $\mathbf{u}_{\rm IL}$ are one unit higher than the empirical convergence orders $p^{l}$ of the local numerical solution $\mathbf{u}_{L}$. This demonstrates the validity of the global error estimates obtained in this work.

\corrtext{Similar to the three previous Examples in Sections~\ref{sec:apps:lin_diss},~\ref{sec:apps:harm_osc} and~\ref{sec:apps:pend}, the system of equations (\ref{eq:bratu_ode}) considered in this Example can be presented in a conservative form, which is characterized by the energy conservation law $E(t) = \mathrm{const}$~\cite{LandauMechanics, GoldsteinMechanics} of the following form
\begin{equation}\label{eq:econs_bratu}
E(t) = \frac{\dot{x}^{2}(t)}{2} - 2\exp(x(t)) = \mathrm{const} \equiv
E(t_{0}) = \frac{\dot{x}^{2}(0)}{2} - 2\exp(x(0)) = -2,
\end{equation}
which represents the integral of motion.}

\corrtext{Fig.~\ref{fig:econs_bratu} shows the dependence of the negative logarithm of the error $-\lg|E(t)-E(t_{0})|$ of the energy conservation law (\ref{eq:econs_bratu}) of the numerical solution of the system (\ref{eq:bratu_ode}) on the argument $t$ for solution at nodes $\mathbf{u}_{n}$, the local solution $\mathbf{u}_{L}(t)$ and the improved local solution $\mathbf{u}_{\rm IL}(t)$, obtained using polynomials with degrees $N = 1$~(\subref{fig:econs_bratu:a1}), $N = 4$~(\subref{fig:econs_bratu:a4}), $N = 12$~(\subref{fig:econs_bratu:a3}), $N = 60$~(\subref{fig:econs_bratu:a4}). The obtained results clearly demonstrate that the energy conservation law (\ref{eq:bratu_ode}) is not strictly satisfied in the numerical solution. However, due to the possibility of achieving a high order $p$, even in the cases of degrees $N = 12$ and $60$, the error in the energy conservation law's fulfillment over the entire solution domain becomes smaller than the characteristic rounding error of double-precision floating-point numbers $\sim 10^{-15}$-$10^{-17}$. The presented results also clearly demonstrate that the accuracy of the fulfillment of the energy conservation law (\ref{eq:econs_bratu}) for the improved local solution $\mathbf{u}_{\rm IL}(t)$ is significantly higher than for the local solution $\mathbf{u}_{L}(t)$. Therefore, it can be concluded that, similar to the three previous Examples in Sections~\ref{sec:apps:lin_diss},~\ref{sec:apps:harm_osc} and~\ref{sec:apps:pend}, despite the dissipative nature of the ADER-DG numerical method with local DG predictor, a sufficiently high degree $N$ can be chosen such that the accuracy of the fulfillment of the energy conservation law will be at or below the characteristic error of representing real numbers as floating-point numbers. The results obtained are qualitatively consistent with the results presented previously in Example in Sections~\ref{sec:apps:pend}.}

The obtained results allowed to conclude that the ADER-DG numerical method with a local DG predictor provides a highly accurate numerical solution to the initial value problem for the nonlinear ODE system (\ref{eq:bratu_ode}) presented in this example. The obtained results are in quite acceptable agreement with the theory developed above. The obtained results show that the initial strict assumption (\ref{eq:lip_cond_f}) is necessary mainly only for the rigorous proof and derivation of the approximation order of the local numerical solution $\mathbf{u}_{L}$, while empirical estimates confirm the obtained result in the case of violation of this assumption. The improved local numerical solution $\mathbf{u}_{\rm IL}$ demonstrates higher accuracy and a higher convergence order compared to the local numerical solution $\mathbf{u}_{L}$.

\section{Features of the numerical method for DAE systems}
\label{sec:daes}

This Section presents a description of the structural features of the ADER-DG numerical method with a local DG predictor for solving DAE systems of the following type:
\begin{equation}\label{eq:ivp_dae_diff_src}
\begin{split}
\frac{d\mathbf{u}}{dt} = \mathbf{F}(\mathbf{u}, \mathbf{v}, t),\quad &t\in\Omega = \left\{t\, |\, t \in [t_{0},\ t_{f}]\right\},\\
\mathbf{0} = \mathbf{G}(\mathbf{u}, \mathbf{v}, t),\quad &\mathbf{u}(t_{0}) = \mathbf{u}_{0},\quad \mathbf{v}(t_{0}) = \mathbf{v}_{0},\\
\end{split}
\end{equation}
where the functions $\mathbf{u}: \Omega \rightarrow \mathbb{R}^{D_{\rm u}}$, $\mathbf{v}: \Omega \rightarrow \mathbb{R}^{D_{\rm v}}$ are desired functions, the functions $\mathbf{F}: \mathbb{R}^{D_{\rm u}}\times\mathbb{R}^{D_{\rm v}}\times\Omega\rightarrow\mathbb{R}^{D_{\rm u}}$, $\mathbf{G}: \mathbb{R}^{D_{\rm u}}\times\mathbb{R}^{D_{\rm v}}\times\Omega\rightarrow\mathbb{R}^{D_{\rm v}}$ are right side functions, which are given. In terms of DAE systems the desired functions $\mathbf{u}$ are a set of differential variables, the desired functions $\mathbf{v}$ are a set of algebraic variables, $t\in\Omega$ is an independent variable. $D_{\rm u}$ is the number of differential variables, $D_{\rm v}$ is the number of algebraic variables, $D = D_{\rm u} + D_{\rm v}$ is the total number of variables. The initial conditions $\mathbf{u}(t_{0}) = \mathbf{u}_{0}$, $\mathbf{v}(t_{0}) = \mathbf{v}_{0}$ are given at point $t_{0}$ and must be consistent with the DAE system (\ref{eq:ivp_dae_diff_src}), however, technically the initial value for the algebraic variables $\mathbf{v}_{0}$ can be formally discarded if its value is not required in solving the problem. The use of the ADER-DG numerical method with a local DG predictor for solving DAE system (\ref{eq:ivp_dae_diff_src}) was proposed in~\cite{ader_dg_dae}, where the applicability of the method was substantiated in detail, primarily empirically.

In this work, the formulation of the ADER-DG numerical method with a local DG predictor for solving DAE systems is justified strictly based on the use of the $\epsilon$-embedding method for the DAE system (\ref{eq:ivp_dae_diff_src}) (see~\cite{Hairer_book_2}), which is widely used to prove the applicability of numerical methods for solving the initial value problem for ODE systems (\ref{eq:ivp_ode_diff_src}) to DAE systems~\cite{Ascher_book_1998, Kunkel_book_2006, Riaza_book_2008, Gerdts_book_2012, handbook_num_anal_daes_2002, Hairer_book_2}.

The initial value problem for an ODE system (\ref{eq:ivp_ode_diff_src}) with a singularly perturbed subsystem, corresponding to the $\epsilon$-embedding of the original DAE system (\ref{eq:ivp_dae_diff_src}), is represented in the form~\cite{Hairer_book_2}:
\begin{equation}\label{eq:ivp_dae_eps_diff_src}
\begin{split}
\frac{d\mathbf{u}^{\epsilon}}{dt} = \mathbf{F}(\mathbf{u}^{\epsilon}, \mathbf{v}^{\epsilon}, t),\quad &t\in\Omega = \left\{t\, |\, t \in [t_{0},\ t_{f}]\right\},\\
\epsilon\frac{d\mathbf{v^{\epsilon}}}{dt} = \mathbf{G}(\mathbf{u}^{\epsilon}, \mathbf{v}^{\epsilon}, t),\quad
&\mathbf{u}^{\epsilon}(t_{0}) = \mathbf{u}_{0},\quad \mathbf{v}^{\epsilon}(t_{0}) = \mathbf{v}_{0},\\
\end{split}
\end{equation}
where the functions $\mathbf{u}^{\epsilon}: \Omega \rightarrow \mathbb{R}^{D_{\rm u}}$, $\mathbf{v}^{\epsilon}: \Omega \rightarrow \mathbb{R}^{D_{\rm v}}$ are desired perturbed functions, which are parameterized by parameter $\epsilon$. The initial value problem for an ODE system (\ref{eq:ivp_dae_eps_diff_src}) is transformed into the original DAE system (\ref{eq:ivp_dae_diff_src}) in the limit $\epsilon \rightarrow 0$. The formulation of the ADER-DG numerical method with a local DG predictor for the ODE system (\ref{eq:ivp_dae_diff_src}) is close to the original formulation presented in Section~\ref{sec:gen_descr} and is based on the integral representation of the original system (\ref{eq:ivp_dae_diff_src}) in the discretization domain $\Omega_{n}$:
\begin{equation}\label{eq:ivp_dae_eps_int_src_by_tau}
\begin{split}
&\mathbf{u}_{n+1}^{\epsilon} = \mathbf{u}_{n}^{\epsilon} +
{\Delta t}_{n} \int\limits_{0}^{1} \mathbf{F}(\mathbf{u}^{\epsilon}(t(\tau)),\, \mathbf{v}^{\epsilon}(t(\tau)),\, t(\tau)) d\tau,\\
&\epsilon\mathbf{v}_{n+1}^{\epsilon} = \epsilon\mathbf{v}_{n}^{\epsilon} +
{\Delta t}_{n} \int\limits_{0}^{1} \mathbf{G}(\mathbf{u}^{\epsilon}(t(\tau)),\, \mathbf{v}^{\epsilon}(t(\tau)),\, t(\tau)) d\tau.
\end{split}
\end{equation}
The integrands in (\ref{eq:ivp_dae_eps_int_src_by_tau}) use a local solution $(\mathbf{q}_{n}^{\epsilon}(\tau),\, \mathbf{r}_{n}^{\epsilon}(\tau))$, represented as an expansion in basis functions $\{\varphi_{p}\}$ (\ref{eq:basis_funcs_def}), as arguments:
\begin{equation}\label{eq:local_sol_dae_eps}
\begin{split}
\mathbf{q}_{n}^{\epsilon}(\tau) = \sum\limits_{p} \hat{\mathbf{q}}_{n, p}^{\epsilon}\varphi_{p}(\tau),\qquad
\mathbf{r}_{n}^{\epsilon}(\tau) = \sum\limits_{p} \hat{\mathbf{r}}_{n, p}^{\epsilon}\varphi_{p}(\tau),
\end{split}
\end{equation}
where $(\{\hat{\mathbf{q}}_{n, p}^{\epsilon}\},\, \{\hat{\mathbf{r}}_{n, p}^{\epsilon}\})$ is the set of expansion coefficients, and is found as a solution to the weak form of the ODE system (\ref{eq:ivp_ode_diff_src}) in $\Omega_{n}$:
\begin{equation}\label{eq:ivp_dae_eps_weak}
\begin{split}
&\int\limits_{0}^{1}\varphi_{p}(\tau)\left[
	\frac{d\mathbf{q}_{n}^{\epsilon}(\tau)}{d\tau} - {\Delta t}_{n}\mathbf{F}(\mathbf{q}_{n}^{\epsilon}(\tau), \mathbf{r}_{n}^{\epsilon}(\tau), t(\tau))
\right]d\tau = 0,\quad\ \mathbf{q}_{n}^{\epsilon}(0) = \mathbf{u}_{n}^{\epsilon},\\
&\int\limits_{0}^{1}\varphi_{p}(\tau)\left[
	\epsilon\frac{d\mathbf{r}_{n}^{\epsilon}(\tau)}{d\tau} - {\Delta t}_{n}\mathbf{G}(\mathbf{q}_{n}^{\epsilon}(\tau), \mathbf{r}_{n}^{\epsilon}(\tau), t(\tau))
\right]d\tau = 0,\quad \mathbf{r}_{n}^{\epsilon}(0) = \mathbf{v}_{n}^{\epsilon},
\end{split}
\end{equation}
where the initial condition is also understood in weak form. The quadrature formula (\ref{eq:gl_rule}) is also used to calculate the integrals, as well as the point-wise evaluation (\ref{eq:point_wise_evaluation}) that follow from it:
\begin{equation}\label{eq:point_wise_evaluation_dae_eps}
\begin{split}
&\mathbf{F}(\mathbf{q}_{n}^{\epsilon}(\tau), \mathbf{r}_{n}^{\epsilon}(\tau), t(\tau)) = \sum\limits_{p = 0}^{N} \hat{\mathbf{F}}_{n, p}^{\epsilon}\varphi_{p}(\tau)\ \mapsto\ 
\sum\limits_{p = 0}^{N} \mathbf{F}(\hat{\mathbf{q}}_{n, p}^{\epsilon}, \hat{\mathbf{r}}_{n, p}, t(\tau_{n}))\varphi_{p}(\tau),\\
&\mathbf{G}(\mathbf{q}_{n}^{\epsilon}(\tau), \mathbf{r}_{n}^{\epsilon}(\tau), t(\tau)) = \sum\limits_{p = 0}^{N} \hat{\mathbf{G}}_{n, p}^{\epsilon}\varphi_{p}(\tau)\ \mapsto\ 
\sum\limits_{p = 0}^{N} \mathbf{G}(\hat{\mathbf{q}}_{n, p}^{\epsilon}, \hat{\mathbf{r}}_{n, p}^{\epsilon}, t(\tau_{n}))\varphi_{p}(\tau).
\end{split}
\end{equation}
The system of nonlinear algebraic equations of a local DG predictor for calculating the expansion coefficients $(\{\hat{\mathbf{q}}_{n, p}^{\epsilon}\},\, \{\hat{\mathbf{r}}_{n, p}^{\epsilon}\})$ of the local solution $(\mathbf{q}_{n}^{\epsilon}(\tau),\, \mathbf{r}_{n}^{\epsilon}(\tau))$ is obtained by substituting expressions (\ref{eq:local_sol_dae_eps}), (\ref{eq:point_wise_evaluation_dae_eps}), with using properties of basis functions (\ref{eq:basis_funcs_def}) and definitions (\ref{eq:kappa_mu_defs}), into expressions (\ref{eq:ivp_dae_eps_weak}) (similar to system of nonlinear algebraic equations (\ref{eq:snae_lstdg})), and in this case takes the following form:
\begin{equation}\label{eq:snae_lstdg_dae_eps}
\begin{split}
&\hat{\mathbf{q}}_{n, p}^{\epsilon} -
{\Delta t}_{n}\sum\limits_{q = 0}^{N} \mathrm{A}_{pq}\mathbf{F}(\hat{\mathbf{q}}_{n, q}^{\epsilon}, \hat{\mathbf{r}}_{n, q}^{\epsilon}, t(\tau_{q})) =
\mathbf{u}_{n}^{\epsilon},\\
&\epsilon\hat{\mathbf{r}}_{n, p}^{\epsilon} - 
{\Delta t}_{n}\sum\limits_{q = 0}^{N} \mathrm{A}_{pq}\mathbf{G}(\hat{\mathbf{q}}_{n, q}^{\epsilon}, \hat{\mathbf{r}}_{n, q}^{\epsilon}, t(\tau_{q})) =
\epsilon\mathbf{v}_{n}^{\epsilon}.
\end{split}
\end{equation}
Calculating the solution at grid nodes $(\mathbf{u}_{n}^{\epsilon},\, \mathbf{v}_{n}^{\epsilon})$ is carried out by the following expressions obtained by applying the quadrature formula to the integral form (\ref{eq:ivp_dae_eps_int_src_by_tau}) of the representation of the original ODE system (\ref{eq:ivp_dae_eps_diff_src}):
\begin{equation}\label{eq:ader_dg_node_sol_dae_eps}
\begin{split}
&\mathbf{u}_{n+1}^{\epsilon} = \mathbf{u}_{n}^{\epsilon} +
{\Delta t}_{n}\sum\limits_{p = 0}^{N} w_{p} \mathbf{F}\left(\hat{\mathbf{q}}_{n, p}^{\epsilon}, \hat{\mathbf{r}}_{n, p}^{\epsilon}, t(\tau_{p})\right),\\
&\epsilon\mathbf{v}_{n+1}^{\epsilon} = \epsilon\mathbf{v}_{n}^{\epsilon} +
{\Delta t}_{n}\sum\limits_{p = 0}^{N} w_{p} \mathbf{G}\left(\hat{\mathbf{q}}_{n, p}^{\epsilon}, \hat{\mathbf{r}}_{n, p}^{\epsilon}, t(\tau_{p})\right).
\end{split}
\end{equation}
The resulting expressions (\ref{eq:snae_lstdg_dae_eps}) and (\ref{eq:ader_dg_node_sol_dae_eps}) define the ADER-DG numerical method with a local DG predictor for the ODE system (\ref{eq:ivp_dae_eps_diff_src}).

By performing transformations similar to those used to obtain expression (\ref{eq:local_sol_rp}), the following relations are obtained for the local solution $(\mathbf{q}_{n}^{\epsilon}(\tau),\, \mathbf{r}_{n}^{\epsilon}(\tau))$ at point $\tau = 1$ (which corresponds to point $t_{n+1}$ by mapping (\ref{eq:tau_mapping})):
\begin{equation}\label{eq:local_sol_dae_eps_q_rp}
\begin{split}
\mathbf{q}_{n}^{\epsilon}(1)
&= \sum\limits_{p = 0}^{N} \hat{\mathbf{q}}_{n, p}^{\epsilon} \varphi_{p}(1)
= \mathbf{u}_{n}^{\epsilon} + {\Delta t}_{n}\sum\limits_{q = 0}^{N}\left[
	\sum\limits_{p = 0}^{N} \mathrm{A}_{pq} \varphi_{p}(1)
\right]\mathbf{F}(\hat{\mathbf{q}}_{n, q}^{\epsilon}, \hat{\mathbf{r}}_{n, q}^{\epsilon}, t(\tau_{q}))\\
&= \mathbf{u}_{n} + {\Delta t}_{n}\sum\limits_{q = 0}^{N}\left[
	\left(\sum\limits_{r = 0}^{N}\left[\mathrm{K}^{-1}\right]_{rq} \varphi_{r}(1)\right) w_{q}
\right]\mathbf{F}(\hat{\mathbf{q}}_{n, q}^{\epsilon}, \hat{\mathbf{r}}_{n, q}^{\epsilon}, t(\tau_{q}))\\
&= \mathbf{u}_{n}^{\epsilon} +
{\Delta t}_{n}\sum\limits_{q = 0}^{N} w_{q}\mathbf{F}(\hat{\mathbf{q}}_{n, q}^{\epsilon}, \hat{\mathbf{r}}_{n, q}^{\epsilon}, t(\tau_{q})) \equiv \mathbf{u}_{n+1}^{\epsilon},
\end{split}
\end{equation}
\begin{equation}\label{eq:local_sol_dae_eps_r_rp}
\begin{split}
\epsilon\mathbf{r}_{n}^{\epsilon}(1)
&= \sum\limits_{p = 0}^{N} \hat{\mathbf{r}}_{n, p}^{\epsilon} \varphi_{p}(1)
= \epsilon\mathbf{v}_{n}^{\epsilon} + {\Delta t}_{n}\sum\limits_{q = 0}^{N}\left[
	\sum\limits_{p = 0}^{N} \mathrm{A}_{pq} \varphi_{p}(1)
\right]\mathbf{G}(\hat{\mathbf{q}}_{n, q}^{\epsilon}, \hat{\mathbf{r}}_{n, q}^{\epsilon}, t(\tau_{q}))\\
&= \epsilon\mathbf{v}_{n}^{\epsilon} + {\Delta t}_{n}\sum\limits_{q = 0}^{N}\left[
	\left(\sum\limits_{r = 0}^{N}\left[\mathrm{K}^{-1}\right]_{rq} \varphi_{r}(1)\right) w_{q}
\right]\mathbf{G}(\hat{\mathbf{q}}_{n, q}^{\epsilon}, \hat{\mathbf{r}}_{n, q}^{\epsilon}, t(\tau_{q}))\\
&= \epsilon\mathbf{v}_{n}^{\epsilon} +
{\Delta t}_{n}\sum\limits_{q = 0}^{N} w_{q}\mathbf{G}(\hat{\mathbf{q}}_{n, q}^{\epsilon}, \hat{\mathbf{r}}_{n, q}^{\epsilon}, t(\tau_{q})) \equiv
\epsilon\mathbf{v}_{n+1}^{\epsilon},
\end{split}
\end{equation}
where the last relation (\ref{eq:local_sol_dae_eps_r_rp}) implies $\mathbf{r}_{n}^{\epsilon}(1) = \mathbf{v}_{n+1}^{\epsilon}$ under condition $\epsilon \neq 0$, and the first expression corresponds to the expected property (\ref{eq:local_sol_rp}) of the local solution $\mathbf{q}_{n}^{\epsilon}(\tau)$ for the differential subsystem of equations within the DAE system (\ref{eq:ivp_dae_diff_src}). To perform the limit $\epsilon \neq 0$ transition from the $\epsilon$-embedded singularly perturbed ODE system (\ref{eq:ivp_dae_eps_diff_src}) to the original DAE system (\ref{eq:ivp_dae_diff_src}), the following notation is introduced:
\begin{equation}\label{eq:limit_eps_defs}
\begin{split}
&\hat{\mathbf{r}}_{n, p} = \lim\limits_{\epsilon \rightarrow 0} \hat{\mathbf{r}}_{n, p}^{\epsilon}, \quad
\mathbf{r}_{n}(\tau) = \lim\limits_{\epsilon \rightarrow 0} \mathbf{r}_{n}^{\epsilon}(\tau), \quad
\mathbf{v}_{n}  = \lim\limits_{\epsilon \rightarrow 0} \mathbf{v}_{n}^{\epsilon}, \quad
\mathbf{v}_{n+1}  = \lim\limits_{\epsilon \rightarrow 0} \mathbf{v}_{n+1}^{\epsilon},\\
&\hat{\mathbf{q}}_{n, p} = \lim\limits_{\epsilon \rightarrow 0} \hat{\mathbf{q}}_{n, p}^{\epsilon}, \quad\!\!
\mathbf{q}_{n}(\tau) = \lim\limits_{\epsilon \rightarrow 0} \mathbf{q}_{n}^{\epsilon}(\tau), \quad\!
\mathbf{u}_{n}  = \lim\limits_{\epsilon \rightarrow 0} \mathbf{u}_{n}^{\epsilon}, \quad\!
\mathbf{u}_{n+1}  = \lim\limits_{\epsilon \rightarrow 0} \mathbf{u}_{n+1}^{\epsilon},
\end{split}
\end{equation}
which correspond to the notation in the formulation of the ADER-DG numerical method with a local DG predictor for the DAE system (\ref{eq:ivp_dae_diff_src}) in~\cite{ader_dg_dae}. By performing the limit $\epsilon \neq 0$ transition in the second subsystem of the system of nonlinear algebraic equations (\ref{eq:snae_lstdg_dae_eps}) of the local DG predictor, the following expression is obtained:
\begin{equation}
\sum\limits_{q = 0}^{N} \mathrm{A}_{pq}\mathbf{G}(\hat{\mathbf{q}}_{n, q}^{\epsilon}, \hat{\mathbf{r}}_{n, q}^{\epsilon}, t(\tau_{q})) =
\frac{\epsilon}{{\Delta t}_{n}}\left[\mathbf{v}_{n}^{\epsilon} - \hat{\mathbf{r}}_{n, p}^{\epsilon}\right] \underset{{\epsilon \rightarrow 0}}{\rightarrow} \mathbf{0},
\end{equation}
from which, and from the non-singularity of matrix $\mathrm{A}$, a new system of nonlinear algebraic equations follows:
\begin{equation}
\sum\limits_{q = 0}^{N} \mathrm{A}_{pq}\mathbf{G}(\hat{\mathbf{q}}_{n, q}, \hat{\mathbf{r}}_{n, q}, t(\tau_{q})) = \mathbf{0},\ \Rightarrow\ 
\mathbf{G}(\hat{\mathbf{q}}_{n, p}, \hat{\mathbf{r}}_{n, p}, t(\tau_{q})) = \mathbf{0},
\end{equation}
which contains the expansion coefficients $(\hat{\mathbf{q}}_{n, q}, \hat{\mathbf{r}}_{n, q})$ defined by expressions (\ref{eq:limit_eps_defs}). The final system of nonlinear algebraic equations of the local DG predictor for the DAE system (\ref{eq:ivp_dae_diff_src}), obtained in limit $\epsilon \neq 0$, takes the following form:
\begin{equation}\label{eq:snae_lstdg_dae}
\begin{split}
&\hat{\mathbf{q}}_{n, p} -
{\Delta t}_{n}\sum\limits_{q = 0}^{N} \mathrm{A}_{pq}\mathbf{F}(\hat{\mathbf{q}}_{n, q}, \hat{\mathbf{r}}_{n, q}, t(\tau_{q})) = \mathbf{u}_{n},\quad
\mathbf{G}(\hat{\mathbf{q}}_{n, p}, \hat{\mathbf{r}}_{n, p}, t(\tau_{q})) = \mathbf{0},
\end{split}
\end{equation}
and the expressions for obtaining the numerical solution $(\mathbf{u}_{n+1},\, \mathbf{v}_{n+1})$ at the next grid node $t_{n+1}$ take the following form:
\begin{equation}\label{eq:ader_dg_node_sol_dae}
\begin{split}
&\mathbf{u}_{n+1} = \mathbf{u}_{n} +
{\Delta t}_{n}\sum\limits_{p = 0}^{N} w_{p} \mathbf{F}\left(\hat{\mathbf{q}}_{n, p}, \hat{\mathbf{r}}_{n, p}, t(\tau_{p})\right),\\
&\mathbf{v}_{n+1} = \mathbf{r}_{n}(1) = \sum\limits_{p} \hat{\mathbf{r}}_{n, p}\varphi_{p}(1),
\end{split}
\end{equation}
which is completely identical to the formulation of the ADER-DG numerical method with a local DG predictor for the DAE system (\ref{eq:ivp_dae_diff_src}) in~\cite{ader_dg_dae}.

Thus, based on the use of $\epsilon$-embedding method for the DAE system, it is rigorously shown that the version of the ADER-DG numerical method with a local DG predictor for the DAE system (\ref{eq:ivp_dae_diff_src}) proposed in~\cite{ader_dg_dae} is the original numerical method for solving the problem initial value for ODE systems (\ref{eq:ivp_ode_diff_src})~\cite{ader_dg_ode_jsc}, which can be applied to solving DAE systems (\ref{eq:ivp_dae_diff_src}) in formulation (\ref{eq:snae_lstdg_dae}) and (\ref{eq:ader_dg_node_sol_dae}). Thus, this explains the results obtained in~\cite{ader_dg_dae}, which are quantitatively some similar to the results in~\cite{ader_dg_ode_jsc}, especially for DAE systems of low index.

It should be noted that the improved local solution proposed in the Subsection~\ref{sec:imp_dg:imp_loc_sol} cannot be created by analogy in the case of the DAE system (\ref{eq:ivp_dae_diff_src}). An improved local solution $\mathbf{q}_{n}^{\rm IL}(\tau)$ for differential variables can be created:
\begin{equation}\label{eq:imp_local_sol_dae}
\begin{split}
\mathbf{q}_{n}^{\rm IL}(\tau) &= \mathbf{u}_{n} + {\Delta t}_{n}\int\limits_{0}^{\tau} \mathbf{F}\left(\mathbf{q}_{n}(\xi), \mathbf{r}_{n}(\xi), t(\xi)\right) d\xi\\
&= \mathbf{u}_{n} + {\Delta t}_{n} \sum\limits_{p = 0}^{N} \mathbf{F}\left(\mathbf{q}_{n, p}, \mathbf{r}_{n, p}, t(\tau_{p})\right) \int\limits_{0}^{\tau} \varphi_{p}(\xi) d\xi.
\end{split}
\end{equation}
However, a similar improved local solution $\mathbf{r}_{n}^{\rm IL}(\tau)$ for algebraic variables $\mathbf{v}$ cannot be created without transforming the original DAE system (for example, by lowering the system index to 0), but this would be a different problem, not the DAE system (\ref{eq:ivp_dae_diff_src}).

\section*{Conclusion}
\addcontentsline{toc}{section}{Conclusion}
\label{concl}

In conclusion, it should be noted that in this paper, an improved local numerical solution for the ADER-DG numerical method with a local DG predictor for solving the initial value problem for a first-order ODE system (\ref{eq:ivp_ode_diff_src}) is proposed. The improved local numerical solution $\mathbf{u}_{\rm IL}$ demonstrates convergence orders $p = N+2$ of one higher than the convergence order $p = N+1$ of the local numerical solution $\mathbf{u}_{L}$ of the original ADER-DG numerical method. The improved local numerical solution has the property of continuity at grid nodes $t_{n}$: $\mathbf{u}_{\rm IL}(t_{n}^{-}) = \mathbf{u}_{\rm IL}(t_{n}^{+})$, while the local numerical solution $\mathbf{u}_{L}$ of the original ADER-DG numerical method was generally discontinuous on the right at grid nodes $t_{n}$: $\mathbf{u}_{L}(t_{n}^{-}) \neq \mathbf{u}_{L}(t_{n}^{+})$. Rigorous proofs of the approximation orders $N+1$ of the local numerical solution $\mathbf{u}_{L}$ of the original ADER-DG method and $N+2$ the improved local numerical solution $\mathbf{u}_{\rm IL}$ are presented. A rigorous derivation of the approximation order of the local numerical solution $\mathbf{u}_{L}$ of the original ADER-DG method is essentially based on the satisfiability of the Lipschitz condition (\ref{eq:lip_cond_f}) for the right-hand side function $\mathbf{F}$ of the ODE system (\ref{eq:ivp_ode_diff_src}). However, a rigorous derivation of the approximation order $N+2$ of the improved local numerical solution $\mathbf{u}_{\rm IL}$ is not related to such a strong requirement.

Obtaining the proposed improved local numerical solution $\mathbf{u}_{\rm IL}$ does not require significant changes to the structure of the ADER-DG numerical method with a local DG predictor. It is calculated based on the local numerical solution $\mathbf{u}_{L}$ and does not have an effect on the numerical solution $\mathbf{u}_{n}$ at the grid nodes. Therefore, all conclusions regarding the convergence orders of the numerical solution $\mathbf{u}_{n}$ at grid nodes, the resulting superconvergence $p = 2N+1$, and the very high stability of the ADER-DG numerical method with a local DG predictor, as established in~\cite{ader_dg_ode_jsc}, remain unchanged. Using the improved local numerical solution $\mathbf{u}_{\rm IL}$ will not require significant additional computational costs, other than calculating the solution values at the tabulation points.

A wide range of applications of the ADER-DG numerical method with a local DG predictor is presented for solving specific initial value problems for ODE systems (\ref{eq:ivp_ode_diff_src}). Computational results were obtained for a wide range of polynomial degrees $1 \leqslant N \leqslant 60$, covering both the range of relatively low polynomial degrees $N$ required in applied computing and the range of very high polynomial degrees $N$ suitable for high-precision calculations. The obtained results provide strong confirmation for the developed rigorous theory. The improved local numerical solution $\mathbf{u}_{\rm IL}$ is shown to exhibit both higher accuracy, compared to the original local numerical solution $\mathbf{u}_{L}$, and improved smoothness and point-wise comparability. Empirical convergence orders $p$ of all individual numerical solutions $\mathbf{u}_{L}$, $\mathbf{u}_{\rm IL}$, $\mathbf{u}_{n}$ were calculated for a wide range of error norms (\ref{eq:eps_un_global_def}), (\ref{eq:eps_ul_global_def}), which well agree with the expected convergence orders equal to the calculated approximation orders $p_{\rm th.}$ (\ref{eq:conv_ords_exp}). It was found that for linear ODE systems, the empirical convergence orders typically exactly match the expected convergence orders $p_{\rm th.}$, while for the numerical solution of nonlinear ODE systems, some small deviations from the expected values $p_{\rm th.}$ occur. A specific example was chosen that violated the Lipschitz condition (\ref{eq:lip_cond_f}) for the right-hand side function $\mathbf{F}$ of the ODE system (\ref{eq:ivp_ode_diff_src}) to clearly demonstrate the necessity of this requirement for applying the ADER-DG numerical method with a local DG predictor. It is shown that the obtained numerical solutions do not exhibit significant qualitative or quantitative deviations from cases where the Lipschitz condition (\ref{eq:lip_cond_f}) for the right-hand side function $\mathbf{F}$ of an ODE system is satisfied. \corrtext{The satisfiability of the energy conservation law in the numerical solution of conservative ODE systems is investigated. It is shown that despite the dissipative nature of the numerical method, a sufficiently high degree can be chosen so that the accuracy of the energy conservation law is at or below the characteristic error of representing real numbers as floating-point numbers.}

The final section of the paper presents a rigorous proof, based on the $\epsilon$-embedding method, of the applicability of the ADER-DG numerical method with a local DG predictor to solving DAE systems (\ref{eq:ivp_dae_diff_src}). It is shown that the proposed version~\cite{ader_dg_dae} of the method can be classified as a classical version of the ADER-DG method with a local DG predictor for solving the initial value problem for a first-order ODE system~\cite{ader_dg_ode_jsc}.

\bmhead{Acknowledgments}
The author would like to thank the anonymous reviewers for their comments and remarks that helped to improve the quality of this paper.
The author would like to thank A.P.~Popova for help in correcting the English text.

\subsection*{Data Availability Statement}
\addcontentsline{toc}{subsection}{Data Availability Statement}
Data will be made available on reasonable request.


\section*{Declarations}
\addcontentsline{toc}{section}{Declarations}

\subsection*{Conflict of interest}
\addcontentsline{toc}{subsection}{Conflict of interest}
The author declares that he has no conflict of interest.

\addcontentsline{toc}{section}{References}


\end{document}